\documentclass{memo-l}

\usepackage{amssymb}
\usepackage{amstext}
\usepackage{amsbsy}
\usepackage{amsxtra}
\usepackage{amscd}
\usepackage{latexsym}
\usepackage{graphicx}

\usepackage[pdfusetitle]{hyperref}
\hypersetup{
  colorlinks,
  allcolors=blue,
  linktoc=all
}

\usepackage{tikz}
\usetikzlibrary{cd}

\usepackage{scalerel}

\usepackage[mathscr]{eucal}

\usepackage{slashed}

\usepackage{ifthen}

\usepackage{bbm}

\usepackage{exscale}

\usepackage[shortcuts]{extdash}

\usepackage{enumitem}
\setlist[enumerate,1]{label={\upshape(\roman*)}}

\usepackage{titlesec}

\usepackage{nonumonpart}

\usepackage{tabularx}

\counterwithout*{chapter}{part} 

\titleformat
{\part} 
[display] 
{\bfseries\Huge} 
{} 
{0ex} 
{
    \rule{\textwidth}{1pt}
    \centering
} 
[
\vspace{-1.5ex}
\rule{\textwidth}{1pt}
] 

\let\oldcontentsline\contentsline
\newcommand{\nopagecontentsline}[3]{\oldcontentsline{#1}{#2}{}}
\newcommand{\TOCstoppagenumber}{%
\addtocontents{toc}{\let\protect\contentsline\protect\nopagecontentsline}%
}
\newcommand{\TOCunstoppagenumber}{%
\addtocontents{toc}{\let\protect\contentsline\protect\oldcontentsline}%
}

\newcounter{savepagenumber}
\newcommand{\partpage}[1]{%
\newpage%
\setcounter{savepagenumber}{\value{page}}
\pagenumbering{gobble}%
\TOCstoppagenumber%
\part*{#1}%
\TOCunstoppagenumber%
\pagenumbering{arabic}%
\setcounter{page}{\value{savepagenumber}}
}

\theoremstyle{plain}
\newtheorem{theorem}{Theorem}[chapter]
\newtheorem{corollary}[theorem]{Corollary}
\newtheorem{lemma}[theorem]{Lemma}
\newtheorem{proposition}[theorem]{Proposition}

\theoremstyle{definition}
\newtheorem{example}[theorem]{Example}
\newtheorem{definition}[theorem]{Definition}

\AtBeginEnvironment{example}{%
  \pushQED{\qed}%
}
\AtEndEnvironment{example}{\popQED\endexample}

\theoremstyle{remark}
\newtheorem{remark}[theorem]{Remark}
\numberwithin{equation}{chapter}

\newtheorem*{genericthm*}{\thistheoremname}
\newcommand{\thistheoremname}{The}
\newcounter{genericthm}

\newenvironment{namedtheorem*}[1]
  {\renewcommand{\thistheoremname}{#1}%
  ~\refstepcounter{genericthm}%
   \begin{genericthm*}}
  {\end{genericthm*}}

\newcommand{\proofstep}[1]{\noindent\textsc{#1}}

\newlist{enumeratethm}{enumerate}{1}
\setlist[enumeratethm]{label=(\roman*), leftmargin=\parindent}

\newcommand{\Op}[1]{\mathsf{#1}} 
\DeclareMathAlphabet{\Set}{U}{eur}{m}{n} 
\newcommand{\Spa}[1]{\mathrm{#1}} 
\newcommand{\Sh}[1]{\mathcal{#1}} 
\newcommand{\Func}[1]{\mathcal{#1}} 
\newcommand{\De}[1]{\emph{#1}} 

\newcommand{\qat}{\mathbb{H}}
\newcommand{\X}{\Spa{X}}
\newcommand{\U}{U}
\newcommand{\V}{V}
\newcommand{\W}{W}
\newcommand{\C}{C}
\newcommand{\B}{B}
\newcommand{\SO}{\Set{O}}
\newcommand{\SU}{\Set{U}}
\newcommand{\SV}{\Set{V}}

\newcommand{\con}{\Spa{C}}
\newcommand{\forms}[2]{\mathcal{A}^{#1}_{#2}}
\newcommand{\pot}[1]{\Spa{H}_{#1}}
\newcommand{\Q}[1]{\mathcal{Q}_{#1}}
\newcommand{\tQ}[1]{\Tilde{\mathcal{Q}}_{#1}}
\newcommand{\ban}[2]{\mathcal{L}^{#1}_{#2}}
\newcommand{\banach}[1]{L^{\!#1}}
\newcommand{\sob}[2]{\mathcal{W}^{#1}_{#2}}
\newcommand{\sobolev}[1]{W^{#1}}
\newcommand{\kodaira}{\mathrm{Kod}}
\newcommand{\leftsp}{\mathrm{Left}}
\newcommand{\prewei}{\mathrm{PreWei}}
\newcommand{\wei}{\mathrm{Wei}}

\newcommand{\loc}{_{\text{\scriptsize\textrm{loc}}}}
\newcommand{\reg}{\text{\scriptsize\textrm{reg}}}
\newcommand{\sing}{\text{\scriptsize\textrm{sing}}}
\newcommand{\ri}{\text{\scriptsize\textrm{right}}}
\newcommand{\pa}{^\ast}
\newcommand{\sd}{^\dagger}

\newcommand{\dirac}{\slashed{\partial}}
\newcommand{\diracJ}{\Op{J}}
\newcommand{\barpartial}{\Bar{\partial}}
\newcommand{\boundary}{\partial}
\newcommand{\hatast}{\mathbin{\hat{\ast}}}
\newcommand{\del}[1]{%
\ifthenelse{\equal{#1}{}}%
{\partial}%
{\partial_{\scaleto{#1\mathstrut}{6pt}}}%
}
\newcommand{\delbar}[1]{%
\ifthenelse{\equal{#1}{}}%
{\barpartial}%
{\barpartial_{\scaleto{#1\mathstrut}{6pt}}}%
}
\newcommand{\unity}{\mathbbm{1}}
\newcommand{\Two}{\mathrm{I\!I}}
\newcommand{\var}{\delta}
\DeclareMathOperator{\Real}{Re}
\DeclareMathOperator{\Imag}{Im}
\newcommand{\dmu}{d\mu}
\renewcommand{\:}{\,{:}\,}

\DeclareMathOperator{\diam}{diam}
\DeclareMathOperator{\Hom}{Hom}
\DeclareMathOperator{\End}{End}
\DeclareMathOperator{\tr}{tr}
\DeclareMathOperator*{\res}{Res}
\DeclareMathOperator*{\ord}{ord}

\DeclareMathOperator{\id}{id}
\DeclareMathOperator{\per}{Per}
\newcommand{\willmore}{\mathscr{W}}
\newcommand{\fundsol}{\Func{F}}

\newcommand{\genus}{\mathrm{g}}

\newcommand{\moebius}{g}
\newcommand{\ci}{\mathrm{i}}
\newcommand{\qi}{\mathrm{i}}
\newcommand{\qj}{\mathrm{j}}
\newcommand{\qk}{\mathrm{k}}

\newcommand{\lz}{(}
\newcommand{\rz}{)^{(0,0)}}
\newcommand{\lp}{(}
\newcommand{\rp}{)^{(1,0)}}
\newcommand{\lh}{(}
\newcommand{\rh}{)^{(1,1)}}
\newcommand{\ls}{(}
\newcommand{\rs}{)^{(0,1)}}

\newcommand{\labelthis}[1]{\stepcounter{equation}\tag{\theequation}\label{#1}}

\makeindex

\begin{document}
\pagenumbering{gobble}

\title{Quaternionic Analysis of Conformal Maps and the Willmore Functional}

\author[R.~Ogilvie]{Ross Ogilvie}
\address{Mathematics Chair III\\Universit\"at Mannheim\\D-68131 Mannheim, Germany}
\email{r.ogilvie@math.uni-mannheim.de}

\author[M.U.~Schmidt]{Martin Ulrich Schmidt}
\address{Mathematics Chair III\\Universit\"at Mannheim\\D-68131 Mannheim, Germany}
\email{schmidt@math.uni-mannheim.de}

\date{\today}




\begin{abstract}
Quaternionic analysis, which describes conformal maps from Riemann surfaces into $\mathbb{R}^3$ or $\mathbb{R}^4$, is extended to weakly conformal maps.
As a consequence we present a new proof that on any compact Riemann surface $\X$ the Willmore functional, the integral of the square of the mean curvature, attains a minimum on the space of smooth conformal maps from $\X$ to $\mathbb{R}^3$ or $\mathbb{R}^4$.
This was first proven by Kuwert and Schätzle under the assumption that the infimum of the Willmore functional is less than $8\pi$.
In this case all conformal maps are unbranched, due to an estimate of Li and Yau.
Rivière removed this restriction by allowing as limits conformal maps with ramification points.
Our approach admits these weakly conformal maps from the very beginning, by extending the quaternionic function theory as developed by Pedit and Pinkall to square-integrable potentials.

In Part~I we carry over the most important properties of holomorphic functions to quaternionic analysis with square-integrable potentials.
We develop a differential operator, the Darboux transformation, that generalizes the 
$\partial$\=/derivative and with respect to which these holomorphic functions are infinitely many times differentiable.
It transforms the Kodaira representation of a weakly conformal map into its Weierstraß representation.

Parts~II and~III treat respectively the existence of a minima on the space of weakly conformal maps with ramification points and the regularity of this minima.
The existence result is a variant of Montel's theorem for quaternionic functions.
Its proof combines the concentration compactness principle with a characterization of possible singularities, namely poles whose orders are bounded in terms of the concentrated energy.
For regularity, we use that an analogue of the Euler-Lagrange equation is fulfilled, which implies smoothness by elliptic bootstrapping.
The difficulty is that the space of weakly conformal mappings with possible ramification points may have singularities at strongly isothermic mappings.
We show that these singularities are relatively mild: the space is the zero set of a smooth real function on a Banach manifold whose Hessian is indefinite.
In particular, the tangent cone spans the tangent space.
\end{abstract}

\maketitle

\hypersetup{allcolors=black}
\tableofcontents
\hypersetup{allcolors=blue}

\frontmatter

\TOCstoppagenumber
\chapter*{Introduction}
\label{chapter:introduction}
\TOCunstoppagenumber

Quaternionic analysis (also called quaternionic function theory) is a powerful tool in the study of conformal maps from a Riemann surface into $4$\=/dimensional Euclidean space.
Just as complex analysis is indispensable in the study of conformal maps between surfaces, 
by identifying maps into $\mathbb{R}^4$ with quaternion-valued functions a similarly rich theory emerges.
Many fundamental results have natural generalizations, such as the notion of holomorphicity, the Riemann-Roch theorem, and the Plücker formula.
The Willmore functional, the integral of the square of the mean curvature of an immersion, is an important conformal invariant and finds a natural home in this framework.
Surfaces in $\mathbb{R}^3$ are obtained as a special case, with $\mathbb{R}^3$ identified as the imaginary quaternions.
The foundations of this theory are due to the work of the ``Berlin group'', principally Franz Pedit and Ulrich Pinkall~\cite{PP,BFLPP,FLPP}. 
The goal of this book is to extend quaternionic analysis to weakly conformal maps in order to prove the existence of minima of the Willmore functional on the space of conformal immersions from a fixed Riemann surface.

Let us sketch the concepts of quaternionic analysis in the case of smooth immersions and provide an example.
A complete treatment is given in Chapter~\ref{chapter:prelim}.
A conformal map is one that preserves angles.
The angle between vectors in the quaternions $\qat = \mathbb{R}^4$ is given by the usual Euclidean inner product.
If we want to consider conformal maps on surfaces, it is natural to work with Riemann surfaces.
This is because biholomorphic functions are conformal and so the natural conformal structure on open subsets of $\mathbb{C}$ is inherited by a Riemann surface.
The complex structure of a Riemann surface can be equivalently encoded by the Hodge star $\ast$ on $1$\=/forms.
Let us consider the conformal map $F$ locally so that we may take the Riemann surface to be $\Omega \subset \mathbb{C}$.
If $dF$ has rank two at some point, then the tangent plane of the image can be described by two quaternions called the left normal $N$ and the right normal $R$:
\[
dF(T_p\X) = \{ x \in \mathbb{H} \mid Nx = xR \}
\]
for $N^2 = R^2 = -1$, where the sign of $N$ is chosen such that $N$ acts on the tangent plane by a $\pi/2$ rotation compatible with induced orientation.
Therefore $N$ encodes a conformal structure for the tangent plane and $F$ is conformal if and only if
\begin{align*}
-\ast dF &= N dF
&
N^2&=-1.
\end{align*}
Equivalently for the right normal
\begin{align*}
-\ast dF &= dFR
&
R^2&=-1.
\end{align*}
We can take this to be the definition of conformality.
The tangent plane lies in the imaginary quaternions exactly when $R = -N$.
Conversely, if all its tangent planes have $R = -N$ then $F$ lies in the imaginary quaternions up to a real additive constant.

\begin{example}[Catenoid]
\label{eg:catenoid}
\index{Catenoid!Admissibility}
Consider $F : \mathbb{C}^\ast \to \qat$ given by
\[
F(z) = \ln(z\Bar{z}) \qi + (\Bar{z} + z^{-1})\qj.
\]
It is easier to recognize this surface in the coordinates $z = \exp (u+\qi v)$, where we have
\[
F(z) 
= 2u \qi + (e^u + e^{-u})e^{-\qi v}\qj
= 2u \qi + 2\cosh u (\cos v\, \qj - \sin v\, \qk).
\]
This is a catenoid in $\Imag \qat$, the surface of revolution generated by the catenary.
Its conformality can be checked by direct calculation, with 
$\left|\frac{\partial F}{\partial u}\right|^2 = 4 + 4\sinh^2u$ and $\left|\frac{\partial F}{\partial u}\right|^2 = 4\cosh^2u$ plainly equal, and the two coordinate tangent vectors orthogonal.
Alternatively, as we calculate in detail in Example~\ref{eg:catenoid derivative}, we can factorize the derivative of $F$ into a special form
\begin{equation}
\label{eq:catenoid dF}
dF
= (1-z^{-1}\qk) \qj dz (1+z^{-1}\qk).
\end{equation}
Now, applying Hodge star gives
\[
-\ast dF
= - (1-z^{-1}\qk) \qj (-\qi dz) (1+z^{-1}\qk).
\]
By comparing the expressions we can read off the only possibilities for $N$ and $R$
\[
N = (1-z^{-1}\qk) (-\qi) (1-z^{-1}\qk)^{-1}, \qquad
R = (1+z^{-1}\qk)^{-1} \qi (1+z^{-1}\qk).
\]
It is easy to check that $N^2 = R^2 = -1$, so $F$ satisfies Equations~\eqref{eq:left normal} and~\eqref{eq:right normal}. Further $R = -N$, which we shall see is characteristic of maps into $\mathbb{R}^3$.
\end{example}
In complex analysis, multiplication by $\qi$ is used to split $d$ by type:
\begin{align}\label{eq partial}
d'&=\tfrac{1}{2}(d+\qi\ast d)&d''&=\tfrac{1}{2}(d-\qi\ast d).
\end{align}
The operator $d''$ is called the holomorphic structure, since holomorphic functions are defined by the condition $d''f = 0$.
In the present situation, we have a complex structure on the trivial right $\qat$\=/line bundle $\X \times \qat$
given by the left multiplication with $N$. 
This action induces an analogous splitting of $d$ acting on $\qat$\=/valued functions:
\begin{align}
\label{eq:holomorphic structure}
d'_N&=\tfrac{1}{2}(d+N\ast d)&d''_N&=\tfrac{1}{2}(d-N\ast d).
\end{align}
The functions $F: \Omega \to \qat$ that lie in the kernel of $d''_N$ are exactly the conformal maps with left normal $N$.
For general Riemann surfaces $\X$, not just subsets $\Omega \subset \mathbb{C}$, we should generalize from $\qat$\=/valued functions to sections of $\qat$\=/line bundles.
Motivated by the properties of $d''_N$, we define (see~\eqref{eq:holo product rule}) a holomorphic structure $D$ to be a (right) $\qat$\=/linear differential operator that satisfies
\[
D(\xi\alpha)=D(\xi)\alpha+\tfrac{1}{2}(\xi d\alpha-N\xi\ast d\alpha)\quad\text{for all }\alpha\in C^\infty(\Omega,\qat).
\]
This is analogous to the definition of a holomorphic structure on a $\mathbb{C}$\=/line bundle. 
In fact, the interplay of the holomorphic structure with the complex structure gives rise to a splitting $D=D^++D^-$ into a commuting and an anti-commuting part of $D$:
\begin{align*}
D^+(\xi)&=\tfrac{1}{2}\big(D(\xi)-ND(N\xi)\big),&D^-(\xi)&=\tfrac{1}{2}\big(D(\xi)+ND(N\xi)\big).
\end{align*}
The left $\mathbb{C}$\=/action and the right $\qat$\=/action identifies a $\mathbb{C}$\=/subbundle $E$: The subspace where left action coincides with right action by $\qi \in \qat$.
The commuting part $D^+$ endows $E$ with a holomorphic structure. 
Conversely, given a holomorphic structure $\delbar{E}$ on $E$, there is an extension to the $\qat$\=/line bundle.
Define another real two-dimensional subbundle $E'$ as the subbundle on which $N$ from the left acts as $-\qi \in \qat$ from the right.
The $\qat$\=/line bundle is the direct sum $E \oplus E'$.
The Berlin group makes $E'$ into a $\mathbb{C}$\=/subbundle by the restriction of the left action (by $N$), so that $E \cong E'$.
The extension is then the sum of two copies of $D^+$.
In this book, we give $E'$ the complex structure of right multiplication by $\qi \in \qat$.
The advantage is that the extension is just the $\qat$\=/linear extension of $\delbar{E}$, simplifying calculations in local coordinates.
Both constructions give the same holomorphic structure $D^+$ on the $\qat$\=/line bundle, which is characterized as the unique holomorphic structure that restricts to $E$ to give back $\delbar{E}$ and preserves the splitting $E\oplus E'$.
To sum up, the structure of conformality is encoded in the notion of a \emph{holomorphic $\qat$\=/line bundle}:

\begin{itemize}
\item[(i)] a $\qat$\=/line bundle with right $\qat$\=/action, 
\item[(ii)] a $\qat$\=/linear complex structure $N$ and 
\item[(iii)] a holomorphic structure $D=D^++D^-$.
\end{itemize}

Before we switch to the non-smooth case let us describe one further structure that we want to extend to our class of weakly conformal maps.
Namely, we want that the Willmore functional extends to this class. 
Classically, it is a functional on the space of immersions $F:\X\to\mathbb{R}^n$ from a compact two-dimensional manifold into the $n$\=/dimensional Euclidean space, defined as the integral over the square of the mean curvature with respect to induced measure $\dmu$:
\[
\willmore(F)=\int_\X H^2\dmu.
\]
\cite{Willmore1965} introduced this quantity with the idea that its minimum could potentially be an interesting invariant on immersed closed surfaces.
Its importance is related to its large symmetry group, which contains all conformal transformations of the ambient space~\cite{Pinkall2024}. Fortunately, this Willmore energy has a simple formula in terms of the anti-commuting part $D^-$ of the holomorphic structure in (iii), which interchanges $E$ and $E'$ and is tensorial:
\[
\tfrac{1}{2}\big(d''_N f - N d''_N(Nf)\big)=\tfrac{1}{4}\big(df-N\ast df+Nd(Nf)+\ast d(Nf)\big)=\tfrac{1}{4}(NdN+\ast dN)f.
\]
In our setup, it is $\mathbb{C}$\=/linear with respect to the complex structures on $E$ and $E'$.
Pedit and Pinkall call the tensor $\frac{1}{4}(NdN+\ast dN)$ the Hopf field. They show that, up to a topological constant, four times the square of its $\banach{2}$\=/norm is equal to the Willmore energy~\cite{PP} (compare also~\cite{Ta1} and Chapter~\ref{chapter:isothermic}).
For admissible maps we shall assume that the density of the Willmore energy is locally integrable. Since $N^2 = -1$ implies $|N|=1$, this is guaranteed by $N\in\sobolev{1,2}\loc(\Omega,\qat)$. Due to the invariance of the Willmore functional with respect to the scaling transformations, it is no accident that the analysis of these admissible maps will turn out to depend on the Sobolev embedding for the critical exponent $2$; this is the scaling invariant case and results in an additional analytical challenge.

Now we are ready to describe our approach to weakly conformal maps. The aim is to include in our considerations all maps that give rise to equivalent forms of the three structures (i)-(iii) above and have a well defined Willmore energy. We will see in Chapter~\ref{chapter:kodaira} that with  modified structures (ii) and (iii) this is true for the class of \emph{admissible} maps defined in below. After this definition we shall explain these modifications of the structures (ii) and (iii).

\begin{definition}
\label{def:weakly conformal local}
A map $F\in\sobolev{1,1}\loc(\Omega,\qat)$ on an open subset $\Omega\subset\mathbb{C}$ is called \emph{weakly conformal} if $dF$ has almost everywhere rank two and there exists a measurable map $N:\Omega\to\qat$ which is almost everywhere a left normal~\eqref{eq:left normal}.
If in addition $N\in\sobolev{1,2}\loc(\Spa{\Omega},\qat)$ then the map $F$ is called \emph{admissible}.
\end{definition}
Admissible maps do not fit into the framework of smooth holomorphic $\qat$\=/line bundles described above: The left normal $N\in\sobolev{1,2}\loc(\Omega,\qat)$ is a measurable function and does not endow all fibers of the $\qat$\=/line bundle with a complex structure.
Consequently, the subset $E = \{ Nx = x\qi \}$ is neither a $\mathbb{C}$\=/line bundle nor a holomorphic $\mathbb{C}$\=/line bundle. 
However, we shall prove in Theorem~\ref{thm:kodaira normal} that admissible maps fit perfectly into the framework of quaternionic analysis, if we replace the non-trivial complex structure given by left multiplication with a left normal $N$ by the trivial complex structure given by left multiplication with the quaternion $\qi$, so that the $\mathbb{C}$\=/subbundle $E$ is well-defined, and replace the holomorphic structure $d''_N$~\eqref{eq:holomorphic structure} induced by the trivial flat connection $d$ by another holomorphic structure $D$. 
The triviality of the left $\mathbb{C}$\=/action allows us to recover the $\qat$\=/bundle as $E_\qat=E\otimes_{\mathbb{C}}\qat$ with the natural left $\mathbb{C}$\=/action. 
Theorem~\ref{thm:kodaira normal} will prove that the holomorphic structure $D$ has the form $D=\delbar{E}-\V$ for the unique extension $\delbar{E}$ to $E_\qat$ discussed above, and with a tensorial part denoted by $D^-=-\V$, which we call a potential. 
Since the Willmore energy is four times the square of the $\banach{2}$\=/norm of this potential, one major subject of this book is the extension of the quaternionic function theory of Pedit and Pinkall to the setting of $\banach{2}\loc$\=/potentials.

\begin{remark}\label{re:weakly conformal}
Let us relate the class of admissible maps to another class of weakly conformal maps, which has been used by Kuwert and Schätzle~\cite{KL, KS} and by Riv\`{e}re~\cite{Ri,Ri2} in their investigation of the Willmore functional.
In~\cite[Section 2]{KL} the class of conformal $\sobolev{2,2}\loc$\=/maps is defined as maps $F\in\sobolev{2,2}\loc(\Omega,\mathbb{R}^n)$ with induced metric $g_{ij}=e^{2u}\delta_{ij}$ with $u\in\banach{\infty}\loc(\Omega)$. As shown in~\cite[Section 2]{KL}, these conditions imply $u\in\sobolev{1,2}\loc(\Omega)$.
By using Wente's inequality we show in Chapter~\ref{chapter:kodaira} that admissible maps in fact belong to $\sobolev{2,2}\loc(\Omega,\qat)\cap\sobolev{1,\infty}\loc(\Omega,\qat)$. For $n=4$ we identify $\mathbb{R}^4\cong\qat$.
Then the orthogonal first derivatives $\frac{\partial F}{\partial x}$ and $\frac{\partial F}{\partial y}$ belong to $\sobolev{1,2}\loc(\Omega,\qat)$ and have length $e^{2u}$.
Therefore the left and right normals,
\begin{align*}
N&=\tfrac{\partial F}{\partial y} \left(\tfrac{\partial F}{\partial x}\right)^{-1},&R&=\left(\tfrac{\partial F}{\partial x}\right)^{-1} \tfrac{\partial F}{\partial y},
\end{align*}
both belong to $\sobolev{1,2}\loc(\Omega,\qat)$.
This shows that these maps are contained in the class of admissible maps.
We shall see in Remark~\ref{rem:weierstrass} that this class of conformal $\sobolev{2,2}\loc$\=/maps is equal to the subclass of admissible maps without roots of $dF$.
For $n=3$ the conformal $\sobolev{2,2}\loc$\=/maps are admissible maps with $R=-N$ and without roots of $dF$.
\end{remark}

The extension of the analysis to weakly conformal maps with ramifications points, allows one to show that admissible maps with bounded Willmore energy have convergent subsequences (compare~\cite{Ri2}).
We shall prove in Theorem~\ref{compactness of conformal mappings}: given a compact Riemann surface any sequence of admissible maps with bounded Willmore energy has a subsequence which converges in an appropriate sense to an admissible map. 
Due to a result of Li and Yau, which we extend to all admissible maps, the bound $\willmore(F)<8\pi$ excludes ramification points of $F$. 
Our analysis shows that the possible singularities of the limit can develop singularities only at points where the Willmore energy concentrates with at least the energy of a round sphere.

The extension to weakly conformal maps with possible roots of the first derivative is another analytical challenge, since the weakly conformal maps are represented by weakly differentiable functions (elements of a Sobolev space), which do not have an obvious notion of roots. 
However, in complex analysis roots are completely described by comparing their growth to the roots of polynomials; this is the definition of the order of a root.
The theory builds out powerful tools for describing roots.
The roots of a complex holomorphic function $f$ are either isolated or $f$ vanishes identically: the strong unique continuation property. 
As a consequence they define a root divisor, or more generally, meromorphic functions define a root and pole divisor. 
Furthermore, meromorphic forms $\omega$ have a residue, whose vanishing is locally equivalent to the existence of a function $f$ with $df=\omega$. 
We shall see that all these tools have a counterpart in quaternionic analysis. 

Let us now state the main theorem. 
The book is organized as a more or less straight path to a complete proof of this theorem:
\begin{namedtheorem*}{Main Theorem}
\label{thm:main theorem}
Let $\X$ be a compact Riemann surface
The Willmore functional attains a minimum on the space of smooth branched conformal maps $F:\X\to\mathbb{R}^4$.
The Willmore functional also attains a minimum on the space of smooth branched conformal maps $F:\X\to\mathbb{R}^3$ with only even order roots of $dF$.
Furthermore, if the Willmore functional of this minimizer in $\mathbb{R}^3$ is smaller than $12\pi$ then it is an immersion.
Finally, in both cases, if the Willmore functional of the minimizer is smaller than $8\pi$ then the minimizer is an immersed embedding.
\end{namedtheorem*}
For compact Riemann surfaces $\X$, on which the infimum of the Willmore functional is smaller than $8\pi$ the statement in the theorem is proven in~\cite{KS}. 
Rivière gave in~\cite{Ri2} the first proof of this statement without the bound on the Willmore functional in all codimensions. Let us mention some earlier results in this direction. The existence of a minimizer on the space of all immersions from a Riemann surface of prescribed genus into the $n$\=/dimensional Euclidean spaces was proven by Simon for genus one~\cite{Si1,Si2}, and by Bauer and Kuwert for all finite genera~\cite{BK}.

The proof is divided between the three parts of the book. 
Part~I provides the necessary tools from quaternionic analysis. 
Part~II proves a version of Montel's theorem in quaternionic analysis: On a compact Riemann surface a sequence of admissible maps with bounded Willmore energy has a subsequence which converges in an appropriate sense to an admissible map. 
Finally, in Part~III, we combine three tools: We first derive a substitute of the Euler-Lagrange equation. Secondly we show that these equations also hold at the singularities of the space of admissible maps. Lastly we apply an elliptic bootstrapping argument in order to prove smoothness and analyticity of the minimizers. 
Let us now explain in detail which statements fit together to a complete proof of this theorem:
\begin{proof}[Proof of the~\ref{thm:main theorem}]
Choose any sequence $(F_m)_{m\in\mathbb{N}}$ of admissible maps $F_m:\X\to\mathbb{R}^n$, for $n\in\{3,4\}$, whose Willmore energies $(\willmore(F_m))_{m\in\mathbb{N}}$ converge to the infimum of the set of Willmore energies of all such admissible maps.
Theorem~\ref{compactness of conformal mappings} shows that $(F_m)_{m\in\mathbb{N}}$ has a convergent subsequence whose limit $F:\X\to\mathbb{R}^n$ is an admissible minimizer.  Furthermore, if $n=3$ and all $dF_m$ have only roots of even order, then this theorem also shows that $dF$ has only roots of even order.
Due to the Theorems~\ref{constrained Willmore 1} and~\ref{constrained Willmore 2} the minimizers are smooth.
Moreover, due to Corollary~\ref{thm:embedding}, admissible maps with Willmore energy less than $8\pi$ are immersed embeddings.
Finally, due to the same corollary, admissible maps with Willmore energy less than $12\pi$ and whose derivatives $dF$ have only roots of even order, are smooth immersions.
This completes the proof of the~\ref{thm:main theorem}.
\end{proof}
We finish this introduction with a detailed summary of the three parts and their chapters. 
Because the precise definitions are yet to be introduced, we take the liberty in this summary to omit the $\banach{2}\loc$\=/potential when we speak of holomorphic sections of $\qat$\=/line bundles.

Part~I deals mostly with the transfer of the concepts and statements of quaternionic analysis from the smooth setting into the setting of $\banach{2}\loc$\=/potentials. 
Chapter~\ref{chapter:prelim} sets many of the conventions used throughout the book, for example the notion of a complex $\qat$\=/bundle, which is a vector bundle with a right $\qat$\=/action and a left $\mathbb{C}$\=/action.
It is for readers who are familiar with differential geometry but not necessarily quaternionic analysis. 
After a short introduction of the quaternions and their linear algebra, we define the exterior derivative and wedge product of differential forms with values in line- and $\qat$\=/vector bundles.
Of particular importance is the notion of a connection on a $\qat$\=/vector bundle. 
Finally, this chapter also introduces weakly differentiable sections on a Riemann surface independent of any further choices (Riemannian metric, cover, etc).

The two subsequent chapters establish many local concepts and properties of quaternionic holomorphic functions, which are well known for complex functions. In Chapter~\ref{chapter:local} these are Cauchy's Integral formula (Theorem~\ref{cauchy formula}), order of roots and poles (Definition~\ref{order of roots} and Lemma~\ref{lem:removable singularity}), strong unique continuation property (Theorem~\ref{strong unique continuation}), isolatedness of roots and one-dimensionality of quotients of holomorphic functions divided by holomorphic functions with a root at some given point (Lemma~\ref{quotient dimension}).

Chapter~\ref{chapter:darboux} is chiefly concerned with Theorem~\ref{thm:darboux}, which shows locally roughly speaking two statements. 
Firstly, the inverse of a quaternionic holomorphic function without roots is also a quaternionic holomorphic function. 
Secondly, the derivative of the quotient of two quaternionic holomorphic functions is a product of the inverse of the denominator quaternionic function with another quaternionic holomorphic function, which is some sense the derivative of the numerator quaternionic function. 
This theorem clarifies locally the relationship between the two representations of admissible maps which are the subject of the two subsequent chapters, namely the Kodaira and the Weierstraß representations. 
This theorem implies that any quaternionic holomorphic function has in some sense infinitely many holomorphic derivatives. 
The proof for $\banach{2}\loc$\=/potentials is much more involved than for potentials in $\bigcup_{2<q<\infty}\banach{q}\loc$.

In Chapter~\ref{chapter:kodaira} we first prove, that on any Riemann surface a left normal $N\in\sobolev{1,2}\loc(\X,\qat)$ gives rise to a holomorphic $\mathbb{C}$\=/line bundle $E$, a $\banach{2}$\=/potential $\V$ acting on the $\qat$\=/line bundle $E_\qat=E\otimes_\mathbb{C}\qat$ which anti-commutes with the obvious left multiplication by $\qi$ on $E_\qat$ and a section of $\upsilon$ of $E_\qat$ in the kernel of $\barpartial-\V$ without roots, such that $N=\upsilon^{-1}\qi\upsilon$. This implies that for any admissible map $F$ on $\X$ with left normal $N$ the section $\upsilon F$ of $E_\qat$ belongs to the kernel of $\barpartial-\V$:
\[
(\barpartial-\V)(\upsilon F)=\big(\tfrac{1}{2}(d-\qi\ast d)-\V\big)(\psi F)=\big((\barpartial-\V)\upsilon\big)F+\psi\tfrac{1}{2}(d-N\ast d)F=0.
\]
In turn this implies that $F \mapsto \upsilon F$ induces an isomorphism from the admissible maps on $\X$ with left normal $N$ onto the space of global sections of $E_\qat$ in the kernel of $\barpartial-\V$ which are not in the one-dimensional subspace $\upsilon\qat$. In particular, any admissible map on a Riemann surface can be represented as the quotient $F=\upsilon^{-1}\phi$ of two holomorphic sections of a holomorphic $\qat$\=/line bundle. 
We call it the Kodaira representation since it generalizes Kodaira's construction of holomorphic maps to projective spaces in terms of linear systems of holomorphic sections of a holomorphic line bundle~\cite[Chapter I, Section 4]{GrHa}. 
We also improve in Chapter~\ref{chapter:kodaira} the regularity of holomorphic sections $\xi$ of $\qat$\=/line bundles, by applying Wente's inequality: they belong to $\sobolev{1,2}\loc\cap\banach{\infty}\loc$.

Chapter~\ref{chapter:weierstrass} now extends the local transformation in Chapter~\ref{chapter:darboux} to a global transformation. In particular the Kodaira representation $F=\upsilon^{-1}\phi$ is transformed to a representation of $dF$ as product of two holomorphic sections of two holomorphic $\qat$\=/line bundles, which are related to each other. 
Pedit and Pinkall called such holomorphic $\qat$\=/line bundles paired. 
We denote the pairing of two holomorphic sections $\chi$ and $\psi$ by $dF=\lp\chi,\psi\rp$. 
This representation of the derivative of an admissible maps is called Weierstraß representation~\cite{Kon2,Ta1,Ta2,Fr2} (a local version was already known by Eisenhardt before the Dirac operator was invented~\cite{Ei}). 
We call this transformation from the Weierstraß data to the Kodaira data the Darboux transformation, because it represents the Darboux transformation for a particular integrable system, namely the Davey-Stewartson equation~\cite{Kon1}. This terminology denotes a general class of transformations associated to integrable systems with Lax operators. For any such integrable system, the corresponding Darboux transformation is bespoke. Generally given a solution of the non-linear differential equation of such an integrable system, the corresponding Darboux transformation produces another solution of the non-linear differential equation from a solution of the associated linear Lax equation (compare~\cite{EK,GHZ}). In geometric applications of integrable systems these transformations preserve some geometric properties of the corresponding geometric objects. In the present situation this point of view is presented in~\cite{BLPP}.

The two following chapters explore the applications of these representations to geometrical problems.
In Chapter~\ref{chapter:3-space} the behavior of the representations under Möbius transformation of the surface $F$ are investigated.
Admissible maps whose derivatives lie in $\Imag \qat$ have a certain symmetry in their Weierstraß representations.
This is used to give a criterion on the Weierstraß representation as to when the admissible map is a map to $\mathbb{R}^3 = \Imag \qat$.
In Chapter~\ref{chapter:isothermic} the differential geometry of admissible maps as a submanifolds is investigated.
We show how to compute common quantities, such as the second fundamental form, mean curvature, and Willmore energy, in terms of the Kodaira and Weierstraß data.
This knowledge is then used to extend the notion of isothermic (and strongly isothermic) to admissible maps.
Common results, such as the characterization of isothermic maps in terms of the existence of a dual map~\cite{Pa} and that totally umbilic surfaces lie in planes or round spheres, are likewise generalized.
Strongly isothermic maps play an important role in Part~III.

In Chapter~\ref{chapter:pluecker} we utilize an iteration of the global Darboux transformation in Chapter~\ref{chapter:weierstrass} for a proof of the Plücker formula.
The Plücker formula was proved in~\cite{FLPP} but here we generalize it to the setting of $\banach{2}\loc$\=/potentials. 
A direct consequence is the estimate of Li and Yau~\cite{LY} for general admissible maps. Finally Chapter~\ref{chapter:riemann roch} proves the Riemann-Roch theorem and Serre duality for holomorphic $\qat$\=/line bundles.

Part~II proves a version of Montel's theorem for quaternionic analysis. 
This part culminates in Theorem~\ref{compactness of conformal mappings}, which states that on a compact Riemann surface a sequence of admissible maps either into $\qat$ or into $\Imag\qat$ with bounded Willmore energy has a subsequence which converges in an appropriate sense to such an admissible map with Willmore energy bounded from above by the limit inferior of the corresponding sequence of Willmore energies.
In other words, the Willmore energy is lower semi-continuous.
The proof uses the concentration compactness principle. In Lemma~\ref{weakly continuous} at the end of Chapter~\ref{chapter:local} we proved that the inverse operators of the holomorphic structures $\delbar{E}-\V$ depend weakly continuous on $\V$, if the concentration of the energy at single points is bounded by a constant related to the Sobolev constant. From this we already know that the space of bounded holomorphic functions is locally compact, if we prevent the energy from massive concentration. The idea is now to use the conformal symmetry to deconcentrate the energy. To control the deconcentration in the limit we need some compactness of the space. This suggests to embed locally the sequence into $\mathbb{P}^1$ and then use the Möbius group of this compact space to deconcentrate the energy.

The organization of Part~II follows the lines of this idea: First we need good control over the global holomorphic structures on $\mathbb{P}^1$. 
This is done in Chapter~\ref{chapter:resolvents} by calculating explicit formulas for the resolvent of the free Dirac operator on $\mathbb{P}^1$. 
Then Chapter~\ref{chapter:weak limits} applies the concentration compactness principle locally, by embedding a sequence of holomorphic functions on a disc into the space of holomorphic functions on $\mathbb{P}^1$. 
After deconcentration of the energy, the invariance of the inverse of the Dirac operator allows for the control of the limit, if the weak limit of the holomorphic structures has a trivial kernel. Otherwise we add some poles and the same number of roots at $\infty\in\mathbb{P}^1$ until this is the case. Here we profit from the Plücker formula, the Riemann-Roch theorem, and Serre duality. 
In Chapter~\ref{chapter:global limits}, on a compact Riemann surface we construct for any sequence of holomorphic sections of a sequence of $\qat$\=/line bundles with bounded degrees and bounded $\banach{2}$\=/norms of the potentials a single Banach space, in which we can embed a subsequence of the original sequence of holomorphic sections. 
This is needed to define convergence of the subsequence. 
Finally Chapter~\ref{chapter:singular holomorphic sheaves} extracts from the additional poles and roots, which we added in Chapter~\ref{chapter:weak limits} to achieve convergence, a concept of a singular holomorphic sheaf, such that the limits are sections of such a sheaf. 
This concept is related to cusp-like singularities of algebraic curves, which is interesting in its own right.

To finish the proof of the~\ref{thm:main theorem} it remains to show in Part~III that all local minima are smooth. 
For this purpose we first define a subspace of a Banach space, whose elements correspond to admissible maps. 
We shall do this for our two representations of admissible maps, namely the Kodaira and the Weierstraß representations, separately. 
First Chapter~\ref{chapter:isothermic kodaira triples} determines all admissible maps into $\qat$ at which the space of Kodaira data is not a submanifold. 
It turns out that this are exactly the strongly isothermic admissible maps. 
This was known for smooth immersions: The subspace of conformal maps in the space of smooth immersions has a singularity exactly at the strongly isothermic surfaces (compare~\cite{BPP}). 
For admissible maps the proof depends on the characterization of the variations $\var\V$ of the Kodaira potential for which there exist variations $(\var\upsilon,\var\phi)$ of both Kodaira sections. 
The variations give rise to a tangent space.

However, only in the not strongly isothermic situation can the implicit function theorem be used to show that this tangent space is the tangent space of a Banach submanifold. 
In the strongly isothermic situation, it is not clear whether all elements of this tangent space arise as the derivative of a differentiable path of admissible maps. 
Those that do form the tangent cone, which in general is a subset of the tangent space. 
The strongly isothermic admissible maps split into two cases.
The first is more or less trivial case: the Kodaira potential is identically zero and the derivative of the Willmore energy vanishes on the tangent space.
In the second, more interesting case, by adopting a strategy of Kuwert and Schätzle~\cite{KS}, we describe the space of admissible maps as the zero set of a real smooth function on a Banach submanifold, whose Hessian is non-degenerate. 
This non-degeneracy follows from an effective local representation of the Hessian at a generic point. 
We achieve this by blowing up the surface to a flat plane. 
As a consequence of the Hessian being non-degenerate we prove that the tangent cone spans the tangent space.

In Chapter~\ref{chapter:constrained 1} we characterize the critical points of the Willmore functional. 
Fortunately, in the strongly isothermic case where the Kodaira potential vanishes, the derivative of the Willmore functional vanishes on the whole tangent space. 
Since in all other cases the tangent cone spans the tangent space, at critical points of the Willmore functional the variation of the Willmore functional vanishes on the tangent space. 
That $\var\willmore=4\pi\deg E+4\var\|\V\|_2^2=8\langle\var\V,\V\rangle$ vanishes together with the characterization of $\var\V$ for tangent vectors $(\var\V,\var\upsilon,\var\phi)$ gives a substitute of the Euler-Lagrange equation.
It also allows us to conclude smoothness by a standard elliptic bootstrapping argument.

We do not know how to describe effectively the admissible maps that take values in $\Imag\qat$ in terms of Kodaira data.
Therefore in the Chapters~\ref{chapter:isothermic weierstrass triples} and~\ref{chapter:constrained weierstrass}, we consider the more complicated space of Weierstraß data and the subspace whose admissible maps take values in $\Imag\qat$. 
In these cases, the characterization of the variation $\var\U$ of the Weierstraß potential that guarantees the existence of variations $(\var\chi,\var\psi)$ of the Weierstraß sections must be supplemented by further conditions on $\var\U$ such that the differential $\var\lp\chi,\psi\rp$ has no periods. 
Remarkably this problem is solved in Lemma~\ref{lem:wei-necessary condition} by using a pairing on the space $\Hom(H_1(\X,\mathbb{Z}),\qat)$. 
Then the other results of the Chapters~\ref{chapter:isothermic kodaira triples} and~\ref{chapter:constrained 1} have a natural counterpart for the spaces of Weierstraß data of admissible maps to $\qat$ and $\Imag\qat$.

We close the book with Chapter~\ref{ch:Flat families}. 
It incorporates into our framework those families of flat connections that were used by the Berlin group to characterize constrained Willmore immersions. 
We shall argue that such families only exists in the case of admissible maps without roots of the derivative. 
Since our approach emphasizes the strength of the setting of admissible maps with possible roots of the derivative, we prefer the description of constrained Willmore surfaces in terms of the analog of the Euler-Lagrange equation in the Theorems~\ref{constrained Willmore 1} and~\ref{constrained Willmore 2}. 
This difference between our approach and the approach of the Berlin group is related to our preference for holomorphic $\qat$\=/line bundles instead of holomorphic $\qat$\=/vector bundles. 
At the beginning of Chapter~\ref{ch:Flat families} we give an argument that the pairing on the space of potentials of a $\qat$\=/line bundle, which we use throughout the book, has no counterpart in the case of vector bundles. 
Hence for vector bundles we use the same pairing as the Berlin group.

For historical reasons let us comment on the relation of the present book to the preprint~\cite{Sch2}. This preprint does not provide enough details for the claims which are made. 
However it contains the main ideas of many statements in Part~I and most of the statements in Part~II. 
In particular, this preprint does not contain any ideas in direction of the Theorem~\ref{thm:kodaira normal} and Lemma~\ref{wente}: 
The maps that can be described in terms of holomorphic sections of $\qat$\=/line bundles with $\banach{2}\loc$\=/potentials, are not further specified and not related to Wente's inequality and to the weakly conformal maps considered in~\cite{KL, KS, Ri}. 
With exception of the substitute of the Euler-Lagrange equation in the Theorems~\ref{constrained Willmore 1} and~\ref{constrained Willmore 2}, this preprint has no overlap with Part~III. 
It overlooked completely the singularities of the space of Kodaira data and Weierstraß data at the strongly isothermic surfaces. 
For this reason the proof of the Main Theorem~8.2 in this preprint is not complete. 
More precisely, there is no argument for the smoothness given in the case of strongly isothermic minimizers.

\mainmatter

\partpage{Quaternionic Function Theory}

\chapter{Preliminaries}
\label{chapter:prelim}

\section*{Quaternions}

Let us start by introducing the quaternions with an emphasis on the geometric aspects.
The quaternions are a four dimensional real vector space $\qat = \{\alpha_0 + \alpha_1 \qi + \alpha_2 \qj + \alpha_3 \qk \}$.
A quaternion has a real part $\Real \alpha = a_0$ and an imaginary part $\Imag \alpha = \alpha_1 \qi + \alpha_2 \qj + \alpha_3 \qk$.
Unlike for complex numbers, the imaginary part of a quaternion is not a real number.
Elements of the subspace $\Imag \qat = \{\alpha_1 \qi + \alpha_2 \qj + \alpha_3 \qk \}$ are called imaginary.
Just as we often identify $\qat$ with $\mathbb{R}^4$, likewise we often identify $\Imag \qat$ with $\mathbb{R}^3$.
The $\qat$\=/conjugate is $\Bar{\alpha} = \Real \alpha - \Imag \alpha$.

Famously the quaternions have an associative but non-commutative multiplication, defined by $\qi^2 = \qj^2 = \qk^2 = \qi \qj \qk = -1$ and that multiplication with $\mathbb{R}$ is commutative.
It is said that Hamilton scratched this relation into the stone of Brougham Bridge, but we are saddened to report we could not find it.
Instead, we offer a simple trick to remember the imaginary multiplication rule.
Draw $\qi,\qj,\qk$ on a directed circle.
Multiplication of two elements gives the third, with a plus sign if they are in the correct direction and a minus sign if they are in the reverse direction.
This is of course the same rule as for the cross product in $\mathbb{R}^3$, and it makes $\qat$ a ring with the unit $1$.

There are a few ways to connect the multiplication to geometry.
Notice that when we square an imaginary quaternion $\alpha$ the mixed terms cancel out because of the anti-commutativity of $\qi,\qj,\qk$ and we are left with $\alpha^2 = - \alpha_1^2 - \alpha_2^2 - \alpha_3^2$.
Extending to $\alpha \in \qat$ we have
\begin{align*}
\alpha\Bar{\alpha} 
&= (\Real \alpha + \Imag \alpha)(\Real \alpha - \Imag \alpha)
= \alpha_0^2 + \alpha_0 \Imag \alpha - \alpha_0 \Imag \alpha - (\Imag \alpha)^2 \\
&= \alpha_0^2 + \alpha_1^2 + \alpha_2^2 + \alpha_3^2
= |\alpha|^2,
\end{align*}
the standard norm on $\mathbb{R}^4$.
Perhaps the most important fact is that the norm is multiplicative: $|\alpha\beta| = |\alpha||\beta|$.
We can prove this using the mutual orthogonality: since multiplication with $\qi, \qj,$ or $\qk$ permutes the orthonormal basis vectors $\{1,\qi,\qj,\qk\}$ (up to a sign) it follows that $\alpha, \alpha\qi, \alpha\qj, \alpha\qk$ are mutually orthogonal.
Therefore
\begin{align*}
|\alpha \beta|^2
= |\alpha \beta_0|^2 + |\alpha \beta_1 \qi|^2 + |\alpha \beta_2\qj|^2 + |\alpha \beta_3\qk|^2
= |\alpha|^2|\beta|^2.
\end{align*}
The unit quaternions (those with norm $1$) are the set $\mathbb{S}^3 \subset \mathbb{R}^4$.
They are a multiplicative subgroup of $\qat$ isomorphic to $\operatorname{SU}(2)$.
The norm shows that every non-zero quaternion has a two-sided inverse, namely $\alpha^{-1} = |\alpha|^{-2} \Bar{\alpha}$.
Therefore the quaternions are a division algebra, also called a non-commutative field.
Finally, because taking inverse reverses the order of multiplication and the norm is multiplicative we see that conjugation also reverses the order of multiplication $\overline{\alpha \beta} = \Bar{\beta}\Bar{\alpha}$.

Via the polarization identity we can derive a nice formula for the standard inner product on $\mathbb{R}^4$:
\begin{align*}
\langle \alpha, \beta \rangle_{\mathbb{R}^4}
&= \frac{1}{2}[ |\alpha + \beta|^2 - |\alpha|^2 - |\beta|^2 ]
= \frac{1}{2}[ \alpha\Bar{\beta} + \beta\Bar{\alpha} ]
= \Real (\Bar{\alpha} \beta)
\end{align*}
A less satisfying way to reach this conclusion is to directly check the formula $\Real (\alpha \Bar{b})$ on the basis $1,\qi,\qj,\qk$ and use bilinearity.
Incidentally we have shown that real part is cyclic with respect to multiplication, $\Real(\alpha \beta) = \Real(\beta \alpha)$.

The reals are orthogonal to the imaginary quaternions, so we focus on understanding orthogonality in the latter.
Generalizing from $\alpha$ and $\alpha\qi$, if $\alpha \in \qat$ and $\beta \in \Imag\qat$ then both $\alpha \beta$ and $\beta \alpha$ are orthogonal to $\alpha$, since for example
\[
\langle \alpha\beta, \alpha \rangle_{\mathbb{R}^4}
= \Real (\overline{\alpha \beta} \alpha) 
= \Real (\alpha \Bar{\alpha} \Bar{\beta}) 
= |\alpha|^2 \Real(\Bar{\beta}) 
= 0.
\]
Further, if $\alpha,\beta \in \Imag\qat$ then they are orthogonal if and only if $\Real (\Bar{\alpha} \beta) = -\Real (\alpha \beta)$ vanishes, ie if $\alpha \beta \in \Imag \qat$ too.
In this case $\{\alpha,\beta,\alpha \beta\}$ is an orthogonal basis of $\Imag\qat$.
It is always a right handed basis, as $\alpha = \qi, \beta = \qj$, and $\alpha \beta = \qk$ shows.


Orthogonality is key to understanding commutation and anti-commutation in the quaternions.
If $\alpha, \beta \in \Imag\qat$ are orthogonal, i.e.\ $\alpha \Bar{\beta} \in \Imag \qat$, then they anti-commute: $\alpha\beta = -\alpha\Bar{\beta} = \overline{ \alpha \Bar{\beta} } = \beta \Bar{\alpha} = - \beta \alpha$.
The converse is also true: $\Real (\Bar{\alpha} \beta) = - \Real(\alpha\beta) = \Real(\beta \alpha) = - \Real (\Bar{\alpha} \beta)$.
Therefore every imaginary quaternion $\alpha$ decomposes $\qat$ into the subspace $\mathbb{R}\{1,\alpha\}$, which commutes with $\alpha$, and the subspace $\mathbb{R}\{1,\alpha\}^\perp = \Imag\qat \cap \alpha^\perp$, which anti-commutes.

We have already calculated that the square of an imaginary quaternion is real and negative.
The converse is also true.
In particular a quaternion is a square root of $-1$ if and only if it is unit and imaginary.
We denote the set of unit imaginary quaternions by $\mathbb{S}^2$.
This follows since $\nu^2 = -1$ iff $\nu^{-1} = -\nu$ iff $|\nu|^2 \nu = - \Bar{\nu}$ (for $\nu$ non-zero).
For any such $\nu$, we have that $\mathbb{R}\{1,\nu\}$ is isomorphic to the complex plane.
Any unit quaternion can be written in polar form as $\alpha = \cos\theta + \sin\theta\, \nu$ for some $\nu \in \mathbb{S}^2$.
Squaring this expression leads to $\alpha^2 = \cos2\theta + \sin2\theta\, \nu$, which lets us see that any non-real quaternion has exactly two square roots, both of which lie in the plane $\mathbb{R}\{1,\nu\}$.

\begin{definition}[Notations for sets]\phantom{boo}
\begin{enumeratethm}
\item
Set of imaginary quaternions $\Imag \qat = \{\alpha_1 \qi + \alpha_2 \qj + \alpha_3 \qk \}$.
\item 
Set of unit quaternions $\mathbb{S}^3 \subset \qat$.
\item 
Set of unit imaginary quaternions $\mathbb{S}^2 = \mathbb{S}^3 \cap \Imag \qat$.
\item
The standard copy of $\mathbb{C}$ in $\qat$ is also called $\qat^+ = \mathbb{R}\{1,\qi\}$.
\item
Its orthogonal complement $\qat^- = (\qat^+)^\perp =  \qj\mathbb{C} = \mathbb{R}\{\qj,\qk\}$.
\item
Additionally, the authors disagree about whether $\mathbb{N}$ includes $0$.
Where it is significant, we prefer $\mathbb{N}_0$ or $\mathbb{N}_+$ to be clear.
\end{enumeratethm}
\end{definition}

\section*{Rotations}

The decomposition $\qat = \mathbb{R}\{1,\alpha\}\oplus \mathbb{R}\{1,\alpha\}^\perp$ makes plain the `double rotation' action of quaternions on themselves.
Write any $\alpha \in \mathbb{S}^3$ in polar form $\alpha = \cos\theta + \sin\theta\, \nu$ for $\nu \in \mathbb{S}^2$.
We know that the subalgebra $\mathbb{R}\{1,\alpha\} = \mathbb{R}\{1,\nu\}$ is isomorphic to $\mathbb{C}$, hence multiplication from either the left or right by $\alpha$ is rotation by $\theta$.
On the other hand, for any $\beta \in \mathbb{R}\{1,\alpha\}^\perp \cap \mathbb{S}^2$ we have an ordered orthonormal basis $\beta, \nu\beta$ of $\mathbb{R}\{1,\alpha\}^\perp$.
Acting $\alpha$ from the left sends $\beta$ to $\cos\theta\,\beta + \sin\theta\,\nu\beta$ and $\nu\beta$ to $-\sin\theta\,\beta + \cos\theta\,\nu\beta$.
We recognize that this is also rotation by $\theta$.
Acting $\alpha$ from the right is rotation by $-\theta$.

Building on this, another operation is $K_\alpha(x) = \alpha x \alpha^{-1}$.
The formula is unchanged by rescaling $\alpha$, so it is typical to take $\alpha$ to be unit-length, whereupon $K_\alpha(x) = \alpha x \Bar{\alpha}$.
This is a group action of $\mathbb{S}^3$ on $\qat$.
If we again write $\alpha$ in polar form then $\Bar{\alpha} = \cos(-\theta) + \sin(-\theta)\, \nu$.
By applying our understanding from the previous paragraph, we deduce that the operation $K_\alpha$ fixes the subspace $\mathbb{R}\{1,\alpha\}$ while rotating the plane $\mathbb{R}\{1,\alpha\}^\perp$ by angle $2\theta$.
In particular, since $K_\alpha$ preserves $\mathbb{R}$ it also preserves $\Imag \qat = \mathbb{R}^3$.
The restriction to $\Imag \qat$ fixes the line $\mathbb{R}\nu$ while rotating the perpendicular plane.
This is the quaternionic representation of rotations in $\mathbb{R}^3$, also known as the double cover of $\operatorname{SO}(3)$ by $\operatorname{SU}(2)$ since $K_{\alpha} = K_{-\alpha}$.
A simple example is that of $-\qk x \qk = x_0 - x_1\qi - x_2\qj + x_3\qk$, where the plane $\mathbb{R}\{1,\qk\}$ is fixed and the plane $\mathbb{R}\{\qi,\qj\}$ is rotated by $2\theta = \pi$.
This is one example of the special case $\alpha \in \mathbb{S}^2$, where $K_\alpha$ fixes the commuting plane and reflects the anti-commuting plane.
Analogously we may define $L_\alpha(x) = \alpha x \alpha$.
In this case it is the plane $\mathbb{R}\{1,\alpha\}$ that is rotated by $2\theta$, whereas the orthogonal plane $\mathbb{R}\{1,\alpha\}^\perp$ is fixed.

We have seen that both left and right multiplication by quaternions gives rotations of $\qat = \mathbb{R}^4$.
A natural challenge therefore is to describe all rotations $\operatorname{SO}(4)$ in terms of quaternions.
Suppose $T \in \operatorname{SO}(4)$ and take a square root $\beta$ of $T(1)$.
Then $L_{\beta^{-1}} \circ T$ gives a rotation that fixes $1$.
By orthogonality, this composition must restrict to give an $\operatorname{SO}(3)$ rotation on the imaginary quaternions.
These are classical rotations, described by an invariant axis $\mathbb{R}\nu$ (for $\nu\in\mathbb{S}^2$) and a rotation angle $2\theta$.
In other words, the composition is nothing other than $K_\alpha$ for $\alpha = \cos\theta + \sin\theta\,\nu$.
Undoing the composition proves that $T(x) = L_{\beta} (K_\alpha(x)) = (\beta\alpha)x(\Bar{\alpha}\beta)$.
This shows we have a surjective group homomorphism from $\mathbb{S}^3\times\mathbb{S}^3$ to $\operatorname{SO}(4)$, and one can calculate that the kernel is $\{(1,1),(-1,-1)\}$.

\section*{Normals}

Let us introduce the idea of the left and right normal of a plane.
In $\mathbb{R}^3$ it is fruitful to describe a plane in terms of its orthogonal line, because lines are simpler than planes.
Given a plane in $\qat$ however the orthogonal complement is another plane.
Moreover only using orthogonality of $\mathbb{R}^4$ misses what makes $\qat$ special.
Instead we encode a plane (through the origin) as a rotation.
Specifically, we identify a plane $E$ with the rotation of $\qat$ that fixes $E$ and negates $E^\perp$.
This rotation has the form $T(x) = Nx(-R)$ for unit quaternions $N,R \in \mathbb{S}^3$.
Our convention differs from that of~\cite[Lemma~2]{BFLPP} by the sign of $R$.
Since $T$ is an involution, we can conclude that $(N^2, R^2)$ is either $(1,1)$ or $(-1,-1)$, due to the double covering of $\operatorname{SO}(4)$.
The possibility $(1,1)$ implies $T = \pm \id_\qat$,  neither of which singles out a plane.
Therefore the plane is described by $N,R\in\mathbb{S}^2$ unique up to the choice of sign $(N,R)$ or $(-N,-R)$.
The plane $E$ is exactly the points fixed by the rotation $T$ and so equivalently described as the quaternions $x$ for which $Nx = xR$.
Moreover, the rotation $x \mapsto Nx$ preserves the plane since it commutes with $T$.
Because $x \mapsto Nx$ applied twice is $-\id_\qat$, $N$ itself must act on the plane by a quarter turn.
The same argument applies to $R$.
If the plane is oriented then the sign of $N$ can be chosen to respect this orientation, in which case $R$ does too.

\begin{definition}
\label{def:left and right normals}
\index{Normal}
Given any real oriented plane $E \subset \qat$, the \emph{left} and \emph{right normals} of the plane are respectively the unique $N,R \in \mathbb{S}^2$ such that $E = \{ x \in \qat \mid N x = x R \}$ and left multiplication by $N$ agrees with the orientation.
\end{definition}

If the plane lies in $\Imag \qat$, then $T$ is one of the special rotations that we have previously discussed, namely $L_N$ for $N$ orthogonal to the plane and $R = -N$.
The name left and right normal comes from this case; in general $N$ and $R$ are not orthogonal to the plane.

\begin{example}[Plane]
\label{eg:plane normals}
\index{Plane}
This example illustrates that $\mathbb{C} \subset \qat$ is the prototypical case of an oriented plane in the quaternions, with left and right normal $\qi$, and all other cases are rotations of this.
Suppose that we have a plane that is the right $\mathbb{C}$\=/span of $(1-\qk)\qj$, that is $E = \{ (1-\qk)\qj h \mid h \in \mathbb{C} \}$.
Clearly multiplication from the right by $\qi$ rotates the plane by a quarter-turn, and let this be its orientation.
Hence $R = \qi$.
The trick to find the left normal is to view our plane as $\mathbb{C} \subset \qat$ rotated by left multiplication by $(1-\qk)\qj$.
Then 
\[
\Big[(1-\qk)\qj\, \qi\, \qj^{-1}(1-\qk)^{-1}\Big]\,(1-\qk)\qj h
= (1-\qk)\qj \, \qi h
= (1-\qk)\qj h \,\qi
\]
expresses the operation of rotating $E$ to $\mathbb{C}$ with $(1-\qk)\qj$, applying $\qi$, and rotating back with $\qj^{-1}(1-\qk)^{-1}$.
It shows that 
\[
N 
= (1-\qk)\qj\, \qi\, \qj^{-1}(1-\qk)^{-1}
= (1-\qk)(-\qi)(1-\qk)^{-1}
\]
is the left normal of $E$.
It is easily seen to obey $N^2 = (1-\qk)(-1)(1-\qk)^{-1} = -1$.
In fact this is the same $N$ as Example~\ref{eg:catenoid} at $z=1$, and we will build on this example through this chapter.
\end{example}

\section*{Linear Algebra}

There are two approaches to explaining linear algebra over the quaternions.
One approach is to treat $\qat$ as sui generis and consider the theory as a special case of modules over a ring.
The other approach is to recognize that the quaternions are an algebra over the reals and so to take vector space theory over $\mathbb{R}$ as one's starting point.
Although the latter has the advantage that linear algebra over the quaternions is indeed quite close to that over the reals, it has the disadvantage of making some constructions seem ad hoc.

\begin{definition}
\label{def: qat vector space}
A $\qat$\=/vector space $\Spa{A}$ is a right $\qat$\=/module.
This means that $\Spa{A}$ is an abelian group with a right $\qat$\=/action: for all $a,b \in \Spa{A}$ and $\gamma,\delta \in \qat$
\begin{enumeratethm}
\item $(a+b)\gamma = a\gamma + b\gamma$ 
\item $a(\gamma + \delta) = a\gamma + a\delta$
\item $a 1 = a$
\item Right action: $a (\gamma \delta) = (a \gamma) \delta$
\end{enumeratethm}
\end{definition}

The prototypical example is $M = \qat^d$ with right multiplication as the scaling action.
Another important example is $M = \qat^d$ with left multiplication by the $\qat$\=/conjugate.
This second example is a right action because 
\[
(\gamma\delta)\cdot x 
= \overline{\gamma\delta}x 
= \Bar{\delta}\Bar{\gamma}x 
= \delta \cdot (\gamma\cdot x).
\]
The adjectives left and right describing an action refer to how the action composes.
Although ordinary left and right multiplication are the motivating examples of left and right actions, this example shows that they are not the only one.

Much of the theory of vector spaces over fields carries over to $\qat$\=/vector spaces without modification, including bases and dimension.
A (right) $\qat$\=/homomorphism $\V$ between two $\qat$\=/vector spaces, also called a (right) $\qat$\=/linear map, is a homomorphism of abelian groups that commutes with the scaling actions $\V(a\gamma) = V(a)\gamma$.
We denote these functions by $\Hom_\qat(\Spa{A},\Spa{B})$ or $\End_\qat(\Spa{A}) = \Hom_\qat(\Spa{A},\Spa{A})$.
They are completely determined by how they act on a basis and thus in finite dimensions can be described as a matrix in the usual way.
Here is an example of an endomorphism acting on an element of $\qat^2$
\[
\begin{pmatrix}
1 & \qi \\ 0 & -\qk
\end{pmatrix}
\begin{pmatrix}
1 \\ \qk  
\end{pmatrix}
= 
\begin{pmatrix}
1 - \qj \\ 1
\end{pmatrix}.
\]
Notice that since the matrix is on the left, it doesn't `get in the way' of scaling from the right.
This is the reason that we prefer right $\qat$\=/vector spaces: it can be laid out just like linear algebra over the reals, with which we are familiar.

An important difference to the commutative case is that scaling by a fixed $\alpha \in \qat$ is not in general an endomorphism of $\Spa{A}$.
The function $a \mapsto a\alpha$ is right $\qat$\=/linear only if $(a\alpha)\gamma = (a\gamma) \alpha$ for all $\gamma \in \qat$, which is only the case if $\alpha \in \mathbb{R}$.
To express the same fact in a different way: if with respect to some basis an endomorphism has the matrix $\alpha I$ such that $\alpha \not\in \mathbb{R}$, it need not be represented by this matrix with respect to another basis.
Therefore $\Hom_\qat(\Spa{A},\Spa{B})$ is a $\mathbb{R}$\=/vector space with pointwise addition and scaling given by the composition with $r\id_\Spa{A}$ for $r \in \mathbb{R}$, but it is not a $\qat$\=/vector space.

This raises the question: given a $\mathbb{R}$\=/vector space, what additional data needs to be given such that it becomes a $\qat$\=/vector space?
It is the choice of two endomorphisms $\qi, \qj \in \End_\mathbb{R}(\Spa{A})$ such that $\qi^2 = \qj^2 = -\id_M$ and $\qi \qj = - \qj \qi$, which is commonly called a \emph{quaternionic structure}.
The difference between a quaternionic structure and Definition~\ref{def: qat vector space} is that in the former the underlying $\mathbb{R}$\=/vector space is already fixed.
This technical point as consequences when considering several actions of $\mathbb{C}$ and $\qat$ on the same space: due to the existence of injective field homomorphisms from $\mathbb{R}$ onto proper subfields of $\mathbb{R}$ 
, it is not guaranteed that the induced actions of $\mathbb{R}$ agree.

\begin{definition}
A \emph{complex structure} $J$ on a $\mathbb{R}$\=/vector space $\Spa{A}$ is an element $J \in \End_\mathbb{R}(\Spa{A})$ such that $J^2 = -\id_\Spa{A}$.
A complex structure on a vector space is equivalent, via $(r + s \qi)\cdot a = (r\id_\Spa{A} + sJ)a$, to a $\mathbb{C}$\=/action on $\Spa{A}$ such that any $r \in \mathbb{R} \subset \mathbb{C}$ acts as $r\id_\Spa{A}$.

Two actions are \emph{compatible} if they commute with one another.
For a complex structure to be compatible with the scaling action of a $\qat$\=/vector space means that $J \in \End_\qat(\Spa{A})$.
We use the adjective \emph{complex} for $\qat$\=/vector spaces with a compatible complex structure.
\end{definition}

Other terms for this structure are a linear complex structure~\cite[II~\S{}1.1]{Aud} or an almost complex structure~\cite[Definition 1.2.1]{Hu}.
On the tangent bundle of a manifold the difference between almost complex structures and complex structures is only relevant for complex manifolds of higher dimension.
The Newlander–Nirenberg theorem~\cite{NN} characterizes when an almost complex structure is induced by a holomorphic atlas.
The condition is always satisfied for Riemann surfaces, so we will not fret over this distinction.

Compatibility of actions is the key to understanding the tensor product of modules.
Recall, a right module $\Spa{A}$ and a left module $\Spa{B}$ over the same ring $S$ can be tensored to construct the abelian group $\Spa{A} \otimes_S \Spa{B}$, whose elements are pairs of the form $a\otimes b$ with equivalences for bi-additivity and $a\gamma \otimes b = a \otimes \gamma b$ for $\gamma \in S$.
But this `uses up' the $S$\=/action.
Suppose we were try to define an action of $S$ on the tensor product by $\delta(a\otimes b) = a \delta \otimes b = a \otimes \delta b$, which at first blush looks well-defined.
But in fact
\begin{equation}
\label{eqn:tensor action definition problem}
\delta (a\gamma \otimes b) = a \gamma \delta \otimes b
\quad\text{compared to}\quad
\delta (a\otimes \gamma b) 
= a \delta \otimes \gamma b
= a \delta \gamma \otimes b,
\end{equation}
shows that it is not well defined if $\gamma,\delta \in S$ are non-commuting.
Hence there is no tensor product for $\qat$\=/vector spaces.
This problem is avoided if the actions on the factor modules are compatible with the $S$\=/action~\cite[II~\S3.4]{Bo}, in which case the above is well-defined.
For example, if an $S$\=/module $\Spa{A}$ has a compatible left $R$\=/action then $\Spa{A} \otimes_S \Spa{B}$ is a left $R$\=/module.
This is how we should understand the familiar tensor product over $\mathbb{R}$ or $\mathbb{C}$: every action of a commutative ring is both a left and right action and is compatible with itself, hence there is no need to match left and right and a `copy' of the action is inherited by the tensor product~\cite[II~\S3.5]{Bo}.
Similarly, compatible actions on $\Spa{B}$ are inherited by the set of $S$\=/homomorphisms $\Hom_S(\Spa{A},\Spa{B})$ through post-composition.
Compatible actions on $\Spa{A}$ are also inherited by $\Hom_S(\Spa{A},\Spa{B})$ through pre-composition, but with the wrinkle that a left action becomes right and vice versa~\cite[II~\S1.14]{Bo}.
To see this, let a left action of $R$ on $\Spa{A}$ be described by $R_\alpha \in \End_S(\Spa{A})$ for each $\gamma \in R$ on $\Spa{A}$:
\[
V \circ R_{\gamma\delta}
= V \circ (R_{\gamma} \circ R_{\delta})
= (V \circ R_{\gamma}) \circ R_{\delta}.
\]

Let us consider tensor products of a complex $\qat$\=/vector space $\Spa{A}$ more concretely.
Any $\mathbb{C}$\=/vector space $E$ can be tensored with a complex $\qat$\=/vector space $\Spa{A}$ to produce $E \otimes_\mathbb{C} \Spa{A}$.
A complex scalar can pass between the two factors: $eh \otimes a = e \otimes ha$ for $h\in\mathbb{C}$.
The right action of $\qat$ is inherited by the tensor product, $(e\otimes a)\gamma = e \otimes a\gamma$.
The complex structure of $E$ is also inherited, $J(e\otimes a) = e\qi \otimes a$; since it is compatible with the action used to create the tensor product there is no problem with this definition.
In fact every complex $\qat$\=/vector space is of this form:
\begin{lemma}
\label{lem:quaternionification}
\index{Underlying bundle}
Every complex $\qat$\=/vector space $\Spa{A}$ can canonically realized as a tensor product $E_\qat := E \otimes_\mathbb{C} \qat$ for $E$ a $\mathbb{C}$\=/vector space and $\qat$ with the complex structure given by left multiplication by $\qi$.
$E$ is called the underlying $\mathbb{C}$\=/vector space.
\end{lemma}
\begin{proof}
The key idea is to find a subspace of $\Spa{A}$ on which the complex structure $J$ agrees with the complex structure given by right multiplication by $\qi$.
First we consider $\Spa{A}$ as a $\mathbb{C}$\=/vector space by restricting the scalars to $\qat^+ \subset \qat$.
Then there is a $\mathbb{C}$\=/basis of $\Spa{A}$ that diagonalizes the complex structure $J$ with eigenvalues $\lambda = \pm \qi$.
If $Ja = a \lambda$ then $J(a\qj) = a \lambda \qj = (a\qj) \Bar{\lambda}$.
So the eigenvectors come in pairs $a, a\qj$ with eigenvalues $\qi, -\qi$ respectively.
Defining the $\mathbb{C}$\=/vector space $E = \{ a \in \Spa{A} \mid Ja = a\qi\}$, we see that $\Spa{A} = E \otimes_\mathbb{C} \qat$.
\end{proof}

\begin{example}
Consider $\Spa{A} = \qat$ as a $\qat$\=/vector space with right multiplication by scalars.
Let the complex structure be left multiplication by $J = (1-\qk)(-\qi)(1-\qk)^{-1}$.
Because it is left multiplication, the complex structure is compatible with scaling by $\qat$.
We saw in Example~\ref{eg:plane normals} that $E = \{ (1-\qk)\qj h \mid h \in \mathbb{C} \}$ was exactly the subset of $\qat$ such that $Ja = a \qi$.
Choose a nonzero element $\mathbf{e} \in E$.
If we write a general element of $\qat$ as $\mathbf{e} \gamma$ for $\gamma \in \qat$ then
$J(\mathbf{e} \gamma) = (J \mathbf{e}) \gamma = \mathbf{e} (\qi \gamma)$.
Relative to the basis $\mathbf{e}$ of $\Spa{A}$ the element $\mathbf{e}\gamma$ is represented by column vector $(\gamma)$, on which scaling by $\alpha \in \qat$ gives $(\gamma\alpha)$ and the action of $J$ gives $(\qi \gamma)$.
\end{example}

The correspondence between complex $\qat$\=/vector spaces and $\mathbb{C}$\=/vector spaces of the above lemma respects tensoring by a $\mathbb{C}$\=/vector space from the left, 
\[
F \otimes_\mathbb{C} E_\qat 
= F \otimes_\mathbb{C} ( E \otimes_\mathbb{C} \qat )
= (F \otimes_\mathbb{C} E) \otimes_\mathbb{C} \qat
= (F \otimes_\mathbb{C} E)_\qat,
\]
since the tensor product is associative.
It is natural to ask which other operations on $\mathbb{C}$\=/vector spaces can be extended to complex $\qat$\=/vector spaces.
From the distributive property of direct sums: $E_\qat \oplus F_\qat = (E\oplus F)_\qat$.

The dual of a vector space is the set of linear maps from it to its field.
For example, $E^{-1} = \Hom_\mathbb{C}(E,\mathbb{C})$ for a $\mathbb{C}$\=/vector space $E$.
In the quaternionic case, we consider $\Hom_\qat(E_\qat, \qat)$.
Since $\qat$ as a right $\qat$\=/vector space also has a compatible left quaternionic structure coming from left multiplication, the set $\Hom_\qat(E_\qat, \qat)$ inherits the left quaternionic structure from $\qat$ and a compatible complex structure from $E_\qat$.
As we wish to work with $\qat$\=/vector spaces, that is to say right $\qat$\=/modules, we work instead with the $\qat$\=/conjugate of the scaling action.

\begin{lemma}
\label{lem:quaternionic dual space}
Define the dual of complex $\qat$\=/vector space $(E_\qat)^{-1}$ as the set $\Hom_\qat(E_\qat, \qat)$ with the right $\qat$\=/action $(V \cdot \gamma)(a) = \Bar{\gamma}V(a)$ and the complex structure $(JV)(a) = V(Ja)$, for $V \in \Hom_\qat(E_\qat, \qat)$, $a \in E_\qat$, and $\gamma \in \qat$.
Compare with~\cite[Example~9, Theorem~1]{BFLPP}.
Then $(E_\qat)^{-1} \cong (E^{-1})_\qat$.
\end{lemma}
\begin{proof}
As in Lemma~\ref{lem:quaternionification}, we should find the elements $V \in (E_\qat)^{-1}$ such that the complex structure acts as multiplication from the right by $\qi$.
Due to $\qat$\=/linearity any homomorphism is determined by how it acts on $E$, so it suffice to consider only $e \in E$.
We compute
\begin{align*}
(JV)(e)
= V(Je)
= V(e \qi)
= V(e) \qi, 
\quad\text{and}\quad
(V\cdot \qi)(e)
= \Bar{\qi}V(e)
= -\qi V(e).
\end{align*}
We see that these are equal if and only if $V$ acts on $E$ to give values in $\qj \mathbb{C}$, which suggests we can write $V$ in a special form.
Let $\Tilde{V} \in \Hom_\mathbb{C}(E,\mathbb{C})$ be $\Tilde{V}(e) = -\qj (V|_E)(e)$.
Then $V$ is the $\qat$\=/linear extension of $\qj\Tilde{V}$ to $E_\qat = E \otimes_\mathbb{C} \qat$.
This gives a bijection between $E^{-1}$ and the underlying $\mathbb{C}$\=/vector space of $(E_\qat)^{-1}$.
Moreover this bijection is an isomorphism of $\mathbb{C}$\=/vector spaces: it is an isomorphism of $\mathbb{R}$\=/vector spaces and if $\Tilde{V} \in E^{-1}$ and $\Tilde{W} = \qi \Tilde{V}$ then
\[
W(e) 
= \qj \Tilde{W}(e)
= \qj \qi\Tilde{V}(e)
= -\qi \qj \Tilde{V}(e)
= -\qi V(e)
= (V \cdot \qi)(e).
\qedhere
\]
\end{proof}

It is somewhat counterintuitive that although $(E_\qat)^{-1} \cong (E^{-1})_\qat$ the correspondence isn't simply restriction and $\qat$\=/linear extension, but rather there is an additional factor of $\qj$ involved.
Indeed, there is even a `copy' of $\Hom_\mathbb{C}(E,\mathbb{C})$ inside of $\Hom_{\qat}(E_\qat,\qat)$: take an element of $\Tilde{V} \in \Hom_\mathbb{C}(E,\mathbb{C})$ and extend it $\qat$\=/linearly to obtain $V \in \Hom_{\qat}(E_\qat,\qat)$ (without a factor of $\qj$).
However the inclusion map of $\Hom_\mathbb{C}(E,\mathbb{C})$ to this subset is conjugate-linear (without the $\qj$ there is a change of sign in the last equation of the above proof).
At root this is because of the conjugation in the scaling action of $E_\qat^{-1}$, which was necessary in order to have a right $\qat$\=/action.

\begin{definition}
\label{def:00 pairing}
\index{Pairing!(0,0)}
Throughout this book we employ various pairings.
To fit with later definitions, we use the symbol $\lz \cdot,\cdot \rz$ for the natural pairing between $\Hom_\qat(E_\qat,\qat)$ and $E_\qat$.
The definitions in the above lemma are
\begin{align*}
\lz V \cdot \gamma, a \delta \rz
=\Bar{\gamma} \lz V, a \rz \delta,
\quad\text{and}\quad
\lz JV, a \rz = \lz V, Ja \rz,
\end{align*}
for $V \in \Hom_\qat(E_\qat,\qat)$, $a \in E_\qat$, and $\gamma,\delta \in \qat$.
\end{definition}

Let us show how to represent calculations with the dual space as matrices.
In linear algebra, given a basis $\{e_l\}$ of a vector space $E$ one almost always works with the dual basis $\{e^\ast_m\}$ of $E^{-1}$ defined by $e^\ast_m(e_l) = \delta_{lm}$.
In particular, with respect to these bases the pairing of $E^{-1}$ and $E$ is that of a row vector matrix multiplied with a column vector.
However this is only a convenient convention.
For quaternionic vector spaces, it is better to choose a basis of $E_\qat^{-1}$ that lies in the underlying $\mathbb{C}$\=/vector space.
One choice is $e^\ast_m(e_l) = \qj \delta_{lm}$.
We call this a $\qj$\=/dual basis.
With respect to this basis
\begin{equation}
\label{eq:00 pairing in basis}
V(a) 
= \left\lz \sum_m e^\ast_m \V_m, \sum_l e_l a_l \right\rz
= \sum_{m,l} \Bar{V}_m \lz e^\ast_m, e_l \rz a_l
= \sum_{l} \Bar{V}_l \qj a_l.
\end{equation}

As noted previously, for two $\qat$\=/vector spaces $E_\qat,E'_\qat$ with complex structures $J,J'$, the set $\Hom_\qat(E_\qat,E'_\qat)$ is not a $\qat$\=/vector space.
But it does inherit two complex structures.
We are particularly interested in those homomorphisms $\Hom_\qat^+(E_\qat,E'_\qat)$ that commute, $J' V = V J$, and those homomorphisms $\Hom_\qat^-(E_\qat,E'_\qat)$ that anti-commute, $J' V = - V J$.
Any homomorphism can be decomposed as 
\[
V = V^+ + V^- = \frac{1}{2}(V - J'VJ) + \frac{1}{2}(V + J'VJ),
\]
giving a direct sum $\Hom_\qat(E_\qat,E'_\qat) = \Hom_\qat^+(E_\qat,E'_\qat)\oplus \Hom_\qat^-(E_\qat,E'_\qat)$.
Commuting homomorphisms are exactly those such that $V(E) \subseteq E'$, since for $e \in E$ we have $V(e) \in E'$ if and only if $J'V(e) = V(e)\qi = V(e\qi) = V(Je)$.
This shows that the restriction of $V$ to $E$ gives $V|_E \in \Hom_\mathbb{C}(E,E')$ and any element of $\Hom_\mathbb{C}(E,E')$ can be extended $\qat$\=/linearly to an element of $\Hom_\qat^+(E_\qat,E'_\qat)$.
In contrast to the proof of Lemma~\ref{lem:quaternionic dual space}, this inclusion is $\mathbb{C}$\=/linear, so we can safely identify $\Hom_\qat^+(E_\qat,E'_\qat)$ and $\Hom_\mathbb{C}(E,E')$.
In the particular case of a one dimensional vector space $E_\qat$ the commuting endomorphisms $\End_\qat^+(E_\qat)$ are identified with $\End_\mathbb{C}(E)$, which is canonically isomorphic to $\mathbb{C}$.

Before we discuss $\Hom_\qat^-(E_\qat,E'_\qat)$ it is useful to recall the conjugate $\bar{E}$ of a $\mathbb{C}$\=/vector space $E$.
This is defined to be the same abelian group as $E$ with the conjugate scaling action.
Because the underlying sets of $E$ and $\bar{E}$ are the same, often a bar is used to distinguish whether an element is to be interpreted as belonging to $E$ or $\bar{E}$.
With this convention, the scaling action can be concisely written as $\bar{e}h = \overline{e\bar{h}}$, for $h \in \mathbb{C}$, $e \in E$, and $\bar{e}, \overline{e\bar{h}} \in \bar{E}$.

Now suppose $V \in \Hom_\qat^-(E_\qat,E'_\qat)$.
Then $J'V(e) = -V(Je) = -V(e\qi) = -V(e)\qi$, so $V(E) \subset \{ e' \otimes \qj \mid e' \in E' \}$.
Observe that the subset $\{ e' \otimes \qj \mid e' \in E' \}$ considered as a $\mathbb{C}$\=/vector space with the action of $\qat^+$ is isomorphic to $\bar{E'}$, with the map being $e' \otimes \qj \mapsto \overline{e'}$.
Therefore $V$ can be identified with a linear map $\tilde{V} : E \to \bar{E'}$, and vice versa, with the formula $V(e \otimes \lambda) = \overline{\tilde{V}(e)} \otimes \qj \lambda$.
We can see how this applies in particular to case of $E' = \bar{E}$. 
The identity map $id_E : E \to E = \overline{\bar{E}}$ gives rise to the anti-commuting `conjugation' map
\begin{equation}
\label{eq:quaternionic bar map}
E_\qat \to \bar{E}_\qat,\qquad
e \otimes \lambda \mapsto \bar{e} \otimes \qj \lambda.
\end{equation}

\section*{Quaternionic bundles}

It is generally the case that what can be defined for vector spaces can be defined for vector bundles fiberwise.
Complex or quaternionic structures on a $\mathbb{R}$\=/vector bundle can be defined as sections of the endomorphism bundle such that on each fiber we have a complex or quaternionic structure respectively.
A (right) $\qat$\=/bundle over a manifold is a $\mathbb{R}$\=/bundle with a right quaternionic structure.
A complex $\qat$\=/bundle has additionally a compatible complex structure.
Then Lemma~\ref{lem:quaternionification} can be generalized to write any complex $\qat$\=/bundle as $E \otimes_\mathbb{C} \qat$ where $E$ is a $\mathbb{C}$\=/vector bundle, $\qat$ stands for the trivial $\qat$\=/vector bundle, and the tensor product is taken fiberwise.
Homomorphism and dual bundles are likewise defined.
We see that $E$ is a $\mathbb{R}$\=/subbundle defined by $J a - a \qi = 0$; an application of the constant rank theorem tells us that if $J$ is smooth then $E$ is a smooth subbundle.

\begin{example}[Catenoid]
\label{eg:catenoid underlying line bundle}
\index{Catenoid!Holomorphic line bundle}
Consider the trivial $\qat$\=/line bundle $\X \times \qat$ over the annulus $\X = \mathbb{C}^\ast$.
Now endow this bundle with the complex structure $J_N$ given by left multiplication $J_N(z, \alpha) = (z, N(z)\alpha)$, for $N(z) = (1-z^{-1}\qk) (-\qi) (1-z^{-1}\qk)^{-1}$ the left normal of the catenoid from Example~\ref{eg:catenoid}.
At the point $z=1$ this formula for $N$ reduces to Example~\ref{eg:plane normals}.
Using an analogous calculation, we see that 
\[
E = \{ (z, (1-z^{-1}\qk)\qj h ) \in \X \times \qat \mid h \in \mathbb{C} \}
\]
is the underlying $\mathbb{C}$\=/bundle.
\end{example}

\begin{remark}
\label{rem:special cover}
Without loss of generality, we can restrict the sort of sets on which we work.
Given a Riemann surface $\X$ we fix a locally finite open cover $\{\SO_l\}$ such that each $\SO_l$ is the contained in the domain of a chart $z_l : \Tilde{\SO}_l \to \Tilde{\Omega}_l \subset \mathbb{C}$ whose image is a ball $\Tilde{\Omega}_l = B(0,R)$ and $\Omega_l = z_l[\SO_l]= B(0,r)$ is itself a ball.
In the case that $\X$ is compact, this cover is necessarily finite.

In particular $\Tilde{\SO}_l$ is simply connected and contractible, ensuring that any $\mathbb{C}$\=/bundle $E$ trivializes over it.
We work with trivializations of $E$ over $\SO_l$ that are the restriction of a trivialization over $\Tilde{\SO}_l$.
Since $\SO_l$ is relatively compact in $\Tilde{\SO}_l$ and $\X$, any smooth functions on $\SO_l$ that are the restriction of smooth functions on $\Tilde{\SO}_l$ are well-behaved on the boundary $\boundary \SO_l$.

Note however that we do not fix the chart at the outset.
At certain places, for example Remark~\ref{rem:cover depends on q}, it is necessary to impose additional conditions on the sets $\SO_l$.
\end{remark}

What is perhaps unfamiliar is how to represent sections of $\qat$\=/bundles in a trivialization.
With respect to a cover described in the above remark, a $\mathbb{C}$\=/bundle $E$ is described by a cocycle $f_{ml} : \SO_l \cap \SO_m \to \mathrm{GL}(\mathbb{C}^d)$.
Then a section $\xi$ of $E$ is described by a collection of functions $\xi_l : \SO_l \to \mathbb{C}^d$ such that $\xi_m =f_{ml}\,\xi_l$ on overlaps.
This motivates
\begin{definition}
\label{def:sections of qat bundle}
A section $\xi$ of $E_\qat$ is described by a collection functions $\xi_l: \SO_l \to \qat^d$ with
\begin{align}
\label{eq:section transformation} 
\xi_m&=f_{ml}\,\xi_l, 
\end{align}
where $f_{ml}$ are the cocycles of the underlying bundle $E$.
The $\xi_l$ are called local representatives of $\xi$.
We will generally use the convention that the local representatives have subscripts and the global object does not.
The complex structure acts by $(J\xi)_l = i \xi_l$ and quaternionic scaling by $(\xi \alpha)_l = \xi_l \alpha$.
\end{definition}

Of course, by forgetting the quaternionic structure $E_\qat$ is a $\mathbb{R}$\=/bundle and so already has a notion of sections.
The relationship to the above definition is that we can consider a section $\xi$ of $E_\qat$ locally as a $\qat$\=/linear combination of a local frame $e=(e_1,\ldots,e_d)$ of $E$, i.e $\xi =e_1\otimes\xi_{l,1}+\cdots+e_d\otimes\xi_{l,d}$ for a function $\xi_l=(\xi_{l,1},\ldots,\xi_{l,d}) : \SO_l \to \qat^d$.
We now know how $J$ acts on $\xi$:
\[
J \left(\sum_{n=1}^d e_n \otimes\xi_{l,n}\right)
= \sum_{n=1}^d Je_n \otimes\xi_{l,n}
= \sum_{n=1}^d e_n \qi\otimes\xi_{l,n}
= \sum_{n=1}^d e_n \otimes \qi\xi_{l,n},
\]
and likewise for $\qat$\=/scaling.
If we were to use a frame of $E_\qat$ that did not lie in $E$, then $J$ would act as a different unit imaginary quaternion, or even as different unit imaginary quaternions on different components of $\xi_l$.
In the above definition we are implicitly using the frame of coordinate basis vectors coming from the local trivialization of $E$.
Unless otherwise noted, it should be assumed that we are working with such a trivialization.

The transformations formulas for sections of other bundles from~\eqref{eq:section transformation}.
A section $\V$ of $\Hom_\qat(E_\qat,E'_\qat)$ transforms according to 
\begin{align}
\label{eq:hom transformation}
\V_m = f'_{ml} \V_l f_{ml}^{-1}.
\end{align}
with $f_{ml}$ the cocycle of $E$ and $f_{ml}'$ the cocycle of $E'$.
Lemma~\ref{lem:quaternionic dual space} tells us that the cocycle of $E_\qat^{-1}$ is $f_{ml}^{-1}$.
But be mindful of the point in the previous paragraph: this is only true if we use the trivialization of the underlying $\mathbb{C}$\=/bundle.
In the case of $E_\qat^{-1}$ this is the $\qj$\=/dual basis with respect to the trivialization of $E$.
We can verify that the pairing of a section $\chi$ of $E_\qat^{-1}$ and $\xi$ of $E_\qat$ results in a $\qat$\=/valued function on $\X$:
\begin{gather*}
\lz \chi, \xi \rz_m
= \Bar{\chi}_m \qj \xi_m
= \overline{f_{ml}^{-1}\chi_l} \qj f_{ml}\xi_l
= \Bar{\chi}_l \Bar{f}_{ml}^{-1}\qj f_{ml}\xi_l
= \Bar{\chi}_l \qj \xi_l
= \lz \chi, \xi \rz_l.
\end{gather*}

\section*{Differential Forms}

Next we wish to explain differential forms in the quaternionic setting.
A $\qat$\=/valued $p$\=/form is a section of the bundle $\bigwedge^p T^\ast\X \otimes_\mathbb{R} \qat$, though since we are working over a Riemann surface $\X$ there is only the zero section for $p > 2$.
The wedge product extends as it does for other associative algebras: 
\[
(\omega \otimes \alpha) \wedge (\eta \otimes \beta)
= (\omega \wedge \eta) \otimes (\alpha \beta).
\]
The wedge product with a $0$\=/form reduces to scalar multiplication, in which case we often write it as a juxtaposition.
Note that unlike the real and complex cases, the wedge product in the quaternionic case is not anti-symmetric.
The exterior derivative for functions is defined in any chart $z_l = x_l + \qi y_l : \SO_l \to \mathbb{C}$ by
\[
d\alpha = dx_l \otimes \frac{\partial \alpha}{\partial x_l} + dy_l \otimes \frac{\partial \alpha}{\partial y_l}.
\]
In particular we have the product rule
\begin{align}
\label{eqn:d product rule}
d(\alpha\beta) = d\alpha \wedge \beta + \alpha \wedge d\beta.
\end{align}
We define the exterior derivative for higher degree forms inductively by the properties $d^2 = 0$ and 
\[
d(\omega \wedge \eta) = d\omega \wedge \eta + (-1)^{\deg \omega} \omega \wedge d\eta.
\]


Because we have the complex scalars $\mathbb{C} = \qat^+$ sitting inside the quaternions, $\qat$\=/valued forms contain the usual complexified differential forms:
\begin{align*}
T^\ast\X \otimes_\mathbb{R} \qat
= T^\ast\X \otimes_\mathbb{R} \mathbb{C} \otimes_\mathbb{C} \qat
= \left( K \oplus \Bar{K} \right) \otimes_\mathbb{C} \qat,
\end{align*}
using $K$ and $\Bar{K}$ for the $\mathbb{C}$\=/line bundle of complex differential forms of type $(1,0)$ and $(0,1)$ respectively.
Observe that if we define a function $z_l : \SO_l \to \qat$ by $x_l + \qi y_l \mapsto x_l + y_l \qi$, where the $\qi$ of the input is the complex unit of the chart and the $\qi$ of the output is the unit imaginary quaternion, then 
\begin{align*}
dz_l &= dx_l \otimes 1 + dy_l \otimes \qi, \\
d\Bar{z}_l &= dx_l \otimes 1 - dy_l \otimes \qi,
\end{align*}
are a nowhere vanishing $1$\=/forms on $\SO_l$ that locally span $K\otimes \qat$ and $\Bar{K}\otimes \qat$ respectively.
If $dz_l$ or $d\Bar{z}_l$ act on any vector field of $X$ the result is valued in $\qat^+$.
This leads us to treat these forms as if they we complex numbers in calculations, such as writing expressions such as $\qj dz_l = d\Bar{z}_l \qj$ and $\qj d\Bar{z}_l = dz_l \qj$.
We can also pick out the type of a form using Hodge star, the natural complex structure of the cotangent bundle of a Riemann surface.
Our convention is $\ast dx_l = dy_l$ and $\ast dy_l = -dx_l$ (as in~\cite[Eq~I.3.7.1]{FK}).
One should be aware that our convention of $\ast$ differs by a sign from the convention of the Berlin group.
In particular, this leads to a negative sign in the definition of conformality, Equation~\eqref{eq:left normal} below, in order to match their definition.
The $(1,0)$\=/forms are sections such that $-\ast \omega = \omega\qi$ and the $(0,1)$\=/forms are sections such that $-\ast \omega = -\omega\qi$.
This can be understood as an application of Lemma~\ref{lem:quaternionification}.

We can split the exterior derivative $d = \partial + \barpartial$ by type.
For example, we can calculate for a function
\begin{align*}
\barpartial \alpha
&= \frac{1}{2} (d\alpha - \qi \ast d\alpha)
= \frac{1}{2} \left( dx_l \otimes \frac{\partial \alpha}{\partial x_l} + dy_l \otimes \frac{\partial \alpha}{\partial y_l} - dy_l \otimes \qi \frac{\partial \alpha}{\partial x_l} + dx_l \otimes \qi \frac{\partial \alpha}{\partial y_l} \right) \\
&= d\Bar{z}_l \otimes \frac{1}{2} \left( \frac{\partial \alpha}{\partial x_l} + \qi \frac{\partial \alpha}{\partial y_l} \right)
=: d\Bar{z}_l \otimes \barpartial_l \alpha, \\
\partial \alpha
&= dz_l \otimes \frac{1}{2} \left( \frac{\partial \alpha}{\partial x_l} - \qi \frac{\partial \alpha}{\partial y_l} \right)
=: dz_l \otimes \partial_l \alpha.
\end{align*}
Here $\partial_l$ and $\barpartial_l$ are the natural extension of Wirtinger derivatives to $\qat$\=/valued functions.
The subscript $l$ refers to the chart $z_l$.
These derivatives inherit a product rule:
\begin{align*}
\barpartial(\alpha\beta)
&= \frac{1}{2} (d(\alpha\beta) - \qi \ast d(\alpha\beta))
= \frac{1}{2} ((d\alpha \wedge \beta + \alpha \wedge d\beta) - \qi \ast(d\alpha \wedge \beta + \alpha \wedge d\beta)) \\
\labelthis{eq:barpartial product rule}
&= \barpartial \alpha \wedge \beta + \frac{1}{2} (\alpha \wedge d\beta  - \qi \alpha \wedge \ast d\beta),
\end{align*}
and a similar formula for $\partial(\alpha\beta)$.
Because the $\qi$ does not necessarily commute with $\alpha$, it it not possible to factor $\alpha$ from the second term.
However in the special case $\alpha = \qj$ it is possible, giving $\barpartial(\qj \beta) = \qj \partial \beta$ and $\partial(\qj \beta) = \qj \barpartial \beta$.

\begin{example}[Catenoid]
\label{eg:catenoid derivative}
\index{Catenoid!Differential}
Let us give the detail promised in Example~\ref{eg:catenoid}.
Consider $F : \mathbb{C}^\ast \to \qat$ given by
\[
F(z) = \ln(z\Bar{z}) \qi + (\Bar{z} + z^{-1})\qj.
\]
Applying the above differentiation rules, we compute $dF$
\begin{align*}
dF
&= (dz\, z^{-1} + d\Bar{z}\, \Bar{z}^{-1})\qi + (d\Bar{z} - dz\, z^{-2}) \qj 
= dz (z^{-1} \qi - z^{-2} \qj) + d\Bar{z} ( \Bar{z}^{-1}\qi + \qj) \\
&= dz\, z^{-1} \qi (1 + z^{-1} \qk) + d\Bar{z}\, \qj ( z^{-1}\qk + 1) 
= (1-z^{-1}\qk) \qj dz (1+z^{-1}\qk).
\end{align*}
This last form is an example of treating $dz$ as a complex value, which makes the factorization possible.
\end{example}

We also introduce differential forms valued in a complex $\qat$\=/bundle $E_\qat$.
These are sections of 
\[
\left(\bigwedge^p T^\ast\X \right) \otimes_\mathbb{R} E \otimes_\mathbb{C} \qat
= \bigwedge^p (T^\ast\X \otimes_\mathbb{R}\mathbb{C}) \otimes_\mathbb{C} E \otimes_\mathbb{C} \qat
= \bigwedge^p (K \oplus \Bar{K}) \otimes_\mathbb{C} E \otimes_\mathbb{C} \qat.
\]
In the special case $E = \mathbb{C}$, the trivial $\mathbb{C}$\=/bundle on $\X$, this definition reduces to that of $\qat$\=/valued forms.
A word of caution: above in the tensor product between the complexified differentials forms and $E$, the complex action on $E$ is the complex structure $J$ and not the scaling action $\qi$.
For this reason, one sometimes sees the use of $J$ as the complex unit of the complexified differential forms, i.e.\ to write $dy\, J \otimes e = dy \otimes Je$ instead of $dy\,\qi \otimes e = dy \otimes Je$.
In the case that we are working in a trivialization of $E$, then $J$ acts as $\qi$ and it is not as necessary to distinguish the multiple complex structures.

Unlike in the complex case, we do not have a tensor product of $\qat$\=/vector spaces, so there is no general analogue of taking the wedge product of an $E$\=/valued form and an $E'$\=/valued form to create an $E\otimes E'$\=/valued form.
However, the wedge product between $E_\qat$\=/valued forms and $\qat$\=/valued forms is defined
\[
(\omega \otimes e \otimes \alpha) \wedge (\eta \otimes \beta)
= (\omega \wedge \eta) \otimes e \otimes (\alpha \beta).
\]
Note that the wedge product of $\qat$\=/valued forms and $E_\qat$\=/valued forms is not defined (the reverse order of the above), because there is no way for $\beta$ to ``pass-though'' $e$ without running into the same problem as~\eqref{eqn:tensor action definition problem}.
Alternatively, we can understand this as the fact that $E_\qat$ does not have a left $\qat$\=/action.

\begin{definition}
\label{def:forms notation}
\index{Differential form}
Let the sheaf of differential forms of degree $p$ valued in $E_\qat$ be denoted $\forms{p}{E}$.
Similarly let the sheaf of differential forms of type $(q,r)$ valued in $E_\qat$ be denoted $\forms{q,r}{E}$.
For $\qat$\=/valued forms we use $\forms{p}{}$ or $\forms{q,r}{}$, omitting the trivial bundle.
\end{definition}

\section*{Connections}

There is not a single notion of a derivative of sections of a line bundle, rather one must specify a connection.
A connection $\nabla$ in the quaternionic setting is a $\qat$\=/linear map from the sections of a $\qat$\=/bundle $M$ to $M$\=/valued $1$\=/forms that obeys the product rule
\begin{align}
\label{eq:connection product rule}
\nabla(\xi\alpha) = (\nabla\xi) \alpha + \xi \wedge d\alpha,
\end{align}
for sections $\xi$ and scalar functions $\alpha : \X \to \qat$.
We see that this equation reduces to Equation~\eqref{eqn:d product rule} in the case that $M$ is the trivial bundle and $\nabla = d$.
If $M = E_\qat$ is a complex $\qat$\=/bundle we can decompose $\nabla$. 
Indeed the $E_\qat$\=/valued $1$\=/forms have two complex structures $J$ and $\ast$.
We can decompose the connection into parts valued in the $(1,0)$- and $(0,1)$\=/forms 
\begin{align}
\label{eq:connection split}
\nabla' &= \tfrac{1}{2}(\nabla + J\ast \nabla) &
\nabla'' &= \tfrac{1}{2}(\nabla - J\ast \nabla),
\end{align}
and further by commutativity with the complex structure $J$
\begin{gather}
\label{eq:connection split 2}
\begin{aligned}
(\nabla')^+\xi &= \tfrac{1}{2}(\nabla'\xi - J \nabla' (J\xi)) &
(\nabla')^-\xi &= \tfrac{1}{2}(\nabla'\xi + J \nabla' (J\xi)), \\
(\nabla'')^+\xi &= \tfrac{1}{2}(\nabla''\xi - J \nabla'' (J\xi)) &
(\nabla'')^-\xi &= \tfrac{1}{2}(\nabla''\xi + J \nabla'' (J\xi)).
\end{aligned}
\end{gather}

Every trivialization induces locally a flat connection that acts like $d$ on the component functions of the sections.
So long as we are working with a trivialization of $E$ as in Definition~\ref{def:sections of qat bundle}, this enables us to write
\begin{align*}
\index{Potential}
(\nabla \xi)_l
&= d \xi_l + dz_l \U_l \xi_l - d\Bar{z}_l \V_l \xi_l \\
&= \partial_l \xi_l + dz_l(\U_l^+ + \U_l^-)\xi_l + \barpartial_l \xi_l - d\Bar{z}_l(\V_l^+ + \V_l^-)\xi_l \\
&= dz_l (\partial_l + \U_l^+ + \U_l^-)\xi_l + d\Bar{z}_l(\barpartial_l - \V_l^+ - \V_l^-)\xi_l.
\end{align*}
The signs for $\U$ and $\V$ are arbitrary, and so we already use the sign convention that is most convenient for the rest of the book.
On overlaps $\SO_l \cap \SO_m$ the transformation rule for sections and forms applied to $\nabla \xi$, for example
\begin{align*}
dz_m (\partial_m + \U_m^+ + \U_m^-) (f_{ml} \xi_l)
&= f_{ml} dz_l (\partial_l + \U_l^+ + \U_l^-)\xi_l,
\end{align*}
leads to
\begin{gather}
\label{eq:potential transformation}
\begin{aligned}
dz_m \U_m^+ &= dz_l (\U_l^+ - f_{ml}^{-1} \partial_l f_{ml} )
\quad\text{and}\quad
dz_m \U_m^- &= dz_l f_{ml} \U_l^- f_{ml}^{-1}, \\
d\Bar{z}_m \V_m^+ &= d\Bar{z}_l (\V_l^+ + f_{ml}^{-1}\barpartial_l f_{ml}) 
\quad\text{and}\quad
d\Bar{z}_m \V_m^- &= d\Bar{z}_l f_{ml} \V_l^- f_{ml}^{-1}.
\end{aligned}
\end{gather}
The commuting part of a $\qat$\=/linear connection restricts to give a $\mathbb{C}$\=/linear connection on the underlying bundle $E$.
On the other hand, the anti-commuting parts $U^-$ and $V^-$ of a $\qat$\=/linear connection are tensorial.
We see from the above transformation formula that they are $\End_\qat(E_\qat)$\=/valued $(1,0)$- and $(0,1)$\=/forms.
In Definition~5.2 and Section~5.3 of~\cite{BFLPP} these are called Hopf fields.
Instead we call any endomorphism valued form a \emph{potential}.

The $(0,1)$\=/part of a connection $\nabla''$ defines a \emph{holomorphic structure}\index{Holomorphic structure} on a complex $\qat$\=/bundle.
Such differential operators $D$ are characterized as $\qat$\=/linear maps from the sections of $E_\qat$ to $E_\qat$\=/valued $(0,1)$\=/forms with the product rule
\begin{align}
\label{eq:holo product rule}
D(\xi \wedge \alpha)
= D(\xi) \wedge \alpha+\tfrac{1}{2}(\xi \wedge d\alpha-J\xi \wedge\ast d\alpha)
\end{align}
which is the analogue of Equation~\eqref{eq:barpartial product rule} and matches~\cite[Definition~4]{BFLPP}.
This same rule extends the operator to higher degree forms, see~\eqref{eq:connection applied to forms} below.

Now suppose, as will be the principal case throughout this book, that $E$ is a holomorphic $\mathbb{C}$\=/bundle.
This means that the trivializations can be chosen such that cocycle $f_{ml}$ consists of holomorphic functions.
Equivalently~\cite[Theorem~2.6.26]{Hu}, there is a $\mathbb{C}$\=/linear holomorphic structure $\delbar{E}$ on $E$.
The operator can be defined in the following way: in any chart $\SO_l$ take $\barpartial_l$ of the coefficients of $\xi$.
This is a well-defined $E$\=/valued $(0,1)$\=/form since
\[
d\Bar{z}_m \barpartial_m\xi_m
= d\Bar{z}_l \barpartial_l(f_{ml}\xi_l)
= d\Bar{z}_l f_{ml} \barpartial_l \xi_l.
\]
Conversely, holomorphic sections are those in the kernel of the holomorphic structure, and a frame of holomorphic sections will give a holomorphic trivialization.
If we consider the implications for a holomorphic structure on $E_\qat$, we see that transformation rule~\eqref{eq:potential transformation} simplifies with $\barpartial_l f_{ml}= 0$ to
\begin{equation}
\label{eq:potential transformation V+}
d\Bar{z}_m \V_m^+ = d\Bar{z}_l \V_l^+ .
\end{equation}
The holomorphic structure on $E$ extends $\qat$\=/linearly to a holomorphic structure on $E_\qat$, which we also denote $\delbar{E}$.
Therefore any holomorphic structure on $E_\qat$ differs from $\delbar{E}$ by a $(0,1)$\=/potential $V$.
If the commuting part $V^+$ vanishes we say that it agrees with the underlying holomorphic structure.
We shall see in Lemma~\ref{lem:local holomorphic structure} that every left normal induces a holomorphic $\mathbb{C}$\=/line bundle $E$ and a holomorphic structure $\delbar{E} - V^-$ on $E_{\qat}$.
Sections of $E_\qat$ lying in the kernel of a holomorphic structure $\delbar{E}-\V$ are called $\V$\=/holomorphic sections.
We will give a more precise definition in the case of weakly differentiable sections in Definition~\ref{def:holomorphic}.
The space of holomorphic sections is invariant under right multiplication with quaternions and therefore a $\qat$\=/vector space.
It does not generally inherit the complex structure, since for any $\V$\=/ holomorphic section $\xi$ we have
\[
(\delbar{E}-\V)(J\xi)
= J (\delbar{E}-\V^+)\xi - \V^- J\xi
= J \V^- \xi - \V^- J\xi,
\]
which vanishes if and only if $\V^- = 0$.

On the other hand, the $(1,0)$ part of a connection, which might be called an anti-holomorphic structure, takes sections of $E_\qat$ to $E_\qat$\=/valued $(1,0)$\=/forms.
Such forms can be understood as sections of $KE_\qat$, which is a holomorphic bundle.
It will be important to us to be able to take these derivatives iteratively, generating sections of $K^2E_\qat$ and so on.
Thus we treat this in the main as a differential operator between sections of holomorphic bundles, rather than as differential forms.
If $E$ is a line bundle, then a trivialization equivalent to a non-vanishing section $\xi$.
The corresponding flat connection has the local expression $\nabla = dz_l( \partial_l - (\partial_l\xi_l)\xi_l^{-1} ) + d\Bar{z}_l( \barpartial_l - (\barpartial_l\xi_l)\xi_l^{-1} )$.
In particular, the $(1,0)$\=/part is an operator independent of the choice of chart.
We will on occasion write $\partial_l - (\partial_l\xi)\xi^{-1}$ for this global operator.

In the same way that the definition of $d$ on functions and the product rule generates the exterior derivative, so too can a connection $\nabla$ be extended to $E_\qat$\=/valued forms of higher degree.
For a section $\xi$ of $\forms{p}{E}$ and $\alpha$ of $\forms{q}{}$ we define
\begin{align*}
\nabla(\xi \wedge \alpha) = \nabla\xi \wedge \alpha + (-1)^{p} \xi \wedge d\alpha.
\end{align*}
Although this is a general definition, it reduces to Equation~\eqref{eq:connection product rule} for $0$\=/forms and the derivative of a $2$\=/form must vanish on a Riemann surface.
Therefore the only additional case is really the derivative of $1$\=/forms, for which the following formulas are more useful:
\begin{gather}
\label{eq:connection applied to forms}
\begin{aligned}
\nabla (dz_l \xi_l)
&= \nabla(\xi_l dx_l + \qi \xi_l dy_l)
= \nabla\xi_l \wedge dx_l + \nabla(\qi \xi_l) \wedge dy_l \\
&= - dz_l \wedge \nabla^+ \xi_l - d\Bar{z}_l \wedge \nabla^- \xi_l, \\ 
\nabla (d\Bar{z}_l \xi_l)
&= - d\Bar{z}_l \wedge \nabla^+ \xi_l - dz_l \wedge \nabla^- \xi_l.
\end{aligned}
\end{gather}

\begin{example}[Catenoid]
\label{eg:catenoid holomorphic trivialization}
\index{Catenoid!Potential}
For any trivialization of a bundle there is a connection $\nabla$ on $\X \times \qat$ given by $\nabla\mathbf{1} = 0$ for $\mathbf{1} : z \mapsto 1$.
We call this the trivial connection, and it is the starting point of the application of quaternionic analysis to immersions.
We give the trivial bundle a complex structure $J_N$ defined by $J_N(z,\alpha) = (z, N\alpha)$ as in Example~\ref{eg:catenoid underlying line bundle}.
We can also express these structures using the section $\mathbf{1}$.
For a section $\xi = \mathbf{1}\xi_1$ the complex structure is $J\xi = \mathbf{1}N\xi_1$ and $\nabla$ acts simply as $\nabla\xi = \mathbf{1} d\xi_1$.


The potentials are of particular interest to us.
The anti-commuting $(0,1)$\=/part of the connection is
\begin{align*}
(\nabla'')^- \xi
&= \tfrac{1}{2}(\nabla'' \xi + J_N \nabla'' (J_N\xi)) \\
&=\tfrac{1}{4}(\nabla \xi - J_N \ast \nabla \xi + J_N \nabla (J_N\xi) - J_N^2 \ast \nabla (J_N\xi)) \\
&=\tfrac{1}{4}(\nabla (\mathbf{1}\xi_1) - J_N \ast \nabla (\mathbf{1}\xi_1) + J_N \nabla (\mathbf{1}N\xi_1) + \ast \nabla (\mathbf{1}N\xi_1)) \\
&=\tfrac{1}{4}(\mathbf{1} d\xi_1 - J_N \mathbf{1} \ast d\xi_1 + J_N \mathbf{1} d(N\xi_1) + \mathbf{1} \ast d(N\xi_1)) \\
&=\mathbf{1}\, \tfrac{1}{4}(d\xi_1 - N \ast d\xi_1 + N d(N\xi_1) + \ast d(N\xi_1)) \\
&=\mathbf{1}\, \tfrac{1}{4}(d\xi_1 - N \ast d\xi_1 + N (dN)\xi_1 - d\xi_1 + (\ast dN)\xi_1 + N \ast d\xi_1) \\
&=\mathbf{1}\, \tfrac{1}{4}(N (dN) + (\ast dN)) \xi_1.
\labelthis{eq:V- from N}
\end{align*}

To be concrete, we can consider $N(z) = (1-z^{-1}\qk) (-\qi) (1-z^{-1}\qk)^{-1}$ from our ongoing example of the catenoid.
A non-vanishing section section of $E$ is $\mathbf{e} = \mathbf{1} (1-z^{-1}\qk)\qj$.
Using the definition of $\nabla$,
\begin{align*}
\nabla \mathbf{e}
&= \mathbf{1} d ((1-z^{-1}\qk)\qj) 
= \mathbf{1} (dz\,z^{-2}\qk)\qj
= \mathbf{1}(1-z^{-1}\qk)\qj (-\qj)(1-z^{-1}\qk)^{-1} dz\,z^{-2}(-\qi) \\
&= \mathbf{e} \Big[ - dz\, z^{-1} - d\Bar{z}\qk z^{-1}\Bar{z} \Big](1+|z|^2)^{-1},
\end{align*}
so that
\begin{align*}
(\nabla \xi)_l
&= d\xi_l - dz\, z^{-1}(1+|z|^2)^{-1}\xi_l - d\Bar{z}\,\qk z^{-1}\Bar{z}(1+|z|^2)^{-1}\xi_l
\end{align*}
with respect to the trivialization given by $\mathbf{e}$.
Hence
\[
\U_l^+ = - z^{-1} (1+|z|^2)^{-1}, \qquad
\U_l^- = \V_l^+ = 0, \qquad
\V_l^- = \qk z^{-1}\Bar{z}(1+|z|^2)^{-1},
\]
taking note of the definitional signs of $\U$ and $\V$.
One can check that Equation~\eqref{eq:V- from N}, which picks out the relevant component algebraically rather than by inspection, gives the same result.
In particular, because $\V_l^+ = 0$, the section $\mathbf{e}$ is a holomorphic section of $E$ with respect to the induced holomorphic structure, $\delbar{E} \mathbf{e} = (\nabla'')^+ \mathbf{e} = 0$.
Hence the bundle $E$ is holomorphically trivial.
\end{example}

\section*{Sobolev spaces and Holomorphicity}

In this book, one of our main objectives is to generalize quaternionic analysis to include holomorphic sections of $\qat$\=/line bundles with $\banach{2}\loc$\=/potentials.
This is because our main theorem deals with Willmore energy, which is the $\banach{2}$\=/norm $\|V^-\|_2$ of the anti-commuting part of a potential $V$.
However the theory is also interesting because 2 is the critical exponent for potentials in the following heuristic sense.
A section $\xi$ should be weakly differentiable, say $\sobolev{1,p}$, so that $\barpartial \xi \in \banach{p}$.
By the Sobolev embedding theorem, we also know that $\xi \in \banach{q}$ for $q^{-1} = p^{-1} - \tfrac{1}{2}$.
Then by Hölder's inequality the potential term $\V\xi$ of $\nabla''\xi = (\barpartial - \V)\xi$ also belongs to $\banach{p}$ if $\V$ belongs to $\banach{2}$.
The use of these theorems in the quaternionic setting will be justified in Lemma~\ref{lem:sobolev regularity}.

For $\mathbb{C}$\=/valued $1$\=/forms $\omega = dzU + d\Bar{z} V$ on $\Omega \subset \mathbb{C}$ there is a natural way to produce a measure~\cite[Eq~I.4.2.1]{FK}, independent of a choice of metric,
\begin{align*}
|\omega|^2 
&= - \tfrac{1}{2} \ast\overline{\omega} \wedge \omega 
= - \tfrac{1}{2} (-\Bar{U}\ci d\Bar{z} + \Bar{V}\ci dz) \wedge (dz U + d\Bar{z}V) \\
&= \Bar{U} \dmu U + \Bar{V}\dmu V
= (|U|^2 + |V|^2)\dmu,
\labelthis{eq:square integrable}
\end{align*}
with volume form $\dmu = dx \wedge dy = \frac{\ci}{2}dz\wedge d\Bar{z}$.
The space of forms has a real positive-definite inner product
\begin{align}\label{inner product}
\langle \omega,\omega'\rangle
&= \frac{1}{2}\int_\Omega\Real\left(\ast\overline{\omega}\wedge \omega'\right)
=\int_\Omega\Real \left( \Bar{U}U' + \Bar{V}V' \right) \dmu 
\end{align}
(compare to the Hermitian product defined in~\cite[Eq~I.4.2.2]{FK}).
We see in this inner product that $(1,0)$- and $(0,1)$\=/forms are orthogonal.
Not only does this same formula extend to $\qat$\=/valued $1$\=/forms, it is well-defined for $\End_\qat(E_\qat)$\=/valued $1$\=/forms in the case that $E_\qat$ is a line bundle.
For example, for the inner product of two $(0,1)$\=/potentials the integrand is indeed a $2$\=/form:
\begin{gather}\begin{aligned}
\label{eq:inner product}
\index{Pairing!Potentials}
\Real\left(\Bar{\V}_m \V'_m\right)\dmu_m
&=\Real\left(\Bar{f}_{ml}^{-1}\Bar{\V}_l\overline{\tfrac{d\Bar{z}_l}{d\Bar{z_m}}f_{ml}}f_{ml}\tfrac{d\Bar{z}_l}{d\Bar{z}_m}\V_lf_{ml}^{-1}\right) \bigl|\tfrac{dz_m}{dz_l}\bigr|^2 \dmu_l \\
&=\Real\left(\Bar{f}_{ml}^{-1}\Bar{\V}_l\big|f_{ml}\big|^2 \V'_lf_{ml}^{-1}\right)\dmu_l
=\Real\left(\Bar{\V}_l \V'_l\right) \dmu_l.
\end{aligned}\end{gather}
The final equality uses the cyclic property of $\Real$.
The calculation for $(1,0)$\=/potentials almost identical.
The inner product induces a norm $\|\V\|_2$
\begin{gather}\label{eq:2-norm}
\|\V\|_2
=\left(\langle \V,\V\rangle\right)^{1/2}
=\left(\int_\X |\V|^2\right)^{1/2}.
\end{gather}
Note the use of single bars to indicate the corresponding measure $|\V|^2_l = |\V_l|^2 \dmu_l$, which is the $\banach{2}$\=/norm of the local representative multiplied by the Lebesgue measure of the coordinate chart.

\begin{definition}
\label{def:potentials}
\index{Holomorphic structure}
For a holomorphic $\mathbb{C}$\=/line bundle $E$ on a Riemann surface $\X$, we call the potentials (locally) \emph{square-integrable} if they are (locally) bounded with respect to $\|\cdot\|_2$.
Let $\pot{E}$ denote the space of locally square-integrable $(0,1)$\=/potentials acting on $E_\qat$.
We do not have a symbol for $(1,0)$\=/potentials.
We decompose $\pot{E}=\pot{E}^+\oplus\pot{E}^-$ into potentials which commute and anti-commute with the complex structure.
For any $V\in\pot{E}$, we use the pair $(E,V)$ to refer to the complex $\qat$\=/line bundle $E_\qat$ with holomorphic structure given by $\delbar{E} - V$.
\end{definition}

Such an inner product and norm could have been expected for the endomorphisms of a $\mathbb{C}$\=/line bundle, because they are canonically isomorphic to the trivial $\mathbb{C}$\=/bundle, but is somewhat surprising in the quaternionic case where we do not have such an identification.
But although the identification of the endomorphisms of $\qat$\=/line bundle with $\qat$ depends on a choice of frame, any two identifications differ by a rotation $\V_m = f_{ml} \V_l f_{ml}^{-1}$.
This property means that it is possible to define a frame independent inner product.
The situation in higher rank is revisited in Chapter~\ref{ch:Flat families}.

Now we turn our attention to defining weakly differentiable sections of line bundles, those that belong to a Sobolev space.
These are kind of sections that are suitable for a definition of holomorphic with respect to a square integral potential.
It is relatively easy to extend the idea of a Schwartz distribution to manifolds.
On occasion in this book we will make use of \emph{currents}, which ``can be considered as a differential form for which the coefficients are distributions''~\cite{dR}.
The normal operations of forms, such as $\wedge$ and $d$ are defined on currents~\cite[pp.~36,46]{dR}\cite[p.~369]{GrHa}.
What is more difficult to define are Lebesgue and Sobolev spaces.
Many works that define Sobolev spaces for manifolds are concerned with Riemannian manifolds, using the metric and the Levi-Civita covariant derivative to give natural generalizations of the Sobolev norms.
A standard reference is~\cite[Chapter~2]{Au} and the counterexamples in~\cite{Heb} are also helpful to understand the issues involved in the definition.
However as we are working with a Riemann surface $\X$, there is not a distinguished metric for us to use.
A survey of possible approaches can be found in~\cite[Remark~7.8]{BB}.
Theorem~10.2.36 in~\cite{Ni} gives a clear statement on compact manifolds that resulting Sobolev spaces are independent of the choices: the set of sections is the independent of the choices, and the identity map is an isomorphism of Banach spaces.
The definitions we give below follow~\cite[Theorem~2.6.2]{Ho}.
We have already adopted in Remark~\ref{rem:special cover} the suggestion of~\cite[p.~82]{Jo} to use a special subatlas in order to simplify the check of chart independence.

\begin{definition}
\label{def:sobolev norm}
\index{Sobolev space}
Consider a cover $\{\SO_l\}$ of $\X$ as in Remark~\ref{rem:special cover}.
We identify functions on $\SO_l$ with functions on $z_l[\SO_l]$ using the chart $z_l$, endows them with the norms
\[
\|f\|_{\sobolev{k,p}(\SO_l,\qat)}
= \sum_{n=0}^k \sum_{m=0}^n \|\partial_l^m\barpartial_l^{n-m}f\|_{\banach{p}(\SO_l,\qat)},
\]
which is equivalent to the usual Sobolev norms.
We say that a section $\xi$ of $E_\qat$ over $U \subset \X$ belongs to $H^0(U,\sob{k,p}{E})$ if $\xi_l$ belongs to $\sobolev{k,p}(U \cap \SO_l,\qat)$ for all $l$.
This defines the sheaves $\sob{k,p}{E}$.
The zeroth order case has the special notation $\ban{p}{E} = \sob{0,p}{E}$.
For differential forms whose representatives belong to the Sobolev spaces we use the notation $\forms{q,r}{}\sob{k,p}{E}$.
\end{definition}

There are two questions arising from this definition, which are dealt with in the lemma below.
The first is the dependence of the sheaf $\sob{k,p}{E}$ on the cover.
The second is that this definition defines does not give a norm on the space of sections.
On a compact Riemann surface $\X$ it is useful to have a global norm.
It is possible to construct a global norm for non-compact Riemann surfaces, but it depends on the choice of charts (as discussed in the listed sources).
However we found that sections on non-compact surfaces with finite global norm were not the appropriate notion for our purposes and so we make no use of them.

\begin{lemma}
\label{lem:locally sobolev welldefined}
The definition of $\sob{k,p}{E}$ is independent of the choice of cover.
If $\X$ is compact then $H^0(\X,\sob{k,p}{E})$ is closed subspace of $\bigoplus_{l=1}^L \sobolev{k,p}(\SO_l,\qat)$, and so is endowed with a global norm making it a Banach space.
and $\pot{E}$ is a Hilbert space with real inner product~\eqref{inner product}.
The equivalence class of these norms are independent of the choice of cover.
\end{lemma}
\begin{proof}
We prepare two rules of calculus.
If $f$ is a smooth function then by the product rule
\[
\|f \xi\|_{\sobolev{k,p}}
\leq C(k)\|f\|_{\sobolev{k,\infty}} \|\xi\|_{\sobolev{k,p}},
\]
where $C(k)$ is a constant depending only on $k$.
A crude estimate for $C(k)$ is $k 2^{k+1}+1$, which comes from bounding every derivative in the product rule by either $\|f\|_{\sobolev{k,\infty}}$ or $\|\xi\|_{\sobolev{k,p}}$.
Similarly, from Faà di Bruno's formula for the higher order chain rule, we have
\[
\|\partial_0^m \barpartial_0^{n} \xi \|_{\banach{p}}
\leq (1+\|z_l\|_{\sobolev{m+n,\infty}})^{m+n} \sum_{\pi \in \text{Part}(m)} \sum_{\pi' \in \text{Part}(n)} \|\partial_l^{|\pi|} \barpartial_l^{|\pi'|} \xi \|_{\banach{p}},
\]
where $|\pi|$ is the number of sets in $\pi$, a partition of $\{1,\ldots,m\}$.

The union of two covers of the type considered in Remark~\ref{rem:special cover} is again such a cover.
Because adding sets makes the condition to belong to $\sob{k,p}{E}$ more restrictive, it suffices to prove that adding sets to a cover does not diminish this set of sections.
Moreover, each additional set $\SO_0$ is relatively compact, and so can be covered by finitely many sets $\SO_l$ of the original cover (arrange for the indices to be $l = 1, \dots, L$).
For the sake of brevity, below we shorten $\SU \cap \SO_0$ to $\SU_0$ and $\SU \cap \SO_0 \cap \SO_l$ to $\SU_l$, and omit $\qat$ from the Sobolev spaces $\sobolev{k,p}(\SO,\qat)$.
Then the norm of the local representative $\xi_0$ of $\xi$ on $\SO_0$ can be estimated by
\begin{align*}
\|&\xi_0\|_{\sobolev{k,p}(\SU_0)}
\leq \sum_{l=1}^L \|f_{0l}\xi_l\|_{\sobolev{k,p}(\SU_l)} 
\leq \sum_{l=1}^L C(k) \|f_{0l}\|_{\sobolev{k,\infty}(\SU_l)} \|\xi_l \|_{\sobolev{k,p}(\SU_l)}.
\end{align*}
The functions $f_{0l}$ as well as all their derivatives are bounded functions on $\SU_l$ due to Remark~\ref{rem:special cover}.
So it suffices to control $\|\xi_l\|_{\sobolev{k,p}(\SU_l,z_0)}$ in terms of $\|\xi_l\|_{\sobolev{k,p}(\SO_l,z_l)}$, using $z_l$ to make the coordinate chart explicit.
Changing the chart changes the measure
\[
\|\xi_l \|_{\banach{p}(\SU_l,z_0)}^p
= \int_{\SU_l} |\xi_l|^p\, \dmu_0
\leq \sup_{\SU_l} \left| \frac{\dmu(z_0)}{\dmu(z_l)} \right| \int_{\SU_l} |\xi_l|^p\, \dmu_l.
\]
Together with the chain rule, $\|\partial_0^m\barpartial_0^n\xi_l\|_{\banach{p}(\SU_l,z_0)}$ is bounded by 
\[
\sup_{\SU_l} \left| \frac{\dmu(z_0)}{\dmu(z_l)} \right|^{1/p}
(1+\|z_l\|_{\sobolev{m+n,\infty}})^{m+n} B_m B_n \| \xi \|_{\sobolev{m+n,p}(\SU_l,z_l)},
\]
where $B_m$ is the number of partitions of $\{1,\ldots,m\}$ (the Bell number).
Again by the nature of the charts the norms of the coordinate change and change of measure are finite.
Finally we have assumed that $\xi_l \in \sobolev{k,p}(U'\cap\SO_l)$ for all $U'$, in particular for $U' = U \cap \SO_0$.
Increasing the order of the Sobolev norm only increases it, hence we have a bound on $\|\xi_l\|_{\sobolev{k,p}(\SU_l,z_0)}$ in terms of $\|\xi_l\|_{\sobolev{k,p}(\SU_l,z_l)}$, as required.

In the case that $\X$ is compact, the cover is finite.
The above estimates now yield an equivalence between norms defined using different covers.
\end{proof}

Since the sets $\SO_l$ are relatively compact in $\X$, one should think of these as sheaves of locally $\sobolev{k,p}$ sections.
This definition also applies to functions on a Riemann surface by considering them as sections of the trivial bundle, although we prefer the notation $\sobolev{k,p}\loc(\X,\qat)$ for this case.
For forms valued in one of these sheaves the notions of a connection and the exterior derivative extend by replacing strong derivatives with their weak versions.
As we have now successfully defined weakly differentiable section, we are ready for the central definition of this work.

\begin{definition}
\label{def:holomorphic}
For fixed $1<p<2$ and for $\V\in\pot{E}$ let $\Q{E,\V}$ denote the kernel of the holomorphic structure $\delbar{E} - \V$ as a subsheaf of $\sob{1,p}{E}$.
On any open subset $\SO\subset\X$ elements of $H^0(\SO,\Q{E,W})$ are called \emph{$\V$\=/holomorphic}.
\end{definition}

In particular, the statement $\xi$\ is $\V$\=/holomorphic, entails the assumption $\xi \in H^0(\SO,\sob{1,p}{E})$ for some $1<p<2$.
We shall see in Theorem~\ref{cauchy formula} that this statement does not depend on the choice of $1<p<2$ and is equivalent to several weaker and stronger statements.

Local regularity results, such as Hölder's inequality, carry over to weakly differentiable sections.
Further, because the sets $z_l[\SO_l]$ are balls, they are bounded and have smooth boundaries, the Sobolev embedding theorem generalizes easily to global sections.
The following lemma summarizes some of the useful results that carry over.

\begin{lemma}
\label{lem:sobolev regularity}
\index{Hölder's inequality}
\index{Sobolev embedding}
\phantom{boo}
\begin{enumeratethm}
\item
\textup{Hölder inequality:}
For $1<p<2$ the potentials $\V\in\pot{E}$ act as continuous linear operators from $H^0(\X,\ban{2p/(2-p)}{E})$ into $H^0(\X,\forms{0,1}{}\ban{p}{E})$.
Compare~\cite[Theorem~III.1~(c)]{RS1}.

\item
\textup{Sobolev Embedding theorem:}
For any $1<p<2$ and $1 \leq q \leq \frac{2p}{2-p}$ the embedding $H^0(\X,\sob{1,p}{E})\hookrightarrow H^0(\X,\ban{q}{E})$ is continuous. 
Compare~\cite[4.12~Theorem]{Ad}.

\item
\textup{Rellich-Kondrachov theorem:} For any $1<p<2$ and $1\leq q<\frac{2p}{2-p}$ the above embedding is moreover compact. 
Compare~\cite[6.3~Theorem]{Ad}.
\end{enumeratethm}
\end{lemma}

Weakly conformal maps were defined in the introduction for subsets of the complex plane.
We now extend the notion to Riemann surfaces.
Since the property is local, most of the difficulty was already overcome in the definition of a weakly differentiable function on a Riemann surface.
As we saw in the proof of Lemma~\ref{lem:locally sobolev welldefined} that our choice of charts means that `measurable' and `almost everywhere' are well-defined.
In particular if $dF$ has almost everywhere rank two, then $N$ in Equation~\eqref{eq:left normal} is unique (as a function defined almost everywhere).
We take Equations~\eqref{eq:left normal} and~\eqref{eq:right normal} as the definition of conformality.

\begin{definition}
\label{def:weakly conformal global}
\index{Weakly conformal map}
\index{Admissible map}
A map $F\in\sobolev{1,1}\loc(\X,\qat)$ on a Riemann surface $\X$ is called \emph{weakly conformal} if $dF$ has almost everywhere rank two and there exists a measurable map $N:\X\to\qat$ that is almost everywhere a left normal,
\begin{align}\label{eq:left normal}
-\ast dF &= N dF
&
N^2&=-1,
\end{align}
or equivalently, if there exists a measurable map $R:\X\to\qat$ that is almost everywhere a right normal,
\begin{align}\label{eq:right normal}
-\ast dF &= dFR
&
R^2&=-1.
\end{align}
If in addition $N\in\sobolev{1,2}\loc(\X,\qat)$ then the map $F$ is called \emph{admissible}.
\end{definition}

\begin{example}[Plane]
\label{eg:plane admissible}
\index{Plane!Admissibility}
A simple case is an admissible map to the complex plane $F : \X \to \mathbb{C}$.
Around any point, take a coordinate $z : \SO \subset \X \to \mathbb{C}$.
We have
\begin{align*}
dF &= dz \partial F + d\Bar{z}\barpartial F \\
-\ast dF 
&=  dz \qi \partial F - d\Bar{z} \qi\barpartial F \\
N dF
&= dz( N^+ \partial F + N^- \barpartial F ) + d\Bar{z}( N^- \partial F + N^+ \barpartial F ).
\end{align*}
By consideration of type and commutativity with $\qi$, we can say that $F$ is weakly conformal on $\SO$ only if
\begin{align*}
\qi \partial F = N^+ \partial F, \quad
0 = N^- \barpartial F, \quad
0 = N^- \partial F, \quad
\text{and}\quad
-\qi \barpartial F = N^+ \barpartial F.
\end{align*}
If neither $\partial F$ nor $\barpartial F$ vanish, then these equations imply $\qi = N^+$ and $-\qi = N^+$, which is impossible.
If both $\partial F$ and $\barpartial F$ vanish, then so does $dF$ and it is not full rank (i.e.\ constant maps are excluded).
Therefore exactly one vanishes: either $F$ or $\Bar{F}$ is holomorphic on $\SO$.
The four equations now imply $N = \qi$ and $N = -\qi$ in the two cases respectively.
In either case then $F$ is analytic.
\end{example}

To finish the chapter, let us consider the integration of $1$\=/forms in $\forms{1}{}\sob{1,p}{E}$.
It is not obvious that the definition of integration in terms of local charts carries over to this case, because the path of integration, being a curve, is a set of measure zero.
We offer two practical approaches.
The first uses additional regularity of the path.
Given a continuously differentiable immersed path $\gamma: (0,1) \to\X$ we can choose charts of $\X$ that map $\gamma$ locally to straight lines in $\Omega\subset\mathbb{R}^2$.
By the Sobolev Embedding theorem~\cite[4.12~Theorem Case C]{Ad} for any $1<p<2$ the restriction of a function in $\sobolev{1,p}\loc(\Omega,\qat)$ to a straight line $\ell$ in $\Omega\subset\mathbb{R}^2$ belongs to $\banach{\frac{p}{2-p}}\loc(\ell)$.
Thus the integral along the line segment is well-defined.
The second approach is specific to closed forms and closed paths, but holds under weaker regularity.
For any continuous cycle $\gamma$, we choose a smooth $1$\=/form $\eta_\gamma$ as in~\cite[II.3.3, p.\ 40]{FK} with the property that
\begin{align}
\label{eq:def line integral}
\int_\X \eta_{\gamma} \wedge \phi&=\int_{\gamma} \phi
\quad\text{for any closed smooth $1$\=/form $\phi$.}
\end{align}
Moreover the support of $\eta_\gamma$ can be chosen to be a tubular neighborhood of the cycle, and in particular compact.
Under this approach, we define $\int_\gamma \phi$ for $\phi$ a closed $1$\=/current to be the integral over $\X$ of $\eta_\gamma\wedge\phi$.
The integral of a closed $\qat$\=/valued $1$\=/current along a cycle only depends on the homology class of that cycle.
Therefore we have the \emph{period map}\index{Period map} of $\omega$ defined by
\begin{align}
\label{eq:def period map}
\per(\omega) \in \Hom(H_1(\X,\mathbb{Z}),\qat), \qquad
\gamma \mapsto \int_\gamma \omega.
\end{align}
The general theory of integration of currents, which goes by the name ``homological integration theory'' can be found in~\cite[Chapter~IV]{dR} and~\cite[Chapter~4]{Federer1996}.

\chapter{Local Behavior of Holomorphic Sections}
\label{chapter:local}
In this chapter, we extend several important results of complex analysis to quaternionic functions: namely Cauchy's integral formula, the notion of the order of a root, the strong unique continuation, and the isolation of the roots.
These results show that $\V$\=/holomorphic sections, even for $\banach{2}\loc$\=/potentials, have a surprising level of regularity.
The analysis depends on the use of the Sobolev embedding for the critical exponent $p=n=2$. 
As a consequence the equations have additional scaling symmetries. 
Not coincidentally the Willmore functional~\eqref{eq:willmore energy kodaira}, a multiple of the square of the $\banach{2}$\=/norm of one of the potentials, is also invariant with respect to the full conformal group of the ambient space.
In spite of the additional challenge caused by the critical exponent, with appropriate modifications most of the strong tools of complex analysis have a counterpart in this approach to quaternionic analysis.
Our decision to always work in a local trivialization in which the complex structure acts as left multiplication with $\qi$, which will be explained in Chapter~\ref{chapter:kodaira}, so that $\qat = \qat^+\oplus\qat^-$, makes the adaption of classical tools to the present setting seem natural.

These results are fundamentally local.
We shall use $\SO\subset\X$ to refer to some open subset over which the bundles are trivial, and $z:\SO\to\Omega \subset \mathbb{C}$ is the corresponding coordinate chart.
For most of this chapter we work in a single coordinate chart, in which case we will omit subscripts such as $\upsilon_l$, $\V_l$, and $\barpartial_l$. 
If, say, $\V\in \banach{2}\loc(\Omega,\qat)$ then around any point of $\Omega$ there is an open subset $\Omega'\subset\Omega$ with compact closure in $\Omega$ such that the restriction of $\V$ to $\Omega'$ belongs to $\banach{2}(\Omega',\qat)$.
We shall always assume, in accordance with Remark~\ref{rem:special cover}, that $\Omega'$ has been chosen with a smooth boundary so as to satisfy the Sobolev Embedding theorem and the Rellich-Kondrachov theorem.

We introduce several operators on Banach spaces of $\qat$\=/valued functions, in order to understand $\barpartial - \V$ and its kernel, the $\V$\=/holomorphic functions.
Dolbeault's Lemma~\cite[Chapter~I Section~D 2.~Lemma]{GuRo} implies that the operator $\Op{I}_{\mathbb{C}}$ with the integral kernel $(z-w)^{-1}\tfrac{1}{\pi} \dmu(w)$ is a right inverse of the operator $\barpartial$.
That is
\[
\Op{I}_{\mathbb{C}}(f) 
:= z \mapsto \int_{\mathbb{C}} f(w) (z-w)^{-1}\tfrac{1}{\pi} \dmu(w).
\]
Due to the Hardy-Littlewood-Sobolev theorem~\cite[Chapter~V.~\S1.2 Theorem~1]{St}, for all $1<p<2$ and $2<q<\infty$ with $\frac{1}{p}=\frac{1}{q}+\frac{1}{2}$, this is a bounded operator from  $\banach{p}(\mathbb{C},\qat)$ into $\banach{q}(\mathbb{C},\qat)$.
The Banach space of such operators is denoted by $\mathcal{L}(\banach{p}(\mathbb{C},\qat),\banach{q}(\mathbb{C},\qat))$.
For any $f\in C^\infty_0(\mathbb{C},\qat)$ the holomorphic function $\Op{I}_\mathbb{C}(\barpartial f)-f$ vanishes identically, since any harmonic function in $\banach{q}(\mathbb{C},\qat)$ has this property (compare~\cite[8.2 Corollary]{ABR}).
Moreover, the restrictions $\Op{I}_{\Omega}$ and $\Op{I}_{\Omega'}$ of $\Op{I}_{\mathbb{C}}$ to $\Omega$ and $\Omega'$, respectively, are bounded operators from $\banach{p}$ into $\banach{q}$.
By considering $C^\infty_0(\Omega,\qat)$ as a subset of $C^\infty_0(\mathbb{C},\qat)$ we obtain
\begin{align}\label{left inverse}
\Op{I}_\Omega(\barpartial f)&=f
\quad\text{for any }f\in C^\infty_0(\Omega,\qat).
\end{align}

On the other hand Hölder's inequality~\cite[Theorem~III.1~(c)]{RS1} implies that the multiplication with $\V|_{\Omega'}\in\banach{2}(\Omega',\qat)$ belongs to $\mathcal{L}(\banach{q}(\Omega',\qat),\banach{p}(\Omega',\qat))$.
Hence $\unity-\Op{I}_{\Omega'}\V|_{\Omega'}$ belongs to $\mathcal{L}(\banach{q}(\Omega',\qat)):=\mathcal{L}(\banach{q}(\Omega',\qat),\banach{q}(\Omega',\qat))$.

We shall now show that this operator also belongs to $\mathcal{L}(\sobolev{1,p}(\Omega'))$. The integral kernel of $\partial\circ\Op{I}_{\Omega'}$ is $-(w-z)^{-2}\frac{1}{\pi} \dmu(w)$.
Due to~\cite[Chapter~II \S4.2 Theorem~3]{St}, for $1<p<\infty$ the corresponding operator belongs to $\mathcal{L}(\banach{p}(\Omega',\qat))$.
Therefore $\Op{I}_{\Omega'}$ belongs for such $1<p<\infty$ to $\mathcal{L}(\banach{p}(\Omega',\qat),\sobolev{1,p}(\Omega',\qat))$.
For the $1<p<2$ and $q=\frac{2p}{2-p}$ as above $\sobolev{1,p}(\Omega',\qat)$ is embedded in $\banach{q}(\Omega',\qat)$ by the Sobolev Embedding theorem~\cite[4.12~Theorem]{Ad}.
Therefore $\unity-\Op{I}_{\Omega'}\V|_{\Omega'}$ belongs to $\mathcal{L}(\sobolev{1,p}(\Omega'))$ and maps $\ker(\barpartial-\V)$ into $\ker(\barpartial)$, whose elements are smooth functions due to Weyl's Lemma~\cite[Theorem~IX.25]{RS2}.
        
For any $\V\in\banach{2}\loc(\Omega,\qat)$, any $1<p<2$ and any point of $\Omega$ we may choose an open neighborhood $\Omega'$ of that point such that 
\begin{align}\label{eq:neumann-convergence}
\| \V\|_{\banach{2}(\Omega',\qat)}<S_p
\quad\text{with}\quad
S_p=\|\Op{I}_{\mathbb{C}}\|_{\mathcal{L}(\banach{p}(\mathbb{C},\qat), \banach{\frac{2p}{2-p}}(\mathbb{C},\qat))}^{-1}
\end{align}
holds.
Due to the bound $\|\Op{I}_{\Omega'}\|_{\mathcal{L}(\banach{p}(\Omega',\qat),\banach{\frac{2p}{2-p}}(\Omega',\qat))}\le S_p^{-1}$ the Neumann series, 
\begin{align}\label{eq:neumann}
\left(\unity-\Op{I}_{\Omega'} \V|_{\Omega'}\right)^{-1}=\sum_{l=0}^{\infty}\left(\Op{I}_{\Omega'}\V|_{\Omega'}\right)^{l},
\end{align}
converges as an operator on $\sobolev{1,p}(\Omega',\qat)$.
In this situation the operator
\begin{equation}\label{eq:resolvent}
\index{Resolvent}
\Op{I}_{\Omega',\V}:=\Op{I}_{\Omega'}\left(\unity-\V|_{\Omega'}\Op{I}_{\Omega'}\right)^{-1}=\left(\unity-\Op{I}_{\Omega'}\V|_{\Omega'}\right)^{-1}\Op{I}_{\Omega'}
\end{equation}
is a right inverse of the operator $\barpartial-\V$ on $\Omega'$. For later use we state a lemma:
\begin{lemma}
\label{lem:fredholm}
On a Riemann surface $\X$ with holomorphic line bundle $E$ there exists for any $1<p<2$ the following exact sequence of morphisms of sheaves:
\begin{align}\label{eq:exact sequence}
0
\to\Q{E,\V}
\hookrightarrow\sob{1,p}{E}
\xrightarrow{\delbar{E}-\V} \forms{0,1}{}\ban{p}{E}
\to 0.
\end{align}
\end{lemma}
\begin{proof}
The operators $\barpartial_l-\V_l:\sobolev{1,p}(\Omega_l,\qat)\to\banach{p}(\Omega_l,\qat)$ on the images of the charts $z_l:\SO_l\to\Omega_l$ induce the morphism of sheaves $\delbar{E}-\V : \sob{1,p}{E}\to\forms{0,1}{}\ban{p}{E}$.
The subsheaf $\Q{E,\V}$ of $\sob{1,p}{E}$ is by Definition~\ref{def:holomorphic} the kernel of this morphism. It is surjective, since on sufficiently small $\Omega'$ the operator $\Op{I}_{\Omega',\V}$~\eqref{eq:resolvent} is a right inverse of $\barpartial-\V:\sobolev{1,p}(\Omega',\qat)\to\banach{p}(\Omega',\qat)$.
\end{proof}
With the help of the operator $\Op{I}_{\Omega',\V}$~\eqref{eq:resolvent} we generalize Cauchy's Integral formula in the following theorem to the quaternionic case.
More precisely, we prove for $\V\in\banach{2}\loc(\Omega,\qat)$ which obey~\eqref{eq:neumann-convergence} on an open subset $\Omega'\subset\Omega$ that $\xi$ is $\V$\=/holomorphic if and only if the following condition holds:
\begin{align}\label{eq:cif}
f\xi=\Op{I}_{\Omega',\V}((\barpartial f)\xi)
\quad\text{a.e. for all } f\in C_0^\infty(\Omega',\mathbb{R}).
\end{align}
In order to explain the analogy to Cauchy's Integral formula we choose a subset $\Omega''$ with compact closure in $\Omega'$ and a function $f\in C_0^\infty(\Omega',\mathbb{R})$ which is equal to $1$\ on $\Omega''$.
Then this formula expresses the values of $\xi$ on $\Omega''$ in terms of the values of $(\barpartial f)\xi$ on $\Omega'\setminus\Omega''$.

Let us now state several conditions, which are equivalent to $\xi$ being $\V$\=/holomorphic:
\begin{theorem}
\label{cauchy formula}
\index{Cauchy's integral formula}
For open $\Omega\subset\mathbb{C}$ and $\V\in\banach{2}\loc(\Omega,\qat)$ the following conditions are equivalent:
\begin{enumeratethm}
\item $\xi$ is $\V$\=/holomorphic, i.e.\ $\xi\in\sobolev{1,p}\loc(\Omega,\qat) \subset \banach{2p/(2-p)}\loc(\Omega,\qat)$ for some $1<p<2$ and $\barpartial \xi=\V\xi$.
\item \textup{(Cauchy's Integral Formula):} $\xi\in\banach{p}\loc(\Omega,\qat)$ for some $1<p<2$ and~\eqref{eq:cif} holds for all open $\Omega'\subset\Omega$ obeying~\eqref{eq:neumann-convergence}.
\item $\xi\in\bigcap_{1<p<2}\sobolev{1,p}\loc(\Omega,\qat) \subset \bigcap_{2<q<\infty}\banach{q}\loc(\Omega,\qat)$ and $\barpartial \xi=\V\xi$.
\item $\xi\in\banach{q}\loc(\Omega,\qat)$ for some $2<q<\infty$ and the $\barpartial$\=/derivative of the distribution induced by $\xi$ is equal to the distribution induced by $\V\xi$.
\end{enumeratethm}
\end{theorem}
\begin{proof}
(i)$\Rightarrow$(ii): For $\V$\=/holomorphic $\xi$ and $f\in C_0^\infty(\Omega',\mathbb{R})$ the product $f\xi$ has compact support in $\Omega'$ and belongs to the closure of $C^\infty_0(\Omega',\qat)$ in $\sobolev{1,p}(\Omega',\qat)$ for some $1<p<2$.
Therefore $\Op{I}_{\Omega'}((\barpartial f)\xi+\V|_{\Omega'} f\xi)=\Op{I}_{\Omega'}(\barpartial(f\xi))=f\xi$ follows from~\eqref{left inverse}.
Now we insert~\eqref{eq:resolvent} and the following calculation shows~(ii):
\begin{align*}\Op{I}_{\Omega',\V}((\barpartial f)\xi)&=(\unity-\Op{I}_{\Omega'} \V|_{\Omega'})^{-1}\Op{I}_{\Omega'}((\barpartial f)\xi)\\&=(\unity-\Op{I}_{\Omega'} \V|_{\Omega'})^{-1}(f\xi-\Op{I}_{\Omega'} \V|_{\Omega'}f\xi)=f\xi.
\end{align*}
(ii)$\Rightarrow$(iii): We consider functions $f \in C_0^\infty(\Omega',\mathbb{R})$ which are constant equal to $1$ on some open subset $\Omega'' \subset \Omega'$ with compact closure in $\Omega'$.
Since $\Op{I}_{\Omega',\V}$ is a right inverse of $\barpartial-\V$ and bounded in $\mathcal{L}(\banach{p}(\Omega',\qat),\sobolev{1,p}(\Omega',\qat))$, Equation~\eqref{eq:cif} implies $\xi\in\ker(\barpartial-\V|_{\Omega''})$ on $\Omega''$ and, by the Sobolev Embedding theorem~\cite[Theorem~4.12]{Ad}, $\xi|_{\Omega''}\in\banach{q}(\Omega'',\qat)$ with $\frac{1}{q}=\frac{1}{p}-\frac{1}{2}$.
Because this is true for all $\Omega''$, it is true locally on $\Omega'$.
Since $\Op{I}_{\Omega',\V}$ maps $\banach{p}(\Omega',\qat)$ into $\sobolev{1,p}(\Omega',\qat)$ for all $1<p<2$, it maps $\banach{2}\loc(\Omega',\qat) \subset \bigcap_{1<p<2}\banach{p}\loc(\Omega',\qat)$ into $\bigcap_{1<p<2}\sobolev{1,p}\loc(\Omega',\qat)$, the latter of which is included in $\bigcap_{2<q<\infty}\banach{q}\loc(\Omega',\qat)$ by the Sobolev Embedding theorem.
Another application of~\eqref{eq:cif} with $f\in C_0^\infty(\Omega'',\mathbb{R})$ implies (iii).

\noindent(iii)$\Rightarrow$(iv) is obvious.

\noindent (iv)$\Rightarrow$(i): (iv) implies that $\barpartial \xi=\V\xi\in\banach{p}\loc(\Omega,\qat)$ with $\frac{1}{p}=\frac{1}{q}+\frac{1}{2}$.
We have seen above that the right inverse $\Op{I}_{\Omega'}$ of $\barpartial$ is bounded in $\mathcal{L}(\banach{p}(\Omega',\qat),\sobolev{1,p}(\Omega',\qat))$ for any open and bounded $\Omega'\subset\Omega$ with smooth $\partial\Omega'$.
So $\barpartial(\xi-\Op{I}_{\Omega'}(\V|_{\Omega'}\xi|_{\Omega'}))=\barpartial \xi-\V\xi=0$.
By Weyl's Lemma~\cite[Theorem~IX.25]{RS2} $\xi-\Op{I}_{\Omega'}(\V|_{\Omega'}\xi|_{\Omega'})$ is smooth on $\Omega'$ and $\xi\in\sobolev{1,p}\loc(\Omega,\qat)$.
Hence (iv) implies (i).
\end{proof}

\begin{example}
\label{eg:ln is V-holo}
Consider the function $\xi(z) = \ln |z|^2$ on $\Omega = B(0, e^{-1})$.
It belongs to $\banach{q}(\Omega,\qat)$ for all $1 \leq q < \infty$ but is not bounded at $z = 0$.
Similarly its weak derivatives,
\[
\frac{\partial\xi}{\partial x}
= \frac{2x}{x^2 + y^2} 
\quad\text{and}\quad
\frac{\partial\xi}{\partial y}
= \frac{2y}{x^2 + y^2},
\]
belong to $\banach{p}(\Omega,\mathbb{C})$ for $1 \leq p < 2$, since
\begin{align*}
\left\| \frac{\partial\xi}{\partial x} \right\|_p^p
&= \int_{\Omega} \left| \frac{2x}{x^2 + y^2} \right|^p\,\dmu
= 2^p \left( \int_0^{2\pi} |\cos \theta|^p \,d\theta \right) \left(\int_{0}^{e^{-1}} r^{1-p}\,dr \right) \\
&= \begin{cases}
2^p\, \Theta_p\, \frac{1}{2-p}e^{p-2} &\text{for } 1 \leq p < 2, \\
\infty &\text{for } 2 \leq p.
\end{cases}
\end{align*}
This shows the function has the regularity of Condition (iii) above, namely $\xi \in \bigcap_{1<p<2}\sobolev{1,p}\loc(\Omega,\qat)$.
Moreover it is a $\V$\=/holomorphic function for
\[
\V 
:= (\barpartial \xi) \xi^{-1}
= (\bar{z} \ln |z|^2 )^{-1},
\]
since this potential belongs to $\banach{2}(\Omega,\mathbb{C})$:
\begin{align*}
\|\V\|_2^2
&= 2\pi \int_0^{e^{-1}} (r \ln r^2 )^{-2} \,r\,dr
= \frac{\pi}{2} \int_0^{e^{-1}} (\ln r )^{-2} \,\frac{dr}{r} 
= \frac{\pi}{2} \int_{-\infty}^{-1}  u^{-2} \,du
= \frac{\pi}{2}.
\qedhere
\end{align*}
\end{example}

For $\V=\V^-\in\pot{E}^-$ we shall see in Chapter~\ref{chapter:kodaira} that $p$ in (i) can be also $1$, that $q$ in (iv) can be also $2$ and that $\bigcap_{1<p<2}\sobolev{1,p}\loc(\Omega,\qat) \subset \bigcap_{2<q<\infty}\banach{q}\loc(\Omega,\qat)$ in (iii) can be replaced by $\sobolev{1,2}\loc(\Omega,\qat) \cap \banach{\infty}\loc(\Omega,\qat)$.
However the above example shows that these improvements are specific to anti-commuting potentials.

We extend on open subsets $\Omega\subset\mathbb{C}$ the definition of a root's order to the quaternionic case. A complication is that unlike in the complex case, if $\xi$ is on $\Omega\setminus\{z_0\}$ $\V$\=/holomorphic then in general $(z-z_0)^{-l}\xi$ is not $\V$\=/holomorphic but $\V_{z_0,l}$\=/holomorphic with
\begin{align}\label{root potential}
\V_{z_0,l}=(z-z_0)^{-l}\V(z-z_0)^l
\quad\text{for }\V\in\banach{2}\loc(\Omega,\qat), z_0\in\Omega\text{ and }l\in\mathbb{Z}.
\end{align}  

\begin{definition}[Order of roots and poles]\label{order of roots}
\index{Root}
\index{Pole}
\index{Meromorphic function}
The \emph{order} of a root $\ord_{z_0}(\xi)$ of a nontrivial $\V$\=/holomorphic function $\xi$ on $\Omega$ is defined as the supremum of the set of all $l\in\mathbb{N}_0$ such that $(z-z_0)^{-l}\xi$ is $\V_{z_0,l}$\=/holomorphic.

A measurable function $\xi$ on $\Omega$ is called \emph{$\V$\=/meromorphic}, if for any $z_0\in\Omega$ there exists an $l\in\mathbb{Z}$, such that the restriction of $(z-z_0)^{-l}\xi$ to some neighborhood of $z_0$ is $\V_{z_0,l}$\=/holomorphic. In this case $\ord_{z_0}(\xi)$ is defined as the supremum of the set of all $l\in\mathbb{Z}$ with this property.
\end{definition}
The following characterization of removable singularities simplifies the determination of $\ord_{z_0}(\xi)$. Consequently, Theorem~\ref{strong unique continuation} will prove the finiteness of $\ord_{z_0}(\xi)$ for non-trivial $\V$\=/holomorphic and $\V$\=/meromorphic $\xi$.
\begin{lemma}\label{lem:removable singularity}
Let $2<q<\infty$ and let $\V\in\banach{2}\loc(\Omega,\qat)$ on an open neighborhood $\Omega\subset\mathbb{C}$ of $z_0$. If $\xi$ is $\V$\=/holomorphic on $\Omega\setminus\{z_0\}$ such that there exists $l\in\mathbb{Z}$ with $(z-z_0)^{-l}\xi\in\banach{q}\loc(\Omega,\qat)$, then $\ord_{z_0}(\xi)$ is the supremum of the subset of $l\in\mathbb{Z}$ with this property. This means that if this supremum $\ord_{z_0}(\xi)$ is non-negative, then $\xi$ is $\V$\=/holomorphic on $\Omega$ and has at $z_0$ a root of order $\ord_{z_0}(\xi)$. Otherwise, i.e.\ if $\ord_{z_0}(\xi)<0$, then $\xi$ is $\V$\=/meromorphic and has a pole at $z_0$ of order $-\ord_{z_0}(\xi)$.
\end{lemma}
The statement of the lemma is true independent of the choice of $2<q<\infty$.
\begin{proof}
For any integer $m$, on the open subset $\Omega\setminus\{z_0\}$, $\xi$ is $\V$\=/holomorphic if and only if $(z-z_0)^{-l}\xi$ is $\V_{z_0,l}$~\eqref{root potential}\=/holomorphic. The bounded operator $\unity-\Op{I}_\Omega\Tilde{\V}$ in $\mathcal{L}(\banach{q}\loc(\Omega,\qat))$ maps $(z-z_0)^{-l}\xi$ into the kernel of $\barpartial$.
Due to~\cite[\S21.9 9.3~Proposition]{Co2}, $z=z_0$ is a removable singularity of any member of $\banach{q}\loc(\Omega,\mathbb{H})$ which is holomorphic on $\Omega\setminus \{z_0\}$.
Hence, on all of $\Omega$, $(z-z_0)^{-\ord_{z_0}(\xi)}\xi$ is $\V_{z_0,\ord_{z_0}(\xi)}$~\eqref{root potential}\=/holomorphic and is therefore an element of $\bigcap_{2<q<\infty}\banach{q}\loc(\Omega,\qat)$ by the previous theorem.
\end{proof}	
If $(z-z_0)^{-l}\xi$ belongs to $\banach{q}\loc(\Omega,\qat)$ with $2<q<\infty$, then for $B(z_0,\epsilon) \subset \Omega' \subset \Omega$ and  $(z-z_0)^{-l}\xi \in \banach{q}(\Omega',\qat)$ due to Hölder's inequality 
\[
\|\xi\|_{\banach{q}(B(z_0,\epsilon),\qat)}\le\|(z-z_0)^l\|_{\banach{\infty}(B(z_0,\epsilon),\qat)}\|(z-z_0)^{-l}\xi\|_{\banach{q}(\Omega',\qat)}=\textbf{O}(\epsilon^l).
\]
We now show that if this holds for all $n\in\mathbb{N}$, then $\xi$ vanishes identically.
This is the so called \De{Strong unique continuation property} which implies that $\ord_{z_0}(\xi)$ is always finite for non-trivial $\xi$.
The proof of the \De{Strong unique continuation property} uses a Carleman inequality (compare with~\cite{Ca} and~\cite[Proposition~1.3]{Wo}).
\begin{proposition}[Carleman Inequality]\label{carleman inequality}
There exists a constant $S_p$ depending only on $1<p<2$, such that the following inequality holds:
\begin{align}\label{eq:carleman}
\left\| |z|^{-l}\xi\right\|_{\frac{2p}{2-p}}&\leq S_p\left\| |z|^{-l}\barpartial \xi\right\|_{p}
\quad\forall l\in\mathbb{Z}
\text{ and }\forall\xi\in C^{\infty}_0(\mathbb{C}\setminus\{0\},\qat).
\end{align}
\end{proposition}
The literature~\cite{Je,Ma,Ki2} deals with the much more difficult higher-dimensional setting and does not treat our case.
In~\cite[Proposition~2.6]{Wo} the analogous but weaker statement about the gradient term of the Laplace operator is proven.
David Jerison pointed out to the authors that these arguments carry over to the Dirac operator $\dirac-\diracJ\V^-=\diracJ(\barpartial-\V^-)$, which will be discussed further in Chapter~\ref{chapter:resolvents}.
\begin{proof}
Dolbeault's Lemma 
generalizes to the case of quaternionic-valued functions $\xi = \xi_1+\qj \xi_2$.
Thus for any smooth $\xi$ with compact support we have the equality
\begin{equation}\label{eq:dolbeault}
\xi(z) 
= \Op{I}_{\mathbb{C}} \barpartial \xi 
=\int_{\mathbb{C}}\frac{1}{z-w}\barpartial \xi(w)\frac{\dmu(w)}{\pi} .
\end{equation}
Due to partial integration for $l\in\mathbb{N}_0$ the integrals $\int_\mathbb{C}z^l\barpartial \xi\, \dmu(z)$ vanish.
Moreover, if the support of $\xi$ does not contain $0$ and therefore is disjoint from a small neighborhood of $0$, then these integrals vanish for all $l\in\mathbb{Z}$.
We conclude that for all $l\in\mathbb{Z}$
\[
\xi(z)=\int_{\mathbb{C}}\left(\frac{z}{w}\right)^l\frac{1}{z-w}\barpartial \xi(w)\frac{\dmu(w)}{\pi}.
\]
This is because the right hand side minus the right hand side of~\eqref{eq:dolbeault} is a finite sum of integrals of the form $z^{-m\pm 1}\int_\mathbb{C}w^l\barpartial \xi(w)\, \dmu(w)$ and therefore vanishes.
Finally, the first factor of the integral kernel, 
\[\left(\frac{|w|z}{|z|w}\right)^l\frac{1}{z-w}\frac{\dmu(w)}{\pi},\]
has absolute value $1$.
The corresponding integral operator maps $|z|^{-l}\barpartial \xi$ onto $|z|^{-l}\xi$ and, due to the Hardy-Littlewood-Sobolev theorem~\cite[Chapter~V. \S1.2 Theorem~1]{St}, it is bounded in $\mathcal{L}(\banach{p}(\mathbb{C},\qat),\banach{\frac{2p}{2-p}}(\mathbb{C},\qat))$ by some constant $S_p$ not depending on $m$.
\end{proof}
Due to a standard argument (e.g.~\cite[Proof of Theorem~5.1.4]{So} and~\cite[Section Carleman Method]{Wo}) this Carleman inequality implies the
\begin{theorem}[Strong Unique Continuation Property]\label{strong unique continuation}
Let $\xi\in \sobolev{1,p}\loc(\Omega,\qat)$ be $\V$\=/holomorphic with $1<p<2$ and $\V\in\banach{2}\loc(\Omega,\qat)$ on an open connected set $0\in\Omega\subset\mathbb{C}$ that has a root of infinite order at $0$.
This means, 
\begin{align}
\label{eq:infinite order}
\left\|\xi \right\|_{\banach{\frac{2p}{2-p}}(B(0,\varepsilon),\qat)}
&\leq\textup{\textbf{O}}(\varepsilon^l)\quad\text{as }\varepsilon \to 0, \text{ for all }l\in \mathbb{N}.
\end{align}
Then $\xi$ vanishes identically on $\Omega$.
\end{theorem}
\begin{proof}
The question is local so we may consider an open subset $0 \in\Omega'\subset\Omega$ with smooth boundary and compact closure in $\Omega$, so that $\V\in \banach{2}(\Omega',\qat)$ holds.
Choose $\varepsilon>0$ such that $\left\|\V\right\|_{\banach{2}(B(z,2\varepsilon),\qat)}\leq 1/(2S_p)$ holds for all $z\in\Omega'$ with the constant $S_p$ of~\eqref{eq:neumann-convergence}.
Let $\rho\in C^{\infty}_0(\mathbb{C},\mathbb{R})$ be $1$ on $B(0,\varepsilon)$ and $0$ on $\mathbb{C}\setminus B(0,2\varepsilon)$.

We now extend Carleman's inequality~\eqref{eq:carleman} to the function $\rho\xi$ that has a root of infinite order at zero.
Let $(\rho_n)_{n\in\mathbb{N}}$ be the sequence in $C^{\infty}(\mathbb{C},\mathbb{R})$ with $\rho_n=(z)=1-\rho(nz)$.
The support of $\rho_{n+1}-\rho_n$ is contained in $B(0,\tfrac{2\varepsilon}{n}) \setminus B(0,\tfrac{\varepsilon}{n+1})$.
On this set, by the bound~\eqref{eq:infinite order}, the $\banach{\frac{2p}{2-p}}$\=/norm of $\xi$ and the $\banach{p}$\=/norm of $\barpartial \xi = \V\xi$ are of order $\textbf{{O}}(n^{-l})$ for all $l\in\mathbb{N}$, and $\left\| |z|^{-l}\right\|_\infty$ is of the order   $\textbf{{O}}(n^l)$.
The norms $\|\rho\|_\infty$, $\|\rho_n\|_\infty$ and $\|\barpartial \rho\|_\infty$ are of order $\textbf{{O}}(n^0)$, and $\|\barpartial \rho_n\|_\infty = \textbf{{O}}(n)$.
This implies that for all $l\in\mathbb{Z}$ the norms $\||z|^{-l}\barpartial(\rho_{n+1}\rho\xi-\rho_{n}\rho\xi)\|_p$ decay in the limit $n\to\infty$ faster than every power of $\frac{1}{n}$.
Hence the sequence $|z|^{-l}\barpartial(\rho_n\rho\xi)$ is a Cauchy sequence in $\banach{p}(B(0,2\varepsilon),\qat)$.
Because it converges to $|z|^{-l}\barpartial(\rho\xi)$ pointwise almost everywhere, it follows from dominated convergence that it converges to $|z|^{-l}\barpartial(\rho\xi)$ in this Banach space.
Therefore Carleman's inequality holds for $\rho\xi$.

It now follows that
\[\left\| |z|^{-l}\rho\xi\right\|_{\frac{2p}{2-p}}\leq S_p\left\| |z|^{-l}\barpartial \rho\xi\right\|_p\leq S_p\left\| |z|^{-l}\rho\barpartial \xi\right\|_p+S_p\left\| |z|^{-l}E\right\|_p.\]
Here $E$ (for error) $=(\barpartial \rho)\xi$ is an $\banach{p}$\=/function supported in $\{z\mid \varepsilon\leq|z|\leq 2\varepsilon\}$.
Using the equality $\barpartial \xi=\V\xi$ and Hölder's inequality 
\begin{align*}
\left\| |z|^{-l}\rho\xi\right\|_{\frac{2p}{2-p}}
&\leq S_p\left\| |z|^{-l}\rho \V\xi\right\|_p+S_p\left\||z|^{-l}E\right\|_p\\
&\leq S_p\left\|\V\right\|_{\banach{2}(B(0,2\varepsilon))}\left\| |z|^{-l}\rho\xi\right\|_{\frac{2p}{2-p}}+S_p\left\| |z|^{-l}E\right\|_p.
\end{align*}
By the choice of $\varepsilon$, we have $S_p\left\|\V\right\|_{\banach{2}(B(0,2\varepsilon))} \leq \tfrac12$.
By subtracting the first term on the right-hand side, we obtain 
\[\left\| |z|^{-l}\rho\xi\right\|_{\frac{2p}{2-p}}\leq 2S_p\left\| |z|^{-l}E\right\|_p.\]
Now comes the crucial observation: $E$ is supported in
$\{z\mid \varepsilon\leq|z|\leq 2\varepsilon\}$, so
\begin{align*}
\left\| |z|^{-l}\rho\xi\right\|_{\frac{2p}{2-p}}\leq 2S_p\varepsilon^{-l}\left\|E\right\|_p
\quad\text{and}
\quad\left\|\bigl(\tfrac{\varepsilon}{|z|}\bigr)^l\rho\xi\right\|_{\frac{2p}{2-p}}\leq 2S_p\left\|E\right\|_p.
\end{align*}
Using the limit $l\to\infty$ we conclude that $\rho\xi$ vanishes on $B(0,\varepsilon)$, and therefore also $\xi$.
In other words, the set $\{z\mid\xi\text{ has a root of infinite order at }z\}$ is open, and in fact contains a ball of fixed radius $\varepsilon$ centered at any of its points.
So this set must be all of $\Omega'$ and the proof is complete.
\end{proof}
The above constructions of quaternionic functions on $\Omega\subset\mathbb{C}$ with prescribed poles can be extended to any Riemann surface $\X$.
A divisor is a locally finite $\mathbb{Z}$\=/linear combination of points of $\X$,~\cite[Section 16.1]{Fo}.
Meromorphic sections of a holomorphic $\mathbb{C}$\=/line bundle may be considered as holomorphic sections of another bundle: We recall the standard notation, for a holomorphic $\mathbb{C}$\=/line bundle $E$ on a Riemann surface $\X$ and a divisor $D$ there is a line bundle $E(D)$ whose sections are the meromorphic sections $\xi$ of $E$ with divisors $(\xi)$ that obey $(\xi)\ge-D$.
\begin{definition}
\label{meromorphic}
Let $E$ be a holomorphic $\mathbb{C}$\=/line bundle on a Riemann surface $\X$, $\V\in\pot{E}$, and let $D$ be a divisor $D$ on $\X$. 
Then $\Q{E(D),\V(D)}$ denotes the subsheaf of the meromorphic sections $\xi$ of $\Q{E,\V}$ on $\X$, such that $\ord_x(\xi)$ defined in Definition~\ref{order of roots} obeys $\ord_x(\xi)\ge-D(x)$ for all $x\in\X$.
%
\end{definition}
Note that for different holomorphic $\mathbb{C}$\=/line bundles $E\ne E'$ the spaces of the corresponding potentials $\pot{E}$ and $\pot{E'}$ are also different. In case of a pair of line bundles $E$ and $E'=E(D)$ a potential $\V\in\pot{E}$ induces an element of $\pot{E(D)}$, which we shall distinguish form $\V$. For this reason we denote the induces potential by $\V(D)$ and not just by $\V$. However, the sections of $\Q{E(D),\V(D)}$ are according to Definition~\ref{order of roots} called $\V$\=/meromorphic. To avoid confusion we shall use the notion of $\V$\=/meromorphic sections $\xi$ only if we cannot specify a divisor $D$, such that $\xi$ is a section of $\Q{E(D),\V(D)}$.

In order to clarify this notation we explain now its usage on open subset $\Omega\subset\mathbb{C}$. 
Like on any Riemann surface let $\unity$ denote on such open subsets $\Omega\subset\mathbb{C}$ the trivial $\mathbb{C}$\=/line bundle $\Omega\times\mathbb{C}$. 
For $z_0\in\Omega$ we set $\Omega_0=\Omega\setminus\{z_0\}$ and choose an open neighborhood $\Omega_1$ of $z_0$ with compact closure in $\Omega$. 
Then the holomorphic line bundle $\unity(lz_0)$ can be described by the cocycle $f_{10}=(z-z_0)^l$ with respect to the cover $\Omega_0\cup\Omega_1$ of $\Omega$. 
Because of uniqueness of analytic continuation, it is sufficient to give the sections of $\unity(lz_0)$ on $\Omega_0$. 
By Riemann's theorem on removable singularities, the sections of $\unity(lz_0)$ are exactly holomorphic functions $\xi$ on $\Omega_0$ such that $(z-z_0)^l\xi$ is bounded on $\Omega_0\cap\Omega_1$. 
Hence for $l<m$ any section of $\unity(lz_0)$ is also a sections of $\unity(mz_0)$. 
On the other hand, for any $l,m\in\mathbb{Z}$, if $\xi$ is a section of $\unity(lz_0)$, then $(z-z_0)^{l-m}\xi$ is a section of $\unity(mz_0)$.

For $\V\in\banach{2}\loc(\Omega,\qat)$ the sheaf $\Q{\unity(lz_0),\V(lz_0)}$ coincides on $\Omega_0$ with $\Q{\unity,\V}$. The stalk at $z_0$ contains all products of $(z-z_0)^{-l}$ with germs $\xi$ in the stalk of the sheaf $\Q{\unity,V_{z_0,-l}}$ at $z_0$ (with $\V_{z_0,-l}$ as in~\eqref{root potential}). As in the case of $\mathbb{C}$\=/line bundles, there is an inclusion of $\Q{\unity(lz_0),\V(lz_0)}$ into $\Q{\unity(mz_0),\V(mz_0)}$ for $l<m$.
By contrast, multiplication with $(z-z_0)^{-l}$ changes the potential to $\V_{z_0,l}$~\eqref{root potential} and instead gives an isomorphism $H^0(\Omega,\Q{\unity(-lz_0),\V(-lz_0)})\cong H^0(\Omega,\Q{\unity,\V_{z_0,l}})$ for any $l\in\mathbb{Z}$.
In particular, $H^0(\Omega,\Q{\unity(-lz_0),0})\cong H^0(\Omega,\Q{\unity,0})$. The above description of the bundle $\unity(lz_0)$ in terms of a cover and a corresponding cocycle induce such a description of the holomorphic $\qat$\=/line bundle $E(D)_\qat$ and the potential $\V\in\pot{E(D)}$. 
Let us now establish some further properties of holomorphic sections of $\qat$\=/line bundles, which are well known in the complex case:
\begin{lemma}
\label{quotient dimension}
\index{Divisor}
Let $E$ be a holomorphic $\mathbb{C}$\=/line bundle on a Riemann surface $\X$ and $\V\in\pot{E}$.
Any non-trivial section $\xi$ of $\Q{E,\V}$ defines a unique effective divisor $D=(\xi)$ on $\X$, such that $\xi$ is a holomorphic section of $\Q{E(-D),\V}$ without roots.
For any pair of divisors $D'\geq D''$ on $\X$ we have
\[\dim_\qat H^1(\X,\Q{E(D''),\V}/\Q{E(D'),\V}) =0,\]
and if $\X$ is compact,
\[\hspace{20mm}\dim_\qat H^0(\X,\Q{E(D''),\V}/\Q{E(D'),\V}) =\deg(D''-D').\hspace{18mm}\qed\]
\end{lemma}
\begin{proof}
The proof that on an open subset $\Omega\subset\mathbb{C}$ the roots of an $\V$\=/holomorphic function $\xi$ with potential $\V\in\banach{2}\loc(\Omega,\qat)$ are isolated, needs some preparation.

By the classical Cauchy's Integral Formula the Taylor polynomial of order $l-1$ depends continuously on the elements of the Bergman space $\Spa{A}^q(\Omega',\qat):=H^0(\Omega',\Q{\unity,0})\cap\banach{q}(\Omega',\qat)$ for $l\geq 1$, $2<q<\infty$, and $\Omega'\subset\Omega$.
Hence the projection onto the remainder term defines a bounded projection operator
\begin{equation}\label{eq:taylor-projector}
P:\Spa{A}^q(\Omega',\qat)\to H^0(\Omega',\Q{\unity(-lz_0),0})\cap\banach{q}(\Omega',\qat) .
\end{equation}
Let $Y$ be the image of the subspace $H^0(\Omega',\Q{\unity(-lz_0),\V(-lz_0)})\cap\banach{q}(\Omega',\qat)$ with respect to the map
\begin{align}\label{eq:op1}
\unity-\Op{I}_{\Omega'} \V|_{\Omega'}&:&H^0(\Omega',\Q{\unity,\V})\cap\banach{q}(\Omega',\qat)&\to\Spa{A}^q(\Omega',\qat) .
\end{align}
We claim that for sufficiently small $\Omega' \subset \Omega$, the restriction of $P$ to $Y$ is an isomorphism.
If $\Omega'$ is chosen such that both $\unity-\Op{I}_{\Omega'}\V|_{\Omega'}$ and $\unity-\Op{I}_{\Omega'}\V_{z_0,l}|_{\Omega'}$ are invertible then the map
\begin{equation} \label{eq:iso}
(\unity-\Op{I}_{\Omega'}\V|_{\Omega'}) \circ (z-z_0)^{l} \circ (\unity-\Op{I}_{\Omega'}\V_{z_0,l}|_{\Omega'})^{-1}
\end{equation}
from $\Spa{A}^q(\Omega',\qat)$ to itself is a small perturbation of $(z-z_0)^l$.
Since $P \circ (z-z_0)^l$ is an isomorphism, so too is $P$ composed with~\eqref{eq:iso}.
By the definition of $Y$, Equation~\eqref{eq:iso} is an isomorphism of $\Spa{A}^q(\Omega',\qat)$ onto $Y$, and therefore $P$ restricted to $Y$ must itself be an isomorphism as claimed.

The following diagram summarizes the relationship between the aforementioned spaces and maps.
Note that we have omitted the domain $\Omega'$ for brevity.
\begin{equation*}
\begin{tikzcd}[row sep = 10mm, column sep = 25mm]
\Q{\unity(-lz_0),0}\cap\banach{q}&\arrow[l,two heads,hook',"(z-z_0)^l"']\Spa{A}^q\arrow[r,two heads,hook,"(\unity-\Op{I}_{\Omega'}\V_{z_0,l}|_{\Omega'})^{-1}"]&\Q{\unity,\V_{z_0,l}}\cap\banach{q}\arrow[d,two heads,hook,"(z-z_0)^l"]\\\arrow[u,"P"]\Spa{A}^q&\arrow[lu,"P|_{Y}"']\arrow[l,hook',"i"']Y&\arrow[l,two heads,hook',"\unity-\Op{I}_{\Omega'}\V|_{\Omega'}"]\Q{\unity(-lz_0),\V(-lz_0)}\cap\banach{q}
\end{tikzcd}
\end{equation*}

For $\V=0=\V_{z_0,l}$ the isomorphisms on the right hand side of both columns and the operator $P|_Y$ simplify to the unity maps of the corresponding spaces.
In this case the diagram is indeed commutative.
For general $\V$ the small triangle on the left hand side is commutative, in contrast to the pentagon on the right hand side, which commutes only approximately: by this we mean that the three operators in~\eqref{eq:iso} composed with $P|_Y$ differ from the isomorphism $(z-z_0)^l$ on the left hand side in the upper column by an operator whose norm is arbitrarily small for sufficiently small $\|\V\|_{\banach{2}(\Omega',\qat)}$.

Because $P$ is a projection, it firstly follows from the claim that the intersection of $Y$ with $\ker(P)$ is trivial.
Secondly, it also follows that $Y+\ker(P) =\Spa{A}^q(\Omega',\qat)$ holds, because $\xi-(P|Y)^{-1}\,P\,\xi$ belongs to $\ker(P)$ for any $\xi\in\Spa{A}^q(\Omega',\qat)$.
Therefore $\Spa{A}^q(\Omega',\qat)$ is isomorphic to $Y \oplus \ker(P)$.
Hence $\mathrm{codim}(Y) = \dim(\ker(P))=l$ because $\ker(P)$ is the space of quaternionic polynomials of degree at most $l-1$.
Equivalently the stalk of the quotient sheaf $\Q{\unity,\V}/\Q{\unity(-lz_0),\V(-lz_0)}$ at $z_0$ has quaternionic dimension $l$.

We now claim that the roots of a non-trivial $\xi\in H^0(\Omega,\Q{\unity,\V})$ are isolated.
Because~\eqref{eq:iso} is an isomorphism, $H^0(\Omega',\Q{\unity(-lz_0),\V(-lz_0)})\cap\banach{q}(\Omega',\qat)$ is a closed subspace of $H^0(\Omega',\Q{\unity,\V})\cap\banach{q}(\Omega',\qat)$.
The closed balls in $\banach{q}(\Omega',\qat)$ are weakly compact (because $\banach{q}(\Omega',\qat)$ is the dual of a Banach space), implying that the norm $\|\,\cdot\,\|_q$ is weakly semi-continuous.
Hence there exists a $\Tilde{\xi}\in H^0(\Omega',\Q{\unity(-lz_0),\V(-lz_0)})\cap\banach{q}(\Omega',\qat)$ with minimal distance $\|\xi-\Tilde{\xi}\|_q$.
Because the operator~\eqref{eq:iso} depends continuously on $z_0\in\Omega'$, the minimal distance $\|\xi-\Tilde{\xi}\|_q$ also depends continuously on $z_0$.
This shows that the complement of the set of roots is open.
For every non-trivial $\V$\=/holomorphic $\xi$ and every $z_0\in\Omega$, the function $(z-z_0)^{-\ord_{z_0}(\xi)}\xi$ is $\V_{z_0,\ord_{z_0}(\xi)}$~\eqref{root potential}\=/holomorphic and has no root at $z_0$.
Hence $\xi$ has no roots in a punctured neighborhood of $z_0$, and thus the roots of $\xi$ are isolated. In particular, any section $\xi\in H^0(\X,\Q{E,\V})$ of a complex $\qat$\=/line bundle with holomorphic structure defined by $\V\in\pot{E}$ has a well defined divisor $D$ with $D(x)=\ord_x(\xi)$ for all $x\in\X$. In particular, $\xi$ is a holomorphic section of $\Q{E(D),\V(D)}$ without roots. The statements about the dimensions of the cohomology groups now follows by standard methods (compare~\cite[Proof of 16.7 Lemma]{Fo}) from the local statement $\dim\big(\Q{\unity,\V}/\Q{\unity(-lz_0),\V(-lz_0)}\big)=l$.
\end{proof}
We shall enhance such similarities between complex and quaternionic analysis in Chapter~\ref{chapter:weierstrass}: There is a well-defined notion of residue (Lemma~\ref{lem:10 pairing properties}) and the order of poles increases by one under differentiation (Lemma~\ref{lem:change of order}).

We finish the chapter with a result which is used in the subsequent chapter and improved in several steps in Part~II.
A natural question is the behavior of the holomorphic structure under changes of the potential.
Due to Cauchy's integral formula~\eqref{eq:cif} we may investigate its inverse instead, the bounded operator $\Op{I}_{\Omega',\V}$.
We have seen that this can be controlled through its Neumann series (compare~\eqref{eq:neumann} and~\eqref{eq:resolvent}):
\[
\Op{I}_{\Omega',\V}
=\sum_{l=0}^\infty\Op{I}_{\Omega'}(\V\Op{I}_{\Omega'})^l.
\]
Specifically the terms 
\begin{equation}\label{eq:op2}
\V\mapsto\Op{I}_{\Omega'}(\V\Op{I}_{\Omega'})^l.
\end{equation}
of the series are bounded operators in $\mathcal{L}(\banach{p}(\Omega',\qat), \banach{q}(\Omega',\qat))$, with $q = \frac{2p}{2-p}$, and the series converges uniformly for $\|\V\|_{\banach{2}(\Omega',\qat)}<S_p$.
The idea of the lemma below is that, if we decrease $q$ and $S_p$ appropriately, then we gain weak continuity in $\V$.

It is clear what it means to decrease $q$, but things a little trickier for $S_p$.
Recall its definition
\[
S_{p}^{-1} = \|\Op{I}_{\mathbb{C}}\|_{\mathcal{L}(\banach{p}(\mathbb{C},\qat),\banach{\frac{2p}{2-p}}(\mathbb{C},\qat))}
\]
from~\eqref{eq:neumann-convergence} for $1<p<2$.
In particular, its dependence on $p$ might be complicated.
Therefore define
\begin{align}\label{def:Sp-}
S_{p,\varepsilon} =\inf\limits_{p'\in[p-\varepsilon,p]}S_{p'}
\quad\text{and}\quad
S_p^-&=\lim_{\epsilon\downarrow0}S_{p,\varepsilon}.
\end{align}
Since the infima $S_{p,\varepsilon}$ are bounded from above ($S_{p,\varepsilon} \leq S_p$) and are increasing for diminishing $\varepsilon$, the limit~\eqref{def:Sp-} must exist. 
$S_p^-$ is positive, because the Marcinkiewicz Interpolation Theorem~\cite[Chapter~4, Theorem~4.13]{BS} yields an upper bound on $S_p^{-1}$ that depends continuously on $p\in(1,2)$.
``Appropriately'' decreasing $S_p$ means considering potentials $\V$ with norm smaller that $S_p^-$.

\begin{lemma}\label{weakly continuous}
Let $1<p<2$, $1<q<\tfrac{2p}{2-p}$, and $S'_p<S_p^-$.
Then the map
\begin{align}\label{eq:integral kernel}
\overline{B(0,S'_p)}\subset\banach{2}(\Omega',\qat)&\to\mathcal{L}(\banach{p}(\Omega',\qat),\banach{q}(\Omega',\qat)),&\V&\mapsto\Op{I}_{\Omega',\V}
\end{align}
is weakly continuous.
By weak continuity we mean the continuity with respect to the weak topology of $\banach{2}(\Omega',\qat)$ and the norm topology of $\mathcal{L}(\banach{p}(\Omega',\qat),\banach{q}(\Omega',\qat))$.
\end{lemma}

We remark that the lemma does not imply that the map $\V \to \Op{I}_{\Omega',\V}$ is weakly continuous on $\V \in \bigcup_{S_{p}' < S_p^-} \overline{B(0,S_{p}')} = B(0, S_p^-)$, because there exists sequences in $B(0, S_p^-)$ that weakly converge to $\V \in B(0, S_p^-)$ but whose norms converge to $S_p^-$.
For example, choose $V_n := S_p^- (1-n^{-1}) \mu(B(0,n^{-1}))^{-0.5} \unity_{B(0,n^{-1})}$, a sequence of potentials supported on ever smaller balls.
This sequence converges weakly to $0$, but $\|\V_n\|_{\banach{2}} = S_p^- (1-n^{-1})$ converges from below to $S_p^-$.

\begin{proof}
As motivated above, by decreasing $q$ we gain some regularity.
More precisely the operators $\Op{I}_{\Omega'}$ become compact and can be approximated by finite rank operators.
The annoying detail is that we must decrease $q$ for each use of $\Op{I}_{\Omega'}$.
Choose $\varepsilon \in (0,p-1)$ sufficiently small, such that $S_p' < S_{p,\varepsilon} = \inf_{p'\in[p-\varepsilon,p]}S_{p'}$.
We are able to accommodate a decrease in $p$ as low as either $p- \varepsilon$ or the Hölder conjugate of $q$.
Therefore choose a decreasing sequence $p_k$ with 
\begin{align*}
p
=p_0 > \cdots > p_k > \cdots > \max\{p-\varepsilon, \tfrac{2q}{2+q}\},
\quad\text{and}\quad
q_{k} := \tfrac{2p_{k+1}}{2-p_{k+1}}.
\end{align*}
The operators $\Op{I}_{\Omega'}$ from $\banach{p_k}(\Omega',\qat)$ to $\sobolev{1,p_k}(\Omega',\qat)$ are continuous, and the further inclusion into $\banach{q_k}(\Omega',\qat)$ is compact, due to the Rellich-Kondrachov theorem~\cite[Theorem~6.3]{Ad} and $q_k < \frac{2p_k}{2-p_k}$.
This makes the composition~\eqref{eq:op2} an operator from $\banach{p_0}(\Omega',\qat)$ to $\banach{q_l}(\Omega',\qat)$.
Since $q_l > q$, with one further embedding it is an operator from $\banach{p}(\Omega',\qat)$ to $\banach{q}(\Omega',\qat)$ for the given $p$ and $q$ in the statement of the lemma.

Due to~\cite[Theorem~II.5.11]{LT}, the Lebesgue spaces $\banach{q}(\Omega',\qat)$ for any $1<q<\infty$ have a Schauder basis.
Consequently they have the approximation property, and these compact operators $\Op{I}_{\Omega'}$ are all norm-limits of finite rank operators (compare~\cite[Section~I.1.a]{LT}).
If we replace all operators $\Op{I}_{\Omega'}$ in~\eqref{eq:op2} by finite rank operators, then the result is weakly continuous.
This shows that~\eqref{eq:op2} is a uniform limit of weakly continuous functions and therefore itself weakly continuous.

Finally, we want to conclude that $\Op{I}_{\Omega',\V}$ is weakly continuous because it is the uniform limit of weakly continuous functions. 
However by changing the spaces, we have changed the norms of the terms, so we need to re-check the uniform convergence of the Neumann series.
The norm of the compact embedding $\sobolev{1,p_k}(\Omega',\qat)\hookrightarrow\banach{q_k}(\Omega',\qat)$ is $\mu(\Omega')^{\frac{1}{q_k}-\frac{2-p_k}{2p_k}}$.
Therefore the norm of the composition~\eqref{eq:op2} from $\banach{p}(\Omega',\qat)$ to $\banach{q}(\Omega',\qat)$ is bounded by the product of the norms
\[
\mu(\Omega')^{\frac{1}{q}-\frac{2-p_l}{2p_l}}S_{p_l}^{-1}
\prod_{k=0}^{l-1} \|\V\|_{\banach{2}(\Omega')} \mu(\Omega')^{\frac{1}{q_k}-\frac{2-p_k}{2p_k}} S_{p_k}^{-1}
\leq \mu(\Omega')^{\frac{1}{q}-\frac{2-p}{2p}} S_{p,\varepsilon}^{-1}
\left( \frac{S_p'}{S_{p,\varepsilon}} \right)^{l-1}.
\]
As an aside, when calculating the norm, we can in effect bring all the embeddings into a single embedding and treat the operators $\Op{I}_{\Omega'}$ as if they were from $\banach{p_k}(\Omega',\qat)$ to $\banach{\frac{2p_k}{2-p_k}}(\Omega',\qat)$.
This same idea will reoccur in the proof of Lemma~\ref{weakly continuous resolvent}.
This is a geometric series with ratio $S_p' S_{p,\varepsilon}^{-1} < 1$, so the Neumann series converges uniformly for $\|\V\|_{\banach{2}(\Omega',\qat)} < S_p'$.
\end{proof}

\chapter{Local Darboux Transformations}
\label{chapter:darboux}
Given a partial differential equation that can be written in Lax form, a Darboux transformation produces a new solution from a known solution and a solution of the related linear Lax equation (compare~\cite{GHZ}).
In our situation, the linear equation $(\barpartial -\V)\phi=0$ is the linear Lax equation  for the Davey-Stewartson equation, which is a non-linear evolution equation for the complex-valued potential $\V$.
Given a particular solution $\upsilon$ without roots of the linear Lax equation $(\barpartial -\V)\upsilon=0$, we will transform any other solution $\phi$ into a solution $\psi$ of the linear Lax equation $(\barpartial-\U)\psi=0$ of the new solution $\U$ of the Davey-Stewartson equation.
Rather than the transformation $\V\mapsto\U$ between the solutions of the Davey-Stewartson equation, our main interest concerns the transformation $(\upsilon,\phi)\mapsto(\chi,\psi)$, where $\chi$ is essentially the inverse of $\upsilon$ and solves $(\barpartial-\qk\Bar{\U}\qk)\chi=0$.
Here both sections of the former pair $(\phi,\upsilon)$ belong to the kernel of the same holomorphic structure $\barpartial-\V$ while the sections of the transformed pair $(\psi,\chi)$ belong to the kernels of two different holomorphic structures $\barpartial-\U$ and $\barpartial-\qk\Bar{\U}\qk$.
They will turn out to be in some sense paired, and their potentials will be denoted by $\U$ and $\U\sd$, respectively.
In this chapter we establish this transformation only locally: on open and bounded subsets $\Omega\subset\mathbb{C}$.
We will address the transformation on Riemann surfaces in Corollary~\ref{cor:global darboux}.

Let us comment on the notation of the local functions $\upsilon$, $\phi$, $\chi$ and $\psi$.
In the Chapters~\ref{chapter:kodaira}-\ref{chapter:weierstrass} we shall establish two different representations of admissible maps in terms of two holomorphic sections of $\qat$\=/line bundles.
The Greek letters $\upsilon$ and $\phi$ are reserved for the notation of the pair of holomorphic sections of the first representation, which is called the Kodaira representation and developed in Chapter~\ref{chapter:kodaira}.
The two other letters $\psi$ and $\chi$ are reserved for the holomorphic sections of the other representation, which is called the Weierstraß representation and developed in Chapter~\ref{chapter:weierstrass}.
Later we will see that the global Darboux Transformation, Corollary~\ref{cor:global darboux}, transforms the Kodaira pair $(\upsilon,\phi)$ into the corresponding Weierstraß pair $(\chi,\psi)$.
Therefore in this chapter the four local functions that are related by the Darboux transformation are also denoted by these symbols.
The Greek letter $\xi$ is reserved for a general local or global section of a $\qat$\=/line bundle.
In this chapter, like the previous chapter, we mostly work in a single coordinate chart and therefore omit subscripts that indicate the chart.

We investigate the local Darboux transformation first for smooth potentials and sections in the following calculation, which serves as a template.
One of the main results of Part~I, Theorem~\ref{thm:darboux}, extends this local Darboux transformation and its inverse to $\banach{2}\loc$\=/potentials.
By applying this theorem iteratively, we obtain an analogue to the statement of classical complex analysis that holomorphic functions are infinitely differentiable.
We end this chapter with Corollary~\ref{weyls lemma}, an analogue of Weyl's Lemma, where we prove that a function $\psi$ is $\U$\=/holomorphic if it gives rise to a closed current.

As motivation, we consider the case of a smooth potential $\V:\Omega\rightarrow\qat$ on an open subset $\Omega\subset\mathbb{C}$.
Let $\phi$ and $\upsilon$ be two $\V$\=/holomorphic functions on $\Omega$.
These functions have the derivatives
\begin{align*}
\index{Potential}
d\upsilon&=(d\Bar{z}\V-dz(\B+\U^-))\upsilon &
d\phi&=(d\Bar{z}\V+dz\partial)\phi
\end{align*}
with $\qat^+$\=/valued potentials $\B$ and $\qat^-$\=/valued potential $\U^-$, defined as $\B+\U^-=-(\partial\upsilon)\upsilon^{-1}$.
If $\upsilon$ has no roots, then it induces the connection $\nabla = dz(\partial+\B+\U^-) + d\Bar{z}(\barpartial-\V)$ on the $\qat$\=/line bundle.
Equivalently, it is the unique connection for which $\nabla \upsilon = 0$.
This connection is flat since $\nabla^2(\upsilon \alpha) = \nabla(\upsilon \wedge d\alpha) = - \upsilon \wedge d^2\alpha = 0$, as for $\mathbb{R}$\=/connections.
The zero curvature equation implies certain relations among the potentials.
By defining $\U^+:=\V^+$ the application of Equation~\eqref{eq:connection applied to forms} gives
\begin{align*}
0 = \nabla^2\xi
&= \nabla [dz(\partial+\B+\U^-)\xi + d\Bar{z}(\barpartial-\V)\xi] \\
&= -dz\wedge \nabla^+(\partial+\B+\U^-)\xi -d\Bar{z}\wedge \nabla^-(\partial+\B+\U^-)\xi  \\
&\qquad - d\Bar{z}\wedge \nabla^+ (\barpartial-\V)\xi - dz\wedge \nabla^- (\barpartial-\V)\xi \\
&= -dz\wedge d\Bar{z} (\barpartial - \V^+)(\partial+\B+\U^-)\xi -d\Bar{z}\wedge dz \U^-(\partial+\B+\U^-)\xi  \\
&\qquad - d\Bar{z}\wedge dz (\partial + B) (\barpartial-\V)\xi - dz\wedge d\Bar{z} (-\V^-) (\barpartial-\V)\xi \\
&= -dz\wedge d\Bar{z} (\barpartial - \U)(\partial+\B+\U^-)\xi + dz \wedge d\Bar{z}(\partial + B + \V^-) (\barpartial-\V)\xi,
\end{align*}
or without the forms
\begin{align}
\label{eq:zero curvature 2}
\index{Zero curvature equation}
(\barpartial-\U)(\partial+\B+\U^-)&=(\partial+\B+\V^-)(\barpartial-\V).
\end{align}
Since $\upsilon$ has no roots the quotient $\upsilon^{-1}\phi=|\upsilon|^{-2}\Bar{\upsilon}\phi$ is well-defined and obeys 
\begin{align}\label{eq:derivative quotient pre}
d(\upsilon^{-1}\phi)&=  \upsilon^{-1}
\left(d\phi-d\upsilon\,\upsilon^{-1}\,\phi\right) 
= \upsilon^{-1}dz\left(\partial+\B+\U^-\right)\phi = \upsilon^{-1} dz \psi  
\end{align}
for 
\begin{align}
\label{eq:psi}
\psi=(\partial+\B+\U^-)\phi.  
\end{align}
As we will explore further in Chapter~\ref{chapter:weierstrass}, the right hand side of Equation~\eqref{eq:derivative quotient pre} can be interpreted as a quadratic expression of two holomorphic sections of two so-called ``paired'' $\qat$\=/line bundles.
More concretely, we saw for the pairing in Definition~\ref{def:00 pairing} that a $\qj$\=/dual basis was appropriate.
Hence we should insert a $\qj$ by defining
\begin{align}
\label{eq:chi}
\chi = (\overline{\qj\upsilon})^{-1},
\end{align}
so that Equation~\eqref{eq:derivative quotient pre} comes to resemble Equation~\eqref{eq:00 pairing in basis}:
\begin{equation}
\label{eq:derivative quotient final}
d(\upsilon^{-1}\phi)=\Bar{\chi}\qj dz\psi  .
\end{equation}
The functions $\psi$ and $\chi$ are moreover holomorphic.
In fact, due to zero curvature equation~\eqref{eq:zero curvature 2}, and since $\phi$ is $\V$\=/holomorphic, we have $(\barpartial-\U)(\partial+\B+\U^-)\phi=0$.
Therefore the first $\psi=(\partial+\B+\U^-)\phi$ is $\U$\=/holomorphic.
For $\chi$ we calculate directly
\begin{gather}\begin{aligned}
d(\overline{\qj \upsilon})^{-1}&=-\qj\overline{d\upsilon\cdot\upsilon^{-1}}\, \Bar{\upsilon}^{-1}=-\qj \, \overline{\left(d\Bar{z}\V-dz(\B+\U^-)\right)} \Bar{\upsilon}^{-1}\\&=\qj(d\Bar{z}(\Bar{\B}+\V^-)-dz(\Bar{\V}^++\U^-))\Bar{\upsilon}^{-1}\\&=(dz(\B-\qj\V^-\qj)-d\Bar{z}(\V^+-\qj\U^-\qj))\chi.\label{eq:dchi}
\end{aligned}\end{gather}
We supplemented the anti-commuting $(1,0)$\=/potential $\U^-$ of $\upsilon$ by the commuting potential $\U^+=\V^+$.
As a consequence the $(0,1)$\=/potential $\U\sd$ of $\chi$ is equal to $-\V^++\qj\U^-\qj=-\U^++\qj\qi\U^-\qi\qj=\qk\Bar{\U}\qk$:
\begin{equation}
\label{eq:inverse derivative}
\barpartial \chi=\barpartial \,(\overline{\qj \upsilon})^{-1}=\qk\Bar{\U}\qk\chi.
\end{equation}
This calculation establishes a one-to-one correspondence between a pair $(\upsilon,\phi)$ of two $\V$\=/holomorphic functions and another pair $(\chi,\psi)$ of a $\qk\Bar{\U}\qk$\=/holomorphic and a $\U$\=/holomorphic function, respectively.
The transformation $(\phi,\upsilon)\mapsto(\psi,\chi)$ is the local Darboux transformation.

We now do some preparatory calculations that will later help us to recover the potentials from the functions.
Equation~\eqref{eq:zero curvature 2} implies the following equation
\begin{align}
\label{eq:zero curvature 1}
d(\V^+d\Bar{z}-\B dz)&=(|\U^-|^2-|\V^-|^2)d\Bar{z}\wedge dz
=2\ci (|\U^-|^2-|\V^-|^2)\dmu.
\end{align}
We also calculate the derivative of $2\ln|\upsilon|$, which must be a $\mathbb{R}$\=/valued $1$\=/form:
\begin{gather}
\label{derivative log}
\begin{aligned}
2d\ln|\upsilon|
&=\frac{\overline{d\upsilon} \upsilon+\Bar{\upsilon}d\upsilon}{\Bar{\upsilon}\upsilon}
=\frac{\Bar{\upsilon}\left(\Bar{\V}dz-(\Bar{\B}+\Bar{\U}^-)d\Bar{z}+d\Bar{z}\V-dz(\B+\U^-)\right)\upsilon}{\Bar{\upsilon}\upsilon}\\
&=dz(\Bar{\V}^+-\B)+d\Bar{z}(\V^+-\Bar{\B}).
\end{aligned}\end{gather}
By setting $\V^+=0$ in the above two equations we obtain
\begin{align}
\label{poisson}
-2\partial\ln|\upsilon|&=\B,
&\barpartial\B&=|\V^-|^2-|\U^-|^2.
\end{align}

In the following discussion concerning this transformation we make use of the \emph{Lorentz spaces}\index{Lorentz spaces} $\banach{p,q}$ (\cite[Chapter~V\S3]{SW},\cite[Chapter~3 Definition~4.1]{BS} and~\cite[1.8.6~Definition]{Zi}).
These rearrangement invariant Banach spaces are an extension of the family of the usual Banach spaces $\banach{p}$ indexed by an additional parameter $1\leq q\leq\infty$.
Their definition follows from two observations.
The $\banach{p}$\=/norm of a function can be calculated in terms of its cumulative distribution $\mu_f(s) = \mu \{ z \in \Omega \mid |f(z)| > s \}$:
\[
\|f\|_p^p = p \int_0^\infty s^p\, \mu_f(s) \,\frac{ds}{s}.
\]
There is also the definition of \emph{weak $\banach{p}$\=/spaces}, which consists of functions bounded with respect to
\[
\sup_{s > 0} s\, \mu_f(s)^{1/p}.
\]
The essential idea is to view these as the $\banach{p}$- and $\banach{\infty}$\=/norms of $s \mu_f(s)^{1/p}$ respectively, and to generalize to the $\banach{q}$\=/norm.
The correct definition is more technical however and we refer the reader to the references.
For $1 <p < \infty$ the Lorentz spaces $\banach{p,p}$ coincide with the usual $\banach{p}$\=/spaces.
Note that for $p=1,\infty$ there are competing definitions of $\banach{p,q}$ due to the existence for two families of quasinorms that are equivalent only for $1 < p < \infty$ (cf.~\cite[Chapter~3 Lemma~4.5]{BS},\cite[1.8.10~Lemma]{Zi}).
We will only use $\banach{1,\infty} = \banach{1}$ and $\banach{\infty,\infty} = \banach{\infty}$.
Like the Lebesgue spaces, on a finite measure space the Lorentz spaces are nested: 
$\banach{p,q}$ is contained in $\banach{p,q'}$ if $q\leq q'$.
Together with $\banach{p,p} = \banach{p}$ this implies $\banach{p,q}$ is contained in $\banach{p',q'}$ if $p>p'$.
In this sense, the Lorentz spaces are a refinement of the Lebesgue hierarchy.

\begin{example}
\label{eg:weak L2}
On $\Omega = B(0,1) \subset \mathbb{C}$ the function $f(z) = z^{-1}$ is $\banach{p}$ for $1 \leq p < 2$ but it is not $\banach{2}$:
\[
\|f\|_2^2
= \int_\Omega |z^{-1}|^2\,\dmu
= 2\pi \int_0^1 r^{-2+1}\,dr
= \infty.
\]
It is however weak $\banach{2}$, since $\mu_f(s) = \mu \{ |z| < s^{-1}\} = \pi s^{-2}$ and
\[
\sup_{s>0} s (\pi s^{-2})^{1/2} = \sqrt{\pi}.
\]
A consequence of this is that $\xi = \Bar{z}^\alpha$ is not $\V$\=/holomorphic (other than $\alpha = 0$) due to the fact that $\V = (\barpartial \xi) \xi^{-1} = \alpha \Bar{z}^{-1}$ is not $\banach{2}$.
\end{example}

In~\cite{O} Hölder's inequality and Young's inequality are generalized to these Lorentz spaces (\cite[Chapter~4 Section~7]{BS} and~\cite[Chapter~2 Section~10]{Zi}):
\index{Hölder's inequality}
\begin{lemma}[Generalised Hölder's inequality]\label{generalized hoelder}
Either for $1/p_1+1/p_2=1/p_3<1$ and $1/q_1+1/q_2\geq1/q_3$ or for $1/p_1+1/p_2=1$, $1/q_1+1/q_2\geq 1$ and $(p_3,q_3)=(1,\infty)$ there exists some constant $C>0$ with
\[\|fg\|_{(p_3,q_3)}\leq C\|f\|_{(p_1,q_1)}\|g\|_{(p_2,q_2)}.\]
\end{lemma}
\begin{lemma}[Generalised Young's inequality]\label{generalized young}
Either for $1/p_1+1/p_2-1=1/p_3>0$ and $1/q_1+1/q_2\geq 1/q_3$ or for $1/p_1+1/p_2=1$, $1/q_1+1/q_2\geq 1$ and $(p_3,q_3)=(\infty,\infty)$ there exists some constant $C>0$ with
\[\|f\ast g\|_{(p_3,q_3)}\leq C\|f\|_{(p_1,q_1)}\|g\|_{(p_2,q_2)}.\]
\end{lemma}
Lemma~\ref{generalized young} implies that $\Op{I}_\Omega$ belongs on any open and bounded subset $\Omega\subset\mathbb{C}$ and for any $1<p<2$ and $q=\frac{2p}{2-p}$ to the following four Banach spaces:
\begin{itemize}
\item $\mathcal{L}(\banach{1}(\Omega,\mathbb{H}),\banach{2,\infty}(\Omega,\mathbb{H}))$
\item $\mathcal{L}(\banach{2,1}(\Omega,\mathbb{H}),\banach{\infty}(\Omega,\mathbb{H}))$
\item $\mathcal{L}(\banach{p}(\Omega,\mathbb{H}),\banach{q,p}(\Omega,\mathbb{H}))$
\item $\mathcal{L}(\banach{p,q}(\Omega,\mathbb{H}),\banach{q}(\Omega,\mathbb{H}))$
\end{itemize}
In the literature the first two statements are condensed to the terminology that the convolution with the potential $z\mapsto|z|^{-1}$ is of restricted weak type $(1,2)$ and $(2,\infty)$ (compare~\cite[Chapter~V\S2]{SW} and~\cite[Remark~5.]{CP}).
In the present situation of open and bounded subsets $\Omega\subset\mathbb{C}$ a direct proof is contained in~\cite[Proof of Theorem~8.3]{KZPS}.
Due to the extension of the Marcinkiewicz Interpolation theorem in Stein and Weiss~\cite[Chapter~V Theorem~3.15]{SW}, the first two statements together imply the latter two.

We now restrict to the case $\V^+=0$ and omit the superscripts of $\V^-$ and $\U^-$.
In the following main theorem of this chapter, we extend the local Darboux transformation to $\banach{2}\loc$\=/potentials.
But before we state this theorem, we remark that the extension to the slightly smaller space of $\banach{2,1}$\=/potentials is relatively straightforward.
Indeed, for $\V\in\banach{2,1}\loc(\Omega,\qat^-)$ the Generalized Hölder's inequality~\ref{generalized hoelder} and Generalized Young's inequality~\ref{generalized young} combined with the arguments similar to those in Chapter~\ref{chapter:local} imply that the $\V$\=/holomorphic functions $\upsilon$ and $\phi$ are continuous and belong to the Sobolev space $\sobolev{1,2}\loc(\Omega,\qat)$.
Hence for $\upsilon$ without roots the potentials belong to $\B\in\banach{2}\loc(\Omega,\qat^+)$ and $\U\in\banach{2}\loc(\Omega,\qat^-)$.
Furthermore, $\upsilon^{-1}$ is continuous and belongs to $\sobolev{1,2}\loc(\Omega,\qat)$.
In this case the statements (1)(a)--(f) in the following Theorem~\ref{thm:darboux} follow from~\eqref{eq:zero curvature 2},~\eqref{eq:inverse derivative} and~\eqref{eq:derivative quotient final}.
\begin{theorem}[Local quaternionic Darboux transformation]\label{thm:darboux}
\index{Darboux transformation}
\phantom{boo}
\begin{enumeratethm}[label={\upshape(\arabic*)}]
\item 
Let $\Omega\subset\mathbb{C}$ be open, $V\in\banach{2}\loc(\Omega,\qat^-)$ and let $\upsilon$ and $\phi$ be $\V$\=/holomorphic with $\upsilon$ having no roots on $\Omega$ (in the sense of Definition~\ref{order of roots}).
Then the following holds:
\begin{enumerate}[label={\upshape(\alph*)}]
\item $\B:=-2\partial\ln|\upsilon| \in \banach{2,\infty}\loc(\Omega,\qat^+)$.
\item $\U:=-(\partial\upsilon)\upsilon^{-1}-\B\in \banach{2}\loc(\Omega,\qat^-)$.
\item $\barpartial\B=|\V|^2-|\U|^2$ as distributions.
\item $\chi:=(\overline{\qj \upsilon})^{-1}$ is $\qk\Bar{\U}\qk$\=/holomorphic.
\item $\psi:=(\partial+\B+\U)\phi$ is $\U$\=/holomorphic.
\item $d(\upsilon^{-1}\phi)=\Bar{\chi}\qj dz\psi$ and $\upsilon^{-1}\phi\in \bigcap_{1<p<2}\sobolev{2,p}\loc(\Omega,\qat)\subset C(\Omega,\qat)$.
\end{enumerate}

\item
Conversely, let $\Omega\subset\mathbb{C}$ be simply connected and open, $\U\in\banach{2}\loc(\Omega,\qat^-)$, $\psi$ be $ \U$\=/holomorphic, and $\chi$ be $\qk\Bar{\U}\qk$\=/holomorphic without roots on $\Omega$.
Then
\begin{enumerate}[label={\upshape(\alph*)}]
\item
$dF:=\Bar{\chi}\qj dz\psi$ is a closed quaternionic-valued current on $\Omega$.
\item
$\V:=\Bar{\chi}^{-1}\,\overline{\partial\chi}-2\barpartial\ln|\chi|\in\banach{2}\loc(\Omega,\qat^-)$.
\item
$\upsilon:=-(\overline{\qj\chi})^{-1}$ and $\phi:=\upsilon F$ are $\V$\=/holomorphic.
\end{enumerate}
\end{enumeratethm}
\end{theorem}
\begin{proof}
The proof of this theorem is divided into seven steps, which include four lemmas and a corollary.

\proofstep{Step A.} 
In Steps~A--D we choose for $\V\in\banach{2}\loc(\Omega,\qat^-)$ small bounded open subsets $\Omega'\subset\Omega''\subset\Omega$ with smooth boundaries in $\Omega''$ and $\Omega$, respectively.
We define $\V'$ as the product of $\V$ with the characteristic function of $\Omega'$.
Due to Chapter~\ref{chapter:local} $\upsilon=(\unity-\Op{I}_{\Omega}\V')^{-1}\alpha$ with $\alpha\in\qat\setminus\{0\}$ is $\V'$\=/holomorphic on $\Omega$ and is holomorphic on $\Omega\setminus\Bar{\Omega}'$.
For this $\upsilon$, we define $\B$ and $\U$ in this step to be weak limits of sequences of smooth potentials and prove $(\partial+\B+\U)\upsilon=0$, i.e.\ the equation (1)(a)--(1)(b) (compare~\eqref{derivative log}).
\begin{lemma}\label{step 2}
For $\V\in\banach{2}\loc(\Omega,\qat^-)$ each point of $\Omega$ has two open neighborhoods $\Omega'\subset\Omega''$ with smooth boundaries and compact closures in $\Omega''$ and $\Omega$, respectively, such that $\upsilon=(\unity-\Op{I}_{\Omega}\V')^{-1}\alpha$ has no root on $\Bar{\Omega}''$ for $\alpha\in\qat\setminus\{0\}$, where $\V'$ is the product of $\V$ with the characteristic function of $\Omega'$.
Furthermore, for some $\B\in\banach{2,\infty}(\Omega'',\qat^+)$ and $\U\in\banach{2}(\Omega'',\qat^-)$ we have $(\partial+\B+\U)\upsilon=0$ on $\Omega''$.
\end{lemma}
\begin{proof}
Without loss of generality we assume that $\Omega$ is bounded.
The statement is independent of the choice of $\alpha\in\qat\setminus\{0\}$, because $\ker(\barpartial -\V')$ is a right $\qat$\=/vector space.
We fix $\alpha=1$, choose $2<q<\infty$ and for each point of $\Omega$ a neighborhood $\Omega''$ with compact closure in $\Omega$.
We now show that $\upsilon$ as defined in the Lemma has no roots on $\Bar{\Omega}''$ if $\|\V'\|_2$ is small enough, which we achieve by choosing a suitable $\Omega'$.

Let $\Spa{A}^q(\Omega,\qat)$ be the Bergman space $H^0(\Omega,\Q{\unity,0})\cap\banach{q}(\Omega,\qat)$ with norm $\|\cdot\|_q$ of $\banach{q}(\Omega,\qat)$.
The continuous function
\begin{align*}
\Bar{\Omega}''&\to\mathbb{R},&z_0&\mapsto\inf\{\|\alpha-(z-z_0)\xi\|_q\mid\xi\in\Spa{A}^q(\Omega,\qat)\}
\end{align*}
has a positive minimal value $C_1>0$.
By distinguishing between the cases $\|(z-z_0)\xi\|_q\le2\|\alpha\|_q$ and  $\|(z-z_0)\xi\|_q>2\|\alpha\|_q$ we conclude
\[\|\alpha-(z-z_0)\xi\|_q\ge\min\left\{\tfrac{C_1}{2\|\alpha\|_q},\tfrac{1}{2}\right\}\|(z-z_0)\xi\|_q.\]
The classical Cauchy's Integral Formula implies $\|\xi\|_q\le C_2\|(z-z_0)\xi\|_q$ for some $C_2>0$ and all $(z_0,\xi)\in\Bar{\Omega}''\times \Spa{A}^q(\Omega,\qat)$.
This implies for all $(z_0,\xi)\in\Bar{\Omega}''\times \Spa{A}^q(\Omega,\qat)$
\begin{align}\label{lower bound}
\|\alpha - (z-z_0)\xi\|_q&\ge\max\left\{C_3\|\xi\|_q,\tfrac{C_1+C_3\|\xi\|_q}{2}\right\},
&C_3&=\tfrac{1}{2C_2}\min\left\{\tfrac{C_1}{\|\alpha\|_q},1\right\}.
\end{align}
For sufficiently small $\|\V'\|_2$ the operators $\unity-\Op{I}_\Omega\V'$ and $\unity-\Op{I}_\Omega\Tilde{\V}'$ with $\Tilde{\V}' =(z-z_0)^{-1}\,\V'\,(z-z_0)$ are on $\banach{q}(\Omega,\qat)$ invertible with inverse operators nearby $\unity$.
By the triangle inequality we have
\begin{multline*}
\|\alpha-(z-z_0)\xi\|_q \leq \| (\unity-\Op{I}_\Omega\V')^{-1}\alpha - (z-z_0)(\unity-\Op{I}_\Omega\Tilde{\V}')^{-1}\xi \|_q+\\
+\|\alpha-(\unity-\Op{I}_\Omega\V')^{-1}\alpha \|_q+\| (z-z_0)\xi-(z-z_0)(\unity-\Op{I}_\Omega\Tilde{\V}')^{-1}\xi \|_q.
\end{multline*}
Therefore the distance between $\upsilon$ and an arbitrary $\V'$\=/holomorphic function that has a root at $z_0$ can be estimated as
\begin{multline*}
\|(\unity-\Op{I}_\Omega\V')^{-1}\alpha-(z-z_0)(\unity-\Op{I}_\Omega\Tilde{\V}')^{-1}\xi\|_q\ge\|\alpha-(z-z_0)\xi\|_q-\\-\left\|\unity-(\unity-\Op{I}_\Omega\V')^{-1}\right\|\cdot\|\alpha\|_q-\left\|\unity-(\unity-\Op{I}_\Omega\Tilde{\V}')^{-1}\right\|\cdot\|\xi\|_q  .
\end{multline*}
If $\left\|\unity-(\unity-\Op{I}_\Omega\V')^{-1}\right\| < \tfrac{C_1}{2\,\|\alpha\|_q}$ and $\left\|\unity - (\unity-\Op{I}_\Omega\qj \Tilde{\V}')^{-1}\right\| < \tfrac{C_3}{2}$ then due to the lower bound~\eqref{lower bound} this distance does not vanish on $(z_0,\xi)\in\Bar{\Omega}''\times \Spa{A}^q(\Omega,\qat)$.
Thus $\upsilon$ has no root on $\Bar{\Omega}''$ for sufficiently small $\Omega'$.

To prove the second part of the lemma we construct $\B$ and $\U$ as limits by approximating $\V'$ by smooth potentials.
Choose a sequence $\V_n$ of smooth potentials with support in $\Omega'$ which converge in $\banach{2}(\Omega,\qat^-)$ to $\V'$.
The corresponding sequence $\upsilon_n=(\unity-\Op{I}_{\Omega}\V_n)^{-1}\alpha$ converges for some $2<q<\infty$ in $\banach{q}(\Omega,\qat)$ to $\upsilon$.
By definition of $\V_n$ the functions $\upsilon_n$ are holomorphic on $\Omega\setminus\Bar{\Omega}'$.
By passing to a subsequence we may assume that $\upsilon_n$ converge uniformly on $\boundary \Omega''$ to $\upsilon$ and have no root in $\Bar{\Omega}''$.
Therefore there exist potentials $\U_n$ and $\B_n$ so that $\B_n+\U_n= -\boundary \upsilon_n\cdot\upsilon_n^{-1}$.
Furthermore, after passing to a subsequence the corresponding $1$\=/forms $\B_n\, dz$ converge in a neighborhood of the smooth boundary $\boundary \Omega''$.
Equation~\eqref{eq:zero curvature 1} implies for $\V^+=0$
\[
d(\B_n dz)
= 2\ci (|\V_n|^2-|\U_n|^2)\dmu.
\]
Since the sequence of measures $\V_n\Bar{\V}_n\dmu$ converges, the sequence $\U_n$ is a bounded sequence in $\banach{2}(\Omega'',\qat^-)$.
Due to the Banach-Alaoglu theorem~\cite[Theorem~IV.21]{RS1}, this sequence $\U_n$ has a weakly convergent subsequence with limit $\U\in\banach{2}(\Omega'',\qat^-)$.
By the Lemma~\ref{generalized young} the integral kernel $(z-w)^{-1}\,\tfrac{1}{\pi}\dmu(w)$ defines a bounded operator $\Op{I}_{\Omega''}\in\mathcal{L}(\banach{2,1}(\Omega'',\qat^-),C_b(\Omega'',\qat^-))$.
Note that we use the same symbol $\Op{I}_{\Omega''}$ to denote the restriction to $\qat^-$\=/valued functions of the analogous operator defined at the beginning of Chapter~\ref{chapter:local}.
Due to~\cite[Chapter~2 Theorem~2.7~and Chapter~4 Corollary~4.8]{BS} the Lorentz space $\banach{2,\infty}(\Omega'',\qat)$ is the dual space of $\banach{2,1}(\Omega'',\qat)$.
Hence $\Op{I}_{\Omega''}$ is also bounded as an operator from the space of signed measures on $\Omega''$ into $\banach{2,\infty}(\Omega'',\qat)$.
Consequently a subsequence of $\Op{I}_{\Omega''}(|\V_n|^2-|\U_n|^2)$ converges weakly in $\banach{2,\infty}(\Omega'',\qat)$.
Since the $\upsilon_n$ are holomorphic on $\Omega\setminus\Bar{\Omega}'$ without roots on $\Bar{\Omega}''\setminus\Bar{\Omega}'$, the inverses $\upsilon_n^{-1}$ are bounded on $\boundary \Omega''$.
By the Maximum Principle and by Montel's theorem a subsequence of the sequence of holomorphic functions $B_n-\Op{I}_{\Omega''}(|\V_n|^2-|\U_n|^2)$ converges uniformly on $\Bar{\Omega}''$.
Hence a subsequence of $\B_n$ converges weakly in $\banach{2,\infty}(\Omega'',\qat^+)$, and the limit $\B$ obeys $\barpartial\B=|\V|^2-|\U|^2$ on $\Omega''$.
Since $\upsilon_n$ converges in $\banach{q}(\Omega'',\qat)$ for some $2<q<\infty$ and since subsequences of $\U_n$ and $\B_n$ both converge weakly, a subsequence of $(\B_n+\U_n)\upsilon_n=-\partial\upsilon_n$ converges weakly in $\banach{p}(\Omega'',\qat)$ to $(\B+\U)\upsilon$ for $1<p<\frac{2q}{2+q}$.
At the beginning of Chapter~\ref{chapter:local} the operator $\Op{I}_\Omega$ on the very left of the Neumann series~\eqref{eq:neumann} of $\upsilon_n$ is shown to be bounded in $\mathcal{L}(\banach{p}(\Omega,\qat),\sobolev{1,p}(\Omega,\qat))$.
Thus $\upsilon_n$ converges in $\sobolev{1,p}(\Omega'',\qat)$, and the limit $\upsilon$ belongs on $\Omega''$ to $\ker(\partial+\B+\U)$.
\end{proof}

\proofstep{Step B.} 
In this step we prove that $\chi=(\overline{\qj\upsilon})^{-1}$ is  $\qk\Bar{\U}\qk$\=/holomorphic, i.e.\ (1)(d).
The essential point requires us to show that $\upsilon^{-1}$ lies in $\banach{q}\loc(\Omega,\qat)$ for some $2<q<\infty$, which is accomplished by controlling $-2\barpartial \partial\ln|\upsilon|\dmu$.
In particular we show that this is a signed measure without point measures.
As we are considering the logarithm of $\upsilon$, it is natural to use the Zygmund space $\banach{}_{\exp}(\Omega,\mathbb{R})$~\cite[Chapter~4 Section~6]{BS}\index{Zygmund space}, which is the space of functions $f$ such that there exists a constant $C(f)$ with $\exp(C|f|)$ integrable.
Its dual is $\banach{}\log\banach{}(\Omega,\mathbb{R})$~\cite[Chapter~4 Theorem~6.5]{BS}, which is the space of functions $f$ such that $|f|\cdot \max\{0,\ln|f|\}$ is integrable.
These are both rearrangement invariant spaces, with norms defined in terms of the non-increasing rearrangement of the function.
Given a function $f$, let $\mu_f$ denote the associated distribution function.
Its non-increasing rearrangement $f^{\ast}$ (\cite[Chapter~II \S3. Chapter~V\S3]{SW},~\cite[Chapter~2 Section~1]{BS} and~\cite[Chapter~1 Section~8]{Zi}) is defined by 
\[f^{\ast}(t) = \inf \{ s \in [0,\infty) \mid \mu_f(s) \leq t \},\]
and its maximal function~\cite[Chapter~2 Definition~3.1]{BS} is defined by
\[f^{\ast\ast}(t) = \frac{1}{t}\int_0^t f^{\ast}(s)\,ds.\]
The norms on $\banach{}_{\exp}(\Omega,\mathbb{R})$ and on $\banach{}\log\banach{}(\Omega,\mathbb{R})$ are given by
\begin{align}\label{eq:normLexpLlogL}
\|g\|_{\banach{}_{\exp}}
&= \sup_{0<t<|\Omega|} \frac{g^{\ast\ast}(t)}{1-\ln(t/|\Omega|)},
&\|g\|_{\banach{}\log\banach{}}
&= \int_{0}^{|\Omega|}g^{\ast\ast}(t)\,dt,
\end{align}
see~\cite[Chapter~4 (6.10) and (6.9)]{BS}.
With these preliminaries covered, we now state a lemma in preparation of the proof of the main statement of this step, Corollary~\ref{step 3}.
\begin{lemma}\label{zygmund estimate}
Let $\Omega$ denote a bounded open subset of $\mathbb{C}$, $d\sigma$ any finite signed Baire measure on $\Omega$~\cite[Chapter~13 Section~5]{Ro2}, and $d\sigma^+$ the positive part of the Hahn decomposition of $d\sigma$~\cite[Chapter~11 Section~5]{Ro2}.
Then there exists a solution $h$ in the Zygmund space $\banach{}_{\exp}(\Omega,\mathbb{R})$ of the equation
\begin{align}\label{eq:liouville}
-\barpartial \partial h&=d\sigma
\quad\text{in the sense of distributions.}
\end{align}
Moreover, if for some $\varepsilon>0$ all $\varepsilon$\=/balls of $\Omega$ have measure smaller than $\pi/q$ with respect to $d\sigma^+$, then $\exp(h)\in\banach{q}\loc(\Omega,\mathbb{R})$.
Conversely, if the positive part $d\sigma^+$ contains a point measure with mass $\geq \pi/q$, then $\exp(h)\not\in\banach{q}\loc(\Omega,\mathbb{R})$.
\end{lemma}
\begin{proof}
The strategy to prove the first statement is to construct a right inverse of $-\barpartial \partial$ in $\mathcal{L}(\banach{}\log\banach{}(\Omega,\mathbb{R}),C_b(\Omega,\mathbb{R}))$, whose dual operator maps the signed measures to $\banach{}_{\exp}(\Omega,\mathbb{R})$.
Due to Dolbeault's Lemma 
the convolution with the function $-\frac{2}{\pi}\ln|z|$ is a right inverse of $-\barpartial \partial$.
We show first that the restriction of this convolution operator to $\Omega$ belongs to $\mathcal{L}(\banach{1}(\Omega,\mathbb{R}),\banach{}_{\exp}(\Omega,\mathbb{R}))$.
Since $\Omega$ is bounded, the claim is equivalent to the analogous statement about the restriction of the convolution with the positive part of $-\tfrac{2}{\pi}\ln |z|$, i.e.\ with the function 
\[f(z)=\begin{cases}-\frac{2}{\pi}\ln|z| & \text{if }|z|<1\\0 & \text{if }|z|\geq 1.\end{cases}\]
Moreover as $f$ is non-negative, $f\ast d\sigma^-$ is non-positive, and therefore $\exp(f\ast d\sigma^-)$ is bounded by $1$.
Therefore the question whether $f\ast d\sigma$ belongs to $\banach{q}\loc(\Omega,\mathbb{R})$ depends only on $d\sigma^+$.
Due to~\cite[Chapter~4 Theorem~6.5]{BS}, $\banach{}_{\exp}(\Omega,\mathbb{R})$ is the dual space of the Zygmund space $\banach{}\log\banach{}(\Omega,\mathbb{R})$, and due to the Riesz Representation Theorem~\cite[Chapter~13 Section~5 Theorem~25]{Ro2} the space of signed measures on $\Omega$ is the dual space of the Banach space $C_b(\Omega,\mathbb{R})$.
Hence to prove the first part of the lemma it is sufficient to show that convolution with $f$ defines a bounded operator from $\banach{}\log\banach{}(\Omega,\mathbb{R})\subset\banach{1}(\Omega,\mathbb{R})$ into $C_b(\Omega,\mathbb{R})\subset\banach{}_{\exp}(\Omega,\mathbb{R})$.

Let us start by computing the distribution function $\mu_f$ of the above function $f$ and its non-increasing rearrangement $f^{\ast}$:
\begin{align*}
\mu_f(s)&=\pi\exp\left(-\pi s\right) &f^{\ast}(t)&=\begin{cases}-\frac{\ln\left(t/\pi\right)}{\pi} &\text{if }0\leq t\leq\pi\\
0 & \text{if } \pi\leq t.\end{cases}
\end{align*}	
For any $g\in\banach{}\log\banach{}(\Omega,\mathbb{R})$ we have a bound on the maximal function of a convolution $h = f \ast g$.
The below calculation begins with Equations~\cite[(1.8.14) and (1.8.15)]{Zi} in the proof of~\cite[1.8.8.~Lemma]{Zi} and is adapted from~\cite[Lemma~1.5]{O}.
\begin{align*}
h^{\ast}(t)
\leq h^{\ast\ast}(t)
&\leq g^{\ast\ast}(t)\int_{f^{\ast}(t)}^{\infty} \mu_{f}(s)\,ds + tg^{\ast\ast}(t)f^{\ast}(t) + \int_t^{\pi} f^{\ast}(s)g^{\ast}(s)\,ds \\ 
&= g^{\ast\ast}(t) \min\left\{\frac{t}{\pi},1\right\} - \int_{t}^{\pi} sg^{\ast\ast}(s)df^{\ast}(s) \\
&= \min\left\{\frac{1}{\pi},\frac{1}{t}\right\} \int_0^t g^{\ast}(s)ds + \frac{1}{\pi} \int_{t}^{\pi} g^{\ast\ast}(s)ds \\
&\leq\frac{1}{\pi}\int_{0}^{t}g^{\ast}(s)\,ds + \frac{1}{\pi}\int_{t}^{|\Omega|}g^{\ast\ast}(s)\,ds \leq\frac{1}{\pi}\|g\|_{\banach{}\log\banach{}}.
\labelthis{eq:convolution bound}
\end{align*}
Consequently $h\in\banach{\infty}(\Omega,\mathbb{R})$ with $\|h\|_\infty\leq\frac{1}{\pi}\|g\|_{\banach{}\log\banach{}}$ and hence convolution with $f$ belongs to $\mathcal{L}(\banach{}\log\banach{}(\Omega,\mathbb{R}),C_b(\Omega,\mathbb{R}))$.
This proves the first part of the lemma.

Due to Weyl's Lemma~\cite[Theorem~IX.25]{RS2} the difference of two solutions $h$ of~\eqref{eq:liouville} is analytic.
Therefore, it suffices to show the second and third statement of the lemma on all $\epsilon$\=/balls in $\Omega$ for some $\epsilon>0$.
Because the norm of the convolution operator $\banach{}\log\banach{}(\Omega,\mathbb{R}) \to C_b(\Omega,\mathbb{R})$ is $\leq \tfrac{1}{\pi}$ by the above calculation, and by~\cite[Chapter~4 (6.12)]{BS} $\|\,\cdot\,\|_{L_{\exp}}$ (see Equation~\eqref{eq:normLexpLlogL}) is bounded by the norm of the dual space of $\banach{}\log\banach{}(\Omega,\mathbb{R})$, the norm of the dual operator from the space of signed measures on $\Omega$ to $\banach{}_{\exp}(\Omega,\mathbb{R})$ is also $\leq \tfrac{1}{\pi}$.
It follows that if the measure of $\Omega$ with respect to a finite positive measure $d\sigma$ is smaller than $\pi/q$, then the  non-increasing rearrangement $h^\ast$ and its maximal function $h^{\ast\ast}$ of solutions of~\eqref{eq:liouville} obey
\begin{align*}
h^{\ast}(t)\leq h^{\ast\ast}(t)&
\leq\frac{\int_\Omega d\sigma}{\pi}\left(1-\ln\frac{t}{|\Omega|}\right)
\quad\text{for }t\in[0,|\Omega|].
\end{align*}
A standard argument~\cite[Chapter~2 Exercise~3]{BS} and the fact that the non-increasing rearrangement $h^{\ast}(t)$ vanishes for $t>|\Omega|$ show that 
\begin{align*}
\int_{\Omega} (\exp(q|h|)-1) \,\dmu &=\int_{0}^{\infty}\left(\exp\left(q h^{\ast}(t)\right)-1\right)\,dt \leq\int_{0}^{|\Omega|}\exp\left(q h^{\ast}(t)\right)\,dt\\&\leq |\Omega|\exp\left(\frac{q\int_\Omega d\sigma}{\pi}\right)\int_{0}^{|\Omega|}\left(\frac{t}{|\Omega|}\right)^{-\frac{q\int_\Omega d\sigma}{\pi}}dt.
\end{align*}
Thus the finiteness of the right-hand-side integral is equivalent to the statement that $\exp|h|$ belongs to $\banach{q}(\Omega,\mathbb{R})$.
To sum up, the exponential $\exp(f\ast d\sigma)$ belongs to $\banach{q}(\Omega,\mathbb{R})$ if $q<\frac{\pi}{\int_\Omega d\sigma^+}$.
This shows the second statement of the lemma.
The third statement follows by direct calculation of the convolution of $f$ with point measures.
\end{proof}
\begin{corollary}\label{step 3}
In the situation of Lemma~\ref{step 2} $\barpartial\B \dmu=(|\V|^2-|\U|^2)\dmu$ is a signed measure on $\Omega''$ without point measures and $\chi=(\overline{\qj\upsilon})^{-1}$ is $\qk\Bar{\U}\qk$\=/holomorphic on $\Omega''$.
\end{corollary}
\begin{proof}
As in the proof of Lemma~\ref{step 2} we may assume that $\Omega$ is bounded and consider the same sequences $\V_n$ and $\upsilon_n=(\unity-\Op{I}_\Omega\V_n)^{-1}\alpha$.
As shown in the proof of this lemma, a subsequence of the sequence of measures $(|\V_n|^2-|\U_n|^2)\dmu$ converges weakly in the dual space of $C_b(\Omega'',\mathbb{R})$ to a signed measure $d\sigma$.
Due to Lemma~\ref{zygmund estimate} the corresponding subsequence $\ln|\upsilon_n|$ differs from a weakly convergent sequence in $\banach{}\log\banach{}(\Omega'',\mathbb{R})$ by a sequence of harmonic functions in $\Omega''$.
Moreover, we have chosen the sequence $\V_n$ in such a way that $\ln|\upsilon_n|$ is bounded on $\boundary \Omega''$.
Therefore, after passing to a subsequence, $\ln|\upsilon_n|$ converges weakly in $\banach{}\log\banach{}(\Omega'',\mathbb{R})$.

Since $\V_n$ converges in $\banach{2}(\Omega'',\qat^-)$, only the weak limit of $|\U_n|^2\,\dmu$ can have point measures.
Thus only the part $d\sigma^-$ of the Hahn decomposition~\cite[Chapter~11 Section~5]{Ro2} $d\sigma=d\sigma^+-d\sigma^-$ can have point measures.

For any $m\in \mathbb{N}_+$ and $\Tilde{\upsilon}_n=(z-z_0)^{-l}\upsilon_n$, the measure $2\partial \barpartial \ln \left(\tfrac{|\Tilde{\upsilon}_n|}{|{\upsilon}_n|}\right) \dmu$ is a point measure of mass $m\pi$ at $z_0$ and we have
\[2\barpartial \partial \ln |\Tilde{\upsilon}_n| = 2\barpartial \partial \ln |\upsilon_n| + m\pi\,\delta(z-z_0).\]
If $d\sigma^-$ contains a point measure at $z=z_0\in\Bar{\Omega}''$ of mass $\ge m\pi$, then adding the point measure $m\pi \,\delta(z-z_0)$ to $d\sigma$ does not create a point measure in the positive part $d\Tilde{\sigma}^+$ of the Hahn decomposition of the weak limit $d\Tilde{\sigma}$ of $d\Tilde{\sigma}_n=2\barpartial\partial\ln|\Tilde{\upsilon}_n|\dmu$.
In particular, for every $z_1 \in \Bar{\Omega}''$ we have $\lim_{\varepsilon\downarrow 0} \Tilde{\sigma}^+(B(z_1,\varepsilon)) = 0$.
Hence we may cover for any $2<q<\infty$ the compact set $\Bar{\Omega}''$ by finitely many balls $B(z_1,\varepsilon_1),\dotsc, B(z_L,\varepsilon_L)$ with $z_l \in \Bar{\Omega}''$ and $\Tilde{\sigma}^+(B(z_l,\varepsilon_l)) < \pi/q$.
Furthermore the same is true for all sufficiently large $n$ if we replace $d\Tilde{\sigma}^+$ by the positive part of $d\Tilde{\sigma}_n$.
There exists $\varepsilon>0$ such that all $B(z,\varepsilon)$ with $z\in\Bar{\Omega}''$ are contained in one of the finitely many balls $B(z_l,\varepsilon_l)$.
Due to Lemma~\ref{zygmund estimate}, for these $n$ the function $\Tilde{\upsilon}_n$ is bounded in $\banach{q}\loc(\Omega'',\qat)$.
Since $\upsilon_n=(z-z_0)^l\Tilde{\upsilon}_n$ converges for some $2<q<\infty$ in $\banach{q}(\Omega'',\qat)$ to $\upsilon$, a subsequence of $\Tilde{\upsilon}_n$ converges weakly to $(z-z_0)^{-l}\upsilon$ in $\banach{q}(\Omega'',\qat)$.
This would imply that $\upsilon$ has a root of order at least $m$ at $z_0$.
Hence, due to our assumptions, the masses of all point measures of $d\sigma^-$ are smaller than $\pi$.

In particular, for any $z\in\Bar{\Omega}''$ the limit $\lim_{\varepsilon \downarrow 0}\sigma^-(B(z,\varepsilon))$ is smaller than $\pi$.
By the same argument as before, for some $2<q<\infty$,  Lemma~\ref{zygmund estimate} implies that for sufficiently large $n$, the $(\overline{\qj\upsilon_n})^{-1}$ are bounded in $\banach{q}\loc(\Omega'',\qat)$.
  
The product of a weak limit of a subsequence of $(\overline{\qj\upsilon_n})^{-1}$ with the norm limit $\overline{\qj\upsilon}$ in $\banach{q}(\Omega'',\qat)$ of the sequence $\overline{\qj\upsilon_n}$ is equal to the weak limit of the products, which is the constant function $1$.
Hence $(\overline{\qj\upsilon})^{-1}$ belongs to $\banach{q}(\Omega'',\qat)$ for the same $2<q<\infty$ as before.
Since $\upsilon$ belongs to $\ker(\barpartial -\V)\cap\ker(\partial+\B+\U)\subset\bigcap_{1<p<2}\sobolev{1,p}(\Omega'',\qat)$, the following map describes the $1$\=/current induced by $d\chi\cdot\chi^{-1}$ on the space of smooth $1$\=/forms $\omega$ with compact support in $\Omega''$:
\begin{align*}
\omega\mapsto&-\int_{\Omega''}\chi\,d(\chi^{-1}\omega) = \int_{\Omega''}\qj\Bar{\upsilon}^{-1}\,d(\Bar{\upsilon}\qj\omega) = \int_{\Omega''}\qj\Bar{\upsilon}^{-1}\left(d\Bar{\upsilon}\qj \wedge\omega+\Bar{\upsilon}\qj d\omega\right)\\
\labelthis{eqn:derivative of inverse}
&=\int_{\Omega''} \qj\,\overline{d\upsilon\cdot\upsilon^{-1}}\qj\wedge\omega-\int_{\Omega''} d\omega= - \int_{\Omega''} \qj\left(\overline{dz(\B+\U)-d\Bar{z}\V}\right)\qj\wedge\omega\\
&=\int_{\Omega''} \left(dz(\B+\qj\Bar{\V}\qj)-d\Bar{z}\qj\Bar{\U}\qj\right)\wedge\omega  =\int_{\Omega''} \left(dz(\B-\qk\Bar{\V}\qk)+d\Bar{z}\qk\Bar{\U}\qk\right)\wedge\omega.
\end{align*}
So $\chi$ fulfills the analogue of the condition of Theorem~\ref{cauchy formula} (iv).
Therefore $\chi$ belongs on $\Omega''$ to $\ker(\barpartial-\qk\Bar{\U}\qk)\subset\bigcap_{2<q<\infty}\banach{q}(\Omega'',\qat)$, and Lemma~\ref{zygmund estimate} implies that $d\sigma^-$ contains no point measures.
\end{proof}
	
\proofstep{Step C.}
Now we show  $\psi=(\partial+\B+\U)\phi$ is $\U$\=/holomorphic and~\eqref{eq:derivative quotient final} on $\Omega'$.
Equation~\eqref{eq:derivative quotient final} implies that the partial derivatives of $\upsilon^{-1}\phi$ belong to $\bigcap_{1<p<2}\sobolev{1,p}\loc(\Omega',\qat)$, and thus $\upsilon^{-1}\phi$ belongs to 
\[\bigcap_{1<p<2}\sobolev{2,p}\loc(\Omega',\qat)\subset\bigcap_{2<q<\infty}\sobolev{1,q}\loc(\Omega',\qat)\subset C(\Omega',\qat),\]
by the Sobolev Embedding theorem.
So Lemma~\ref{step 4} implies (1)(e) and (1)(f).
\begin{lemma}\label{step 4}
In the situation of Lemma~\ref{step 2} and Corollary~\ref{step 3} let $\phi$ be $\V$\=/holomorphic on $\Omega''$.
Then the derivative $(\partial+\B+\U)\phi$ is $\U$\=/holomorphic on $\Omega'$ and~\eqref{eq:derivative quotient final} holds with $\upsilon=(\unity-\Op{I}_{\Omega}\V')^{-1}\alpha$, $\alpha\in\qat\setminus\{0\}$.
\end{lemma}
\begin{proof}
We use again the sequence of smooth potentials $\V_n\in\banach{2}(\Omega,\qat^-)$ with limit $\V'$ and the corresponding sequence of functions $\upsilon_n$ on $\Omega$.
Let $\U_n$ and $\B_n$ denote the corresponding sequences with $\B_n+\U_n=-(\partial\upsilon_n)\upsilon_n^{-1}$.
The sequence $\upsilon_n$ converges in $\banach{q}(\Omega,\qat)$ and is holomorphic on $\Omega\setminus\Bar{\Omega}'$.
Therefore it converges on compact subsets of $\Omega\setminus\Bar{\Omega}'$ with all derivatives to $\upsilon$.
In the limit $\|\V'\|_2\to 0$ the sequence converges with all derivatives on compact subsets of $\Omega\setminus\Bar{\Omega}'$ to $q$.
As in the proof of Lemma~\ref{step 2} we choose $\Omega''$ with compact closure in $\Omega$ containing $\Bar{\Omega}'$ and with smooth boundary $\partial\Omega''$.
The integrals of $\B dz$ along this boundary $\partial\Omega''$ converge for $\|\V'\|_2\to0$ to zero.
For sufficiently small $\Omega'$ and $\Omega''$ the sequence $\U_n$ stays in the domain of the map in~\eqref{eq:integral kernel} and $\Op{I}_{\Omega'',\U_n}$ converges by Lemma~\ref{weakly continuous} in $\mathcal{L}(\banach{p}(\Omega'',\qat),\banach{q}(\Omega'',\qat))$.

Since we assume that $\phi$ belongs on $\Omega''$ to $\ker(\barpartial -\V)$ and since $\Omega'$ has compact closure in $\Omega''$, it follows that $\phi$ belongs to $\bigcap_{1<p<2} \sobolev{1,p}(\Omega',\qat)$ rather than only to $\bigcap_{1<p<2} \sobolev{1,p}\loc(\Omega',\qat)$.

For the remainder of this proof, we restrict all functions to $\Omega'$.
We apply Chapter~\ref{chapter:local} and obtain $\phi$ for any $1<p<2$ as the limit of
\begin{align*}
\phi_n=(\unity-\Op{I}_{\Omega'}\V_n)^{-1}\alpha'
\quad\text{with}\quad
\alpha'=(\unity-\Op{I}_{\Omega'}\V)
\phi\in\sobolev{1,p}(\Omega',\qat).
\end{align*}
Clearly $(\upsilon_n,\phi_n)$ is a pair of $\V_n$\=/holomorphic functions which converges to $(\upsilon,\phi)$.
Furthermore, the bounded sequence $(\partial+\B_n+\U_n)\phi_n$ belongs to the kernel of $\barpartial-\U_n$ and a subsequence converges weakly in $\banach{p}(\Omega',\qat)$ to $(\partial+\B+\U)\phi$.

In order to prove that this function is $\U$\=/holomorphic it suffices to verify Cauchy's Integral Formula Theorem~\ref{cauchy formula}~(ii) with any test function $f\in C_0^\infty(\Omega',\mathbb{R})$ that is $1$ on some subset $\Omega'''\subset\Omega'$ with compact closure in $\Omega'$.
The convergent sequence of compact operators $\Op{I}_{\Omega',\U_n}$ maps the bounded sequence $(\partial+\B_n+\U_n)((\barpartial f)\phi_n)$ to a convergent sequence in $\banach{q}\loc(\Omega',\qat)$ with limit $(\partial+\B+\U)\phi$ for any $2<q<\tfrac{2p}{2-p}$.
Hence the limit belongs on such $\Omega'''$ to the image of a function with support in $\Omega'\setminus\Bar{\Omega}'''$ with respect to the operator $\Op{I}_{\Omega',\U^-}$ and belongs on $\Omega'''$ to the kernel of $\barpartial -\qj{U}$.
Due to Step~B for sufficiently small $\Omega'$ the sequence $\ln|\upsilon_n|$ is bounded in $\banach{}_{\exp}(\Omega',\mathbb{R})$ by any $\epsilon>0$.
Hence $\upsilon_n^{-1}$ and $\upsilon_n^{-1}\phi_n$ converge in $\banach{q}(\Omega',\qat)$ weakly to $\upsilon^{-1}$ and $\upsilon^{-1}\phi$, respectively.
Now~\eqref{eq:derivative quotient final} follows from the corresponding formulas for $\upsilon_n^{-1}\phi_n$.
\end{proof}
	
\proofstep{Step D.}
In this step we show that the continuous function $\upsilon^{-1}\phi$ has the same roots as $\phi$.
\begin{lemma}\label{step 5}
\index{Root}
In the situation of Lemma~\ref{step 2}, Corollary~\ref{step 3} and Lemma~\ref{step 4}, for any $\V$\=/holomorphic $\phi\in\banach{q}(\Omega',\qat)$, the roots of the continuous function $\upsilon^{-1}\phi$ coincide with the roots of $\phi$ in the sense of Definition~\ref{order of roots}.
\end{lemma}
\begin{proof}
Due to Lemma~\ref{quotient dimension} for any $z_0\in \Omega'$ the quotient of $\ker(\barpartial -\V)$ by the subspace of functions with a root at $z_0$ is a one-dimensional $\qat$\=/right-linear space.
Since $\upsilon$ has no roots on $\Omega'$, there thus exists a unique $\Tilde{\alpha}\in\qat$ such that $\Tilde{\phi}=\phi-\upsilon \Tilde{\alpha}$ has in the sense of Definition~\ref{order of roots} a root at $z_0$, i.e.\ $(z-z_0)^{-1}\Tilde{\phi}\in \banach{q}(\Omega',\qat)$.
Note that $\phi$ has a root at $z_0$ if and only if $\Tilde{\alpha}=0$ holds.
We need to show that $\Tilde{\alpha}=0$ holds if and only if the continuous function $\upsilon^{-1}\phi$ has a root at $z_0$.
In fact, we will show that $\Tilde{\alpha}$ is the value of $\upsilon^{-1}\phi$ at $z_0$.
If $(\upsilon^{-1}\phi)(z_0)\neq \Tilde{\alpha}$, then the continuous function $\upsilon^{-1}\Tilde{\phi}=\upsilon^{-1}\phi-\Tilde{\alpha}$ would not have a root at $z_0$ and there would exist an open ball $B(z_0,\epsilon) \subset \Omega'$ on which $\upsilon^{-1}\Tilde{\phi}$ does not have a root.
Then $(\upsilon^{-1}\Tilde{\phi})^{-1}$ would belong to $C(B(z_0,\epsilon),\qat)$.
This would imply that $(z-z_0)^{-1}\upsilon=(z-z_0)^{-1}\Tilde{\phi}(\upsilon^{-1}\Tilde{\phi})^{-1}$ belongs to $\banach{q}\loc(B(z_0,\epsilon),\qat)$ and therefore $\upsilon$ would have a root at $z=z_0$, which contradicts the assumptions.
\end{proof}

\proofstep{Step E.}
In Steps~A--C we proved the statements (1)(a)-(1)(b) and (1)(d)-(1)(f) only in the special case $\upsilon=(\unity-\Op{I}_{\Omega}\V')^{-1}\alpha$ with $\alpha\in\qat\setminus\{0\}$.
These arguments showed that any $\V$\=/holomorphic $\phi$ is the product $\phi=\upsilon(\upsilon^{-1}\phi)$ of this special $\upsilon$ with the continuous function $\upsilon^{-1}\phi\in\bigcap_{2<q<\infty}\sobolev{1,q}\loc(\Omega',\qat)$.
We distinguish the general case $(\Tilde{\phi},\Tilde{\upsilon})$ from the special case $(\phi,\upsilon)$ of Steps~A--C.
In order to prove the statements (1)(a)-(1)(b) and (1)(d)-(1)(f) in the general case, we apply the results of Steps~A--C first with the special $\upsilon$ and $\phi=\Tilde{\upsilon}$ and then with $\upsilon$ and $\phi=\Tilde{\phi}$.
  
If $\phi=\Tilde{\upsilon}$ has no root, then by Lemma~\ref{step 5} $(\upsilon^{-1}\phi)^{-1}$ is continuous without roots and belongs by Step~C to $\bigcap_{2<q<\infty}\sobolev{1,q}\loc(\Omega',\qat)$.
We calculate
\begin{align*}
d\phi&=d\upsilon(\upsilon^{-1}\phi)+\upsilon\,d(\upsilon^{-1}\phi),&d\phi^{-1}&=d\left((\upsilon^{-1}\phi)^{-1}\right)\upsilon^{-1}+(\upsilon^{-1}\phi)^{-1}\,d(\upsilon^{-1}).
\end{align*}
Due to the equation $(d\phi)\phi^{-1}-(d\upsilon)\upsilon^{-1}=\upsilon\left(d(\upsilon^{-1}\phi)\right)\phi^{-1}$ the corresponding potentials $\Tilde{\B}$ and $\Tilde{\U}$ with $\Tilde{\upsilon}=\phi\in\ker(\partial+\Tilde{\B}+\Tilde{U})$ obey
\begin{align}\label{estimate difference}
\Tilde{\B}-\B+\Tilde{\U}-\U&\in\bigcap_{1<p<2}\sobolev{1,p}\loc(\Omega',\qat)\subset\bigcap_{2<q<\infty}\banach{q}\loc(\Omega',\qat).
\end{align}
Hence $\Tilde{\upsilon}$ obeys (1)(a)--(1)(b) with respect to the  $\Tilde{B}\in\banach{2,\infty}\loc(\Omega',\qat)$ defined above.
Now the analogue of~\eqref{eqn:derivative of inverse} shows (1)(d) for $\Tilde{\upsilon}$.

If $\Tilde{\upsilon}=\phi$ has no roots, then by Step~D $\upsilon^{-1}\phi$ has no roots.
We write $\Tilde{\upsilon}^{-1}\Tilde{\phi}=(\upsilon^{-1}\phi)^{-1}\upsilon^{-1}\Tilde{\phi}$ and derive~(1)(f) from Step~C for $\upsilon$ and $\phi$ and $\upsilon$ and $\Tilde{\phi}$:
\begin{align*}
d(\Tilde{\upsilon}^{-1}\Tilde{\phi})&=-(\upsilon^{-1}\phi)^{-1}d(\upsilon^{-1}\phi)(\upsilon^{-1}\phi)^{-1}\upsilon^{-1}\Tilde{\phi}+(\upsilon^{-1}\phi)^{-1}d(\upsilon^{-1}\Tilde{\phi})\\&=-(\upsilon^{-1}\phi)^{-1}\Bar{\chi}\qj dz\left(\left((\partial+\B+\U)\phi\right)\phi^{-1}\Tilde{\phi}-(\partial+\B+\U)\Tilde{\phi}\right)\\&=-\phi^{-1}dz\left((\B-\Tilde{\B}+\U-\Tilde{\U})\Tilde{\phi}-(\partial+\B+\U)\Tilde{\phi}\right)\\&=-\Tilde{\upsilon}^{-1}\qj\qj dz(\partial+\Tilde{\B}+\Tilde{\U})\Tilde{\phi}=\Bar{\Tilde{\chi}}\qj dz\Tilde{\psi}.
\end{align*}
For the proof of (1)(e) we first obtain for $\Tilde{\phi}\in\ker(\barpartial -\V)$ from Step~C and~\eqref{estimate difference}
\[\Tilde{\psi}=(\partial+\Tilde{\B}+\Tilde{\U})\Tilde{\phi}=(\partial+\B+\U)\Tilde{\phi}+(\Tilde{\B}-\B+\Tilde{\U}-\U)\Tilde{\phi}\in\bigcap_{1<p<2}\sobolev{1,p}\loc(\Omega',\qat).\]
So (1)(e) follows from the analogue of~\eqref{eqn:derivative of inverse} and the following calculation
\begin{align*}
0&=\Tilde{\upsilon}d(d(\Tilde{\upsilon}^{-1}\Tilde{\phi})) = -\qj \overline{\Tilde{\chi}}^{-1} d(\overline{\Tilde{\chi}} \qj dz \Tilde{\psi})
=-\qj\left(\overline{d\Tilde{\chi} \cdot \Tilde{\chi}^{-1}} \right) \qj \wedge dz \Tilde{\psi} + d(dz \Tilde{\psi}) \\
& =-\qj\left( \overline{-\qj \left(\overline{dz(\Tilde{\B}+\Tilde{\U})-d\Bar{z}\V}\right) \qj} \right) \qj \wedge dz \Tilde{\psi} -dz\wedge d\Bar{z}\barpartial \Tilde{\psi} \\
& = (dz(\Tilde{\B}+\Tilde{\U})-d\Bar{z}\V) \wedge dz\Tilde{\psi}-dz\wedge d\Bar{z}\barpartial \Tilde{\psi} 
= - dz\wedge d\Bar{z} (\barpartial-\Tilde{\U})\Tilde{\psi}  .
\end{align*}
So the statements (1)(a)--(1)(b) and (1)(d)-(1)(f) follow from this application of Steps~A--D to a general $\Tilde{\upsilon}=\phi$ without roots.

\proofstep{Step F.}
We consider the inverse transformation $(\psi,\chi)\mapsto(\phi,\upsilon)$ and prove (2)(a)--(2)(c).
For this purpose we interchange $\upsilon\in\ker(\barpartial-\V)$ with $\chi\in\ker(\barpartial-\qk\Bar{\U}\qk)$.
The arguments in Steps~A--E demonstrating $\upsilon\in\ker(\partial+\B+\U)$ with unique $\B\in\banach{2,\infty}\loc(\Omega,\qat^+)$ and $\V\in \banach{2}\loc(\Omega,\qat^-)$ apply to $\chi$, and show $\chi\in\ker(\partial-\B+\qk\Bar{\V}\qk)$ again with unique $\B\in\banach{2,\infty}\loc(\Omega,\qat^+)$ and $\V\in\banach{2}\loc(\Omega,\qat^-)$.
We then have 
\begin{align*}
d(dF) 
& = d(\Bar{\chi} \qj dz \psi) 
= \overline{d\chi} \wedge \qj dz \psi + \Bar{\chi}\qj d(dz \psi)
= \overline{d\Bar{z}\barpartial\chi} \wedge \qj dz \psi - \Bar{\chi}\qj dz \wedge d\Bar{z}\barpartial \psi \\
&= \overline{d\Bar{z}\qk \Bar{\U} \qk \chi} \wedge \qj dz \psi - \Bar{\chi}\qj dz \wedge d\Bar{z}\U \psi 
= \Bar{\chi} \qk \U \qk \qj d\Bar{z} \wedge dz \psi - \Bar{\chi}\qj dz \wedge d\Bar{z}\U \psi \\
&= \Bar{\chi} \qk \U \qi dz \wedge d\Bar{z} \psi - \Bar{\chi}\qj dz \wedge d\Bar{z}\U \psi 
= 0.
\end{align*}
In the last step we used that $\qi dz\wedge d\Bar{z}$ is real, so commutes with $\U$.
This proves (2)(a).
The arguments in Steps~B and~E demonstrating $\chi\in\ker(\barpartial-\qk\Bar{\U}\qk)$ show $\upsilon=-(\overline{\qj\chi})^{-1}\in\ker(\barpartial -\V)$ and (2)(b).
Finally, $\phi=\upsilon F$ has the derivative $d\phi=(d\upsilon)F+\upsilon dF=(d\upsilon)F+dz\psi$ and belongs also to $\ker(\barpartial -\V)$.
This proves (2)(c).
		
\proofstep{Step G.}
Finally we prove the remaining equation (1)(c), i.e.\ the second equation of~\eqref{poisson}.
Again, we first consider $\upsilon=(\unity-\Op{I}_{\Omega}\V')^{-1}\alpha$ with $\alpha\in\qat\setminus\{0\}$.
In Step~C we proved the second equation of~\eqref{poisson} with the right hand side replaced by the weak limit of $(|\V_n|^2-|\U_n|^2)\dmu$ as a signed measure.
For any measurable $M\subset\Omega'$ the closed balls in $\banach{2}(M,\mathbb{C})$ are weakly compact and therefore weakly closed.
Therefore the $\banach{2}$\=/norm is weakly semi-continuous.
Since $\V_n$ converges in $\banach{2}(\Omega',\qat^-)$ to $\V$ and $\U_n$ weakly to $\U$, we get
\begin{equation}\label{poisson2}
\int_M 2\barpartial \partial\ln|\upsilon|\,\dmu \ge\int_M \bigl(|\U|^2-|\V|^2\bigr)\,\dmu 
\quad\text{for every measurable } M\subset\Omega''.
\end{equation}
As in Step~E, we write general $\Tilde{\upsilon}$ without roots as $\Tilde{\upsilon}=\upsilon(\upsilon^{-1}\phi)$ with $\upsilon^{-1}\phi\in\sobolev{2,p}\loc(\Omega,\qat)\subset C(\Omega,\qat)$.
This implies that~\eqref{poisson2} holds for all $\V$\=/holomorphic $\Tilde{\upsilon}$ without roots.
In the proof of (2)(c) in Step~F we applied the same line of arguments to $\chi=(\overline{\qj \upsilon})^{-1}$ and $\qk\Bar{\U}\qk$ instead of $\upsilon$ and $\V$, which shows the reverse inequality of~\eqref{poisson2}.
Together this proves the second equation of~\eqref{poisson}, which completes the proof of Theorem~\ref{thm:darboux}.
\end{proof}
An immediate consequence is the following corollary:
\begin{corollary}[Weyl's Lemma of quaternionic function theory]
\label{weyls lemma}
\index{Weyl's lemma}
Let $\U$ be a potential in $\banach{2}\loc(\Omega,\qat^-)$ on an open subset $\Omega\subset\mathbb{C}$ and $\chi$ a $\qk\Bar{\U}\qk$\=/holomorphic function without roots on $\Omega$.
Then a function $\xi\in\banach{p}\loc(\Omega,\qat)$ with $1<p<2$ is $\U$\=/holomorphic if and only if $\Bar{\chi}\qj dz\xi$ is a closed current on $\Omega$.
\end{corollary}
\begin{proof}
By Theorem~\ref{thm:darboux}~(2)(a) $\Bar{\chi}\qj dz\xi$ is a closed current for $\U$\=/holomorphic $\xi$.

Conversely, if $dF=\Bar{\chi}\qj dz\xi$ is a closed current, let $F\in\bigcap_{1<r<p}\sobolev{1,r}(\Omega',\qat)$ be an anti-derivative on a simply connected subset of $\Omega'\subset\Omega$.
Theorem~\ref{thm:darboux}~(2)(c) implies that $\upsilon=-(\overline{\qj\chi})^{-1}$ is $\V$\=/holomorphic with $\V\in\banach{2}\loc(\Omega,\qat^-)$, and Theorem~\ref{thm:darboux}~(1)(a)-(b) implies $d\upsilon=\big(d\Bar{z}\V-dz(\B+\U)\big)\upsilon$ with $\B\in\banach{2,\infty}\loc(\Omega,\qat^+)$ and $\U\in\banach{2}\loc(\Omega,\qat^-)$. We conclude that $\phi=\upsilon F$ is also $\V$\=/holomorphic and $\xi=(\partial+\B+\U)\phi$:
\[
d(\upsilon F)=(d\upsilon)F+\upsilon dF=\big(d\Bar{z}\V-dz(\B+\U)\big)\phi+\upsilon\chi\qj dz\xi=d\Bar{z}\V\phi+dz\big(\xi-(\B+\U)\phi\big).
\]
Finally Theorem~\ref{thm:darboux}~(1)(d)-(1)(e) implies that $\xi$ is $\U$\=/holomorphic.
\end{proof}
We remark that at least locally we may iterate the transformation $\phi\mapsto\psi=(\partial+\B+\U)\phi$ infinitely many times.
In fact, for any $\alpha\in\qat\setminus\{0\}$ on a small domain $\Omega'$ the function $\Tilde{\upsilon}=(\unity-\Op{I}_{\Omega'}\U|_{\Omega'})^{-1}\alpha$ is $\U$\=/holomorphic without roots on $\Omega'$.
Hence the Darboux transformation can be applied to the pair of $\U$\=/holomorphic functions $\Tilde{\upsilon}$ and $\psi$ (playing the role of $\upsilon$ and $\phi$ respectively).
This application shows that there exists $\C\in\banach{2,\infty}\loc(\Omega,\qat^+)$ and $\W\in\banach{2}\loc(\Omega',\qat^-)$ such that $\C+\W=(\partial\Tilde{\upsilon})\Tilde{\upsilon}^{-1}$.
Furthermore, the function $(\barpartial+\C+\W)\psi$ is $\W$\=/holomorphic.
Consequently we may take locally infinitely many such derivatives of $\phi$.
This is an quaternionic function theoretical version of the classical statement that holomorphic functions have infinitely many derivatives.
We shall use finitely many of such global iterations in Chapter~\ref{chapter:pluecker} in order to prove the Plücker formula.

\chapter{Kodaira Representation}
\label{chapter:kodaira}

This chapter marks the pivot from mostly local questions to mostly global ones.
We apply the theory developed in the preceding chapters to represent an admissible map $F: \X \to \qat$ as a quotient of two holomorphic sections of a holomorphic $\qat$\=/line bundle in Theorem~\ref{thm:kodaira}.
We call this the Kodaira representation of $F$.
More precisely in Theorem~\ref{thm:kodaira normal} we prove a one-to-one correspondence between left normals~\eqref{eq:left normal} with $N\in\sobolev{1,2}(\Omega,\qat)$ and triples $(E,\V,\upsilon)$ of a holomorphic $\mathbb{C}$\=/line bundle $E$ on $\X$, a potential $\V \in \pot{E}^-$, and a $\V$\=/holomorphic section $\upsilon$ without roots.
Then $F = \upsilon^{-1}\phi$ for some $\V$\=/holomorphic section $\phi$.
Remarkably, this correspondence allows us to parlay Wente's inequality for $N$ in Lemma~\ref{wente} to an improvement of the regularity of the results of Chapters~\ref{chapter:local} and~\ref{chapter:darboux}.
In particular $\V$\=/holomorphic sections belong to $\banach{\infty}\loc(\Omega,\qat) \cap \sobolev{1,2}\loc(\Omega,\qat)$.

As a first step, we examine the objects of Theorem~\ref{thm:kodaira normal} locally in the smooth case.
This explanation, like Example~\ref{eg:catenoid holomorphic trivialization}, aims to familiarize the reader with how to write down a holomorphic $\qat$\=/line bundle whose complex structure is simple (left multiplication by $\qi$) but whose connection may be complicated.
For any left normal $N$ on $\Omega \subset \mathbb{C}$ we show that there exists a corresponding triple $(E,\V,\upsilon)$.
A left normal $N$ endows the trivial bundle $\X \times \qat$ with a complex structure $J_N$.
And the trivial bundle of course has the trivial connection $\nabla$.
The underlying $\mathbb{C}$\=/bundle $E$ is endowed by $(\nabla'')^+$ with a holomorphic structure $\delbar{E}$.
The difference $\nabla'' - \delbar{E} = (\nabla'')^-$ is described by an anti-commuting potential $\V$.
Let $\mathbf{e}$ be a non-vanishing holomorphic section of $E$ on a possibly smaller set $\Omega' \subset \Omega$.
That $\mathbf{e}$ is a section of $E$ means that $J_N\mathbf{e} = \mathbf{e}\qi$ and that it is holomorphic means $\delbar{E} \mathbf{e} = 0$.
For other sections of $E_\qat$
\[
J_N(\mathbf{e}\alpha) 
= \mathbf{e}\qi\alpha, \qquad
\nabla''(\mathbf{e}\alpha)
= \mathbf{e}(d'' - d\Bar{z}\V) \alpha.
\]
But we do not yet have a correspondence because we have not defined $\upsilon$.
The missing information is essentially the section $\mathbf{1}$.
It gives not only the trivialization $\Omega \times \qat$ but also determines $\nabla$ through $\nabla(\mathbf{1}\alpha) 
= \mathbf{1}d\alpha$.
In terms of the trivializing section $\mathbf{e}$ the section $\mathbf{1}$ should be described as $\mathbf{1} = \mathbf{e}\upsilon$ for a necessarily non-vanishing $\upsilon : \Omega \to \qat$. 
That $\mathbf{e}$ belongs to $E$ implies
\begin{equation*}
\mathbf{e}\qi \upsilon = J_N (\mathbf{e}\upsilon)
= J_N\mathbf{1} = \mathbf{1}N = \mathbf{e} \upsilon N
\qquad\Rightarrow\qquad 
\qi \upsilon = \upsilon N.
\end{equation*}
Additionally, rearranging yields $N = \upsilon^{-1}\qi \upsilon$, showing how to recover $N$.
Likewise, the fact that $\mathbf{e}$ is holomorphic is equivalent to
\begin{align*}
\mathbf{e}d''\upsilon
= (\nabla'')^+(\mathbf{e}\upsilon)
&= (\nabla'')^+\mathbf{1}
= \tfrac{1}{2}(\nabla''\mathbf{1} - J_N \nabla''(J_N\mathbf{1}))
= 0 - \tfrac{1}{2}J_N \nabla''(\mathbf{1}N) \\
&= - \tfrac{1}{4} J_N (\mathbf{1}dN - J_N\mathbf{1} \ast dN) 
= - \mathbf{1}\tfrac{1}{4}(N dN + \ast dN) \\
&= - \mathbf{e}\upsilon\tfrac{1}{4}(N dN + \ast dN).
\end{align*}
This equation can be written in a simpler form, since
\begin{equation}
\label{eq:N potential is 01}
\begin{aligned}
\upsilon\tfrac{1}{4}(N dN + \ast dN)
= \upsilon\tfrac{1}{4}\left(N dx \frac{\partial N}{\partial x} + N dy \frac{\partial N}{\partial y} + dy \frac{\partial N}{\partial x} - dx \frac{\partial N}{\partial y}\right)\\
= \upsilon\tfrac{1}{4}(dx - N dy)\left(N \frac{\partial N}{\partial x} - \frac{\partial N}{\partial y}\right)
= d\bar{z}\,\upsilon\tfrac{1}{4}\left(N \frac{\partial N}{\partial x} - \frac{\partial N}{\partial y}\right),
\end{aligned}
\end{equation}
using $\upsilon (dx - N dy) = (dx - \qi dy) \upsilon$.
In Equation~\eqref{eq:V- from N} we computed $(\nabla'')^-(\mathbf{1}\xi_1) = \mathbf{1}\tfrac{1}{4}(NdN+\ast dN)\xi_1$.
Applying this formula with $\xi_1 = \upsilon^{-1}$ gives a formula for $\V$:
\begin{align*}
-\mathbf{e}\,d\bar{z}\,\V
= (\nabla'')^-\mathbf{e}
= \mathbf{e}\,\upsilon\tfrac{1}{4}(NdN+\ast dN)\upsilon^{-1} 
= \mathbf{e} \,d\bar{z}\,\upsilon\tfrac{1}{4}\left(N \frac{\partial N}{\partial x} - \frac{\partial N}{\partial y}\right)\upsilon^{-1}.
\end{align*}

In the case that $N\in\sobolev{1,2}(\Omega,\qat)$ is not smooth, this geometric picture breaks down because there is no clear way to define the holomorphic $\mathbb{C}$\=/bundle $E$; since $N$ is only measurable we do not have a complex structure $J_N \in \End_\qat(\qat)$, only a (linear) complex structure on almost all fibers.
However the correspondence is maintained by allowing a non-smooth pair $(\V,\upsilon)$ and by defining $E$ differently, specifically allowing potentials $\V$ in $\banach{2}$ and $\V$\=/holomorphic sections $\upsilon$ in $\sobolev{1,p}$ for some $1<p<2$.
In the following lemma we address the local existence question: given the left normal $N$ do there exist the desired functions $\upsilon$ and $\V$?

\begin{lemma}\label{lem:local holomorphic structure}
Let $N\in\sobolev{1,2}\loc(\Omega,\qat)$ obey $N^2=-1$ on an open set $\Omega\subset\mathbb{C}$.
Then there exists on an open neighborhood $\Omega'$ of any point in $\Omega$ a quaternionic function $\upsilon \in \sobolev{1,p}(\Omega',\qat)$ for some $1<p<2$ with the following properties:
\begin{align}\label{eq:intertwiner 2}
\qi\upsilon&=\upsilon N,
& 
d''\upsilon + \upsilon\, \tfrac{1}{4}(NdN+\ast dN) & = 0
\end{align}
Moreover there is a solution such that $\upsilon$ is $\V$\=/holomorphic for $\V \in \banach{2}(\Omega',\qat^-)$ and has no root on $\Omega'$ in the sense of Definition~\ref{order of roots}.
\end{lemma}
\begin{proof}
The following projection $P$ projects quaternionic-valued functions onto functions with the symmetry $\qi \upsilon = \upsilon N$:
\[
P: \xi \mapsto P\xi = \tfrac{1}{2}(\xi-\qi \xi N).
\]
Let us examine how derivatives behave under symmetrization.
For $\xi \in \sobolev{1,p}\loc(\Omega, \qat)$ we compute
\begin{align*}
 d''(P\xi) - P(d''\xi) &=-\tfrac{1}{4}(\qi \xi dN+\xi\ast dN)\\
 & = -\tfrac{1}{4}(P\xi)(NdN+\ast dN)+ P\left(\xi\tfrac{1}{4}(NdN -\ast dN)\right), 
\end{align*}
which is equivalent to 
\begin{align}
 \label{eq:doubleprimeprojection}
d''(P\xi) + (P\xi)\tfrac{1}{4}(NdN+\ast dN) 
= P\left(d''\xi + \xi \tfrac{1}{4}(NdN-\ast dN)\right).
\end{align}
Any function $\xi$ for which the right-hand side of~\eqref{eq:doubleprimeprojection} vanishes gives a solution $\upsilon=P\xi$ of~\eqref{eq:intertwiner 2}.
As in Equation~\eqref{eq:N potential is 01}, the left hand side is a $(0,1)$\=/form.
Let's expand out the right hand side:
\begin{align*}
P&\left(d''\xi + \xi \tfrac{1}{4}(NdN-\ast dN)\right) \\
&= P\left(d\bar{z}\barpartial\xi + dx\,\xi \tfrac{1}{4} \left(N \tfrac{\partial N}{\partial x} + \tfrac{\partial N}{\partial y}\right)\right) + dy\,P\left(\xi \tfrac{1}{4} \left( N \tfrac{\partial N}{\partial y} - \tfrac{\partial N}{\partial x} \right) \right) \\
&= P\left(d\bar{z}\barpartial\xi + dx\,\xi \tfrac{1}{4} \left(N \tfrac{\partial N}{\partial x} + \tfrac{\partial N}{\partial y}\right)\right) + \qi dy\,P\left( \qi\xi \tfrac{1}{4} \left( N \tfrac{\partial N}{\partial x} + \tfrac{\partial N}{\partial y} \right)N \right) \\
&= d\bar{z}\, P \left( \barpartial\xi + \xi \tfrac{1}{4} \left( N\tfrac{\partial N}{\partial x} + \tfrac{\partial N}{\partial y} \right) \right),
\end{align*}
where we have used that $N$ anti-commutes with $\tfrac{\partial N}{\partial x}$ and $\tfrac{\partial N}{\partial y}$ (since $N^2 = -1$ implies $N dN = -dN N$), that $P(\qi\xi) = \qi P\xi$, and that $P(\qi \xi N) = - P\xi$.
In particular, it vanishes if
\begin{align}
\barpartial \xi+\frac{1}{4}\xi\left(N\frac{\partial N}{\partial x} + \frac{\partial N}{\partial y}\right) = 0.
\label{eq:intertwiner 4}
\end{align}

To show the existence of non-trivial solutions $\xi$ of~\eqref{eq:intertwiner 4}, we transfer the arguments at the beginning of Chapter~\ref{chapter:local} to potentials that are multiplied from the right.
Let us use the notation $\frac{1}{4}\left(N\frac{\partial N}{\partial x}+\frac{\partial N}{\partial y}\right)_{\mathrm{r}}$ for right multiplication with this potential.
Since $N$ has absolute value $|N|=1$, it in fact belongs to
$\sobolev{1,2}\loc(\Omega,\qat) \cap \banach{\infty}(\Omega,\qat)$.
Hence the potential $\frac{1}{4}\left(N\frac{\partial N}{\partial x}+\frac{\partial N}{\partial y}\right)$ belongs to $\banach{2}\loc(\Omega,\qat)$.
Then for sufficiently small $\Omega'\subset\Omega$ the Neumann series,
\[
\left(\unity+\Op{I}_{\Omega'}\tfrac{1}{4}\left(N\tfrac{\partial N}{\partial x}+\tfrac{\partial N}{\partial y}\right)_{\mathrm{r}}\right)^{-1}=\sum_{l=0}^{\infty}\left(-\Op{I}_{\Omega'}\tfrac{1}{4}\left(N\tfrac{\partial N}{\partial x}+\tfrac{\partial N}{\partial y}\right)_{\mathrm{r}}\right)^l,
\]
converges as an operator on $\banach{q}(\Omega',\qat)$, which maps the elements in the kernel of $\barpartial$ onto solutions of~\eqref{eq:intertwiner 4}.
Hence there exist locally non-trivial solutions $\xi$ of~\eqref{eq:intertwiner 4}.

Next we show that the solution $\xi$ of~\eqref{eq:intertwiner 4} can be chosen such that $\upsilon=P\xi$ is a non-trivial solution of~\eqref{eq:intertwiner 2}.
Let us assume to the contrary that all solutions of~\eqref{eq:intertwiner 4} belong to $\ker(P)$.
For any $z_0$ in the Lebesgue set of $N$~\cite[Chapter~I \S1.8]{St} we have
\begin{align*}
\lim_{R\downarrow 0}\tfrac{1}{R^2\pi}\big\|N-N(z_0)\big\|_{\banach{1}(B(z_0,R))}&=0,&
\lim_{R\downarrow 0}\tfrac{1}{4}\big\|N\tfrac{\partial N}{\partial x}-\tfrac{\partial N}{\partial y}\big\|_{\banach{2}(B(z_0,R))}&=0.
\end{align*}
Because of $N(z_0)^2=-1$, the solutions of $\alpha =-\qi \alpha N(z_0)$ form a linear subspace of $r\in\qat$ of real dimension $2$.
We choose one with $|\alpha|=1$ and define
\begin{align*}
\xi&=\left(\unity+\Op{I}_{B(z_0,R)}\tfrac{1}{4}\left(N\tfrac{\partial N}{\partial x}-\tfrac{\partial N}{\partial y}\right)_{\mathrm{r}}\right)^{-1}\alpha.
\end{align*}
For some $2<q<\infty$ and sufficiently small $R>0$ we then obtain
\begin{align*}
\|N-N(z_0)\|_{\banach{2}(B(z_0,R))}
&\leq \tfrac{1}{2}R^2\pi, \\
\|\xi - \alpha\|_{\banach{q}(B(z_0,R))}
&\leq\tfrac{1}{2}\|\alpha\|_{\banach{q}(B(z_0,R))} =\tfrac{1}{2}(R^2\pi)^{\frac{1}{q}}.
\end{align*}
We estimate with the Hölder inequality for $\frac{1}{q'}=1-\frac{1}{q}$
\[
\|\xi-\alpha\|_{\banach{1}(B(z_0,R))}\le\|1\|_{\banach{q'}(B(z_0,R))}\|\xi-\alpha\|_{\banach{q}(B(z_0,R))}\le\tfrac{1}{2}R^2\pi.
\]
We assumed $P\xi=0$, i.e.\ $\xi=\qi\xi N$, $\alpha=-\qi \alpha N(z_0)$ and $|\alpha|=1$.
This gives
\begin{align*}
2R^2\pi
&= 2\|\alpha\|_{\banach{1}(B(z_0,R))}
=\left\|(\xi-\qi\xi N)-(\alpha-\qi \alpha N(z_0))\right\|_{\banach{1}(B(z_0,R))} \\
&\leq \|\xi-\alpha\|_{\banach{1}(B(z_0,R))}+\|\qi(\xi-\alpha)N\|_{\banach{1}(B(z_0,R))}+\|\qi \alpha(N-N(z_0))\|_{\banach{1}(B(z_0,R))} \\ 
&= 2\|\xi-\alpha\|_{\banach{1}(B(z_0,R))}+\|N-N(z_0)\|_{\banach{1}(B(z_0,R))}
\leq\left(\tfrac{2}{2}+\tfrac{1}{2}\right)R^2\pi=\tfrac{3}{2}R^2\pi.
\end{align*}
This contradicts the assumption that all solutions of~\eqref{eq:intertwiner 4} are contained in the kernel of $P$ and shows that there exist non-trivial solutions of~\eqref{eq:intertwiner 2}.

Fix now a non-trivial solution $\upsilon$ of~\eqref{eq:intertwiner 2} and define
\[
- \V 
= \begin{cases}
\upsilon\tfrac{1}{4}\left(N \frac{\partial N}{\partial x} - \frac{\partial N}{\partial y}\right)\upsilon^{-1}
&\text{for } \upsilon \neq 0, \\
0 &\text{for } \upsilon = 0.
\end{cases}
\]
The second case of the definition is necessary here because we do not yet know that $\upsilon$ is non-zero.
Observe that the norm of $\V$ at any point is either $\frac{1}{4}\left|N\frac{\partial N}{\partial x} - \frac{\partial N}{\partial y}\right|$ or $0$, so $\V \in \banach{2}(\Omega',\qat)$.
Where $\upsilon$ is non-zero
\[
-\qi\V
= \upsilon N \tfrac{1}{4}\left(N \frac{\partial N}{\partial x} - \frac{\partial N}{\partial y}\right)\upsilon^{-1}
= - \upsilon \tfrac{1}{4}\left(N \frac{\partial N}{\partial x} - \frac{\partial N}{\partial y}\right)N\upsilon^{-1}
= \V \qi,
\]
which shows that $\V$ is valued in $\qat^-$.
Whether or not $\upsilon$ is zero, Equation~\eqref{eq:N potential is 01} shows that 
\[
\upsilon\, \tfrac{1}{4}(NdN+\ast dN)
= - d\bar{z} \V \upsilon,
\]
and so $(\barpartial - \V)\upsilon = 0$.
Therefore $\upsilon$ is $\V$\=/holomorphic.

Finally, to produce a solution that has no roots, we may have to modify $\upsilon$ and $\V$.
According to Lemma~\ref{quotient dimension} any roots of $\upsilon$ are discrete and of finite order.
Let $f$ be a holomorphic function with poles at the roots of $\upsilon$ of equal order; such a function exists due to the solubility of Mittag-Leffler distribution on non-compact Riemann surfaces (see~\cite[26.3~Theorem]{Fo}).
The product $f\upsilon$ still solves~\eqref{eq:intertwiner 2}, is $f\V f^{-1}$\=/holomorphic, and has no roots on $\Omega'$.
\end{proof}
The function $\upsilon$ acts as an `intertwiner', allowing us to move back and forth across the correspondence between sections of $\X\times \qat$ and $E_\qat$, even in the non-smooth case.
For example, consider the section $\mathbf{1}\alpha = \mathbf{e}\upsilon\alpha$.
If $N$ is smooth and defines a complex structure $J_N$ we have $\nabla''(\mathbf{1}\alpha) = \mathbf{1} \tfrac{1}{2}(d\alpha - N \ast d\alpha)$ from Equation~\eqref{eq:holo product rule}.
But even in the non-smooth case the analogous formula follows from Equation~\eqref{eq:intertwiner 2}:
\begin{align*}
(\delbar{E} - d\Bar{z}\V)(\mathbf{e}\upsilon\alpha)
&= \delbar{E} (\mathbf{e}\upsilon)\alpha + \tfrac{1}{2} (\mathbf{e}\upsilon d\alpha - J \mathbf{e}\upsilon \ast d\alpha) 
- d\Bar{z}\V \upsilon\alpha \\
&= (\mathbf{e}d''\upsilon - d\Bar{z}\V \upsilon)\alpha + \mathbf{e} \tfrac{1}{2} (\upsilon d\alpha - \qi\upsilon \ast d\alpha) \\
\labelthis{eq:intertwiner 1}
&= \mathbf{e}\upsilon \tfrac{1}{2} (d\alpha - N \ast d\alpha).
\end{align*}

As the last part of the proof of the lemma shows, the $\upsilon$ is not unique.
Importantly any two local solutions differ by left multiplication with a complex function that is moreover holomorphic.
Therefore these different local solutions fit together to give the cocycle of a holomorphic line bundle $E$.
Hence the line bundle $E$ in the correspondence is indirectly defined.
This is the key idea that allows us in general to move from results about local holomorphic functions to global sections.
In particular, until now we have left off the subscripts of local objects, but from this point on we are more fastidious about the distinction.

\begin{theorem}
\label{thm:kodaira normal}
\index{Normal}
On a Riemann surface $\X$ to every left normal~\eqref{eq:left normal} $N\in\sobolev{1,2}\loc(\X,\qat)$ we can associate a triple $(E,\V,\upsilon)$, consisting of a holomorphic $\mathbb{C}$\=/line bundle $E$, a $\banach{2}\loc$\=/potential $\V\in\pot{E}^-$ and a $\V$\=/holomorphic section $\upsilon \in H^0(\X,\Q{E,\V})$ without roots, such that $N=\upsilon^{-1}\qi\upsilon$.

The line bundle $E$ is uniquely determined up to isomorphism class.
Fixing $E$ the associated triples have the form $(E,g \V g^{-1},g\upsilon)$ for any holomorphic function $g : \X \to \mathbb{C}^\ast$.
\end{theorem}
\begin{proof}
Let $N\in\sobolev{1,2}\loc(\X,\qat)$.
Due to Lemma~\ref{lem:local holomorphic structure} there exists on each member of an open cover of $\X$ a square-integrable potential $\V$ and a solution $\upsilon$ of~\eqref{eq:intertwiner 2} without roots.
Let $(\V_l,\upsilon_l)$ and $(\V_m,\upsilon_m)$ be two such pairs on the open sets $\SO_l$ and $\SO_m$, respectively.
It is always possible to shrink the members of the cover $\{\SO_l\}$ such that they comport with Remark~\ref{rem:special cover}.
By Theorem~\ref{thm:darboux}(d) we know that $\chi_l = \overline{(\qj \upsilon_l)}^{-1}$ is $\sobolev{1,p}$ for all $1<p<2$, so therefore so too is $\upsilon_l^{-1}$ since they have the same length.
Using the product rule, $f_{ml} = \upsilon_m\upsilon_l^{-1}$ and $f_{ml}^{-1}=\upsilon_l\upsilon_m^{-1}$ belong to $\bigcap_{1<p<2}\sobolev{1,p}\loc(\SO_l\cap\SO_m,\qat)$.
The first equation of~\eqref{eq:intertwiner 2} implies that $f_{ml}$ is $\mathbb{C}$\=/valued, while the second equation applied to $\upsilon_m = f_{ml}\upsilon_l$ yields
\begin{align*}
0
= d\Bar{z}_l\barpartial_l(f_{ml}\upsilon_l) + (f_{ml}\upsilon_l)\, \tfrac{1}{4}(NdN+\ast dN)
= d\Bar{z}_l(\barpartial_l f_{ml})\upsilon_l + 0,
\end{align*}
and similarly $0 = d\Bar{z}_l(\barpartial_l f_{ml}^{-1}) \upsilon_m$.
Hence $f_{ml}$ and $f_{ml}^{-1}$ are holomorphic.
By this process we obtain a holomorphic cocycle $f_{ml}$ of a holomorphic $\mathbb{C}$\=/line bundle $E$ on $\X$.
Inspection of the formula defining $\V$ in the proof of Lemma~\ref{lem:local holomorphic structure} shows that $\V_m = f_{ml}\V_lf_{ml}^{-1}$.
The local pairs $(\V_l,\upsilon_l)$ define a global potential $\V\in\pot{E}^-$ and a $\V$\=/holomorphic section $\upsilon$ without roots of the holomorphic $\qat$\=/line bundle $E_\qat$.
The first equation of~\eqref{eq:intertwiner 2} implies $N=\upsilon^{-1}\qi\upsilon$.
		
Vice versa, let $(E,\V,\upsilon)$ be a triple as described in the theorem.
Due to Theorem~\ref{thm:darboux} $\upsilon^{-1}$ is measurable and $N=\upsilon^{-1}\qi\upsilon$ belongs to $\bigcap_{1<p<2}\sobolev{1,p}\loc(\X,\qat)$ with 
$dN=\upsilon^{-1} (\qi (d\upsilon)\upsilon^{-1} - (d\upsilon)\upsilon^{-1} \qi) \upsilon = \upsilon^{-1}[\qi,(d\upsilon)\upsilon^{-1}]\upsilon$.
Now Theorem~\ref{thm:darboux} implies $dN = 2\upsilon^{-1} \qi (d\Bar{z} \V- dz \U)\upsilon$ for some $\U\in\banach{2}\loc(\X,\qat^-)$, since the $\qat^+$ part of $(d\upsilon)\upsilon^{-1}$ drops out.
Therefore $N \in\sobolev{1,2}\loc(\X,\qat)$.

Suppose $(E,\V,\upsilon)$ and $(\Tilde{E},\Tilde{\V},\Tilde{\upsilon})$ are both associated to $N$.
Consider the quotient $g_l = \Tilde{\upsilon}_l\upsilon_l^{-1}$ locally.
Again, $g_l$ must be non-vanishing and holomorphic.
Observe
\[
g_m 
= \Tilde{\upsilon}_m\upsilon_m^{-1}
= \Tilde{f}_{ml}\Tilde{\upsilon}_l\upsilon_l^{-1}f_{ml}^{-1}
= \Tilde{f}_{ml} g_l f_{ml}^{-1},
\]
which shows that $g$ is a section of $\Tilde{E}\otimes E^{-1}$.
Hence $\Tilde{E}$ and $E$ are isomorphic holomorphic $\mathbb{C}$\=/line bundles.
Fix $E$.
Now $g_m = g_l$, so we have a global holomorphic function $g$.
Again by the above reasoning $\Tilde{\V}_l = g_l \V_l g_l^{-1}$.
Conversely, if $g : \X \to \mathbb{C}^\ast$ is non-vanishing and holomorphic, then $g\upsilon$ is $g\V g^{-1}$\=/holomorphic and has no roots.
\end{proof}

Before we prove the main theorem of this chapter, Theorem~\ref{thm:kodaira}, which extends the above theorem to admissible maps, we show that this correspondence between left normals and holomorphic bundles can be used to strengthen the results of Chapters~\ref{chapter:local} and~\ref{chapter:darboux}.
In~\cite{We} the regularity of constant mean curvature surfaces, conformal maps $F : \Omega \to \mathbb{R}^3$ satisfying the non-linear PDE $\triangle F = H \frac{\partial F}{\partial x} \times \frac{\partial F}{\partial y}$, was found to be better than what might be expected from classical theory.
The technique was extracted and generalized, and regularity results of the PDE $-\triangle u = \det dA$ for $A : \Omega \to \mathbb{R}^2$ that use this idea have come to be known as Wente's inequality. 

\begin{lemma}[Wente's inequality]
\label{wente}
\index{Wente's inequality}
On an open $\Omega\subset\mathbb{C}$ any $\V$\=/holomorphic function $\upsilon$ with $V\in\banach{2}\loc(\Omega,\qat^-)$ belongs to $\upsilon\in\banach{\infty}\loc(\Omega,\qat)\cap\sobolev{1,2}\loc(\Omega,\qat)$.
Furthermore, $|\upsilon|$ is continuous with the same roots as $\upsilon$ in Definition~\ref{order of roots}.
\end{lemma}
\begin{proof}
Let $\upsilon$ be $\V$\=/holomorphic with $\V\in\banach{2}\loc(\Omega,\qat^-)$.
Due to Lemma~\ref{quotient dimension} the roots of $\upsilon$ are isolated and we may assume locally that $\upsilon$ has no roots.
Due to Theorem~\ref{thm:kodaira normal} $N=\upsilon^{-1}\qi\upsilon$ belongs to $\sobolev{1,2}\loc(\Omega,\qat)$.

$N$ is a map from $\Omega$ into the $2$\=/sphere $S^2$ of purely imaginary unit quaternions.
Thus we can regard $dN = 2\upsilon^{-1} \qi (d\Bar{z} \V- dz \U)\upsilon$ as a $1$\=/form with values in $TS^2$.
The determinant $\det(dN)$ is equal to the determinant of the  $\qat^-$\=/valued $1$\=/form $2\qi (d\Bar{z}\V-dz\U)$ because determinants are invariant under changes of coordinates.
Letting $\U = \U_2 \qj + \U_3 \qk$ and likewise for $\V$, we can write this $1$\=/form, a linear map from $T\Omega$ to $\qat^-$, as a matrix:
\begin{align*}
2\qi (d\Bar{z}\V-dz\U)
&= 2 dx ( \V_2\qk - \V_3\qj - \U_2\qk + \U_3\qj )
+ 2 dy ( \V_2\qj + \V_3\qk + \U_2\qj + \U_3\qk ) \\
&= 2\begin{pmatrix}
-\V_3 + \U_3 & \V_2 + \U_2 \\
\V_2 - \U_2 & \V_3 + \U_3
\end{pmatrix}
\end{align*}
A direct calculation shows $\det(dN)=4(|\U|^2-|\V|^2)$.
By Theorem~\ref{thm:darboux}~(1)(a) and (1)(c), 
\[-\triangle \ln |\upsilon| = -\tfrac{1}{4} \barpartial \partial \ln |\upsilon| = \tfrac18 \barpartial\B=-\tfrac18 (|\U|^2-|\V|^2)= -\tfrac{1}{32} \det(dN).\]
Because $N\in\sobolev{1,2}\loc(\Omega,\qat)$, Wente's inequality~\cite[Theorem~3.1.2]{Hel} implies that on any open ball $\Set{B}$ with compact closure in $\Omega$ the unique solution $f$ of the Dirichlet problem $-\triangle f = -\tfrac{1}{32}\det(dN)$ with the boundary condition $f|_{\boundary \Set{B}}=0$ belongs to $f \in C(\Bar{\Set{B}})\cap\sobolev{1,2}(\Set{B})$.
Clearly $\ln|\upsilon|$ differs from this solution by a harmonic distribution.
Due to Weyl's Lemma 
any harmonic distribution is a harmonic function.
Therefore $\ln|\upsilon|$ is continuous and hence locally bounded.
This implies $\upsilon\in\banach{\infty}\loc(\Omega,\qat)$ and with $\barpartial \upsilon=\V\upsilon\in\banach{2}\loc(\Omega,\qat)$ and the arguments of Chapter~\ref{chapter:local} also $\upsilon\in\sobolev{1,2}\loc(\Omega,\qat)$.
This proves also that the continuous function $|\upsilon|$ has no zeros if $\upsilon$ has no roots in the sense of Definition~\ref{order of roots}.
Vice versa, if $\upsilon$ has a root then the same arguments as above apply to the holomorphic section $(z-z_0)^{-\ord_{z_0}(\upsilon)}\upsilon$ without roots.
Therefore in this case the continuous function $|\upsilon|=|z-z_0|^{\ord_{z_0}(\upsilon)}|(z-z_0)^{-\ord_{z_0}(\upsilon)}\upsilon|$ has a zero at $z_0$.  
\end{proof}
This lemma improves several statements.
First we strengthen Theorem~\ref{cauchy formula} and Corollary~\ref{weyls lemma}:
\begin{corollary}
\label{cor:weyls lemma 2}
\index{Cauchy's integral formula}
\index{Weyl's lemma}
For open $\Omega\subset\mathbb{C}$ and $\U\in\banach{2}\loc(\Omega,\qat^-)$ the following conditions are equivalent:
\begin{enumeratethm}
\item \textup{(Cauchy's Integral Formula):} $\xi\in\banach{p}\loc(\Omega,\qat)$ for some $1<p<2$ and~\eqref{eq:cif} holds for all open $\Omega'\subset\Omega$ obeying~\eqref{eq:neumann-convergence} for $\U$ instead of $\V$.
\item $\xi\in\sobolev{1,2}\loc(\Omega,\qat)\cap\banach{\infty}\loc(\Omega,\qat)$ and $\barpartial \xi=\U\xi$.
\item $\xi\in\sobolev{1,1}\loc(\Omega,\qat) \subset \banach{2}\loc(\Omega,\qat)$ and $\barpartial \xi=\U\xi$.
\item $\xi\in\banach{2}\loc(\Omega,\qat)$ and the $\barpartial$\=/derivative of the distribution induced by $\xi$ is equal to the distribution induced by $\U\xi$.
\item \textup{(Weyl's Lemma of quaternionic function theory II):} $\xi\in\banach{1}\loc(\Omega,\qat)$ and $\Bar{\chi}\qj dz\xi$ is a closed current for some $\qk\Bar{\U}\qk$\=/holomorphic $\chi$ without roots on $\Omega$.
\end{enumeratethm}
\end{corollary}
\begin{proof}
(i)$\Rightarrow$(ii) follows from Theorem~\ref{cauchy formula} and Lemma~\ref{wente}.

\noindent(ii)$\Rightarrow$(iii)$\Rightarrow$(iv) are obvious.

\noindent(iv)$\Rightarrow$(v): On $\Omega$, any $\qk\Bar{\U}\qk$\=/holomorphic  $\chi$ belongs to $\sobolev{1,2}\loc(\Omega,\qat)\cap\banach{\infty}\loc(\Omega,\qat)$ by Lemma~\ref{wente}.
For $\xi$ as in (iv) and $f\in C_0^\infty(\Omega,\mathbb{R})$ the following calculation is valid since all integrands are products of elements in $C_0^\infty(\Omega,\mathbb{R})$ and $\banach{1}\loc(\Omega,\qat)$:
\begin{gather*}
\int_\Omega df\wedge\Bar{\chi}\qj dz\xi
=\int_\Omega d\left(f\Bar{\chi}\qj dz\xi\right) - \int_\Omega fd\left(\Bar{\chi}\qj dz\xi\right)
= - \int_\Omega f\overline{d\left(dz\chi\right)}\qj\xi- \int_\Omega f\Bar{\chi}\qj d(dz\xi)\\
=\int_\Omega f\overline{dz\wedge d\Bar{z}(\barpartial \chi)}\qj\xi+ \int_\Omega f\Bar{\chi}\qj dz\wedge d\Bar{z}\barpartial \xi
= \int_\Omega f\Bar{\chi}\left(\qk\U\qk d\Bar{z}\wedge dz\qj + \qj dz\wedge d\Bar{z}\U\right)\xi
=0.
\end{gather*}
(v)$\Rightarrow$(i): Due to Wente's inequality, $\chi$ belongs to $\banach{\infty}\loc(\Omega,\qat)$ and consequently the closed current $\Bar{\chi}\qj dz\xi$ belongs to $\banach{1}\loc(\Omega,\qat)$.
On a simply connected bounded $\Omega$ there exists a unique 
\begin{align*}
F\in\sobolev{1,1}\loc(\Omega,\qat)
\quad\text{with}\quad
dF=\Bar{\chi}\qj dz\xi
\quad\text{and}\quad
\int_\Omega F\,\dmu=0.
\end{align*}
Due to the Theorem~\ref{thm:darboux}~(2)(b)-(2)(c) $\upsilon=-(\overline{\qj\chi})^{-1}$ is $\V$\=/holomorphic with $\V\in\banach{2}\loc(\Omega,\qat^-)$, and $\phi=\upsilon F\in\sobolev{1,1}\loc(\Omega,\qat)\subset\banach{2}\loc(\Omega,\qat)$ solves $(\barpartial -\V)\phi=0$ weakly.
On small subsets $\Omega'\subset\Omega$ the operator $\unity-\Op{I}_{\Omega}\V$ is invertible and $\Tilde{\upsilon}=(\unity-\Op{I}_{\Omega}\V)^{-1}r$ with $r\in\qat\setminus\{0\}$ is $\V$\=/holomorphic without roots.
By Lemma~\ref{wente} $\Tilde{\upsilon}$ belongs to $\banach{\infty}\loc(\Omega',\qat)$ and $\Bar{\Tilde{\upsilon}}\qj dz\phi$ is a closed current on $\Omega'$.
Corollary~\ref{weyls lemma} implies that $\phi$ is $\V$\=/holomorphic.
Due to Theorem~\ref{thm:darboux}~(1)(f) $\xi$ is equal to $(\partial+\B+\U)\phi$, and by Theorem~\ref{thm:darboux}~(1)(e) it is $\U$\=/holomorphic on $\Omega'$.
Now (i) follows from Theorem~\ref{cauchy formula}.
\end{proof}
Furthermore, as applied to the Darboux transformation, Lemma~\ref{wente} yields:
\begin{corollary}
\label{cor:darboux improved}
Let $\upsilon$ and $\phi$ be $\V$\=/holomorphic sections with the assumptions of Theorem~\ref{thm:darboux}.
The functions defined therein have the following improvements to their regularity:
\begin{enumerate}
\item[\textup{(1)(a$'$)}] $\B\in\banach{2}\loc(\Omega,\qat^+)$ instead of $B\in\banach{2,\infty}\loc(\Omega,\qat^+)$,
\item[\textup{(1)(d$'$)}] $\chi\in\banach{\infty}\loc(\Omega,\qat)\cap\sobolev{1,2}\loc(\Omega,\qat)$ instead of $\chi\in\bigcap_{1<p<2}\sobolev{1,p}\loc(\Omega,\qat)$.
\item[\textup{(1)(e$'$)}] $\psi\in\banach{\infty}\loc(\Omega,\qat)\cap\sobolev{1,2}\loc(\Omega,\qat)$ instead of $\psi\in\bigcap_{1<p<2}\sobolev{1,p}\loc(\Omega,\qat)$.
\item[\textup{(1)(f$'$)}] $F = \upsilon^{-1}\phi\in\sobolev{2,2}\loc(\Omega,\qat)\cap\sobolev{1,\infty}\loc(\Omega,\qat) \subset C^{0,1}(\Omega,\qat)$ instead of $\upsilon^{-1}\phi\in\bigcap_{1<p<2}\sobolev{2,p}\loc(\Omega,\qat)$.
\end{enumerate}
\end{corollary}

These improvements mean that in the Definitions~\ref{def:weakly conformal local} and~\ref{def:weakly conformal global} of admissible maps it suffices to assume $F\in\sobolev{1,1}\loc$.
Without applying Wente's inequality we would have had to define admissible to be a more limited class of maps, namely $F\in\sobolev{1,p}\loc$ for some $1<p<2$ (compare $F$ in the proofs of Corollaries~\ref{weyls lemma} and~\ref{cor:weyls lemma 2}) in order to have the same level of regularity of the sections.
It was for this reason that we interposed Wente's inequality and the discussion of regularity between Theorems~\ref{thm:kodaira normal} and~\ref{thm:kodaira}.

So then we come to the main theorem of this chapter.
In Theorem~\ref{thm:kodaira normal} we constructed a correspondence between left normals and triples $(E,\V,\upsilon)$.
To construct a correspondence to admissible maps $F$ we need additional data to distinguish between maps with the same normal.
Just as $\upsilon$ is the section of $E_\qat$ associated to $\mathbf{1}$ of $\X \times \qat$, the idea is to conceive of $F$ as a section of $\X \times \qat$ associated to a section $\phi$ of $E_\qat$.
Since $E_\qat$ is a line bundle and $\upsilon$ is a section without roots, it is possible to write any section as $\phi = \upsilon F$.
Compare with~\cite[Example, p.~395]{PP}. 
We exclude the case that $\phi$ is a constant multiple of $\upsilon$, since then the map $F = \upsilon^{-1}\phi$ is constant and not weakly conformal.

\begin{theorem}[Quaternionic Kodaira representation]
\label{thm:kodaira}
\index{Kodaira representation}
On a Riemann surface $\X$ to every admissible map $F:\X\to\qat$ we can associate a quadruple $(E,\V,\upsilon,\phi)$, consisting of a holomorphic $\mathbb{C}$\=/line bundle $E$, a potential $\V\in\pot{E}^-$ and two global linear independent $\V$\=/holomorphic sections $\upsilon$ and $\phi$, such that $\upsilon$ has no roots and $F=\upsilon^{-1}\phi$.
The left normal $N$ of $F$ corresponds in Theorem~\ref{thm:kodaira normal} to $(E,\V,\upsilon)$.

The line bundle $E$ is uniquely determined up to isomorphism class.
Fixing $E$, the associated quadruples have the form $(E,g \V g^{-1},g\upsilon, g\phi)$ for any holomorphic function $g : \X \to \mathbb{C}^\ast$.
\end{theorem}
\begin{proof} Due to Theorem~\ref{thm:kodaira normal} it remains to prove that for a given left normal $N\in\sobolev{1,2}\loc(\X,\qat)$ with corresponding triple $(E,\V,\upsilon)$ the admissible maps $F$ with this left normal are in one-to-one correspondence to global $\V$\=/holomorphic sections $\phi$.
Due to~\eqref{eq:intertwiner 1} the section $\phi=\upsilon F$ is $\V$\=/holomorphic if and only if
\begin{align*}
0 
= (\delbar{E} - \V)\phi
= \upsilon \tfrac{1}{2}(dF - N \ast dF)
\qquad \Leftrightarrow \qquad
dF = N \ast dF.
\end{align*}
In other words, if and only if $F$ has $N$ as its left normal, Equation~\eqref{eq:left normal}.
Therefore $F \mapsto \phi$ yields the one-to-one correspondence.
\end{proof}

\begin{example}[Plane]
\label{eg:plane kodaira}
\index{Plane!Kodaira representation}
A simple case is an admissible map to a plane, $F(z) = \alpha f(z) \beta$ for a $\mathbb{C}$\=/valued function $f$ and constants $\alpha,\beta \in \mathbb{S}^3$.
We know from Example~\ref{eg:plane admissible} that $f$ is analytic.
Assume first that $f$ is holomorphic.
Then $dF = \alpha\, dz f'(z) \beta$ implies that $N = \alpha\qi\alpha^{-1}$ and $E = \X \times \alpha \mathbb{C}$.
We read off that $\upsilon = \alpha^{-1}$ is a solution to Equation~\eqref{eq:intertwiner 2}, whose potential is $\V = 0$.
Lastly $\phi$ is also $0$\=/holomorphic for $\phi = \upsilon F = \alpha^{-1} \alpha f(z) \beta = f(z) \beta$.

If $f$ is anti-holomorphic, then write $F(z) = (-\alpha\qj)\qj f(z) \beta$.
The same line of reasoning gives a Kodaira representation $\V = 0$, $\upsilon = \qj \alpha^{-1}$, and $\phi = \qj f(z) \beta$.
\end{example}

\begin{example}[Catenoid]
\label{eg:catenoid kodaira}
\index{Catenoid!Kodaira representation}
We can also find the Kodaira representation of the catenoid.
Recall the work we have already done in Example~\ref{eg:catenoid holomorphic trivialization}.
There we found that $\mathbf{e} = \mathbf{1} (1-z^{-1}\qk)\qj$ was a holomorphic section of $E$.
By definition
\[
\mathbf{1} = \mathbf{e}\upsilon
\qquad\Rightarrow\qquad
\upsilon 
= -\qj (1-z^{-1}\qk)^{-1},
\]
which is $\V$\=/holomorphic for
\[
\V = \qk z^{-1}\Bar{z}(1+|z|^2)^{-1}.
\]
And $\mathbf{e}\phi = \mathbf{e}\upsilon F$ completes the data.
\end{example}

\begin{example}[Round Sphere]
\label{eg:sphere kodaira}
\index{Round sphere!Kodaira representation}
Let us give an example for a Riemann surface that is not a subset of $\mathbb{C}$: the sphere $\X = \mathbb{P}^1$.
We use the standard coordinates $z_1 = z$ and $z_2 = z^{-1}$.
The simplest conformal immersion is the embedding as the round sphere in $\Imag\qat$, which is essentially stereographic projection:
\[
F(z) = \frac{1}{|z|^2 + 1}\left( (|z|^2 - 1)\qi - 2\qj z\right).
\]
We see, for example, that $F(0) = -\qi$ and $F(\infty) = \qi$ are the south and north poles, where as the equator $|z| = 1$ is mapped to $F(z) = -\qj z \in \qj \mathbb{C}$.
We use the same trick as in Example~\ref{eg:catenoid} of writing the derivative in a special form
\begin{equation}
\label{eq:sphere dF}
\frac{1}{2}(|z|^2 + 1)^2\, dF
= (-\qj z - \qi) \qj dz (\qj z + \qi),
\end{equation}
which again allows us to read off
\[
N 
= (-\qj z - \qi)(-\qi)(-\qj z - \qi)^{-1}
=(-\qj - \qi z_2)(-\qi)(-\qj - \qi z_2)^{-1}.
\]
We see that $N(0) = -\qi$ and $N(\infty) = \qi$ are outward pointing normals.
The $\mathbb{C}$\=/line bundle $E$ is the subbundle 
\begin{multline*}
\{ (z, (\qj z + \qi)\qj h) \in \mathbb{P}^1\setminus\{\infty\} \times \qat \mid h \in \mathbb{C} \} \\
\cup \;\{ (z,(\qj + \qi z_2)\qj h) \in \mathbb{P}^1\setminus\{0\} \times \qat \mid h \in \mathbb{C} \},
\end{multline*}
where the fibers agree for $z \in \mathbb{C}\setminus\{0\}$.
We could proceed as in Example~\ref{eg:catenoid holomorphic trivialization}, but instead let's solve Equation~\eqref{eq:intertwiner 2}.
First we calculate the parts of potential term:
\begin{align*}
dN
&= (-\qi + N)\qj dz (-\qj z - \qi)^{-1} \\
NdN
&= -(N\qi + 1)\qj dz (-\qj z - \qi)^{-1} \\
\ast dN
&= (1 + N\qi) \qj dz (-\qj z - \qi)^{-1},
\end{align*}
so the second equation of~\eqref{eq:intertwiner 2} reduces to $\barpartial \upsilon_1 = 0$.
Hence $\V_1 \equiv 0$.
We have used the subscript $1$ here to indicate that we are working in the $z_1$\=/coordinate chart.
Functions that have the symmetry are all of the form $\upsilon_1(z) = f(z)\qj(-\qj z - \qi)^{-1}$ for $f(z) \in \mathbf{C}$. 
Then the PDE becomes
\begin{gather*}
0
= \barpartial (f \,\qj(-\qj z - \qi)^{-1}) 
= \barpartial f\, \qj(-\qj z - \qi)^{-1} + f \qj \partial (-\qj z - \qi)^{-1} \\
\Leftrightarrow\qquad
\barpartial f = f (1 + |z|^2)^{-1}z,
\end{gather*}
which has the solution $f = g(z)(1 + |z|^2)$ for $g$ holomorphic.
Taking $g \equiv -1$ gives a nowhere vanishing solution 
\[
\upsilon_1 
= -(1 + |z|^2) \qj(-\qj z - \qi)^{-1}
= -\qj(\qj z + \qi)
= z + \qk.
\]

We can repeat this procedure in the $z_2$ chart.
Again we find that $N dN + \ast dN \equiv 0$, so the PDE to solve is $\barpartial \upsilon_2 = 0$ and $\V_2 \equiv 0$.
This time the symmetry forces $\upsilon_2(z) = f(z) \qj (-\qj - \qi z_2)^{-1}$, which leads to the same equation for $f$ as above.
Thus we obtain a solution $\upsilon_2 = 1 + \qk \Bar{z}_2$.
As expected from the proof of Theorem~\ref{thm:kodaira normal}, the two solutions differ by a $\mathbb{C}$\=/cocycle $\upsilon_2 = z_2 \upsilon_1$.
This cocycle $f_{21} = z_2$ constructs the bundle $E_\qat$ in which these two local sections define a global section $\upsilon$ without roots.
As an aside, the cocycle shows that $\deg E = 1$.
The potential is $\V \equiv 0$.

To complete the Kodaira representation of the round sphere, we need the second section:
\begin{align*}
\phi_1
= \upsilon_1 F
&= \qi z - \qj,
&
\phi_2
= \upsilon_2 F
&= \qi  - \qj\Bar{z}_2.
\end{align*}
This is also $\V$\=/holomorphic, 
The Kodaira representation of $F$ is $F = \upsilon_1^{-1}\phi_1 = (z + \qk)^{-1}\qi (z + \qk)$.
This pleasing form was not apparent from the initial definition of $F$.
\end{example}


In later chapters we will consider variations of admissible maps through variations of their Kodaira representations.
In order to include variations of the holomorphic structure of the underlying $\mathbb{C}$\=/line bundle, we shall extend the statement in Theorem~\ref{thm:kodaira}.
Specifically, we will show that $F=\upsilon^{-1}\phi$ is admissible even in the case where the sections $\upsilon,\phi$ are $\V$\=/holomorphic with respect to a potential not necessarily in $\pot{E}^-$.
This is done in the remainder of this chapter by gauging the part $A\in\pot{E}^+$ of $\V$ away.
The following lemma sets up the proof in Corollary~\ref{gauge holomorphic structure}.
\begin{lemma}\label{bar solution}
On a bounded open subset $\Omega\subset\mathbb{C}$ all $g\in\banach{2,\infty}(\Omega)$ have a solution $h$ of $\barpartial h=g$, whose exponential $\exp(\pm h)$ belongs to $\banach{q}\loc(\Omega)$ if $1\leq q<\frac{2\sqrt{\pi}}{\|g\|_{2,\infty}}$.
\end{lemma}
\begin{proof} 
The proof is similar to the proof of Lemma~\ref{zygmund estimate}.
Due to Dolbeault's lemma 
the convolution with $f(z)=\frac{1}{\pi z}$ is a right inverse of $\barpartial$, so $h=f\ast g$ is a solution of $\barpartial h = g$.
The distribution function associated to $f$ and its non-increasing rearrangement are equal to (\cite[Chapter~II \S3.
Chapter~V\S3]{SW},~\cite[Chapter~2 Section~1]{BS} and~\cite[Chapter~1 Section~8]{Zi}):
\begin{align*}
\mu_f(s)&=\frac{1}{\pi s^2}&f^{\ast}(t)&=\frac{1}{\sqrt{\pi t}}.
\end{align*}
By definition of the norm $\|\cdot\|_{2,\infty}$ the maximal function $g^{\ast\ast}$ of the non-increasing rearrangement $g^\ast$ with $g\in\banach{2,\infty}(\Omega)$ obeys
\begin{align*}
g^{\ast\ast}(t)&\leq\frac{\|g\|_{2,\infty}}{\sqrt{t}}&g^{\ast\ast}(t)=\frac{1}{t}\int_0^T g^\ast(s)\,ds.
\end{align*}
Making use of the first estimate in~\eqref{eq:convolution bound} we obtain
\begin{align*}
h^{\ast}(t)\leq h^{\ast\ast}(t)&\leq g^{\ast\ast}(t)\int_{\frac{1}{\sqrt{\pi t}}}^{\infty}\frac{ds}{\pi s^2}+\frac{1}{2\sqrt{\pi}}\int_{t}^{|\Omega|}g^{\ast\ast}(s)\frac{ds}{\sqrt{s}}\\&\leq\frac{\|g\|_{2,\infty}}{\sqrt{\pi}}\left(1+\frac{\ln|\Omega|-\ln(t)}{2}\right).
\end{align*}
This implies that $h$ belongs to $\banach{}_{\exp}(\Omega)$~\cite[Chapter~4.6]{BS}.
Due to a standard argument~\cite[Chapter~2 Exercise~3]{BS} $\exp|h|$ belongs to $\banach{q}(\Omega)$ if the integral
\[\int_0^{|\Omega|}\exp(qh^\ast(t))\,dt \leq C\int_0^{|\Omega|}t^{-\frac{q\|g\|_{2,\infty}}{2\sqrt{\pi}}}\,dt.\]
is finite.
This completes the proof.
\end{proof}
\begin{corollary}
\label{gauge holomorphic structure}
On any Riemann surface $\X$ for any pair $(E,\V)$ of holomorphic $\mathbb{C}$\=/line bundle $E$ and potentials $\V=\V^++\V^-\in\pot{E}$, there exists another pair $(E',\V')$ with $\deg(E') = \deg(E)$, $\V'\in\pot{E'}^-$ and $|\V'|^2=|\V^-|^2$ as defined in~\eqref{eq:square integrable}, such that the sheaves $\Q{E,\V}$ and $\Q{E',\V'}$ are isomorphic on $\X$.
Moreover, for given holomorphic line bundles $E$ and $E'$ with $\deg(E)=\deg(E')$ and $\V'\in\pot{E'}^-$ there exists $\V\in \pot{E}$ such that $\Q{E,\V}$ is isomorphic to $\Q{E',\V'}$.
\end{corollary}
\begin{proof}
We fix a pair $2<q'<q<\infty$.
By Theorem~\ref{cauchy formula}~(iv) the sheaf $\Q{E,\V}$ is the kernel of $\barpartial-\V^-$ as a subsheaf either of $\ban{q}{E}$ or of $\ban{q'}{E}$ defined in Definition~\ref{def:holomorphic}.
Cover $\X$ with coordinate charts $z_l:\SO_l\to\Omega_l$ such that each $\SO_l$ is relatively compact in $\X$ and $\Omega_l$ is bounded, as per Remark~\ref{rem:special cover}.
The holomorphic structure acts on $\SO_l$ as $\barpartial_l-\V^+_l-\V^-_l$ with $\V^\pm_l\circ z_l^{-1}\in\banach{2}(\Omega_l,\qat^\pm)$.
By choosing the $\SO_l$ small enough, we may achieve additionally $\|\V^+_l\circ z_l^{-1}\|_{2,\infty}\le2\sqrt{\pi}(\frac{1}{q'}-\frac{1}{q})$.
Due to Lemma~\ref{bar solution} there exists on $\SO_l$ a solution of $\barpartial_l h_l=-\V^+_l$ with $\exp(\pm h_l\circ z_l^{-1})\in\banach{qq'/(q-q')}(\Omega_l,\mathbb{C})$.
The arguments in the first two paragraphs of Chapter~\ref{chapter:local} concerning $\Op{I}_{\Omega'}$ show $h_l\circ z_l^{-1}\in\sobolev{1,2}(\Omega_l,\mathbb{C})$.

Due to Equation~\eqref{eq:potential transformation V+} we have $\V^+_m\,d\Bar{z}_m =\V^+_l \,d\Bar{z}_l$ on the intersection $\SO_l\cap\SO_m$ of two such sets.
Therefore the function $h_m-h_l$ and thus also the exponential $\exp(h_m-h_l)$ is holomorphic and is the cocycle of a holomorphic $\mathbb{C}$\=/line bundle $E''$ on $\X$.
Consider $E'=E''\otimes E$ with potential $\V'\in \pot{E'}^-$ locally defined by $\V'_l=\exp(-h_l)\V^-_l\exp(h_l)=\V^-\exp(h_l-\Bar{h}_l)$.
Because $E''$ is defined in terms of the exponential of a holomorphic cycle, we have $\deg(E'')=0$ and therefore $\deg(E') = \deg(E)$.
Since $h_l-\Bar{h}_l$ is imaginary we deduce $|\V'|^2=|\V^-|^2$ as defined in~\eqref{eq:square integrable}.
The left multiplication by $\exp(h_l)$ defines a morphism of sheaves from $\ban{q}{E}$ into $\ban{q'}{E'}$ and the left multiplication by $\exp(-h_l)$ from $\ban{q}{E'}$ into $\ban{q'}{E}$.
Due to 
\begin{align*}
\begin{aligned}(\barpartial_l-\V'_l)\exp(h_l)& =\exp(h_l)(\barpartial_l-\V^+-\V^-)=\exp(h_l)(\barpartial_l-\V) \\(\barpartial_l-\V_l)\exp(-h_l) &=(\barpartial_l-\V^+_l-\V^-_l)\exp(-h_l)=\exp(-h_l)(\barpartial_l-\V'_l),\end{aligned}
\end{align*}
these morphisms map the kernels of $\barpartial-\V$ and $\barpartial -\V'$ onto each other.
Theorem~\ref{cauchy formula}~(iv) shows again that $\Q{E',\V'}$ is the kernel of $\barpartial -\V'$ either as a subsheaf of $\ban{q}{E'}$ or $\ban{q'}{E'}$.
Hence the above morphisms map the subsheaves $\Q{E,\V}$ and $\Q{E',\V'}$ isomorphically onto each other.

To prove the second statement, we represent the quotient $E'' = E' \otimes E^{-1}$ of two holomorphic $\mathbb{C}$\=/line bundles by the holomorphic cocycle $(g_{ml})$ with respect to the open cover $\{\SO_l\}$.
Additionally we now assume that the $\SO_l$ are simply connected.
With a smooth partition of unity subordinate to $\{\SO_l\}$ we then write $g_{ml} = \exp(h_m-h_l)$ with smooth functions $h_l$ on $\SO_l$, see~\cite[Theorem~12.6]{Fo}.
The local $1$\=/forms $A_l d\Bar{z}_l$ with $\V^+_l = -\barpartial_l h_l$ fit together to a global $1$\=/form on $\X$ and therefore define a potential $\V^+\in \pot{E}^+$.
Moreover, let $\V'\in \pot{E'}^-$ be represented by the local functions $\V'_l$.
Then the local functions $\V^-_l =\V_l'\,\exp(\Bar{h}_l - h_l)$ define a smooth potential $\V^-\in \pot{E}^-$ and $\exp(-h_l)$ induces an isomorphism from $\Q{E',\V'}$ to $\Q{E,\V^++\V^-}=\Q{E,\V}$ as in the first part of the proof.
\end{proof}

\begin{remark}
\label{rem:cover depends on q}
In the sequel we shall apply this corollary in several situations in order to transform for a given holomorphic $\mathbb{C}$\=/line bundle on a Riemann surface $\X$ and $\V\in\pot{E}$, a section $\xi\in H^0(\X,\Q{E,\V})$ and one or several other sections of $H^0(\X,\ban{q}{E})$ into sections of $H^0(\X,\ban{q'}{E'})$ where the holomorphic $\mathbb{C}$\=/line bundle $E'$ is chosen in such a way that $\xi'\in H^0(\X,\Q{E',\V'})$ for some $\V'\in\pot{E'}^-$. This applications needs in general $2<q'<q$. Furthermore, the cover $\{\SO_l\}$ of $\X$, on whose members we construct in the proof the local functions $h_l$ will in general depend on the exponents $2<q'<q$. For this reason any such application contains first a choice of the exponents $q'<q$ and then a choice of a corresponding cover $\{\SO_l\}$ of $\X$. Furthermore, if we repeat such an application we need to further decrease $q''<q'<q$ and to refine the cover $\{\SO_l\}$.

These applications will be necessary in order to transform the given holomorphic sections into situations where we can apply either the global Kodaira representation, Theorem~\ref{thm:kodaira}, or the global Darboux transformation, Corollary~\ref{cor:global darboux}, or the quaternionic Weierstraß representation, Theorem~\ref{thm:weierstrass}. 
However, the statements which we prove by applying these theorems will in general not depend on the choice of the exponents $2<q'<q$, since due to Theorem~\ref{cauchy formula}(iv) the property on $\xi\in H^0(\X,\ban{q'}{E})$ to be $\V'$\=/holomorphic with respect to some potential $\V'\in\pot{E}$ does not depend on the exponents $q'$. 
In such situations we shall assume that the reader knows to first choose the exponents $2<q'<q$ and afterwards choose an appropriate cover $\{\SO_l\}$ of $\X$, as described in the proof of Corollary~\ref{gauge holomorphic structure}. So we shall comment on these choices, only if the truth of the statements depends on the correct choice of the exponents $q'<q$ and the cover $\{\SO_l\}$ of $\X$.
\end{remark}

\chapter{Weierstraß Representation}
\label{chapter:weierstrass}
In this chapter we merge the local transformations of Chapter~\ref{chapter:darboux} to a global transformation on a Riemann surface $\X$.
This gives us another way to describe admissible maps, namely the Weierstraß representation.
An important advantage of the Weierstraß representation over the Kodaira representation is that there is a simple condition to determine whether the admissible map lies in $\Imag\qat$, which will be developed in the following chapter.
We show that the inverse of the Darboux transformation has essentially the same form, and use these in combination to the Darboux transformation to a tuple of sections.

Central to the Weierstraß representation is the concept of paired holomorphic $\qat$\=/line bundles~\cite[Theorem~4.2]{PP}.
If $E$ is the holomorphic $\mathbb{C}$\=/line bundle of a Kodaira representation, then the paired $\qat$\=/line bundles of the Weierstraß representation are $E^{-1}_\qat$ and $KE_\qat$.
In general, the tensor product of the underlying $\mathbb{C}$\=/line bundles should be the holomorphic cotangent bundle $K$.
This makes it possible to pair their sections to produce a $1$\=/form.

\begin{definition}
\label{def:10 pairing}
\index{Pairing!(1,0), Paired line bundles}
Take any sections $\chi$ of $E^{-1}_\qat$ and $\psi$ of $KE_\qat$, represented by local functions $\chi_l:\SO_l\to\qat$ and $\psi_l:\SO_l\to\qat$ with respect to a cover as in Remark~\ref{rem:special cover}.
We define the pairing $\lp\chi,\psi\rp$ to be the $1$\=/form on $\X$ with local representatives on $\SO_l$
\begin{align*}
\lp\chi,\psi\rp_l = \Bar{\chi}_l\qj dz_l\psi_l.
\end{align*}
\end{definition}
The holomorphic cotangent bundle $K$ corresponds to the cocycle $\tfrac{dz_l}{dz_m}$ with respect to this cover.
Consequently the local representatives of the sections obey (compare~\eqref{eq:section transformation})
\begin{align*}
\chi_m&=f_{ml}^{-1}\chi_l,
&\psi_m&=f_{ml}\tfrac{dz_l}{dz_m}\psi_l
\end{align*}
and the local $1$\=/forms $\Bar{\chi}_l\qj dz_l\psi_l$ on $\SO_l$ indeed define a global $1$\=/form $\lp\chi,\psi\rp$ on $\X$:
\begin{equation*}
\Bar{\chi}_m \qj dz_m \psi_m = \Bar{\chi}_l \Bar{f}_{ml}^{-1} \qj \tfrac{dz_m}{dz_l} dz_l f_{ml} \tfrac{dz_l}{dz_m} \psi_l = \Bar{\chi}_l \qj dz_l \psi_l.
\end{equation*}
We summarize the properties of this $(1,0)$\=/form in the following Lemma:
\begin{lemma}\phantom{boo}
\label{lem:10 pairing properties}
\begin{enumeratethm}
\item
For any function $g:\X\to\mathbb{C}$ the pairing obeys
\begin{align*}
\lp g\chi,\psi\rp
= \lp\chi,g\psi\rp.
\end{align*}

\item 
For any function $h:\X\to\qat$ the pairing obeys
\begin{align*}
\lp\chi h,\psi\rp=\Bar{h}\lp\chi,\psi\rp
\quad\text{and}\quad
\lp\chi,\psi h\rp=\lp\chi,\psi\rp h.
\end{align*}

\item\hspace{30mm}
$\displaystyle{\overline{\lp\chi,\psi\rp}=-\lp\psi,\chi\rp.}$

\item
For any $1<p<2$ and any $\chi\in H^0(\X,\sob{1,p}{E^{-1}})$ and $\psi\in H^0(\X,\sob{1,\frac{2p}{3p-2}}{KE})$ the coefficients of $\lp\chi,\psi\rp$ belong to $\sobolev{1,1}\loc$.

\item
The local $1$\=/forms $dz_l\U_l$ with $\U_l\in\banach{2}\loc(\SO_l,\qat^-)$ represent on the open sets of the cover $\{\SO_l\}$ the potential of the $(1,0)$ part of a connection on $E$ which anti-commutes with the complex structure, if and only if $\U_l$ represents a global potential $\U\in\pot{KE}^-$

\item
The local functions $\U_l\in\banach{2}\loc(\SO_l,\qat)$ on the open sets of the cover $\{\SO_l\}$ represent a global potential $\U\in\pot{KE}$ if and only if the local functions $(\U\sd)_l = \qk\Bar{\U}_l\qk$ represent a global potential in $\pot{E^{-1}}$ which is denoted by $\U\sd$.

\item
For any triple $(\U,\chi,\psi)\in\pot{KE}\times H^0(\X,\Q{E^{-1},\U\sd})\times H^0(\X,\Q{KE,\U})$ the $1$\=/form $\lp\chi,\psi\rp$ is closed.

\item 
Let $\Omega\subset\mathbb{C}$ an open neighborhood of $z_0$ and $(\U,\chi,\psi)$ as in (vii) on $\X=\Omega\setminus\{z_0\}$. 
The residue $\res_{z_0}\lp\chi,\psi\rp$ is defined as $\frac{1}{2\pi\qi}$ times the integral of $\lp\chi,\psi\rp$ along a closed smooth path in $\X$ surrounding $z_0$ once in anti-clockwise order and does not depend on the choice of the path.
\end{enumeratethm}
\end{lemma}
\begin{proof}
The statements~(i)-(iii) follow algebraically from the local definition.
The proof of~(iv) follows from Lemma~\ref{lem:sobolev regularity}. 
Specifically, the Sobolev embedding theorem says that $H^0(\X,\sob{1,p}{E^{-1}})\subset H^0(\X,\ban{\frac{2p}{2-p}}{E^{-1}})$ and $H^0(\X,\sob{\frac{2p}{3p-2}}{KE})\subset H^0(\X,\ban{\frac{p}{p-1}}{KE})$.
Since the inverse exponents add up to $\frac{1}{p}+\frac{p-1}{p}=1=\frac{2-p}{2p}+\frac{3p-2}{2p}$, the Sobolev product rule and Hölder's inequality imply (iv).

To prove~(v) let $dz_l \U_l$ represent an anti-commuting $(1,0)$\=/potential on $E$, as in Chapter~\ref{chapter:darboux}.
It gives rise to a $(0,1)$\=/potential on $KE$.
The representatives
\begin{align*}
d\Bar{z}_m \U_m
= d\Bar{z}_l \overline{\tfrac{dz_m}{dz_l}} \tfrac{dz_l}{dk_m} f_{ml} \U_l f_{ml}^{-1}
= d\Bar{z}_l \tfrac{dz_l}{dz_m} f_{ml} \U_l \left(\tfrac{dz_l}{dz_m}f_{ml}\right)^{-1}
\end{align*}
transform according to Equation~\eqref{eq:hom transformation}, as required.

For (vi), due to~\eqref{eq:inverse derivative}, the local expressions $\qk\Bar{\U}_l\qk$ define the $(1,0)$\=/potential $\U\sd$ of the section $\chi$ of $E^{-1}_\qat$.
To show that it is indeed a global potential in $\pot{E^{-1}}$ we calculate
\begin{align*}
d\Bar{z}_m(\U\sd)_m = d\Bar{z}_l \overline{\tfrac{dz_m}{dz_l}} \,\qk\,\overline{\tfrac{dz_l}{dk_m} f_{ml} \U_l f_{ml}^{-1}}\,\qk
= d\Bar{z}_l \overline{\tfrac{dz_m}{dz_l}} \,f_{ml}^{-1} \qk \tfrac{dz_l}{dz_m} \Bar{\U}_l \qk f_{ml}
= d\Bar{z}_l f_{ml}^{-1} \qk\Bar{\U}_l\qk f_{ml}.
\end{align*}
And in fact these are equivalences, in the sense that either of the potentials $\U$ and $\U\sd$ determine the other.
Thus a holomorphic structure naturally induces a holomorphic structure on the paired $\qat$\=/line bundle, which proves~(vi).

(vii) follows from Theorem~\ref{thm:darboux}(2)(a).

(viii) 
As discussed above Equation~\eqref{eq:def line integral}, the integrals of $\lp\chi,\psi\rp$ along continuously differentiable immersed cycles $\gamma$ are well defined. 
Due (vii) this integral only depends on homotopy class of $\gamma$. 
This proves~(viii).
\end{proof}

We now give the definition of paired bundles: two holomorphic $\qat$\=/line bundles are \emph{paired} if their underlying bundles tensor to give $K$ and the potentials of their holomorphic structures obey (vi).
Note that either potential determines the other.

From the first property, if we change the line bundle $E$ by multiplying its cocycles by a non-vanishing complex-valued function $g$ then the $(1,0)$ pairing is invariant, since $\lp g^{-1}\chi,g\psi\rp = \lp\chi,\psi\rp$.
We use the superscript $(1,0)$ for this pairing because in some sense it is a $(1,0)$\=/form, although not in the sense of Definition~\ref{def:forms notation}: instead of a complex structure acting on $\lp \chi,\psi\rp$ from the left it acts `internally' with $-\ast \lp \chi,\psi\rp = \lp \chi,\qi \psi\rp$.
Related, we can view the pairing as the $(1,0)$\=/derivative of $\phi$, since locally the pairing has the same formula as~\ref{thm:darboux}(2)(a) with $\psi=(\partial+\B+\U)\phi$.
This is made more precise in the following corollary, which merges the local Darboux transformation to a global transformation:
\begin{corollary}[Global quaternionic Darboux transformation]\label{cor:global darboux}
\index{Darboux transformation}
Let $\X$ be a Riemann surface, $E$ a holomorphic $\mathbb{C}$\=/line bundle on $\X$ and $\V\in\pot{E}^-$ and $\upsilon\in H^0(\X,\Q{E,\V})$ a section without roots.
The local Darboux transformation in Theorem~\ref{thm:darboux} shows that $\chi:=(\overline{\qj\upsilon})^{-1}\in H^0(\X,\Q{E^{-1},\U\sd})$ for $\U\sd\in\pot{E^{-1}}^-$.
Moreover, another $\phi\in H^0(\X,\Q{E,\V})$ is transformed into a paired $\psi\in H^0(\X,\Q{KE,\U})$ with
\begin{equation}\label{eq:weierstrass}
dF=d(\upsilon^{-1}\phi)=\lp\chi,\psi\rp.
\end{equation}
We call $(E,\V,\upsilon,\phi)\mapsto(E,\U,\chi,\psi)$ the \emph{global Darboux transformation}.
If additionally $\X$ is compact, then the global $\banach{2}$\=/norms~\eqref{eq:square integrable} fulfill 
\begin{align}
\label{eq:hopf norms}
\|\U\|_2^2-\|\V\|_2^2=\|\U\sd\|_2^2-\|\V\|_2^2&=\pi\deg(E).
\end{align}
\end{corollary}
\begin{proof}
As standard, we work with the cover from Remark~\ref{rem:special cover}.
The first part follows just by applying Theorem~\ref{thm:darboux} on each chart.
The local sections on $\SO_l$ are defined as the composition of $z_l$ with the corresponding function on $\Omega_l$.
For example, the local representatives $\chi_l = (\overline{\qj\upsilon_l})^{-1}$ define a section $\chi$ of $E^{-1}_\qat$, since:
\begin{align*}
\chi_m
&= (\overline{\qj f_{ml} \upsilon_l})^{-1} 
= \qj \overline{f_{ml}}^{-1} \overline{\upsilon_l}^{-1} = f_{ml}^{-1} \chi_l.
\end{align*}
It has no roots because none of the local representatives do.
The fact that the local functions $\U_l$ and $\U\sd_l$ transform so as to define potentials $\U\in\pot{KE}^-$ and $\U\sd\in\pot{E^{-1}}^-$ on $\X$ are (v) and (vi) of Lemma~\ref{lem:10 pairing properties}.
So Theorem~\ref{thm:darboux}(d) tells us that $\chi$ is a global holomorphic section of $\Q{E^{-1},\U\sd}$.
Likewise $\psi$ is well-defined, and global holomorphic section of $\Q{KE,\U}$ due to Theorem~\ref{thm:darboux}(e).
Finally  Theorem~\ref{thm:darboux}~(1)(f) gives~\eqref{eq:weierstrass}.

For the second part of the corollary, where we assume the compactness of $\X$, we want to integrate Theorem~\ref{thm:darboux}~(1)(c) over $\X$.
As the potentials locally obey $d(\B_l dz_l)=2\ci (|\V_l|^2-|\U_l|^2)\,\dmu$, you might be tempted to conclude from Stokes' theorem that the integral is zero.
But $\B_l dz_l$ does not define a global $1$\=/form, instead it follows from Theorem~\ref{thm:darboux}~(1)(a) that
\[
dz_m \B_m 
= dz_m \partial_m(-2\ln |\upsilon_m| )
= dz_l\B_l - dz_m f_{ml}^{-1}\partial_m f_{ml}
= dz_l\B_l - f_{ml}^{-1}df_{ml}.
\]
Observe that this is Equation~\ref{eq:potential transformation}, so $\B_l$ defines the commuting part of a connection.
The logarithmic derivative $f_{ml}^{-1}df_{ml}$ of the cocycle $f_{ml}$ of $E$ defines a cocycle of holomorphic $1$\=/forms on $\X$, whose residue (compare~\cite[\S17.1]{Fo}) is equal to $\deg(E)$.
Therefore the following equation proves~\eqref{eq:hopf norms}:
\[
\frac{1}{2\ci}\int_\X d(\B_l dz_l)
= \pi (- \deg(E)).\qedhere
\]
\end{proof}

As we have stated above, we think of the Darboux transformation $\psi = (\partial + \B + \U)\phi$ as a holomorphic derivative.
In this lemma we show that its behavior with respect to roots and poles of a $\V$\=/holomorphic section is similar to that of $\partial$ in complex analysis.

\begin{lemma}
\label{lem:change of order}
\index{Pole}
Let $\Omega\subset\mathbb{C}$ be a bounded open neighborhood of $z_0$, $\V\in\banach{2}(\Omega,\qat^-)$, let $\upsilon$ be $\V$\=/holomorphic without roots on $\Omega
$ and $\U\in\banach{2}\loc(\Omega,\qat^-)$ such that $\chi=(\overline{\qj\upsilon})^{-1}$ is $\U\sd$\=/holomorphic on $\Omega$. Then the following statements hold:
\begin{enumeratethm}
\item For a meromorphic section $\psi$ of $\Q{KE,\U}$ with first order pole at $z_0$ the residue $\res_{z_0}\lp\chi,\psi\rp$ is not zero.
\item For all $\phi\in H^0(\Omega\setminus\{z_0\},\Q{E,\V})$ with finite $\ord_{z_0}(\phi)\in\mathbb{Z}$ the global Darboux transformed $\psi=(\partial-(\partial\upsilon)\upsilon^{-1})\phi$ defined on $\X=\Omega\setminus\{z_0\}$ in Corollary~\ref{cor:global darboux}, extends to a meromorphic section of $\Q{KE,\U}$ on $\Omega$ such that
\begin{equation*}
\ord_{z_0}(\psi)\begin{cases}=\ord_{z_0}(\phi)-1&\text{if }\ord_{z_0}(\phi)\ne0.\\\ge0&\text{if }\ord_{z_0}(\phi)=0.\end{cases}
\end{equation*}
\item The Darboux transformation $\phi\mapsto\psi=\big((\partial-\partial\upsilon)\upsilon^{-1}\big)\phi$ defines a surjective map from the stalk of meromorphic sections of $\Q{E,\V}$ at $z_0$ onto the space of all germs $\psi$ of meromorphic sections of $\Q{KE,\U}$ at $z_0$ with $\res_{z_0}\lp\chi,\psi\rp=0$.
\end{enumeratethm}
\end{lemma}
\begin{proof}
The residue $\res_{z_0}\lp\chi,\psi\rp$ coincides with $\res_{z_0}\lp\chi',\psi'\rp$ if $\chi'-\chi$ is a germ of $\Q{E^{-1},\U\sd}$ at $z_0$ with $\ord_{z_0}(\chi'-\chi)>0$ and $\psi'-\psi$ is a germ of $\Q{KE,\U}$ at $z_0$. 
Hence the statement in (i) is independent of the choice of the $\U\sd$\=/holomorphic germ $\chi$ with $\ord_{z_0}(\chi)=0$ and the germ $\psi$ of a meromorphic section of $\Q{KE,\U}$ with $\ord_{z_0}(\psi)=-1$. 
Due to Lemma~\ref{step 2}, for $\alpha\in\qat\setminus\{0\}$ and sufficiently small $r>0$ the elements $\chi_r=(\unity-\Op{I}_{\Omega}\U\sd|_{B(z_0,r)})^{-1}\alpha$ and $\psi_r=\frac{1}{z-z_0}\big(\unity-\Op{I}_{\Omega}\frac{\overline{z-z_0}}{z-z_0}\U|_{B(z_0,r)}\big)^{-1}\alpha$ have these properties on $B(z_0,r)$. 
So it suffices to prove the statement for these $\chi_r$ and $\psi_r$ for small $r>0$. 
By the definition in Lemma~\ref{lem:10 pairing properties}(viii) the residue $\res_{z_0}\lp\chi_r,\psi_r\rp$ is equal to the integral of $\frac{1}{2\pi\qi}\lp\chi_r,\psi_r\rp$ along any closed smooth cycle in $\Omega\setminus\{z_0\}$, which surrounds $z_0$ once in anti-clockwise order.
In particular, the closed smooth cycle may be fixed independent of $r$.
Since $\U\sd|_{B(z_0,r)}$ and $\U|_{B(z_0,r)}$ converge to $0$ in $\banach{2}(\Omega,\qat^-)$ for $r\downarrow 0$, the functions $\chi_r$ and $(z-z_0)\psi_r$ converge in $\banach{\frac{2p}{2-p}}(\Omega,\qat)$ for any $1<p<2$ to the constant $\alpha$. 
Since they are $\U\sd|_{B(z_0,r)}$ and $\frac{\overline{z-z_0}}{z-z_0}\U|_{B(z_0,r)}$\=/holomorphic, respectively, they also converge in $\sobolev{1,p}(\Omega,\qat)$ to $\alpha$. 
The arguments in the proof of Lemma~\ref{lem:10 pairing properties}(viii) show that the residue $\res_{z_0}\lp\chi_r,\psi_r\rp$ depends continuously on $r$ and converges in the limit $r\downarrow0$ to $\res_{z_0}\Bar{\alpha}\qj dz(z-z_0)^{-1}\alpha=\Bar{\alpha}\qj\alpha\ne0$. 
This shows $\res_{z_0}\lp\chi_r,\psi_r\rp\ne0$ for small $r>0$ and completes the proof of~(i).

The statement~(ii) is local and does not depend on the open neighborhood $\Omega$ of $z_0$. We first give the proof for $\ord_{z_0}(\phi)=0$. In this case $\ord_{z_0}(\psi)\ge0$ since $\psi$ is $\U$\=/holomorphic if $\phi$ is $\V$\=/holomorphic.

For $\ord_{z_0}(\phi)>0$ we may choose $\Omega$ such that $\phi\in H^0(\Omega,\Q{E,\V})$. Due to Lemma~\ref{quotient dimension} there exists for any $L\in\mathbb{N}$ an open neighborhood $\Omega$ of $z_0$ such that for any $1\le l\le L$ there exist $\phi\in H^0(\Omega,\Q{E,\V})$ with $\ord_{z_0}(\phi)=l$. In Corollary~\ref{cor:global darboux} we construct the operator
\begin{align}\label{darboux transform}
\partial-(\partial\upsilon)\upsilon^{-1}&:H^0(\X,\Q{E,\V})\to H^0(\X,\Q{KE,\U}),&\phi\mapsto\psi=\big(\partial-(\partial\upsilon)\upsilon^{-1}\big)\phi
\end{align}
It has the kernel $\upsilon\qat$. Hence the restriction to $\{\phi\in H^0(\Omega,\Q{E,\V})\mid1\le\ord_{z_0}(\phi)\}$ is an isomorphism. Now we prove inductively the following inequality
\begin{align}
\label{eq:order inequality}
\ord_{z_0}(\phi)-1&\le\ord_{z_0}(\psi).
\end{align}
For $\ord_{z_0}(\phi)=1$ this is clear. Let us assume that this inequality holds holds for $\phi$ with $\ord_{z_0}(\phi)=l$ and let $\phi'\in H^0(\Omega,\Q{E,\V})$ have $\ord_{z_0}(\phi')=l+1$. We shrink $\Omega$ such that $\phi$ has no other root on $\Omega$ besides $z_0$. In this case we may conceive $\phi$ and $\phi'$ as holomorphic sections of another holomorphic $\qat$\=/line bundle on $\Omega$ such that $\phi$ has no roots. Consequently Theorem~\ref{thm:darboux} applies and yields $\ord_{z_0}\big(d(\phi^{-1}\phi')\big)\ge0$. We conclude that~\eqref{eq:order inequality} holds for $\phi'$ and $\psi'=\big(\partial-(\partial\upsilon)\upsilon^{-1}\big)\phi'$:
\begin{align*}
\ord_{z_0}(\psi')&=\ord_{z_0}\big(\lp\chi,\psi'\rp\big)=\ord_{z_0}\big(d\big(\upsilon^{-1}\phi'\big)\big)=\ord_{z_0}\big(d\big(\upsilon^{-1}\phi\phi^{-1}\phi'\big)\big)\\&=\ord_{z_0}\left(\lp\chi,\psi\rp\big(\phi^{-1}\phi'\big)+\big(\upsilon^{-1}\phi\big)d\big(\phi^{-1}\phi'\big)\right)\\
&\ge\min\Big\{\ord_{z_0}\Big(\lp\chi,\psi\rp\big(\phi^{-1}\phi'\big)\Big),\ord_{z_0}\Big(\big(\upsilon^{-1}\phi\big)d\big(\phi^{-1}\phi'\big)\Big)\Big\}\\
&\ge\min\big\{\ord_{z_0}(\psi)+1,l\big\}=l=\ord_{z_0}(\phi')-1.
\end{align*}
This finishes the proof of~\eqref{eq:order inequality} for $\ord_{z_0}(\phi)>0$. Next we prove that the inequality~\eqref{eq:order inequality} for $\ord_{z_0}(\phi)>0$ implies in the same cases equality in~\eqref{eq:order inequality}. In fact, this inequality shows for $l\in\mathbb{Z}^+$ that~\eqref{darboux transform} maps $\{\phi\in H^0(\Omega,\Q{E,\V})\mid l\le\ord_{z_0}l\}$ into $\{\psi\in H^0(\Omega,\Q{KE,\U})\mid l-1\le\ord_{z_0}(\psi)\}$. Due to Lemma~\ref{quotient dimension} for any $l$ and any sufficiently small $\Omega$, the quotients of $\{\phi\in H^0(\Omega,\Q{E,\V})\mid 1\le\ord_{z_0}(\phi)\}$ and $H^0(\Omega,\Q{KE,\V})$ by these subspaces, respectively, are $(l-1)$\=/dimensional. Hence the inequality~\eqref{eq:order inequality} implies that~\eqref{darboux transform} induces for appropriate $\Omega$ a surjective map
\begin{multline*}
\big\{\phi\in H^0(\Omega,\Q{E,\V})\mid 1\le\ord_{z_0}(\phi)\big\}\Big/\big\{\phi\in H^0(\Omega,\Q{E,\V})\mid l+1\le\ord_{z_0}(\phi)\big\}\\
\to H^0(\Omega,\Q{KE,\V})\Big/\big\{\psi\in H^0(\Omega,\Q{KE,\U})\mid l\le\ord_{z_0}(\psi)\big\}.
\end{multline*}
Since both spaces have for all $l=1,\ldots,L$ the dimension $l$ the corresponding maps are isomorphisms. Hence equality must hold in~\eqref{eq:order inequality} for $\ord_{z_0}(\phi)\in\{1,\ldots,L\}$. This finishes the proof of (ii) for $\ord_{z_0}(\phi)>0$.

For $l=-\ord_{z_0}(\phi)>0$ we again prove the inequality~\eqref{eq:order inequality}.
As a first step, we show that $\psi$ is meromorphic with a pole of order $\le l+1$.
We claim that this is equivalent to the statement that $\lp\chi',\psi\rp$ is a closed current in $\forms{1}{}\banach{1}\loc(\Omega,\qat)$ for $\chi'\in H^0(\Omega,\Q{E^{-1},\U\sd})$ with $\ord_{z_0}(\chi')=l+1$ and no other roots on $\Omega$. To prove the claim we first note such $\chi'$ exist by Lemma~\ref{quotient dimension} on sufficiently small $\Omega$.
Using Lemma~\ref{lem:10 pairing properties}(i) we write  $\lp \chi', \psi\rp = \lp (z-z_0)^{-l-1}\chi', (z-z_0)^{l+1}\psi\rp$. Since $(z-z_0)^{-l-1}\chi'$ is $\U\big((l+1)z_0\big)$\=/holomorphic, due to the local Darboux transformation Theorem~\ref{thm:darboux} and Wente's inequality~\ref{wente} its inverse belongs to $\banach{\infty}\loc(\Omega,\qat)$ and $\lp\chi',\psi\rp\in\forms{1}{}\banach{1}\loc(\Omega,\qat)$ implies $(z-z_0)^{l+1}\psi\in\banach{1}\loc(\Omega,\qat)$. Now Corollary~\ref{cor:weyls lemma 2}(v) shows the claim: $\ord_{z_0}(\psi)\ge-l-1$ if and only if $\lp\chi',\psi\rp$ is a closed current.

Corollary~\ref{cor:global darboux} applied to $G=\chi^{-1}\chi'$ gives $dG=\big\lp(\overline{\qj\chi})^{-1},\big(\partial-(\partial\chi)\chi^{-1}\big)\chi'\big\rp=-\big\lp\upsilon,\big(\partial-(\partial\chi)\chi^{-1}\big)\chi'\big\rp$ with $\big(\partial-(\partial\chi)\chi^{-1}\big)\chi'\in H^0(\Omega,\Q{KE^{-1},\V\sd})$.
Furthermore, $\ord_{z_0}(\chi')=l+1\in\mathbb{N}$ implies $\ord_{z_0}\big(\big(\partial-(\partial\chi)\chi^{-1}\big)\chi'\big)=l$ by the forgoing arguments. Now we calculate with Lemma~\ref{lem:10 pairing properties}(i)-(ii)
\begin{align*}
\lp\chi',\psi\rp&=\lp\chi G,\psi\rp=\Bar{G}\lp\chi,\psi\rp=d\big(\Bar{G}(\upsilon^{-1}\phi)\big)-d\Bar{G}(\upsilon^{-1}\phi)\\&=d\big(\Bar{G}(\upsilon^{-1}\phi)\big)+\big\lp\big(\partial-(\partial\chi)\chi^{-1}\big)\chi',\upsilon\big\rp\big(\upsilon^{-1}\phi\big)\\&=d\big(\Bar{G}(\upsilon^{-1}\phi)\big)+\big\lp\big(\partial-(\partial\chi)\chi^{-1}\big)\chi',\phi\big\rp\\&=d\big(\Bar{G}(\upsilon^{-1}\phi)\big)+\big\lp(z-z_0)^{(-l)}\big(\partial-(\partial\chi)\chi^{-1}\big)\chi',(z-z_0)^l\phi\big\rp.
\end{align*}
Due to Wente's inequality~\ref{wente} and the product rule of Sobolev functions, the function $\Bar{G}(\upsilon^{-1}\phi)=\Bar{\chi}'\Bar{\chi}^{-1}\upsilon^{-1}\phi=-\Bar{\chi}'\qj\phi=-\overline{(z-z_0)^{-l-1}\chi'}\qj(z-z_0)^{l+1}\phi$ belongs to $\sobolev{1,2}\loc(\Omega,\qat)\cap\banach{\infty}\loc(\Omega,\qat)$. 
It's exterior derivative belongs to $\banach{2}\loc(\Omega,\qat)$.
Therefore the first term on the right hand side is a closed current. 
The second term is the exterior derivative of an admissible map, because we have a pairing of paired holomorphic sections, Lemma~\ref{lem:10 pairing properties}(vii). 
Therefore $\lp\chi',\psi\rp$ is a closed current in $\forms{1}{}\banach{1}\loc(\Omega,\qat)$. 
This proves~\eqref{eq:order inequality} for $\ord_{z_0}(\phi)<0$.
Observe additionally that the inequality implies that $\psi$ extends to a meromorphic section of $\Q{KE,\U}$. 

Before we complete the case $\ord_{z_0}(\phi)<0$ we need to prove that $\ord_{z_0}(\psi)\ne-1$ holds in the situation of~(ii). 
This follows from~(i) and the statement that $\res_{z_0}\lp\chi,\psi\rp=0$ for $\psi$ in the image of~\eqref{darboux transform}. 
This latter statement is indeed true, since the integral along smooth closed cycles in $\Omega\setminus\{z_0\}$ of the derivative $\lp\chi,\psi\rp$ of the admissible map $\upsilon^{-1}\psi\in\sobolev{2,p}\loc(\Omega\setminus\{z_0\},\qat)\subset C(\Omega\setminus\{z_0\},\qat)$ vanishes.

To finish the proof of~(ii) for $\ord_{z_0}(\phi)<0$ we remark that the inequality~\eqref{eq:order inequality} together with $\ord_{z_0}(\psi)\ne-1$ implies that~\eqref{darboux transform} induces an injective map
\begin{multline*}
\big\{\phi\in H^0(\X,\Q{E,\V})\mid-l\le\ord_{z_0}(\phi)\big\}\Big/\big\{\phi\in H^0(\X,\Q{E,\V})\mid 0\le\ord_{z_0}(\phi)\big\}\to\\
\big\{\psi\in H^0(\X,\Q{KE,\V})\mid-l-1\le\ord_{z_0}(\psi)\big\}\Big/\big\{\psi\in H^0(\X,\Q{KE,\U})\mid -1\le\ord_{z_0}(\psi)\big\}.
\end{multline*}
Since both spaces have for $l\in\{1,\ldots,L\}$ the dimension $l$ the corresponding maps are isomorphisms. Therefore equality must hold in~\eqref{eq:order inequality} for $-\ord_{z_0}(\phi)\in\{1,\ldots,L\}$. This finishes the proof of~(ii).


Part~(iii) contains the statement that $\res_{z_0}\lp\chi,\psi\rp$ vanishes for $\psi$ in the image of~\eqref{darboux transform}. This statement has been shown in the proof of~(ii). In order to compete the proof of~(iii) it remains to show that for all germs $\psi$ of meromorphic sections of $\Q{KE,\U}$ at $z_0$ with vanishing $\res_{z_0}\lp\chi,\psi\rp$ there exists a germ of a meromorphic section $\phi$ of $\Q{E,\V}$ at $z_0$ with $\psi=\big(\partial-(\partial\upsilon)\upsilon^{-1}\big)\phi$. Due to the local Darboux transform, Theorem~\ref{thm:darboux}, this is true for $\U$\=/holomorphic $\psi$. Since by Lemma~\ref{quotient dimension} the quotient of meromorphic germs $\psi$ at $z_0$ with poles of order at most $l+1$ divided by germs $\psi$ with $\ord_{z_0}(\psi)\ge-1$ has the same quaternionic dimension $l$ as the quotient of meromorphic germs $\phi$ at $z_0$ with poles of order at most $l$ divided by holomorphic germs $\phi$, this follows from~(ii).
\end{proof}
With the Corollary~\ref{cor:global darboux}, applying the global Darboux transformation to the quaternionic Kodaira representation $(E,\V,\upsilon,\phi)$ of an admissible map gives us the data $(E,\U,\chi,\psi)$.
This is called the quaternionic Weierstraß representation.
The following theorem addresses to what extent such a quadruple is in correspondence with admissible maps.
Compare with~\cite[Theorem~4.3]{PP}.

\begin{theorem}[Quaternionic Weierstraß Representation]
\label{thm:weierstrass}
\index{Weierstraß representation}
For any admissible map $F:\X\to\qat$ on a Riemann surface $\X$ there exist a holomorphic $\mathbb{C}$\=/line bundle $E$, a potential $\U\in\pot{KE}^-$, $\chi\in H^0(\X,\Q{E^{-1},\U\sd})$ and $\psi\in H^0(\X,\Q{KE,\U})$ with $dF =\lp\chi,\psi\rp$ as in~\eqref{eq:weierstrass}.
There is a unique choice of $E$ up to isomorphism such that $\chi$ has no roots. 
For any $E$ in this isomorphism class the possible quadruples are $(E,g\U g^{-1}, g^{-1}\chi, g\psi)$ for a holomorphic function $g : \X \to \mathbb{C}$ without roots.

Vice versa, for any such quadruple $(E,\U,\chi,\psi)$ any anti-derivative $F$ of $dF=\lp\chi,\psi\rp$ defines an admissible map on the universal covering $\Tilde{\X}$ of $\X$.
\end{theorem}
\begin{proof}
The Quaternionic Kodaira representation in Theorem~\ref{thm:kodaira} establishes a one-to-one correspondence between quadruples $(E,\V,\upsilon,\phi)$ and admissible maps.
Further, the global Darboux transformation Corollary~\ref{cor:global darboux} transforms $(E,\V,\upsilon,\phi)$ into $(E,\U,\chi,\psi)$.
Since the local Darboux transformation has an inverse, such data is unique up to adding a constant multiple $\upsilon \beta$ to $\phi$, which corresponds to translating $F$ by $\beta \in \qat$.
Any root of $dF =\lp\chi, \psi\rp$ must appear in either $\chi$ or $\psi$.
We have seen that $\lp g^{-1}\chi,g\psi\rp = \lp \chi,\psi\rp$, for holomorphic $g$.
This is the only manner that quadruples $(E,\U,\chi,\psi)$ corresponding to the same $dF$ can differ from each other.
If $g$ has no roots, then this is just changing $E$ without changing its isomorphism class.
If $g$ has roots then we can move roots between $\chi$ and $\psi$, so that $\chi$ has no roots.
If $dF$ has no roots, then the quadruple $(E,\U,\chi,\psi)$ is unique up to isomorphism.
%

So it remains to recover $F$ from its Weierstraß representation.
For all such maps $F$ with~\eqref{eq:weierstrass} the left normal is locally equal to $N=\upsilon^{-1}\qi\upsilon=-\Bar{\chi}\qi\Bar{\chi}^{-1}$.
When $\chi$ has a root let $D=(\chi)$ be the divisor of $\chi$ (see Lemma~\ref{quotient dimension}).
Then $\chi$ is a section without roots of $E^{-1}(D)_\qat$  and $\psi$ is a section with roots of $KE(D)_\qat$.
By Corollary~\ref{cor:global darboux}, $\upsilon = -(\overline{\qj \chi})^{-1}$ is a $\V$\=/holomorphic section of $E(D)_\qat$ with potential $\V\in \pot{E(D)}^-$.
Due to Theorem~\ref{thm:kodaira normal} $N = \upsilon^{-1} \qj \upsilon$ belongs to $\sobolev{1,2}\loc (\X, \qat)$ and due to the Theorems~\ref{thm:darboux} and~\ref{thm:kodaira} any anti-derivative $F$ of $dF = \lp\chi, \psi\rp$ is admissible.
By Lemma~\ref{lem:10 pairing properties} this form is closed and has an anti-derivative on the universal covering $\Tilde{\X}$ on $\X$.
\end{proof}
\begin{remark}
\label{rem:weierstrass}
Due to Wente's inequality Lemma~\ref{wente} the Quaternionic Weierstraß representation implies that all admissible maps belong to $\sobolev{2,2}\loc$ and the conformal factor is the continuous function $|\chi|\cdot|\psi|$ which vanishes only at roots of $\chi$ or $\psi$.
We call such points \emph{branch points}.
In particular, an admissible map belongs to the class of conformal $\sobolev{2,2}\loc$\=/maps (see Remark~\ref{re:weakly conformal}) if and only if $dF$ has no roots.
\end{remark}

We should explain the connection to the classical Weierstraß-Enneper representation for minimal surfaces.
Recall that for a holomorphic function $f$ and a meromorphic function $g$ such that $fg^2$ is holomorphic then
\begin{align*}
z \mapsto \frac{1}{2} \Real \int_{z_0}^z \left( f (1-g^2), \qi f(1+g^2), 2fg \right) \; dw \; \in \mathbb{R}^3
\end{align*}
is a minimal surface~\cite[p.~395]{Sp}.
Of relevance is the following formula~\cite[p.~397]{Sp} for the normal of this surface
\[
N(z) = \frac{1}{1 + |g(z)|^2} \Big( 2 \Real g(z), 2 \Imag g(z), |g(z)|^2 - 1 \Big) \in \mathbb{S}^2.
\]
We see that $g$ is just the normal of the surface, where $\mathbb{C}\cup\{\infty\}$ has been identified with $\mathbb{S}^2$ via stereographic projection.
Hence the classical representation, like our representation, is foremost a parameterization of the normal of the surface.
A further explanation of this point of view can be found in~\cite[\S2]{Ta3}.

\begin{example}[Plane]
\label{eg:plane weierstrass}
\index{Plane!Weierstraß representation}
Continuing from Example~\ref{eg:plane kodaira}, we have an admissible map into the plane $F(z) = \alpha f(z) \beta$ for $f$ a $\mathbb{C}$\=/valued holomorphic function.
Its Kodaira data is $\V = 0, \upsilon = \alpha^{-1}, \phi = f(z) \beta$.
Actually, all the potentials vanish since $d\upsilon \equiv 0$.
The Weierstraß data is $\chi = \qj \Bar{\upsilon}^{-1} = \qj \Bar{\alpha}$ and $\psi = \partial (f\beta) = f'\beta$. 
In particular we see that $\psi$ is just derivative of $f$ in the complex sense, which demonstrates our interpretation of the Darboux transformation as a generalization of the derivative of a complex function to the quaternionic setting.
The Weierstraß representation of $F$ is likewise a repackaging of the holomorphic derivative of $f$, since $dF = \alpha\Bar{\qj}\,\qj dz\, f' \beta = \alpha\, dz\, f' \beta$.
\end{example}

\begin{example}[Catenoid]
\label{eg:catenoid weierstrass}
\index{Catenoid!Weierstraß representation}
Let us compute the Weierstraß representation for the catenoid.
We continue from Example~\ref{eg:catenoid kodaira}.
We have
\[
\chi 
= (\overline{\qj \upsilon})^{-1}
= \overline{(1-z^{-1}\qk)}
= 1+z^{-1}\qk.
\]
In fact, using the potentials from Example~\ref{eg:catenoid holomorphic trivialization} we have
$\B = -z^{-1} (1 + |z|^2)^{-1}$ and $\U = 0$ in the present notation.
We could deduce the other section $\psi$ from $\chi$ and $dF$ (Equation~\eqref{eq:catenoid dF}), but let us instead carry out the Darboux transform from scratch.
The calculation is rather messy however, so we try to avoid denominators.
First we simplify $\phi = \upsilon F$:
\begin{align*}
\qj (1+|z|^2) \phi
&= \Bar{z}\ln(z\Bar{z}) (z\qi + \qj) + (1+|z|^2)(\Bar{z}\qj - \qi)
\end{align*}
Now we differentiate both sides with $\barpartial$ so that a term with $\partial \phi$ is produced:
\begin{align*}
\barpartial (\text{LHS})
&= \barpartial \left[\qj (1+ z\Bar{z}) \phi \right]
= \qj \Bar{z} \phi + \qj (1+ z\Bar{z}) \partial \phi, \\
\barpartial (\text{RHS})
&= \ln(z\Bar{z}) (z\qi + \qj) 
+ \Bar{z}\Bar{z}^{-1} (z\qi + \qj) 
+ 0
+ z(\Bar{z}\qj - \qi)
+ (1+|z|^2)\qj \\
&= 
\ln(z\Bar{z}) (z\qi + \qj) 
+ 2(1 + z\Bar{z})\qj.
\end{align*}
Thus
\begin{align*}
\qj(1 + |z|^2)(\partial + B + U)\phi
&= \barpartial \left[\qj (1+ z\Bar{z}) \phi \right] - \qj z \phi - \qj z^{-1}\phi
= (1 + |z|^2)( \qj + \Bar{z}^{-1}\qi),
\end{align*}
and we conclude that
\[
\psi 
= (\partial + B + U)\phi
= -\qj ( \qj + \Bar{z}^{-1}\qi)
= 1 + z^{-1}\qk.
\]
Both $\chi$ and $\psi$ are $0$\=/holomorphic, which is connected to the fact that the catenoid is a minimal surface.
As we will see below in Equation~\eqref{eq:willmore energy weierstrass}, the $\banach{2}$\=/norm of the Weierstraß potential $U$ is Willmore energy of the admissible map, the integral of square of the mean curvature.
In particular, $\U \equiv 0$ if and only if $F$ is a minimal surface.
The fact that $\chi$ and $\psi$ are equal is also not a coincidence.
This is possible because the surface lies in $\Imag\qat$, see Lemma~\ref{lem:admissible imaginary} below.
\end{example}


We close this chapter by supplementing the global Darboux transformation $(\upsilon,\phi)\mapsto(\chi,\psi)$ in a natural way by another transformation $(\chi,\psi)\mapsto(\zeta,\eta)$.
This other transformation is from the Weierstraß representation of $\Bar{F}$ to its Kodaira representation, and so is in some sense an inverse to the Darboux transform.
Furthermore, we also explain how we may generalize these two transformations $(\upsilon,\phi)\mapsto(\chi,\psi)\mapsto(\zeta,\eta)$ to act on $d$ linear independent global holomorphic sections of various holomorphic $\qat$\=/line bundles, and then how the concatenated transformation may be iterated $d$ times.
This will be an important method in Chapter~\ref{chapter:pluecker}, as well as essential to understand the relationship between the special forms of the Kodaira and Weierstraß potentials of constrained Willmore surfaces, Proposition~\ref{prop:constraint transformation}. 

Consider then a $\qat$\=/vector subspace $H$ of the finite dimensional space of global holomorphic sections of a holomorphic $\qat$\=/line bundle.
We shall associate to $H$ another subspace of dimension $\dim H-1=:d-1$.
More precisely, first we choose an element $\xi\in H\setminus\{0\}$ whose vanishing order at any $x\in\X$ is the minimum of the vanishing orders at $x$ of all elements of $H\setminus\{0\}$.
Lemma~\ref{lem:minimal-order} proves the existence such a $\xi$.
Afterwards we change the underlying holomorphic $\mathbb{C}$\=/line bundle $E$ if necessary, such that $\xi$ has no roots as a section of $E_\qat$.
Hence it defines a unique flat connection given locally by $dz_l(\partial_l-(\partial_l\xi_l)\xi_l^{-1}) + d\Bar{z}_l(\barpartial_l-\V_l)$.
Due to the global quaternionic Darboux Transformation, the $(1,0)$\=/part defines a differential operator
\begin{align}\label{10 differential operator}
\partial_l-(\partial_l\xi)\xi^{-1}&:&H^0(\X,\Q{E,\V})\to H^0(\X,\Q{KE,\U}).
\end{align}
Now we associate to $H_0 = H$ and $\xi$ the image $H_1$ with respect to~\eqref{10 differential operator}.
The dimension of this image is $d-1$ since the kernel of~\eqref{10 differential operator} is spanned by $\xi$.
We may now choose a new $\xi$ from $H_1 \setminus \{0\}$ and iterate.
For any given $d$\=/dimensional $H$ we apply this transformation $d$\=/many times.
In this way we obtain the spaces $H_0=H,H_1,\ldots,H_d=\{0\}$ of dimension $\dim H_l=d-l$.
They depend on the choices of the corresponding $\xi$'s.
We might choose a basis of $H$ such that for any $l$ the $l$\=/th element of this basis is mapped by the composition of the first $l-1$  maps~\eqref{10 differential operator} to the corresponding element of $H_l$ without roots.
The last $d-l$ elements of the basis are then mapped by this composition to a basis of $H_l$, whereas the first $l-1$ element are mapped to zero.
In this way we endow all the spaces $H_l$ with a basis.
This choice is not unique: adding to the $l$th basis element of $H_0$ multiples of the first $l-1$ basis elements does not change its image in $H_l$.

In a second step we supplement the sequence $H_0=H,H_1,\ldots,H_d=\{0\}$ by another sequence $H\pa_0=\{0\}, H\pa_1,\ldots,H\pa_d$ such that any $H\pa_l$ is a $l$\=/dimensional space of global holomorphic sections which are paired with the holomorphic sections in $H_l$.
Moreover, the construction is defined in such a way, that the reversed sequence $H\pa_d$, \ldots, $H\pa_0=\{0\}$ is obtained by a sequence of transformations as described above.
To explain the details we use the notation of the global Darboux transformation:
\begin{lemma}\label{lem:extended darboux}
\index{Darboux transformation}
On a possibly non-compact Riemann surface $\X$, let $(E,\V,\upsilon,\phi)\to(E,\U,\chi,\psi)$ be a Global quaternionic Darboux transformation as in Corollary~\ref{cor:global darboux}.
Furthermore, let the closed $1$\=/form $\lp\zeta,\upsilon\rp$ of some $\zeta\in H^0(\X,\Q{KE^{-1},\V\sd})$ has no periods.
Then there exists $\eta\in H^0(\X,\Q{E^{-1},\U\sd})$ such that $\big(\delbar{E^{-1}}-(\delbar{E^{-1}}\chi)\chi^{-1}\big)\eta=\zeta$.
This $\eta$ is unique up to an additive right multiple of $\chi$.
Further, whenever $\lp\zeta,\phi\rp$ has no periods, neither does $\lp\eta,\psi\rp$.
\end{lemma}
\begin{proof}
Since $\lp\zeta,\upsilon\rp$ has no periods, there exist an admissible map $G$ on $\X$ with $dG=\lp\zeta,\upsilon\rp$.
Taking the conjugate gives $d\Bar{G}=\lp-\upsilon,\zeta\rp$ so that the Kodaira section without roots representing $\Bar{G}$ is $\chi = -(\overline{\qj (-\upsilon)})^{-1}$, as per Corollary~\ref{cor:global darboux}.
This yields $\eta=\chi\Bar{G}\in H^0(\X,\Q{E^{-1},\U\sd})$ and $(\delbar{E^{-1}}-(\delbar{E^{-1}}\chi)\chi^{-1})\eta=\zeta$.
The function $G$ is unique up to an additive constant and as such $\eta$ is unique up the addition of $\chi\qat$.
It remains to show that the closed $1$\=/form $\lp\eta,\psi\rp$ has no periods.
Let $F$ denote as usual the admissible map $F=\upsilon^{-1}\phi$ with $dF=\lp\chi,\psi\rp$.
Since $d(GF)$ and $\lp\zeta,\phi\rp$ have no periods, Lemma~\ref{lem:10 pairing properties}(ii) finishes the proof:
\begin{align*}
\lp\eta,\psi\rp&=\lp\chi\Bar{G},\psi\rp=G\lp\chi,\psi\rp=GdF=d(GF)-(dG)F\\&=d(GF)-\lp\zeta,\upsilon\rp F=d(GF)-\lp\zeta,\upsilon F\rp=d(GF)-\lp\zeta,\phi\rp.\qedhere
\end{align*}
\end{proof}
Let us now inductively construct the spaces $H\pa_0=\{0\}$, $H\pa_1$, \ldots, $H\pa_d$.
We assume that for some $0\le l<d$ the spaces $H\pa_0=\{0\}$, \dots, $H\pa_l$ are constructed, and that the $1$\=/forms $\lp\zeta,\phi\rp$ have no periods for all $(\zeta,\phi)\in H\pa_l\times H_l$.
Let $\Q{E,\V}$ denote the sheaf, such that $H_l\subset H^0(\X,\Q{E,\V})$, let $\upsilon$ be the element of $H$ without roots and let~\eqref{10 differential operator} denote the operator from $H_l$ onto $H_{l+1}$.
Then $H\pa_{l+1}$ is spanned by $\chi=(\overline{\qj\upsilon})^{-1}$ and the elements $\eta$ constructed in the Lemma, where $\zeta$ runs through a basis of $H\pa_l$.
Note that the space $H\pa_{l+1}$ does not depend on additive right multiples of $\chi$ since $\chi$ is already a basis vector.
This gives the definition of $H\pa_{l+1}$.

It remains to close the induction.
Any element of $H\pa_{l+1}$ is the sum of a multiple of $\chi$ and an element $\eta$ from the above lemma.
Any element of $H_{l+1}$ is the image of $\phi \in H_l$ under~\eqref{10 differential operator}.
The final clause of the lemma tells us that the $1$\=/forms $\lp\eta,\psi\rp $ have no periods.
The vanishing of the periods of $1$\=/forms $\lp\chi,\psi\rp$ follows from the global Darboux transform.
Therefore the $\lp\cdot,\cdot\rp$ pairing of any element of $H\pa_{l+1}$ and $H_{l+1}$ has no periods.
Therefore the induction hypothesis is again fulfilled and the inductive construction of $H\pa_0=\{0\},\ldots,H\pa_d$ is finished.
Due to the uniqueness of $\eta$ statement in Lemma~\ref{lem:extended darboux} the second sequence of spaces $H\pa_0,\ldots, H\pa_d=\{0\}$ is completely determined by the first sequence $H_0,\ldots,H_d=\{0\}$.

We now explain the symmetry between the sequence $(H_0,H\pa_0),\dots,(H_d,H\pa_d)$ and the sequence $(H\pa_d,H_d),\dots,(H\pa_0,H_0)$.
The transformation from $H_l$ to $H_{l+1}$ is a Darboux transformation using the flat connection induced by an element $\upsilon$ without roots in $H_l$.
The corresponding element $\chi = (\overline{\qj \upsilon})^{-1}$ in $H\pa_{l+1}$ induces the flat connection given locally by $dz_m(\partial_m-(\partial_m\chi_m)\chi_m^{-1})+ d\Bar{z}_m(\barpartial_m-\U\sd_m)$.
As we see in the lemma, the construction of $\eta \in H\pa_{l+1}$ is like an inverse Darboux transformation of $\zeta \in H\pa_l$, since $(\partial_m-(\partial_m\chi)\chi^{-1})\eta = \zeta$.
Since $\chi$ spans the kernel of the $(1,0)$ part $(\partial_m-(\partial_m\chi)\chi^{-1})$, this operator indeed maps $H\pa_{l+1}$ onto $H\pa_l$.

Let us explain in detail the case $d=2$, where the second transformation is in some sense the inverse of the global Darboux transformation of the admissible map $\Bar{F}$.
The Kodaira representation describes an admissible map $F = \upsilon^{-1}\phi$ in terms of two holomorphic sections $\upsilon,\phi\in H^0(\X,\Q{E,\V})$ and the Weierstraß representation $dF=\lp\chi,\psi\rp$.
The Kodaira sections span a two-dimensional linear subspace $H_0\subset H^0(\X,\Q{E,\V})$.
By definition $\upsilon$ has no roots.
We therefore choose $\upsilon$ as the element which defines the flat connection and the corresponding $(1,0)$ part.
Consequently we associate to $H_0$ the one-dimensional space $H_1$ spanned by $\psi=\big(\partial_m-(\partial_m\upsilon)\upsilon^{-1}\big)\phi$.
The other space $H\pa_1$ is spanned by $\chi=\big(\overline{\qj\upsilon})^{-1}$.
For the second transformation we first replace $\Q{KE,\U}$ by $\Q{KE(-D),\U(-D)}$ such that $\psi$ is a section without roots.
Here $D$ is the root divisor of $\psi$.
The corresponding flat connection has $(0,1)$\=/part $\delbar{}-\U(-D)$.
Let $\zeta=\big(\overline{\qj\psi}\big)^{-1}$ denote the element of $H\pa_2$ without roots.
By applying Lemma~\ref{lem:extended darboux} to the triple $(\psi,\chi,F)$ instead of $(\zeta,\upsilon,G)$ we construct $\eta=\zeta\Bar{F}$.
We read off that $(\zeta,\eta)$ is the Kodaira pair and $(-\psi,\chi)$ the Weierstraß pair of $\Bar{F}$.
In this sense the second transformation $(\chi,\psi)\mapsto(\zeta,\eta)$ is the inverse of the Darboux transformation from the Kodaira pair to the Weierstraß pair of $\Bar{F}$.

\chapter{Admissible Maps into \texorpdfstring{$\mathbb{R}^3$}{R3}}
\label{chapter:3-space}

In this chapter we explore how the Kodaira and Weierstraß representations change under various transformations: covering maps and compositions, reflections and Möbius transformations.
The most important result, Lemma~\ref{lem:admissible imaginary}, gives a condition on the Weierstraß representation that determines whether or not the admissible map is into $\Imag\qat$, which we identify with $\mathbb{R}^3$.
If additionally the roots of $dF$ are all of even order then the Weierstraß representation simplifies.

Let us begin with two easy transformations.
Let $F : \X \to \qat$ be an admissible map with Kodaira representation $(E,\V,\upsilon,\phi)$ and Weierstraß representation $(E, \U, \chi, \psi)$.
Now consider the quadruple $((KE)^{-1}, \U\sd, \psi, \chi)$.
Lemma~\ref{lem:01 pairing properties}(vi) implies it also meets the conditions of Theorem~\ref{thm:weierstrass}, that it is a Weierstraß representation.
In fact Lemma~\ref{lem:01 pairing properties}(iii) even tells us the corresponding admissible map:
\begin{align*}
\lp \psi, \chi \rp
= - \overline{\lp \chi, \psi \rp}
= - \overline{dF}
= d(-\Bar{F}).
\end{align*}
This map $-\Bar{F}$ is the reflection of $F$ across $\Imag\qat$.
It was in fact discussed at the end of the previous chapter, where it was explained that its Kodaira representation has non-vanishing section $- (\overline{\qj\psi})^{-1}$ (but the line bundle must be chosen so that $\psi$ has no roots).
Similarly, one can manipulate the Weierstraß representation of $F$ to derive that of the point reflection $-F$:
it has Kodaira representation $-F = \upsilon^{-1}(-\phi)$ and Weierstraß representation $d(-F) = \lp\chi, -\psi\rp$.

\emph{Möbius transformations}\index{Möbius transformation} arise by extending this second idea.
Consider linear combinations of the Kodaira sections $\phi$ and $\upsilon$:
\begin{align}
\label{eq:quaternionic moebius 1}
\phi&\mapsto\Tilde{\phi}=\phi\alpha+\upsilon\beta\quad\upsilon\mapsto\Tilde{\upsilon}=\phi\gamma+\upsilon\delta,\quad\text{with}\quad\alpha,\beta,\gamma,\delta\in\qat,
\end{align}
By linearity, the new sections $\Tilde{\upsilon},\Tilde{\phi}$ are also $\V$\=/holomorphic sections of $E_\qat$.
The resulting admissible map $\Tilde{F}$ is a quaternionic Möbius transformation of the original:
\begin{align}
\label{eq:quaternionic moebius 2}
F&\mapsto\Tilde{F}=\Tilde{\upsilon}^{-1}\Tilde{\phi}=(\upsilon F\gamma+\upsilon\delta)^{-1}(\upsilon F\alpha+\upsilon\beta)=(F\gamma+\delta)^{-1}(F\alpha+\beta).
\end{align}
An affine transformation corresponds to $\gamma=0$ and is bijective if and only if $\alpha$ and $\delta$ are not vanishing.
We see from the above formula that the surface is rotated $\delta^{-1}F\alpha$ and subsequently translated by $\delta^{-1}\beta$.
The non-affine transformations, $\gamma\neq 0$, are slightly more complicated.
This time it is a bijective transformation if and only if $\delta\gamma^{-1}\alpha-\beta$ does not vanish.
Any point of the surface with $F = -\delta\gamma^{-1}$ is sent to infinity, in which case we no longer have a map from the Riemann surface $\X$ into $\qat$.
The reverse situation is considered in the examples below.

How does a Möbius transformation affect the Weierstraß representation?\index{Weierstraß representation!Möbius transformation} 
From the Darboux transformation we have
\begin{equation}
\label{eq:mobius chi}
\chi \mapsto \Tilde{\chi}
= \overline{(\qj \Tilde{\upsilon})}^{-1}
= \overline{(\qj \upsilon(F\gamma + \delta))}^{-1}
= \chi\overline{(F\gamma + \delta)}^{-1}.
\end{equation}
For the other part of the Darboux transformation, we need the other potentials too:
\begin{align*}
-(\Tilde{\B}_l + \Tilde{\U}_l)\Tilde{\upsilon}_l
&= \partial_l\Tilde{\upsilon}_l
= \partial_l\phi_l \gamma + \partial_l\upsilon_l \delta 
= \psi_l\gamma -(\B_l + \U_l)\phi_l\gamma -(\B_l + \U_l)\upsilon_l\delta \\
&=\psi_l\gamma -(\B_l + \U_l)\Tilde{\upsilon}_l .
\end{align*}
From this we read off that 
\begin{align}
\label{eq:mobius potentials}
\Tilde{\B}_l = \B_l - (\psi_l\gamma\Tilde{\upsilon}_l^{-1})^+,
\quad
\Tilde{\U}_l = \U_l - (\psi_l\gamma\Tilde{\upsilon}_l^{-1})^-,
\quad\text{as well as}\quad
\Tilde{\V}_l = \V_l.
\end{align}
Hence
\begin{align*}
\psi_l \mapsto \Tilde{\psi}_l
&= (\partial_l + \Tilde{\B}_l + \Tilde{\U}_l)\Tilde{\phi}_l
= (\partial_l + \B_l + \U_l - \psi_l\gamma\Tilde{\upsilon}_l^{-1})\Tilde{\phi}_l \\
&= (\partial_l + \B_l + \U_l)(\phi_l\alpha + \upsilon_l\beta) - \psi_l\gamma\Tilde{\upsilon}_l^{-1}\Tilde{\phi}_l \\
&= \psi_l(\alpha - \gamma\Tilde{F}).
\labelthis{eq:mobius psi}
\end{align*}
In the case of an affine transformation the corresponding transformations of the Weierstraß data simplify to $\chi \mapsto \Tilde{\chi} = \chi \Bar{\delta}^{-1}$ and $\psi \mapsto \Tilde{\psi} = \psi\alpha$.

\begin{example}[Round Sphere]
\label{eg:sphere weierstrass}
\index{Round sphere!Weierstraß representation}
Consider the plane $\X = \mathbb{C}$ embedded in $\Imag \qat$ as $F(z) = \qj z$.
This is a simple case of Example~\ref{eg:plane weierstrass}: the potentials are zero, $\upsilon = - \qj$, $\phi = z$, $\chi = 1$, and $\psi = \partial z = 1$.
In terms of Lemma~\ref{lem:admissible imaginary} we have $g \equiv 1$.

Now apply the Möbius transformation with $\alpha = \delta = \qi$ and $\beta = \gamma = 1$.
The new immersion is
\begin{align*}
\Tilde{F} 
&= (\qj z + \qi)^{-1} (\qj z \qi + 1)
= (|z|^2+1)^{-1}(-\qj z - \qi)(\qj z \qi + 1) \\
&= (|z|^2+1)^{-1}((|z|^2-1)\qi - 2\qj z),
\end{align*}
which is the round sphere of Example~\ref{eg:sphere kodaira}.
Formulas~\eqref{eq:quaternionic moebius 1} give the same Kodaira representation as Example~\ref{eg:sphere kodaira}:
\[
\Tilde{\upsilon}
= z\, 1 -\qj\, \qi = z + \qk,
\qquad
\Tilde{\phi}
= z\,\qi - \qj\,1
= z\qi - \qj.
\]
For the sake of completeness, these formulas also give us a Weierstraß representation of the sphere, compare Equation~\eqref{eq:sphere dF},
\begin{align*}
\Tilde{\chi}
&= \chi \overline{(F\gamma + \delta)}^{-1}
= 1 \overline{(\qj z + \qi)}^{-1}
= -(\qj z + \qi)^{-1}, \\
\Tilde{\psi}
&= \psi(\alpha - \gamma\Tilde{F})
= 1\left(\qi - (\qj z + \qi)^{-1} (\qj z \qi + 1)\right)
= -2(\qj z + \qi)^{-1}, \\
\Tilde{\U}
&= \U - (\psi\gamma\Tilde{\upsilon})^-
= 0 - ( 1 \cdot 1 \cdot (z + \qk)^{-1} )^-
= (|z|^2 + 1)^{-1}\qk.
\qedhere
\end{align*}
\end{example}

\begin{example}[Inverted Catenoid]
\label{eg:inverted catenoid}
\index{Catenoid!Inverted catenoid}
Let apply the same Möbius transformation $\alpha = \delta = \qi$ and $\beta = \gamma = 1$ to the catenoid (Example~\ref{eg:catenoid}).
The new surface is 
\begin{align*}
\Tilde{F}(z)
&= \left[ F + \qi \right]^{-1} \left[ F\qi + 1 \right] 
= \left[ z F + z\qi \right]^{-1} \left[ zF\qi + z \right] \\
&= \left[ z(1 + \ln(z\Bar{z})) \qi + (z\Bar{z} + 1)\qj \right]^{-1} (-\qi) \left[ z(1 -\ln(z\Bar{z}))\qi + (z\Bar{z} + 1)\qj \right].
\end{align*}
The presence of the $z$ factor is sufficient to dominate the singularity of $F$ at $z=0$.
Specifically, $zF = z \ln(z\Bar{z}) \qi + (z\Bar{z} + 1)\qj$ is continuous with $\lim_{z \to 0} zF = \qj$ but it is not (strongly) differentiable.
Hence $\Tilde{F}$ extends to a map $\Tilde{F} : \mathbb{P}^1 \to \qat$ with $\Tilde{F}(0) = \Tilde{F}(\infty) = \qi$, called the \emph{compactified inverted catenoid}.
\begin{figure}[t]
\centering
\includegraphics[width=0.8\linewidth]{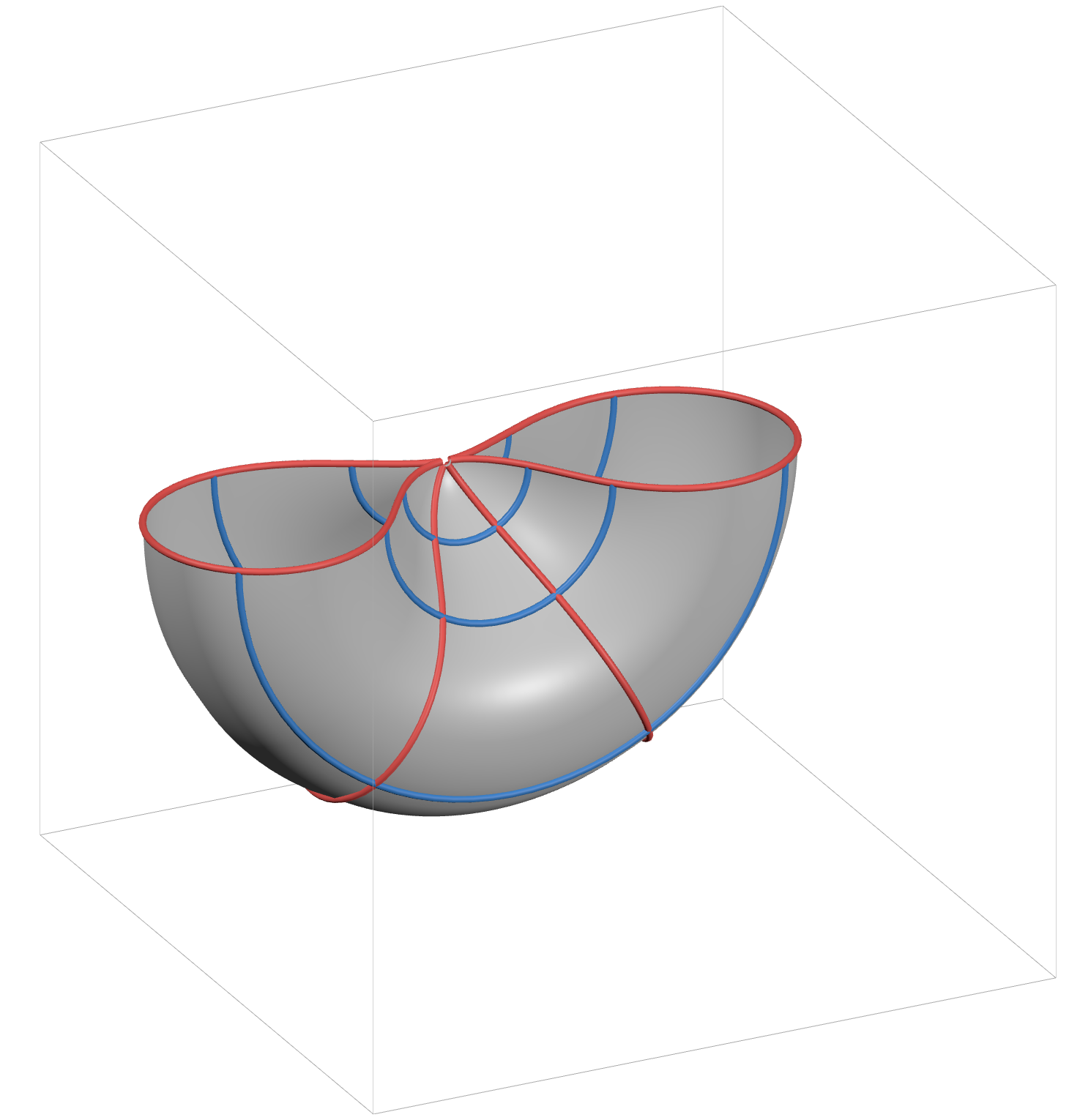}
\caption{A cut-away diagram of an inverted catenoid.
It is a surface of revolution of the teardrop-shaped red curve.
The point in the center is $\qi$, which is the point added in the compactification.}
\end{figure}

In fact this extended map is admissible, Definition~\ref{def:weakly conformal global}, as we shall now demonstrate.
At points of $\mathbb{C}^\ast$ we already know from the Möbius transformation that it is admissible, because it has a Kodaira representation on neighborhoods of these points (Example~\ref{eg:catenoid kodaira}).
So it remains to prove admissibility at the newly added points.
A direct approach would be to check that $\Tilde{N}$ belongs to $\sobolev{1,2}$ in a neighborhood of the new points, perhaps using the formula
\[
\Tilde{N} 
= \Tilde{\upsilon}^{-1}\qi\Tilde{\upsilon}
= (F\gamma + \delta)^{-1}N (F\gamma + \delta).
\]
If one pursues this approach, in checking that $d\Tilde{N} \in \banach{2}$ it suffices to show that the potentials $\B,\U,\V$ as well as $\psi\alpha\Tilde{\upsilon}^{-1}$ extend as $\banach{2}$ functions to the additional points.

But let us give an explanation more in keeping with the present theme and show that it has Kodaira data.
First we handle the extension to $z = 0$.
The potential $\Tilde{\V} = \V = dz\, \qk z^{-1}\Bar{z}(1 + |z|^2)^{-1}$ cannot be extended to $z=0$ continuously, but any extension is an $\banach{2}$ function on $\mathbb{C}$.
For $z \neq 0$, we have the Kodaira representation from the Möbius transformation
\begin{align*}
\Tilde{\upsilon}
&= \upsilon (F + \qi)
= -\qj (z - \qk)^{-1} (zF + z\qi), \\
\Tilde{\phi}
&= \upsilon (F \qi + 1)
= -\qj (z - \qk)^{-1} (zF\qi + z).
\end{align*}
The sections extend continuously with $\Tilde{\upsilon}(0) = - \qk$ and $\Tilde{\phi}(0) = -\qj$.
We argued in Example~\ref{eg:ln is V-holo} that $\ln z\Bar{z}$ is $\sobolev{1,p}$ for all $1 < p < 2$.
It follows that both sections are weakly differentiable with $\barpartial\Tilde{\upsilon} = \V \Tilde{\upsilon}$ and $\barpartial\Tilde{\phi} = \V \Tilde{\phi}$.
In particular, they are $\V$\=/holomorphic at $0$.
We can conclude that $\upsilon$ does not have a root at $z=0$ in the sense of Definition~\ref{order of roots} from Lemma~\ref{wente} and $|\upsilon(0)| \neq 0$.
Therefore Theorem~\ref{thm:kodaira} tells us that $\Tilde{F}$ is admissible.

The proof at $z = \infty$ is similar.
Rename $z$ to $z_1$ and use the coordinate $z_2 = z_1^{-1}$.
Instead of $zF$, it is convenient to work with 
\[
\Bar{z}_2 F(z_2^{-1})
= - \Bar{z}_2 \ln (z_2\Bar{z}_2)\qi + (z_2\Bar{z}_2 + 1) \qj,
\]
which has the value $\qj$ in the limit $z_2 \to 0$.
For $z_2 \neq 0$ we have
\begin{align*}
\Tilde{\upsilon}_1
&= -z_2^{-1} \qj (1 - \Bar{z}_2 \qk)^{-1} \Bar{z}_2(F + \qi), \\
\Tilde{\phi}_1
&= -z_2^{-1} \qj (1 - \Bar{z}_2 \qk)^{-1} \Bar{z}_2(F\qi + 1).
\end{align*}
In order to hold $\Tilde{\upsilon}$ bounded at $z_2 = 0$, we see that $f_{21}$ should be linear.
Let $E$ be the line bundle with holomorphic cocycle $f_{21} = z_2$. Note $\deg(E) = 1$.
Then
\begin{align*}
\Tilde{\upsilon}_2
&= f_{21}\Tilde{\upsilon}_1
= -\qj (1 - \Bar{z}_2 \qk)^{-1} \Bar{z}_2(F + \qi), \\
\Tilde{\phi}_2
&= f_{21}\Tilde{\phi}_1
= - \qj (1 - \Bar{z}_2 \qk)^{-1} \Bar{z}_2(F\qi + 1).
\end{align*}
Just as for the $z_1 = 0$ case, these sections extend continuously to $z_2 = 0$ with $\Tilde{\upsilon}_2(0) = 1$ and $\Tilde{\phi}_2(0) = \qi$.
These are $\V_2$\=/holomorphic for 
\[
d\Bar{z}_2 \V_2 
= z_2\, d\Bar{z}_1\V_1\, z_2^{-1}
= -d\Bar{z}_2\, \qk z_2^{-1}\Bar{z}_2 (1 + |z_2|^{2})^{-1}.
\]
This is square-integrable.
As it is represented by Kodaira data, the map $\Tilde{F}$ is admissible.
This example was the source of Example~\ref{eg:ln is V-holo}.
\end{example}

Möbius transformations and reflection account for all angle preserving transformations of $\qat$.
We can also try to understand how changing the Riemann surface $\X$ changes the representations.
If we have a holomorphic or anti-holomorphic map between Riemann surfaces, then the composition with an admissible map will also be admissible.
We deal with holomorphic maps in Example~\ref{eg:kodaira covering}.
For anti-holomorphic maps we instead investigate the conjugate Riemann surface $\Bar{\X}$ in Example~\ref{eg:kodaira conjugate}.

\begin{example}[Holomorphic coverings]
\label{eg:kodaira covering}
Suppose that $F : \X \to \qat$ is an admissible map with Kodaira representation $(E,\V,\upsilon,\phi)$.
Now suppose that we have a holomorphic map $f : \Tilde{\X} \to \X$.
Then $\Tilde{F} = F \circ f : \Tilde{\X} \to \qat$ is an admissible map.
From this point of view, we can see the admissible map $\Tilde{F}(z) = \alpha f(z) \beta$ from Example~\ref{eg:plane kodaira} as the composition of $f : \X \to \mathbb{C}$ and the embedded plane $F: \mathbb{C} \to \qat, z \mapsto \alpha z \beta$.

Let us describe its Kodaira data $(\Tilde{E},\Tilde{\V},\Tilde{\upsilon},\Tilde{\phi})$.
Essentially, one just pulls back the data along $f$.
The pullback bundle $\Tilde{E} = f^\ast E$ is a holomorphic line bundle on $\Tilde{\X}$.
Certainly $\Tilde{\upsilon} = \upsilon \circ f$ and $\Tilde{\phi} = \phi \circ f$ are sections of $\Tilde{E}_\qat$.
And $\Tilde{F} = (\upsilon^{-1}\phi) \circ f = \Tilde{\upsilon}^{-1}\Tilde{\phi}$.
So the point of difficulty is to understand the potential $\Tilde{\V}$ and why the sections are $\Tilde{\V}$\=/holomorphic.

By shrinking if necessary, we use a cover of $\Tilde{\X}$ by charts $w_l : \Tilde{\SO}_l \to \mathbb{C}$ and of $\X$ by $z_l : \SO_l \to \mathbb{C}$ such that $f[\Tilde{\SO}_l] \subset \SO_l$.
Define 
\[
\Tilde{\V}_l = \overline{\frac{\partial f}{\partial w_l}} \V_l \circ f.
\]
We see that this transforms as a $(0,1)$\=/potential of $f^\ast E$ (compare~\eqref{eq:potential transformation}), which has cocycle $f_{ml}\circ f$:
\begin{align*}
d\Bar{w}_m \Tilde{\V}_m
&= d\Bar{w}_m \overline{\frac{\partial f}{\partial w_m}} \V_m \circ f
= d\Bar{w}_l \overline{\frac{\partial f}{\partial w_l}} (f_{ml}\circ f)(\V_l \circ f)(f_{ml}\circ f)^{-1} \\
&= d\Bar{w}_l (f_{ml}\circ f)\Tilde{\V}_l(f_{ml}\circ f)^{-1}.
\end{align*}
Finally, it is now easy to see that the sections are $\Tilde{\V}$\=/holomorphic:
\[
\frac{\partial}{\partial \Bar{w}_l} \Tilde{\upsilon}_l
= \frac{\partial f}{\partial \Bar{w}_l}\frac{\partial \upsilon_l}{\partial z_l} \circ f + \frac{\partial \Bar{f}}{\partial \Bar{w}_l}\frac{\partial \upsilon_l}{\partial \Bar{z}_l} \circ f
= 0 + \overline{\frac{\partial f}{\partial w_l}} (\V_l \circ f) (\upsilon_l \circ f)
= \Tilde{\V}_l \,\Tilde{\upsilon}_l,
\]
since $f$ is holomorphic.
The same calculation holds for $\phi$.

The Weierstraß representation can be handled with the same ideas.
One section is simple: $\Tilde{\chi} = (\overline{\qj \Tilde{\upsilon}})^{-1} = \chi \circ f$ is a section of $f^\ast E^{-1}_\qat$.
However, $\Tilde{\psi}$ must be a section of $K_{\Tilde{\X}}\Tilde{E}$ and it is not generally the case that the pullback of the canonical bundle of $\X$ gives the canonical bundle of $\Tilde{\X}$.
This provides a concise reason of why we needed a corrective factor in the above definition of $\Tilde{\V}$.
The remaining formulas are forced by the chain rule:
\[
\Tilde{\psi}_l 
= \frac{\partial f}{\partial w_l} \psi_l \circ f, 
\qquad
\Tilde{\B}_l = 
\frac{\partial f}{\partial w_l} \B_l \circ f,
\qquad
\Tilde{\U}_l = 
\frac{\partial f}{\partial w_l} \U_l \circ f.
\qedhere
\]
\end{example}

\begin{example}[Conjugate Riemann Surface]
\label{eg:kodaira conjugate}
For any Riemann surface $\X$, there is another Riemann surface $\Bar{\X}$ which consists of the same underlying smooth manifold, but with opposite orientation, called the conjugate Riemann surface to $\X$.
For any chart $z_l : \SO_l \subset \X \to \Omega_l \subset \mathbb{C}$ of $\X$ the conjugate $\Bar{z}_l : \SO_l \subset \X \to \Bar{\Omega}_l \subset \mathbb{C}$ is a chart of $\Bar{\X}$.
Let $\rho : \Bar{\X} \to \X$ be the identity map, which is anti-holomorphic and hence conformal.
Given any admissible map $F : \X \to \qat$, we can consider $\Tilde{F} = F \circ \rho : \Bar{X} \to \qat$, which must also be admissible.
Because we have changed the orientation, the normals of $\Tilde{F}$ are negated compared to $F$.

A natural question is whether there is a simple description of its Kodaira data $(\Tilde{E},\Tilde{\V},\Tilde{\upsilon},\Tilde{\phi})$ in terms of the Kodaira data $(E,\V,\upsilon,\phi)$ of $F$.
As you may have inferred from the very fact we are raising this question, the answer is yes.
The second part of Example~\ref{eg:plane kodaira}, where $f$ is anti-holomorphic, serves as our template.
We have $\Tilde{E} = \Bar{E}$, the bundle whose cocycle is $\Bar{f}_{ml}$, where $f_{ml}$ is the cocycle of $E$.
This is a holomorphic line bundle on $\Bar{X}$ with holomorphic structure $\partial_{\Bar{E}}$.
Then we can define $\Tilde{\upsilon}_l = \qj \upsilon_l \circ \rho$ and $\Tilde{\phi}_l = \qj \phi_l \circ \rho$.
They are indeed sections of $\Tilde{E}$ since
\[
\Tilde{\upsilon}_m 
= \qj f_{ml} \upsilon_l \circ \rho
= \Bar{f}_{ml} \qj \upsilon_m \circ \rho
= \Bar{f}_{ml} \Tilde{\upsilon}_l.
\]
In their quotient the $\qj$ cancels out, i.e.\ $\Tilde{\upsilon}^{-1}\Tilde{\phi} = (\upsilon\circ\rho)^{-1}(\phi\circ\rho) = F \circ \rho = \Tilde{F}$.
It remains to see that these sections are $\Tilde{V}$\=/holomorphic for some potential.
We compute
\[
(\delbar{\Bar{E}})_l \Tilde{\upsilon}_l
= \qj (\partial_{\Bar{E}})_l (\upsilon_l \circ \rho)
= \qj (\barpartial_l \upsilon_l)\circ \rho
= \qj (\V_l \upsilon_l)\circ \rho
= - \qj (\V_l\circ \rho) \qj \Tilde{\upsilon}_l.
\]
Hence define $\Tilde{\V}_l = \qj (\V_l\circ \rho) \qj^{-1}$.
It only remains to check that this transforms according to~\eqref{eq:potential transformation} as an anti-commuting potential of $\Bar{E}_\qat$:
\[
\Tilde{\V}_m 
= - \qj \V_m \qj
= - \qj f_{ml} \V_l f_{ml}^{-1} \qj
= - \Bar{f}_{ml} \qj \V_l \qj \Bar{f}_{ml}^{-1}
= \Bar{f}_{ml} \Tilde{\V}_l \Bar{f}_{ml}^{-1}.
\]
Because $\qj$ is a constant, it does not produce any additional terms in the derivatives.
The Weierstraß data is $(\Tilde{E}, \qj(\U \circ \rho)\qj^{-1},\qj \chi\circ \rho, \qj \psi\circ\rho)$.
\end{example}

Our last example of this section is not really a transformation.
Rather it considers the situation that two admissible maps `overlap' with one another.
A typical example is Example~\ref{eg:kodaira covering} of a composition; clearly the image of the composition is a subset of the image of the original map.
The follow example should be considered as a converse to that example.

\begin{example}[Factorisation of admissible maps]
\label{eg:kodaira factorization}
We can generalize Example~\ref{eg:plane admissible}.
Let $F : \X \to \qat$ be an injective admissible map that is bi-Lipschitz continuous onto its image.
Suppose that $\Tilde{F} : \Tilde{\X} \to F[\X]$ is also admissible.
Then we have an induced map $f : \Tilde{\X} \to \X$, also Lipschitz.
Because the image of $\Tilde{F}$ is contained in the image of $F$, their tangent planes must coincide, and up to sign their normals agree.
Without loss of generality, we assume that in fact the normals are equal (if not, apply Example~\ref{eg:kodaira conjugate} to $\X$).

We claim that $f$ is a holomorphic map between Riemann surfaces.
This is a local claim, so we can assume that $\Tilde{\X}$ is contained in a single chart $w$ and likewise $\X$ is contained in a single chart $z$.
The definition of weakly conformal implies (a la Example~\ref{eg:plane admissible}) that
\[
N^+ \frac{\partial F}{\partial z} + N^- \frac{\partial F}{\partial \Bar{z}} 
= \qi \frac{\partial F}{\partial z}, \qquad
N^- \frac{\partial F}{\partial z} + N^+ \frac{\partial F}{\partial \Bar{z}} 
= - \qi \frac{\partial F}{\partial z}.
\]
Then 
\begin{align*}
N d\Tilde{F} 
&= N df \frac{\partial F}{\partial z} + N d\Bar{f} \frac{\partial F}{\partial \Bar{z}} 
= df \left[ N^+ \frac{\partial F}{\partial z} + N^- \frac{\partial F}{\partial \Bar{z}} \right] + d\Bar{f} \left[ N^- \frac{\partial F}{\partial z} + N^+ \frac{\partial F}{\partial \Bar{z}} \right] \\
&= df \qi \frac{\partial F}{\partial z} - d\Bar{f} \qi \frac{\partial F}{\partial \Bar{z}} 
= dw \left[ \frac{\partial f}{\partial w} \qi \frac{\partial F}{\partial z} - \frac{\partial \Bar{f}}{\partial w} \qi \frac{\partial F}{\partial \Bar{z}} \right] + d\Bar{w} \left[ \frac{\partial f}{\partial \Bar{w}} \qi \frac{\partial F}{\partial z} - \frac{\partial \Bar{f}}{\partial \Bar{w}} \qi \frac{\partial F}{\partial \Bar{z}} \right].
\end{align*}
On the other hand
\begin{align*}
-\ast d\Tilde{F}
= dw \qi \left[ \frac{\partial f}{\partial w} \frac{\partial F}{\partial z} + \frac{\partial \Bar{f}}{\partial w} \frac{\partial F}{\partial \Bar{z}} \right] - d\Bar{w} \qi \left[ \frac{\partial f}{\partial \Bar{w}} \frac{\partial F}{\partial z} + \frac{\partial \Bar{f}}{\partial \Bar{w}} \frac{\partial F}{\partial \Bar{z}} \right].
\end{align*}
These two expressions are equal ($\Tilde{F}$ is weakly conformal) if and only if 
\[
\frac{\partial \Bar{f}}{\partial w} \frac{\partial F}{\partial \Bar{z}} = 0 
\quad\text{and}\quad
\frac{\partial f}{\partial \Bar{w}} \frac{\partial F}{\partial z} = 0.
\]
Since $F$ cannot be constant, at least one of the derivatives of $F$ must be non-zero.
This forces $\barpartial f = 0$, proving the claim that $f$ is holomorphic.

It follows that $\Tilde{F}$ has the same regularity as $F$.
In particular, admissible maps into the round sphere or plane must be analytic functions.
\end{example}

\section*{Maps into $\Imag \qat = \mathbb{R}^3$}

An important advantage of the Weierstraß representation compared to the Kodaira representation is that one can single out immersions into $\mathbb{R}^3$.
If an immersion $F$ is valued in $\Imag\qat$ then so is its derivative $dF$ (and the converse is true up to a constant of integration), so consider immersions with $\overline{dF} = - dF$.
At any point we ask ``Given $\chi$, what are the possible values of $\psi$?'' 
The key is to write $\psi$ as a multiple of the non-zero $\chi$:
\begin{align*}
\left\{\psi\in\qat\mid\overline{\Bar{\chi}\qj dz\psi} =-\Bar{\chi}\qj dz\psi\right\} 
&= \left\{ \psi = g\chi \mid g \in \qat,\, \Bar{g}\qj dz = \qj dz g \right\} 
= \left\{\psi = g\chi\mid g\in\mathbb{C}\right\}.
\end{align*}
Hence an admissible map $F:\X\to\qat$ with a Weierstraß representation~\eqref{eq:weierstrass} takes values in the imaginary quaternions if and only if the local representatives $\chi_l$ and $\psi_l$ of the sections obey $\psi_l = g_l\chi_l$ for a measurable function $g_l:\SO_l\to\mathbb{C}$.
Assuming $\chi_l$ is $\qk\Bar{\U}_l\qk$\=/holomorphic, the function $g_l\chi_l$ obeys
\[
(\barpartial_l -\U_l)(g_l\chi_l)
=\left(\barpartial_l g_l + g_l\qk\Bar{\U}_l\qk - \U_l g_l \right)\chi_l,
\]
which is zero (that is, $\psi_l$ is $\U_l$\=/holomorphic) if and only if 
\begin{equation}
\label{eq:imag dF properties of g}
\barpartial_l g_l =2\U_l^+g_l
\quad\text{and}\quad 
g_l\qk\Bar{\U}_l^-\qk = \U_l^-g_l.
\end{equation}
Both of these conditions can be simplified.
This second condition may also be written as $\U_l g_l\qk =- g_l\qk\Bar{\U}_l= \overline{\U_l g_l \qk}$, showing that the phase of the potential $\U_l$ is determined by $g_l$.
For the first condition, by using Corollary~\ref{gauge holomorphic structure} we may replace 
\begin{itemize}
\item 
the holomorphic line bundle $E$ by another line bundle $E'$ of the same degree, 
\item the potential $\U$ by a potential $\U'$ with $(\U')^+ = 0$, such that $\Q{E',\U'}$ is isomorphic to $\Q{E,\U}$ and $|\U'|=|\U^-|$, 
\item and the local representatives of the sections $\chi$ and $\psi$ by local representative $\chi'_l$ and $\psi'_l$ of a global sections of $\Q{E,\qk\Bar{\U}'\qk}$ and $\Q{KE',\U'}$, respectively.
\end{itemize}
Changing the line bundle $E'$ does not change $dF$, and so again $\psi'_l=g'_l\chi'_l$ for some complex functions $g'_l$.
However the first condition of Equation~\eqref{eq:imag dF properties of g} now reads $\barpartial_l g'_l = 0$ since $\U^+ = 0$.
Hence the $g'_l$ are holomorphic.
Having carried out this change of line bundle, from now on we shall omit the primes.
We can understand the local holomorphic functions $g_l$ in a global way.
Since they transform the section $\chi$ of $E^{-1}_\qat$ into the sections $\psi$ of $KE_\qat$ they define a holomorphic section $g$ of $\Hom_\qat^+(E^{-1}_\qat, KE_\qat) = \Hom_\mathbb{C}(E^{-1}, KE) = KE^2$.

If $\chi$ has no roots, then the roots of the $g$ are the roots of $\psi$.
In particular, on compact $\X$ the line bundle $E$ and the cocycle $g$ are uniquely determined by the root divisor $(dF)$ of $dF$.
If $g$ has no roots then $E$ is a spin bundle $E^2 = K^{-1} = T\X$.
More generally, on the simply connected $\SO_l$ there exist holomorphic functions $f_l:\SO_l\to\mathbb{C}$ with $g_l=f_l^2$ if and only if all roots of $\psi$ are of even order.
In this case, after moving half the roots from $\psi$ to $\chi$ both sections coincide $\Tilde{\chi}_l=f_l\chi_l=f_l^{-1}\psi_l=\Tilde{\psi}_l$ and are holomorphic sections of $\Q{\Tilde{E},\Tilde{\U}}$.
Here $\Tilde{E}$ is a spin bundle $\Tilde{E}^2=T\X$ and the potential $\Tilde{\U}=f^{-1}\U f\in\pot{\Tilde{E}^-}$ obeys
\[
(f\Tilde{\U}f^{-1})\sd=f^{-2}f\Tilde{\U}f^{-1}f^2\quad\Longleftrightarrow\quad\qk\Tilde{\U}\qk=\Tilde{\U}.
\]
This means that the local representatives of $\Tilde{\U}$ are $\mathbb{R}\qk$\=/valued functions.
To summarize
\begin{lemma}
\label{lem:admissible imaginary}
\index{Weierstraß representation!In $\mathbb{R}^3$}
If $F$ is an admissible map on a Riemann surface $\X$ into $\Imag\qat$, then there is a line bundle $E$ such that the corresponding quadruple $(E,\U,\chi,\psi)$ in Theorem~\ref{thm:weierstrass} obeys
\begin{align}\label{eq:imaginary 1}
\psi = g\chi, \quad
\U = g\U\sd g^{-1},
\quad\text{for some holomorphic section $g$ of $KE^2$}.
\end{align}
Additionally the following conditions can be fulfilled if and only if all roots of $dF$ are of even order:
\begin{align}\label{eq:imaginary 2}
\hspace{25mm}
g&=1,&E^2&=T\X,&\psi&=\chi, &\U&=\U\sd.
\hspace{24mm}\qed
\end{align}
\end{lemma}

Let us examine how Möbius transformations interact with Lemma~\ref{lem:admissible imaginary}.
The transformation~\eqref{eq:quaternionic moebius 2} preserves the imaginary quaternions if and only if
\begin{align*}
(F\gamma +\delta)^{-1}(F\alpha+\beta)=-(\Bar{\alpha}F-\Bar{\beta})(\Bar{\gamma} F-\Bar{\delta})^{-1}
\quad\text{whenever }\Bar{F}=-F.
\end{align*}
This implies 
\[
F(\alpha\Bar{\gamma}+\gamma\Bar{\alpha})F + F(-\alpha\Bar{\delta} + \delta\Bar{\alpha}) + (\beta\Bar{\gamma} - \gamma\Bar{\beta})F - (\beta\Bar{\delta} + \delta\Bar{\beta}) = 0
\]
for all imaginary $F$.
We claim that each of the brackets must vanish separately.
Taking $F = 0$ shows that the last bracket must vanish.
Then comparing $F$ and $-F$ shows that the first bracket must vanish.
This leaves an equation of the form $AF + FB = 0$.
If neither $A$ nor $B$ vanishes, they must have the same length, and then we know from Definition~\ref{def:left and right normals} that this equation has a $2$\=/dimensional space of solutions, by interpreting $A$ and $-B$ as left and right normals.
But we need it to hold for all $F \in \Imag \qat$.
Hence at least one of them vanishes, in which case the other does too.
This proves the claim.
Conversely, if each of the four brackets vanishes then certainly the equation holds for all $F$.
Therefore a Möbius transformation preserves the imaginary quaternions if and only if
\begin{align}
\label{eq:imaginary preserved}
\alpha\Bar{\gamma} = -\gamma\Bar{\alpha}, &&
\alpha\Bar{\delta} = \delta\Bar{\alpha}, &&
\beta\Bar{\gamma} = \gamma\Bar{\beta}, &&
\beta\Bar{\delta} = - \delta\Bar{\beta}.
\end{align}

We saw in the argument preceding Lemma~\ref{lem:admissible imaginary} how to modify the Weierstraß data into the form of Equation~\eqref{eq:imaginary 1}.
It is not automatic therefore that the transformed Weierstraß data $(E,\Tilde{U},\Tilde{\chi},\Tilde{\psi})$ will also have the special form.
But in fact it does, as we will now verify.
We compute
\begin{equation}
\label{eq:mobius weierstrass imag cond}
\begin{aligned}
\Tilde{\psi}
&= \psi(\alpha - \gamma\Tilde{F})
= g\chi(\alpha + \gamma\overline{\Tilde{F}})
= g\chi(\alpha\overline{(F\gamma+\delta)} + \gamma\overline{(F\alpha + \beta)})\overline{(F\gamma+\delta)}^{-1} \\
&= g\chi(\alpha\Bar{\delta} + \gamma\Bar{\beta})\overline{(F\gamma+\delta)}^{-1}
= (\alpha\Bar{\delta} + \gamma\Bar{\beta})g\Tilde{\chi}
=: \Tilde{g}\Tilde{\chi},
\end{aligned}
\end{equation}
using that $\alpha\Bar{\delta} + \gamma\Bar{\beta}$ is real.
In particular $\Tilde{g}$ is holomorphic section of $KE^2$.
Note, that rescaling the $\alpha,\beta,\gamma,\delta$ by a real number does not affect the Möbius transformation.
In particular, we can normalize $\alpha\Bar{\delta} + \gamma\Bar{\beta} = 1$, so that Equation~\eqref{eq:imaginary 2} is preserved.
In Example~\ref{eg:sphere weierstrass}, we had sections with $\psi = \chi$ but after the Möbius transformation $\Tilde{\psi} = 2 \Tilde{\chi}$ in accordance with Equation~\eqref{eq:mobius weierstrass imag cond}.
If we had wanted to obtain a Weierstraß representation for which the sections remained equal, then we should have normalized the coefficients of the Möbius transformation: $\alpha = \delta = \qi/\sqrt{2}$ and $\beta = \gamma = 1/\sqrt{2}$.
This amounts to redistributing the factor of $2$ between the two sections.

\chapter{Willmore Energy and Isothermic Maps}
\label{chapter:isothermic}

In this chapter we apply the Kodaira and Weierstraß representations to extend some aspects of classical differential geometry to admissible maps.
In particular, we produce formulas for the mean curvature~\eqref{eq:mean curvature} and the Willmore energy~\eqref{eq:willmore energy weierstrass} in terms of the Weierstraß data.
Then we investigate the class of isothermic surfaces, giving a classical definition in terms of curvature lines, and characterizations in terms of dual maps and holomorphic quadratic differentials.
An important notion for Part~III, which we introduce here, is that of a strongly isothermic surface.
In Theorem~\ref{thm:isothermic local characterization} we give a proof purely in the language of quaternionic analysis that a totally umbilic surface is contained in a plane or round sphere.

We begin with the normals.
The Kodaira representation and Theorem~\ref{thm:kodaira normal} gave us the formula $N = \upsilon^{-1}\qi\upsilon$ for the left normal.
The Weierstraß representation gives us formulas for both normals.
As we did in Example~\ref{eg:catenoid}, from the local form $dF = \Bar{\chi} \qj dz \psi$ and Equations~\eqref{eq:left normal} and~\eqref{eq:right normal} that $-\ast dF = N dF = dF R$ we can read off the left and right normals:
\begin{align}
\label{eq:weierstrass-NR}
\index{Normal}
N 
&= \Bar{\chi}(-\qi)\Bar{\chi}^{-1} 
= \overline{\chi^{-1}\qi\chi}
= - \chi^{-1}\qi\chi,
&
R &= \psi^{-1}\qi\psi,
\end{align}
using that $\Bar{N} = -N$ to remove the conjugation from the first formula.
In the second formula, $\psi$ may have roots.
However, we can change the line bundle $E$ to create a section $\psi' = g\psi$ of $KE'_\qat$ without roots, as in the proof of Theorem~\ref{thm:weierstrass}.
Since $g$ is complex, $\psi'^{-1}\qi\psi' = \psi^{-1}\qi\psi$.
Thus these formulas for the normals are still well-defined at roots.

\begin{example}
\label{eg: - bar F}
We observed following Definition~\ref{def:left and right normals} that for maps into $\Imag\qat$ the left and right normals are the negatives of one another.
We can see this in the formulas, because such an admissible map has a Weierstraß representation of the form in Lemma~\ref{lem:admissible imaginary}, from which it follows that
\[
R 
= \psi^{-1}\qi\psi
= (g \chi)^{-1}\qi g \chi
= \chi^{-1}\qi \chi
= -N.
\]

The effect of a Möbius transformation on the normals follows from Equations~\eqref{eq:mobius chi} and~\eqref{eq:mobius psi} but is in general rather complicated.
For an affine transformation we have $\Tilde{\chi} = \chi\Bar{\delta}^{-1}$ and $\Tilde{\psi} = \psi\alpha$ so
\begin{align*}
\Tilde{N}
&= - (\chi\Bar{\delta}^{-1})^{-1} \qi (\chi\Bar{\delta}^{-1})
= \Bar{\delta} N \Bar{\delta}^{-1}, &
\Tilde{R}
&= (\psi\alpha)^{-1}\qi (\psi\alpha)
= \alpha^{-1} R \alpha,
\end{align*}
which are both rotations of $\mathbb{S}^2$.

Another example is that of $-\Bar{F}$.
From its Weierstraß representation $d(-\Bar{F}) = \lp \psi, \chi \rp$ we have
\begin{align*}
N_{-\Bar{F}} 
&= - \psi^{-1} \qi \psi = - R_F, &
R_{-\Bar{F}} 
&= \chi^{-1} \qi \chi = - N_F.
\end{align*}
In words, the left normal of $-\Bar{F}$ is equal to the negative of right normal of $F$ and the right normal of $-\Bar{F}$ is equal to the negative of the left normal of $F$.
\end{example}

In the reverse direction, the quaternionic Weierstraß representation yields the following solution to the problem of whether two given maps $N,R\in\sobolev{1,2}(\X,\qat)$ with $N^2=-1=R^2$ are the left and right normals of an admissible map on the universal covering $\Tilde{X}$.
Apply Theorem~\ref{thm:kodaira normal} to $-N$ and $R$ to get the corresponding triples $(E^{-1},\U\sd,\chi)$ and $(E',\U',\psi')$ respectively, where neither section $\chi$ nor $\psi'$ has roots.
Because $N$ and $R$ are arbitrary, there is not necessarily any relationship between $E^{-1}$ and $E'$.
In particular, there is no reason that they should be paired bundles.
There exists an admissible map on the universal covering $\Tilde{\X}$ with left normal $N$ and right normal $R$ if and only if the holomorphic $\mathbb{C}$\=/line bundle $KE \otimes E'^{-1}$ has a non-trivial holomorphic section $g$ and the two potentials $g\U'g^{-1}$ and $\U$ in $\pot{KE}^-$ agree.
In this case $\psi = g\psi'$ is a $\U$\=/holomorphic section of $KE_\qat$ and $(\U, \chi, \psi)$ is the Weierstra{\ss} triple of an admissible map on $\Tilde{\X}$.
This means that $\lp\chi,\psi\rp$ defines a $\qat$\=/valued closed $1$\=/form which is the derivative of an admissible map on $\Tilde{\X}$ with left normal $N$ and right normal $R$.
The roots of $\lp\chi,\psi\rp = \lp\chi,g\psi'\rp$ are the same as the roots of the holomorphic section $g$ of $KE \otimes E'^{-1}$.

It is well-understood under which conditions the Gauss map determines an immersion into $\mathbb{R}^n$.
In particular, Theorem~2.5 and the paragraph that precedes it in~\cite{Hoffman1983} tells us that non-minimal surfaces are uniquely determined by their Gauss map up to similarity and minimal surfaces with the same Gauss map are in one-to-one correspondence with holomorphic functions.
We can easily reproduce these results using quaternionic analysis.
\begin{example}[Minimal Surfaces]
\index{Minimal surfaces}
The Gauss map is the map $p \mapsto T_p F[\X] \in \mathrm{Gr}(2,n)$, and in our setting this is exactly the information of the left and right normals.
If two admissible maps $F,F' : \X \to \qat$ with Weierstraß data $(E,\U,\chi,\psi)$ and $(E',\U',\chi',\psi')$ have the same left and right normals, then $\qi \chi'\chi^{-1} = \chi'\chi^{-1} \qi$ shows that $\chi' = g \chi$ for a holomorphic section $g$ of $E'\otimes E^{-1}$.
Since neither $\chi$ nor $\chi'$ vanish, neither does $g$ and therefore the two bundles $E$ and $E'$ are isomorphic.
We may assume that $E' = E$, so $g$ is a function $\X \to \mathbb{C}^\ast$.
As stated in Theorem~\ref{thm:weierstrass}, there is some freedom in the Weierstraß representation.
We can use this freedom to choose $\chi' = \chi$, which also implies $\U' = \U$.

Up to this point, we have only used the equality of the left normals of $F$ and $F'$.
The same reasoning can be applied to the right normal, giving a $\mathbb{C}$\=/valued holomorphic function $h$ that obeys $\psi' = h\psi$ and $\U = \U' = h \U h^{-1}$ on the complement of the roots of $\psi$ and $\psi'$.
If $\U$ is not zero, because $\U \in \pot{KE}^-$ it must be that $h$ is real valued.
But the only real valued holomorphic functions are constants.
This allows us to conclude that $dF' = \lp \chi, h \psi \rp = dF\, h$.
In other words, $F$ and $F'$ agree up to scaling and translation.
However, if $\U$ vanishes then the admissible maps with the same normals as $F$ are of the form $dF' = \lp \chi, h \psi\rp$ for some meromorphic function $h$ with $\ord_x h + \ord_x \psi \geq 0$.
Because these have the same line bundle $E$, Theorem~\ref{thm:weierstrass} tells us that they are distinct admissible maps.
That $\U$ vanishing corresponds to minimal surfaces is Equation~\eqref{eq:mean curvature} below.

For $h \in \mathbb{S}^1$ these new admissible maps are the associated family of the minimal surface.
For example, consider the catenoid.
If we insert $h \in \mathbb{S}^1$ then $dF' = \lp \chi, h \psi\rp$ is no longer an exact differential:
\[
\int_{|z| = 1} dF'
= \int_0^{2\pi} (\overline{1 + e^{-\qi\theta}\qk})\,\qj d\theta\qi e^{\qi\theta} h (1 + e^{-\qi\theta}\qk)
= - 4\pi\qi \Imag h.
\]
These associated minimal surfaces belong to the catenoid-helicoid family, with the helicoid corresponding to $h = \pm \qi$.
In general we can take the Weierstraß data of the new surface to be $(E',\U', \chi', \psi') = (E, 0, \chi, h \psi)$.
Then $\upsilon' = \upsilon$ and the other potentials are also unchanged $\B' = \B$, $\V' = \V$.
The only difficulty is to find a $\V$\=/holomorphic section $\phi'$ of $E$ with $(\partial_l + B_l + 0)\phi' = \psi' = h \psi$; the obvious guess $h \phi$ is instead holomorphic with respect to $h \V h^{-1} = h^2 \V$.
However if we can find a solution $\phi'$ for $h'=\qi$, then a linear combination gives the general solution $\phi'' = \phi\cos \theta + \phi'\sin\theta$ for $h'' = \cos \theta + \qi \sin\theta \in \mathbb{S}^1$:
\begin{gather*}
(\partial_l + B_l + 0)\phi''
= \psi\cos \theta + \qi\psi \sin\theta = h'' \psi = \psi'', \\
F'' 
= (\upsilon'')^{-1}\phi'' = F \cos\theta + F' \sin\theta. 
\qedhere
\end{gather*}
\end{example}

The relationship between the normals and the Weierstraß representation opens the path to express the classical geometric information of an immersed surface in the language of quaternionic analysis.
We should note that the following calculations only make sense for immersions, or at least away from the branch points.
As in Remark~\ref{rem:weierstrass} branch points are the roots of $dF$ or equivalently $\psi$.
At these points the conformal factor $|\chi|\cdot|\psi|$ is zero, and likewise the below formulas for the tangent and perpendicular planes contain only the zero vector.
This is to be expected: at these non-immersive points there is not a correspondence between $T\X$ and $TF[\X]$.
With this in mind, from the formula for the left normal $N = \upsilon_l^{-1}\qi\upsilon_l$ it follows
\begin{align}
\label{eq:dN potentials}
dN 
&= - \upsilon_l^{-1} d\upsilon_l N + \upsilon_l^{-1} \qi d\upsilon_l
= 2 \upsilon_l^{-1}\qi ( - dz_l \U_l + d\Bar{z}_l \V_l) \upsilon_l.
\end{align}
Likewise
\begin{align}
\label{eq:dR potentials}
dR
= - \psi_l^{-1} d\psi_l R + \psi_l^{-1} \qi d\psi_l
= 2 \psi_l^{-1} \qi ( dz_l \V'_l + d\Bar{z}_l \U_l) \psi_l,
\end{align}
where $\V'_l = ((\partial_l\psi_l)\psi_l^{-1})^-$.
These formulas give a geometric interpretation of the potentials $\U$ and $\V$; in essence, the potentials determine the derivative of the normals.

Next we describe the \emph{second fundamental form} in terms of the left and right normals.
As we commented upon in Chapter~\ref{chapter:prelim} and despite their name, the left and right normals are not perpendicular to the immersion.
Instead they describe the tangent plane $\{ \Bar{\chi}_l \qj \alpha \psi_l \mid \alpha \in \mathbb{C}\}$ and the perpendicular plane $\{ \Bar{\chi}_l \alpha \psi_l \mid \alpha \in \mathbb{C}\}$.
In any chart $z_l$ let us consider the Weingarten operator $A^{\alpha} : T\X \to T\X$ associated to the perpendicular vector field $\Xi^\alpha : z_l \mapsto \Bar{\chi}_l(z) \alpha(z) \psi_l(z)$ for a complex-valued function $\alpha : \SO_l \to \mathbb{C}$.
Following the definition in~\cite[Equation~(2.3)]{Pal}, it is the negative of the tangential component of the derivative of the perpendicular vector field:
\begin{align*}
A^{\alpha}
&= -\frac{1}{2} \Big[d(\Bar{\chi}_l \alpha \psi_l) - N d(\Bar{\chi}_l \alpha \psi_l)R \Big] \\
&= - \frac{1}{2} \Big[ \Big( \overline{d\chi}_l - N \overline{d\chi}_l \qi \Big) \qj\alpha\psi_l 
+ \Bar{\chi}_l \Big( d\alpha + \qi d\alpha \qi \Big) \psi_l
+ \Bar{\chi}_l \alpha \Big( d\psi_l + \qi d\psi_l R \Big) \Big]\\
&= \phantom{-}\frac{1}{2} \Big[ \overline{dN}\, \Bar{\chi}_l \qi\alpha\psi_l + 0 + \Bar{\chi}_l \qi \alpha \psi_l\, dR\Big] 
= \frac{1}{2} \Big[ \overline{dN}\, \Xi^{\qi\alpha} + \Xi^{\qi\alpha}\, dR\Big].
\end{align*}
In particular, it does not depend on the derivative of $\alpha$.
We might further simplify by the insertion of Equations~\eqref{eq:dN potentials} and~\eqref{eq:dR potentials}:
\begin{align}
\label{eq:weingarten operator}
A^{\alpha}
&= \Bar{\chi}_l \qj \big[ dz_l \,2\Real (\qj \alpha \U_l) 
+ d\Bar{z}_l ( \qj \alpha \V'_l - \V_l\qj\alpha ) \big]\psi_l.
\end{align}
This expression is a little muddled: as a $\qat$\=/valued differential form it acts on $T\X$ to produce something in $\qat$.
But notice that the expression in the square bracket is purely $\mathbb{C}$\=/valued, so overall we indeed have a tangent vector in $TF[\X] \cong T\X$.
It is perhaps more familiar to think of the Weingarten operator in real terms, in which case with respect to the real coordinates $z_l = x_l + \qi y_l$ it has the matrix
\begin{align}
\label{eq:weingarten matrix}
\begin{pmatrix}
2\Real (\qj \alpha \U_l) + \Real((\V'_l - \V_l)\qj\alpha) & - \Imag((\V'_l - \V_l)\qj\alpha) \\
- \Imag((\V'_l - \V_l)\qj\alpha) & 2\Real (\qj \alpha \U_l) - \Real((\V'_l - \V_l)\qj\alpha)
\end{pmatrix}.
\end{align}

The second fundamental form of two tangent vectors is the projection of $\nabla_X Y$ to the normal bundle.
Its inner product with $\Xi^\alpha$ equal to the inner product of $X$ with $A^\alpha(Y)$, as written below.
Note that the perpendicular vector field $\Xi^\alpha$ need not be unit length.
\begin{align*}
\index{Second fundamental form}
\Two^\alpha(X,Y)
&= \left\langle dF(X), A^\alpha(Y) \right\rangle_\qat
= \Real \Big[\overline{\Bar{\chi_l} \qj dz_l(X)\psi_l} A^\alpha(Y) \Big] \\
&= |\chi_l|^2|\psi_l|^2 \Real \Big[d\Bar{z}_l \otimes dz_l \,2\Real (\qj \alpha \U_l) + d\Bar{z}_l \otimes d\Bar{z}_l ( \qj \alpha \V'_l - \V_l\qj\alpha ) \Big].
\end{align*}
These values completely determine the second fundamental form.
We can take a local orthonormal frame for the normal bundle, such as $\alpha_1 = |\chi|^{-1}|\psi|^{-1}$ and $\alpha_2 = \qi \alpha_1$.
Then the second fundamental form is $\Two = \Two^{\alpha_1}\Xi^{\alpha_1} + \Two^{\alpha_2}\Xi^{\alpha_2}$.

The mean curvature vector is the half the trace of the second fundamental form.
We can give a formula for it by evaluating the second fundamental form on an orthonormal basis of the tangent space.
Such a basis is $\{X_1,X_2\} = \{|\chi|^{-1}|\psi|^{-1}\partial_x,|\chi|^{-1}|\psi|^{-1}\partial_y\}$.
Since $\Two^\alpha(X_1,X_1) + \Two^\alpha(X_2,X_2)
= 4\Real (\qj \alpha \U_l)$, we have
\begin{align}
\label{eq:mean curvature}
\index{Mean curvature}
H 
&= \frac{1}{2} \Big[ 4\Real (\qj \alpha_1 \U_l)\Xi^{\alpha_1} + 4\Real (\qj \alpha_2 \U_l)\Xi^{\alpha_2} \Big]
= 2 \chi_l^{-1} \, \qj \U_l\, \Bar{\psi}_l^{-1}.
\end{align}
This formula contains only the Weierstraß data.
Similar to Equation~\eqref{eq:weierstrass-NR}, the transformation behavior of the sections and potential ensure that $H$ is independent of the trivialization $z_l$.
We see in particular that $\U = 0$ corresponds to minimal surfaces.
The following formulas, which are~\cite[Proposition~8]{BFLPP} when the sign convention of $\ast$, $N$ and $R$ are accounted for, will also be useful later:
\begin{align}\label{eq:mean curvature 2}
2\Bar{H}dF&=\ast dR+RdR,&2dF\Bar{H}&=\ast dN-NdN.
\end{align}

\begin{example}[Round Sphere]
\index{Round sphere!Mean curvature}
Since the plane and catenoid have vanishing mean curvature, the non-trivial example to consider is the round sphere.
Using the Weierstraß representation from Example~\ref{eg:sphere weierstrass} we have the inward pointing unit normal of the sphere
\begin{align*}
H 
&= 2 (\qj z + \qi) \qj (|z|^2 + 1)^{-1}\qk \tfrac{1}{2}(\overline{\qj z + \qi})
= - (z + \qk)^{-1} \qi (z + \qk) 
= -F.
\qedhere
\end{align*}
\end{example}

From the mean curvature we can derive a useful formula for the Willmore energy.
The \emph{Willmore energy}\index{Willmore energy} is the integral of the square of the mean curvature with respect to the induced area measure of the surface:
\begin{equation}
\label{eq:willmore energy weierstrass}
\willmore(F)
= \int_\X |H|^2 |dF|^2
= \int_\X 4 |\chi_l|^{-2} |\U_l|^2 |\psi_l|^{-2} \cdot |\chi_l|^2 |dz|^2 |\psi_l|^2 
= 4\|\U\|_2^2.
\end{equation}
This formula is well-known~\cite{Ta1,Ta2,PP}.
For compact Riemann surfaces $\X$ we can also compute the Willmore energy in terms of the Kodaira representation $(E, \V, \upsilon, \phi)$ for $\V \in \pot{E}^-$, since Equation~\ref{eq:hopf norms} implies
\begin{equation}
\label{eq:willmore energy kodaira}
\willmore(F) = 4\pi\deg(E) + 4\|\V\|_2^2. 
\end{equation}
Because the line bundle and Kodaira potential are invariant, Willmore energy is a Möbius invariant for compact Riemann surfaces.

\begin{example}
\label{eg:willmore energy}
\index{Willmore energy}
\index{Plane!Willmore energy}
\index{Catenoid!Willmore energy}
\index{Round sphere!Willmore energy}
Let us compute the Willmore energy for all our examples thus far.
First consider the plane and the catenoid.
From Examples~\ref{eg:plane weierstrass} and~\ref{eg:catenoid weierstrass} we know that $\U \equiv 0$, so by Equation~\eqref{eq:willmore energy weierstrass} their Willmore energies are both zero.

On the other hand, Equation~\eqref{eq:willmore energy kodaira} is more useful for the example of the round sphere.
We saw in Example~\ref{eg:sphere kodaira} that it has $\V \equiv 0$ and $\deg(E) = 1$.
Hence its Willmore energy is $4\pi$.

There is a slight generalization of this case that is still interesting.
Consider an admissible map from a compact Riemann surface $\X$ to the round sphere $\mathbb{S}^2$.
Due to Example~\ref{eg:kodaira factorization} we know that this factors through a branched covering $f : \X \to \mathbb{P}^1$, its Kodaira potential is zero, and the line bundle is $f^\ast E$ for $E$ the line bundle of the Kodaira representation of the round sphere.
The Willmore energy is $4\pi\deg(f)$.

We can also compute the Willmore energy of the compactified inverted catenoid.
We will use~\eqref{eq:willmore energy kodaira}, in particular the degree of the underlying bundle is $1$ and it has the same Kodaira potential as the catenoid $\V = dz\, \qk z^{-1}\Bar{z}(1 + |z|^2)^{-1}$.
The value of $\|\V\|_2$ is not changed by ignoring the point $z = \infty$.
Hence we calculate
\begin{align*}
\|\V\|_2^2
&= \int_{\mathbb{C}} \Big|\qk z^{-1}\Bar{z}(1 + |z|^2)^{-1}\Big|^2\,\dmu(z)
= \int_{\mathbb{C}} (1 + r^2)^{-2}\,r\, dr\, d\theta
= \pi.
\end{align*}
In sum then $\willmore(\Tilde{F}) = 4\pi\cdot 1 + 4\cdot \pi = 8\pi$.

The Willmore energy is also known for certain classes of surfaces.
For example, Bryant has shown in~\cite[\S5]{Bryant1984} that immersions of the sphere into $\mathbb{S}^3$ that are critical points of the Willmore energy function have energies that lie in $4\pi\mathbb{N}$.
We examine the critical points of the Willmore energy in Chapters~\ref{chapter:constrained 1} and~\ref{chapter:constrained weierstrass}.
\end{example}


Isothermic surfaces are a well-studied class of surfaces.
They include surfaces of constant mean curvature as well as surfaces of revolution.
Moreover isothermic surfaces which can be conceived of as integrable systems~\cite{Burstall1997a,Cieslinski1995}.
Classically, a surface in $\mathbb{R}^3$ is called isothermic if it has a conformal parameterization by lines of curvature.
This means that, away from the umbilic points (points where the principal curvatures are equal), every point has a coordinate neighborhood $z_l$ in which the Weingarten operator is diagonalized.
In higher codimension $n-2$ the Weingarten operator and the lines of curvature depend on the choice of perpendicular vector field $\Xi$, and so to be isothermic means that there is a coordinate that simultaneously diagonalizes all Weingarten operators.
This is not as impossible as it sounds: because the Weingarten operator is linear in $\Xi$, it amounts to simultaneously diagonalizing the $n-2$ Weingarten operators corresponding to a frame of the normal bundle.

For surfaces in $\qat$ we already have a formula for the Weingarten operator, Equation~\eqref{eq:weingarten matrix}.
We see that $A^\alpha$ is diagonal if and only $\Imag ((\V'_l - \V_l)\qj\alpha) = 0$.
We can take a frame of the normal bundle $\Xi^1, \Xi^\qi$, i.e.\ $\alpha = 1, \qi$.
Then $A^1$ and $A^\qi$ are simultaneously diagonal in this coordinate if and only if $\V'_l = \V_l$, where $\V'_l = ((\partial_l\psi_l)\psi_l^{-1})^-$.
This condition is well-defined even if the potentials are only $\banach{2}$.

\begin{definition}
\label{def:isothermic}
\index{Isothermicity}
An admissible map $F : \X \to \qat$ is called \emph{isothermic} if it has a cover of coordinate neighborhoods $z_l : \SO_l \to \mathbb{C}$ such that $\V'_l = \V_l$, where $d\Bar{z}_l\V_l$ is the Kodaira potential in local coordinates and similarly $\V'_l = ((\partial_l\psi_l)\psi_l^{-1})^-$.
\end{definition}

We should not (yet) try to interpret this equality globally, since $\V \in \pot{E}^-$ but $\V'$ is the anti-commuting $(1,0)$ part of a connection on $KE_\qat$.
By Lemma~\ref{lem:10 pairing properties}(v) $\V'_l$ also defines a potential in $\pot{(KE)^{-1}}^-$.
We can at least remove the explicit use of coordinates locally.
Define a holomorphic quadratic differential form on $\SO_l$ by $q_l = -\qi dz_l^2$.
Then $q_l \V' q_l^{-1}$ is a potential in $\pot{KE^{-1}}^-$ and $\V = - q_l \V' q_l^{-1}$.
Conversely, if 
\begin{equation}
\label{eq:isothermic quad diff}
q_l \in H^0(\SO_l, K^2)  
\quad\text{and}\quad
\V = - q_l \V' q_l^{-1},
\end{equation}
then the coordinate $z_l$ for which $q_l = -\qi dz_l^2$ is an isothermal coordinate for $F$.
This is essentially~\cite[Remark~18]{BPP}.
The factor of $-\qi$ in the definition of $q$ is of course conventional, but with this convention $\Real q = dx_l\otimes dy_l + dy_l\otimes dx_l$.
In other words the $x_l$ and $y_l$ directions are the null directions of $\Real q$, allowing us to elide the coordinate $z_l$.

\begin{example}[Round Sphere]
\label{eg:sphere isothermic}
\index{Round sphere!Isothermicity}
Let us consider the round sphere, continuing from Example~\ref{eg:sphere weierstrass}.
We already know that $\V = 0$, so it remains to find $\V'$.
From $\psi = -2(\qj z + \qi)^{-1}$ we compute
\begin{align*}
(d \psi)\psi^{-1}
&= 2(\qj z + \qi)^{-1} \qj dz (\qj z + \qi)^{-1} \cdot -\tfrac{1}{2} (\qj z + \qi) 
= -(\qj z + \qi)^{-1} \qj dz \\
&= \frac{1}{|z|^2+1}(\qj z + \qi) \qj dz 
= \frac{1}{|z|^2+1}(- dz\,\Bar{z} + d\Bar{z}\,\qk).
\end{align*}
$\V'$ is the anti-commuting part $(1,0)$ part of this expression, which we see is zero.
We conclude that the usual coordinate $z$ is an isothermal coordinate on $\mathbb{C} \subset \mathbb{P}^1$.
\end{example}

We should discuss the difficulties with umbilic points.
For smooth immersions, an umbilic point is one where for each Weingarten operator the two principal curvatures are equal.
This implies that~\eqref{eq:weingarten matrix} is a multiple of the identity matrix, which forces $\Real ((\V'_l - \V_l)\qj\alpha) = \Imag ((\V'_l - \V_l)\qj\alpha) = 0$.
Not only does this imply that $\V'_l = \V_l$ at this point, but since this must hold for any coordinate $z_l$, and $\V'_l$ and $\V_l$ transform differently, they are each separately zero.
Thus we might define an umbilic point of an admissible map to be a point at which $\V$ and $\V'$ vanish.
If the potentials are only functions in $\banach{2}$, the problem is that this condition is not well-defined.
The condition of totally umbilic is well-defined however, since then $\V$ and $\V'$ must vanish on the whole surface.
We see from the above example that the round sphere is totally umbilic in this sense.
In any case, we see that points for which $\V$ and $\V'$ vanish are not an obstacle for the definition of isothermic in the sense of Definition~\ref{def:isothermic}.

We recall the characterization of Palmer in~\cite[Theorem~I]{Pal} of isothermic surfaces in arbitrary codimension: A smooth conformal immersion $F:\X \to \mathbb{R}^n$ on a Riemann surface $\X$ is \emph{isothermic} if and only if locally there exists conformal immersions $G_l: \SO_l \subset \X \to \mathbb{R}^n$ with parallel tangent planes and opposite orientation, called a local dual map to $F$.
For $n\in\{3,4\}$ we describe a tangent plane with a left and right normal, where the sign of the normals are chosen exactly to describe the orientation of the plane as per Definition~\ref{def:left and right normals}.
Hence a parallel plane with opposite orientation is described by normals with the opposite signs.
Let us extend the definition of dual map to admissible maps, in preparation for the generalization of Palmer's result.

\begin{definition}
\label{def:isothermic admissible}
\label{def:dual map}
\index{Isothermicity!Dual map}
Two admissible maps $F, G: \X \to \qat$ are called \emph{dual} (also known as Christoffel dual) to one another if left and right normals of $G$ are the negatives of the left and right normals of $F$ respectively.
If on $\SO_l \subset \X$ there exists an admissible map $G_l: \SO_l \to\qat$ dual to $F|_{\SO_l}$, we call $G_l$ a local dual map to $F$.
If we have an admissible map $F:\X\to\Imag(\qat)$ with a dual map $G$, then it is possible to translate $G$ (which doesn't change its left and right normals) such that it is also map into $\Imag(\qat)$.
\end{definition}

\begin{example}[Catenoid]
\label{eg:catenoid isothermic}
\index{Catenoid!Isothermicity}
The catenoid is an isothermic surface in $\Imag \qat$.
Already in Example~\ref{eg:catenoid} we found the left and right normals of the catenoid, which are the negatives of one another.
By elementary manipulations
\begin{align*}
N_{\text{cat}}
&= (1- z^{-1}\qk)(-\qi)(1- z^{-1}\qk)^{-1}
= - (\bar{z}- \qk)\qi(\Bar{z} - \qk)^{-1} \\
&= - (\overline{z + \qk})\qi(z + \qk) (1 + |z|^2)^{-1} 
= - (z + \qk)^{-1}\qi(z + \qk) 
= - N_{\text{sphere}},
\end{align*}
where we recognize the left normal of the round sphere from Example~\ref{eg:sphere kodaira}.
The round sphere is also a surface in $\Imag\qat$, so 
\[
R_{\text{cat}}
= - N_{\text{cat}}
= N_{\text{sphere}}
= - R_{\text{sphere}}.
\]
Hence the catenoid is dual to the round sphere restricted to $z \in \mathbb{C}^\ast$.
\end{example}

If there exists an admissible map $G : \X \to \qat$ dual to $F$, we can investigate its Kodaira and Weierstraß representations.
Let $(E,\V,\upsilon,\phi)$ and $(E,\U,\chi,\psi)$ be the Kodaira and Weierstraß representations of $F$.
The left normal of $G$ is $-\upsilon^{-1}\qi\upsilon$ by Equation~\eqref{eq:weierstrass-NR}.
Therefore $G$ has a Weierstraß representation of the form $dG = \lp\upsilon,\phi\pa\rp$ for some $\V\sd$\=/holomorphic section $\phi\pa$ of $KE^{-1}_\qat$.
Note that the symbol $\ast$ in $\phi\pa$ is not an operator; it only denotes another section, similar to a prime.
We use this notation to be compatible with Remark~\ref{rem:pa is not an operator}.
We can read off its Weierstraß representation that the right normal of $G$ is $(\phi\pa)^{-1}\qi\phi\pa$.
This is, per the definition of a dual map, the negative of the right normal of $F$ and the left normal of $-\Bar{F}$ by Example~\ref{eg: - bar F}.
It follows that we have a Kodaira representation $-\Bar{F}=(\phi\pa)^{-1}\upsilon\pa$ on the complement of the roots of $\phi\pa$ for some section $\upsilon\pa$ of $KE^{-1}_\qat$ that is also $\V\sd$\=/holomorphic.
By applying Lemma~\ref{quotient dimension} to $\phi\pa$ we see that $\ord_z \phi\pa \leq \ord_z \upsilon\pa$ for all $z \in \X$.
Defining $\upsilon\pa$ to be zero at the roots of $\phi\pa$, it extends holomorphically to all of $\X$.
Locally the two Kodaira representations imply $(\phi\pa_l)^{-1}\upsilon\pa_l = - \bar{F} = - \overline{\upsilon_l^{-1}\phi_l}$.
The converse is also true.
\begin{lemma}
\label{lem:dual map correspondence}
Let $F: \X \to \qat$ be an admissible map.
Its dual maps $G : \X \to \qat$ are in correspondence (up to translations) with pairs $(\upsilon\pa,\phi\pa) \in H^0(\X,\Q{KE^{-1},\V\sd})^{\times 2}$ that satisfy
\begin{equation}
\label{eq:dual map condition}
\upsilon\pa_l \Bar{\upsilon}_l + \phi\pa_l \Bar{\phi}_l = 0
\end{equation}
and for which $dG = \lp\upsilon,\phi\pa\rp$ is an exact $1$\=/form.
\end{lemma}
\begin{proof}
One direction of implication has already been proven above.
Suppose conversely that a pair of sections $(\upsilon\pa,\phi\pa)\in H^0(\X,\Q{KE^{-1},\V\sd})^{\times 2}$ satisfies~\eqref{eq:dual map condition}.
Since $\phi\pa$ is a section of $KE^{-1}_\qat$, it is possible to pair it with $\upsilon$, and because it is $\V\sd$\=/holomorphic, the resulting form $\lp\upsilon,\phi\pa\rp$ is a closed $\qat$\=/valued $1$\=/form on $\X$.
However, there is no reason to suppose that this form is exact.
But if it is, then we have an admissible map $G : \X \to \qat$ with Weierstraß data $dG=\lp\upsilon,\phi\pa\rp$, unique up to translation.
It remains to show that $G$ is indeed a dual map to $F$.
The left normal of $G$ is $-\upsilon^{-1}\qi\upsilon$ by Equation~\eqref{eq:weierstrass-NR}, which is the negative of the left normal of $F$.
Rearranging~\eqref{eq:dual map condition} gives $(\phi\pa)^{-1}\upsilon\pa = - \overline{\upsilon^{-1}\phi} = -\Bar{F}$ on the complement of the roots of $\phi\pa$.
Since the roots are discrete points, this holds almost everywhere on $\X$.
It follows that the right normal $(\phi\pa)^{-1}\qi\phi\pa$ of $G$ is the left normal of $-\Bar{F}$, and by Example~\ref{eg: - bar F} again the left normal of $-\Bar{F}$ is the negative of the right normal of $F$.
Hence the right normal of $G$ is the negative of the right normal of $F$.
This shows that $G$ is a dual map to $F$.
\end{proof}

\begin{example}[Catenoid]
\label{eg:catenoid isothermic kernel}
\index{Catenoid!Isothermicity}
Consider the catenoid.
We have $\V = \qk z^{-1}\Bar{z} (1 + |z|^2)^{-1}$ and $\V\sd = \qk z \Bar{z}^{-1} (1 + |z|^2)^{-1}$.
Observe that $\V\sd = z^{-2} \V z^2$.
In particular, if $\xi$ is $\V$\=/holomorphic then $z^{-2}\xi$ is $\V\sd$\=/holomorphic.
Using this observation, and the fact that $- \Bar{F} = F$ since $F \in \Imag\qat$, we see
\begin{align*}
-\lp z^{-2} \phi, \upsilon \rp - \lp z^{-2}\upsilon, \phi \rp
&= - 2\qk z^{-2} \Big[ \upsilon F\, \Bar{\upsilon} + \upsilon\, (\overline{\upsilon F}) \Big]
= 0.
\end{align*}
Thus $(\upsilon\pa,\phi\pa) = (z^{-2} \phi, z^{-2} \upsilon)\in H^0(\X,\Q{KE^{-1},\V\sd})^{\times 2}$ satisfy~\eqref{eq:dual map condition}.
\end{example}

We are now ready to generalize to admissible maps the characterization of isothermic surfaces by the existence of local dual maps.
First assume that there exist local dual maps $G_l$ to $F$ on a cover of $\X$.
We may assume $\SO_l$ to be discs, so that exactness is guaranteed and the local dual maps $G_l$ are in correspondence with pairs of $\V\sd$\=/holomorphic functions $(\upsilon\pa_l,\phi\pa_l)$ that satisfy~\eqref{eq:dual map condition}.
We have the two Weierstraß representations $dF = \lp \chi, \psi\rp$ and $dG_l = \lp \upsilon_l, \phi\pa_l\rp$. 
By comparing right normals $(\phi\pa_l)^{-1}\qi\phi\pa_l = - \psi^{-1}\qi\psi$, per the definition of dual map.
For elegance, let $\upsilon'_l := (\overline{\qj \psi_l})^{-1}$ on the complement of the roots of $\psi$ in $\SO_l$, which is a $(-\V'_l)\sd$\=/holomorphic section of $(KE)^{-1}_\qat$.
Then the relation of the normals reads $(\phi\pa_l)^{-1}\qi\phi\pa_l = (\upsilon'_l)^{-1}\qi\upsilon'_l$.
Hence there exists a section $q_l$ of $H^0(\SO_l,K^2)$ with $\phi\pa_l = q_l \upsilon'_l$.
Now $\phi\pa_l$ is $\V\sd$\=/holomorphic if and only if
\begin{gather*}
0
= (\barpartial_l - \V\sd_l)(q_l\upsilon'_l) 
= (\barpartial_l q_l - q_l(\V'_l)\sd -\V\sd_l q_l)\upsilon'_l \\
\Leftrightarrow \qquad
\barpartial_l q_l = 0 
\quad\text{and}\quad 
\V_l = - q_l^{-1} \V'_l q_l
\end{gather*}
This is exactly~\eqref{eq:isothermic quad diff}, which we know to be equivalent to isothermic.
This proves:

\begin{lemma}
\label{lem:isothermic local duals}
An admissible map $F : \X \to \qat$ is isothermic if and only if there exists a cover $\{\SO_l\}$ of $\X$ and admissible maps $G_l: \SO_l \to\qat$ dual to $F|_{\SO_l}$.
\qed
\end{lemma}

Let us now indulge in an exploration of some of the properties of dual maps.
As might be expected given their relation to isothermicity, for simply connected $\X$, the existence of a dual map is a Möbius invariant.
Let $\Tilde{\phi}=\phi\alpha+\upsilon\beta$ and $\Tilde{\upsilon}=\phi\gamma+\upsilon\delta$ be a Möbius transformation from $F$ to $\Tilde{F}$.
If $\Tilde{F}$ has a dual, let $(\Tilde{\upsilon}\pa,\Tilde{\phi}\pa)\in H^0(\X,\Q{KE^{-1},\V\sd})^{\times 2}$ be the corresponding pair.
Rearranging $\Tilde{\upsilon}\pa \overline{\Tilde{\upsilon}} + \Tilde{\phi}\pa \overline{\Tilde{\phi}} = 0$ yields
\[
(\Tilde{\phi}\pa \Bar{\beta} + \Tilde{\upsilon}\pa\Bar{\delta})\Bar{\upsilon} 
+ (\Tilde{\phi}\pa \Bar{\alpha} + \Tilde{\upsilon}\pa\Bar{\gamma}) \Bar{\phi}
= 0,
\]
so the pair $(\upsilon\pa,\phi\pa) 
= (\Tilde{\phi}\pa \Bar{\alpha} + \Tilde{\upsilon}\pa\Bar{\gamma}, \Tilde{\phi}\pa \Bar{\beta} + \Tilde{\upsilon}\pa\Bar{\delta})$ is also $\V\sd$\=/holomorphic and satisfies~\eqref{eq:dual map condition} with respect to $F$.
The exactness condition does not necessarily carry over for general Riemann surfaces however, which is why we have assumed that $\X$ is simply connected.
Therefore $F$ also has a dual map.

A natural question is whether the dual map is unique.
In a trivial sense the answer is no: we have already mentioned that a dual map may be translated without changing its left and right normals.
Similarly we can scale a dual map to produce another.
In the correspondence with $(\upsilon\pa,\phi\pa)$ this is multiplying both sections by a real scalar.
So a more precise question is whether the dual map is unique up to scaling and translation.
In the following lemma, we show that admissible maps with more than one dual map are very special: they are contained in a plane or a round sphere.
This mirrors the result for smooth isothermic immersions, where it implies that the immersion is totally umbilic.
As is well-known for surfaces $\mathbb{R}^3$, but also true in higher dimension~\cite[Chapter~7,~Theorem~26]{Sp}, connected totally umbilic surfaces are contained in a plane or a round sphere.

\begin{theorem}
\label{thm:isothermic local characterization}
Let $F: \mathbb{D} \to \qat$ be an admissible map.
The following are equivalent
\begin{enumeratethm}
\item 
$F$ has more than one independent dual map (that are not related by scaling or translation).
\item
$F$ has a dual map and its Kodaira potential vanishes $\V = 0$.
\item 
$F$ is totally umbilic.
\item
$F[\mathbb{D}]$ is contained in a $2$\=/plane or a round $2$\=/sphere.
\end{enumeratethm}
In this case, the dual maps are minimal surfaces and are in correspondence with the holomorphic functions on $\mathbb{D}$.
\end{theorem}
\begin{proof}
First we show that (i) implies (ii).
We operationalize the existence of two independent dual maps as pairs $(\upsilon\pa_l,\phi\pa_l) \in H^0(\mathbb{D},\Q{KE^{-1},\V\sd})^{\times 2}$ for $l=1,2$ that are $\mathbb{R}$\=/linearly independent.
On the complement of the roots of $\phi\pa_1$ and $\phi\pa_2$, Theorem~\ref{thm:kodaira} shows that the left normal of $-\Bar{F}$ is equal to $(\phi\pa)^{-1}_2\qi\phi\pa_2=(\phi\pa)^{-1}_1\qi\phi\pa_1$.
Rearranging gives $\qi\phi\pa_2(\phi\pa_1)^{-1}=\phi\pa_2(\phi\pa_1)^{-1} \qi$, which is only possible if this product commutes with $\qi$.
Hence we define a $\mathbb{C}$\=/valued function $f$ by $\phi\pa_2 = f \phi\pa_1$.
Since both are $\V\sd$\=/holomorphic this implies
\begin{align*}
0 = (\barpartial - \V\sd)\phi\pa_2 = (\barpartial f + f\V\sd - \V\sd f) \phi\pa_1.
\end{align*}
The equation forces both $\barpartial f = 0$ and $\V\sd = 0$ on the complement of the roots of $\phi\pa_1$.
Since the roots of $\phi\pa_1$ are isolated, we can set $\V = 0$ everywhere.

For the converse implication and the final statement of the lemma, suppose (ii), i.e.\ $\V = 0$ and $(\upsilon\pa_1, \phi\pa_1) \in H^0(\mathbb{D},\Q{KE^{-1},\V\sd})^{\times 2}$ with~\eqref{eq:dual map condition}.
From the Weierstraß representation of the dual map $dG=\lp\upsilon,\phi\pa_1\rp$ we see that the Weierstraß potential of $G$ vanishes.
In other words, $G$ is a minimal surface.
But we have discussed above that any minimal surface belongs to a large family of minimal surfaces all of which have the same left and right normals.
Hence every member of this family is a dual map to $F$, proving (i).

In the situation of this lemma, the Riemann surface is $\X = \mathbb{D}$, so we can give an even better description of this family.
The pair $(f\upsilon\pa_1,f\phi\pa_1)$ belongs to $H^0(\X,\Q{KE^{-1},\V\sd})^{\times 2}$ and obeys~\eqref{eq:dual map condition} for any meromorphic function $f$ so long as the product with the sections remains holomorphic.
Because $\upsilon$ has no roots, Equation~\eqref{eq:dual map condition} tells us that $\ord_z \upsilon\pa_1 \geq \ord_z \phi\pa_1$ for all $z\in\mathbb{D}$.
We can find an $f$ without roots and with poles of order equal to the order of the roots of $\phi\pa_1$.
Hence we have a pair $(\upsilon\pa, \phi\pa)$ with~\eqref{eq:dual map condition} and $\phi\pa$ has no roots, and a Kodaira representation $-\Bar{F} = (\phi\pa)^{-1}\upsilon\pa$ on all of $\mathbb{D}$.
Any dual map corresponds to a pair $(f \upsilon\pa, f \phi\pa)$ for $f$ holomorphic.

As discussed above, an admissible map is totally umbilic if both the Kodaira potential $\V$ and $\V' = ((\partial \psi)\psi^{-1})^-$ vanish everywhere.
That $\V$ vanishes is already a condition of (ii), so in order to prove (ii) implies (iii) we should show that additionally $\V' = 0$.
Because we have a Kodaira representation $-\Bar{F} = (\phi\pa)^{-1}\upsilon\pa$, its left normal $(\phi\pa)^{-1}\qi\phi\pa$ is the negative of the right normal $\psi^{-1}\qi\psi$ of $F$, i.e.\ $\psi(\phi\pa)^{-1}\qi = - \qi\psi(\phi\pa)^{-1}$.
This time the quotient anti-commutes with $\qi$ (compare to $f$ above), so we conclude that there is a $\mathbb{C}$\=/valued function $g$ with $\psi = g \qj \phi\pa$.
Differentiating $\psi = g \qj \phi\pa$ shows
\[
\partial \psi 
= \partial g \qj \phi\pa + g \qj \barpartial \phi\pa
= (\partial g) g^{-1} \psi + 0.
\]
In particular the anti-commuting part $\V'$ is zero.
Hence $F$ is totally umbilic.

We now prove that (iii) implies (iv), which is the most difficult of the implications.
Because $\V = 0$ we know that $\upsilon$ is analytic, and because it is non-vanishing the other potentials, defined by $\partial \upsilon = -(\B + \U)\upsilon$, are also analytic.
Hence $U$ and $B$ may be differentiated.
Applying $\barpartial$ yields
\begin{align*}
0 
&= \barpartial\B \upsilon + \B\barpartial \upsilon + \barpartial \U \upsilon + \U\partial \upsilon 
= (\barpartial\B + \barpartial \U - \U\B - \U^2) \upsilon 
\end{align*}
Hence $\barpartial\B = \U^2$ and $\barpartial \U = \U\B$ (compare to the introduction of Chapter~\ref{chapter:darboux}).
Similarly $\V' = 0$ implies $\partial \psi = C \psi$ for a $\qat^+$\=/potential $C$, per the definition of $\V'$, in addition to $\barpartial \psi = \U \psi$.
The compatibility between these two derivatives is
\begin{align*}
\partial\barpartial \psi 
&= \partial \U \psi + \U \barpartial \psi
= \partial \U \psi + \U^2 \psi \\
\barpartial\partial \psi 
&= \barpartial C \psi + C \barpartial \psi
= \barpartial C \psi + C \U\psi.
\end{align*}
Therefore $\partial \U = C\U$ and $\barpartial C = \U^2$.

We now put together these facts about the potentials to argue that $F$ is contained in a plane or round sphere.
In the present situation, Equation~\eqref{eq:dN potentials} says $dN = -2\upsilon^{-1} dz \qi \U \upsilon$ and Equation~\eqref{eq:dR potentials} says $dR = 2\psi^{-1} d\Bar{z} \qi \U \psi$.
Since $\U$ is analytic on $\mathbb{D}$ either $\U \equiv 0$ or its roots are isolated.
In the former case, then both normals are constant and we have a map to a plane.
So assuming it's not identically vanishing define a function $\alpha$ through $\psi = -2\qi \U\upsilon \alpha$ on the complement of the roots of $\U$.
The geometric meaning of $\alpha$ will become apparent shortly.
Differentiating both sides yields
\begin{gather*}
- 2\qi \U \upsilon d\alpha
= \left[ - (dz C + d\Bar{z} \U) 2\qi \U 
+ 2\qi (dz C\U + d\Bar{z} \U B)
+ 2\qi \U dz ( -\B - \U) \right] \upsilon \alpha 
= 0.
\end{gather*}
Therefore $\alpha$ is in fact constant.
Now from the Weierstraß representation
\begin{align*}
dF
= \Bar{\chi}\qj dz \psi
= -\upsilon^{-1} dz \, 2\qi \U\upsilon \alpha
= dN\, \alpha.
\end{align*}
Integrating gives $F = N \alpha + \beta$, showing that $F$ is contained in a round $2$\=/sphere of radius $|\alpha|$ centered at $\beta$.
This argument is rather directly an adaption of the standard argument that totally umbilic surfaces are contained in a round sphere or plane.


Finally, we prove that (iv) implies (ii).
Suppose that we have an admissible map $\Tilde{F}$ to a $2$\=/plane or a $2$\=/sphere.
By Möbius transformation, we assume without loss of generality that $\Tilde{F} : \mathbb{D} \to \mathbb{S}^2$.
From Example~\ref{eg:kodaira factorization} we know that $\Tilde{F}$ factors as a holomorphic map $f : \mathbb{D} \to \mathbb{P}^1$ composed with the map to the round sphere $F : \mathbb{P}^1 \to \mathbb{S}^1$.
Combining Examples~\ref{eg:sphere kodaira} and~\ref{eg:kodaira covering} we even have its Kodaira data
\[
\Tilde{\upsilon}_1(z)
= f(z) + \qk,
\qquad
\Tilde{\phi}_1(z)
= \qi \Tilde{\upsilon}_1(z)
= \qi (f(z) + \qk)
\]
on the complement of $f^{-1}[\infty] \subset \mathbb{D}$, and around the points of $\mathbb{D}$ that $f$ maps to $\infty \in \mathbb{P}^1$ we would change to a different trivialization of the bundle $\Tilde{E} = f^\ast E$.
However, we can modify this Kodaira representation of $\Tilde{F}$, essentially using the fact that $f^\ast E$ is trivial over $\mathbb{D}$ even though $E$ is not necessarily trivial over $\mathbb{P}^1$.
Indeed, there exists a holomorphic function $g : \mathbb{D} \to \mathbb{C}$ which has roots at the poles of $f$ to the same order.
Then $fg$ is a non-vanishing holomorphic function and we have the Kodaira representation on all of $\mathbb{D}$
\[
\Tilde{\upsilon}(z) 
= g(z)f(z) + g(z)\qk, \qquad
\Tilde{\phi}(z) 
= \qi \Tilde{\upsilon}(z)
= \qi (g(z)f(z) + g(z)\qk).
\]
But now $(\upsilon\pa, \phi\pa) = (\Tilde{\phi}, \Tilde{\upsilon})$ is a solution to~\eqref{eq:dual map condition} since
\[
- \phi\pa \overline{\Tilde{\phi}}
= \Tilde{\upsilon} \overline{\Tilde{\upsilon}} \qi
= \qi \Tilde{\upsilon} \overline{\Tilde{\upsilon}}
= \Tilde{\phi} \overline{\Tilde{\upsilon}}
= \upsilon\pa \overline{\Tilde{\upsilon}}.
\qedhere
\]
\end{proof}

We now turn our investigation to the difference between having a dual map and being isothermic.
Clearly, if $F : \X \to \qat$ has a dual map $G : \X \to \qat$ then the latter's restrictions give local dual maps $G_l = G|_{\SO_l}$.
The question is whether all isothermic surfaces have a (global) dual map.
Here is an argument that almost succeeds to show this reverse implication.
Suppose that $F$ is isothermic but not contained in a plane or round sphere.
Then by the above lemma, each local dual map is unique (up to scaling and translation).
This gives a unique continuation, allowing us to glue the maps on the overlaps $\SO_l \cap \SO_m$ and assemble the local dual maps $G_l$ into a map $G: \Tilde{\X}\to \qat$.
Let us give a little more detail.
On $\SO_l$ we have $dG_l = \lp \upsilon_l, \phi\pa_l \rp$.
Since $G_l$ is determined up to real scaling and translation, $\phi\pa_l$ is unique up to real scaling.
This means, on $\SO_l \cap \SO_m$ we can assume that $G_m = G_l$ without loss of generality.
Now we see that on the overlap
\[
\overline{\upsilon_m}\qj dz_m \phi\pa_m
= dG_m
= dG_l
= \overline{\upsilon_l}\qj dz_l \phi\pa_l
\quad \Rightarrow \quad
\phi\pa_m = f_{ml}^{-1} \frac{dz_l}{dz_m} \phi\pa_l,
\]
which is to say that these local functions combine to a section of $KE^{-1}_\qat$ over $\SO_l \cup \SO_m$.
Continuing in this way, we construct on $\Tilde{\X}$ a $(\V\sd\circ p)$\=/holomorphic section $\phi\pa$ of $p^\ast (KE^{-1}_\qat)$ with $dG = \lp \upsilon \circ p, \phi\pa \rp$.
Note the appearance of the universal cover here: if $\X$ is not simply connected, then there may be monodromy in this operation.
To be precise, let $\Lambda : \Tilde{\X} \to \Tilde{\X}$ be a deck transformation.
The possible translation and scaling monodromy of $G$ means exactly that there exist constants $\alpha \in \mathbb{R}$ and $\beta \in \qat$ such that $G \circ \Lambda = G \alpha + \beta$ and $\phi\pa \circ \Lambda = \phi\pa \alpha$.
Then $G$ descend to a map on $\X$, in which case it is dual to $F$, exactly if the monodromy is trivial.

\begin{definition}
\label{def:strongly isothermic}
\index{Isothermicity!Strongly isothermic}
An isothermic surface $F : \X \to \qat$ is called \emph{strongly isothermic} if there exists an admissible map $G : \tilde{\X} \to \qat$ on the universal cover $p : \tilde{\X} \to \X$ such that $F \circ p$ and $G$ are dual maps and any deck transformation $\Lambda : \Tilde{\X} \to \Tilde{\X}$ acts on $G$ as a translation $G \circ \Lambda = G + \beta_{\Lambda}$ for $\beta_{\Lambda} \in \qat$.
\end{definition}

We can understand strongly isothermic as an intermediate condition on $F$; it is stronger than isothermic but weaker that having a dual map.
By allowing translation monodromy but excluding scaling monodromy, it means that $dG$ descends to a differential form on $\X$, and in particular it has a Weierstraß representation with data which is single valued on $\X$.
We see that these $G$ are exactly the result of dropping the exactness condition in Lemma~\ref{lem:dual map correspondence}.
For simply connected Riemann surface, the exactness condition was already superfluous, hence in this case there is no difference between isothermic and strongly isothermic.
As this will be important in Chapters~\ref{chapter:isothermic kodaira triples} and~\ref{chapter:isothermic weierstrass triples} we restate the lemma for this case.

\begin{lemma}
\label{lem:strongly isothermic correspondence}
Let $F: \X \to \qat$ be an admissible map.
It is strongly isothermic if and only if there is a pair $(\upsilon\pa,\phi\pa) \in H^0(\X,\Q{KE^{-1},\V\sd})^{\times 2}$ that satisfies~\eqref{eq:dual map condition}.
The Weierstraß representation of the dual map $G$ on $\tilde{X}$ is $dG = \lp\upsilon,\phi\pa\rp$.
If $\X$ is simply connected, all isothermic admissible maps are strongly isothermic.
\qed
\end{lemma}

The name strongly isothermic comes from the appendix of~\cite{BPP}, where strongly isothermic surfaces in $\mathbb{R}^3$ are defined as those that have a non-trivial holomorphic quadratic differential $q \in H^0(K^2)$ such that the null-directions of $\Real(q)$ are principal curvature directions.
These two notions are equivalent.
We have already seen this equivalence locally.
In the case of a strongly isothermic map in the sense of Definition~\ref{def:strongly isothermic}, the holomorphic quadratic form $q$ is defined globally through $\phi\pa = q \upsilon'$ for $\upsilon' = (\overline{\qj \psi})^{-1}$.
Conversely, the existence of such a non-trivial holomorphic quadratic differential $q \in H^0(K^2)$ defines $\phi\pa = q \upsilon'$ such that the necessary $\upsilon\pa$ also exists.

For compact Riemann surfaces $\X$ this gives a constraint on $F$ to be strongly isothermic, and therefore also on the existence of a dual map.
Namely, there must exist a non-trivial holomorphic section $q$ of $K^2$, whose degree is $4\genus-4$ where $\genus$ denotes the genus of $\X$.
More precisely, the root divisor of $q$ is the sum of the root divisors of $\psi$ and $\phi\pa$, which are respectively the root divisors of $dF = \lp \chi, \psi\rp$ and $dG = \lp \upsilon, \phi\pa\rp$.
In the case $\X = \mathbb{P}^1$, there can be no strongly isothermic admissible maps, as the number of roots cannot be negative.
We can use Example~\ref{eg:sphere isothermic}, which shows that the round sphere is isothermic, to illustrate this obstacle.
In that example, the isothermicity gave us a the holomorphic quadratic differential $q = -\qi dz^2$ on $\mathbb{C}\subset \mathbb{P}^1$.
However, if we attempt to extend this to the whole surface then it has a fourth order pole.
In the case of a compact Riemann surface of genus $1$ that is strongly isothermic, its derivative $dF$ cannot have roots and so it must be an immersion.
Willmore and constrained Willmore strongly isothermic tori have been investigated~\cite{Thomsen1924,Richter1997thesis,Burstall2002b,Bohle2012a}.

\chapter{Plücker Formula}
\label{chapter:pluecker}

In this chapter we prove a version of the Plücker formula for holomorphic $\qat$\=/line bundles that bounds the $\banach{2}$\=/norm of the potential in terms of algebraic data, 
namely the genus of the Riemann surface, the degree of the holomorphic $\mathbb{C}$\=/line bundle, the number of holomorphic sections, and the order of their roots.
This reproduces a result of~\cite{FLPP}.
The main application, Corollary~\ref{thm:embedding}, is a lower bound on the Willmore energy in terms of the roots of the sections of the Weierstraß representations, which implies a well-known result of Li and Yau~\cite{LY} that any immersion of a compact Riemann surface with Willmore energy lower than $8\pi$ is injective.
Throughout this chapter, let $E$ be a holomorphic $\mathbb{C}$\=/line bundle on a compact Riemann surface $\X$ and $H$ a finite-dimensional $\qat$\=/linear subspace of the right-$\qat$ vector space $H^0(\X,\Q{E,\V})$ of $\V$\=/holomorphic sections with $\V\in\pot{E}^-$.

At any point $x\in\X$ we show that there exists $d:=\dim_\qat H$ many distinct numbers that occur as orders of sections of $H$ at $x$.
Since any choice of sections $\upsilon_1,\ldots,\upsilon_d$ with these orders $\ord_x(\upsilon_l)$ form a basis, these numbers are unique.
Choose any basis $\upsilon_1,\dotsc,\upsilon_d$ of $H$ and arrange it so that $\ord_x \upsilon_1=\min\{\ord_x \upsilon_1,\ldots,\ord_x \upsilon_d\}=\min\{\ord_x \upsilon\mid\upsilon \in H\}$.
If $z$ is a local coordinate at $x$, then $z^{-\ord_x(\upsilon_1)}\upsilon_l$ are sections of $\Q{E(-\ord_x(\upsilon_1)),\V(\frac{z}{\Bar{z}})^{\ord_x(\upsilon_1)}}$ and $\upsilon_1^{-1}\upsilon_l=(z^{-\ord_x(\upsilon_1)}\upsilon_1)^{-1}z^{-\ord_x(\upsilon_1)}\upsilon_l$ admissible and continuous nearby $x$ for $l=1,\dotsc,d$.
So we replace $\upsilon_l$ by $\upsilon_l - \upsilon_1\cdot (\upsilon_1^{-1}\upsilon_l)(x)$ for $l=2,\ldots,d$ and obtain $\ord_x(\upsilon_l) > \ord_x(\upsilon_1)$ by Lemma~\ref{step 5}.
Next we repeat the same procedure for $\upsilon_2,\dotsc,\upsilon_d$.
Inductively we get a basis $\upsilon_1,\dotsc,\upsilon_d$ such that 
\begin{gather}\label{sequence of orders}
\ord_{x}\!_{1} H<\cdots<\ord_{x}\!_d H 
\quad\text{for}\quad
\ord_{x}\!_{l} H := \ord_x(\upsilon_l),\;l=1,\dotsc,d.
\end{gather}
We shall see that this sequence differs only at finitely many points from $0,\ldots,d-1$.
The order of $H$ is defined as~\cite[Definition~4.2.]{FLPP}:
\[\ord H=\sum_{x\in\X}\left(\ord_{x}\!_{1} H-0\right)+\cdots+\left(\ord_{x}\!_d H-d+1\right).\]
For smooth $\V$ the following estimate is proven in~\cite[Corollary~4.8.]{FLPP}:
\begin{align}\label{eq:pluecker}
\frac{1}{\pi}\|\V\|_2^2\geq d\left((d-1)(1-\genus)-\deg(E)\right)+\ord H.
\end{align}
Here $\genus$ denotes the genus of $\X$ and $\|\V\|_2^2$ is defined in~\eqref{eq:square integrable}.
For square-integrable $\V$ this inequality follows from a $d$\=/fold application of Corollary~\ref{cor:global darboux}.
First we choose a member $\upsilon$ of the linear system $H$ with roots of least order at all points of $\X$.
Let us prove that this is always possible:
\begin{lemma}
\label{lem:minimal-order}
\index{Root}
For a holomorphic $\mathbb{C}$\=/line bundle $E$ and $\V\in\pot{E}^-$ on a compact Riemann surface $\X$ any non-trivial, finite-dimensional $\qat$\=/right linear subspace $H\subset H^0(\X,\Q{E,\V})$ contains $\upsilon$ with $\ord_x(\upsilon)=\min\{\ord_x(\phi)\mid\phi\in H\}$ simultaneously for all $x\in\X$.
\end{lemma}
\begin{proof} 
We prove the statement by induction on $d:= \dim_\qat H$.
In the case $d=1$, the statement is true because all elements of $H\setminus \{0\}$ have the same divisor.
We now suppose that the statement of the Lemma holds for some $d \geq 1$.
A subspace $H$ with $\dim_\qat H = d+1$ can be written as the direct sum $H\pa \oplus \phi \qat$ with $\phi \in H\setminus \{0\}$.
Choose $\upsilon \in H\pa$ with $\ord_x(\upsilon)=\min\{\ord_x(\xi)\mid\xi\in H\pa\}$ for all $x\in\X$.
Let $D$ be the pointwise minimum of all divisors $(\xi)$ with $\xi \in H$.
Observe that $D$ is in fact the pointwise minimum of $(\upsilon)$ and $(\phi)$, since $(\upsilon\qat + \phi\qat) \geq \min\{(\upsilon), (\phi)\}$.

The $\V$\=/holomorphic sections $\upsilon$ and $\phi$ induce $\Tilde{\upsilon},\Tilde{\phi}\in H^0(\X,\Q{E(-D),\V(-D)})$ without common roots.
Choose an open neighborhood $\SO$ with smooth boundary that contains all roots of $\Tilde{\upsilon}$, but no roots of $\Tilde{\phi}$.
Then $\Tilde{\upsilon}^{-1}\Tilde{\phi}$ is admissible on the complement of the roots of $\Tilde{\upsilon}$ and bounded on the compact set $\X\setminus \SO$ by $\epsilon^{-1}$ for some $\epsilon>0$.
We claim that the image $F[\SO]$ of the admissible map $F=\Tilde{\phi}^{-1}\Tilde{\upsilon}$ on $\SO$ does not contain $B(0,\epsilon)\subset\qat$.
This follows by an argument of Kaufman~\cite[Proof of Theorem~1]{Ka}, which shows that the Hausdorff dimension of the image $F[\SO]$ is at most 2.
Let us recall this argument: By the Sobolev Embedding theorem $F$ belongs to the Hölder space $\sobolev{2,\frac{2q}{2+q}}(\SO,\qat) \subset \sobolev{1,q}(\SO,\qat)\subset C^{0,1-\frac{2}{q}}(\Bar{\SO},\qat)$ with $2<q<\infty$.
Indeed the proof of~\cite[Lemma~4.28]{Ad} gives on $\Omega\subset\mathbb{C}$ with Lipschitz boundary the bound
\[
\sup_{x\ne y\in\Omega}\frac{|u(x)-u(y)|}{|x-y|^{1-\frac{2}{q}}}\le C(q,\Omega)\|\nabla u\|_{\banach{q}(\Omega,\mathbb{R})}
\quad\text{for all}\quad 
u\in C^\infty_0(\Omega,\mathbb{R}),
\]
with some constant $C(q,\Omega)$ depending on $q$ and $\Omega$.
If we cover $\SO$ by small essentially disjoint cubes $Q_n$ of side $r_n$, then this yields the following estimate:
\[\diam(F[Q_n])\le C(q)r_n^{1-\frac{2}{q}}\|\nabla F\|_{\banach{q}(Q_n,\qat)},\]
with $C(q)$ depending only on $q$.
With the Hölder inequality we obtain
\[\sum_n\diam^2(F[Q_n])\le\left(\sum_n(C(q))^{\frac{2q}{q-2}}r_n^2\right)^{1-\frac{2}{q}}\left(\sum_n\|\nabla F\|_{\banach{q}(Q_n,\qat)}^q\right)^{\frac{2}{q}} < \infty.\]
The finiteness of the sum on the left hand side (compare~\cite[1.4.1 Definition]{Zi}) implies that the Hausdorff dimension of $F[\SO]$ is at most $2$, so our claim that there exist $\alpha\in B(0,\epsilon)\setminus F[\SO]$ is true.
Then $\Tilde{\upsilon}-\Tilde{\phi}\alpha$ has no root on $\SO$ by Lemma~\ref{step 5}.
On $\X\setminus \SO$ it is equal to $(1-\Tilde{\phi}\alpha\Tilde{\upsilon}^{-1})\Tilde{\upsilon}$ which does not vanish by $|\Tilde{\phi}\alpha\Tilde{\upsilon}^{-1}|=|\alpha|\cdot|\Tilde{\phi}\Tilde{\upsilon}^{-1}|<1$.
So $\upsilon-\phi \alpha\in H$ has the required properties.
\end{proof}
\begin{remark}
The arguments of the proof reflects the difference between the dimensions of the domain $\X$ and the codomain $\qat$ of admissible maps.
For differentiable functions, the lemma follows from Sard's theorem.
From this point of view, Kaufman's estimate is a version of Sard's theorem for Sobolev functions.
For holomorphic complex functions, the dimensions of domain and codomain are equal, which is why an analogous statement does not hold.
\end{remark}
We transform the linear system $H\subset H^0(\X,\Q{E,\V})$ into a linear system $\Tilde{H}\subset H^0(\X,\Q{E(-D),\V(-D)})$ that contains an element $\Tilde{\upsilon}$ without roots using Lemma~\ref{lem:minimal-order}.
Now we can apply the Darboux transform~\ref{cor:global darboux} to pairs $\Tilde{\upsilon}$ and $\phi$, where $\phi$ can be any element of $\Tilde{H}$.
The map $\phi \mapsto \psi$ transforms $\Tilde{H}$ into another linear system $H_1\subset H^0(\X,\Q{KE(-(\upsilon)),\V_1})$ with square-integrable potential $\V_1\in\pot{KE(-(\upsilon))}^-$ satisfying (compare~\eqref{eq:hopf norms})
\[\frac{\|\V_1\|_2^2-\|\V\|_2^2}{\pi}=\deg(E(-(\upsilon))) =\deg(E)-\deg(\upsilon)=\deg(E)-\sum_{x\in \X} \ord_x\!_1 H.\]
The kernel of the map $\phi\mapsto \psi$ is spanned by $\Tilde{\upsilon}$, so $\dim(H_1)=d-1$.
Due to Lemma~\ref{lem:change of order}, for any $\phi\in H^0(\X\Q{E(-D),\V(-D)})$ with a root at $x\in\X$ the transformed $\psi_m=(\partial_m+\B_m+\U_m)\phi$ has a root at $x$ of order $\ord_x(\psi)=\ord_x(\phi)-1$.

For $d>1$ we repeat the whole procedure: We choose by Lemma~\ref{lem:minimal-order} a section $\upsilon_1\in H_1$ with minimal order of roots and replace $H_1$ by the corresponding linear system $\Tilde{H}_1$ in $H^0(\X,\Q{KE(-(\upsilon)-(\upsilon_1)),\V_1(-(\upsilon_1))})$.
The corresponding element $\Tilde{\upsilon}_1\in\Tilde{H}_1$ has no roots, and the Darboux transformation $\phi \mapsto \psi$ from Corollary~\ref{cor:global darboux} with $\Tilde{\upsilon}_1$ transforms $\Tilde{H}_1$ into a linear system $H_2$ in $H^0(\X,\Q{K^2E(-(\upsilon)-(\upsilon_1)), \V_2})$ with some potential $\V_2 \in \pot{K^2E(-(\upsilon)-(\upsilon_1))}^-$.
By Lemma~\ref{lem:change of order}(i) we have
\[\deg(\upsilon_1)=\sum_{x\in \X} (\ord_x\!_2 H-\ord_x\!_1 H-1)=\sum_{x\in \X}(\ord_x\!_2H-1)-\deg(\upsilon),\]
and the sum on the right hand side is finite.
Therefore gives~\eqref{eq:hopf norms}
\begin{multline*}
\frac{\|\V_2 \|_2^2-\|\V_1\|_2^2}{\pi}=\deg(KE(-(\upsilon)-(\upsilon_1))\\=\deg(E)+2\genus-2-\deg(\upsilon)-\deg(\upsilon_1)
=\deg(E)+2\genus-2-\sum_{x\in \X}(\ord_x\!_2 H-1).
\end{multline*}
We repeat such applications of Lemma~\ref{lem:minimal-order}  and  Corollary~\ref{cor:global darboux} until we end up with a trivial linear system with potential $\V_d$.
In each step the elements of the linear system $H_l$ are differentiated and mapped to a linear system $H_{l+1}$ of one dimension less.
Due to Lemma~\ref{lem:change of order}(i), the degree of the divisor $(\upsilon_l)$ of the corresponding element $\upsilon_l\in H_l$ constructed in Lemma~\ref{lem:minimal-order} is equal to
\[\deg(\upsilon_l)=\sum_{x\in\X}(\ord_x\!_{l+1}H-\ord_x\!_l H-1)=\sum_{x\in\X}(\ord_x\!_{l+1}H-l)-\sum_{x\in\X}(\ord_x\!_l H-l+1).\]
In particular, $\ord H$ is equal to the following finite sum and therefore finite:
\[\ord H=\deg(\upsilon)+\deg(\upsilon_1)+\cdots+\deg(\upsilon_{d-1}).\]
Furthermore, the sum over the degrees on the right hand sides of~\eqref{eq:hopf norms} is equal to $\deg(E)d-\ord H+(2\genus-2)\sum_{l=0}^{d-1}l$.
Consequently the negative of the sum of the equations~\eqref{eq:hopf norms} for all these Darboux transformations gives
\begin{align*}
\frac{\|\V\|_2^2-\|\V_d\|_2^2}{\pi}&=-\deg(E)d+\ord H-\sum_{l=0}^{d-1}l\deg(K))\\&=d\left((d-1)(1-\genus)-\deg(E)\right)+\ord H.
\end{align*}
This implies
\begin{corollary}[Quaternionic Plücker Formula]
\label{pluecker formula}
\index{Plücker formula}
For a holomorphic $\mathbb{C}$\=/line bundle $E$ on a compact Riemann surface and $\V\in\pot{E}^-$ a linear system $H\subset H^0(\X,\Q{E,\V})$ satisfies~\eqref{eq:pluecker}.
\qed
\end{corollary}
\begin{corollary}\cite[Theorem~4.12.]{BFLPP}\label{dimension bound}
For a  holomorphic $\mathbb{C}$\=/line bundle $E$ on a compact Riemann surface and $\V\in\pot{E}^-$ we have
\begin{align*}
\frac{\|\V\|_2^2}{\pi}&\geq\begin{cases}d((d-1)(1-\genus)-\deg(E))&\text{ for }\deg(E)<(d-1)(1-2\genus)\\\left\lceil\frac{(d-1+\genus-\deg(E))^2-\genus^2}{4\genus}\right\rceil&\text{ for }(d-1)(1-2\genus)\le\deg(E)<d-1.\end{cases}
\end{align*}
\end{corollary}
\begin{proof}
For $\genus=0$ the upper case is the Plücker formula and the lower case is void.
For $\genus\geq 1$, given $x\in\X$ and $H=H^0(\X,\Q{E,\V})$ we consider the sequence $\upsilon_1,\dotsc,\upsilon_d$ in~\eqref{sequence of orders}.
For $1 \leq l\leq d$ let $H_l = \upsilon_{d-l+1}\qat \oplus \dotsc \oplus \upsilon_d\qat$.
This is a $l$\=/dimensional subspace of $H$.
Every section $\xi \in H_l \setminus \{0\}$ has a root of order $\geq d-l$ at $x$, and therefore $\ord H_l \geq l(d-l)$ holds.
Due to the Plücker formula we have for all $l\in\{1,\ldots,d\}$ the inequality
\[\frac{\|\V\|_2^2}{\pi}\geq l(d-1-\deg(E)-(l-1)\genus):=p(l)\]
The quadratic polynomial $p(l)$ has a maximum at $l_\text{\scriptsize max}=\frac{d-1-\deg(E)}{2\genus} + \tfrac12$.
If $l_\text{\scriptsize max} \leq \tfrac12$, then $p(d) < \dotsc < p(1) \leq p(0)=0$, and therefore we do not obtain a non-trivial bound.
On the other hand, if $l_\text{\scriptsize max} > \tfrac12$, equivalently $\deg(E) < d-1$, then the interval $[l_\text{\scriptsize max}-\frac12,l_\text{\scriptsize max}+\frac12]$ contains a positive integer $l_\text{\scriptsize opt}$.
If $l_\text{\scriptsize max}+\frac12\le d$, equivalently $\deg(E)\ge(d-1)(1-2\genus)$, then $l_\text{\scriptsize opt} \in \{1,\ldots,d\}$ with $p(l_\text{\scriptsize opt}) = \lceil p(l_\text{\scriptsize max} \pm \tfrac12) \rceil$ because $p(l)$ is a monic quadratic polynomial.
Therefore
\[
\frac{\|\V\|_2^2}{\pi}\ge\left\lceil\frac{(d-1+\genus-\deg(E))^2-\genus^2}{4\genus}\right\rceil\text{ for }(d-1)(1-2\genus)\le\deg(E)<(d-1).
\]
For $l_\text{\scriptsize max}+\frac12>d$, equivalently $\deg(E)<(d-1)(1-2\genus)$, the maximum of $p(l)$ on $\{1,\ldots,d\}$ is $p(d)=d((d-1)(1-\genus)-\deg(E))$.
\end{proof}

Finally we generalize to admissible maps $F:\X\to\qat$ a result of Li and Yau~\cite{LY} about a lower bound of the Willmore functional~\eqref{eq:willmore energy weierstrass}.

\begin{corollary}
\label{thm:embedding}
\index{Willmore energy}
\index{Underlying bundle}
If $F:\X\to\qat$ is admissible on a compact Riemann surface $\X$ with local Weierstraß representation~\eqref{eq:weierstrass}, then the Willmore functional obeys
\begin{align*}
\willmore(F)&\geq4\pi\sum_{x\in F^{-1}(\{r\})}\left(1+\ord_x(\chi)+\ord_x(\psi)\right)
\quad\text{for all }r\in\qat.
\end{align*}
In particular $\willmore(F)\geq4\pi$.
For $\willmore(F)<8\pi$, the map $F$ is injective and $dF$ has no roots.
Furthermore, for maps obeying~\eqref{eq:imaginary 2} and $\willmore(F)<12\pi$, again $dF$ has no roots and $E^2=T\X$.
In general the degree of the holomorphic line bundle $E$ for the Kodaira and the Weierstraß representation is bounded by
\[
2-2g-\frac{\willmore(F)}{4\pi}\le\deg(E)\le\frac{\willmore(F)}{4\pi}.
\]
\end{corollary}
\begin{proof}
In the proof of Theorem~\ref{thm:weierstrass} we constructed two $\V$\=/holomorphic sections $\upsilon$ and $\phi=\upsilon F$ of a holomorphic $\qat$\=/line bundle $E_\qat$, whose potential $\V$ obeys~\eqref{eq:hopf norms}.
Now~\eqref{eq:pluecker} implies  $\willmore(F)=4\|\U\|_2^2 = 4\|\V\|_2^2+4\pi \deg(E) \geq4\pi\ord(H)$ for any one-dimensional subspace $H$ of $\upsilon\qat\oplus\phi\qat$.
Hence for all $r\in\qat$ the number of roots of $\phi-\upsilon r=\upsilon(F-r)$ counted with multiplicities is a lower bound of $\frac{1}{4\pi}\willmore(F)$.
By Lemma~\ref{lem:change of order}(i) this order is the sum on the right-hand side of the inequality in the corollary.
If $\willmore(F)<8\pi$ then this sum is equal to $1$ for every $r\in \qat$, hence $F$ is injective and $dF=\lp\chi,\psi\rp$ has no roots.
If $F$ obeys~\eqref{eq:imaginary 2}, then by Lemma~\ref{lem:admissible imaginary} $dF$ can only have roots of even order, and therefore no roots for $\willmore(F)<12\pi$.

We prove the last statement separately for the Kodaira and the Weierstraß representation.
For the Kodaira representation we have $\willmore(F)=4\|\V\|_2^2+4\pi\deg(E)$ and $\dim H^0(\X,\Q{E,\V})\ge 2$.
The first formula implies the upper bound for $\deg(E)$ and the Plücker formula the lower bound:
\[
\frac{\willmore(F)}{4\pi}\ge2(1-\genus-\deg(E))+\deg(E)=2-2\genus-\deg(E).
\]
In terms of the Weierstraß potential the Willmore functional is equal to $\willmore(F)=4\|\U\|_2^2$ and both spaces $H^0(\X,\Q{E^{-1},\U\sd})$ and $H^0(\X,\Q{KE,\U})$ are at least one-dimen\-sional.
The Plücker formula for the first space implies the upper bound on $\deg(E)$ and for the second space the lower bound.
\end{proof}

\begin{example}[Holomorphic coverings]
Let $\X$ be a branched cover of $\mathbb{P}^1$ of degree $d$.
Then, as in Example~\ref{eg:kodaira covering}, we can compose the covering map $p$ with the round sphere map $F : \mathbb{P}^1 \to \qat$ of Examples~\ref{eg:sphere kodaira} and~\ref{eg:willmore energy}.
This gives us an admissible map $\Tilde{F} = F \circ p : \X \to \qat$ with Willmore energy $4\pi d$.
On the other hand, we know that $\deg E = 1$ and so the pullbacked line bundle $\Tilde{E} = p^\ast E$ will have degree $d$.
The bound in Corollary~\ref{thm:embedding} reads $2 - 2g - d \leq d \leq d$.
This shows the sharpness of the upper bound.

Likewise, if we now reflect $\Tilde{F}$ as in Example~\ref{eg: - bar F}, since the round sphere lies in $\mathbb{R}^3$, it is unchanged, and in particular its Willmore energy is unchanged.
But we produce a new Weierstraß representation $-d\overline{\Tilde{F}} = \lp \psi, \chi \rp$ with line bundle $(K_\X\, p^\ast E)^{-1}$.
Because $\psi$ has roots at the ramification points of the covering map, this is not the form of the Weierstraß representation we usually consider in this book, however the above corollary applies to this case.
The bound in Corollary~\ref{thm:embedding} reads $2 - 2g - d \leq -(2g-2) - d \leq d$, showing sharpness from below.
\end{example}

\begin{example}[Inverted Catenoid]
\label{eg:inverted catenoid not injective}
\index{Catenoid!Inverted catenoid}
Recall the compactified inverted catenoid from Example~\ref{eg:inverted catenoid} and that it's Willmore energy is $8\pi$ from Example~\ref{eg:willmore energy}.
In particular, it maps both $z = 0$ and $\infty$ to $\qi \in \qat$, so $\tilde{F}$ is not injective.
This shows the $8\pi$ bound for injectivity in Corollary~\ref{thm:embedding} is sharp.
\end{example}

\chapter{Serre Duality and Riemann-Roch Theorem}
\label{chapter:riemann roch}
In this chapter, for compact Riemann surfaces we prove Serre duality and the Riemann-Roch theorem for all holomorphic $\qat$\=/line bundles with square-integrable potentials.
For this purpose, for any holomorphic $\mathbb{C}$\=/line bundle $E$ we introduce a Hermitian pairing between $(1,0)$\=/forms with values in $E_\qat$ and sections of $KE^{-1}_\qat$.
This pairing is defined in such a way that the two holomorphic structures $\delbar{KE}-\V$ and $\delbar{E^{-1}}-\V\sd$ of the two paired holomorphic line bundles are the adjoint of one other.
If we conceive of these holomorphic structures with the aid of elliptic theory as Fredholm operators between appropriate Banach spaces, then these pairings directly imply Serre duality.
Afterwards the Riemann-Roch theorem follows by standard arguments from sheaf theory.
As an application, in the situation that $\X$ is an elliptic curve and the Willmore energy is less than $8\pi$, we further restrict the possible values $\deg E\in\{-1,0,1\}$ in Corollary~\ref{thm:embedding} to $\deg E=0$

We consider only compact Riemann surfaces $\X$ in this chapter.
Recall our standard set-up as in Remark~\ref{rem:special cover} of an open cover $\SO_1\cup\ldots\cup\SO_L$ of $\X$ with charts $z_l:\SO_l\to\Omega_l\subset\mathbb{C}$ and smooth boundaries $\boundary \Omega_l$.
Hence any holomorphic $\mathbb{C}$\=/line bundle $E$ is described by the holomorphic cocycle $f_{ml}$ with respect to this cover.

Due to Lemma~\ref{lem:10 pairing properties}(vi), the holomorphic structure $\delbar{E} -\V$ is in the sense of Chapter~\ref{chapter:weierstrass} the paired holomorphic structure of $\delbar{KE^{-1}} - \V\sd$ on the $\qat$\=/line bundle $KE^{-1}_\qat$.
We shall show now that the corresponding operators are in a sense the adjoints of each other.
We supplement the two pairings $\lz\cdot,\cdot\rz$ and $\lp\cdot,\cdot\rp$ by another sesquilinear pairing $\lh\xi,\omega\rh$ with values in the $2$\=/forms between sections $\xi$ of $KE^{-1}_\qat$ and $(0,1)$\=/forms $\omega$ with values in $E_\qat$.
Its integral $\int_\X\lh\xi,\omega\rh$ is a pairing between $\xi\in H^0(\X,\ban{\frac{p}{p-1}}{KE^{-1}})$ and $\omega\in H^0(\X,\forms{0,1}{}\ban{p}{E})$.
\begin{definition}
\label{def:11 pairing}
\index{Pairing!(1,1)}
Take any sections $\xi$ of $KE^{-1}_\qat$ and $\omega$ of $\forms{0,1}{E}$, represented by $\xi_l$ and $d\Bar{z}_l\omega_l$.
We define the pairing $\lh\xi,\omega\rh$ to be the $\qat$\=/valued $2$\=/form on $\X$ with local representatives 
\begin{align*}
\lh\xi,\omega\rh_l&=-\Bar{\xi}_l\qj dz_l \wedge d\Bar{z}_l\omega_l=\Bar{\xi}_l\qj d\Bar{z}_l\wedge dz_l\omega_l.
\end{align*}
\end{definition}
It should behave therefore somewhat like a $(1,1)$\=/form, hence its name.
Its true utility however is that it gives a non-degenerate pairing between spaces of sections.
We should first check that this formula actually defines a global $2$\=/form.
On $\SO_l\cap\SO_m$ it has indeed the correct transformation behavior:
\begin{align*}
-\Bar{\xi}_m\qj dz_m\wedge d\Bar{z}_m\omega_m
&=-\overline{\tfrac{dz_l}{dz_m}f^{-1}_{ml}\xi_l}\qj dz_m\wedge\tfrac{d\Bar{z}_m}{d\Bar{z}_l}f_{ml}d\Bar{z}_l\omega_l
=-\Bar{\xi}_l\qj dz_l\wedge d\Bar{z}_l\omega_l.
\end{align*}
This form is closely related to the exterior derivative of $\lp\cdot,\cdot\rp$ from Definition~\ref{def:10 pairing}.
For differentiable sections $\xi$ of $KE^{-1}_\qat$ and $\upsilon$ of $E$ we have
\begin{gather}
\label{eq:d of 01 pairing}
\begin{aligned}
d\lp\xi,\upsilon\rp_l
&=d(\Bar{\xi}_l\qj dz_l\upsilon_l)=\overline{d\xi_l}\wedge\qj dz_l\upsilon_l-\Bar{\xi}_l\qj dz_l\wedge d\upsilon_l\\
&=\overline{(dz_l\partial_l\xi_l+d\Bar{z}_l\barpartial_l\xi_l)}\wedge\qj dz_l\upsilon_l-\Bar{\xi}_l\qj dz_l\wedge d\Bar{z}_l\barpartial_l\upsilon_l\\
&=\overline{\barpartial_l\xi_l}\qj d\Bar{z}_l\wedge dz_l\upsilon_l+\lh\xi,\delbar{E}\upsilon\rh_l\\
&=-\overline{\Bar{\upsilon}_l\qj d\Bar{z}_l\wedge dz_l\qj\barpartial_l\xi_l}+\lh\xi,\delbar{E}\upsilon\rh_l\\
&=-\overline{\lh\upsilon,\delbar{KE^{-1}}\xi\rh_l}+\lh\xi,\delbar{E}\upsilon\rh_l.
\end{aligned}
\end{gather}
Here we used that $\qj d\Bar{z}_l\wedge dz_l=-2\qk \dmu_l$ is a $\Imag\qat$\=/valued $2$\=/form.
Furthermore, the $2$\=/form-valued pairing $\lh\upsilon,\delbar{KE^{-1}}\xi\rh$ is obtained by the application of Definition~\ref{def:11 pairing} to the line bundle $KE^{-1}$ instead of $E$ and $K(KE^{-1})^{-1}=E$ instead of $KE^{-1}$.
\begin{lemma}\phantom{boo}\label{lem:11 pairing properties}
\begin{enumeratethm}
\item
For functions $g:\X\to\mathbb{C}$ the pairing obeys
\begin{align*}
\lh g\xi,\omega\rh &= \lh\xi,g\omega\rh.
\end{align*}
\item 
For functions $h:\X\to\qat$ the pairing obeys
\begin{gather*}
\lh\xi h,\omega\rh=\Bar{h}\lh\xi,\omega\rh
\quad\text{and}\quad
\lh\xi,\omega h\rh=\lh\xi,\omega\rh h.
\end{gather*}

\item
The pairing $(\xi,\omega)\mapsto\int_\X\lh\xi,\omega\rh$ is a non-degenerate sesquilinear pairing between $\xi\in H^0(\X,\ban{\frac{p}{p-1}}{KE^{-1}})$ and $\omega\in H^0(\X,\forms{0,1}{}\ban{p}{E})$.

\item
If $\xi\in H^0(\X,\sob{1,\frac{2p}{3p-2}}{KE^{-1}})$, $\upsilon\in H^0(\X,\sob{1,p}{E})$ and $\V\in\pot{E}$, then (compare~\cite[Theorem~4.2]{PP})
\begin{gather*}
d\lp\xi,\upsilon\rp=-\overline{\lh\upsilon,\delbar{KE^{-1}}\xi\rh}+\lh\xi,\delbar{E}\upsilon\rh
\quad\text{and}\quad
\overline{\lh\upsilon,\V\sd\xi\rh}=\lh\xi,\V\upsilon\rh.
\end{gather*}
\end{enumeratethm}
\end{lemma}
\begin{proof}
Similar to Lemma~\ref{lem:10 pairing properties}, the first two properties follow algebraically from the definition.
Thus (iii) has the behavior of a sesquilinear pairing (conjugate-linear in the first argument, linear in the second); we just need to show it is finite and non-degenerate.
The dual space of $\banach{p}(\Omega_l,\qat)$ is $\banach{p/(p-1)}(\Omega_l,\qat)$,
so the product is integrable on $\SO_l$.
Since $\X$ is compact and covered by finitely many sets $\SO_l$, the integrals of the $2$\=/forms are finite.
The non-degeneracy of these pairings follows because $d\Bar{z}_l\wedge dz_l$ is equal to $2\ci\dmu$ on $\Omega_l$, for $\dmu$ the real Lebesgue measure.

That leaves (iv) to be proven.
The potential $\V\sd$ is defined in Lemma~\ref{lem:10 pairing properties}(vi), where its well-definedness as a global potential was verified.
By the Sobolev Embedding theorem~\ref{lem:sobolev regularity} the sheaves $\sob{1,p}{E}$ and $\sob{1,\frac{2p}{3p-2}}{KE^{-1}}$ are subsheaves of $\ban{\frac{2p}{2-p}}{E}$ and $\ban{\frac{p}{p-1}}{KE^{-1}}$, respectively.
So all terms are locally represented by functions in $\banach{1}\loc(\Omega_l,\qat)$.
Equation~\eqref{eq:d of 01 pairing} shows the equation on the left hand side.
The equation on the right hand side follows by inserting $\qj dz_l\wedge d\Bar{z}_l=2\qk \dmu_l$ locally on $\SO_l$:
\begin{gather*}
\hspace{-5mm}\overline{\lh\upsilon,\V\sd\xi\rh}_l=-2\overline{\Bar{\upsilon}_l\qk \dmu_l\qk\Bar{\V}_l\qk\xi_l}=2\Bar{\xi}_l\qk\V_l\qk^2 \dmu_l\upsilon_l=-2\Bar{\xi}_l\qk \dmu_l\V_l\upsilon_l=\lh\xi,\V\upsilon\rh_l.\qedhere
\end{gather*}
\end{proof}
The pairing in (iii) should be thought of as a generalization of a Hermitian pairing, whose arguments belong to different spaces.
This could be formalized in the language of \emph{dual systems}, but we will not do so.
In combination with (iv), we should think of the operators
\begin{align}
\label{fredholm 1}
\delbar{E}-\V&:&H^0(\X,\sob{1,p}{E})&\to H^0(\X,\forms{0,1}{}\ban{p}{E}),\\
\label{fredholm 2}
\delbar{KE^{-1}}-\V\sd&:&H^0(\X,\ban{p/(p-1)}{KE^{-1}})&\to H^0(\X,\forms{0,1}{}\sob{-1,p/(p-1)}{KE^{-1}}),
\end{align}
as adjoint to one another.
This is because the integral of $d\lp\chi,\psi\rp$ over $\X$, which is compact, is zero and so for sections meeting the assumptions of (iv) we have
\begin{align}
\label{eq:integration by parts}
\int_\X \overline{\lh\upsilon,(\delbar{KE^{-1}}-\V\sd)\xi\rh}
= \int_\X \lh\xi,(\delbar{E}-\V)\upsilon\rh.
\end{align}
We also refer to this as ``integration by parts.''
From~\eqref{fredholm 2} we see that the adjoint is a paired holomorphic structure on the $\qat$\=/line bundle $KE^{-1}_\qat$.
Furthermore, the relation in Lemma~\ref{lem:11 pairing properties}(iv) is a quantitative version of the statement~(vii) in Lemma~\ref{lem:10 pairing properties}. 

We observe at this point that if we have paired holomorphic bundles then the simultaneous application of Corollary~\ref{gauge holomorphic structure} results again in paired bundles.
If we decompose $\V=\V^++\V^-$ into $\V^+\in\pot{E}^+$ and $\V^-\in\pot{E}^-$, then we transform in Corollary~\ref{gauge holomorphic structure} the pair $(E,\V)$ into another pair $(E',\V')$, such that $\V'\in\pot{E'}^-$ and $|\V'|^2=|\V^-|^2$ in the sense of~\eqref{eq:square integrable}.
Furthermore the sheaf $\Q{E,\V}$ is isomorphic to $\Q{E',\V'}$ and $\Q{KE^{-1},\V\sd}$ is isomorphic to the paired sheaf $\Q{KE'^{-1},{\V'}\sd}$, since $(\V^+)\sd=-\V^+$ and the cocycle induced by $-\V^+$ is the inverse of the cocycle induced by $\V^+$.

Since the sheaf $\sob{1,p}{E}$ is invariant under multiplication with a smooth partition of unity standard arguments~\cite[Theorem~12.6.]{Fo} show that the first cohomology group of this sheaf is trivial.
Consequently, the long exact cohomology sequence~\cite[\S15.]{Fo} corresponding to the short exact sequence~\eqref{eq:exact sequence},
\[
0\to H^0(\X,\Q{E,\V})\hookrightarrow H^0(\X,\sob{1,p}{E})\xrightarrow{\delbar{E}-\V}H^0(\X,\forms{0,1}{}\ban{p}{E})\to H^1(\X,\Q{E,\V})\to 0,
\]
shows that the cokernel of~\eqref{fredholm 1} is isomorphic to $H^1(\X,\Q{E,\V})$.

Since the pairing in Lemma~\ref{lem:11 pairing properties}(iii) is non-degenerate, this cokernel is also dual to the kernel of~\eqref{fredholm 2}.
Since $\frac{p}{p-1}>2$ the equivalence of~(i) and (iv) in Theorem~\ref{cauchy formula} identifies this kernel of~\eqref{fredholm 2} with $H^0(\X,\Q{KE^{-1},\V\sd})$.
The isomorphisms $H^0(\X,\Q{E,\V})\cong H^0(\X,\Q{E',\V'})$ and  $H^0(\X,\Q{KE^{-1},\V\sd})\cong H^0(\X,\Q{KE'^{-1},{\V'}\sd})$ obtained by Corollary~\ref{gauge holomorphic structure} show with Corollary~\ref{dimension bound} that these four vector spaces are finite-dimensional.
So~\eqref{fredholm 1} is a Fredholm operator with kernel $H^0(\X,\Q{E,\V})$ and with cokernel dual to $H^0(\X,\Q{KE^{-1},\V\sd})$.
This proves (compare~\cite[\S8.-\S9.]{Na})
\begin{theorem}[Serre Duality]
\label{Serre duality}
\index{Serre duality}
Let $E$ be a holomorphic $\mathbb{C}$\=/line bundle on a compact Riemann surface $\X$ and $\V\in\pot{E}$.
Then $\V\sd\in\pot{KE^{-1}}$ and the \v{C}ech cohomology group $H^1(\X,\Q{E,\V})$ is dual to $H^0(\X,\Q{KE^{-1},\V\sd})$.
\qed
\end{theorem}
\begin{theorem}[Riemann-Roch Theorem]
\label{riemann roch}
\index{Riemann-Roch theorem}
Let $E$ be a holomorphic $\mathbb{C}$\=/line bundle on a compact Riemann surface $\X$ and $\V\in\pot{E}$.
Then the quaternionic dimensions of the corresponding \v{C}ech cohomology groups are finite and obey
\[\dim_{\qat}H^0(\X,\Q{E,\V})-\dim_{\qat}H^1(\X,\Q{E,\V})=1-\genus+\deg(E).\]
\end{theorem}
\begin{proof}
We already showed that~\eqref{fredholm 1} is a Fredholm operator depending continuously on the potential $\V\in\pot{E}$.
Since the index depends continuously on the operator~\cite[Chapter~VII Theorem~4.]{Pa} it is equal to the index of the corresponding Fredholm operator with trivial potential $\V=0$.
Due to the classical Riemann-Roch theorem this index is equal to $1-\genus+\deg(E)$.
\end{proof}

In Chapter~\ref{chapter:singular holomorphic sheaves} we utilize the following description of $H^1(\X,\Q{E,\V})$ together with an explicit formula for the Serre duality between elements of this space with elements of $H^0(\X,\Q{KE^{-1},\V\sd})$. 
Let $E$ be a holomorphic $\mathbb{C}$\=/line bundle on a compact Riemann surface $\X$ with finitely many pairwise distinct points $x_1,\ldots,x_L$ and $\V\in\pot{E}$. 
Moreover, let $\Sh{M}_{E,\V}$ denote the sheaf of $\V$\=/meromorphic sections of $E_\qat$ and $\Q{E,\V,x_1,\ldots,x_L}\subset\Sh{M}_{E,\V,x_1,\ldots,x_L}$ the cartesian product of the stalks of the sheaves $\Q{E,\V}$ and $\Sh{M}_{E,\V}$ at the points $x_1,\ldots,x_L$, respectively. 
Finally, let $\iota$ denote
\begin{align}\label{eq:embedding iota}
\iota:H^0(\X,\Sh{M}_{E,\V})\cap H^0(\X\setminus x_1,\ldots,x_L\},\Q{E,\V})&\hookrightarrow\Sh{M}_{E,\V,x_1,\ldots,x_L}
\end{align}
the natural embedding of the global sections of $\Q{E,\V}$ with possible poles of finite order at the points $x_1,\ldots,x_L$ into $\Sh{M}_{E,\V,x_1,\ldots,x_L}$. It maps any global section to the cartesian product of the corresponding germs at the points $x_1,\ldots,x_L$. Now we can state and prove the following lemma (compare~\cite[17.3~Theorem]{Fo}):
\begin{lemma}\label{lem:serre duality in terms of residue}
Let on a compact Riemann surface $\X$ the following data be given: finitely many pairwise distinct points $x_1,\dots,x_L\in\X$, a holomorphic $\mathbb{C}$\=/line bundle $E$ and $\V\in\pot{E}$. 
Furthermore, let $\Q{E,\V,x_1,\ldots,x_L}\subset\Sh{M}_{E,\V,x_1,\ldots,x_L}$ denote the cartesian products of the corresponding stalks and $\iota$ the embedding~\eqref{eq:embedding iota} as defined above. 
Then $H^1(\X,\Q{E,\V})$ is naturally isomorphic to the following quotient space
\begin{align}\label{eq:quotient space}
\Sh{M}_{E,\V,x_1,\ldots,x_L}\!\Big/\hspace{-1mm}\Big(\!\Q{E,\V,x_1,\ldots,x_L}\hspace{-1mm}+\!\iota\big[H^0(\X,\Sh{M}_{E,\V})\!\cap\!H^0(\X\!\setminus\!\{x_1,\ldots,x_L\},\Q{E,\V})\big]\!\Big)\!.\hspace{-1.2mm}
\end{align}
Furthermore, for any $\eta$ in this quotient space and any $\xi\in H^0(\X,\Q{KE^{-1},\V\sd})$ the Serre duality between $\eta$ and $\xi$ is equal to $2\pi\qi\sum_{l=1}^L\res_{x_l}\lp\eta_l,\xi\rp$.
\end{lemma}
\begin{proof}
We first construct an embedding of the quotient space~\eqref{eq:quotient space} into $H^1(\X,\Q{E,\V})$. Afterwards we prove the formula for the Serre duality and utilize this formula in order to prove the surjectivity of the embedding.

Let $\eta=(\eta_1,\ldots,\eta_L)\in\Sh{M}_{E,\V,x_1,\ldots,x_L}$ be given. 
In the sequel we make statements about the indices $\{1,\ldots,L\}$ which should hold for any $l\in\{1,\ldots,L\}$. 
We choose small disjoint open discs $\SU_l$ with smooth boundaries such that $\eta_l$ extends to a section of $\Q{E,\V}$ on $\SU_l\setminus\{x_l\}$ and supplement $\SU_1,\ldots,\SU_L$ by $\SU_0=\X\setminus\{x_1,\ldots,x_L\}$ to an open cover of $\X$. 
Any three of its sets are disjoint. 
Hence $\eta$ defines a cocycle of $\Q{E,\V}$ with respect to this cover, and therefore by~\cite[12.5 Definition]{Fo} an element of $H^1(\X,\Q{E,\V})$. 
Due to~\cite[12.3 Lemma]{Fo} this element of $H^1(\X,\Q{E,\V})$ does not depend on the choice of $\SU_1,\ldots,\SU_L$. 
The coboundaries are the elements of the denominator in~\eqref{eq:quotient space}. 
Therefore this construction induces a map of\eqref{eq:quotient space} to $H^1(\X,\Q{E,\V})$.
If $\eta$ is mapped to zero in $H^1(\X,\Q{E,\V})$ then in the above construction we know it is a coboundary $\eta = \delta^\ast(f)$ due to~\cite[12.4 Lemma]{Fo}.
In particular, at $x_l$ we have $f_0 = f_l - \eta_l$, showing that $f_0$ belongs to $H^0(\X,\Sh{M}_{E,\V})\!\cap\!H^0(\X\!\setminus\!\{x_1,\ldots,x_L\},\Q{E,\V})$, and therefore $\eta_l = f_l - f_0$ belongs to the denominator of~\eqref{eq:quotient space}.
In summary, we have an injective map from~\eqref{eq:quotient space} to $H^1(\X,\Q{E,\V})$.

In order to prove the formula for the Serre duality between the given $\eta$ we first construct an element of $H^0(\X,\forms{0,1}{}\ban{p}{E})$ that maps to $\eta$ under the connecting homomorphism of~\eqref{eq:exact sequence}.
Choose a $h_l\in C^\infty_0(\SU_l,[0,1])$ that is identically equal to $1$ on some neighborhood of $x_l$. 
We diminish $\SU_0$ such that each $h_l$ vanishes on $\SU_0\cap\SU_l$.
Due to~\cite[12.3 Lemma and 12.4 Lemma]{Fo} the restriction of $\eta$ to this new cover represents the same element of $H^1(\X,\Q{E,\V})$. 
Consider $(1-h_l)\eta_l|_{\SU_l\setminus\{x_l\}}$ and extend it by zero at $x_l$.
Together with $0$ on $\SU_0$ this gives a cochain in $C^0(\mathcal{U},\sob{1,p}{E})$.
Finally, the application of $\delbar{E}-\V$ gives a global section of $\forms{0,1}{}\ban{p}{E}$.
By definition of the connecting homomorphism this element of $H^0(\X,\forms{0,1}{}\ban{p}{E})$ is mapped to $\eta$.
Due to Lemma~\ref{lem:11 pairing properties}(iv) the Serre duality of $\eta$ and $\xi$ is equal to
\[
\sum_{l=1}^L\int\limits_{\SU_l}\hspace{-1.5mm}\big\lh(\delbar{E}-\V)\big((1-h_l)\eta_l\big),\xi\big\rh
=\sum_{l=1}^L\int\limits_{\SU_l}\hspace{-1.5mm}d\big\lp(1-h_l)\eta_l,\xi\rp
=2\pi\qi\sum_{l=1}^L\res_{x_l}\lp\eta_l,\xi\rp.
\]
The last equality follows by Stokes theorem, since it is equal the integral of $\lp(1-h_l)\eta_l,\xi\rp = \lp\eta_l,\xi\rp$ along the boundary of the disc $\SU_l$ (since $\boundary\SU_l = \SU_0$, where $h_l$ vanishes), and by Lemma~\ref{lem:10 pairing properties}(viii) equal to $2\pi\qi\res_{x_l}\lp\eta_l,\xi\rp$.

To finish the proof we need to show that the embedding from the quotient~\eqref{eq:quotient space} into $H^1(\X,\Q{E,\V})$ is surjective. By Lemma~\ref{quotient dimension}, for any $\xi\in H^0(\X,\Q{KE^{-1},\V\sd})$ we might choose $\eta_1\in\Sh{M}_{E,\V,x_1}$ such that $\lp\eta_1,\xi\rp$ has a first order pole at $x_1$. By Lemma~\ref{lem:change of order}(i) the pairing of $\eta=(\eta_1,0,\ldots,0)\in\Sh{M}_{E,\V,x_1,\ldots,x_L}$ with $\xi$ does not vanish. This implies that the orthogonal complement in $H^0(\X,\Q{KE^{-1},\V\sd})$ of the image of the map from the quotient~\eqref{eq:quotient space} into $H^1(\X,\Q{E,\V})$ with respect to the pairing in Lemma~\ref{lem:serre duality in terms of residue} is trivial. Hence by Serre duality this map is surjective.
\end{proof}

\begin{corollary}
\label{cor:bundle degree vanished}
For all admissible maps $F:\X\to\qat$ on an elliptic curve $\X$ with $\willmore(F)<8\pi$ the degree of the holomorphic line bundle $E$ for the Kodaira and the Weierstra{\ss} representations vanishes.
\end{corollary}
\begin{proof}
The bounds in Corollary~\ref{thm:embedding} give $\deg(E)\in\{-1,0,1\}$ for $\genus=1$ and $\willmore(F)<8\pi$.
We first exclude $\deg(E)=1$ by contradiction for the Kodaira representation.
From $\dim H^0(\X,\Q{E,\V})\ge 2$, the Riemann-Roch theorem, and Serre duality, we obtain $\dim H^1(\X,\Q{E,\V})=\dim H^0(\X,\Q{KE^{-1},\V\sd})\ge 1$ with $\deg(KE^{-1})=-1$.
However, the Plücker formula then implies the contradiction $8\pi \leq 4\|\V\|_2^2+4\pi\deg(E) = \willmore(F) < 8\pi$.

By  Corollary~\ref{thm:embedding} $dF$ and the holomorphic sections $\chi$ and $\psi$ of the Weierstraß representation have no roots.
So the holomorphic line bundle $E$ is unique and the same as in the Kodaira representation.
Interchanging $\chi$ and $\psi$ in $dF=\lp\chi,\psi\rp$ yields the Weierstraß representation of $-\Bar{F}$ (Example~\ref{eg: - bar F}).
For $\deg(E)=-1$ the Kodaira representation of $-\Bar{F}$ gives two holomorphic sections of $(KE)^{-1}_\qat$ with $\deg((KE)^{-1})=1$ and $\willmore(-\Bar{F})=\willmore(F)$.
So the arguments for $\deg(E)=1$ exclude also this case.
Hence $\deg(E)=0$ for the Kodaira and the Weierstraß representations.
\end{proof}

\begin{example}[Inverted Catenoid]
\index{Catenoid!Inverted catenoid}
Recall the compactified inverted catenoid from Example~\ref{eg:inverted catenoid}.
Because $\Tilde{F}(0) = \Tilde{F}(\infty) = \qi \in \qat$, we can interpret $\Tilde{F}$ as a map from the sphere with $z = 0,\infty$ identified.
This singular curve is the nodal curve $y^2 = x^2(x+1)$, the stable limit in the space of elliptic curves.
Morally, we should be able to perturb the inverted catenoid to an admissible map on an elliptic curve (we do not consider limits of Riemann surfaces in this book, although we will consider limits of admissible maps on a fixed Riemann surface in the next part).
The above corollary suggests that any such perturbation is necessarily Willmore energy increasing, since $\deg(E) = 1$.
\end{example}

\partpage{Limits of Sequences of Admissible Maps}

\chapter{Resolvent of the Dirac Operator on \texorpdfstring{$\mathbb{P}^1$}{P1}}
\label{chapter:resolvents}
Lemma~\ref{weakly continuous} at the end of Chapter~\ref{chapter:local} shows that the map $\V\mapsto(\barpartial-\V)^{-1}$ is weakly continuous, if $\|\V\|_2$ is bounded by some given positive constant. 
In this chapter we improve this result:
Theorem~\ref{inverse holomorphic structure} states that on the compact Riemann surface $\mathbb{P}^1$ the former map is for all $\epsilon>0$ weakly continuous, if the restrictions to all $\epsilon$ balls $B(x,\epsilon)$ of $\mathbb{P}^1$ have norms $\|\V|_{B(x,\epsilon)}\|_2$ bounded by some given positive constant.
Any weakly convergent sequence $(\V_n)_{n\in\mathbb{N}} \rightharpoonup \V$ of potentials on $\mathbb{P}^1$ is bounded and has a subsequence such that the sequence of measures $|\V_n|^2$ on $\mathbb{P}^1$ converges weakly. 
The combination of Theorem~\ref{inverse holomorphic structure} and Lemma~\ref{bounded point measures} implies that $((\barpartial-\V_n)^{-1})_{n\in\mathbb{N}}$ converges, under two conditions: The holomorphic structure $\barpartial-\V$ is invertible and the measures of all points of $\mathbb{P}^1$ with respect to the weak limit of the sequence of measures $(|\V_n|^2)_{n\in\mathbb{N}}$ are smaller than a given positive constant. 
So the weak continuity of $\V\mapsto(\barpartial-\V)^{-1}$ can only fail in the case of massive concentration of energy. 
In the subsequent chapter we shall determine the optimal value of this threshold:
The weak continuity can only fail if there is a concentration at single points of energy exceeding the Willmore energy of a round sphere.

The proof of these results use techniques which are well known for the Laplace operators~\cite[Chapter~V.~\S3]{St}.
Namely, the resolvent $(\lambda\unity-\triangle)^{-1}$ of the free Laplace operator on $\mathbb{R}$ has for $\lambda>0$ an integral kernel, which has at the diagonal a singularity not depending on $\lambda$ and decays away from the diagonal like $\exp(-C\lambda|x-x'|)$ with $C>0$.
We shall show an analogous statement for $\barpartial-\V$ on $\mathbb{P}^1$.

To enact our strategy we will relate $\barpartial-\V$ to the Dirac operator.
The proper setting for a Dirac operator\index{Dirac operator} is a $\mathrm{Spin}$\=/bundle of an oriented manifold with Riemannian metric $g$~\cite[\S3.2]{Fr1}.
It is a first order differential operator from the bundle to itself.
On a Riemann surface $\X$ spin bundles are realized in~\cite[Section~5]{FLPP} as complex $\qat$\=/bundles $E_\qat$ whose underlying holomorphic $\mathbb{C}$\=/bundle obeys $E^2\cong K$.
Moreover the Dirac operator is the composition of the holomorphic structure $\delbar{E} : H^0(\X,\sob{1,p}{E})\to H^0(\X,\ban{p}{\bar{K}E})$ (here interpreting it as a map between sections, rather than a map to forms) with a bundle morphism that we will denote by $\diracJ$.
On $\X$, a Riemannian metric is equivalent to a Hermitian metric, which additionally induces a Hermitian metric on $E$.
This yields an isomorphism $\bar{K}E \cong \bar{E}^2E \cong \bar{E}\mathbb{C} = \bar{E}$, where $\bar{E}E \cong \mathbb{C}$ is induced by the Hermitian metric.
For quaternionic bundles, we have then $\bar{K}E_\qat \cong \bar{E}_\qat$, which we further compose with the `conjugation' map $\bar{E}_\qat \to E_\qat$ of Equation~\eqref{eq:quaternionic bar map}.
We call this composition $\diracJ : \bar{K}E_\qat \to E_\qat$.
The content of~\cite[Lemma~5.3]{FLPP} is that the (free) Dirac operator, often denoted by $\dirac$, is $\diracJ \delbar{E} : H^0(\X,\sob{1,p}{E})\to H^0(\X,\ban{p}{E})$.

We now specialize, for the remainder of the chapter, to $\X = \mathbb{P}^1$ and the unique line bundle $E$ of degree $-1$.
For the moment it suffices to equip $\mathbb{P}^1$ with any conformal metric.
As we will see in the lemma below, the free Dirac operator $\diracJ\delbar{E}$ is an unbounded closed operator on $H^0(\mathbb{P}^1,\ban{p}{E})$ with compact resolvent.
Due to Lemma~\ref{lem:fredholm}, the addition of a potential gives a unbounded closed operator of $H^0(\mathbb{P}^1,\ban{p}{E})$ on the same domain $H^0(\mathbb{P}^1,\sob{1,p}{E})$ for any $1<p<2$.
Because these are unbounded operators, it is preferable to work with their inverses.

\begin{definition}
\label{def:resolvent}
For $\ci\lambda$ in the resolvent set, we define the \emph{resolvent}\index{Resolvent}
\[
\Op{R}(\V,\ci\lambda)
:=\big(\iota \circ (\qi\lambda)\ri-\diracJ(\delbar{E}-\V)\big)^{-1}
\]
from $H^0(\mathbb{P}^1,\ban{p}{E})$ to $H^0(\mathbb{P}^1,\sob{1,p}{E})$.
Here $(\qi\lambda)\ri$ denotes right multiplication by $\qi\lambda\in\qat^+$ on $E_\qat$ and $\iota$ is the compact embedding $H^0(\mathbb{P}^1,\sob{1,p}{E})\hookrightarrow H^0(\mathbb{P}^1,\ban{q}{E})$.
\end{definition}

We relate different resolvents to each other by the following variant of the first resolvent formula~\cite[Theorem~VI.5]{RS1}:
\begin{align*}
\Op{R}(\W,\ci\mu)-\Op{R}(\V,\ci\lambda)
&= \Op{R}(\V,\ci\lambda)\big((\qi\lambda-\qi\mu)\ri + \diracJ\V - \diracJ\W\big)\Op{R}(\W,\ci\mu) \\
\Op{R}(\W,\ci\mu)
&= \Op{R}(\V,\ci\lambda) \left[\unity +\big((\qi\lambda-\qi\mu)\ri + \diracJ\V - \diracJ\W\big) \Op{R}(\W,\ci\mu) \right] \\
\Op{R}(\V,\ci\lambda)
&= \Op{R}(\W,\ci\mu) \left[\unity + \big((\qi\lambda-\ci\mu)\ri + \diracJ\V-\diracJ\W\big) \Op{R}(\W,\ci\mu) \right]^{-1}\hspace{-3mm}.
\labelthis{eq:first resolvent formula}
\end{align*}
This means, and indeed it will form a central part of our strategy, that it is often sufficient to understand the \emph{free resolvent} $\Op{R}(0,\ci\lambda)$.
In the following lemma we establish its existence.

\begin{lemma}\label{dirac resovent}
For $1<p<\infty$ and the holomorphic line bundle $E$ on $\mathbb{P}^1$ of degree $-1$ the free Dirac operator $\diracJ\delbar{E}$ is invertible in $\mathcal{L}(H^0(\mathbb{P}^1,\sob{1,p}{E}),H^0(\mathbb{P}^1,\ban{p}{E}))$.

Further, $\iota\circ(\qi\lambda)\ri-\diracJ\delbar{E}$ of $\mathcal{L}(H^0(\mathbb{P}^1,\sob{1,p}{E}),H^0(\mathbb{P}^1,\ban{p}{E}))$ is invertible for $\lambda\in\mathbb{R}$.
The compact operator $\iota\circ\Op{R}(0,\ci\lambda)\in\mathcal{L}(H^0(\mathbb{P}^1,\ban{p}{E}))$ is diagonalizable with eigenvalues and eigenfunctions not depending on $1<p<\infty$.
\end{lemma}
\begin{proof}
The operator $\diracJ$ is invertible since it is essentially multiplication with the conformal factor, which is a positive smooth function on a compact space and so its reciprocal is also a bounded smooth function. 
Hence the first statement follows if we show that $\delbar{E}$ is an isomorphism from $H^0(\mathbb{P}^1,\sob{1,p}{E})$ onto $H^0(\mathbb{P}^1,\forms{0,1}{}\ban{p}{E})$.
The argument in the paragraph between~\eqref{left inverse} and~\eqref{eq:neumann-convergence} which shows that $\Op{I}_{\Omega'}$ belongs to $\mathcal{L}(\banach{p}(\Omega',\qat),\sobolev{1,p}(\Omega',\qat))$ is valid for all $1<p<\infty$.
Hence for $\V=0$ the Lemma~\ref{lem:fredholm} holds for all such $1<p<\infty$, and $\delbar{E}$ is a Fredholm operator.
The kernel of this operator is $H^0(\mathbb{P}^1,\Q{E,0})$ which is trivial since $\deg(E)<0$. As explained in the proof of Serre duality~\ref{Serre duality}, the dimension of the cokernel is $\dim H^1(\mathbb{P}^1,\Q{E,0})$.
Hence, due to the Riemann-Roch Theorem~\ref{riemann roch} the index of the Fredholm operator is $1-\genus+\deg{E}=0$.
This proves the first statement.

Since $\iota\circ(\qi\lambda)\ri$ is compact for all $\lambda\in\mathbb{C}$, the operator $\iota\circ(\qi\lambda)\ri-\diracJ\delbar{E}$ is a Fredholm operator of index zero (compare \cite[Theorem~4.4.2~(i)]{BSa}).
\index{Index!Fredholm operator} 
Hence the second statement follows if the kernel is trivial. 
Moreover, the first statement implies that $\diracJ\delbar{E}\circ\iota$ on the domain $H^0(\mathbb{P}^1,\sob{1,p}{E})$ is an unbounded closed operator of $H^0(\mathbb{P}^1,\ban{p}{E})$ with compact resolvent.
Furthermore, due to~\cite[Theorem~6.1.15]{BSa}, the complement of the resolvent set contains only eigenvalues with finite-dimensional generalized eigenspaces and $0$ is not such an eigenvalue. 
For $p=2$ the Dirac operator is self adjoint (compare~\cite[Section~4.1 and 4.2.]{Fr1}). 
In this case $\diracJ\delbar{E}\circ\iota$ is diagonalizable and the resolvent set contains in addition to $\{0\}$ also $\ci\mathbb{R}$. 
Now the second and the third statement follow, if for any $1<p<\infty$ any generalized eigenfunction of $\diracJ\delbar{E}\circ\iota$ in $H^0(\mathbb{P}^1,\sob{1,p}{E})$ in fact belongs to $\cap_{1<p<\infty}H^0(\mathbb{P}^1,\sob{1,p}{E})$.

We consider for $\lambda\in\mathbb{C}$ elements $\xi\in H^0(\mathbb{P}^1,\sob{1,p}{E})$ in the kernel of $\iota\circ(\qi\lambda)\ri-\diracJ\delbar{E}$. 
Since $0$ belongs to the resolvent set of $\diracJ\delbar{E}$, we may assume $\lambda\ne0$. 
In this case the (weakly differentiable) section $\xi$ is a fixed point of $-(\diracJ\delbar{E})^{-1}(\qi\lambda)\ri$, and thus iteratively a fixed point of $\big(-(\diracJ\delbar{E})^{-1}(\qi\lambda)\ri\big)^k$ for any $k \in \mathbb{N}$. 
The first statement implies that $-(\diracJ\delbar{E})^{-1}(\qi\lambda)\ri$ improves the regularity. 
Due to the Sobolev embedding, Lemma~\ref{lem:sobolev regularity}, the space $H^0(\mathbb{P}^1,\sob{1,p}{E})$ is embedded for $1<p<2$ into $H^0(\mathbb{P}^1,\ban{2p/(2-p)}{E})$ and for $2<p<\infty$ into the continuous sections of $E_\qat$. 
Therefore $\xi\in\cap_{1<p<\infty}H^0(\mathbb{P}^1,\sob{1,p}{E})$. 
Since the generalized eigenfunctions belong to the kernel of higher powers of $\iota\circ(\qi\lambda)\ri-\diracJ\delbar{E}$, appropriate differences of such generalized eigenfunctions are fixed points of $-(\diracJ\delbar{E})^{-1}(\qi\lambda)\ri$. 
Inductively we conclude that all generalized eigenfunctions belong to $\cap_{1<p<\infty}H^0(\mathbb{P}^1,\sob{1,p}{E})$. 
This completes the proof.
\end{proof}


In the proof of Lemma~\ref{weakly continuous resolvent} we shall find an integral kernel for $\Op{R}(0,\ci\lambda)$ which allows to control the behavior of $\Op{R}(0,\ci\lambda)$ for large real $\lambda$.
In the case of the Laplace operator on $\mathbb{R}^n$ these integral kernels can be determined explicitly by using the invariance with respect to translations and rotations. 
In order to imitate this procedure on $\mathbb{P}^1$ we need a transitive action of the isometry group and a description of its action on $H^0(\mathbb{P}^1,\ban{p}{E})$.
To this end, we realize $\mathbb{P}^1$ as the homogeneous Riemannian manifold $\mathbb{S}^3/\mathbb{S}^1$, where $\mathbb{S}^3$ are the unit quaternions and $\mathbb{S}^1 = \{ e^{\ci\varphi} \mid \varphi\in\mathbb{R}/2\pi\mathbb{Z}\big\} \subset \qat^+$.
The canonical projection $\Op{P}$, called the Hopf fibration, has the simple expression $\alpha = \alpha_1 + \qj \alpha_2 \mapsto [\alpha_1\:\alpha_2]$ in homogeneous coordinates.
The standard metric on $\mathbb{S}^3 \subset \qat$ descends to the Fubini-Study metric on $\mathbb{P}^1$.
It can also be obtained as the metric of a sphere of radius $\frac{1}{2}$ in $\mathbb{R}^3$.
The unique spin bundle $E$ arises in this picture as the tautological bundle.
That is, to each point $\Op{P}(\alpha) \in \mathbb{P}^1$ we assign the complex line $\alpha \mathbb{C} \subset \qat$ to be its fiber.
In particular, by identifying $\qat \cong \mathbb{C}^2$ and restricting the standard Hermitian product, $E$ inherits a Hermitian metric $h_E$.
This identification makes left multiplication by unit quaternions $\mathbb{S}^3$ into the group $\operatorname{SU}(2)$, hence the action on $E$ is an isometry.

It is useful to put explicit formulas to the above geometric ideas.
Let us cover $\mathbb{P}^1$ by $\SO_1=\{[1\:y]\mid y\in \mathbb{C}\}$ and $\SO_2=\{[x\:1]\mid x\in \mathbb{C}\}$ with the charts $z_1:\SO_1\to \mathbb{C},[x\:y]\mapsto\frac{y}{x}$ and $z_2:\SO_2\to \mathbb{C},[x\:y]\mapsto\frac{x}{y}$.
Here $[x\:y]$ denotes `homogeneous coordinates', the equivalence class of the relation $(x,y)\sim(x',y')\Leftrightarrow xy'=x'y$.
The Fubini-Study metric is locally given by $(1+|z_l|^2)^{-2} \Real(dz_l\otimes d\bar{z}_l)$ (see~\cite[p.~31]{GrHa}).
The measure on $\mathbb{P}^1$ induced by the Fubini-Study metric is
\[
\dmu_{\mathbb{P}^1}(z_l)
= (1 + |z_l|^2)^{-2} \dmu(z_l).
\]
Here $\dmu(z_l)$ denotes the Lebesgue measure on $z_l\in\mathbb{C}$.
The unique spin bundle $E$ has cocycle $f_{21} = \qi z_2^{-1}$ and inherits the Hermitian metric given locally by 
\[
h_E(\xi,\xi')(z_l)
= (1 + |z_l|^2) \overline{\xi_l(z_l)} \xi'_l(z_l),
\]
which extends to a Hermitian pairing on $E_\qat$.
This means the length of a section
\[
|\xi|_E(z_l) = \sqrt{h_E(\xi,\xi)(z_l)}
\]
is a $\mathbb{R}$\=/valued function on $\mathbb{P}^1$.
Finally, we can express the free Dirac operator in these charts by
\begin{align}
\label{eq:Dirac in chart}
(\diracJ \delbar{E} \xi)_l 
= \qj (1 + |z_l|^2) \barpartial_l \xi_l.
\end{align}
Indeed we can check explicitly that $\diracJ \delbar{E} \xi$ is again a section of $E_\qat$:
\begin{align*}
(\diracJ \delbar{E} \xi)_2 
&= \qj (1+ |z_2|^2)\barpartial_2 \xi_2
= \qj (1+ z_2\bar{z}_2)\overline{\tfrac{dz_1}{dz_2}} \barpartial_1 (f_{21} \xi_1) \\
&= \qj (1+ z_2\Bar{z}_2) (-\Bar{z}_2^{-2})\qi z_2^{-1} \barpartial_1 \xi_1
= \qi z_2^{-1} \qj (z_1\Bar{z}_1 + 1) \barpartial_1 \xi_1
= f_{21}(\diracJ \delbar{E} \xi)_1 .
\end{align*}
For the description of the measurable section $\xi$ of $E_\qat$ up to null set, we can neglect the single point $[0:1]$ in $\mathbb{P}^1\setminus\SO_1$.
So the function $\xi_1$ on $\SO_1$ completely determines $\xi$.
In the sequel we will identify the objects on $\Hat{\SO}_1$ with the global objects and omit the subscripts.

The element $\alpha \in \mathbb{S}^3$ acts on $\mathbb{P}^1$ as
\[
\alpha [x\:y]
\sim (\alpha_1 + \qj \alpha_2) (x + \qj y)
= (\alpha_1 x - \Bar{\alpha}_2y) + \qj(\Bar{\alpha}_1 y + \alpha_2 x)
\sim [\alpha_1 x - \Bar{\alpha}_2y \: \Bar{\alpha}_1 y + \alpha_2 x ].
\]
This is action is a double cover of a subgroup of the Möbius group, the group of biholomorphic transformations of $\mathbb{P}^1$, as both $\alpha$ and $-\alpha$ act on a point $z$ identically.
In particular, in the $z_l$\=/coordinates
\[
\alpha z_1
= \frac{\Bar{\alpha}_1 z_1 + \alpha_2}{\alpha_1 - \Bar{\alpha}_2z_1},
\qquad
\alpha z_2
= \frac{\alpha_1z_2 - \Bar{\alpha}_2}{\Bar{\alpha}_1 + \alpha_2 z_2}.
\]
Given an action on the domain of a function, we have the standard left action on functions $(\alpha \cdot f)(z) = f(\alpha^{-1}z)$.
In general, this is not sufficient to define an action on a section of a bundle, because the two sides of the equation would belong to different fibers of the bundle.
However, the action of $\mathbb{S}^3$ on the tautological bundle $E$, via its action on $\qat$, does induce an action on the sections of $E_\qat$ (compare with the representation $\mathcal{P}^{1,-1}$ of the non-unitary series in~\cite[Chapter~II.\S4]{Kn}, or~\cite[\S2.2]{Rouviere2014}).

\begin{lemma}\label{moebius action}
The group $\mathbb{S}^3$ acts on sections of $E_\qat$ by
\begin{align}
\label{eq:SU2 action on sections}
(\alpha\cdot\xi)(z_1)
:= (\Bar{\alpha}_1 + \Bar{\alpha}_2 z_1)^{-1}\xi\left(\alpha^{-1} z_1\right)
=(\Bar{\alpha}_1 + \Bar{\alpha}_2 z_1)^{-1}\xi\left(\frac{\alpha_1z_1-\alpha_2}{\Bar{\alpha}_1 + \Bar{\alpha}_2z_1}\right).
\end{align}
Here $(\Bar{\alpha}_1 + \Bar{\alpha}_2 z_1)^{-1}$ denotes the left multiplication by this complex function of $z_1$.
This preserves the Hermitian pairing and length in the sense that
\begin{align}
h_E(\alpha\cdot\xi, \alpha\cdot\xi')(z_1)
= h_E(\xi,\xi')(\alpha^{-1}z_1), \text{ and }
|\alpha\cdot\xi|_E(z_1)
= |\xi|_E(\alpha^{-1}z_1).
\labelthis{eq:SU2 isometry section}
\end{align}
Further, the measure $\mu_{\mathbb{P}^1}$ is invariant under these transformations, and the Dirac operator commutes with them.
\end{lemma}
\begin{proof}
As discussed about, it is sufficient to describe a section except for a null set, so we may ignore that~\eqref{eq:SU2 action on sections} is not defined at $z_1 = \Bar{\alpha}_1/\Bar{\alpha}_2$.
If one wishes to give a well-defined formula at this point, the second chart $\SO_2$ must be used.
The crux of the calculation is that
\[
1 + |\alpha^{-1}z|^2
= \frac{|\Bar{\alpha}_1 + \Bar{\alpha}_2 z|^2 + |\Bar{\alpha}_1 + \Bar{\alpha}_2 z|^2}{|\Bar{\alpha}_1 + \Bar{\alpha}_2 z|^{2}}
= (1+|z|^2) |\Bar{\alpha}_1 + \Bar{\alpha}_2 z|^{-2}.
\]
It then follows that
\begin{align*}
h_E(\alpha\cdot\xi, \alpha\cdot\xi')(z_1)
&= (1+|z_1|^2) \overline{\xi_1(\alpha^{-1}z_1)} |\Bar{\alpha}_1 + \Bar{\alpha}_2 z_1|^{-2}\xi_1'(\alpha^{-1}z_1) \\
&= (1 + |\alpha^{-1}z_1|^2) \overline{\xi_1(\alpha^{-1}z_1)} \xi_1'(\alpha^{-1}z_1) 
= h_E(\xi,\xi')(\alpha^{-1}z_1).
\end{align*}
and consequently
\[
|\alpha\cdot\xi|(z_1)
= \sqrt{h_E(\xi,\xi)(\alpha^{-1}z_1)}
= |\xi|(\alpha^{-1}z_1).
\]
It also shows that the invariance of the measure
\begin{align*}
\int_{\mathbb{P}^1} f(\alpha z_1)\,\dmu_{\mathbb{P}^1}(z_1)
&= \int_{\mathbb{C}} f(\alpha z_1) (1 + |z_1|^2)^{-2}\, \dmu(z_1) \\
&= \int_{\mathbb{C}} f(w_1) (1 + |\alpha^{-1} w_1|^2)^{-2}\, \frac{\ci}{2} d(\alpha^{-1} w_1) \wedge \overline{d(\alpha^{-1} w_1)} \\
&= \int_{\mathbb{C}} f(w_1) (1 + |\alpha^{-1} w_1|^2)^{-2}\, \frac{\ci}{2} \frac{dw_1}{(\Bar{\alpha}_1 + \Bar{\alpha}_2 w_1)^2} \wedge \overline{\frac{dw_1}{(\Bar{\alpha}_1 + \Bar{\alpha}_2 w_1)^2}} \\
&= \int_{\mathbb{P}^1} f(w_1) \,\dmu_{\mathbb{P}^1}(w_1),
\end{align*}
and commutation with the Dirac operator
\begin{align*}
\diracJ \delbar{E} (\alpha \cdot \xi)(z_1)
&= \qj (1 + |z_1|^2) \barpartial \left( (\Bar{\alpha}_1 + \Bar{\alpha}_2 z_1)^{-1}\xi_1\left(\frac{\alpha_1z_1-\alpha_2}{\Bar{\alpha}_1 + \Bar{\alpha}_2z_1}\right) \right) \\
&= \qj (1 + |z_1|^2) (\Bar{\alpha}_1 + \Bar{\alpha}_2 z_1)^{-1} \barpartial \left(\frac{\Bar{\alpha}_1\Bar{z}_1-\Bar{\alpha}_2}{\alpha_1 + \alpha_2\Bar{z}_1}\right) (\barpartial\xi_1)(\alpha^{-1}z_1) \\
&= \qj (1 + |z_1|^2) (\Bar{\alpha}_1 + \Bar{\alpha}_2 z_1)^{-1} (\alpha_1 + \alpha_2\Bar{z}_1)^{-2} (\barpartial\xi_1)(\alpha^{-1}z_1) \\
&= (\Bar{\alpha}_1 + \Bar{\alpha}_2 z_1)^{-1} \qj (1 + |\alpha^{-1}z_1|^2) (\barpartial\xi_1)(\alpha^{-1}z_1) \\
&= \alpha \cdot (\diracJ \delbar{E} \xi)(z_1).
\qedhere
\end{align*}
\end{proof}

In light of the above properties, we deviate form our general strategy, embodied in Definition~\ref{def:sobolev norm}, which avoids the use of Riemannian metrics.
More precisely, we endow functions $\banach{p}(\mathbb{P}^1,\qat)$ with the norm with respect to the induced measure $\mu_{\mathbb{P}^1}$ and the sections $H^0(\mathbb{P}^1,\ban{p}{E})$ with norm
\begin{align*}
\labelthis{eq:Lpnorminvariance}
\index{Sobolev space}
\|\xi\|_p^p
&= \|\,|\xi|_E\,\|_{\banach{p}(\mathbb{P}^1,\qat)}^p
=\int_{\mathbb{P}^1}|\xi|^p \,\dmu_{\mathbb{P}^1},
\end{align*}
which when expressed in charts is
\begin{equation}
\label{eq:Lpnorminvariance2}
\|\xi\|_p^p
= \int_{\mathbb{C}}|\xi_1(z_1)|^p(1+|z_1|^2)^{\frac{p}{2}-2} \,\dmu(z_1)
= \int_{\mathbb{C}}|\xi_2|^p(1+|z_2|^2)^{\frac{p}{2}-2} \,\dmu(z_2).
\end{equation}
Further, it is conceptually helpful to recognize that $H^0(\mathbb{P}^1,\ban{2}{E})$ has the Hermitian pairing
\[
\langle \xi, \xi' \rangle
= \int_{\mathbb{P}^1} h_E(\xi,\xi')\,\dmu_{\mathbb{P}^1},
\]
although we will not directly use it.
On the subsets with compact closures in $\SO_l$ the weight $(1+z_l\Bar{z}_l)^{\frac{p}{2}-2}$ is bounded from above and from below by strictly positive numbers.
Hence the norm~\eqref{eq:Lpnorminvariance} is equivalent to the norm of $H^0(\mathbb{P}^1,\ban{p}{E})$ constructed in Definition~\ref{def:sobolev norm}.
The advantage of this definition, of course, is that it makes $\alpha \mapsto \alpha\cdot\xi$ an isometry of the Banach space $H^0(\mathbb{P}^1,\ban{p}{E})$ and of the Hilbert space $H^0(\mathbb{P}^1,\ban{2}{E})$, due to~Lemma~\ref{moebius action}.

To derive an integral kernel for $\Op{R}(0,\ci\lambda)$, we first find one for the inverse of the composition 
\begin{align*}
\big((\qi\lambda)\ri-\diracJ\delbar{E}\big)&\circ\big((-\qi\lambda)\ri-\diracJ\delbar{E}\big) \\
&= (\qi\lambda)\ri\circ(-\qi\lambda)\ri-\diracJ\delbar{E}\circ(-\qi\lambda)\ri-(\qi\lambda)\ri\circ\diracJ\delbar{E}+(\diracJ\delbar{E})^2 \\
&=\lambda^2+\big(\diracJ\delbar{E}\big)^2,
\labelthis{dirac square}
\end{align*}
for $\lambda\in\mathbb{R}$.
Unlike the Dirac operator, this composition commutes with the complex structure. 
In fact it bears a strong resemblance to the Laplace operator and has a well-behaved spectral theory~\cite[\S3.3; p~87.~Remark]{Fr1}.
We exploit this in the following lemma to find the integral kernel of its inverse.

\begin{lemma}\label{resolvent square dirac}
For $\lambda\in\mathbb{R}$ the inverse of~\eqref{dirac square} acts on $\xi\in H^0(\mathbb{P}^1,\ban{2}{E})$ as
\begin{align}\label{integral kernel 2}
\left(\big(\lambda^2+\big(\diracJ\delbar{E}\big)^2\big)^{-1}\xi\right)(z)
&=\int_{\mathbb{C}}\Func{K}_{\lambda}(z,w)\xi(w)\frac{d\mu(w)}{1+|w|^2}.
\end{align}
The kernel $\Func{K}_{\lambda}$ has the explicit formula
\begin{equation}
\label{explicit formula}
\begin{gathered}
\Func{K}_{\lambda}(z,w)
= \frac{1}{\sqrt{1+|z|^2}} \alpha_z\cdot\big(\fundsol_{\lambda}\big(|w|\big)\big)
=\frac{1}{1+z\Bar{w}}\fundsol_{\lambda}\left( \left|\frac{w-z}{1+\Bar{z}w}\right| \right), \text{ where }
\\
\alpha_z = (1+|z|^2)^{-1}(1 + \qj z) \in \mathbb{S}^3, 
\text{ and } \\
\fundsol_{\lambda}(r)
= \frac{2}{\pi} \int_{r}^{\infty}\frac{\sinh(\lambda(\pi-2\arctan(s)))}{\sinh(\lambda\pi)}\, \frac{\sqrt{1+r^2}}{\sqrt{s^2-r^2}}  \,\frac{ds}{1+s^2}.
\end{gathered}
\end{equation}
\end{lemma}
\begin{proof}
We note that~\eqref{dirac square} is an unbounded closed operator of $H^0(\mathbb{P}^1,\ban{2}{E})$, whose inverse is the composition $\Op{R}(0,-\ci\lambda)\circ\Op{R}(0,\ci\lambda)$ of two compact resolvents constructed in Lemma~\ref{dirac resovent}.
Furthermore, $\xi(0)$ is well defined, since for $\V=0$ the arguments in the proof of Lemma~\ref{lem:sobolev regularity}~(ii) extend for $1<p<2$ to $H^0(\X,\sob{2,p}{E})$ and show that this space is embedded into the continuous section of $E_\qat$.
We first construct $\fundsol_{\lambda}$ on $r\in[0,\infty)$ such that any $\xi$ in the domain of~\eqref{dirac square} obeys:
\begin{align}
\label{eq:fundamental solution}
\xi(0)
=\int_{\mathbb{C}}\fundsol_{\lambda}(|w|)\, \Big(\big(\lambda^2+(\diracJ\delbar{E})^2\big)\xi\Big)(w)\frac{d\mu(w)}{1+|w|^2}.
\end{align}
Consider the right hand side of this formula.
By changing to polar coordinates $w = re^{\ci\varphi}$, we obtain
\[
\int_{0}^r \fundsol_{\lambda}(r) \left( \int_0^{2\pi} \Big(\big(\lambda^2+(\diracJ\delbar{E})^2\big)\xi\Big)(re^{\ci\varphi}) \,d\varphi \right) \frac{r\,dr}{1+r^2}.
\]
We interpret this as saying that only the average of $\big(\lambda^2+(\diracJ\delbar{E})^2\big)\xi$ over the compact subgroup $\mathbb{S}^1 = \{ e^{\ci\varphi} \mid \varphi\in\mathbb{R}/2\pi\mathbb{Z}\big\} \subset\mathbb{S}^3$ is necessary for the calculation of~\eqref{eq:fundamental solution}.
Therefore define the averaging operator
\begin{gather*}
\Op{A}\xi(|z|)
=\frac{1}{2\pi}\int_{0}^{2\pi} \xi(|z|e^{\ci\varphi})\,d\varphi
=\frac{1}{4\pi}\int_0^{4\pi} e^{\ci\varphi/2} (e^{\ci\varphi/2} \cdot\xi )(|z|)\,d\varphi.
\end{gather*}
The averaging operator $\Op{A}$ projects $H^0(\mathbb{P}^1,\ban{2}{E})$ onto the closed subspace $\Spa{A}$ that consists of sections whose representative on $\SO_1$ depends only on $r=|z|$. 
The second expression for $\Op{A}$ shows that it is composed of essentially three operators: an integration over a compact space, left multiplication with a complex number, and an $\mathbb{S}^3$\=/action.
But $\lambda^2+\big(\diracJ\delbar{E}\big)^2$ commutes with all three, and so therefore also with $\Op{A}$.
The integral in~\eqref{eq:fundamental solution} now reduces to
\[
\int_{0}^r \fundsol_{\lambda}(r) \Big(\big(\lambda^2+(\diracJ\delbar{E})^2\big)\Op{A}\xi\Big)(r) \frac{2\pi r\,dr}{1+r^2}.
\]
By a direct calculation $\lambda^2+(\diracJ\delbar{E})^2$ acts on elements of $\Spa{A}$ as
\begin{equation}
\label{eq:green-vacuum-operator-1}
\lambda^2-\left(\frac{1+r^2}{2}\right)^2\left(\frac{d^2}{dr^2}+\frac{1}{r}\frac{d}{dr}\right)-\frac{1+r^2}{2}r\frac{d}{dr}.
\end{equation}
Due to Lemma~\ref{dirac resovent}, this operator is an unbounded closed operator with compact resolvent on $\Spa{A}$.
The substitution $x=\frac{r^2-1}{r^2+1}$ makes the coefficients of~\eqref{eq:green-vacuum-operator-1} polynomial.
It transforms the measure $\frac{2\pi r\,dr}{1+r^2}$ on $r\in[0,\infty)$ into $\frac{\pi}{1-x}dx$ on $x\in[-1,1)$ and~\eqref{eq:green-vacuum-operator-1} into
\begin{equation}
\label{eq:green-vacuum-operator-3}
\lambda^2-(1-x^2)\frac{d^2}{dx^2}-(1-x)\frac{d}{dx}.
\end{equation}
This operator is also an unbounded closed operator with compact resolvent of the $\banach{2}$\=/space on $[-1,1)$ with measure $\frac{\pi}{1-x}dx$.
Thus~\eqref{eq:fundamental solution} holds for all $\xi$ in the domain of~\eqref{dirac square} if and only if all $\Tilde{\xi}$ in the domain of~\eqref{eq:green-vacuum-operator-3} obey
\begin{align*}
\Tilde{\xi}(0)
&=\int_0^\pi \fundsol_{\lambda}(r(x))\, \left( \lambda^2-(1-x^2)\frac{d^2}{dx^2}-(1-x)\frac{d}{dx} \right)\Tilde{\xi}(x)\,\frac{\pi}{1-x}dx.
\end{align*}

The operator~\eqref{eq:green-vacuum-operator-3} is almost, but not quite, in standard form.
The multiplication with $(1-x)$ is an isomorphism from the $\banach{2}$\=/functions on $x\in[-1,1]$ with measure $(1-x)\pi dx$ onto the $\banach{2}$\=/functions on $[-1,1)$ with measure $\frac{\pi}{1-x} dx$. 
Due to~\cite[\S10.8 (14)]{MOT}, the Jacobi polynomials $(P_{n-1}^{(1,0)})_{n\in\mathbb{N}}$ are a complete set of eigenfunctions with eigenvalues $-n^2$ of the following operator:
\[
(1-x)^{-1}\left((1-x^2)\frac{d^2}{dx^2}+(1-x)\frac{d}{dx}\right)(1-x)
=(1-x^2)\frac{d^2}{dx^2}-(1+3x)\frac{d}{dx}-1,
\]
This conjugation of~\eqref{eq:green-vacuum-operator-3} is indeed self-adjoint in the $\banach{2}$\=/space on $[-1,1]$ with measure $(1-x)\pi dx$, as is~\eqref{eq:green-vacuum-operator-3} in the $\banach{2}$\=/space on $[-1,1)$ with measure $\frac{\pi}{1-x}dx$.

Therefore we can use the Jacobi polynomials to write a basis of eigenfunctions to solve the integral equation.
Let $h_n$ denote $\int_{-1}^1(P_n^{(1,0)}(x))^2(1-x)dx$. Since $(P_n^{(1,0)})_{n\in\mathbb{N}_0}$ is an orthogonal basis, formally the following equation holds for all $g\in\mathbb{R}[x]$:
\begin{align}\label{delta}
\hspace{-1mm}\int_{-1}^1\hspace{-3mm}g(x)\delta(x\!+\!1)\frac{\pi dx}{1\!-\!x} &\!=\!g(-1)&\hspace{-2.5mm}\text{with }
\delta(1\!+\!x)&\!=\!\sum_{n=1}^\infty\frac{2P_{n-1}^{(1,0)}(-1)(1\!-\!x)P_{n-1}^{(1,0)}(x)}{\pi h_{n-1}}.
\end{align}
Let us show that $\delta(x+1)\frac{\pi dx}{1-x}$ is the delta distribution at $x=-1$, not just formally. We insert the following equalities listed in~\cite[\S10.8. (3),(4) and (13)]{MOT}:
\begin{align*}
P_n^{(1,0)}(-1)&=(-1)^nP_n^{(0,1)}(1)=(-1)^n,&h_n&
=\frac{2^2\Gamma(n+2)\Gamma(n+1)}{(2n+2)n!\Gamma(n+2)}=\frac{2}{n+1},
\end{align*}
in the following sum, which is formally transformed by~\eqref{eq:green-vacuum-operator-3} into $\delta(1+x)$:
\begin{equation}\label{greens function 1}
\sum_{n=1}^\infty\frac{2P_{n-1}^{(1,0)}(-1)(1-x)P_{n-1}^{(1,0)}(x)}{\pi h_{n-1}(\lambda^2+n^2)}
=-\sum_{n=1}^\infty\frac{(-1)^nn}{\pi(\lambda^2+n^2)}(1-x)P_{n-1}^{(1,0)}(x).
\end{equation}
This sum is a $\banach{2}$\=/function on $[-1,1)$ with measure $\frac{\pi dx}{1-x}$, since $(\sqrt{\frac{n}{\pi}}(1-x)P_{n-1}^{(1,0)})_{n\in\mathbb{N}}$ is an orthonormal basis of this Hilbert space.
The self-adjointness of the closed unbounded operator~\eqref{eq:green-vacuum-operator-3} on this Hilbert space implies for any $g$ in the domain of this operator, that the integral of the product of the image of $g$ under~\eqref{eq:green-vacuum-operator-3} with~\eqref{greens function 1} times $\frac{\pi dx}{1-x}$ yields $g(-1)$.
Consequently~\eqref{greens function 1} is the fundamental solution at $x=-1$ of~\eqref{eq:green-vacuum-operator-3}.

In order to derive the explicit formula~\eqref{explicit formula} for this fundamental solution $\fundsol_{\lambda}$ we modify the derivation of Mehler's integral for Legendre polynomials and derive a similar integral representation for the Jacobi polynomials $P_{n-1}^{(1,0)}(x)$.
The definition of Mehler's integral can be found in~\cite[\S10.10~(43)]{MOT} and we will follow the derivation of~\cite[\S4.8]{Sz}.  
The starting point is the replacement of \cite[(4.8.1)]{Sz} by~\cite[(4.4.6)]{Sz}:
\[
P_{n-1}^{(1,0)}(x)=\frac{1}{2\pi\ci}\oint\left(\frac{t^2-1}{2t-2x}\right)^{n-1}\left(\frac{1-t}{1-x}\right)\frac{dt}{t-x}\quad\text{for}\quad x\in(-1,1).
\]
Here the path surrounds $t=x$ in the complex plane once in the anti-clockwise direction. 
We choose the same path as in \cite[(4.8.4)]{Sz}, a circle of radius $\sqrt{\frac{1-x}{2}}$ around $t=1$.
This circle is parameterized in a somewhat unusual manner.
The circle is divided into two arcs by the unit circle~\cite[(4.8.5)]{Sz}.
On the outside arc $\varphi\in(-\arccos(x),\arccos(x))$ and $t=e^{\ci\varphi}+e^{\ci\frac{\varphi}{2}}\sqrt{2\cos\varphi-2x}$; 
On the inside arc $\varphi\in(-\arccos(x),\arccos(x))$ also, but traversed from high to low, and $t=e^{\ci\varphi}-e^{\ci\frac{\varphi}{2}}\sqrt{2\cos\varphi-2x}$.
Here $\sqrt{2\cos\varphi-2x}$ is the unique positive square root. 
The advantage of this choice is that for both arcs we have $\frac{t^2-1}{2-2x}=e^{\ci\varphi}$. 
We insert the formula for $\frac{dt}{t-x}$ in \cite[p.~86]{Sz} and obtain for $2\pi(1-x)P_{n-1}^{(1,0)}(x)$ 
$$\int\limits_{-\arccos(x)}^{\arccos(x)}\hspace{-3mm}e^{\ci n\varphi}\left(\frac{e^{-\ci\frac{\varphi}{2}}-e^{\ci\frac{\varphi}{2}}}{\sqrt{2\cos(\varphi)-2x}}-1\right)d\varphi+\hspace{-3mm}\int\limits_{\arccos(x)}^{-\arccos(x)}\hspace{-3mm}e^{\ci n\varphi}\left(\frac{e^{-\ci\frac{\varphi}{2}} - e^{\ci\frac{\varphi}{2}}}{-\sqrt{2\cos(\varphi)-2x}}-1\right)d\varphi.$$
We insert $x=-\cos(2y)=\cos(\pi-2y)$ and obtain a variant of Mehler's formula:
\begin{align*}
\pi&(1+\cos(2y)) P_{n-1}^{(1,0)}(-\cos(2y))
=-\int\limits_{2y-\pi}^{\pi-2y}\frac{e^{\ci n\varphi}2\ci\sin(\frac{\varphi}{2})\,d\varphi}{\sqrt{2\cos(\varphi)+2\cos(2y)}} \\
&=-\int\limits_{2y-\pi}^0\frac{4\sin(n\varphi)\sin(\frac{\varphi}{2})\,d\varphi}{\sqrt{2\cos(\varphi)\!+\!2\cos(2y)}}
=-(-1)^n\int\limits_{2y}^{\pi}\frac{4\sin(n\varphi)\cos(\frac{\varphi}{2})\,d\varphi}{\sqrt{2\cos(2y)\!-\!2\cos(\varphi)}} \\
&=-(-1)^n \int\limits_{y}^{\pi/2}\frac{8\sin(2n\varphi)\cos(\varphi)\,d\varphi}{\sqrt{2\cos(2y)\!-\!2\cos(2\varphi)}}
=-(-1)^n \int\limits_{y}^{\pi/2}\frac{4\sin(2n\varphi)\cos(\varphi)\,d\varphi}{\sqrt{\sin^2(\varphi)\!-\!\sin^2(y)}},
\end{align*}
where in the final line we have rescaled $\varphi$ but kept the same variable name.
Inserting this into~\eqref{greens function 1} at $x=-\cos(2y)$ yields
\begin{align*}
\sum_{n=1}^\infty\frac{(-1)^{2n} n}{\pi^2(\lambda^2+n^2)}
\int\limits_{y}^{\pi/2}\frac{4\sin(2n\varphi)\cos(\varphi)\,d\varphi}{\sqrt{\sin^2(\varphi)\!-\!\sin^2(y)}}
= \int\limits_{y}^{\pi/2} \sum_{n=1}^\infty\frac{4n\sin(2n\varphi)}{\pi^2(\lambda^2+n^2)} \frac{\cos(\varphi)\,d\varphi}{\sqrt{\sin^2(\varphi)\!-\!\sin^2(y)}} .
\end{align*}
The sum can be evaluated with the following Fourier series on $\varphi\in[0,\pi/2]$:
\begin{align*}
2\pi&\frac{\sinh(\lambda(\pi-2\varphi))}{\sinh(\lambda\pi)}
=\frac{2\pi e^{-\lambda(2\varphi-\pi)}}{2\sinh(\lambda\pi)} + \frac{2\pi e^{\lambda(2\varphi-\pi)}}{2\sinh(-\lambda\pi)}
=\sum_{n\in\mathbb{Z}}\left(\frac{e^{2\ci n\varphi}}{\lambda+\ci n}-\frac{e^{2\ci n\varphi}}{\lambda-\ci n}\right) 
\\
&=\sum_{n=1}^\infty\left(\frac{e^{2\ci n\varphi}\!-\!e^{-2\ci n\varphi}}{\lambda+\ci n}\!-\!\frac{e^{\ci n\varphi}\!-\!e^{-\ci n\varphi}}{\lambda-\ci n}\right)
=\sum_{n=1}^\infty\frac{-2\ci n(e^{2\ci n\varphi}\!-\!e^{-2\ci n\varphi})}{\lambda^2+n^2} 
\\
&=\sum\limits_{n=1}^\infty\frac{4n\sin(2n\varphi)}{\lambda^2+n^2}.
\end{align*}
The substitution $x=-\cos(2y)$ may not seem the natural choice at first, but it is geometrically significant, since $y$ is the radial coordinate of $\mathbb{P}^1$ with distance taken with respect to the Fubini-Study metric.
However, we generally prefer the coordinate $r$ for reasons of ease, and so let us rewrite the above formulas in $r$.
From
\[
r^2 
= \frac{1+x}{1-x}
= \frac{1 - \cos(2y)}{1 + \cos(2y)}
= \frac{\sin^2(y)}{\cos^2(y)}
= \tan^2(y),
\]
we see that $y = \arctan(r)$. 
Likewise we make the substitution $\varphi = \arctan(s)$.
Then, making use of $\sin(\arctan(s)) = \frac{s}{\sqrt{1+s^2}}$ and $\cos(\arctan(s)) = \frac{1}{\sqrt{1+s^2}}$, we can write the fundamental solution as
\begin{align*}
\frac{2}{\pi} & \int_{y}^{\pi/2} \frac{\sinh(\lambda(\pi-2\varphi))}{\sinh(\lambda\pi)} \frac{\cos(\varphi)\,d\varphi}{\sqrt{\sin^2(\varphi)\!-\!\sin^2(y)}} \\
&= \frac{2}{\pi}\int_{r}^{\infty}\frac{\sinh(\lambda(\pi-2\varphi(r))) \frac{1}{\sqrt{1+s^2}} \,\frac{ds}{1+s^2}}{\sinh(\lambda\pi)\sqrt{\frac{s^2}{1+s^2}-\frac{r^2}{1+r^2}}} \\
&= \frac{2}{\pi}\sqrt{1+r^2}\int_{r}^{\infty}\frac{\sinh(\lambda(\pi-2\varphi(s)))}{\sinh(\lambda\pi)}\, \frac{1}{\sqrt{s^2(1+r^2)-r^2(1+s^2)}}  \,\frac{ds}{1+s^2} \\
&= \frac{2}{\pi}\sqrt{1+r^2}\int_{r}^{\infty}\frac{\sinh(\lambda(\pi-2\varphi(s)))}{\sinh(\lambda\pi)}\, \frac{1}{\sqrt{s^2-r^2}}  \,\frac{ds}{1+s^2}.
\end{align*}
This is exactly the formula given for $\fundsol_\lambda(r)$ in~\eqref{explicit formula}.
Although this derivation seems formidable, we may be certain of its correctness by comparison to~\cite[(6)]{Gasper2006}, which reaches the equivalent formula~\eqref{eq:fund sol int by parts} using a different method.

We have succeeded in proving~\eqref{eq:fundamental solution}, which is the special case of the theorem with $z = 0$, such that $\Func{K}_{\lambda}(0,w) = \fundsol_{\lambda}(|w|)$.
To prove the general case, we again make use of the invariance of the Dirac operator with respect to $\mathbb{S}^3$.
Indeed, notice that~\eqref{eq:fundamental solution} can be written as an inner product on $H^0(\ban{2}{E})$:
\[
\Big(\big(\lambda^2+(\diracJ\delbar{E})^2\big)^{-1}\xi\Big)(0)
= \int_{\mathbb{P}^1} h_E\Big( \fundsol_{\lambda}\big(|w|\big), \xi(w)\Big) \, \dmu_{\mathbb{P}^1}(w).
\]
For any point $z \in \mathbb{P}^1$, choose a point $\alpha_z \in \mathbb{S}^3$ with $\Op{P}(\alpha) = z$.
An obvious choice is $\alpha_z = (1 + |z|^2)^{-1/2}(1 + \qj z)$.
Then
\begin{align*}
\Big(\big(\lambda^2&+\big(\diracJ\delbar{E}\big)^2\big)^{-1}\xi\Big)(z)
= \frac{1}{\sqrt{1+|z|^2}}\Big(\big(\lambda^2+\big(\diracJ\delbar{E}\big)^2\big)^{-1}(\alpha_z^{-1}\cdot\xi)\Big)(0)\\
&= \frac{1}{\sqrt{1+|z|^2}} \int_{\mathbb{P}^1} h_E\Big( \fundsol_{\lambda}\big(|w|\big), (\alpha_z^{-1}\cdot\xi)(w)\Big) \, \dmu_{\mathbb{P}^1}(w) \\
&= \frac{1}{\sqrt{1+|z|^2}} \int_{\mathbb{P}^1} h_E\Big( \alpha_z\cdot\big(\fundsol_{\lambda}\big(|w|\big)\big), \xi(w)\Big) \, \dmu_{\mathbb{P}^1}(w),
\end{align*}
where we have used Lemma~\ref{moebius action} to pass $\alpha_z^{-1}$ through the operator, as well as to transform the integral by the action of $\alpha_z$.
Finally, we calculate
\[
\alpha_z\cdot\big(\fundsol_{\lambda}\big(|w|\big)\big)
= \frac{\sqrt{1+|z|^2}}{1+\Bar{z}w} \fundsol_{\lambda}\big(|\alpha_z^{-1} w|\big)
= \frac{\sqrt{1+|z|^2}}{1+\Bar{z}w} \fundsol_{\lambda}\left(\left|\frac{w-z}{1+\Bar{z}w}\right|\right).
\]
The general formula now follows (don't forget the conjugation of the first argument of $h_E$!)
\end{proof}

In the theory of integral kernels, convolution kernels play a distinguished role.
Suppose that $G$ is a compact Lie group and $H$ a closed subgroup (hence a compact Lie group in its own right).
Then given Haar measures $\mu_G$ on $G$ and $\mu_H$ on $H$, the homogeneous space $G/H$ has a $G$\=/invariant measure.
One should scale these measures so that $\mu_G(G) = \mu_{G/H}(G/H)\mu_H(H)$.
For us, taking the usual unit sphere measures on $\mathbb{S}^3$ (total measure $2\pi^2$) and $\mathbb{S}^1$ (total measure $2\pi$), the corresponding measure $\mu_{G/H}$ is exactly the Fubini-Study measure (total measure $\pi$, one quarter of the unit $2$\=/sphere $4\pi$).
Functions on $G$ obey
\[
\int_G f\,\dmu_G
= \int_{G/H} \left( \int_H f(\alpha \beta) \,\dmu_h(\beta) \right) \dmu_{G/H}(\alpha H).
\]
Reiter calls this the extended Weil formula~\cite[(3.4.11)]{Reiter2000}; it might reasonably be thought of as a generalization of Fubini's theorem.
In particular, for a function on $G/H$ and its pullback to $G$ under the canonical projection $\Op{P} : G \to G/H$, their norms differ by a constant:
\begin{equation}
\label{eq:p norms on homogeneous space}
\|f\|_{\banach{p}(G/H)} = \mu_H(H)^{-1/p}\|f \circ \Op{P}\|_{\banach{p}(G)}.
\end{equation}
This fact is easily understood in terms of distribution functions, since for $U \subset G/H$ we have $\mu_G(\Op{P}^{-1}[U]) = \mu_H(H)\mu_{G/H}(U)$.
This idea will be used in the lemma below to extend this formula to Lorentz norms.

The convolution of functions $f,g$ on $G$ is defined in~\cite[(20.10)~Theorem]{Hewitt1979} as
\[
f \ast g (\alpha)
= \int_G f(\beta)g(\beta^{-1}\alpha)\,\dmu_G(\beta)
= \int_G f(\alpha\beta^{-1})g(\beta)\,\dmu_G(\beta).
\]
The equality of these two formula holds if the group $G$ is unimodular, which all compact groups are.
If $g$ is the pull back of a function on $G/H$, which is the case exactly if $g$ is right $H$\=/invariant, then so too is $f\ast g$, as is readily seen from the first integral above.
As in Euclidean space, we have Young's inequality for convolutions.
But in the following lemma, we extend this to a generalized Young's inequality on homogeneous spaces.

\begin{lemma}\label{lem:generalized convolution}
For functions $f,g \in \banach{1}(G/H)$ we define their convolution to be
\[
f \hatast g (\alpha)
= \mu_H(H)^{-1} \int_{G} f(\Op{P}(\beta)) g(\Op{P}(\beta^{-1}\alpha)) \,\dmu_{G}(\beta)
= \mu_H(H)^{-1} (f\circ\Op{P})\ast(g\circ\Op{P})(\alpha),
\]
where $\Op{P}: G \to G/H$ is the canonical projection.
The convolution descends to a function in $L^1(G/H)$.
For $r^{-1} + 1 = p^{-1} + q^{-1}$ and $c^{-1} \leq a^{-1} + b^{-1}$ there is a constant $C$, dependent only on the exponents, such that
\[
\|f \hatast g\|_{\banach{r,c}(G/H)}
\leq C \|f\|_{\banach{p,a}(G/H)} \|g\|_{\banach{q,b}(G/H)}.
\]
\end{lemma}
\begin{proof}
Before we discuss the convolution, let us extend~\eqref{eq:p norms on homogeneous space} to Lorentz norms.
These norms are defined in terms of decreasing rearrangements\cite[Def~2.1.5]{BS} and maximal functions~\cite[Def~2.3.1]{BS}:
\begin{align*}
(f\circ\Op{P})^\ast(t)
&= \inf \left\{ s \mid \mu_G\left\{ |f\circ\Op{P}| > s\right\} \leq t \right\} \\
&= \inf \left\{ s \mid \mu_H(H)\mu_{G/H}\left\{ |f| > s\right\} \leq t \right\}
= f^\ast\left( \mu_H(H)^{-1}t \right) \\
(f\circ\Op{P})^{\ast\ast}(t)
&= \frac{1}{t} \int_0^t (f\circ \Op{P})^\ast(s)\,ds
= \frac{1}{t} \int_0^t f^\ast\left( \mu_H(H)^{-1}s \right)\,ds \\
&= \frac{1}{t} \int_0^{t/\mu_H(H)} f^\ast(s')\,\mu_H(H)ds' 
= f^{\ast\ast}\left( \mu_H(H)^{-1}t \right)
\end{align*}
Therefore the Lorentz norms~\cite[Def~4.4.4]{BS} of $G$ and $G/H$ are related by 
\begin{align*}
\|f\circ \Op{P}\|_{(p,q)}
&= \left[ \int_0^\infty \left[ t^{1/p} f^{\ast\ast}\left( \mu_H(H)^{-1}t \right) \right]^q \frac{dt}{t} \right]^{1/q} \\
&= \left[ \int_0^\infty \left[ (\mu_H(H)t)^{1/p} f^{\ast\ast}(t') \right]^q \frac{dt'}{t'} \right]^{1/q} \\
&= \mu_H(H)^{1/p} \|f\|_{(p,q)},
\end{align*}
for $q< \infty$ and for $q=\infty$
\begin{align*}
\|f\circ \Op{P}\|_{(p,\infty)}
&= \sup_{t} t^{1/p} (f\circ\Op{P})^{\ast\ast}(t)
= \sup_{t} t^{1/p} f^{\ast\ast}\left( \mu_H(H)^{-1}t \right) \\
&= \sup_{t'} (\mu_H(H)t')^{1/p} f^{\ast\ast}(t') 
= \mu_H(H)^{1/p} \|f\|_{(p,\infty)}.
\end{align*}

Now, consider the definition of convolution on $G/H$ above: it is the convolution of the pullbacks of $f$ and $g$ to $G$, normalized by the measure of $H$ (cf.~\cite[\S2.1~Example~1]{Rouviere2014}).
As we have already observed, since $g \circ \Op{P}$ is $H$\=/invariant, so too is the convolution.
The usual Young's inequality for convolution on groups~\cite[(20.18)~Theorem]{Hewitt1979} shows that the result is an $\banach{1}$ function on $G$.
Hence this convolution does indeed belong to $\banach{1}(G/H)$.

In fact, the usual Young's inequality is exactly what we need to show that convolution on groups is a convolution operator in the sense of Bennett and Sharpley~\cite[Theorem~4.7.6]{BS}.
Hence the generalized Young's inequality holds for the convolution of functions on groups.
Indeed, it holds not only for the exponent $c^{-1} = a^{-1} + b^{-1}$, but also for all greater $c$, due to~\cite[Prop~4.4.2]{BS} about the nesting of Lorentz spaces (compare~\cite[2.10.1~Theorem]{Zi}).
But now this result follows easily for the homogeneous space with the same constant.
\begin{align*}
\|f \hatast g \|_{\banach{r,c}(G/H)}
&= \mu_H(H)^{-1} \|(f\circ \Op{P}) \ast (g\circ\Op{P}) \|_{\banach{r,c}(G/H)} \\
&= \mu_H(H)^{-1-r^{-1}} \|(f\circ \Op{P}) \ast (g\circ\Op{P}) \|_{\banach{r,c}(G)} \\
&\leq \mu_H(H)^{-1-r^{-1}} C \|f\circ \Op{P}\|_{\banach{p,a}(G)} \|g\circ\Op{P}\|_{\banach{q,b}(G)} \\
&= \mu_H(H)^{-1-r^{-1}+p^{-1}+q^{-1}} C \|f\|_{\banach{p,a}(G/H)} \|g\|_{\banach{q,b}(G/H)} \\
&= \mu_H(H)^{0} C \|f\|_{\banach{p,a}(G/H)} \|g\|_{\banach{q,b}(G/H)}.
\qedhere
\end{align*}
\end{proof}

The astute reader may have already seen that the action of $\alpha_z$ in~\eqref{explicit formula} makes the integral kernel one of convolution type.
It is important therefore to control the norm of $\Func{K}_{\lambda}(0,w) = \fundsol_{\lambda}(|w|)$.
We prepare some needed estimates in the following lemma.
\begin{lemma}
\label{green}
For $\lambda\in\mathbb{R}$ and $y\in(0,\pi)$ the following estimates hold:
\begin{enumeratethm}
\setlength{\itemsep}{0.5em}
\item 
$\displaystyle{0<\fundsol_{\lambda}(r)\leq \tfrac{1}{\pi} e^{- 2\lambda \varphi(r)} \ln\left( \frac{\sqrt{1+r^2}+1}{\sqrt{1+r^2}-1} \right)}$.
\item 
$\displaystyle{0<r\,|\lambda|\,\fundsol_{\lambda}(r)\leq \tfrac{1}{2e}}$.
\item
Define $\displaystyle{b(\lambda,y) := \frac{\lambda y\cosh(\lambda(\pi-y)) + \sinh(\lambda(\pi-y))}{\sinh(\lambda\pi)}}$.
Then $0 < b(|\lambda|,y) \leq 1$, $b$ is decreasing in $y$, and $\lim_{\lambda \to \infty} e^{\lambda\varepsilon/2} b(\lambda,\varepsilon) = 0$ for any $0 < \varepsilon < \pi/2$.
\item 
$\displaystyle{0<- r (1+r^2) \fundsol_{\lambda}'(r) \leq \tfrac{2}{\pi} b(|\lambda|,\arctan r) }$,
\item 
The function $z \mapsto \sqrt{1+|z|^2}|z|^{-1}$ belongs to $\banach{2,\infty}(\mathbb{P}^1,\qat)$.
\end{enumeratethm}
\end{lemma}
Let us discuss these inequalities, to build our intuition.
Inequality (i) shows, after rationalizing the denominator (giving $r^{2}$) and splitting the $\ln$ into a difference, that $\fundsol_{\lambda}$ has (at worst) a logarithmic singularity at $r=0$ for fixed $\lambda$.
So naturally, $|\lambda|\fundsol_{\lambda}$ also has a logarithmic singularity at $r=0$ for fixed $\lambda$.
But the ``strength'' of the singularity clearly increases with $\lambda$.
To bound it uniformly in $\lambda$, we must allow a singularity of a worse order, namely $r^{-1}$.
This is the idea of (ii).
Similar to (i), inequality (iii) shows that the derivative $\fundsol'_{\lambda}$ has a singularity (at worst) of the order of $r^{-1}$, which is to be expected.
However, away from the singularities at $r=0$, exponentially decay in $|\lambda|$ dominates in (i) and (iii).
More specifically, as shall be needed in the proof of Lemma~\ref{weakly continuous resolvent}, for all $0 < \varepsilon < \pi/2$ both the function and its derivative converge to zero in $|\lambda| \to \infty$ uniformly for $r \in [\tan\varepsilon, \infty)$.
The uniformity comes from the fact that the bounds are decreasing functions of $r$, so it suffices to show the convergence to zero at $r = \tan\varepsilon$ as $|\lambda| \to \infty$.
\begin{proof}
Note that $\sinh(\lambda(\pi-t))$ and $\sinh(\lambda\pi)$ are both odd functions of $\lambda$, so $\fundsol_{\lambda} = \fundsol_{|\lambda|}$ and it suffices to consider only positive values of $\lambda$.
The integrand of the fundamental solution is positive, so $0 <\fundsol_{\lambda}(r)$ holds.
We see in the formula~\eqref{explicit formula} for $\fundsol_\lambda$ that it has essentially two parts: one factor composed of hyperbolic trigonometric functions that contains $\lambda$ and another that is algebraic and free from $\lambda$.
Our modus operandi is to estimate away the first factor and integrate what remains.
This rather tight bound simplifies the first factor:
\begin{equation}
\label{eq:sinh estimate}
\frac{\sinh (\lambda(\pi-2\varphi))}{\sinh \lambda\pi}
= \frac{1- e^{-2\lambda(\pi-2\varphi)}}{1 - e^{-2\lambda\pi}} e^{- 2 \lambda \varphi}
\leq e^{- 2\lambda \varphi}.
\end{equation}
So we will try to prove (i) and (ii) using the following formula:
\[
\fundsol_{\lambda}(r)
\leq \frac{2}{\pi} \int_{r}^{\infty} e^{- 2\lambda \varphi(s)}\, \frac{\sqrt{1+r^2}}{\sqrt{s^2-r^2}}  \,\frac{ds}{1+s^2}.
\]
For the sake of brevity we have continued with the notation $\varphi(s) = \arctan(s)$ from the previous lemma.
The simplest estimate is to notice that the exponential is decreasing in $s$ and therefore has its maximum at $s=r$.
If we take this bound, the remaining integral is then elementary.
Because we will use it again below, define
\[
L(r,s) := \int \frac{\sqrt{1+r^2}}{\sqrt{s^2-r^2}}  \,\frac{ds}{1+s^2} 
= \frac{1}{2}\ln\left( \frac{s\sqrt{1+r^2}+\sqrt{s^2-r^2}}{s\sqrt{1+r^2}-\sqrt{s^2-r^2}} \right).
\]
Note that $L(r,r) = \frac{1}{2}\ln 1 = 0$.
Then
\begin{equation*}
\fundsol_{\lambda}(r)
\leq \frac{2}{\pi} e^{- 2\lambda \varphi(r)} L(r,s) \Big|_{s=r}^{s=\infty}
= \frac{1}{\pi} e^{- 2\lambda \varphi(r)}\ln\left( \frac{\sqrt{1+r^2}+1}{\sqrt{1+r^2}-1} \right).
\end{equation*}
This is exactly (i).

For (ii), we want to bound $\lambda \fundsol_\lambda$ independently of $\lambda$.
So consider the function $\lambda e^{- 2\lambda \varphi(s)}$ as a function of $\lambda$.
It is zero both at $\lambda = 0$ and in the limit as $\lambda \to \infty$.
Differentiating we find it has a maximum at $\lambda = (2\varphi(s))^{-1}$, where it takes the value $(2e\varphi(s))^{-1}$.
Therefore we have the following uniform bound
\[
\lambda \fundsol_{\lambda}(r)
\leq \frac{2}{\pi}\int_{r}^{\infty} (2e\varphi(s))^{-1}\, \frac{\sqrt{1+r^2}}{\sqrt{s^2-r^2}}  \,\frac{ds}{1+s^2}.
\]
As the rest of the integrand is algebraic, it is useful to approximate $\varphi(s) = \arctan(s)$ by an algebraic function too.
The ratio of this function to $s / \sqrt{1+s^2}$ is monotonic, increasing from $1$ at $s=0$ to $\pi/2$ as $s \to \infty$.
In other words
\[
\frac{s}{\sqrt{1+s^2}} \leq \arctan(s) \leq \frac{\pi}{2}\frac{s}{\sqrt{1+s^2}}.
\]
Hence
\begin{align*}
\pi e\, \lambda \fundsol_{\lambda}(r)
&\leq \int_{r}^{\infty} \frac{\sqrt{1+r^2}\, ds}{s\sqrt{(s^2-r^2)(1+s^2)}} 
= -\frac{\sqrt{1+r^2}}{r} \left.\arctan\left(r \sqrt{\frac{1+s^2}{s^2-r^2}}\right) \right|_{r}^{\infty} \\
&= \sqrt{1+r^{-2}} \left[ - \arctan(r) + \frac{\pi}{2} \right] 
= \sqrt{1+r^{-2}} \arctan\left(r^{-1}\right) \\
&\leq \sqrt{1+r^{-2}}\, \frac{\pi}{2}\frac{r^{-1}}{\sqrt{1+r^{-2}}} 
= \frac{\pi}{2} r^{-1},
\end{align*}
where we have again used our algebraic bound on $\arctan$ at the end.
Thus we have shown (ii).

Next we move on to estimates (iii) and (iv), which concern the derivative of $\fundsol_\lambda$.
We cannot differentiate with the Leibniz rule the right hand side of~\eqref{explicit formula} in its current form, because the integrand is not continuous at the endpoint.
To overcome this, we will integrate by parts once. 
$\frac{\pi}{2}\fundsol_{\lambda}(r)$ is equal to
\begin{align*}
\left. \frac{\sinh (\lambda(\pi-2\varphi(s)))}{\sinh(\lambda\pi)} L(r,s) \right|_{s=r}^{s=\infty} 
&- \int_r^\infty \frac{\cosh(\lambda(\pi-2\varphi(s)))}{\sinh(\lambda\pi)}\left(\frac{-2\lambda}{1+s^2}\right) L(r,s) \,ds \\
&= \int_r^\infty \frac{2\lambda\cosh(\lambda(\pi-2\varphi(s)))}{\sinh(\lambda\pi)} L(r,s) \frac{ds}{1+s^2}.
\labelthis{eq:fund sol int by parts}
\end{align*}
Now the integrand is continuous and vanishing at $s=r$, and differentiating yields
\begin{align*}
\frac{\pi}{2}\fundsol_{\lambda}'(r)
&= - 0 
+ \int_r^\infty \frac{2\lambda\cosh(\lambda(\pi-2\varphi(s)))}{\sinh(\lambda\pi)} \frac{\partial L}{\partial r} \frac{ds}{1+s^2} \\
&= - \frac{1}{r\sqrt{1+r^2}} \int_r^\infty \frac{2\lambda\cosh(\lambda(\pi-2\varphi(s)))}{\sinh(\lambda\pi)} \frac{s}{\sqrt{s^2-r^2}} \frac{ds}{1+s^2}, 
\end{align*}
where
\begin{align*}
2\frac{\partial L}{\partial r}
&= \frac{s\frac{r}{\sqrt{1+r^2}} + \frac{-r}{\sqrt{s^2-r^2}}}{s\sqrt{1+r^2}+\sqrt{s^2-r^2}} - \frac{s\frac{r}{\sqrt{1+r^2}} - \frac{-r}{\sqrt{s^2-r^2}}}{s\sqrt{1+r^2}-\sqrt{s^2-r^2}} \\
&= \frac{r}{\sqrt{1+r^2}\sqrt{s^2-r^2}} \frac{-2s(s^2-r^2) - 2s(1+r^2)}{s^2(1+r^2)-(s^2-r^2)}
= \frac{-2s}{r\sqrt{1+r^2}\sqrt{s^2-r^2}} .
\end{align*}

We see that this formula for $\fundsol_{\lambda}'(r)$ already has some of the same factors that appear in (iv).
Clearly the integrand is positive, which makes $-\fundsol_{\lambda}'(r)$ greater or equal to $0$.
To obtain an upper bound, similar to (i), we might use that $\cosh(\lambda(\pi-2\varphi(s)))$ is decreasing in $s$, so that $\cosh(\lambda(\pi-2\varphi(s))) \leq \cosh(\lambda(\pi-2\varphi(r)))$.
This bound alone, however, is too crude when $r$ is small because it grows with $\lambda$.
Therefore we need to consider the `remainder' of this overestimation:
\begin{multline*}
- \frac{\pi}{2} r\sqrt{1+r^2} \fundsol_{\lambda}'(r)
= \frac{2\lambda\cosh(\lambda(\pi-2\varphi(r)))}{\sinh(\lambda\pi)} \int_r^\infty \frac{s}{\sqrt{s^2-r^2}} \frac{ds}{1+s^2} \\
- \int_r^\infty \frac{2\lambda\cosh(\lambda(\pi-2\varphi(r))) - 2\lambda\cosh(\lambda(\pi-2\varphi(s)))}{\sinh(\lambda\pi)} \frac{s}{\sqrt{s^2-r^2}} \frac{ds}{1+s^2} .
\end{multline*}
The integral in the leading term evaluates to
\[
\int_r^\infty \frac{s}{\sqrt{s^2-r^2}} \frac{ds}{1+s^2}
= \left.\frac{1}{\sqrt{1+r^2}}\arctan\left( \sqrt{\frac{s^2-r^2}{1+r^2}} \right) \right|_r^\infty
= \frac{\pi}{2} \frac{1}{\sqrt{1+r^2}}.
\]
For the remainder term, we bound it from below (since it appears with a negative sign) using $s/\sqrt{s^2 - r^2} \geq 1 \geq 1/\sqrt{1+r^2}$:
\begin{align*}
\int_r^\infty &\frac{2\lambda\cosh(\lambda(\pi-2\varphi(r))) - 2\lambda\cosh(\lambda(\pi-2\varphi(s)))}{\sinh(\lambda\pi)} \frac{s}{\sqrt{s^2-r^2}} \frac{ds}{1+s^2} \\
&\geq \frac{1}{\sqrt{1+r^2}}\int_r^\infty \frac{2\lambda\cosh(\lambda(\pi-2\varphi(r))) - 2\lambda\cosh(\lambda(\pi-2\varphi(s)))}{\sinh(\lambda\pi)} \frac{ds}{1+s^2} \\
&= \left.\frac{2\lambda\cosh(\lambda(\pi-2\varphi(r)))}{\sqrt{1+r^2}\sinh(\lambda\pi)} \arctan(s)\right|_r^\infty + \left. \frac{\sinh(\lambda(\pi-2\varphi(s)))}{\sqrt{1+r^2}\sinh(\lambda\pi)} \right|_r^\infty \\
&= \frac{2\lambda\cosh(\lambda(\pi-2\varphi(r)))}{\sqrt{1+r^2}\sinh(\lambda\pi)} \left(\frac{\pi}{2} - \arctan(r)\right) - \frac{\sinh(\lambda(\pi-2\varphi(r)))}{\sqrt{1+r^2}\sinh(\lambda\pi)} .
\end{align*}
Putting these two parts back together produces (iv), motivating the definition of $b(\lambda,y)$:
\begin{multline*}
- \frac{\pi}{2} r(1+r^2) \fundsol_{\lambda}'(r)
\leq \frac{2\lambda\cosh(\lambda(\pi-y))}{\sinh(\lambda\pi)} \left[\frac{\pi}{2} - \left(\frac{\pi}{2}-\varphi(r)\right)\right] 
+ \frac{\sinh(\lambda(\pi-y))}{\sinh(\lambda\pi)} \\
= \frac{2\lambda \varphi(r)\cosh(\lambda(\pi-2\varphi(r))) + \sinh(\lambda(\pi-2\varphi(r)))}{\sinh(\lambda\pi)} 
= b(\lambda,\varphi(r)).
\end{multline*}

We now demonstrate the features of $b(\lambda,y)$ listed in (iii).
For $y \in (0,\pi)$ the function is manifestly positive.
Its derivative,
\begin{align*}
\frac{\partial b}{\partial y}
&= \frac{- \lambda^2 y\sinh(\lambda(\pi-y))}{\sinh(\lambda\pi)},
\end{align*}
is strictly negative on this interval and hence $b$ is decreasing in $y$.
Since we can evaluate $b(\lambda,0) = 1$, this shows that $0 < b \leq 1$.
Lastly, observe that $b(\lambda,\varepsilon) \geq b(\lambda,y)$ for $y \in (\varepsilon,\pi)$ since it is a decreasing function.
So the uniform convergence to zero follows from
\begin{align*}
\lim_{\lambda \to \infty} e^{\lambda\varepsilon/2} b(\lambda,\varepsilon)
&= \lim_{\lambda \to \infty} e^{\lambda\varepsilon/2} \frac{\lambda \varepsilon [e^{\lambda(\pi-\varepsilon)} + e^{-\lambda(\pi-\varepsilon)}] + [e^{\lambda(\pi-\varepsilon)} - e^{-\lambda(\pi-\varepsilon)}]}{e^{\lambda\pi} - e^{-\lambda\pi}} \\
&= \lim_{\lambda \to \infty} \frac{\lambda \varepsilon [e^{-\lambda\varepsilon/2} + e^{-\lambda(2\pi-3\varepsilon/2)}] + [e^{-\lambda\varepsilon/2} - e^{-\lambda(2\pi-3\varepsilon/2)}]}{1 - e^{-2\lambda\pi}}
= 0.
\end{align*}

Only (v) remains. 
Similar to Example~\ref{eg:weak L2} the function $f(z) = \sqrt{1+|z|^2} |z|^{-1}$ is not square-integrable on $\mathbb{P}^1$, but it does belong to $\banach{2,\infty}$.
First we calculate the area of the ball $B(0,R)$ with respect to the Fubini-Study metric:
\begin{align*}
\mu_{\mathbb{P}^1}(B(0,R))
&= \int_0^{R} \frac{2\pi r\, dr}{(1+r^2)^2}
= \left.\frac{-\pi}{1+r^2} \right|_0^{R}
= \pi \left( 1 - \frac{1}{1+R^2} \right)
= \frac{\pi R^2}{1+R^2}.
\end{align*}
In particular, we see that $\mathbb{P}^1$ itself has area $\pi$, as previously noted.
The function $f$ is always greater than or equal to $1$.
Thus 
\[
\sup_{0 < t \leq 1} t \,\mu_{\mathbb{P}^1}\left\{ f(z) > t \right\}^{1/2}
= \sup_{0 < t \leq 1} t \,\mu_{\mathbb{P}^1}(\mathbb{P}^1)^{1/2}
= \sup_{0 < t \leq 1} t \sqrt{\pi}
= \sqrt{\pi}.
\]
On the other hand, for $t>1$:
\[
f(z) = \sqrt{|z|^{-2} + 1} > t
\quad\Leftrightarrow\quad
|z| \leq \left(t^2 - 1 \right)^{-1/2}
=: R_t.
\]
So then 
\begin{align*}
\sup_{t > 1} t \,\mu_{\mathbb{P}^1}\left\{ f(z) > t \right\}^{1/2}
= \sup_{t > 1} t \,\mu_{\mathbb{P}^1}(B(0,R_t))^{1/2}
= \sup_{t > 1} t \left(\pi t^{-2}\right)^{1/2}
= \sqrt{\pi}.
\end{align*}
Together this shows that $f(z)$ belongs to $\banach{2,\infty}(\mathbb{P}^1,\qat)$.
\end{proof}

Let us now return to $\Op{R}(0,\ci\lambda)$. We decompose this resolvent into the sum $\Op{R}(0,\ci\lambda)=\Op{R}_{1}(0,\ci\lambda)+\Op{R}_{2}(0,\ci\lambda)$ of operators which act on $\xi\in H^0(\mathbb{P}^1,\ban{2}{E})$ as the composition of first applying the inverse of~\eqref{dirac square} in Lemma~\ref{resolvent square dirac} and then $-(\qi\lambda)\ri$ and $-\diracJ\delbar{E}$, respectively:
\begin{gather}\label{integral kernel of resolvent 1}
\begin{aligned}
\big(\Op{R}_{1}(0,\ci\lambda)\xi\big)(z)=&-\int_\mathbb{C}\Func{K}_{\lambda}(z,z')\xi(z')\qi\lambda\frac{d\mu(z')}{1+|z'|^2},\\
\big(\Op{R}_{2}(0,\ci\lambda)\xi\big)(z)=&-\int_\mathbb{C}\diracJ\barpartial\Func{K}_{\lambda}(z,z')\xi(z')\frac{d\mu(z')}{1+|z'|^2}.
\end{aligned}
\end{gather}
Here $\barpartial$ is the derivative with respect to $\Bar{z}$. 
We further decompose both operators $\Op{R}_{j}(0,\ci\lambda)$ for $j=1,2$ and their sum $\Op{R}(0,\ci\lambda)$:
\begin{gather}\begin{aligned}\label{eq:near distant}
\Op{R}_{j}(0,\ci\lambda)&=\Op{R}_{j,\varepsilon\text{\scriptsize\textrm{-near}}}(0,\ci\lambda)+\Op{R}_{j,\varepsilon\text{\scriptsize\textrm{-far}}}(0,\ci\lambda),\\
\Op{R}(0,\ci\lambda)&=\Op{R}_{\varepsilon\text{\scriptsize\textrm{-near}}}(0,\ci\lambda)\hspace{2mm}+\Op{R}_{\varepsilon\text{\scriptsize\textrm{-far}}}(0,\ci\lambda),
\end{aligned}\end{gather}
where the integrands in~\eqref{integral kernel of resolvent 1} are either multiplied with the indicator function $\unity_{\{(z,z')\mid d_{\mathbb{P}^1}(z,z')<\varepsilon\}}$ or with the indicator function of $\unity_{\{(z,z')\mid d_{\mathbb{P}^1}(z,z')\ge\varepsilon\}}$.

The next lemma investigates the limit of $\Op{R}(0,\ci\lambda)$ for $\lambda\to\pm\infty$. We shall see that a suitable blowup is an analogous resolvent of the free Euclidean Dirac operator. Let us first introduce the scaling operators whose limits yields the blow up. For any $\lambda\in\mathbb{R}\setminus\{0\}$ and $1<p<\infty$ we introduce the following operators
\begin{align}\label{def isometries}
\Psi_{\lambda,p}&:\banach{p}(\mathbb{C},\qat)\to H^0(\mathbb{P}^1,\ban{p}{E}),
&\big(\Psi_{\lambda,p}\xi\big)_1(z)&=|\lambda|^{\frac{2}{p}}\xi(|\lambda|z)\big(1+|z|^2\big)^{\frac{2}{p}-\frac{1}{2}},
\end{align}
where the representative on $\SO_1$ (denoted with subscript $1$) suffices to define the section.
A direct calculation shows that $\Psi_{\lambda,p}$ is an bijective isometry with inverse
\begin{align}\label{inverse isometries}
\Psi_{\lambda,p}^{-1}&:H^0(\mathbb{P}^1,\ban{p}{E})\to\banach{p}(\mathbb{C},\qat),
&\big(\Psi_{\lambda,p}^{-1}\xi\big)(z)&=|\lambda|^{-\frac{2}{p}}\xi_1\big(\tfrac{z}{|\lambda|}\big)\big(1\!+\!\big|\tfrac{z}{\lambda}\big|^2\big)^{\frac{1}{2}-\frac{2}{p}}.\hspace{-2mm}
\end{align}
\begin{lemma}\label{lambda to infty}
For any $1<p<2$ the operators $\Psi_{\lambda,\frac{2p}{2-p}}^{-1}\Op{R}(0,\ci\lambda)\Psi_{\lambda,p}$ converge in $\mathcal{L}(\banach{p}(\mathbb{C},\qat),\banach{\frac{2p}{2-p}}(\mathbb{C},\qat))$ in the limit $\lambda\to\infty$ to the resolvent $\big((\ci)\ri-\qj\Bar{\partial}\big)^{-1}=\Op{I}_{\mathbb{C}}\big((\ci)\ri\,\Op{I}_{\mathbb{C}}-\qj\big)^{-1}$ of the free Euclidean Dirac operator $\qj\Bar{\partial}$ on $\mathbb{C}$. 
Additionally, $\Op{R}_{|\lambda|^{-\frac{1}{2}}\text{\scriptsize\textrm{-far}}}(0,\ci\lambda)$ converges to zero in $\mathcal{L}(H^0(\mathbb{P}^1,\ban{p}{E}),H^0(\mathbb{P}^1,\ban{\frac{2p}{2-p}}{E}))$.
\end{lemma}

\begin{proof}
Let us first show how to write the length of the section $\Op{R}_{1}(0,\ci\lambda)\xi$ as a convolution.
We take the length of~\eqref{integral kernel 2} and calculate:
\begin{align*}
|\Op{R}_{1}(0,\ci\lambda)\xi|_E(z)
&= \sqrt{1+|z|^2}\left| \int_{\mathbb{C}} \frac{1}{\sqrt{1+|z|^2}} \alpha_z\cdot\big(\fundsol_{\lambda}\big(|w|\big)\big)\xi(w)\frac{d\mu(w)}{1+|w|^2} \ci\lambda\right| \\
&\leq \int_{\mathbb{C}}  |\alpha_z\cdot\big(\fundsol_{\lambda}\big(|w|\big)\big)|\, |\xi(w)|\frac{d\mu(w)}{1+|w|^2} |\lambda| \\
&= \int_{\mathbb{C}} |\lambda|\, |\alpha_z\cdot\fundsol_{\lambda}|_E\big(|w|\big)\, |\xi|_E(w)\;d\mu_{\mathbb{P}^1}(w) \\
&= \int_{\mathbb{P}^1} |\lambda|\, \left| \fundsol_{\lambda} \right|_E (|\alpha_z^{-1}w|)\, |\xi|_E(w)\, \dmu_{\mathbb{P}^1}(w) \\
&= \left| \lambda \fundsol_{\lambda} \right|_E \hatast |\xi|_E.
\end{align*}
For $\Op{R}_{2}(0,\ci\lambda)$ we first need to give an explicit formula for its integral kernel:
\begin{align*}
\diracJ\barpartial\Func{K}_{\lambda}(z,w)
&= \qj (1+|z|^2) \frac{1}{1+z\Bar{w}} \barpartial \left[\fundsol_{\lambda}\left( \left|\frac{w-z}{1+\Bar{z}w}\right| \right) \right] \\
&= \qj \frac{1+|z|^2}{1+z\Bar{w}} \fundsol_{\lambda}'\left( \left|\frac{w-z}{1+\Bar{z}w}\right| \right) \frac{1}{2} \frac{w-z}{|w-z|} \frac{-(1+|w|^2)}{|1+\Bar{z}w|(1+\Bar{z}w)} \\
&= -\frac{1}{2} \qj \frac{(1+|z|^2)(1+|w|^2)}{|1+z\Bar{w}|^3} \fundsol_{\lambda}'\left( \left|\frac{w-z}{1+\Bar{z}w}\right| \right) \frac{w-z}{|w-z|} \\
&= -\frac{1}{2} \frac{1}{\sqrt{1+|z|^2}} \qj \alpha_z \cdot \left( (1+|w|^2) \fundsol_{\lambda}'\left( \left|w\right| \right) \frac{w}{|w|} \right).
\end{align*}
Thus $\Op{R}_{2}(0,\ci\lambda)\xi$ has an expression that is entirely similar to $\Op{R}_{1}(0,\ci\lambda)\xi$,
\begin{align*}
\Op{R}_{2}(0,\ci\lambda)\xi(z)
&= \frac{1}{2} \qj \int_\mathbb{C} \frac{1}{\sqrt{1+|z|^2}} \alpha_z \cdot \left( (1+|w|^2) \fundsol_{\lambda}'\left( \left|w\right| \right) \frac{w}{|w|} \right) \xi(w)\frac{d\mu(w)}{1+|w|^2},
\end{align*}
and by analogous working we have
\begin{align*}
|\Op{R}_{2}(0,\ci\lambda)\xi|_E(z)
\leq \frac{1}{2} \left| (1+|w|^2) \fundsol_{\lambda}'\left( \left|w\right| \right) \right|_E \hatast |\xi|_E.
\end{align*}
There is a trick to simplify the calculation.
We observe
\begin{equation}
\label{eq:supremum trick convolution}
\begin{aligned}
\left| \lambda \fundsol_{\lambda} \right|_E
&= \left| w \lambda \fundsol_{\lambda} w^{-1}\right|_E
\leq \left\| w \lambda \fundsol_{\lambda} \right\|_{\banach{\infty}(\mathbb{P}^1,\qat)} \,\left|w^{-1}\right|_E, 
\text{ and} \\
\left| (1+|w|^2) \fundsol_{\lambda}'\left( \left|w\right| \right) \right|_E
&= \left\| w(1+|w|^2) \fundsol_{\lambda}'\left( \left|w\right| \right) \right\|_{\banach{\infty}(\mathbb{P}^1,\qat)} \,\left|w^{-1}\right|_E.
\end{aligned}
\end{equation}
Due to Lemma~\ref{green}~(ii) and~(iv), the suprema are bounded by constants independent of $\lambda$.
Now, applying Lemma~\ref{lem:generalized convolution}, we see that the norms of $\Op{R}_{1}(0,\ci\lambda)$ and $\Op{R}_{2}(0,\ci\lambda)$ in $\mathcal{L}(H^0(\mathbb{P}^1,\ban{p}{E}),H^0(\mathbb{P}^1,\ban{\frac{2p}{2-p}}{E}))$ by the norm of $|w^{-1}|_E$ in the Lorentz space $\banach{2,\infty}(\mathbb{P}^1,\qat)$, since $\frac{2-p}{2p}+1=\tfrac{1}{p}+\tfrac{1}{2}$ and $\frac{2-p}{2p}\le \tfrac{1}{p} + \tfrac{1}{\infty}$.
That $|w^{-1}|_E$ belongs to $\banach{2,\infty}(\mathbb{P}^1,\qat)$ is exactly Lemma~\ref{green}~(v).
For the near and far components, since they are the same kernels but multiplied by a function that depends on the length with respect to the Fubini-Study metric, when the lengths of these operators are written as convolutions, the above two functions are multiplied by indicator function dependent on the distance from $0 \in \mathbb{P}^1$.

At this point we can already address the second statement of the lemma. 
With respect to the trick above of taking the supremum norm, for the far part with $\epsilon = |\lambda|^{-\frac{1}{2}}$ it suffices to take the supremum over $\mathbb{P}^1 \setminus B_{\mathbb{P}^1}(0,|\lambda|^{-\frac{1}{2}})$.
Since $|\lambda|e^{-|\lambda|^{-\frac{1}{2}}|\lambda|}$ converges to zero in the limit $\lambda\to\infty$, the estimates in Lemma~\ref{green}~(i), (iii), and~(iv) together imply that the supremum norms in~\eqref{eq:supremum trick convolution} converge to zero in this limit. 
This proves the second statement.

For the proof of the first statement we insert~\eqref{def isometries}-\eqref{inverse isometries} into~\eqref{integral kernel of resolvent 1}:
\begin{gather*}
\begin{aligned}
\big(\Psi_{\lambda,\frac{2p}{2-p}}^{-1}\Op{R}_{1}(0,\ci\lambda)\Psi_{\lambda,p}\xi\big)(z)
=&-\int_\mathbb{C} \left[\Func{K}_{\lambda}\big(\tfrac{z}{|\lambda|},\tfrac{z'}{|\lambda|}\big)\frac{\lambda}{|\lambda|}\frac{(1+|\frac{z}{\lambda}|^2)^{\frac{3}{2}-\frac{2}{p}}}{(1+|\frac{z'}{\lambda}|^2)^{\frac{3}{2}-\frac{2}{p}}} \right] \xi(z')\qi\, d\mu(z'),\\
\big(\Psi_{\lambda,\frac{2p}{2-p}}^{-1}\Op{R}_{2}(0,\ci\lambda)\Psi_{\lambda,p}\xi\big)(z)=&-\int_\mathbb{C} \left[\diracJ\barpartial\Func{K}_{\lambda}\big(\tfrac{z}{|\lambda|},\tfrac{z'}{|\lambda|}\big)\frac{1}{|\lambda|}\frac{(1+|\frac{z}{\lambda}|^2)^{\frac{3}{2}\!-\!\frac{2}{p}}}{\big(1+|\frac{z'}{\lambda}|^2)^{\frac{3}{2}\!-\!\frac{2}{p}}} \right] \xi(z')\,d\mu(z').
\end{aligned}
\end{gather*}
Here the exponent $-1$ of $|\lambda|$ is the sum of the corresponding exponents $\frac{2}{p}-\frac{2(2-p)}{2p}=1$ of $\Psi_{\lambda,p}$ and $\Psi_{\lambda,\frac{2p}{2-p}}$ minus the exponent $2$ in the transformation of measures $d\mu(z')=d\mu(\frac{z}{|\lambda|})=|\lambda|^{-2}d\mu(z)$ in~\eqref{integral kernel of resolvent 1}. Moreover, the exponent $\frac{2}{p}-\frac{3}{2}$ of $1+|\frac{z'}{\lambda}|^2$ is the sum of the negative of the corresponding exponent of $\Psi_{\lambda,p}$ and $-1$ coming form the weight in the measure $\frac{d\mu(z')}{1+|z'|^2}$ in~\eqref{integral kernel 2}. Finally the exponent $\frac{3}{2}-\frac{2}{p}$ of $1+|z|^2$ is the corresponding exponent $\frac{1}{2}-\frac{2(2-p)}{2p}$ of $\Psi_{\lambda,\frac{2p}{2-p}}$.

We wish to take the limit of each of these expressions as $\lambda \to \infty$.
Due to the compactness of $\mathbb{P}^1$, it suffices to prove the first statement on the subspace of all $\xi$ with support in an open ball.
Further, because of the invariance with respect to the action of $\mathbb{S}^3$ in Lemma~\ref{moebius action}, without loss of generality we take this ball to lie in $\SO_1$.
We may pass the limit through the integral, due to the bounds on the integrals above that are independent of $\lambda$.
We see then that the additional weights with powers of $1+|\frac{z}{\lambda}|^2$ and $1+|\frac{z'}{\lambda}|^2$ converge in the limit $\lambda\to\infty$ uniformly to $1$ and can be neglected.
We aim to control the limit of
\begin{align*}
\Func{K}_{\lambda}\big(\tfrac{z}{\lambda},\tfrac{w}{\lambda}\big)
\!&= \frac{1}{1+ \lambda^{-2} z\Bar{w}}\fundsol_{\lambda}\left( \lambda^{-1} \left|\frac{w-z}{1+\lambda^{-2}\Bar{z}w}\right| \right).
\end{align*}
On any compact subset of $z\in\mathbb{C}$ the expression $2|\lambda|\arctan(\frac{|z|}{|\lambda|})$ converges in the limit $\lambda\to\infty$ uniformly to $2|z|$.
Therefore we insert $r=\frac{r'}{\lambda}$ in~\eqref{explicit formula} and take the limit $\lambda\to\infty$:
\begin{align*}
\lim\limits_{\lambda\to\infty}\fundsol_\lambda(\lambda^{-1}r')&=\frac{2}{\pi}\lim\limits_{\lambda\to\infty}\int_{\frac{r'}{|\lambda|}}^{\infty}\frac{\sinh(|\lambda|(\pi-2\arctan(s)))}{\sinh(|\lambda|\pi)}\sqrt{\frac{1+(\frac{r'}{\lambda})^2}{s^2-(\frac{r'}{\lambda})^2}}ds\\
&=\frac{2}{\pi}\int_{r'}^{\infty}\lim\limits_{\lambda\to\infty}\frac{\sinh(|\lambda|(\pi-2|\lambda|\arctan(\frac{s'}{|\lambda|})))}{\sinh(|\lambda|\pi)}\sqrt{\frac{1+(\frac{r'}{\lambda})^2}{(\frac{s'}{|\lambda|})^2-(\frac{r'}{\lambda})^2}}\frac{ds'}{|\lambda|}\\
&=\frac{2}{\pi}\int_{r'}^\infty\frac{e^{-2s'}ds'}{\sqrt{s'^2-r'^2}}
=\frac{1}{2} \left(\frac{2}{\pi}\int_{0}^{r'} + \frac{2}{\pi}\int_{r'}^\infty\right) e^{-t}e^{-\frac{r'^2}{t}}\frac{dt}{t}.
\end{align*}
The last equality follows from two related substitutions
\begin{gather*}
t\!+\!\frac{r'^2}{t}=\frac{t^2\!+\!r'^2}{t}=2s'\quad\Leftrightarrow\quad(t-s')^2=s'^2-r'^2\quad\Leftrightarrow\quad t=s'\pm\sqrt{s'^2-r'^2}\\\text{with}\quad dt=\Big(1\pm\frac{s'}{\sqrt{s'^2-r'^2}}\Big)ds'=\frac{\pm s'+\sqrt{s'^2-r'^2}}{\sqrt{s'^2-r'^2}}ds'=\frac{\pm tds'}{\sqrt{s'^2-r'^2}}.
\end{gather*}
The substitution with the positive square root transforms the domain of integration $s'\in(r',\infty)$ to $t\in(r',\infty)$, where as the substitution with the negative square root results in $t\in(0,r')$ but oriented in the opposite direction. 
Although either substitution would give the desired integrand, the advantage of averaging both is that $r'$ no longer appears in the domain of integration.

We use $\fundsol_\infty(r')$ to denote this limit of $\fundsol_\lambda(\lambda^{-1}r')$. 
The limit of $\lambda^{-1}\diracJ\barpartial\Func{K}_{\lambda}\big(\tfrac{z}{\lambda},\tfrac{w}{\lambda}\big)$ as $\lambda \to \infty$ follows from that of $\Func{K}_{\lambda}$.
In summary, the limit of the blow-up gives the operator
\begin{align}\label{eucildean dirac}
 -\int_{\mathbb{C}} \fundsol_{\infty}(|z-w|)\xi(w)\ci + \qj\barpartial\fundsol_{\infty}(|z-w|)\xi(w) \;d\mu(w).
\end{align}
We will now show that this operator is indeed a resolvent of the free Dirac operator on $\mathbb{C}$, as claimed.
In analogy to the resolvent of the Dirac operator on $\mathbb{P}^1$ we determine the resolvent $\big((\ci)\ri-\qj\Bar{\partial}\big)^{-1}=\big(-(\ci)\ri-\qj\barpartial\big)\big(1+(\qj\barpartial)^2\big)^{-1}$ of the Dirac operator $\qj\barpartial$ on $\mathbb{C}$ in terms of the inverse of $1+(\qj\barpartial)^2=1-\frac{1}{4}\triangle$. Due to~\cite[Chapter~V~\S3.1 (26)]{St} this inverse is the convolution on $\banach{p}(\mathbb{C},\qat)$ with the function
\[
\fundsol_{\mathbb{C}}(|z|)=\frac{4}{4\pi\Gamma(1)}\int_0^\infty e^{-\frac{4\pi|z|^2}{\delta}}e^{-\frac{\delta}{4\pi}}\frac{d\delta}{\delta}=\frac{1}{\pi}\int_0^\infty e^{-t}e^{-\frac{|z|^2}{t}}\frac{dt}{t}.
\]
Here we replaced $|x|^2$ by $4|z|^2$ such that $\triangle_x=\partial\barpartial$ and the additional factor $4$ in the denominator changes the measure $d\mu(z)$ of the convolution into $d\mu(x)=4d\mu(z)$. Consequently the resolvent $\big((\ci)\ri\!-\!\qj\Bar{\partial}\big)^{-1}$ of Euclidean free Dirac operator on $\mathbb{C}$ is exactly\eqref{eucildean dirac}.
\end{proof}

For any $1<p<2$ the operator~\eqref{eucildean dirac} belongs to $\mathcal{L}(\banach{p}(\mathbb{C},\qat),\banach{\frac{2p}{2-p}}(\mathbb{C},\qat))$.
Let $S'_p$ denote the inverse of the norm of this operator.
In analogy to~\eqref{def:Sp-} we define
\begin{align}\label{def:Sp-'}
S'_{p,\varepsilon} =\inf\limits_{p'\in[p-\varepsilon,p]}S'_{p'}
\quad\text{and}\quad
{S'_p}^-&=\lim_{\epsilon\downarrow0}S'_{p,\varepsilon}.
\end{align}
The following corollary prepares the improvement of Lemma~\ref{weakly continuous} in Theorem~\ref{inverse holomorphic structure}:
\begin{corollary}\label{weakly continuous resolvent}
For every $1<p<2$ the constant ${S'_p}^->0$ in~\eqref{def:Sp-'} has the following property: For every $C'_p<{S'_p}^-$ and $\epsilon>0$ there exists $0<\delta<1-p$ such that $\diracJ\V\Op{R}(0,\ci\lambda)$ has in the operator space $\mathcal{L}(H^0(\mathbb{P}^1,\ban{p'}{E}),H^0(\mathbb{P}^1,\ban{p''}{E}))$ norm smaller than $1 - (1-C'_p({S'_p}^-)^{-1})/4$ for all $p-\delta<p''\le p'\le p$, $\lambda\in(-\infty,-\frac{1}{\delta})\cup(\frac{1}{\delta},\infty)$ and all $\V$ in
\begin{gather}\label{fubini study sets}
\{\V\in\pot{E}\mid\|\V|_{B_{\mathbb{P}^1}(x,\epsilon)}\|_2 \le C'_p\text{ for all }x\in\mathbb{P}^1\}.
\end{gather}
Here $\V|_{B_{\mathbb{P}^1}(x,\epsilon)}$ denotes the restriction to $B_{\mathbb{P}^1}(x,\epsilon)=\{y\in\mathbb{P}^1\mid d_{\mathbb{P}^1}(x,y)<\epsilon\}$.
\end{corollary}

\begin{proof}
The operator $\diracJ\V$ on $\SO_1$ acts as pointwise multiplication with the corresponding function $z \mapsto \qj (1+|z|^2)\V(z)$.
As we have used many times, the operator norm is bounded using Hölder's inequality.
Using for $q'=\frac{2p'}{2-p'}$ and the norms on $H^0(\mathbb{P}^1,\ban{p'}{E})$:
\[
\|\diracJ\V \xi\|_{p'}
= \|(1+|z|^2)\V(z) |\xi|_E \|_{\banach{p'}(\mathbb{P}^1,\qat)}
\leq \|(1+|z|^2)\V(z)\|_{\banach{2}(\mathbb{P}^1,\qat)} \|\xi \|_{q'}
\]
and
\[
\|(1+|z|^2)\V(z)\|_{\banach{2}(\mathbb{P}^1,\qat)}
= \int_{\mathbb{C}} |\qj(1+|z|^2)\V|^2 \,\dmu_{\mathbb{P}^1}(z)
= \int_{\mathbb{C}} |\V|^2 \dmu(z)
= \|\V\|_{\banach{2}(\mathbb{C},\qat)}^2,
\]
which is the usual definition of the norm of $\V\in\pot{E}$ in Definition~\ref{def:potentials}.
Since $\mathbb{P}^1$ is covered by $N_\epsilon$\=/many balls of radius $\epsilon$, the square of the norms $\|\V\|_2^2$ of the elements of~\eqref{fubini study sets} are bounded by $N_\epsilon (C'_p)^2$.
With this explanation in mind, we aim to control $\diracJ\V\Op{R}(0,\ci\lambda)$ in $\mathcal{L}(H^0(\mathbb{P}^1,\ban{p'}{E}),H^0(\mathbb{P}^1,\ban{p''}{E}))$ using the foregoing Lemma~\ref{lambda to infty} and $\|V\|_{\banach{2}(\mathbb{C},\qat)}$.

For given $C'_p<{S'_p}^-$ and $\epsilon>0$ we now give three conditions on $\delta$; each condition corresponds to an approximation of $\|\diracJ\V\Op{R}(0,\ci\lambda)\|$ in $\mathcal{L}(H^0(\mathbb{P}^1,\ban{p'}{E}),H^0(\mathbb{P}^1,\ban{p''}{E}))$.
Because $C'_p({S'_p}^-)^{-1}<1$ we have a small amount of wiggle room to play with.
As long as each of these three approximations results in an increase of $\|\diracJ\V\Op{R}(0,\ci\lambda)\|$ by an amount smaller than $({S'_p}^--C'_p)/(4{S'_p}^-)$ the inequality $\|\diracJ\V\Op{R}(0,\ci\lambda)\|<1 - (1-C'_p({S'_p}^-)^{-1})/4$ follows.
Hence it suffices to show that each of these approximations has an arbitrary small error, if $\delta$ is sufficiently small.
Rather an explicitly calculating the threshold $\delta$, we show that the approximations of $\Op{R}(0,\ci\lambda)$ converge in the limit $\delta\downarrow0$.
Since the global $\banach{2}$-norms $\|\V\|_2$ of the elements of~\eqref{fubini study sets} are bounded, this implies in all three cases also an convergence of the norms of $\diracJ\V\Op{R}(0,\ci\lambda)$.
Now we explain these three approximations separately:

\begin{enumerate}
\item[1.] Approximate $\Op{R}(0,\ci\lambda)$ by $\Op{R}_{\epsilon-\text{near}}(0,\ci\lambda)$: \\
Finally we approximate $\Op{R}(0,\ci\lambda)$ by $\Op{R}_{\epsilon-\text{near}}(0,\ci\lambda)$, which is implied due the second statement of the foregoing Lemma~\ref{lambda to infty}, since $\lambda^{-1/2}$ will fall below $\epsilon$ in the limit $\lambda \to \infty$.
As a consequence we may take pointwise multiplication by $\diracJ\V$ into the integral kernels in~\eqref{integral kernel of resolvent 1}.
Then the result is unchanged if we replace $\diracJ\V$ by $\unity_{\{(z,z')\mid d_{\mathbb{P}^1}(z,z')<\varepsilon\}}\diracJ\V$.
Hence in this case the $\banach{2}$-norms of the potentials $\V$ in~\eqref{fubini study sets} are effectively bounded by $C'_p$.
\item[2.] Approximate $\Op{R}(0,\ci\lambda)$ by the Euclidean free Dirac operator:\\
Due to the first statement in the foregoing Lemma~\ref{lambda to infty}, we can approximate the norm of $\Op{R}(0,\ci\lambda)$ in $\mathcal{L}(H^0(\mathbb{P}^1,\ban{p'}{E}),H^0(\mathbb{P}^1,\ban{q''}{E}))$, with $q''=\frac{2p''}{2-p''}$, by the norm of the resolvent~\eqref{eucildean dirac} of the free Dirac operator on $\mathbb{C}$ in $\mathcal{L}(\banach{p'}(\mathbb{C},\qat),\banach{q''}(\mathbb{C},\qat))$.
Moreover, for $\delta<p-1$ this convergence is uniform $p'$ and $p''$ for $p-\delta<p''\le p'\le p$, since the constants of the estimates in the interpolation Lemma~\ref{lem:generalized convolution} are uniformly bounded for these exponents.
\item[3.] Approximate the norm of the Euclidean free Dirac operator independent of $p'$ and $p''$:\\
The second upper bound on $\delta$ allows to replace $p-\delta<p''\le p'\le p$ by $p''=p=p'$.
By definition of ${S'_p}^-$ the difference $S'_{p'}-{S'_p}^-$ has an arbitrary small upper bound for sufficiently small $\delta$.
Furthermore, by a standard application of Hölder's inequality the norms of the embeddings $H^0(\mathbb{P}^1,\ban{p'}{E})\hookrightarrow H^0(\mathbb{P}^1,\ban{p''}{E})$ converge to $1$ for $p''\uparrow p'$.
Therefore the difference of the norms of the resolvent~\eqref{eucildean dirac} of the Euclidean free Dirac operator on $\mathbb{C}$ in $\mathcal{L}(\banach{p}(\mathbb{C},\qat),\banach{\frac{2p}{2-p}}(\mathbb{C},\qat))$ and $\mathcal{L}(\banach{p'}(\mathbb{C},\qat),\banach{\frac{2p''}{2-p''}}(\mathbb{C},\qat))$ become arbitrary small, if $\delta$ is sufficiently small.
\end{enumerate}
Each of these limits implies the existence of a value of $\delta$, such that for all smaller $\delta$ the error of the approximation is less than $({S'_p}^--C'_p)/(4{S'_p}^-)$.
Therefore if $\delta$ is smaller than all three thresholds, then all approximations are sufficiently good.
The statement of the corollary follows.
\end{proof}

Let us now derive a simple criterion (stated in Lemma~\ref{bounded point measures}) that describes when the tail of a sequence of potentials $(\V_n)_{n\in\mathbb{N}}$ in $\pot{E}$ on $\mathbb{P}^1$, whose sequence of measures $|\V_n|^2$ converge weakly to $\sigma$, is contained in one of the weakly compact sets described in Corollary~\ref{weakly continuous resolvent}.
This criterion holds if $\sigma$ does not contain point masses larger or equal to the square of ${S'_p}^-$~\eqref{def:Sp-'}.
Due to $\sigma$\=/additivity the measure of any point $z$ with respect to $\sigma$ is the limit of the sequence of measures of $B(z,\frac{1}{n})=\{z\}\cup\bigcup_{m\ge n}\big(B(z,\frac{1}{m})\setminus B(z,\frac{1}{m+1})\big)$, and if this is non-zero we say that the measure has a point mass there.
Hence we may cover $\mathbb{P}^1$ by open sets whose measures with respect to $\sigma$ are smaller than $({S'_p}^-)^2$.
Due to the compactness of $\mathbb{P}^1$ this open cover has a finite subcover $\mathfrak{U}$.
Without loss of generality, we may suppose that $\mathbb{P}^1 \not\in \mathfrak{U}$.
Then the function 
\[
\mathbb{P}^1 \to \mathbb{R}_+\,, z \mapsto \max_{\SU\in\mathfrak{U}} \mathrm{dist}(z, \mathbb{P}^1 \setminus \SU)
\]
is continuous and therefore achieves a positive minimum $2\varepsilon$.
Thus every disc $B(z,2\varepsilon)$ in $\mathbb{P}^1$ is contained in some $\SU \in \mathfrak{U}$, and then $B(z,\varepsilon) \subset \{ z\in \SU \mid \mathrm{dist}(z,\mathbb{P}^1 \setminus \SU) > \varepsilon\} =: \SU'$ holds.
Let $f_\SU: \mathbb{P}^1 \to [0,1]$ be a continuous function with $\mathrm{supp}(f_\SU) \subset \SU$ and $f_\SU|_{\SU'} = 1$.
This implies $\int_\X f_\SU\,d\sigma < \int_\SU d\sigma <({S'_p}^-)^2$.
By weak convergence we also have $\int_\X f_\SU \,|\V_n|^2 <({S'_p}^-)^2$ for all but finitely many $n$.
Because every ball $B(z,\varepsilon)$ is contained in $\SU'$ for some $\SU\in\mathfrak{U}$, it follows that there exist $C'_p<{S'_p}^-$ and, for all $n$ sufficiently large and for all $z\in\mathbb{P}^1$ that $\|\V_n|_{B(z,\varepsilon)}\|_2\leq C'_p$.
We summarize the discussion in the following lemma: %
\begin{lemma}\label{bounded point measures}
Let $(\V_n)_{n\in\mathbb{N}}$ be on $\mathbb{P}^1$ a sequence in $\pot{E}$ whose sequence of measures $(|\V_n|^2)_{n\in\mathbb{N}}$ converges weakly to $\sigma$.
If $\sigma\big(\{x\}\big)$ is smaller than the square of ${S'_p}^-$~\eqref{def:Sp-'} for all $x\in\mathbb{P}^1$, then all but finitely many $\V_n$ belong for some $C'_p<{S'_p}^-$ and some $\epsilon>0$ to the set~\eqref{fubini study sets}.\qed
\end{lemma}
Conversely, if a weak limit of a sequence of finite Baire measures $|\V_n|^2$ contains a point mass greater or equal to the square of ${S'_p}^-$ at some point $z_0$, then
\[
\limsup\limits_{n\to\infty} \|\V_n|_{B(z,\varepsilon)}\|_2\geq{S'_p}^-\quad\text{for every}\quad\varepsilon>0.
\]

Now we can state and prove the main theorem of this chapter.
\begin{theorem}\label{inverse holomorphic structure}
For $p$ specified below, let $(\V_n)_{n\in\mathbb{N}}$ be on $\mathbb{P}^1$ a weakly convergent sequence in $\pot{E}$ with weak limit $\V$ whose sequence of measures $(|\V_n|^2)_{n\in\mathbb{N}}$ converges weakly to $\sigma$.
If in addition $\sigma\big(\{x\}\big)$ is smaller than the square of ${S'_p}^-$~\eqref{def:Sp-'} for all $x\in\mathbb{P}^1$ and $\dim H^0(\mathbb{P}^1,\Q{E,\V})=0$, then the following statements hold:
\begin{enumeratethm}
\item for $1<p<2$ and $1<q<\frac{2p}{2-p}$ the sequence $\big(\Op{R}(\V_n,0)\big)_{n\in\mathbb{N}}$ converges in $\mathcal{L}(H^0(\mathbb{P}^1,\ban{p}{E}),H^0(\mathbb{P}^1,\ban{q}{E}))$ to $\Op{R}(\V,0)$.
\item For $p=\frac{4}{3}$ and large $n$ the operators $(\unity-\V_n\Op{I}_{\mathbb{C}})$ from $\banach{\frac{4}{3}}(\mathbb{C},\qat)$ to $\banach{\frac{4}{3}}(\mathbb{C},\qat)$ are invertible with uniformly bounded inverse.
\end{enumeratethm}
\end{theorem}
\begin{proof}
By the forgoing Lemma~\ref{bounded point measures} all but finitely many $\V_n$ belong to the set~\eqref{fubini study sets} for some $C'_p<{S'_p}^-$ and some $\epsilon>0$.
Hence due to Corollary~\ref{weakly continuous resolvent} there exists $\delta>0$ such that the following Neumann series converges converges for all $\lambda\in(-\infty,-\frac{1}{\delta})\cup(\frac{1}{\delta},\infty)$ and large $n$, and by~\eqref{eq:first resolvent formula} the limit is the resolvent
\[
\Op{R}(\V_n,\ci\lambda)
= \Op{R}(0,\ci\lambda) \left[\unity +  \diracJ\V_n \Op{R}(0,\ci\lambda) \right]^{-1}
= \Op{R}(0,\ci\lambda)\sum_{l=0}^{\infty}\big[-\diracJ\V_n\Op{R}(0,\ci\lambda)\big]^{l}.
\]
Now the same arguments as in Lemma~\ref{weakly continuous} show that for the given $1<p<2$ and $1<q<\frac{2p}{2-p}$ and the same $\lambda\in(-\infty,-\frac{1}{\delta})\cup(\frac{1}{\delta},\infty)$ the map $\V'\mapsto\Op{R}(\V',\ci\lambda)$ from the set~\eqref{fubini study sets} to $\mathcal{L}(H^0(\mathbb{P}^1,\ban{p}{E}),H^0(\mathbb{P}^1,\ban{q}{E}))$ is continuous with respect to the weak topology on~\eqref{fubini study sets} and the norm topology on $\mathcal{L}(H^0(\mathbb{P}^1,\ban{p}{E}),H^0(\mathbb{P}^1,\ban{q}{E}))$: Here we use that $\Op{R}(0,\ci\lambda)$ is a compact operator in $\mathcal{L}(H^0(\mathbb{P}^1,\ban{p'}{E}),H^0(\mathbb{P}^1,\ban{q''}{E}))$ for $1<p'\le p$ and $1<q''<\frac{2p'}{2-p''}$ by Lemma~\ref{dirac resovent}.
Moreover, due to~\cite[Theorem~II.5.11]{LT}, these Banach spaces $H^0(\mathbb{P}^1,\ban{q''}{E})$ have a Schauder basis and the approximation property.
Therefore and all compact operators into $H^0(\mathbb{P}^1,\ban{q''}{E})$ are norm-limits of finite rank operators (compare~\cite[Section~I.1.a]{LT}).
If all resolvents in the Neumann series are approximated by finite rank operators, then all terms of the Neumann series are approximated by functions which are weakly continuous with respect to $\V'$ and which have values in the compact operators from $H^0(\mathbb{P}^1,\ban{p}{E})$ into $H^0(\mathbb{P}^1,\ban{q}{E})$.
The uniform limit of weakly continuous functions on the weakly compact sets~\eqref{fubini study sets} is again a weakly continuous function~\cite[Theorem~IV.8]{RS1}.
In particular, the sequence $(\Op{R}(\V_n,\ci\lambda))_{n\in\mathbb{N}}$ converges in $\mathcal{L}(H^0(\mathbb{P}^1,\ban{p}{E}))$ for all $\lambda\in(-\infty,-\frac{1}{\delta})\cup(\frac{1}{\delta},\infty)$ to $\Op{R}(\V,\ci\lambda)$.
As a preparation for Serre duality~\ref{Serre duality} we proved that $\delbar{E}-\V$ is a Fredholm operator from $H^0(\mathbb{P}^1,\sob{1,p}{E})$ to $H^0(\mathbb{P}^1,\forms{0,1}{}\ban{p}{E})$, whose index vanishes by the Riemann-Roch theorem~\ref{riemann roch}.
Hence this operator is invertible.
Due to the first resolvent formula~\eqref{eq:first resolvent formula}, $\frac{1}{\ci\lambda}$ belongs to the resolvent set of $\Op{R}(\V,\ci\lambda)$, so for sufficiently large $n$ the operator $\unity\!-\!\ci\lambda \Op{R}(\V_n,\ci\lambda)$ is invertible and its inverse converges in $\mathcal{L}(H^0(\mathbb{P}^1,\ban{p}{E}))$.
Moreover, by two applications of~\eqref{eq:first resolvent formula}, the sequence of operators
\begin{align*}
\left(\unity+\diracJ\V_n\Op{R}(0,0)\right)^{-1}
&=\diracJ\delbar{E}\Op{R}(\V_n,0)\\
&=\diracJ\delbar{E}\Op{R}(\V_n,\ci\lambda) \left[\unity\!-\!\ci\lambda \Op{R}(\V_n,\ci\lambda) \right]^{-1} \\
&=\diracJ\delbar{E}\Op{R}(0,\ci\lambda) \left[\unity\!-\!\diracJ \V_n \Op{R}(0,\ci\lambda) \right]^{-1} \left[\unity\!-\!\ci\lambda \Op{R}(\V_n,\ci\lambda) \right]^{-1}\\
&=\left[\ci\lambda\Op{R}(0,\ci\lambda)\!-\!\unity\right]\left[\unity\!-\!\diracJ \V_n \Op{R}(0,\ci\lambda) \right]^{-1} \left[\unity\!-\!\ci\lambda \Op{R}(\V_n,\ci\lambda) \right]^{-1}\\
\end{align*}
converges to $(\unity+\diracJ\V\Op{R}(0,0))^{-1}$.
Another application of~\eqref{eq:first resolvent formula} yields $\Op{R}(\V,0)=\Op{R}(0,0)(\unity+\diracJ\V\Op{R}(0,0))^{-1}$, which implies the statement~(i).
To prove the second statement, we first explain what makes the exponent $\tfrac{4}{3}$, and its dual exponent $4$, special.
The exponent $4$ has the property that $\xi \mapsto \xi_1$, which is to say taking the representative on $\SO_1$, conceived as an operator from $H^0(\mathbb{P}^1,\ban{4}{E})$ to $L^4(\mathbb{C},\qat)$, is an isometry (compare~\eqref{eq:Lpnorminvariance2}). 
Similarly, the map $f \mapsto \qj(1+|z|^2)f(z)$ can be understand as defining a section of $E$ by its representative on $\SO_1$, and is an isometry from $L^{\frac{4}{3}}(\mathbb{C},\qat)$ to $H^0(\mathbb{P}^1,\ban{\frac{4}{3}}{E})$.
We might call this map $\diracJ$.
Thus we have an isometric intertwining
\[
\begin{tikzcd}
H^0(\mathbb{P}^1,\ban{4}{E}) \arrow[r, "\xi_1"] \arrow[d, "\diracJ\delbar{E}", swap] & L^4(\mathbb{C},\qat) \arrow[d, "\barpartial"] \\
H^0(\mathbb{P}^1,\ban{\frac{4}{3}}{E}) & L^{\frac{4}{3}}(\mathbb{C},\qat) \arrow[l, "\diracJ"]
\end{tikzcd}
\qquad
\begin{tikzcd}
H^0(\mathbb{P}^1,\ban{4}{E}) \arrow[r, "\xi_1"] & L^4(\mathbb{C},\qat) \\
H^0(\mathbb{P}^1,\ban{\frac{4}{3}}{E}) \arrow[u, "{- \Op{R}(0,0)}"] & L^{\frac{4}{3}}(\mathbb{C},\qat) 
\arrow[l, "\diracJ"] \arrow[u, "{\Op{I}_{\mathbb{C}}}", swap]
\end{tikzcd}
\]
The diagram on the left holds in the sense of unbounded linear operators, but serves to illustrate the truth of the diagram on the right.
Note the minus sign in $\Op{R}(0,0)$ arising from Definition~\ref{def:resolvent}.
We have an analogous diagram for the action of a potential.
Finally, we calculate
\begin{align*}
(\unity+\diracJ\V_n\Op{R}(0,0))^{-1}
&= (\unity - \diracJ\V_n\Op{I}_{\mathbb{C}}\diracJ^{-1})^{-1}
= \diracJ (\unity - \V_n\Op{I}_{\mathbb{C}})^{-1} \diracJ^{-1}.
\end{align*}
This shows that $\diracJ$ isometrically intertwines $(\unity+\diracJ\V_n\Op{R}(0,0))^{-1}$ on $\banach{\frac{4}{3}}(\mathbb{P}^1,\mathbb{H})$ and $(\unity-\V_n\Op{I}_{\mathbb{C}})^{-1}$ on $\banach{\frac{4}{3}}(\mathbb{C},\mathbb{H})$.
This proves the Theorem.
\end{proof}

\chapter{Local Limits of Sequences of Holomorphic Structures}
\label{chapter:weak limits}
In this chapter we shall further improve our understanding of the failure of the weak continuity of the map $\V\mapsto(\barpartial-\V)^{-1}$. In the preceding chapter we established on $\mathbb{P}^1$ that such a failure can only happen if the energy concentrates at single points with an energy above a positive threshold. 
In this chapter we shall prove in Proposition~\ref{disc limits} three statements:
\begin{enumerate}
\item The weak continuity can fail locally on a disc only if such a concentration at single points exceeds a positive constant.
\item In the case of such massive concentration of energy a blow up is used to establish the following criterion:
The failure of weak continuity happens only if at the point of concentrated energy the blow up yields a holomorphic structure on the exceptional fiber $\mathbb{P}^1$ with non-trivial kernel. 
Hence the positive constant in~(i) is the Willmore energy of a round sphere, $4\pi$.
\item If we increase the degree of the underlying complex line bundle by additional poles at the points of concentration of energy, then the failure of the weak continuity is repaired.
Here the orders of the poles are chosen to be the concentrated energy divided by $4\pi$, rounded down to the nearest integer.
At the end of this part we shall see in Chapter~\ref{chapter:singular holomorphic sheaves}, that this means that the corresponding holomorphic sections can evolve only poles of such orders at the points of massive concentration of energy.
\end{enumerate}
The proof of this proposition is divided into three steps.
In Step~A the sequence of potentials is decomposed into a `singular' part whose limit is a discrete measure and a regular part that does not contain `large' point masses.
Steps~B and~C then deal with the convergence of the resolvents near large concentrations of mass in the two different cases whether the criterion in~(ii) is fulfilled or not.
There is then a corollary that controls the limits of related holomorphic structures, namely the $\mathbb{C}$\=/conjugate and Serre dual.

The key idea of the proof of Proposition~\ref{disc limits} is to embed the disc into $\mathbb{P}^1$.
The advantage of $\mathbb{P}^1$ over the disc is that the former is a compact space, which means that the kernels of the holomorphic structure is finite dimensional.
On the other hand, the conformal group of $\mathbb{P}^1$ is still large enough to prevent the concentration of energy, a technique that is known as `concentration--compactness'.
Throughout this chapter, $z$ refers to the coordinate on $\mathbb{D}$ or $\mathbb{P}^1$ and $z_l$ refers to a selection of points (not different charts).
Further, we denote the operator $\Op{I}_{\Omega,\V}$ from~\eqref{eq:resolvent} instead as $\Op{I}_{\Omega}(\V)$, as otherwise the subscripts become too crowded.

\begin{proposition}
\label{disc limits}
\index{Potential}
Let $(\U_n)_{n\in\mathbb{N}}$ be a weakly convergent sequence of potentials in $\banach{2}(\mathbb{D},\qat^-)$ with weak limit $\U$ such that the corresponding sequence of measures $|\U_n|^2$ converges weakly to some measure $\sigma$ whose restriction to the annulus $\mathbb{D}\setminus B(0,\frac{1}{2})$ obeys the condition of Lemma~\ref{bounded point measures} for $p=\frac{4}{3}$.
Then there exists a sequence of polynomials $(P'_n)_{n\in\mathbb{N}}$ of fixed degree with the following properties:
\begin{enumeratethm}
\item The sequences $(P'_n)_{n\in\mathbb{N}}$ and $(P_n'^{-1})_{n\in\mathbb{N}}$ converge uniformly on compact subsets of the annulus $\{z\in\mathbb{D}\mid\frac{1}{2}<|z|<1\}$.
\item For any $1<q<4$ the sequence $(\Op{I}_{\mathbb{D}}(P'_n\U_n P_n'^{-1}))_{n\in\mathbb{N}}$ converges in the normed space of bounded operators from $\banach{\frac{4}{3}}(\mathbb{D},\qat)$ to $\banach{q}(\mathbb{D},\qat)$.
\end{enumeratethm}
\end{proposition}
We prepare the proof of this proposition in three steps.

\proofstep{Step A:} \textit{The decomposition of the sequence of potentials.}
We describe a procedure how to decompose a subsequence of the sequence of potentials $\U_n$.
This decomposition depends only on the measures $|\U_n|^2$.
Its purpose is to prepare the investigation of the weak discontinuities of the sequence of holomorphic structures in Steps~B and~C.
Due to Lemma~\ref{bounded point measures} the only possible causes of discontinuities are point masses contained in $\sigma$ which have mass greater or equal to some positive constant.
Because $\|\U_n\|_{\banach{2}(\mathbb{D},\qat^-)}$ is bounded, there are only finitely many such points.

This decomposition is based on a sequence $(\{ \Set{D}_{n,l} \mid l\in L\})_{n\in \mathbb{N}}$ of sets of finitely many disks contained in $\mathbb{D}$ where the index set $L$ is independent of $n$.
For each $n\in \mathbb{N}$, any two of these disks are either disjoint or one of them is a proper subset of the other.
The strict inclusion relation defines a strict partial ordering $<$ on the index set $L$ which is also independent of $n \in \mathbb{N}$.
For each $n$, these disks decompose $\mathbb{D}$ into the disjoint union
\begin{equation}\label{eq:decomposition domains}\begin{gathered}
\mathbb{D} = \Set{D}_{\reg,n} \;\cup\; \bigcup_{l\in L} \Set{D}_{\sing,n,l}
\quad\text{with the sets}\\
\Set{D}_{\reg,n}=\mathbb{D} \setminus \bigcup_{l\in L} \Set{D}_{n,l}
\quad\text{and}\quad
\Set{D}_{\sing,n,l} = \Set{D}_{n,l} \setminus \bigcup_{l' < l} \Set{D}_{n,l'} \;.
\end{gathered}\end{equation}
This decomposition induces a decomposition of the potential $\U_n$ into a sum
\begin{equation}\label{eq:decomposition}\begin{gathered}
\U_n=\U_{\reg,n}+\sum_{l\in L} \U_{\sing,n,l}
\quad\text{with}\\
\U_{\reg,n}= \U_n \, \unity_{\Set{D}_{\text{reg},n}}\quad\text{and}\quad
\U_{\sing,n,l} = \U_n \, \unity_{\Set{D}_{\text{sing},n,l}} \;,
\end{gathered}\end{equation}
where $\unity_X$ denotes the indicator function of $X$. 
In particular, the decomposition of $U_n$ does not change if the sets $\Set{D}_{n,l}$ are changed by sets of Lebesgue measure zero.

By conceiving the potentials on $\mathbb{D}$ as potentials on $\mathbb{P}^1$ vanishing outside of $\mathbb{D}$, the Möbius group of $\mathbb{P}^1$ acts on the space of these potentials.
We shall use this action to improve the convergence properties of the potentials $\U_n$ so as to meet the hypotheses of Lemma~\ref{bounded point measures}.
The action of the Möbius group $SL(2,\mathbb{C})/\mathbb{Z}_2$ on $\mathbb{P}^1$
\[
SL(2,\mathbb{C})\ni\begin{pmatrix}
a & b\\
c & d
\end{pmatrix}:\mathbb{P}^1\to\mathbb{P}^1,z\mapsto\frac{az+b}{cz+d},
\]
induces a representation of the Möbius group on the Hilbert space $\banach{2}(\mathbb{C})$ of square-integrable potentials of the spin bundle on $\mathbb{P}^1$.
Due to the global Iwasawa decomposition~\cite[Chapter~VI Theorem~5.1]{Heg}, the Möbius group splits as a product $KAN$, where $K=\mathrm{PSU}(2)$ is the isometry group of $\mathbb{P}^1$, $A$ consists of the real scaling transformations $z\mapsto\exp(t)z$ with $t\in \mathbb{R}$, and $N$ is the group of translations $z\mapsto z-z_0$ with $z_0\in \mathbb{C}$.
We have already seen in Lemma~\ref{moebius action} that the action of $K$ is by isometries.

The desired properties for the sequences of disks $(\Set{D}_{n,l})_{n\in \mathbb{N}}$ and sequences of corresponding Möbius transformations $(g_{n,l})_{n\in \mathbb{N}}$ from the semidirect product $AN$, both for all $l\in L$,  are gathered in the following conditions.
Recall that the sequences of disks induce a decomposition of the potentials $\U_n$ by Equation~\eqref{eq:decomposition}.
\begin{description}
\item[Decomposition (i)] The sequences $(\U_{\reg,n})_{n\in\mathbb{N}}$ and $((\moebius_{n,l}^{-1})^{\ast}\U_{\sing,n,l})_{n\in\mathbb{N}}$ with $l\in L$ obey the assumption of Lemma~\ref{bounded point measures}.
\item[Decomposition (ii)] For every $l\in L$, the weak limit of the measures $|\U_{\sing,n,l}|^2$ is a discrete measure of mass greater or equal to $S_{\frac{4}{3}}^{-2}$.
\item[Decomposition (iii)] For every maximal element $l\in L$, the weak limit of the measures $\sum_{l' \leq l}|\U_{\sing,n,l'}|^2$ is equal to a measure supported at some $z_l \in \mathbb{D}$ whose mass is greater or equal to $S_{\frac{4}{3}}^{-2}$ and equal to $\sigma(\{z_l\})$.
\item[Decomposition (iv)]
For every non-maximal element $l\in L$, there exists  a (necessarily unique) minimum $l''$ of $\{l' \in L \mid l < l'\}$, and the weak limit of the measures $\sum_{l' \leq  l}|(\moebius_{n,l''}^{-1})^{\ast} \U_{\sing,n,l'}|^2$ is equal to a measure supported at some $z_l\in \mathbb{P}^1$ whose mass is greater or equal to $S_{\frac{4}{3}}^{-2}$ and equal to the measure of $\{z_l\}$ with respect to the weak limit of the measures $|(g_{n,l''}^{-1})^* \U_{\sing,n,l''}|^2$.
\item[Decomposition (v)]
There exists an $N_l$ such that $\bigcap_{n>N_l} g_{n,l}^{-1} [\Set{D}_{n,l}]$ contains an open set $\SO_l \subset \mathbb{P}^1$.
\end{description}

We construct a decomposition with these properties for a subsequence of the $\U_n$ by applying two lemmas in alternation.
\begin{lemma}\label{lem:split-off-disks}
For any sequence of square-integrable potentials $(\U_n)_{n\in \mathbb{N}}$ on $\mathbb{P}^1$ such that the measures $|\U_n|^2$ converges weakly to the finite measure $\sigma$, let $z_1,\ldots,z_M$ be the elements of the finite set $\{z\in\mathbb{P}^1\mid\sigma(\{z\})\ge S_{\frac{4}{3}}^{-2}\}$.
Then there exist sequences of positive radii $(r_{n,1},\ldots,r_{n,M})_{n\in\mathbb{N}}$ converging to $0$ with the following properties:
\begin{enumeratethm}
\item For any $n\in\mathbb{N}$ the balls $B(z_1,r_{n,1}),\ldots,B(z_m,r_{n,M})$ are disjoint.
\item For any $m \in \{1,\dotsc,M\}$ the sequence of measures $(|\unity_{B(z_m,r_{n,m})}\U_n|^2)_{n\in\mathbb{N}}$ converges weakly to the single point measure at $z_m$ with mass $\sigma(\{z_m\})$.
\end{enumeratethm}	
\end{lemma}
\begin{proof}
For small $r>0$ the balls $B(z_{m},r)$ are pairwise disjoint and have measures greater or equal to $S_{\frac{4}{3}}^{-2}$ with respect to $\sigma$.
Hence we may choose monotonically decreasing sequences of positive radii with the desired properties.
\end{proof}
\begin{lemma}\label{existence Moebius}
For any bounded sequence $(\U_n)_{n\in\mathbb{N}}$ of square-integrable potentials on $\mathbb{P}^1$ there exists a sequence of Möbius transformations $(\moebius_n)_{n\in\mathbb{N}}$ in the subgroup $AN$ such that a subsequence of $(\moebius_n^{-1})^{\ast}\U_n$ either obeys the assumption of Lemma~\ref{bounded point measures} or the subsequence of measures $|(\moebius_n^{-1})^{\ast}\U_n|^2$ has a weak limit with at least two point masses greater or equal to $S_{\frac{4}{3}}^{-2}$.
\end{lemma}
\begin{proof}
For any square-integrable potential $\U$ on $\mathbb{P}^1$, let $r_\U$ denote the function on $\mathbb{P}^1$ which associates to every $x\in\mathbb{P}^1$ the maximal radius $r \leq \pi$ of the balls $B(x,r)$ having measure not larger than $S_{\frac{4}{3}}^{-2}$ with respect to the measure $|\U|^2$.
If $d(x,x')$ denotes the distance between $x,x'\in\mathbb{P}^1$, then $B(x,r_U(x))\subset B(x',r_U(x')+d(x,x'))$ and $B(x',r_U(x'))\subset B(x,r_u(x)+d(x,x'))$ implies $|r_U(x)-r_U(x')|\le d(x,x')$.
Hence this function is continuous and attains the minimum $r_{\min}(\U)$ at some points of $\mathbb{P}^1$.

We claim that there exists a Möbius transformation $\moebius$ maximizing the function $h \mapsto r_{\min}((h^{-1})^{\ast}\U)$ on the Möbius group.
Let $h_n$ be a maximizing sequence of Möbius transformations, i.e.\ the limit of the sequence $r_{\min}((h_n^{-1})^{\ast}\U)$ is equal to the supremum of the image of the former function.
If $h$ belongs to the subgroup $K$ of isometries of $\mathbb{P}^1$, then $r_{\min}((h^{-1})^{\ast}\U)=r_{\min}(\U)$.
Therefore the sequence $h_n$ may be chosen in the semidirect product $AN$ of the real scaling transformations $z\mapsto\exp(t)z$ with the translations $z\mapsto z-z_0$.
Now we show that if the values of $t_n\in\mathbb{R}$ or of $z_n\in\mathbb{C}$ corresponding to a sequence $h_n:z\mapsto e^{t_n}(z-z_n)$ of such Möbius transformations are not bounded, then there exist arbitrarily small balls $B(x,\varepsilon)\subset\mathbb{P}^1$ with $\overline{\lim}\|(h_n^{-1})^{\ast}\U|_{B(x,\varepsilon)}\|_2=\|\U\|_2$. For the proof we observe that a subsequence of $(t_n,z_n)_{n\in\mathbb{N}}$ belongs to one of the following cases:
\begin{enumerate}
\item If $\mathbb{R} \ni t_n \to \infty$ and if $z_n$ converges in $\mathbb{C}$, then the limit of the pre-images $h_n^{-1}[B(\infty,\epsilon)]$ contains the complement of the one-point-set $\{\lim z_n\}$.
\item If $t_n$ is bounded from below and $|z_n|\to\infty$, then the limit of the pre-images $h_n^{-1}[B(\infty,\epsilon)]$ contains the complement of the one-point-set $\{\infty\}$
\item If $t_n$ is bounded from above and $e^{t_n}|z_n|\to\infty$, then the limit of the pre-images $h_n^{-1}[B(\infty,\epsilon)]$ contains the complement of the one-point-set $\{\infty\}$.
\item If $t_n\to-\infty$ and $e^{t_n}z_n$ converges in $\mathbb{C}$, then the limit of the pre-images $h_n^{-1}[B(\lim e^{t_n}z_n,\epsilon)]$ contains the complement of the one-point-set $\{\infty\}$.
\end{enumerate}  
Hence if $\|\U\|_2>S_{\frac{4}{3}}^{-1}$, then the maximizing sequence of Möbius transformations in $AN$ is bounded and therefore has a convergent subsequence.
In this case the continuity of $h \mapsto r_{\min}((h^{-1})^{\ast}\U)$ implies the claim.
On the other hand, if $\|\U\|_2\leq S_{\frac{4}{3}}^{-1}$, then $r_{\min}((h^{-1})^{\ast}\U) = \pi$ independent of the choice of $h$, and therefore the claim also follows with an arbitrary choice of $g$.  

Let $\moebius$ be a Möbius transformation maximizing $r_{\min}((\moebius^{-1})^{\ast}\U)$.
If the maximum is smaller than $\frac{\pi}{4}$, then  $r_{(\moebius^{-1})^{\ast}\U}$ attains the maximum at two different points with distance larger than $\frac{\pi}{4}$.
Indeed, were there only one maximum, or were the pairwise distance between all such maxima less that $\tfrac{\pi}{4}$, then all such balls are wholly contained in the hemisphere obtained by enlarging the radius of one of the balls to $\tfrac{\pi}{2}$.
In this case, there exists a dilation $h$ of $\mathbb{P}^1$ which increases the radii of the balls with radius $r_{\min}((\moebius^{-1})^{\ast}\U)$ contained in this hemisphere.
The measure of a ball $B$ with respect to the measure $|(\moebius^{-1})^{\ast}\U|^2$ is the same as the measure of $h[B]$ with respect to the measure $|((h\moebius)^{-1})^{\ast}\U|^2$.
Therefore $g$ cannot maximize $r_{\min}((\moebius^{-1})^{\ast}\U)$.
	
In conclusion, for any bounded sequences of potentials $(\U_n)_{n\in\mathbb{N}}$, the corresponding sequence $(\moebius_n)_{n\in\mathbb{N}}$ of Möbius transformations which maximize the minimum $r_{\min}((\moebius_n^{-1})^{\ast}\U_n)$ has the desired properties.
In fact, either a subsequence of $r_{\min}((\moebius_n^{-1})^{\ast}\U_n)$ is bounded from below, and then the assumption of Lemma~\ref{bounded point measures} is satisfied, or else $r_{\min}((\moebius_n^{-1})^{\ast}\U_n)$ converges to zero, and then the corresponding maxima of the sequence of functions $r_{(\moebius_n^{-1})^{\ast}\U_n}$ has at least two different accumulation points on $\mathbb{P}^1$.
Hence a weak limit of the sequence of measures $|(\moebius_n^{-1})^{\ast}\U_n|^2 \dmu(z)$ contains two point masses greater or equal to $S_{\frac{4}{3}}^{-2}$.
\end{proof}

\proofstep{Culmination of Step A.} 
We now describe the alternating application of the Lemmas~\ref{lem:split-off-disks} and~\ref{existence Moebius}.
At first Lemma~\ref{lem:split-off-disks} is applied once to $(\U_n)_{n\in\mathbb{N}}$, conceived as sequence of potentials on $\mathbb{P}^1$ vanishing outside of $\mathbb{D}$, and afterwards for each $l\in L$ at most once to the sequence $(\moebius_{n,l}^{-1})^{\ast}(\unity_{\Set{D}_{n,l}}\U_n)$.
These applications of Lemma~\ref{lem:split-off-disks} add $M$ elements to $L$, denoted by $l'$, which are in the first application maximal elements with respect to $<$, and in the later applications less than $l$.
The new sequences of disks have in the first application $(\Set{D}_{n,l'})_{n\in\mathbb{N}}=(B(z_m,r_{n,m}))_{n\in\mathbb{N}}$ centres inside $\overline{B(0,\frac{1}{2})}$ and are contained in the latter applications $(\Set{D}_{n,l'})_{n\in\mathbb{N}}=(\moebius_{n,l}[B(z_m,r_{n,m})])_{n\in\mathbb{N}}$ in $(\Set{D}_{n,l})_{n\in\mathbb{N}}$, where $m$ runs in both cases through $m=1,\ldots,M$.
For every new element $l'\in L$ an application of Lemma~\ref{existence Moebius} to the sequence $(\unity_{\Set{D}_{n,l'}}\,\U_n)_{n\in\mathbb{N}}$ firstly yields the corresponding sequence of Möbius transformations $g_{n,l'}$ and secondly defines a subsequence such that $|(\moebius_{n,l'}^{-1})^{\ast}(\unity_{\Set{D}_{n,l'}}\U_n)|^2 \dmu(z)$ converges.
Afterwards we continue with the alternating application of Lemma~\ref{lem:split-off-disks} followed by Lemma~\ref{existence Moebius} to all those new sequences $(\moebius_{n,l'}^{-1})^{\ast}(\unity_{\Set{D}_{n,l'}}\U_n)$ whose weak limits of the measures $|(\moebius_{n,l'}^{-1})^{\ast}(\unity_{\Set{D}_{n,l'}}\U_n)|^2$ still contain point masses greater or equal to $S_{\frac{4}{3}}^{-2}$.
Any application of Lemma~\ref{lem:split-off-disks} which follows an application of Lemma~\ref{existence Moebius} constructs inside the sequence of disks $(\Set{D}_{n,l})_{n\in\mathbb{N}}$ at least two new sequences of disks $(\Set{D}_{n,l'})_{n\in\mathbb{N}}$ which all obey
\begin{equation}
\label{eq:iteration-terminator}
\lim_{n\to\infty}\|\U_n\|_{\banach{2}(\Set{D}_{n,l'})}^2\le\lim_{n\to\infty} \|\U_n\|_{\banach{2}(\Set{D}_{n,l})}^2-S_{\frac{4}{3}}^{-2}.
\end{equation}
The iterated construction of disks inside disks terminates when for every new sequence of disks $(\Set{D}_{n,l'})_{n\in\mathbb{N}}$, the weak limit of $(|(\moebius_{n,l'}^{-1})^{\ast}(\unity_{\Set{D}_{n,l'}}\U_n)|^2\dmu(z))_{n\in\mathbb{N}}$ does not contain a point mass greater than $S_{\frac{4}{3}}^{-2}$.
Due to~\eqref{eq:iteration-terminator} this happens after finitely many alternating applications of the Lemmas~\ref{lem:split-off-disks} and~\ref{existence Moebius}.

We explained above how the sequences of disks $\{(\Set{D}_{n,l})_{n\in\mathbb{N}}\mid l\in L\}$ constructed in this way decompose $\mathbb{D}$ into sequences of disjoint domains~\eqref{eq:decomposition domains} and potentials $\U_{\reg,n}$ and $\U_{\sing,n,l}$~\eqref{eq:decomposition}.
By construction, for each $l\in L$ the weak limit of the sequence of measures $(|(\moebius_{n,l}^{-1})^{\ast}\U_{\sing,n,l}|^2)_{n\in\mathbb{N}}$ contains no point masses greater than $S_{\frac{4}{3}}^{-2}$.
This shows that the resulting decomposition has property \emph{Decomposition~(i)}.
The other properties \emph{Decomposition~(ii)--(iv)} follow from the construction and the two properties in Lemma~\ref{lem:split-off-disks}. %

Finally we show the property \emph{Decomposition (v)}.
We know that the volume of the discs $g_{n,l}^{-1}[\Set{D}_{n,l}]$ have lower bound.
If this were not the case, the weak limits of measures $(|(\moebius_{n,l}^{-1})^{\ast}\U_{\sing,n,l}|^2)_{n\in\mathbb{N}}$ would converge to a measure supported on a single point of mass larger than $S_{\frac{4}{3}}^{-2}$, which would contradict Lemma~\ref{existence Moebius}.
By the compactness of $\mathbb{P}^1$, after passing to a subsequence of $(g_{n,l}^{-1}[\Set{D}_{n,l}])_{n\in\mathbb{N}}$ the centres converge.
Hence there is a neighborhood $\SO_l$ of the limit of the centres which is contained in the discs for large enough $n$.%

\proofstep{Step B:} \textit{Limit of the local resolvents near a singular point whose blown-up holomorphic structure has a trivial kernel.} 
In this step we prove that for all $l\in L$ the sequences of potentials $(\U_{\sing,n,l})_{n\in\mathbb{N}}$ do not contribute to the limit of $(\Op{I}_{\mathbb{D}}(\U_n))_{n\in\mathbb{N}}$ if the weak limit of sequence of the blown up holomorphic structures with potentials $((\moebius_{n,l}^{-1})^{\ast}\U_{\sing,n,l})_{n\in\mathbb{N}}$ is a holomorphic structure on $\mathbb{P}^1$ with trivial kernel.
In fact, if for maximal $l\in L$ the assumption of the following lemma is fulfilled and additionally $\unity-\U_{\reg,n}\Op{I}_{\mathbb{C}}(0)$ converges to an invertible operator, then a subsequence of $\unity-\U_{\reg,n}\Op{I}_{\mathbb{C}}(\U_{\sing,n,l})$ is invertible and we obtain
\begin{align*}
\lim_{n\to\infty}\Op{I}_{\mathbb{C}}(\U_{\reg,n}+\U_{\sing,n,l})&=\lim_{n\to\infty}\Op{I}_{\mathbb{C}}(0)\left(\unity-\U_{\reg,n}\bigl(\Op{I}_{\mathbb{C}}(\U_{\sing,n,l})\bigr)\right)^{-1}\\&=\lim_{n\to\infty}\Op{I}_{\mathbb{C}}(0)\left(\unity-\U_{\reg,n}\Op{I}_{\mathbb{C}}(0)\right)^{-1}=\lim_{n\to\infty}\Op{I}_{\mathbb{C}}(\U_{\reg,n}).
\end{align*}
Analogously, if for non-maximal $l\in L$ and the index $l''$ in \emph{Decomposition (iv)} the assumption of the lemma is fulfilled and in addition $\unity-(\moebius_{n,l''}^{-1})^{\ast}  \U_{\sing,n,l''}\Op{I}_{\mathbb{C}}(0)$ converges to an invertible operator, then a subsequence of
\[
\unity-(\moebius_{n,l''}^{-1})^{\ast}  \U_{\sing,n,l''} \Op{I}_{\mathbb{C}}((\moebius_{n,l''}^{-1})^{\ast} \U_{\sing,n,l})
\]
is invertible and we obtain
\[
\lim_{n\to\infty}\Op{I}_{\mathbb{C}}\bigl((\moebius_{n,l''}^{-1})^{\ast}(\U_{\sing,n,l''}+\U_{\sing,n,l})\bigr)=\lim_{n\to\infty}\Op{I}_{\mathbb{C}}\bigl((\moebius_{n,l''}^{-1})^{\ast}\U_{\sing,n,l''}\bigr).
\]
In the subsequent step we will show how to meet the additional assumptions.
\begin{lemma}\label{trivial kernel}
If for some $l\in L$ the sequence of potentials $((\moebius_{n,l}^{-1})^{\ast}\U_{\sing,n,l})_{n\in\mathbb{N}}$ converges weakly to some holomorphic structure on $\mathbb{P}^1$ which has a trivial kernel, then the following sequence converges to zero with respect to the norm topology of operators on $\banach{\frac{4}{3}}(\mathbb{C},\mathbb{H})$:
\begin{align*}
& \text{if $l$ is maximal in $L$:} & \U_{\reg,n}\bigl(\Op{I}_{\mathbb{C}}(\U_{\sing,n,l}) & -\Op{I}_{\mathbb{C}}(0)\bigr) 
\\& \text{if $l$ is not maximal in $L$:} &(\moebius_{n,l''}^{-1})^{\ast}  \U_{\sing,n,l''} \bigl( \Op{I}_{\mathbb{C}}((\moebius_{n,l''}^{-1})^{\ast} \U_{\sing,n,l})  & - \Op{I}_{\mathbb{C}}(0) \bigr)
\;,
\end{align*}
with the index $l''$ introduced in \emph{Decomposition (iv)}.
\end{lemma}
\begin{proof}
We give the proof only for the case of maximal $l\in L$; the other case is proven by analogous arguments.
For $1<p<2$, every Möbius transformation $h$ induces an isometry
\begin{equation*}
\pi_p(h):\banach{p}(\mathbb{C},\mathbb{H})\to\banach{p}(\mathbb{C},\mathbb{H}),\;f\mapsto \Tilde{f} \quad\text{where}\quad \Tilde{f}(z)=f(h^{-1}z)\left|\frac{d(h^{-1}z)}{dz}\right|^{\frac{2}{p}}.
\end{equation*}
A direct calculation shows that the inverse $\Op{I}_{\mathbb{C}}(0)$ of the free holomorphic structure, considered as an operator from $\banach{p}(\mathbb{C},\mathbb{H})$ into $\banach{\frac{2p}{2-p}}(\mathbb{C},\mathbb{H})$ with $1<p<2$, is invariant under the semidirect product $AN$ of the scaling transformations ($z\mapsto\exp(t)z$ with $t\in\mathbb{R}$) and the translations ($z\mapsto z+z_0$ with $z_0\in\mathbb{C}$).
Therefore $\Op{I}_\mathbb{C}(0)$ is invariant under
$\moebius_{n,l}$:
\[
\pi_{\frac{2p}{2-p}}(\moebius_{n,l})\Op{I}_{\mathbb{C}}(0)\pi_p(\moebius_{n,l}^{-1})=\Op{I}_{\mathbb{C}}(0).
\]

Due to the assumption on the weak limit and Theorem~\ref{inverse holomorphic structure}, the sequence of operators $(\unity-\U_{\sing,l}\Op{I}_\mathbb{C}(0))^{-1}$ on $\banach{\frac{4}{3}}(\mathbb{C},\qat)$ is bounded.
Next we claim that the following sequence of operators on $\banach{\frac{4}{3}}(\mathbb{C},\mathbb{H})$ converges to zero:
\begin{align}\notag
\U_{\reg,n}\Op{I}_{\mathbb{C}}(\U_{\sing,n,l})-\U_{\reg,n}\Op{I}_{\mathbb{C}}(0)
&=\U_{\reg,n}\Op{I}_{\mathbb{C}}(0)\U_{\sing,n,l}\Op{I}_{\mathbb{C}}(\U_{\sing,n,l})\\
\label{eq:potential-Ic}
&=\U_{\reg,n}\Op{I}_{\mathbb{C}}(0)\U_{\sing,n,l}\Op{I}_{\mathbb{C}}(0)
\left(\unity-\U_{\sing,n,l}\Op{I}_{\mathbb{C}}(0)\right)^{-1}.
\end{align}
In fact, due to condition \emph{Decomposition~(iii)}, the $\banach{2}$\=/norms of the restrictions of the potentials $\U_{\reg,n}$ to small balls $B(x_l,\varepsilon)$ around the singular point $x_l$ converge to zero for $\varepsilon\to 0$.
Therefore it suffices to show the analogous claim where the potentials $\U_{\reg,n}$ is replaced in\eqref{eq:potential-Ic} by $\U_{\reg,n}\cdot \unity_{\mathbb{P}^1 \setminus B(x_l,\varepsilon)}$ with sufficiently small $\varepsilon>0$.
In doing so, the integral kernel of $\Op{I}_{\mathbb{C}}(0)$ in between the two potentials $\U_{\reg,n}$ and $\U_{\sing,n,l}$ in\eqref{eq:potential-Ic} is evaluated on the cartesian product of the supports of $\U_{\sing,n,l}$ with $\mathbb{P}^1\setminus B(x_l,\epsilon)$.
For large $n$ the supports of $\U_{\sing,n,l}$ are contained in arbitrary small balls $B(x_l,\varepsilon')\subset B(x_l,\varepsilon)$.
For $0<\varepsilon'<\varepsilon$, the integral kernel of $\Op{I}_{\mathbb{C}}(0)$ from Chapter~\ref{chapter:local} is the limit
\[
\frac{\dmu(z)}{\pi(z-w)}
=-\frac{1}{\pi}\sum_{m=0}^\infty\frac{(w-x_l)^m\dmu(z)}{(z-x_l)^{m+1}}
\quad\text{on}\quad
(w,z)\in B(x_l,\epsilon')\times\big(\mathbb{P}^1\setminus B(x_l,\epsilon)\big).
\]
This series defines for $1<p<\frac{4}{3}$ a series of finite rank operators from $\banach{p}(B(x_l,\varepsilon'),\mathbb{H})$ to $\banach{4}(\mathbb{C}\setminus B(x_l,\varepsilon),\mathbb{H})$.
Since $\U_{\sing,n,l}$ is a bounded operator from $\banach{\frac{p}{2+p}}(B(x_l,\varepsilon'),\mathbb{H})$ to $\banach{p}(B(x_l,\varepsilon'),\mathbb{H})$, we conceive the operator $\Op{I}_{\mathbb{C}}(0)$ occurring in~\eqref{eq:potential-Ic} on the left hand side of the bounded operators $\left(\unity-\U_{\sing,n,l}\Op{I}_{\mathbb{C}}(0)\right)^{-1}$ as an operator from $\banach{\frac{4}{3}}(\mathbb{C},\mathbb{H})$ into $\banach{\frac{p}{p+2}}(B(x_l,\varepsilon'),\mathbb{H})$.
This operator is the composition of the bounded operator $\Op{I}_{\mathbb{C}}(0)$ from $\banach{\frac{4}{3}}(\mathbb{C},\mathbb{H})$ to $\sobolev{1,\frac{4}{3}}(B(x_l,\varepsilon'),\qat)$ (see Chapter~\ref{chapter:local}) with the embedding $\sobolev{1,\frac{4}{3}}(B(x_l,\varepsilon'),\qat)\hookrightarrow\banach{\frac{p}{p+2}}(B(x_l,\varepsilon'),\mathbb{H})$ which is compact by the Rellich-Kondrachov theorem~\cite[Theorem~6.3]{Ad}.
Consequently this composition is a limit of finite rank operators.
Since now the operators $\Op{I}_{\mathbb{C}}(0)$ on both sides of $\U_{\sing,n,l}$ in~\eqref{eq:potential-Ic} are uniform limits of finite rank operators, we may replace $\U_{\sing,n,l}$ in the limit $n\to\infty$ by a weak limit.
This weak limit is zero because $\U_{\sing,n,l}$ is bounded and the support shrinks for $n\to\infty$ to $\{x_l\}$.
This proves the Lemma.
\end{proof}
\proofstep{Step C:} \textit{Limit of the local resolvents near a singular point whose blown-up holomorphic structure has a non-trivial kernel.}
In this case we add a sequence of effective divisors $D_n$ with support contained in the supports of $\U_{\sing,n,l}$, such that the corresponding transformed sequences of holomorphic structures with potentials $(\moebius_{n,l}^{-1})^{\ast}\U_{\sing,n,l}$ on $\mathbb{P}^1$ have trivial kernels.  
\begin{lemma}\label{transformed holomorphic structures}
For any holomorphic $\qat$\=/line bundle $E_\qat$ on $\mathbb{P}^1$ whose sheaf $\Q{E,\V}$ has non-trivial sections, there exists a unique $d\in \mathbb{N}_+$ such that
\begin{align}\label{eq:trafo-holo-str:H0d}
\dim H^0\big(\mathbb{P}^1\!,\Q{E(-d\infty),\V(-d\infty)}\big)&\!=\!0,&
\dim H^0\big(\mathbb{P}^1\!,\Q{E((1-d)\infty),\V((1-d)\infty)}\big)&\!=\!1.
\end{align}
Moreover, for every non-empty open subset $\SO\subset\mathbb{C}$, there exists an effective divisor $D$ of degree $d-\deg(E)-1$ whose support is contained in $\SO$ and such that
\begin{align*}
\dim H^0\big(\mathbb{P}^1,\Q{E(D+(m-d)\infty),\V(D+(m-d)\infty)}\big)&=m&\text{for all }m\in\mathbb{N}_0.
\end{align*}
\end{lemma}
\begin{proof}
Due to the Plücker formula (Corollary~\ref{pluecker formula}), for large $d\in \mathbb{N}_+$ we have $\dim H^0\left(\mathbb{P}^1,\Q{E(-d\infty),\V(-d\infty)}\right)=0$.
This dimension can increase by at most $1$ when $d$ is diminished by $1$. Hence the minimal $d$ with $\dim H^0\left(\mathbb{P}^1,\Q{E(-d\infty),\V(-d\infty)}\right)=0$ is the only $d$ that satisfies~\eqref{eq:trafo-holo-str:H0d}.
The Riemann-Roch Theorem~\ref{riemann roch} yields
\begin{align*}
\dim H^0\big(\mathbb{P}^1\!,\Q{E(-d\infty),\V(-d\infty)}\big)&=
\deg(E)+1-d+ \dim H^1\big(\mathbb{P}^1\!,\Q{E(-d\infty),\V(-d\infty)}\big)\\
&\geq\deg(E)+1-d.
\end{align*}
By definition of $d$ this implies $\deg(E)\leq d-1$. Moreover $\deg(E)=d-1$ holds if and only if $\dim H^1(\mathbb{P}^1,\Q{E(-d\infty),\V(-d\infty)})=0$.
By Serre duality~\ref{Serre duality} this is equivalent to
\begin{align*}
\dim H^1\big(\mathbb{P}^1,\Q{E((m-d)\infty),\V((m-d)\infty)}\big)=0
\quad\text{for all }m\in\mathbb{N}_0.
\end{align*}
Finally, due to the Riemann-Roch theorem, this in turn is equivalent to
\begin{align*}
\dim H^0\big(\mathbb{P}^1,\Q{E((m-d)\infty),\V((m-d)\infty)}\big)=m
\quad\text{for all }m\in\mathbb{N}_0.
\end{align*}
So for $\deg(E)=d-1$ the statement of the lemma holds with $D=0$.

We now consider the cases $\deg(E)<d-1$.
We claim that in this case there exists an element $z\in \SO$ such that $\dim H^0(\mathbb{P}^1,\Q{E(z-d\infty),\V(z-d\infty)})=0$.
Let us assume to the contrary that we have $\dim H^0(\mathbb{P}^1,\Q{E(z-d\infty),\V(z-d\infty)})=1$ for all $z\in\SO$.
Consequently $H^0(\mathbb{P}^1,\Q{E(w_1+\ldots+w_m-d\infty),\V(w_1+\cdots+w_m-d\infty)})$ contains $m$ linear independent elements with poles at $w_1,\ldots,w_m$, respectively.
Serre duality and the Plücker formula imply $\dim H^1(\mathbb{P}^1,\Q{E(z_1+\cdots+z_m-d\infty),\V(z_1+\cdots+z_m-d\infty)})=0$ for large $m$.
Consequently, due to the Riemann-Roch theorem, we obtain
\[
m\leq \dim H^0\left(\mathbb{P}^1,\Q{E(z_1+\cdots+z_m-d\infty),\V(z_1+\cdots+z_m-d\infty)}\right)=1+\deg(E)+m-d,
\]
which contradicts $\deg(E)<d-1$.
This proves the claim.

We now iteratively replace $E$ by $E(z)$, where $z\in \SO$ is chosen according to the claim.
This replacement increases $\deg(E)$ by $1$ without changing $d$.
Therefore after $d-1-\deg(E)$ such replacements, we arrive at $E(D)$ with some effective divisor $D$, which satisfies $\deg(E(D)) = d-1$ and therefore has the desired properties.
\end{proof}
We shall use the foregoing lemma in the following manner: For any line bundle $E$ this lemma gives us a divisor so that the $\deg (E)+\deg(D)=d-1$.
We then add further points at infinity so $E(D-d\infty)$ has degree $-1$ and hence must be the spin bundle on $\mathbb{P}^1$.
Setting $n=0$ in the lemma then shows that the holomorphic structure of $\V(D-d\infty)$ has a trivial kernel.

\begin{proof}[Proof of Proposition~\ref{disc limits}]
We now combine the three steps we have prepared above to complete the proof of Proposition~\ref{disc limits}.

After performing Step~A we choose a minimal element $l\in L$.
We know by Lemma~\ref{bounded point measures} that the resolvents of $(g_{n,l}^{-1})^*(\U_{\sing,n,l})_{n\in\mathbb{N}}$ converge to the resolvent of the weak limit.
So if the weak limit has non-trivial kernel then Lemma~\ref{transformed holomorphic structures} defines a divisor with support in $\SO_l$.
We define $D_{n,l}$ as $g_{n,l}^{-1}$ applied to this divisor.
For non-minimal $l\in L$, we first change the potential $\U_{\sing,n,l}$ by the divisor $\sum_{l' < l}D_{n,l'}$ and then apply the same procedure.
Inductively for all $l\in L$, the weak limit the modified potentials $\U_{\sing,n,l}(\sum_{l' \leq l}D_{n,l'})$ gives a holomorphic structure with trivial kernel.

Finally the transformed holomorphic $\qat$\=/line bundle with the sequence of potentials $(\U_n(\sum_{l \in L} D_{n,l}))_{n\in\mathbb{N}}$ converges to the holomorphic $\qat$\=/line bundle whose potentials is the weak limit of $(\U_{\reg,n}(\sum_{l \in L} D_{n,l}))_{n\in\mathbb{N}}$.
If this limit has non-trivial holomorphic sections, then we apply Lemma~\ref{transformed holomorphic structures} with $\SO=\{z\in\mathbb{C}\mid 2<|z|<3\}$ so that the transformed holomorphic $\qat$\=/line bundles converge to a holomorphic $\qat$\=/line bundle on $\mathbb{P}^1$ without holomorphic sections.

We collect all these changes as a sequence $(P'_n)_{n\in\mathbb{N}}$ of complex polynomials whose roots combine the divisors from both the $\U_{\sing}$ and $\U_{\reg}$ modifications.
Since the roots of these polynomials are contained in compact subsets of $\mathbb{C}$, by passing to a subsequence we may assume that these sequences $(P'_n)_{n\in\mathbb{N}}$ converge on $\mathbb{C}$.
This concludes the proof of Proposition~\ref{disc limits}.
\end{proof}

In the following corollary we see that the same arguments apply to three other sequences of resolvents, to be used in Chapter~\ref{chapter:singular holomorphic sheaves}.
\begin{corollary}
\label{cor:other_resolvents_converge}
Let $(\U_n)_{n\in\mathbb{N}}$ be a sequence of potentials in $\banach{2}(\mathbb{D},\qat^-)$ fulfilling the conditions of Proposition~\ref{disc limits}.
Then there exists, besides the sequence of polynomials $(P'_n)_{n\in\mathbb{N}}$, another sequence of polynomials $(P''_n)_{n\in\mathbb{N}}$ which obeys condition~(i) in Proposition~\ref{disc limits} such that the following sequences of operators also converge in the space of bounded operators from $\banach{\frac{4}{3}}(\mathbb{D},\qat)$ to $\banach{q}(\mathbb{D},\qat)$:
\begin{align*}
\Big(\Op{I}_{\mathbb{D}}({P'_n}^{-1}\U_n\sd P'_n) \Big)_{n\in\mathbb{N}},\quad
\Big(\Op{I}_{\mathbb{D}}(P''_n\U_n\sd{P''_n}^{-1})\Big)_{n\in\mathbb{N}},\quad
\text{and}\quad
\Big(\Op{I}_{\mathbb{D}}({P''_n}^{-1}\U_n P''_n)\Big)_{n\in\mathbb{N}}.
\end{align*}
\end{corollary}
\begin{proof}
The proof is a three-fold modification of the arguments of Proposition~\ref{disc limits}.
To prove of the convergence of the first and second sequence of operators we apply the arguments of Step~A to the sequence of potentials $(\U\sd_n)_{n\in\mathbb{N}}$ instead of $(\U_n)_{n\in\mathbb{N}}$.
Since the decomposition only depends on the sequence of measures $|\U_n|^2=|\Bar{\U}_n|^2$ this application results in the same decomposition.
For a potential $\V\in\pot{E}$ acting on the spin bundle $E$ of $\mathbb{P}^1$, namely the bundle with $\deg(E)=-1$ and $E \cong KE^{-1}$, then Serre duality~\ref{Serre duality} and Riemann--Roch Theorem~\ref{riemann roch} implies $\dim H^0(\mathbb{P}^1,\Q{E,\V})=\dim H^0(\mathbb{P}^1,\Q{E,\V\sd})$.
Thus in Step~B the weak limit of the blown up holomorphic structures with potentials $((\moebius_{n,l}^{-1})^{\ast}\qk\Bar{\U}_{\sing,n,l}\qk)_{n\in\mathbb{N}}$ has a trivial kernel if and only if the same is true for the corresponding weak limit of the blown up holomorphic structures with potentials $((\moebius_{n,l}^{-1})^{\ast}\U_{\sing,n,l})_{n\in\mathbb{N}}$.
Furthermore, from Lemma~\ref{eq:trafo-holo-str:H0d} we know there is a divisor $D'$ with degree $d$ such that
\begin{align*}
\dim H^0\left(\mathbb{P}^1,\Q{E(D'+(m-d)\infty),\V(D'+(n-d)\infty)}\right)= \max\{0, m\}
\quad\text{for all } m\in\mathbb{Z}.
\end{align*}
Applying Serre duality and Riemann-Roch then
\begin{align*}
\dim H^0\left(\mathbb{P}^1,\Q{E(-D'+(m+d)\infty),\V(D'-(m+d)\infty)\sd}\right)&= m + \max\{0, -m\} \\
&= \max\{0, m\}
\end{align*}
for all $m \in\mathbb{Z}$.
With these modifications the arguments of the proof of Proposition~\ref{disc limits} show the convergence of the first sequence of operators.

In the proof of the convergence of the second sequence we apply the arguments of Steps~B and~C to the decomposition of the sequence of potentials $(\U\sd_n)_{n\in\mathbb{N}}$.
This yields a sequence of polynomials $(P''_n)_{n\in\mathbb{N}}$ and proves the convergence of the second sequence of operators.

The relationship between the third sequence of operators and the second is the same as the relationship between the first sequence and the sequence of Proposition~\ref{disc limits}.
The convergence of the third sequence follows.
\end{proof}

\chapter{Global Limits of Sequences of Holomorphic Structures}
\label{chapter:global limits}
In this chapter we shall globalize the local understanding, achieved in the foregoing chapter, of limits of sequences of holomorphic structures in the case of massive concentration of energy.
In contrast to the local approach, the global point of view has to take into account the degrees of freedom of non-isomorphic holomorphic $\mathbb{C}$\=/line bundles.
This data consists of the degree of the line bundle, which is discrete, and the finite dimensional Jacobian, which parametrizes the holomorphic structures.
Since for different underlying $\mathbb{C}$\=/line bundles both holomorphic structures and different holomorphic sections thereof belong to different spaces, the notion of limits has no obvious meaning.
At first sight, one might hope that this only concerns the discrete degree of freedom of the underlying $\mathbb{C}$\=/line bundle. This would be easy to fix: we bound the degree and pass to a subsequence of constant degree.
However, in the foregoing section the control of the limit needs in some cases to change the underlying holomorphic $\mathbb{C}$\=/line bundle in such a way that the original sequence of sections has additional roots of finite order nearby the points of massive concentration of energy and converge to a section without such additional roots.
Moreover, the exact location of the sequences of roots of the sequence of sections is determined by a suitable sequence of blow ups of the holomorphic structures nearby these points.
Therefore we need to define limits of sequences $(\xi_n)_{n\in\mathbb{N}}$ of holomorphic sections of sequences $\qat$\=/line bundles, whose underlying holomorphic $\mathbb{C}$\=/line bundles $(E_n)_{n\in\mathbb{N}}$ are not isomorphic.

For this purpose a subsequence of the spaces of global holomorphic sections $H^0(\X,\Q{E_n,\U_n})$ is embedded into a single Banach space.
This follows from the four Lemmas in this chapter:
First we select in Lemma~\ref{global resolvent} a suitable background holomorphic $\mathbb{C}$\=/line bundle $E$, and represent the members $(E_n)_{n\in\mathbb{N}}$ of the sequence of holomorphic line bundles as $E_n=E(D_n)$ with sequences of divisors $D_n$.
In Lemma~\ref{lem:embed sections} we cover $\X$ by finitely many pairwise disjoint open discs $\Set{U}_1,\ldots,\Set{U}_M$ and an additional open subset $\Set{U}_0$.
Let $\Set{U}$ denote the union $\Set{U}=\Set{U}_1\cup\ldots\cup\Set{U}_M$ of the discs.
Then we show that the restriction maps $H^0(\X,\Q{E(D),\U})\hookrightarrow H^0(\Set{U}\cap\Set{U}_0,\sob{1,p}{E})$ are embeddings, if the support of $D$ is disjoint from $\Bar{\Set{U}}_0$ and if $\|\U|_{\Set{\U}_0}\|_2$ is sufficiently small.
The idea behind these conditions has a geometric flavor: We decompose the difficult to control influence of the divisor and the potential into finitely many local situations on topologically trivial discs.
In addition we separate this influence from the global non-trivial topology.
A combination of Cauchy's integral formula with the results of the foregoing chapter allows us to recover the values on the discs of any $\xi\in H^0(\X,\Q{E(D),\U})$ from the restriction $\xi|_{\Set{U}\cap\Set{U}_0}$.
Similarly, due to the choice of $E$ and the bound on $\|\U|_{\Set{U}_0}\|$ this restriction also determines the values on the remainder set $\Set{U}_0$.
In other words, $\xi$ is globally determined by its restriction.

The two final Lemmas~\ref{covering} and~\ref{singular holomorphic convergence} show that under the assumption that $\deg(D_n)$ and $\|\U_n\|_2$ are bounded, any sequence of pairs $(D_n,\U_n)_{n\in\mathbb{N}}$ with $\U_n\in\pot{\X,E(D_n)}$ has a subsequence to which the forgoing lemmas apply.
This allows us to achieve the main aim of this section, namely to prove the following quaternionic function theoretical version of Montel's theorem: A sequence $\xi_n\in H^0(\X,\Q{E_n,\U_n})$ has a subsequence which converges to a global meromorphic section of $\Q{E,\U}$ if $\deg(E_n)$, $\|\U_n\|_2$ and $\|\xi_n\|$ are bounded.
Here $E$ is an accumulation point of $(E_n)_{n\in\mathbb{N}}$ and $\U$ a weak accumulation point of $(\U_n)_{n\in\mathbb{N}}$.
The statement of this main result, Theorem~\ref{global limits}, is postponed until the next chapter, where we will be in a position to specify the poles of the limit $\xi=\lim_{n\to\infty}\xi_n$.

In the following lemma we establish for a suitable fixed $\mathbb{C}$\=/line bundle $E$ the weak continuity of the map from a sufficiently small ball of $\pot{E}$ to the inverse operators of the corresponding holomorphic structures.
This is a necessary technical result required by the subsequent lemmas.
\begin{lemma}\label{global resolvent}
\index{Resolvent}
On a compact Riemann surface $\X$ there exists a holomorphic $\mathbb{C}$\=/line bundle $E$ on a compact Riemann surface $\X$ such that $E$ and $KE^{-1}$ have no holomorphic sections.
For any such line bundle $E$ and $1<p<2$ there exists a $C_p>0$, such that $\delbar{E}-\V:H^0(\X,\sob{1,p}{E})\to H^0(\X,\forms{0,1}{}\ban{p}{E})$ is invertible for any $\V\in\pot{E}$ with $\|\V\|_2\le C_p$.
Furthermore, for any $1<q<\frac{2p}{2-p}$ the map $\V\mapsto(\delbar{E}-\V)^{-1}$ is weakly continuous from $\{\V\in\pot{E}\mid\|\V\|_2\leq C_p\}$ into the compact operators $H^0(\X,\forms{0,1}{}\ban{p}{E})\to H^0(\X,\sob{1,p}{E})\hookrightarrow H^0(\X,\ban{q}{E})$.
\end{lemma}
\begin{proof}
Let us first show that on any compact Riemann surface $\X$ there exists a divisor $D$ such that the corresponding holomorphic $\mathbb{C}$\=/line bundle $E$ obeys the first condition.
By Serre duality and the Riemann-Roch theorem $\deg(D)$ equals $\genus-1$, and conversely any divisor $D$ of degree $\genus-1$ with $H^0(\X,\Sh{O}_D)=0$ has this property.
For $\genus=0$ the unique holomorphic $\mathbb{C}$\=/line bundle of degree $-1$ has this property.
For $\genus\ge1$ such $D$ are mapped by the Abel Jacobi map onto the complement of the $\theta$ divisor in the Jacobian (compare~\cite[Chapter~IV.3]{FK}) translated by some constant.
So the set of equivalence classes of such holomorphic $\mathbb{C}$\=/line bundles $E$ is an open and dense subset of the component of the Picard group that contains all holomorphic line bundles of degree $\genus-1$.

In the proof of Serre duality, Theorem~\ref{Serre duality}, we proved that $\delbar{E}-\V$ defines for any $\V\in\pot{E}$ a Fredholm operator from $H^0(\X,\sob{1,p}{E})$ to $H^0(\forms{0,1}{}\ban{p}{E})$ with kernel $H^0(\X,\Q{E,\V})$ and cokernel isomorphic to $H^1(\X,\Q{E,\V})$.
Hence, in the present situation with $\V=0$, the operator $\delbar{E}$ is an isomorphism from $H^0(\X,\sob{1,p}{E})$ onto $H^0(\X,\forms{1,0}{}\ban{p}{E})$ and has a bounded inverse.
Let $C_p>0$ be smaller than the norm of this inverse operator.
As in Lemma~\ref{lem:sobolev regularity} the multiplication with a potential $\V \in \pot{E}$ is a bounded operator from $H^0(\X,\sob{1,p}{E})$ into $H^0(\X,\forms{0,1}{}\ban{p}{E})$.
For all potentials $\V$ with $\|\V\|_2\leq C_p$
the Neumann series of $(\delbar{E} - \V)^{-1}$ converges in the space of operators from $H^0(\X,\forms{0,1}{}\ban{p}{E})$ to $H^0(\X,\sob{1,p}{E})$.

Finally we will use a similar argument as in Lemma~\ref{weakly continuous resolvent} to conclude that the inverse of the holomorphic structure depends weakly continuously on the potential.
Due to Rellich-Kondrachov theorem~\ref{lem:sobolev regularity} the embedding $H^0(\X,\sob{1,p}{E})$ into $H^0(\X,\ban{q}{E})$ is compact for all $q<\frac{2p}{2-p}$.
Due to~\cite[Theorem~II.5.11]{LT}, the Banach spaces $\banach{q}(\Omega,\mathbb{H})$ have a Schauder basis.
Consequently they have the \emph{approximation property}, namely that all compact operators into $\banach{q}(\Omega,\mathbb{H})$ are norm-limits of finite rank operators (compare~\cite[Section~I.1.a]{LT}).
Hence all terms in the Neumann series are norm-limits of weakly continuous functions from the set of all potentials $\V$ with $\|\V\|_2\leq C_p$ into the compact operators from $H^0(\X,\forms{0,1}{}\ban{p}{E})$ into $H^0(\X,\ban{q}{E})$.
Since this set is weakly compact, the uniform limit of weakly continuous functions is again a weakly continuous function~\cite[Theorem~IV.8]{RS1}.
\end{proof}
\begin{lemma}\label{lem:embed sections}
For a compact Riemann surface $\X$ fix a line bundle $E$ satisfying the conditions of the above Lemma~\ref{global resolvent}.
Choose a cover of $\X$ by a finite number of disjoint open discs $\SU_1,\ldots,\SU_M$ and an additional open subset $\SU_0 \subset \X$.
Let $\SU = \bigcup_{m=1}^M \SU_m$, so that $\SU \cap \SU_0$ is the region of overlap.

Suppose that we have a set of holomorphic $\qat$\=/line bundles where
\begin{enumerate}
\item 
each underlying $\mathbb{C}$\=/line bundle has the form $E(D)$ such that the support of the divisor $D$ lies in the complement of $\overline{\SU_0}$.
\item 
each potential $\U\in\pot{E(D)}^{-}$ restricted to $\SU_0$ has $\banach{2}$\=/norm not larger than the constant $C_p$ from Lemma~\ref{global resolvent} for some $1<p<2$.
\end{enumerate}
The restriction map from $H^0(\X,\Q{E(D),\U})$ to $H^0(\SU\cap\SU_0,\sob{1,p}{E})$ is a well-defined embedding.
For any $C>0$ consider all restriction maps arising from a divisor $D$ satisfying (i) and $|\deg D|\le C$ and a potential $\U$ satisfying (ii) and $\|\U\|_2\le C$.
Let $\Spa{C}_C$ be the union of the images of all such restriction maps.
Then any bounded subset of $\Spa{C}_C$ is relatively compact in $H^0(\SU\cap\SU_0,\sob{1,p}{E})$.
\end{lemma}
\begin{proof}
In order to simplify the geometry we assume that $\SU_0$ is the union of $(\X\setminus\SU)\cup\SV$ with the disjoint union $\SV=\SV_1\cup\ldots\cup\SV_M$ of open annuli $\SV_m\subset\SU_m$ whose outer boundary coincides with $\boundary\SU_m$ and whose inner boundary is a component of $\boundary\SU_0$.
Any open neighborhood of $\X\setminus\SU$ contains an open subset $\SU_0$ with this property, and the statement for $\Set{U}_0$ implies the statement for any open set $\supset\SU_0$.

Fix $D$ and $\U$ obeying (i) and (ii) respectively.
If we restrict the line bundles $E(D)$ and $E$ to $\SV$ then they are isomorphic, because the support of $D$ is contained in $\X\setminus\overline{\SU_0} \supset\SU\setminus\overline{\SV}$.
Hence the restriction map $H^0(\X,\Q{E(D),\U})\to H^0(\SU\cap\SU_0,\sob{1,p}{E})$ is well-defined.
We now prove that this map is an embedding.

For any $m=1,\ldots,M$ we choose two disjoint open sub-annuli $\Set{A}_m$ and $\Set{B}_m$ of $\SV_m$ such that the outside boundary of $\Set{A}_m$ is again $\boundary\SU_m\subset\boundary\SV_m$ and the inside boundary of $\Set{B}_m$ is the inside boundary of $\SV_m$.
Let $\Set{A}$ and $\Set{B}$ denote the following disjoint unions of annuli:
$\Set{A}:= \bigcup_{m=1}^M \Set{A}_m$ and $\Set{B}:=\bigcup_{m=1}^M \Set{B}_m$.

Choose a smooth real cut-off function $f$ on $\X$, which is equal to $1$ on $\Set{U}_0\setminus\Set{B}\supset\X\setminus\Set{U}$, equal to zero on $\SU\setminus\Set{A}$, and whose derivative has compact support in $\Set{B}$.
In particular, $f$ has compact support in $\Set{U}_0$.
For any $\xi\in H^0(\X,\Q{E(D),\U})$, we consider $(\barpartial f)\xi=(\delbar{E}-\U)(f\xi)$ as a global section of $\forms{0,1}{}\ban{p}{E}$ because it is zero outside of $\Set{B}\subset\SU_0$.
By Lemma~\ref{global resolvent}, we can apply the inverse of $\delbar{E} - \U$ to $(\barpartial f)\xi$ to get $f\xi$.
In particular, $\xi|_{\SU_0\setminus\Set{B}}=f\xi|_{\SU_0\setminus\Set{B}}\in H^0(\SU_0\setminus\Set{B},\sob{1,p}{E})$ is determined by $\xi|_{\Set{B}}\in H^0(\Set{B},\ban{p}{E})$.

On the other hand, for every $m=1,\dots,M$, we choose a real cut-off function $f_m$ that is $1$ on $\SU_m\setminus\Set{A}_m$ and has compact support on $\SU_m$.
Since $\SU_m$ is biholomorphic to $\mathbb{D}$ we may trivialize $E(D)$ by a polynomials.
There exists unique polynomials $q_{m}, r_{m}$ whose zeroes are contained in $\mathbb{D}$ such that $\frac{q_{m}}{r_{m}} \xi$ restricted to $\SU_m$ is a section of the spin bundle on $\mathbb{P}^1$.

We apply Proposition~\ref{disc limits} to the constant sequence of holomorphic sections $(\frac{q_{m}}{r_{m}} \xi)_{n\in\mathbb{N}}$ with potential $\frac{q_{m}}{r_{m}}\U\frac{r_{m}}{q_{m}}$ on the disc $\SU_m$.
Because in this situation there is no concentration of the measure associated to $\U$, the effective divisor $D_{m}$ obtained from the proposition is trivial.
The proposition tells us that the resolvent of the holomorphic structure is an operator from $H^0(\Set{A}_m,\ban{4/3}{E})$ into $H^0(\SU_m\setminus\Set{A}_m,\ban{q}{E})$ with $1<q<4$.
If we combine this for $m=1,\ldots,M$, due to the quaternionic version of Cauchy's integral formula~\ref{cauchy formula}~(ii) $\xi$ as a section of $H^0(\Set{A},\ban{4/3}{E})$ determines $\xi$ in $H^0(\Set{B},\ban{q}{E})$.
Moreover the $H^0(\X,\sob{1,\frac{4}{3}}{})$\=/norm is bounded by the $H^0(\X,\ban{q}{})$\=/norm.
This shows that the restriction map is an embedding.

To finish the proof, we must show that for given $C>0$ any bounded sequence $\xi_n\in\Spa{C}_C$ has a convergent subsequence.
We denote by $D_n$ the divisors and $\U_n\in\pot{E(D_n)}^{-}$ the potentials obeying properties (i) and (ii) such that the sections $\xi_n$ lie in $H^0(\X, \Q{E(D_n),\U_n})$ respectively.
Since the degrees of the divisors are bounded, there exists a subsequence of constant degree.
By passing to a subsequence we may assume that the divisors have constant degree and converge.
For each $k$ and $m$ between $1$ and $M$, let $\frac{q_{n,m}}{r_{n,m}}$ be a rational function that trivializes the line bundle $E(D_n)$ on the disc $\SU_m$.
For any $m=1,\dots,M$ we apply Proposition~\ref{disc limits} to the sections $\frac{q_{n,m}}{r_{n,m}}\xi_n|_{\SU_m}$ which are holomorphic with respect to the holomorphic structure $\frac{q_{n,m}}{r_{n,m}}\U_n\frac{r_{n,m}}{q_{n,m}}$.
We obtain from the proposition the polynomials $p_{n,m}$ with corresponding divisor $D'_{n,m}$.
The support of this divisor lies outside $\SU_0$ because it is contained in the set of points where the square of the $\banach{2}$\=/norm of the potential concentrates with energy at least $\pi$, and by property (ii) and the previous lemma, no such points lie in $\SU_0$.
In the case that the limit of $(D_{n}+ \sum_{m}D'_{n,m})_{n\in\mathbb{N}}$ has a point on the boundary of $\SU_0$, we may shrink $\SU_0$ slightly to avoid this.
After passing to a subsequence we may suppose that the polynomials $p_{n,m}, q_{n,m}, r_{n,m}$ and their inverses are uniformly convergent on $\SU\cap\SU_0$.

For each $m=1,\dots,M$, by the same argument with $\Set{A}_m$ and $\Set{B}_m$ above, we know that the $H^0(\SU_m\cap\SU_0,\sob{1,\frac{4}{3}}{})$\=/norms of the sections $\frac{p_{n,m}q_{n,m}}{r_{n,m}}\xi_n$ are bounded uniformly in terms of their $H^0(\SU_m\cap\SU_0,\ban{q}{})$\=/norms.
Because the polynomials and their inverses are uniformly convergent, the $H^0(\SU_m\cap\SU_0,\sob{1,\frac{4}{3}}{})$\=/norms of $\xi_n$ are uniformly bounded also in terms of their $H^0(\SU_m\cap\SU_0,\ban{q}{})$\=/norms.
Rellich-Kondrachov theorem~\cite[Theorem~6.3]{Ad} implies that any sequence of non-trivial eigenfunctions whose $H^0(\SU_m\cap\SU_0,\ban{q}{})$\=/norms are bounded has a convergent subsequence, and that the limit is non-trivial.
By passing to a weakly convergent subsequence of $\U_n$ the resolvents in Lemma~\ref{global resolvent} converge.
Thus the limit of $\xi_n|_{\SU\cap\SU_0}$ extends to a holomorphic section of $\Q{E(D+D'),\U}$ where $D+D'$ is the limit of $D_{n}+ \sum_{m}D'_{n,m}$ and $\U$ is the weak limit of $\U_n(\sum_{m}D'_{n,m})$.
\end{proof}
\begin{lemma}\label{covering}
Let $\X$ be a compact Riemann surface with finitely many marked points $x_1,\ldots,x_K$, $\sigma$ a finite Baire measure on $\X$ and $\epsilon>0$.
Then $\X$ is covered by finitely many open discs $\SU_1,\ldots,\SU_M$ and an equal number of open annuli $\SV_1\supset\boundary \SU_1,\ldots,\SV_M\supset\boundary \SU_M$.
The discs are biholomorphic to $\mathbb{D}$ and pairwise disjoint.
The annuli obey $\sigma(\SV_1\cup\ldots\cup\SV_M) < \epsilon$ and their closures do not contain marked points.
\end{lemma}
\begin{proof}
Our cover $\X=\SO_1\cup\ldots\cup\SO_L$ from Remark~\ref{rem:special cover} is a finite collection of open subsets such that the biholomorphic maps $\SO_l\simeq\mathbb{D}$ continuously extend to homeomorphisms $\Bar{\SO}_l\simeq\Bar{\mathbb{D}}$.
We may further ensure that each marked point is contained in only one member of the cover and that any $\SO_l$ is not contained in $\SO_{l+1}\cup\ldots\cup\SO_L$.
We will now inductively cover the compact sets
\[
\emptyset \subset \X\setminus\left(\SO_2\cup\ldots\cup\SO_L\right)\subset\X\setminus\left(\SO_3\cup\ldots\cup\SO_L\right)\subset\ldots\subset\X\setminus\SO_L\subset\X
\]
by open pairwise disjoint discs and an equal number of corresponding open annuli, each containing the boundary of the corresponding disc.

Assume that $\X\setminus\left(\SO_l\cup\ldots\cup\SO_L\right)$ is covered by open pairwise disjoint discs $\SU_1,\ldots,\SU_{M'}$ and some annuli $\SV_1\supset\boundary \SU_1,\ldots,\SV_{M'}\supset\boundary \SU_{M'}$ all contained in $\SO_1\cup\ldots\cup\SO_{l-1}$.
Assume further that none $\SU_m$ is wholly contained in $\SO_l\cup\ldots\cup\SO_L$ and that the annuli do not contain any marked point and obey $\sigma(\SV_m) < \epsilon 2^{-m}$.

Let first show that $\boundary \SO_l\cup(\SO_l\cap\Bar{\SU}_m)$ is connected for any $m=1,\ldots,M'$.
For empty $\SO_l\cap\Bar{\SU}_m$ this follows from $\Bar{\SO}_l\sim\Bar{\mathbb{D}}$ and $\boundary \SO_l\simeq\boundary \mathbb{D}$.
In the remaining case $\SO_l\cap\Bar{\SU}_m\ne\emptyset$ we conclude that $\boundary \SO_l\cap\Bar{\SU}_m\ne\emptyset$, because otherwise $\Bar{\SU}_m=(\Bar{\SU}_m\cap\SO_l)\cup\Bar{\SU}_m\setminus\Bar{\SO}_l$ would be a disjoint union of two non empty open subsets and contradicts the connectedness of the closure of the connected $\SU_m$ (see~\cite[Chapter~3 Theorem~1.4]{Mun}).
Hence $\boundary \SO_l\cup\Bar{\SU}_m$ is as a non-disjoint union of two connected sets connected (see~\cite[Chapter 3 Theorem~1.3.]{Mun}).
The subset $\boundary \SO_l$ is closed in $\boundary \SO_l\cup(\SO_l\cap\Bar{\SU}_m)$.
So any decomposition into two disjoint closed subsets $\boundary \SO_l\cup(\SO_l\cap\Bar{\SU}_{M'})=\Set{A}\cup\Set{B}$ induces such a decomposition $\boundary \mathbb{D}\simeq\boundary \SO_l=(\Set{A}\cap\boundary \SO_l)\cup(\Set{B}\cap\boundary \SO_l)$.
Consequently $\boundary \SO_l$ is contained in one set which we now call $\Set{A}$ and disjoint from $\Set{B}$.
Since $\Set{A}$ and $\Set{B}$ are closed in $\boundary \SO_l\cup(\SO_l\cap\Bar{\SU}_m)$ which in turn is closed in $\boundary \SO_l\cup\Bar{\SU}_m$, the connected set $\boundary \SO_l\cup\Bar{\SU}_m$ is the union of the disjoint closed subsets $\Set{A}\cup(\Bar{\SU}_m\setminus\SO_n)$ and $\Set{B}$.
Thus $\Set{B}=\emptyset$ since $\Set{A}$ contains $\boundary \SO_l$.
This implies that $\boundary \SO_l\cup(\SO_l\cap\Bar{\SU}_m)$ is connected as well as the union $\boundary \SO_l\cup\bigcup_{m=1}^{M'}(\SO_l\cap\Bar{\SU}_m)$ of connected sets with non-empty intersection.

Now we claim that every non-empty connected component $\Set{C}$ of $\SO_l\setminus(\Bar{\SU}_1\cup\ldots\cup\Bar{\SU}_{M'})$ is simply connected.
Since bounded, open and connected subsets of $\mathbb{C}$ are simply connected if their complement is connected (compare~\cite[Chapter~4 Section~4.2]{Ah}), it suffices to show that $\Bar{\SO}_l\setminus\Set{C}$ is connected.
The connected components $\Set{C}$ of the open and locally connected set $\SO_l\setminus(\Bar{\SU}_1\cup\ldots\cup\Bar{\SU}_{M'})$ are open (see~\cite[Chapter~3 Theorem~4.1.]{Mun}).
Their closures $\Bar{\Set{C}}$ in $\Bar{\SO}_l$ are also connected.
If such a closure $\Bar{\Set{C}}$ is disjoint from $\boundary \SO_l\cup\bigcup_{m=1}^{M'}(\SO_l\cap\Bar{\SU}_m)$, then $\SO_l\setminus(\Bar{\SU}_1\cup\ldots\cup\Bar{\SU}_{M'})\supset\Bar{\Set{C}}=\Set{C}$ is open and closed in $\Bar{\SO}_l$ and therefore either empty or equal to $\Bar{\SO}_l$ which is impossible.
So $\Bar{\Set{C}}\cup\boundary \SO_l\cup\bigcup_{m=1}^{M'}(\SO_l\cap\Bar{\SU}_m)$ is a non-disjoint union of connected sets and connected.
The union of the closures of all connected components of $\SO_l\setminus(\Bar{\SU}_1\cup\ldots\cup\Bar{\SU}_{M'})$ with the exception of one given $\Set{C}$ and of $\boundary \SO_l\cup\bigcup_{m=1}^{M'}(\SO_l\cap\Bar{\SU}_m)$ is a union of connected sets with non-empty intersection and connected.
It is $\Bar{\SO}_l\setminus\Set{C}$ and the claim follows.

Therefore, due to the Riemann mapping theorem, any connected component $\Set{C}$ of $\SO_l\setminus(\Bar{\SU}_1\cup\ldots\cup\Bar{\SU}_{M'})$ is biholomorphic to $\mathbb{D}$.
Consider the inclusion
\[
\X\setminus\left(\SO_{l+1}\cup\!\ldots\!\cup\SO_L\cup\Bar{\SU}_1\cup\!\ldots\!\cup\Bar{\SU}_{M'}\cup\SV_1\!\ldots\!\cup\SV_{M'}\right)\subset \SO_l\setminus\left(\Bar{\SU}_1\cup\!\ldots\!\cup\Bar{\SU}_{M'}\right).
\]
We denote the smaller compact set by $K$.
We know that it is covered by finitely many connected components of the larger set.
All marked points in $\SO_l$ belong to $\Set{K}$, since all marked points are contained only in one of the sets $\SO_1,\ldots,\SO_L$.
Each such connected component $\Set{C}$ that is not disjoint from $\Set{K}$ is a tubular neighborhood of the compact subset $\Set{K}\cap\Set{C}$.
Since $\Set{C}$ itself is a disc, we decompose $\Set{C}$ into an open annulus and a compact disc containing $\Set{K}\cap\Set{C}$.
Because an open annulus contains infinitely many disjoint open annuli with pairwise disjoint closures, there exists one such annulus $\SV_{M'+1}\subset\SO_l$ which obeys $\sigma(\SV_{M'+1}) < \epsilon 2^{-M'-1}$ and whose closure contains none of the marked points in $\SO_l$.
Finally choose any disc $\SU_{M'+1}\subset\SO_l$ with $\boundary \SU_{M'+1}\subset\SV_{M'+1}$.
It must contain $\Set{K}\cap\Set{C}$ and $\SU_{M'+1}\cup\SV_{M'+1}$.

We apply this to every connected component $\Set{C}$ of $\SO_l\setminus(\Bar{\SU}_1\cup\ldots\cup\Bar{\SU}_{M'})$ which is not disjoint from $\Set{K}$.
The resulting discs and annuli in $\SO_l$ together with $\SU_1,\ldots,\SU_{M'}$ and $\SV_1,\ldots,\SV_{M'}$ cover $\X\setminus\left(\SO_{l+1}\cup\ldots\cup\SO_L\right)$.
Moreover, no disc is contained in $\SO_{l+1}\cup\ldots\cup\SO_L$.
After doing so for $l=1,\ldots,L-1$ the open discs $\SU_1,\ldots,\SU_{M}$ and annuli $\SV_1,\ldots,\SV_{M}$ have the desired properties.
\end{proof}
Finally we come to the main lemma of this chapter.
We show that Lemma~\ref{lem:embed sections} can be applied to any sequence of holomorphic $\qat$\=/line bundles whose degree and potentials are bounded.
Then we can consider the question of convergence of a sequence of sections because we can embed the sections in a common Banach space.
Furthermore, since bounded sequences are contained in compact sets, they always have convergent subsequences.
\begin{lemma}\label{singular holomorphic convergence}
Let $(E_n,\U_n)_{n\in\mathbb{N}}$ be a sequence of pairs of line bundles and potentials $\U_n\in\pot{E_n}^-$ on a compact Riemann surface $\X$.
Assume both that $|\deg E_n|$ and $\|\U_n\|_2$ are bounded.
Then there exists a line bundle $E$, a cover $\SU_0,\ldots,\SU_M$ and a sequence of divisors $(D_n)_{n\in\mathbb{N}}$ with $E_n\cong E(D_n)$ such that the conditions of Lemma~\ref{lem:embed sections} hold for all elements of a subsequence of $(E(D_n),\U_n)_{n\in\mathbb{N}}$.
\end{lemma}
\begin{proof}
Choose a line bundle $E$ that obeys the assumptions of Lemma~\ref{global resolvent} and let $C_{\frac{4}{3}}$ be the corresponding constant for $p=\frac{4}{3}$ constructed by that lemma.
The sequence $(E^{-1}\otimes E_n)_{n\in\mathbb{N}}$ of holomorphic line bundles with bounded degrees has a subsequence of constant degree $d$.
We may represent the sequence $(E^{-1}\otimes E_n)_{n\in\mathbb{N}}$ by a sequence of divisors $(D_n)_{n\in\mathbb{N}}$ (compare~\cite[Section~29.18.]{Fo}).
In particular $E_n= E(D_n)$.
Due to Riemann-Roch theorem~\cite[Theorem~16.9.]{Fo} every holomorphic line bundle of degree greater or equal to the genus $\genus$ of $X$ has at least one holomorphic section.
Therefore the divisors $D_n$ are linearly equivalent to differences of two divisors whose degrees are $\genus+|d|$ and $\genus$, which we may assume are effective divisors.
By the compactness of $\X$, any sequence of effective divisors of constant degree has convergent subsequence.
Therefore a subsequence of $(D_n)_{n\in\mathbb{N}}$ converges to $D$ with corresponding line bundle $E(D)$.

Due to the Banach-Alaoglu theorem~\cite[Theorem~IV.21]{RS1} and the Riesz Representation theorem~\cite[Chapter~13 Section~5]{Ro2} the sequence of bounded finite Baire measures $|\U_n |^2$ on $\X$ has a subsequence that converges to a finite Baire measure $\sigma$.
With respect to $\sigma$ only finitely many points can have mass larger than the constant $S_{\frac{4}{3}}^{-2}$ of Lemma~\ref{weakly continuous resolvent}.
We combine these points with the support of the divisor $D$ to finitely many marked points $x_1,\ldots,x_K$ of $\X$.
The preceding lemma, Lemma~\ref{covering}, when applied with $\varepsilon = C_{\frac{4}{3}}$ constructs a cover of $\X$ consisting of discs $\SU_1,\dots,\SU_M$ and a `remainder' set $\SU_0 := \SV_1 \cup\dots\cup \SV_M$.
Since the divisors $D_n$ converge to $D$ whose support is contained in the complement of $\overline{\SU_0}$, after passing to a subsequence the same is true of the divisors $D_n$.
So the divisors $D_n$ obey condition~(i) of Lemma~\ref{lem:embed sections}.
Similarly, we know that $\sigma(\SU_0) < \varepsilon$ so after passing to a further subsequence we know that $\| \U_n |_{\SU_0} \|_{\banach{2}}$ is also smaller that $\varepsilon$.
Thus $(E(D_n),\U_n)$ fulfil both conditions of Lemma~\ref{lem:embed sections}.
This proves the lemma.
\end{proof}

\chapter{Singular Holomorphic Sheaves}
\label{chapter:singular holomorphic sheaves}
In this chapter we introduce a type of sheaf called a singular holomorphic sheaf.
This concept is a natural extension of the concept of holomorphic structures on $\qat$\=/line bundles.
In fact, this is the culmination of Part~II: Given a sequence of sections $(\xi_n)_{n \in \mathbb{N}}$ of holomorphic structures $(E_n, \V_n)_{n \in \mathbb{N}}$ such that $\deg E_n$, $\| \V_n \|_2$, and $\| \xi_n \|$ (in the sense of Lemma~\ref{lem:embed sections}) are all bounded, Theorem~\ref{global limits} shows that it has a subsequence converging to a section of a singular holomorphic sheaf.
A direct consequence is the existence part of the proof of the~\ref{thm:main theorem}.

We begin this chapter with the definition of singular holomorphic sheaves and their duals, and provide a characterization in Lemma~\ref{singular divisor} in terms of divisors.
Theorem~\ref{singular holomorphic convergence 2} builds on Proposition~\ref{disc limits} to give the local version of Theorem~\ref{global limits}.
At the end of the chapter we also extend Riemann-Roch and Serre duality to singular holomorphic sheaves.

To motivate the following definition, let us consider an example from the theory of complex functions on a Riemann surface.
\begin{example}
Let $\X = \mathbb{C}$ be the line, whose holomorphic functions are $\Sh{O} = \Sh{O}(\X)$.
Take $\Sh{S}$ to be the sheaf of functions such that $\Sh{S}_x = \Sh{O}_x$ for all $x\neq 0$ and at the origin $\Sh{S}_0$ consists of all functions with power series of the form $a_{-1}z^{-1} + a_1z + a_2z^2 + \dots$, i.e.\ the constant term is missing.
By definition their stalks only disagree at the origin, where we see that the intersection of $\Sh{O}_0 \cap \Sh{S}_0 = z\mathbb{C}\{z\}$, and hence its codimension in both $\Sh{O}_0$ and $\Sh{S}_0$ is one.
This singular sheaf corresponds to the cusp curve $\X' = \{y^2 = x^3\}$ because if we normalize the cusp through the parameterization $y = z^3, x = z^2$ then we see that $\Sh{S} = z^{-1} \Sh{O}(\X')$.
In particular $\Sh{S}$ is locally free over $\Sh{O}(\X')$.
\end{example}

\begin{definition}
\label{def:singular sheaf}
\index{Singular holomorphic sheaf}
A \emph{singular holomorphic sheaf} $\Sh{S}$ on a Riemann surface $\X$ is a subsheaf of right $\qat$\=/modules of the meromorphic sections of a $\qat$\=/line bundle $E_\qat$ with potential $\V\in\pot{E}$ with the following properties:
\begin{enumeratethm}
\item The set of all $x$ with unequal stalks $\Sh{S}_x\ne\Q{E,\V,x}$ is discrete in $\X$.
\item For all $x\in X$ both codimensions of the intersection $\Sh{S}_x\cap\Q{E,\V,x}$ in the stalks $\Sh{S}_x$ and $\Q{E,\V,x}$ are finite and equal.
\end{enumeratethm}
On a compact Riemann surface $\X$ we define $\deg(\Sh{S})=\deg(E)$.
We refer to $\Q{E,\V}$ as the underlying sheaf.

For any such sheaf $\Sh{S}$ the Serre dual sheaf $\Sh{S}^\ast$ (see Chapter~\ref{chapter:riemann roch}) is defined as the subsheaf of meromorphic sections of $\Q{KE^{-1},\V\sd}$ such that for all $x\in\X$ the stalk $\Sh{S}^\ast_x$ contains all germs $\xi$ of meromorphic sections of $\Q{KE^{-1},\V\sd}$ for which the residue $\res_x\lp\chi,\xi\rp$ at $x$ (defined in Lemma~\ref{lem:10 pairing properties}(vii)) vanishes for all $\chi\in\Sh{S}_x$.
\end{definition}

In the case that a sheaf contains the holomorphic functions, we can simply define its degree through divisors as $\sum_x \dim (\Sh{S}_x/\Sh{O}_x)$.
Generally for locally free sheaves $\Sh{S}' \subset \Sh{S}$ the difference of their degrees is $\sum_x \dim (\Sh{S}_x/\Sh{S}'_x)$.
The idea of condition (ii) is that although $\Sh{S}$ does not contain $\Q{E,\V}$, the number of missing holomorphic functions has been added back as meromorphic functions:
\begin{align*}
\deg \Sh{S}
&= \sum_x \dim (\Sh{S}_x / \Sh{S}_x\cap\Q{E,\V,x}) + \deg (\Sh{S}\cap\Q{E,\V}) \\
&= \sum_x \dim (\Q{E,\V,x} / \Sh{S}_x\cap\Q{E,\V,x}) + \deg (\Sh{S}\cap\Q{E,\V})
= \deg \Q{E,\V}.
\end{align*}
\begin{lemma}\label{dual sheaf}
\index{Serre duality}
The Serre dual sheaf $\Sh{S}^\ast$ of a singular holomorphic sheaf $\Sh{S}$ (as defined above) is also a singular holomorphic sheaf whose underlying sheaf is the Serre dual sheaf of the underlying sheaf of $\Sh{S}$.
On a compact Riemann surface $\X$ we have $\deg(\Sh{S}^\ast)=2\genus-2-\deg(\Sh{S})$ where $\genus$ denotes the genus of $\X$.
\end{lemma}
\begin{proof}
In this proof $\chi$ denotes a meromorphic section of $\Q{E,\V}$ and $\xi$ a meromorphic section of $\Q{KE^{-1}\V\sd}$.
Moreover, we denote by $\ord_x(\chi)$ and $\ord_x(\xi)$ the vanishing order at $x\in\X$ of these meromorphic sections with poles counted negatively, respectively.
The residue of $\lp\chi,\xi\rp$ at $x$ defines a non-degenerate pairing between such $\chi$ and $\xi$.
If $\Sh{S}_x=\Q{E,\V,x}$, then $\Sh{S}^\ast_x=\Q{KE^{-1},\V\sd,x}$, and $\Sh{S}^\ast$ obeys condition~(i).
For $x\in\X$ with $\Sh{S}_x\ne\Q{E,\V,x}$ we define $n_x$ as $n_x:=-\min\{\ord_x(\chi)\mid\chi\in\Sh{S}_x\setminus\{0\}\}$.
Consequently
\begin{align}\label{codimension}
\{\chi\mid\ord_x(\chi)\ge -n_x\}&\supset\Sh{S}_x\setminus\{0\},&
\{\xi\mid\ord_x(\xi)\ge n_x\}&\subset\Sh{S}^\ast_x.
\end{align}
Condition~(ii) on $\Sh{S}_x$ states that the number of linear independent elements of $\Q{E,\V,x}$ which do not belong to $\Sh{S}_x$ is equal to the number of linear independent elements of $\Sh{S}_x$ which do not belong to $\Q{E,\V,x}$.
This condition is equivalent to the codimension of the first inclusion~\eqref{codimension} being equal to the codimension of the inclusion $\{\chi\mid\ord_x(\chi)\ge-n_x\}\supset\Q{E,\V,x}=n_x$.
Analogously condition~(ii) on $\Sh{S}^\ast_x$ with underlying sheaf $\Q{KE^{-1},\V\sd}$ is equivalent to the codimension of the second inclusion~\eqref{codimension} being equal to the codimension of the inclusion $\{\xi\mid\ord_x(\xi)\ge n_x\}\subset\Q{KE^{-1}\V\sd,x}$ which is also equal $n_x$.
Since the elements of $\Sh{S}_x$ gives relations on the elements of $\Sh{S}^\ast_x$ both codimensions are equal.
So $\Sh{S}^\ast$ obeys condition~(ii) with underlying sheaf $\Q{KE^{-1},\V\sd,x}$.

The last statement follows from $\deg KE^{-1}=2\genus-2-\deg E$.
\end{proof}
\begin{lemma}\label{singular divisor}
\index{Divisor}
Let $E$ be a holomorphic $\mathbb{C}$\=/line bundle and $\V\in\pot{E}$ on a Riemann surface $\X$ that is not necessarily compact.

For any singular holomorphic sheaf $\Sh{S}$ with underlying sheaf $\Q{E,\V}$ there exists a divisor $D = \sum_{x\in\X} k_x x$ on $\X$ such that $\Q{E(-D),\V(-D)}\subseteq \Sh{S} \subseteq \Q{E(D),\V(D)}$ and for all $x\in\X$ we have that the codimensions of $\Q{E(-D),\V(-D),x} \subseteq \Sh{S}_x$ and $\Sh{S}_x \subseteq \Q{E(D),\V(D),x}$ are both $k_x$.
For this $D$ the same properties also hold for $\Sh{S}_x^*$.

Conversely, if there exists a divisor and a sheaf $\Sh{S}$ such that $\Q{E(-D),\V(-D)} \subseteq \Sh{S} \subseteq \Q{E(D),\V(D)}$ and for all $x\in\X$ the codimensions of $\Q{E(-D),\V(-D),x} \subseteq \Sh{S}_x$ and $\Sh{S}_x \subseteq \Q{E(D),\V(D),x}$ are both $k_x$, then $\Sh{S}$ is a singular holomorphic sheaf with underlying sheaf $\Q{E,\V}$.
\end{lemma}
\begin{proof}
Choose an $x\in\X$ such that  $\Sh{S}_x\ne\Q{E,\V,x}$.
Due to Lemma~\ref{quotient dimension}, the conditions
\begin{equation}
\label{eq:sing_quotient_cond}
\begin{gathered}
\Q{E(-kx),\V(-kx),x}\subseteq \Sh{S}_x\subseteq\Q{E(nx),\V(nx),x},\\
\dim_{\qat}\Sh{S}_x/\Q{E(-kx),\V(-kx),x}=k
\quad\text{and}\quad
\dim_{\qat}\Q{E(kx),\V(kx),x}/\Sh{S}_x=k.
\end{gathered}
\end{equation}
imply that Definition~\ref{def:singular sheaf}~(ii) holds.

Conversely, if Definition~\ref{def:singular sheaf}~(ii) holds then we choose linear independent elements of $\Sh{S}_x$ and $\Q{E,\V,x}$ that generate their respective quotients by the intersection $\Sh{S}_x \cap \Q{E,\V,x}$.
Due to Lemma~\ref{quotient dimension} the maximum order of poles of sections of $\Sh{S}_x$ and the maximum order of roots of sections of $\Q{E,\V,x} \setminus \Sh{S}_x \cap \Q{E,\V,x}$ are well defined and the conditions~\eqref{eq:sing_quotient_cond} are satisfied with $k$ as the larger of the two.
\end{proof}
In general there are two different minimal effective divisors $D'=\sum_{x\in\X} k'_x x$ and $D''=\sum_{x\in\X} k''_x x$ that fulfil $\Q{E(-D''),\V(-D'')}\subseteq\Sh{S} \subseteq \Q{E(D'),\V(D')}$.
Due to Lemma~\ref{quotient dimension} the conditions~\eqref{eq:sing_quotient_cond} in the lemma are equivalent to the conditions
\begin{gather}\label{separate divisors}
\begin{aligned}
\dim_{\qat}\Sh{S}_x/\Q{E(-D''),\V(-D''),x}&=k''_x\\
\dim_{\qat}\Q{E(D'),\V(D'),x}/\Sh{S}_x&=k'_x
\end{aligned}\quad\text{for all }x\in\X.
\end{gather}
 Then $D=\max\{D',D''\}$ is the minimal divisor fitting Lemma~\ref{singular divisor}.

By the definition of a singular holomorphic sheaf the supports of both $D'$ and $D''$ consists of the points where the stalks are unequal in the sense of Definition~\ref{def:singular sheaf}~(i).
For the Serre dual sheaf $\Sh{S}^\ast$ these divisors $D'$ and $D''$ of minimal degree are interchanged.
Before we show that the limits in Proposition~\ref{disc limits} are indeed singular holomorphic sheaves at the end of this chapter let us now extend Serre duality and the Riemann-Roch theorem to these sheaves.

In both proofs we shall use for a given singular holomorphic sheaf $\Sh{S}$ with underlying sheaf $\Q{E,\V}$ the long exact sequences of the following short exact sequences of sheaves, where $\Q{E,\V}$ denotes the divisor $D$ constructed in Lemma~\ref{singular divisor}:
\begin{eqnarray}\label{short exact 1}
&\begin{tikzcd}[column sep=small, row sep=small]
0 \arrow[r] & \Q{E(-D),\V(-D)} \arrow[r, hook] & \Sh{S} \arrow[r, two heads] & \Sh{S} / \Q{E(-D),\V(-D)} \arrow[r] & 0
\end{tikzcd}\\
&\begin{tikzcd}[column sep=small, row sep=small]
0 \arrow[r] & \Sh{S}^\ast \arrow[r, hook] & \Q{K(E(-D))^{-1},\V(-D)\sd} \arrow[r, two heads] & \Q{K(E(-D))^{-1},\V(-D)\sd} / \Sh{S}^\ast \arrow[r] & 0
\end{tikzcd}\label{short exact 2}
\end{eqnarray}
The support of both quotient sheaves $\Sh{S}/\Q{E(-D),\V(-D)}$ and $\Q{K(E(-D))^{-1},\V(-D)\sd}/\Sh{S}^\ast$ is contained in the discrete set $\{x\in\X\mid\Sh{S}_x\ne\Q{E,\V}\}$. Thus by~\cite[Lemma~16.7]{Fo} their first cohomology groups vanish.

\begin{theorem}[Riemann-Roch Theorem]
\label{thm:riemann roch singular}
\index{Riemann-Roch theorem}
Let $\Sh{S}$ be the sheaf of holomorphic sections of a singular holomorphic sheaf in a compact Riemann surface $\X$ of genus $\genus$.
Then the finite dimensions of the Cohomology groups obey
\[
\dim H^0(\X,\Sh{S})-\dim H^1(\X,\Sh{S})=1-\genus+\deg(\Sh{S}).
\]
\end{theorem}

\begin{proof} 
The alternating sum of the dimensions of the long exact cohomology f the short exact sequence~\eqref{short exact 1} vanish Cite. This implies that $\dim H^0(\X,\Sh{S})-\dim H^1(\X,\Sh{S})$ is equal to
\[
\dim H^0(\X,\Q{E(-D),\V(-D)})+H^0(\X,\Sh{S}/\Q{E(-D),\V(-D)})-\dim H^1(\X,\Q{E(-D),\V(-D)})
\]
Due to the Equation~\eqref{separate divisors} and the Riemann-Roch Theorem~\ref{riemann roch} this is equal to $\deg(D)+\deg(E(-D))=\deg(E)=\deg(\Sh{S})$.
\end{proof}
In order to prepare the proof of Serre duality of singular holomorphic sheaves we first extend the description of the first cohomology groups in Lemma~\ref{lem:serre duality in terms of residue} to such sheaves. Let $\Sh{S}$ be a given singular holomorphic sheaf $\Sh{S}$ on a compact Riemann surface $\X$ and let $x_1,\ldots,x_L$ denote the corresponding points $x\in\X$ with $\Sh{S}_x\ne\Q{E,\V,x}$. As in Lemma~\ref{lem:serre duality in terms of residue} let $\Sh{M}_{E,\V}$ denote the $\V$\=/meromorphic sections of $E_\qat$ and $\Sh{S}_{x_1,\ldots,x_L}\subset\Sh{M}_{E,\V,x_1,\ldots,x_L}$ the cartesian products of the stalks of the sheaves $\Sh{S}$ and $\Sh{M}_{E,\V}$ at the points $x_1,\ldots,x_L$, respectively. Here $\iota$ denotes the natural embedding
\begin{align}\label{eq:embedding iota 2}
\iota:&H^0(\X,\Sh{M}_{E,\V})\cap H^0(\X\setminus\{x_1,\ldots,x_L\},\Sh{S})\hookrightarrow\Sh{M}_{E,\V,x_1,\ldots,x_L},
\end{align}
which maps any global section onto the cartesian product of its germs at $x_1,\ldots,x_L$.
\begin{lemma}\label{lem:serre duality in terms of residue 2}
Let on a compact Riemann surface $\X$ the following data be given: a singular holomorphic sheaf with underlying sheaf $\Q{E,\V}$ and all points $x_1,\ldots,x_L$ of the set $\{x\in X\mid\Sh{S}_x\ne\Q{E,\V,x}\}$. Furthermore, let $\Sh{S}_{x_1,\ldots,x_L}\subset\Sh{M}_{E,\V,x_1,\ldots,x_L}$ denote the cartesian products of the corresponding stalks of the sheaves $\Sh{S}$ and $\Sh{M}_{E,\V}$, respectively, and $\iota$ the embedding~\eqref{eq:embedding iota 2} as defined above. Then $H^1(\X,\Q{E,\V})$ is naturally isomorphic to the following quotient space
\begin{align}\label{eq:quotient space 2}
\Sh{M}_{E,\V,x_1,\ldots,x_L}\Big/\Big(\Sh{S}_{x_1,\ldots,x_L}+\iota\big[H^0(\X,\Sh{M}_{E,\V})\cap\!H^0(\X\setminus x_1,\ldots,x_L\},\Sh{S})\big]\Big).
\end{align}
\end{lemma}
\begin{proof}
For any $\eta=(\eta_1,\ldots,\eta_L)$ in this quotient space we choose small disjoint open neighborhoods $\SU_1,\ldots,\SU_L$ of $x_1,\ldots,x_L$ and supplement them by $\SU_0=\X\setminus\{x_1,\ldots,x_L\}$ to an open cover $\mathcal{U}$ of $\X$. Any three of its sets are disjoints. Since the restriction of each $\eta_l$ to $\SU_l\cap\SU_0=\SU_l\setminus\{x_l\}$ belongs to $H^0(\SU_l\cap\SU_0,\Sh{S})$ they define a cocycle in $C^1(\mathcal{U},\Sh{S})$. This cocyle is a coboundary, if and only if $\eta$ belongs to the space in the denominator of~\eqref{eq:quotient space 2}. By~\cite[12.5 Definition]{Fo} such cocycles define an element of $H^1(\X,\Sh{S})$. Due to~\cite[12.3 Lemma and 12.4 Lemma]{Fo} the corresponding map from the quotient space~\eqref{eq:quotient space 2} to $H^1(\X,\Sh{S})$ is an embedding.

Now we claim that in the long exact cohomology sequence of the exact sequence in~\eqref{short exact 1} we may replace $H^1(\X,\Sh{S})$ by the quotient~\eqref{eq:quotient space}, such that the new sequence is still exact. This implies that $H^1(\X,\Sh{S})$ and the quotient~\eqref{eq:quotient space 2} have the same dimensions, which in turn shows the surjectivity of the former embedding.

To prove the claim we first replace $H^1(\X,\Q{E(-D),V(-D)})$ by the isomorphic quotient space in Lemma~\ref{lem:serre duality in terms of residue}.
As explained in the proof of this Lemma, this quotient space is the first cohomology group of the corresponding cochain complex with respect to covers of the form $\X=\SU_0\cup\ldots\cup\SU_L$.
Here $\SU_1,\ldots,\SU_L$ are appropriate disjoint open neighborhoods of the given $x_1,\ldots,x_L$ and $\SU_0=\X\setminus\{x_1,\ldots,x_L\}$.
In contrast to the quotient space which is isomorphic to $H^1(\X,\Q{E(-D),V(-D)})$, the corresponding cochain complex depends on the choice of this cover, which in turn depends on the representatives of the involved germs.
With the possible exception of the connecting homomorphism, all maps of the long exact cohomology sequence of any short exact sequence of sheaves are defined as maps of the corresponding cochain complexes of a single fixed cover of $\X$.
The construction of the connecting homomorphism in~\cite[\S15.11]{Fo} is also a map between the corresponding cochain complexes with a single fixed cover.
In the proof of the exactness of the long cohomology sequence in~\cite[15.12~Theorem]{Fo} only one argument (denoted by g) uses several covers.
This argument gives the proof that the image of $H^1(\X,\Q{E(-D),\V(-D)})\to H^1(\X,\Sh{S})$  contains the kernel of $H^1(\X,\Sh{S})\to H^1(\X,\Sh{S}/\Q{E(-D),\V(-D)})$.
Below~\eqref{short exact 1}-\eqref{short exact 2} we argued for the vanishing of $H^1(\X,\Q{E(D),\V(D)}/\Sh{S})$.
Hence the long cohomology sequence with $H^1(\X,\Sh{S})$ replaced by the quotient~\eqref{eq:quotient space 2} and $H^1(\X,\Q{E(D),\V(D)})$ replaced by the isomorphic quotient space in Lemma~\ref{lem:serre duality in terms of residue} is also exact, if the induced map from the quotient representing $H^1(\X,\Q{E(-D),\V(-D)})$ into~\eqref{eq:quotient space 2} is surjective.
The natural isomorphism from the meromorphic sections of $\Q{E(-D),\V(-D)}$ onto the meromorphic sections of $\Sh{S}$ embeds the space in the denominator of the quotient representing $H^1(\X,\Q{E(-D),\V(-D)})$ into in the denominator of~\eqref{eq:quotient space 2}.
This implies the surjectivity of the map from the quotient space which represents $H^1(\X,,\Q{E(-D),\V(-D)})$ onto~\eqref{eq:quotient space 2}, and therefore also that~\eqref{eq:quotient space 2} is isomorphic to $H^1(\X,\Sh{S})$.
\end{proof}
The formula for the Serre duality in Lemma~\ref{lem:serre duality in terms of residue} has an obvious extension to a pairing between the elements $\eta$ of~\eqref{eq:quotient space 2} and $\xi\in H^0(\X,\Sh{S}^\ast)$.
\begin{definition}
\label{def:pairing singular}
\index{Pairing!Serre dual sheaf}
We define a pairing between elements of $H^1(\X,\Sh{S})$ and elements of $H^0(\X,\Sh{S}^\ast)$ by representing the elements of $H^1(\X,\Sh{S})$ by the corresponding elements of the quotient space~\eqref{eq:quotient space 2} defined in Lemma~\ref{lem:serre duality in terms of residue 2}.
The pairing between an element of $H^1(\X,\Sh{S})$ represented by $\eta$ and $\xi\in H^0(\X,\Sh{S}^\ast)$ is defined as
\[
\sum_{l=1}^L\res_{x_l}\lp\eta_l,\xi\rp.
\]
\end{definition}
\begin{theorem}[Serre duality]
\label{thm:serre duality singular}
\index{Serre duality}
Let $\Sh{S}$ be the sheaf of holomorphic sections of a singular holomorphic sheaf on a compact Riemann surface $\X$.
Then the pairing in Definition~\ref{def:pairing singular} is non-degenerate.
\end{theorem}
\begin{proof}
We dualize the left part of the long exact cohomology sequence induced by~\eqref{short exact 2} and arrange it together with the right part of the long exact cohomology sequence induced by~\eqref{short exact 1} to the following diagram:
\begin{equation*}
\begin{tikzcd}[column sep=small]
H^0(\X, \Sh{S}/\Q{E(-D),\V(-D)}) \arrow[d, "f_3", hook, two heads] \arrow[r] & H^1(\X, \Q{E(-D),\V(-D)}) \arrow[d, "f_4", hook, two heads] \arrow[r, two heads] & H^1(\X, \Sh{S}) \arrow[d, "f_5"] \arrow[r]& 0\\
\!H^0(\X,\!\Q{K(E(\!-\!D)\!)^{-1}\!,\!\V(\!-\!D)\sd\!}/\!\Sh{S}^\ast)^\ast\!\arrow[r] &\!H^0(\X,\!\Q{K(E(\!-\!D)\!)^{-1}\!,\!\V(\!-D)\sd}\!)^\ast\! \arrow[r, two heads] & H^0(\X,\!\Sh{S}^\ast\!)^\ast\!\arrow[r]&0
\end{tikzcd}
\end{equation*}
We represent both first cohomology groups in the upper row as the corresponding quotient spaces in Lemma~\ref{lem:serre duality in terms of residue} and~\ref{lem:serre duality in terms of residue 2}, respectively.
Therefore all spaces in the upper row are either cartesian products of stalks at the points $x_1,\ldots,x_L$ or quotients of such cartesian products divided by some subspace.
Before we define the vertical maps we explain the horizontal maps in the upper row as maps between the corresponding cartesian products of stalks at $x_1,\ldots,x_L$.

We claim that the connecting homomorphism in the upper row is induced by the embedding $\Sh{S}_{x_1,\ldots,x_L}\hookrightarrow\Sh{M}_{E,\V,x_1,\ldots,x_L}$ between these cartesian product both defined in Lemma~\ref{lem:serre duality in terms of residue 2}.
Note that this embedding indeed induces a map from $H^0(\X,\Sh{S}/\Q{E(-D),\V(-D)})$ into the quotient which represents $H^1(\X,\Q{E(-D),\V(-D)})$ in Lemma~\ref{lem:serre duality in terms of residue}, since this quotient is $\Sh{M}_{E,\V,x_1,\ldots,x_L}$ divided by a subspace containing $\Q{E(-D),\V(-D),x_1,\ldots,x_L}$.
The quotient $\Sh{S}_{x_1,\ldots,x_L}/\Q{E(-D),\V(-D),x_1,\ldots,x_L}$ is the domain of the connecting homomorphism.
We may represent any element in this domain by sections $(f_1,\ldots,f_L)$ of $\Sh{S}$ on disjoint open neighborhoods $\Set{U}_1,\ldots,\Set{U}_L$ of $x_1,\ldots,x_L$, respectively.
They define a cocycle in $C^1(\mathcal{U},\Q{E(-D),\V(-D)})$ with respect to the cover $\mathcal{U}$ of $\X$ by the open sets $\X\setminus\{x_1,\ldots,x_L\},\Set{U}_1,\ldots,\Set{U}_L$, which is for all $l=1,\ldots,L$ on $\X\setminus\{x_1,\ldots,x_L\}\cap\Set{U}_l$ equal to $f_l$.
This cocycle represents by definition of the connecting homomorphism~\cite[Section~15.11]{Fo} the corresponding image in $H^1(\X,\Q{E(-D),\V(-D)})$.
Now the claim follows.

The second map in the upper row is induced by the embedding of sheaves $\Q{E(-D),\V(-D)}\hookrightarrow\Sh{S}$, which maps the corresponding cartesian product of stalks $\Q{E(-D),\V(-D),x_1,\ldots,x_L}$ into $\Sh{S}_{x_1,\ldots,x_L}$.
Moreover, it also maps the space in the denominator of the quotient which represents $H^1(\X,\Q{E(-D),\V(-D)})$ in Lemma~\ref{lem:serre duality in terms of residue} onto a sub-space of the denominator of the quotient that represents $H^1(\X,\Sh{S})$ in Lemma~\ref{lem:serre duality in terms of residue 2}.
Hence the embedding $\Q{E(-D),\V(-D),x_1,\ldots,x_L}\hookrightarrow\Sh{S}_{x_1,\ldots,x_L}$ indeed induces the map between these two quotient spaces, which represents the second map in the upper row.

Let us now explain the three vertical maps which are all induced by pairings.
For $f_3$, observe that $\Sh{S}/\Q{E(-D),\V(-D)}$ and $\Q{K(E(\!-\!D)\!)^{-1}\!,\!\V(\!-\!D)\sd\!}/\!\Sh{S}^\ast$ are skyscraper sheaves supported at points $x$ with $\Sh{S}_x\ne\Q{E,\V,x}$ which we denote by $x_1,\ldots,x_L$ as in the preceding lemma.
Hence the space of sections is the sum of the stalks at these points.
The pairing is the sum over $x=x_1,\ldots,x_L$ of the following pairings
\[
\res_x\lp\cdot,\cdot\rp:(\Q{E(D),\V(D)}/\Sh{S})_x\times(\Sh{S}^\ast/\Q{KE^{-1}(-D),\V(D)\sd})_x\to\qat.
\]
By definition of $\Sh{S}^\ast$ these pairings are non-degenerate.
Hence $f_3$ is an isomorphism.
The second vertical map $f_4$ is induced by the Serre duality for the Serre dual sheaves $\Q{E(-D),V(-D)}$ and $\Q{K(E(-D))^{-1},\V(-D)\sd}$.
Therefore $f_4$ is also an isomorphism.
Finally the third vertical map $f_5$ is induced by the pairing defined in Definition~\eqref{def:pairing singular}.
If we use the representation of $H^1(\X,\Q{E(-D),V(-D)})$ in Lemma~\ref{lem:serre duality in terms of residue} with the corresponding formula for the Serre duality, then all three involved pairings are written as sums over the residues at the points $x_1,\ldots,x_L$.
The first two maps in the lower row are just the dual maps of the natural maps
\[H^0(\X,\Sh{S}\pa)\hookrightarrow H^0(\X,\Q{K(E(-D))^{-1},V(-D)\sd})\twoheadrightarrow H^0(\X,\Q{K(E(-D))^{-1},V(-D)\sd}/\Sh{S}\pa),\]
The given descriptions of all maps implies the commutativity of the diagram.
By the five lemma $f_5$ inherits the property to be an isomorphism from $f_3$ and $f_4$.
\end{proof}

\begin{theorem}\label{singular holomorphic convergence 2}
Let $(\U_n)_{n\in\mathbb{N}}$ in $\banach{2}(\mathbb{D},\qat^-)$ be a sequence fulfilling the condition of Proposition~\ref{disc limits} so that we have the corresponding polynomials $(P'_n)_{n\in\mathbb{N}}$.
Let $D'$ be the effective divisor of the roots of $P' := \lim P'_n$ which lie in $\overline{B(0,\frac{1}{2})}$.
Then after passing to a subsequence of $(\U_n)_{n\in\mathbb{N}}$ there exists a unique singular holomorphic sheaf $\Sh{S} \hookrightarrow \Q{\unity(D'), \U(D')}$ on $\mathbb{D}$ such that for any $\frac{1}{2}<r<1$:
\begin{enumeratethm}
\item\label{condition 1 on singular structure} for any sequence of unit-length sections $\xi_n \in H^0(B(0,r),\Q{\unity,\U_n})$, if $P'_n \xi_n$ converges in $\sobolev{1,p}(B(0,r),\qat)$ then $\lim P'_n \xi_n = P'\xi$ for some non-trivial section $\xi$ in $H^0(B(0,r),\Sh{S})$.

\item\label{condition 2 on singular structure} for any $\xi \in H^0(\mathbb{D},\Sh{S})$ the restriction $P'\xi|_{B(0,r)}$ is the limit of the sequence $P'_n \xi_n$ in $\sobolev{1,p}(B(0,r),\qat)$ for some sequence $\xi_n \in H^0(B(0,r),\Q{\unity,\U_n})$.
\end{enumeratethm}
\end{theorem}
\begin{proof}
Let the limits of $(P'_n)_{n\in\mathbb{N}}$ and $(P''_n)_{n\in\mathbb{N}}$ be respectively $P'$ and $P''$, and further let their roots be represented by the effective divisors $D'$ and $D''$.
By the construction in the proof of Proposition~\ref{disc limits}, the support of these divisors lies in $\overline{B(0,\frac{1}{2})} \cup \overline{B(0,4)}\setminus B(0,3)$.
First observe that the convergence of the third sequence of operators in the foregoing corollary and Cauchy's Integral Formula compare the equivalence of Theorem~\ref{cauchy formula}~(i) and (ii) implies that any limit of a sequence $(\xi_n)_n=(P''_n\phi_n)_{n\in\mathbb{N}}$ with $\phi_n\in H^0(B(0,r),\Q{\unity,{P''_n}^{-1}\U_n P''_n})$ is of the form $P''\phi$ with $\phi \in H^0(B(0,r),\Q{\unity,{P''_n}^{-1}\U P''_n})$.
Therefore, due to condition~(i), $\Sh{S}$ contains $\Q{\unity(-D''),\U(-D'')}$ as a subsheaf.

We now construct inductively additional sections in $H^0(\mathbb{D},\Sh{S})$ that induce a basis of $H^0(\mathbb{D},\Sh{S}/\Q{\unity(-D''),\U(-D'')})$.
Note that
\begin{align}\label{dim quotient}
  \dim H^0(\mathbb{D},\Q{\unity,\U_n}/\Q{\unity(-D''_n),\U_n(-D''_n)})&=\deg D''_n
\end{align}
by Lemma~\ref{quotient dimension}, so we expect that we will need $\deg D''$ additional sections to generate $\Sh{S}$.
We will describe these additional sections as the restrictions of elements of $H^0(\mathbb{P}^1,\Q{E(D' + d\infty),\U(D')})$ with $d = \deg D''$ and the spin bundle $\deg(E)=-1$ on $\mathbb{P}^1$.
For this purpose we recall from the beginning of Chapter~\ref{chapter:resolvents} that elements $\V\in\banach{2}(\mathbb{C},\qat)$ may be conceived as potentials in $\pot{E}$ and extend $\U_n$ from $\mathbb{D}$ to $\mathbb{P}^1$ by zero.
Note that $\U_n(d\infty) = \U_n$ since it vanishes in a neighborhood of infinity.




Choose any divisor $D$ with
\begin{gather}\label{induction divisor}
-D'' + d\infty \leq -D \leq d\infty
\end{gather}
and the corresponding sequence of divisors $D_n$ with constant degree and
\begin{align}\label{divisor sequence}
-D''_n + d\infty \leq -D_n \leq d\infty
\quad\text{with}\quad
\lim_{n\to\infty}D_n=D.
\end{align}
The roots of the polynomials were chosen so that by Lemma~\ref{transformed holomorphic structures} the limiting holomorphic structure has no sections on $\mathbb{P}^1$.
As in that Lemma, the inclusion $H^0(\mathbb{P}^1,\Q{E(D_n),\qk\Bar{\U}_n(D_n)\qk}) \subseteq H^0(\mathbb{P}^1,\Q{E(D''_n-d\infty),\qk\Bar{\U}_n(D''_n)\qk})$ implies similarly
\[\dim H^1(\mathbb{P}^1,\Q{E(-D_n),\U_n(-D_n)}) = 0 = \dim H^1(\mathbb{P}^1,\Q{E(-D),\U(-D)}).\]
In particular by Riemann-Roch 
\begin{align}
\label{eq:dim holo sections}
\dim H^0(\mathbb{P}^1,\Q{E(-D_n),\U_n(-D_n)})
= - \deg D
= \dim H^0(\mathbb{P}^1,\Q{E(-D),\U_n(-D)}).
\end{align}

We argue by induction over the degree of $-D$, whereby at each step we add a point to the divisor $-D$.
For the first step of the induction, consider a divisor $-D$ as in~\eqref{induction divisor} above with $-\deg D = 1$ and the corresponding sequence $-D_n$ in~\eqref{divisor sequence}.
Also for this first step, let us set $\Spa{B}_n = 0$ for all $n$ and $\Spa{B}=0$.
By Equation~\eqref{eq:dim holo sections}, for each $n$ there is a non-trivial, unit-length element
\begin{align}\label{induction bases}
\xi_n&\in H^0(\mathbb{P}^1,\Q{E(-D_n),\U_n(-D_n)}) \setminus \Spa{B}_n.
\end{align}
We know from Lemma~\ref{lem:embed sections} that the spaces $H^0(\mathbb{P}^1,\Q{E(-D_n),\U_n(-D_n)})$ are all embedded in $\Spa{A} := H^0(\SU_0,\Q{E,\U})$ and additionally that the sequence is contained in a compact set.
Therefore it has a convergent subsequence with
\begin{align}\label{limits induction}
\xi := \lim_{n\to\infty}\xi_n&\in H^0(\mathbb{P}^1,\Q{E(D' + d\infty),\U(D')}) \setminus \Spa{B}.
\end{align}
By (i) the restriction of $\xi$ to $\mathbb{D}$ is a non-trivial section of
\[
\xi|_{\mathbb{D}}\in H^0(\mathbb{D},\Q{E(D'),\U(D')}/\Q{E(-D''),\U(-D'')}).
\]

In subsequent iterations, the space $\Spa{A}$ remains the same, but we take $\Spa{B}_n$ to be the span of all the elements $\xi_n$ from Equation~\eqref{induction bases} constructed in previous steps.
The space $\Spa{B}$ is likewise the span of previous $\xi$'s, which by inductive hypothesis are linearly independent.
Now choose a new sequence of sections $\xi_n$ according to Equation~\eqref{induction bases}.
By applying Lemma~\ref{quotient}, we may assume that the new sequence is bounded and has no accumulation points in $\Spa{B}$, since that lemma only modifies $\xi_n$ by the addition of something in $\Spa{B}_n \subset H^0(\mathbb{P}^1,\Q{E(-D_n),\U_n(-D_n)})$ and rescaling.
Again by applying Lemma~\ref{lem:embed sections} we know that there is a subsequence converging to $\xi \not\in \Spa{B}$.
This guarantees that the dimension of $\Spa{B}$ increases at every step by one, and that at the end we have $\deg D''$ additional linearly independent sections.
This completes the induction.

We define $\Sh{S}$ on $\mathbb{D}$ as the unique subsheaf of $\Q{E(D'),\U(D')}$  which contains $\Q{E(-D''),\U(-D'')}$ as a subsheaf such that $H^0(\mathbb{D},\Sh{S}/\Q{E(-D''),\U(-D'')})$ is generated by the additional sections from the induction.
The divisors $D'$ and $D''$ we are exactly as in the situation of conditions~\eqref{separate divisors}.
In particular, $\dim_{\qat}\Sh{S}_x/\Q{E(-D''),\V(-D''),x}=k''_x$ holds by construction and the second condition holds since Lemma~\ref{quotient dimension} yields
\[
\dim_{\qat}\Q{E(D'),\V(D'),x}/\Q{E(-D''),\V(-D''),x}=k'_x + k''_x.
\]
Therefore $\Sh{S}$ is a singular holomorphic sheaf.

We now check that the singular holomorphic sheaf $\Sh{S}$ actually obeys both conditions~(i) and~(ii) of this theorem.

In the situation of~(i), consider first $\xi_n\in H^0(B(0,r), \Q{E(-D''_n),\U_n(-D''_n)})$.
The convergence of the third sequence of resolvents in Corollary~\ref{cor:other_resolvents_converge} and Cauchy's integral formula~\ref{cauchy formula} together implies in this case
\begin{align}\label{limit with roots}
\lim_{n\to\infty}\xi_n\in H^0(B(0,r), \Q{E(-D''),\U(-D'')}) \subset H^0(B(0,r), \Sh{S}).
\end{align}
For the remaining cases in the situation~(i) we consider
\[
[\xi_n] \in H^0(B(0,r), \Q{E,\U_n}/\Q{E(-D''_n),\U_n(-D''_n)}).
\]
It suffices to prove that $[\xi_n]$ converges to an element
\begin{align}\label{limit bases}
\lim_{n\to\infty}[\xi_n]&\in H^0(B(0,r), \Sh{S}/ \Q{E(-D''),\U(-D'')}).
\end{align}
The linear independent elements~\eqref{induction bases} modified by Lemma~\ref{quotient} build a basis of
\[
H^0(\mathbb{P}^1, \Q{E(d\infty),\U_n}/\Q{E(-D''_n+d\infty),\U_n(-D'')})
\]
since their number coincides with the dimension.
For $\frac{1}{2}<r<1$ they restrict to bases of $H^0(B(0,r), \Q{E,\U_n}/\Q{E(-D''_n),\U_n(-D''_n)})$, because the supports of the quotients $\Q{E,\U_n}/\Q{E(-D''_n),\U_n(-D''_n)}$ are contained in $B(0,1)$.
By construction, the limits~\eqref{limits induction} of this sequence of bases together with $\Q{E(-D''),\U(-D'')}$ generates $\Sh{S}$ and form a basis of $H^0(B(0,r), \Sh{S}/\Q{E(-D''),\U(-D'')})$.
This shows~\eqref{limit bases} and (i).

In the situation of (ii), again consider first $\xi \in H^0(\mathbb{D}, \Q{E(-D''),\U(-D'')})$ with $\frac{1}{2}<r<1$.
We choose a smooth function $f : \mathbb{D} \to [0,1]$ with compact support and which is $1$ on $B(0,r)$ and define
\[
{P''_n}^{-1}\xi_n :=\left.\Op{I}_{\mathbb{D}}({P''_n}^{-1}\U_n P''_n)((\barpartial f){P''}^{-1}\xi)\right|_{B(0,r)}.
\]
Since the resolvents are the right inverses of $\barpartial -\U_n$ we see that ${P''_n}^{-1}\xi_n$ is in the kernel of $\barpartial - {P''_n}^{-1}\U_n P''_n$ and thus $\xi_n$ is $\U_n$\=/holomorphic.
Due to Corollary~\ref{cor:other_resolvents_converge} ${P''_n}^{-1}\xi_n$ converges to ${P''}^{-1}\xi$.
This shows (ii) for the given $\xi$.
The remaining elements of $H^0(\mathbb{D}, \Sh{S}/\Q{E(-D''), \U(-D'')})$ were constructed above as limits~\eqref{limits induction} of sequences in $H^0(\mathbb{P}^1,\Q{E(d\infty),\U_n}) \supseteq H^0(B(0,r),\Q{E,\U_n})$.
This shows (ii).

To prove the uniqueness of $\Sh{S}$ we remark that the third sequence of resolvents in Corollary~\ref{cor:other_resolvents_converge} and all the inductively constructed elements~\eqref{induction bases} converge, by the choice of the subsequences.
So any singular holomorphic sheaf $\Sh{S}$ which obeys (i) and (ii) contains $\Q{E(-D''),\U(-D'')}$ as a subsheaf and has the restrictions of the limits~\eqref{limits induction} as sections on $\mathbb{D}$.
By construction these $\deg{D''}$ many limits are linearly independent elements of $H^0(\mathbb{D},\Q{E(D'),\U(D')})$.
So together with condition~\eqref{separate divisors} this uniquely determines a singular holomorphic sheaf with underlying sheaf $\Q{E,\U}$.
\end{proof}
\begin{lemma}\label{quotient}
Consider a Banach space $\Spa{A}$ and a sequence of closed subspaces $\Spa{B}_n$ with constant finite dimension $K$.
Let $\Spa{B}$ be a $K$\=/dimensional subspace such that each of its elements is the limit of a sequence $(b_n)_{n\in\mathbb{N}}$ with $b_n \in \Spa{B}_n$.
For any sequence $(a_n)_{n\in\mathbb{N}}$ with $a_n \in \Spa{A}\setminus\Spa{B}_n$ by Heine-Borel there exists a sequence $(b_n)_{n\in\mathbb{N}}$ with $b_n\in \Spa{B}_n$ such that $\|a_n + b_n\| = \|a_n \|_{\Spa{A}/\Spa{B}_n} > 0$ with the norms $\|x\|_{\Spa{A}/\Spa{Y}} := \inf \{ \|x + y\| \mid y \in \Spa{Y} \}$.
Then the normalized sequence $a'_n := \|a_n + b_n\|^{-1}(a_n + b_n)$ has $\liminf \| a'_n \|_{\Spa{A}/\Spa{B}} > 0$.
\end{lemma}
\begin{proof}
Suppose $\liminf  \| a'_n \|_{\Spa{A}/\Spa{B}} = 0$.
Then by passing to a subsequence there exists a sequence $(b'_n)_{n\in\mathbb{N}}$ in $\Spa{B}$ such that $\|a'_n-b'_n\| \to 0$.
From $\|b'_n\| \leq \|b'_n - a'_n\| + \|a'_n\|$ we see that $(b'_n)_{n\in\mathbb{N}}$ is a bounded sequence, and after passing to a further subsequence we have that $ \lim a'_n = \lim b'_n = b' \in \Spa{B}$.
Let $b''_n \in \Spa{B}_n$ with $\lim b''_n =: b'$.
Then $\|a'_n \|_{\Spa{A}/\Spa{B}_n} \leq \|a'_n - b''_n\| \to 0$.
This is a contradiction, since $\|a'_n \|_{\Spa{A}/\Spa{B}_n} = 1 = \|a'_n \|$ by the choice of $b_n$.
\end{proof}

\begin{theorem}
\label{global limits}
\index{Constrained Willmore}
On a compact Riemann surface $\X$ a sequence $(E_n,\U_n)_{n\in\mathbb{N}}$ of holomorphic $\mathbb{C}$\=/line bundles with bounded $|\deg E_n|$ together with potentials $\U_n\in\pot{E_n}^-$ with bounded norms $\|\U_n\|_2$ has a subsequence such that
\begin{enumeratethm}
\item limit $E := \lim\limits_{n\to\infty}E_n$ and weak limit $\U:=\lim\limits_{n\to\infty}\U_n\in\pot{E}^-$ exist.
\item any bounded sequence $(\xi_n)_{n\in\mathbb{N}}$ of images with respect to the embeddings 
\[
H^0(\X,\Q{E_n,\U_n})\hookrightarrow H^0(\SU\cap\SU_0, \sob{1,p}{E})
\]
constructed in Lemmas~\ref{lem:embed sections}--\ref{singular holomorphic convergence} has a convergent subsequence and the limit extends to a section of a singular holomorphic sheaf $\Sh{S}$ on $\X$ with underlying sheaf $\Q{E,\U}$.
\end{enumeratethm}
\end{theorem}
\begin{proof}
(i) has been proven in the proof of Lemma~\ref{singular holomorphic convergence}.

For the proof of (ii) let $\X=\SU_0\cup\SU=\SU_0\cup\SU_1\ldots\cup\SU_M$ be the cover whose existence was is proven in Lemma~\ref{covering}.
By Lemma~\ref{singular holomorphic convergence} and Lemma~\ref{lem:embed sections} the sequence $(\xi_n)_{n\in\mathbb{N}}$ in (ii) has a convergent subsequence.
For any $m=1,\ldots,M$ Theorem~\ref{singular holomorphic convergence 2} constructs on $\SU_m$ a subsequence and a singular holomorphic sheaf $\Sh{S}$ which obeys both conditions therein.
So there exists a subsequence for which this holds on every $\SU_m$ simultaneously.
The compact Riemann surface $\X$ is also covered by $\SU_0$ and slightly smaller discs contained in $\SU_1,\ldots,\SU_M$.
Then Theorem~\ref{singular holomorphic convergence 2}(i) shows that the limit of the subsequence is a section of $\Sh{S}$ on these smaller discs.
On $\SU_0$ Lemma~\ref{global resolvent} and the quaternionic version of Cauchy's integral formula~\ref{cauchy formula} imply that the limit is a section of $H^0(\SU_0, \Q{E,\U}) = H^0(\SU_0, \Sh{S})$.
Since the limit belongs to $\Sh{S}$ on each set of the cover, this proves (ii).
\end{proof}
In the following theorem we reformulate the preceding theorem in terms of admissible maps.

\begin{theorem}\label{compactness of conformal mappings}
Let $n\in\{3,4\}$.
On a compact Riemann surface $\X$ a sequence of admissible maps to $\mathbb{R}^n$ with bounded Willmore functional has a subsequence whose composition with a sequence of conformal transformations of $\Spa{S}^n\supset\mathbb{R}^n$ converges in $\bigcap_{1<p<2}\sobolev{2,p}\loc(\X\setminus\{x_1,\ldots,x_L\},\mathbb{R}^n)$ to an admissible map $F$ with Willmore functional bounded by the total mass $\mu_\willmore(\X)$ of the weak limit $\mu_\willmore$ of the Willmore densities.
The exceptional points $x_1,\ldots,x_L$ have mass $\mu_\willmore(\{x_l\})\geq 4\pi$.
Around every exceptional point $x_l$ there is a small ball $B(x_l,\epsilon)\subset\X$, local coordinate $z_l$, and order $n_l\in \mathbb{N}_0$ such that $|z_l|^{2n_l}dF$ converges in $\bigcap_{1<p<2}\sobolev{1,p}(B(x_l,\epsilon),\mathbb{R}^n)$.
The order $n_l$ is bounded by
$n_l\leq\frac{1}{4\pi}\mu_\willmore(\{x_l\})$.
\end{theorem}
\begin{proof}
Due to the quaternionic Weierstraß representation~\ref{thm:weierstrass} all conformal mappings are represented by two non-trivial sections of two paired holomorphic $\qat$\=/line bundles.
Let $\chi_n$ and $\xi_n$ be the sequences of paired sections corresponding to a minimizing sequence of the Willmore functional on the space of conformal immersions of a compact Riemann surface $\X$ into $\mathbb{R}^3$ or $\mathbb{R}^4$.
The Plücker formula (Corollary~\ref{pluecker formula}) implies that the degrees of the corresponding holomorphic $\qat$\=/line bundles are bounded from below (compare Corollary~\ref{thm:embedding}).
Since these two line bundles are paired, the degrees are also bounded from above.
Therefore Theorem~\ref{global limits} implies that both sequences after renormalising have convergent subsequences with non-trivial limits.
But it might happen that the corresponding limits have additional poles.
The corresponding Weierstraß representations describe admissible maps from $\X$ into $\Spa{S}^4\supset\mathbb{R}^4$ or $\Spa{S}^3\supset\mathbb{R}^3$.
The conformal transformations of $\Spa{S}^4\supset\mathbb{R}^4$ acts as quaternionic Möbius transformations on the corresponding holomorphic sections $\upsilon$ and $\phi$ constructed in the proof of Theorem~\ref{thm:weierstrass} by Möbius transformation~\eqref{eq:quaternionic moebius 1}.
The subgroup, which preserves the imaginary quaternions and preserve or reverse their orientation is described at the end of Chapter~\ref{chapter:weierstrass}.
Observe that in Theorem~\ref{global limits} we implicitly use translations and rotations of the immersions corresponding to rescalings of the two holomorphic sections of the two paired holomorphic $\qat$\=/line bundles.
If we use more general transformations~\eqref{eq:quaternionic moebius 2}, we may always achieve that the limit stays inside of $\mathbb{R}^3$ or $\mathbb{R}^4$.
So Theorem~\ref{global limits} implies the claim.
\end{proof}

\partpage{Isothermic and Constrained Willmore Admissible Maps}

\chapter{Strongly Isothermic Kodaira Triples}
\label{chapter:isothermic kodaira triples}

The subspace of conformal immersions from a compact Riemann surface into $\mathbb{R}^n$ inside the space of immersions has a singularity if and only if the immersion defines a strongly isothermic surface (compare~\cite{BPP}).
In this chapter we carry over this classical characterization of singularities to admissible maps into $\mathbb{R}^4\cong\qat$.
More precisely, Definition~\ref{def:kodaira} identifies the set of admissible maps with a subspace of a Banach space using the quaternionic Kodaira representation (in Chapter~\ref{chapter:isothermic weierstrass triples} we use the Weierstraß representation for the same purpose).
We prove in Theorem~\ref{thm:kod submanifold} that the subspace of Kodaira triples is a Banach manifold away from the strongly isothermic admissible maps.
On the other hand we see in Lemma~\ref{lem:isothermic1} that at strongly isothermic admissible maps the space of Kodaira triples has a singularity, but that the singularity is not too bad.
Indeed, we shall see that the space of admissible maps is locally the zero set of a real smooth function on a Banach space.
At the strongly isothermic Kodaira triples, the derivative of this function vanishes, and its Hessian is indefinite.
In particular, the tangent cone along the space of Kodaira triples differs from, but spans, the tangent space.

Now we introduce the space of admissible maps from $\X$ into $\qat$, as described by their Kodaira representations $(E,\V,\upsilon,\phi)$.
We know that the triple $(E,\V,\upsilon)$ uniquely determines the left normal of $F$, which motivates the space $\leftsp(\X,E)$.

\begin{definition}
\label{def:kodaira}
\index{Kodaira representation}
Let $E$ be a holomorphic $\mathbb{C}$\=/line bundle on a Riemann surface $\X$ and $1\leq p \leq 2$.
Then $\leftsp(\X,E)$ denotes the subspace of pairs
\[
\{(\V,\upsilon)\in\pot{E}\times H^0(\X,\sob{1,p}{E})\mid \upsilon\in H^0(\X,\Q{E,\V})\setminus \{0\}, \upsilon\text{ has no roots}\}.
\]
Further let $\kodaira(\X,E)$ denote the subspace of triples
\[
\{(\V,\upsilon,\phi)\in\pot{E}\times H^0(\X,\sob{1,p}{E})^{\times2}\mid \upsilon,\phi\in H^0(\X,\Q{E,\V})\setminus \{0\}, \upsilon\text{ has no roots}\}.
\]
\end{definition}

In the above definition, we have fixed the line bundle $E$ but note that the potential $\V \in \pot{E}$ is not necessarily anti-commuting.
This makes it possible to describe all admissible maps whose underlying bundles have the same degree as a subset of a single Banach space.
Even for a fixed line bundle, the Kodaira potentials of this form are not necessarily related by conjugation, unlike Kodaira data $(E',\V',\upsilon',\phi')$ with $\V' \in \pot{E'}^-$.
If we decompose $\V=\V^+ + \V^-$, then the extra freedom concerns $\V^+$.
Indeed, we use Corollary~\ref{gauge holomorphic structure} to transform $(\V,\upsilon,\phi) \in \kodaira(\X,E)$ into $(\V',\upsilon',\phi')\in\kodaira(\X,E')$ with $\deg(E') = \deg(E)$ and $|\V'| = |\V^-|$.
In this transformation, $\V^+$ completely determines $E'$, and the map from the infinite dimensional space $\V^+\in \pot{E}^+$ to the corresponding compact component of the finite-dimensional Picard group of $\X$ is by Corollary~\ref{gauge holomorphic structure} surjective.
Conversely, one can choose $\V^+\in \pot{E}^+$ freely in the pre-image of $E'$, and this choice uniquely fixes $(E, \V, \upsilon, \phi)$.
Therefore to describe all admissible maps, it is sufficient to consider the disjoint union of $\kodaira(\X,E_n)$ for any sequence $(E_n)$ of holomorphic line bundles $E_n$ on $\X$ with $\deg(E_n)=n$.
In essence we have partitioned the space of admissible maps by degree of the underlying bundle.
If we want to consider admissible maps whose Willmore energy $\willmore(F)$ is below a certain threshold, Corollary~\ref{thm:embedding} shows that $\deg(E)$ is bounded by
\[|\deg(E)+\genus-1|\le\frac{\willmore(F)}{4\pi}+\genus-1,\]
where $\genus>0$ is the genus of $\X$. 
In this case we would only need to take finitely many line bundles $E_n$ into account.

While this simplification can be performed at any point of $\kodaira(\X,E)$ it is not possible to perform it simultaneously at multiple points.
Even if $(\V',\upsilon',\phi')\in\kodaira(\X,E')$ with $\V' \in \pot{E}^-$, the commuting part of $\V$ will not be zero for nearby points.
Since locally these transformations are given by left multiplication with an invertible complex function the roots of $(\upsilon',\phi')$ are the same as the roots of $(\upsilon,\phi)$, and the quotients obey $\upsilon^{-1}\phi=\upsilon'^{-1}\phi'$.


Being subsets of Banach spaces, it is natural to ask whether $\leftsp(\X,E)$ and $\kodaira(\X,E)$ are themselves Banach spaces.

\begin{theorem}\label{thm:kod submanifold}
On a compact Riemann surface $\X$, let $F$ be an admissible map with line bundle $E$.
Then $\leftsp(\X,E)$ is a real Banach submanifold of $\pot{E}\times H^0(\X,\sob{1,p}{E})$.
Further if $F=\upsilon^{-1}\phi$ is not strongly isothermic then $\kodaira(\X,E)$ is locally at $(\V,\upsilon,\phi)$ a real Banach submanifold of $\pot{E}\times H^0(\X,\sob{1,p}{E})^{\times2}$.
\end{theorem}

Some preparation is necessary before we will be able to prove this theorem.
We will prove two lemmas which both refer to the concept of complemented subspaces of Banach spaces, Lemmas~\ref{lem:implicit-variant} and~\ref{lem:complemented2}.
We will also introduce in Definition~\ref{def:01 pairing} a pairing $\ls\cdot,\cdot\rs$ that produces a potential in $\pot{E}$ from sections of paired bundles $E_\qat$ and $KE^{-1}_\qat$.
Its properties are collected in Lemma~\ref{lem:01 pairing properties}.
However we start with a revision of the concept of strongly isothermic that was introduced in Definition~\ref{def:isothermic admissible}.

Recall Lemma~\ref{lem:strongly isothermic correspondence}.
It states that an admissible map is strongly isothermic if there exist a pair $(\upsilon\pa,\phi\pa) \in \Q{KE^{-1},\V\sd}$ where $\V \in \pot{E}^-$ is the Kodaira potential.
But as discussed above, in this chapter we allow $\V$ to have a non-zero commuting part.
If we apply Corollary~\ref{gauge holomorphic structure} to gauge away $\V^+$, how are $(\upsilon\pa,\phi\pa)$ are transformed?
Since the parts $\V^+$ and $(\V\sd)^+$ of the decompositions $\V=\V^++\V^-$ and $\V\sd=-\V^++(\V^-)\sd$ are the negatives of each other the corresponding cocycles in Corollary~\ref{gauge holomorphic structure} are the inverse of each other.
Hence $(\upsilon\pa,\phi\pa)$ are transformed into a ${\V'}\sd$\=/holomorphic pair $((\upsilon\pa)',(\phi\pa)')$ of the $\qat$\=/line bundle $(KE'^{-1})_\qat$.
Moreover $(\upsilon\pa,\phi\pa)$ satisfies~\eqref{eq:dual map condition} locally, if and only if $((\upsilon\pa)',(\phi\pa)')$ satisfies the condition for the Kodaira triple $(\V',\upsilon',\phi')\in\kodaira(\X,E')$.
In short, the conclusions of our earlier investigation into dual maps carry over to the present situation.
In the main though, we will use our freedom to make $V^+ = 0$ at the strongly isothermic surface, so the results of Chapter~\ref{chapter:isothermic} may be applied directly.

Now we prove the two lemmas.
Recall that a closed linear subspace $\Spa{K}$ of a Banach space $\Spa{A}$ is called \emph{complemented} if there exists a closed linear subspace $\Spa{K}'$ of $\Spa{A}$ such that $\Spa{K}\cap\Spa{K}'=\{0\}$ and $\Spa{K}+\Spa{K}'=\Spa{A}$ holds.
The first lemma establishes a variant of the implicit function theorem, and the second lemma simplifies the verification of the hypotheses of the first lemma.
\begin{lemma}\label{lem:implicit-variant}
Let $\Spa{A}$ and $\Spa{B}$ be Banach spaces, $f : U \to \Spa{B}$ a $C^1$\=/map on an open subset $U \subset \Spa{A}$ and $a_0\in U$.
Suppose that $f'(a_0)$ is surjective and the linear subspace $\ker f'(a_0)\subset\Spa{A}$ is complemented.
Then there exist open neighborhoods $U' \subset U$ of $a_0$ and $V \subset \Spa{B}$ of $f(a_0)$, an open subset $W \subset \ker f'(a_0)$ and a $C^1$\=/diffeomorphism $g: W \times V \to U'$ such that $f(g(w,b))=b$ for all $(w,b)\in W\times V$.
\end{lemma}
\begin{proof}
By assumption, there exists a complemented subspace $\Spa{A}'$ of $\ker f'(a_0)$, such that the following map is an isomorphism of Banach spaces:
\begin{align*}
\Phi: \ker f'(a_0)\times\Spa{A}'&\to\Spa{A},&(w,a')\mapsto a=w+a',
\end{align*}
 where $\ker f'(a_0) \times \Spa{A}'$ is the cartesian product of Banach spaces.
The map
\begin{align*}
h: \Phi^{-1}[U]&\to \ker f'(a_0) \times \Spa{B},&(w,a')&\mapsto(w,\Tilde{f}(w,a')) \;.
\end{align*}
with $\Tilde{f} = f \circ \Phi|_{\Phi^{-1}[U]}$ has at $\Phi^{-1}(a_0) = (w_0,a'_0)$ a derivative of block form
\[
h'(w_0,a'_0)=\begin{pmatrix} \unity_{\ker f'(a_0)} & \tfrac{\partial \Tilde{f}}{\partial w}(w_0,a'_0) \\ 0 & \tfrac{\partial \Tilde{f}}{\partial a'}(w_0,a'_0) \end{pmatrix}\;.
\]
\begin{tikzpicture}
\draw (-5+0.5,-2.5) node{$\ker f'(a_0)\times \Spa{A}'$};
\draw[<->] (-5-2,-1.5) -- (-5+2,-1.5) node[anchor=north]{$\ker f'(a_0)$};
\draw[<->] (-5-1,-2) -- (-5-1,2) node[anchor=west]{$\Spa{A}'$};
\fill (-6+1.5,0+0.5) circle(0.1) node[anchor=south]{$(w_0,a'_0)$};
\draw[dashed] (-6+1.5,0+0.5) circle [x radius=1.5, y radius=1.1, rotate=45];
\draw (-6+1.5+0.7,0+0.5+1) node[anchor=south west]{$\Phi^{-1}[U]$};

\draw (0.5,-2.5) node{$\Spa{A}$};
\draw[<->] (-2,-1.5) -- (3,-1.5) node[anchor=north]{$\ker f'(a_0)$};
\draw[line width=2] (1.5-1,-1.49) -- (1.5+1,-1.49) node[anchor=south east]{$W$};
\draw[<->] (-0.5,-2) -- (-1.5,2) node[anchor=west]{$\Spa{A}'$};
\fill (0+1,0+0.5) circle(0.1) node[anchor=south]{$a_0$};
\draw[dashed] (0+1,0+0.5) circle(1.4);
\draw (0+0.3,0+0.5+1) node[anchor=south east]{$U$};

\foreach \x in {-0.7,-0.3,0,0.5,0.8}
\draw (0,\x) sin (1,0.5+\x) cos (2,1+\x);
\draw (2,1.8) node[anchor=west]{$f(a)=b$};

\draw[dashed] (0.3,0.9) sin (1,0.5+0.7) cos (1.6,1+0.4);
\draw[dashed] (0.6,-0.5) sin (1,-0.3) cos (1.9,0);
\draw[dashed] (0.6,-0.5) -- (0.3,0.9);
\draw[dashed] (1.6,1+0.4) -- (1.9,0);
\draw (0.7,-0.6) node[anchor=west]{$U'$};

\draw (4.6,-2.5) node{$\Spa{B}$};
\draw[<->] (4.6,-2) -- (4.6,2);
\fill (4.6,0+0.5) circle (0.1) node[anchor=west]{$f(a_0)$};
\draw[line width=2] (4.6,0.5-1) -- (4.6,0.5+1) node[anchor=west]{$V$};

\draw[->] (-3,0) to node[above]{$\Phi$} (-2,0);
\draw[->] (3,0) to node[above]{$f$} (3+1,0);
\end{tikzpicture}

The inclusion $\Spa{A}' \hookrightarrow \Spa{A}$ induces a bijection $\Spa{A}' \to \Spa{A}/\ker f'(a_0)$, and $f'(a_0)$ induces a bijection $\Spa{A}/\ker f'(a_0) \to \Spa{B}$.
Hence $f'(a_0)|_{\Spa{A}'}:\Spa{A}'\to\Spa{B}$ is also bijective.
This shows that $\frac{\partial\Tilde{f}}{\partial a'}$ is bijective and has a continuous inverse by the Open Mapping Theorem~\cite[Theorem~III.10]{RS1}.
Hence $h'(w_0,a'_0)$ is invertible.
By the Inverse Function Theorem~\cite[Supplementary Material, Corollary~1]{RS1}, there exist open neighborhoods $U' \subset U$ of $a_0$, $V \subset \Spa{B}$ of $f(a_0)$ and $W \subset \ker f'(a_0)$ of $w_0$, such that $h \circ \Phi^{-1}$ is a $C^1$\=/diffeomorphism from $U'$ to $W \times V$.
The inverse map $g$ has the desired properties.
\end{proof}
%


\begin{lemma}\label{lem:complemented2}
Let $\Spa{A}$ and $\Spa{B}$ be Banach spaces with the projection $p_{\Spa{A}}: \Spa{A} \times \Spa{B}\to\Spa{A}$.
Then a closed linear subspace $\Spa{K} \subset \Spa{A} \times \Spa{B}$ is complemented, if $p_{\Spa{A}}|_{\Spa{K}}: \Spa{K} \to \Spa{A}$ is a Fredholm operator.
\end{lemma}
\begin{proof}
By abuse of notation, we identify the subspaces $\Spa{A}\times \{0\}$ and $\{0\} \times \Spa{B}$ of $\Spa{A}\times \Spa{B}$ with $\Spa{A}$ and $\Spa{B}$, respectively, and we consider also the natural projection $p_{\Spa{B}}: \Spa{A}\times \Spa{B}\to\Spa{B}$.
The finite-dimensional kernel $\Spa{K}' = \Spa{K}\cap \Spa{B}$ of $p_{\Spa{A}}|_\Spa{K}$ is a subspace both of $\Spa{B}$ and of $\Spa{K}$, and therefore has a complement $\Spa{B}'$ in $\Spa{B}$ and a complement $\Spa{K}''$ in $\Spa{K}$. $p_{\Spa{A}}|_{\Spa{K}''}$ is a bijective linear map onto a linear subspace $\Spa{A}'$ of $\Spa{A}$.
Because of $\Spa{A}' = \Imag(p_{\Spa{A}}|_{\Spa{K}})$, this space has finite codimension in $\Spa{A}$, and therefore has a complement $\Spa{A}''$ in $\Spa{A}$.
By the open mapping theorem, $p_{\Spa{A}}|_{\Spa{K}''}: \Spa{K}'' \to \Spa{A}'$ is an isomorphism of Banach spaces.
This implies that $\Spa{K}''$ is the graph of the continuous linear map $g = p_{\Spa{B}} \circ (p_{\Spa{A}}|_{\Spa{K}''})^{-1}: \Spa{A}' \to \Spa{B}$.
We now extend $g$ to the unique linear map $g: \Spa{A} \cong \Spa{A}' \times \Spa{A}'' \to\Spa{B}$ with $g|\Spa{A}'' = 0$.
The two maps 
$ \Spa{A}\times \Spa{B} \to \Spa{A}\times \Spa{B}$ given in block form by
\begin{align*}
\begin{pmatrix} \unity_{\Spa{A}} & 0 \\ g & \unity_{\Spa{B}} \end{pmatrix}
\quad\text{and}\quad
\begin{pmatrix} \unity_{\Spa{A}} & 0 \\ -g & \unity_{\Spa{B}} \end{pmatrix} \;
\end{align*}
are the inverse of each other.
The first maps $\Spa{A}' \times \Spa{K}'$ onto $\Spa{K}$.
Because $\Spa{A}' \times \Spa{K}'$ is complemented in $\Spa{A}\times \Spa{B}$ by $\Spa{A}'' \times \Spa{B}'$, also $\Spa{K}$ is complemented in $\Spa{A} \times \Spa{B}$.
\end{proof}
Because we want to show that $\leftsp(\X,E)$ is a real Banach submanifold of $\pot{E}\times H^0(\X,\sob{1,p}{E})$, it useful to work with a real non-degenerate $\pot{E}$\=/valued pairing that is related to the inner product on $\pot{E}$.
We denote this pairing by $\ls\cdot,\cdot\rs$ since its values are $(0,1)$\=/potentials, but it does not fit directly into the family of the three other pairings $\lz\cdot,\cdot\rz$, $\lp\cdot,\cdot\rp$ and $\lh\cdot,\cdot\rh$.
We first give the definition and then demonstrate the relationship to the inner product of $\pot{E}$
\begin{definition}
\label{def:01 pairing}
\index{Pairing!(0,1)}
As is standard, assume $\X$ is covered by open sets $\SO_l$ with charts $z_l:\SO_l\to\Omega_l$ such that $E$ is described by the holomorphic cocycle $f_{ml}$ with respect to this cover.
We define a bilinear map $\ls\cdot,\cdot\rs$ from $H^0(\X,\ban{2p/(2-p)}{KE^{-1}})\times H^0(\X,\ban{p/(p-1)}{E})$ to $\pot{E}$
through its local representatives
\begin{equation*}
\ls\upsilon\pa,\upsilon\rs_l=2\qk\upsilon\pa_l\Bar{\upsilon}_l.
\end{equation*}
This gives a global well-defined element of $\pot{E}$ (see Equation~\eqref{eq:potential transformation}) since
\begin{gather*}
\qk\upsilon\pa_m\Bar{\upsilon}_m
=\qk f_{ml}^{-1}\tfrac{dz_m}{dz_l}\upsilon\pa_l\overline{f_{ml}\upsilon_l}
=f_{ml}|f_{ml}|^{-2}\tfrac{d\Bar{z}_m}{d\Bar{z}_l}\qk\upsilon\pa_l\Bar{\upsilon}_l|f_{ml}|^2f_{ml}^{-1}
=f_{ml}\tfrac{d\Bar{z}_m}{d\Bar{z}_l}\qk\upsilon\pa_l\Bar{\upsilon}_l f_{ml}^{-1},
\end{gather*}
and they are square-integrable by Hölder's inequality.
\end{definition}
\begin{remark}
\label{rem:pa is not an operator}
The symbol $\ast$ on $\upsilon\pa$ is \emph{not} an operator.
In the sequel there are many instances of pairings of two holomorphic sections of paired $\qat$\=/line bundles.
To avoid the proliferation of variable names we use $\upsilon\pa$ as the variable name of the section that is paired with $\upsilon$.
The above is one example.
The pairing $\lp \chi,\psi \rp$ is an exception to this notation, since both $\chi$ and $\psi$ themselves will be paired with other sections later on.
\end{remark} 
The inner product $\langle\,\cdot\,,\,\cdot\,\rangle$ defined by Equation~\eqref{inner product} induces an isomorphism from the real Hilbert space $\pot{E}$ onto its dual. 
Therefore the key property~(iv) in the following lemma uniquely determines the pairing $\ls\cdot,\cdot\rs$:

\begin{lemma}
\label{lem:01 pairing properties}
Let $(\upsilon,\upsilon\pa)\in H^0(\X,\ban{\frac{p}{p-1}}{E})\times H^0(\X,\ban{\frac{2p}{2-p}}{KE^{-1}})$
\begin{enumeratethm}
\item
For any function $g:\X\to\mathbb{C}$ the pairing obeys
\begin{align*}
\ls g\upsilon\pa,\upsilon\rs=\Bar{g}\ls\upsilon\pa,\upsilon\rs
\quad\text{and}\quad
\ls\upsilon\pa,g\upsilon\rs=\ls\upsilon\pa,\upsilon\rs\Bar{g}.
\end{align*}
\item
For any function $h:\X\to\qat$ the pairing obeys
\begin{align*}
\ls\upsilon\pa h,\upsilon\rs&=\ls\upsilon\pa,\upsilon \Bar{h}\rs.
\end{align*}
\item
$\hspace{32mm}\displaystyle{\ls\upsilon,\upsilon\pa\rs=(\ls\upsilon\pa,\upsilon\rs)\sd.}$\vspace{2mm}
\item
For $\V\in\pot{E}$ with compact support the pairing obeys
\[\int_{\X}\Real\lh\upsilon\pa,\V\upsilon\rh=\langle\ls\upsilon\pa,\upsilon\rs,\V\rangle.\]
\end{enumeratethm}
\end{lemma}
\begin{proof}
Similar to the Lemmas~\ref{lem:10 pairing properties} and~\ref{lem:11 pairing properties}, the first three properties follow algebraically from the definition.

(iv): On $\SO_l$ the product of the local functions of $\V$, $\upsilon$ and $\Bar{\upsilon}\pa$ is integrable, since the inverse exponents sum up to $\frac{1}{2}+\frac{p-1}{p}+\frac{2-p}{2p}=1$.
Both sides are integrals of such products in different orders.
Because the real part of a product of quaternions is cyclic, the integrands of both sides coincide (compare Equation~\eqref{eq:inner product}):
\begin{align*}
\Real\lh\upsilon\pa, \V\upsilon\rh_l
&= \Real\bigl(\Bar{\upsilon}\pa_l\qj d\Bar{z}_l\wedge dz_l \V_l\upsilon_l\bigr)
=\Real\bigl( \upsilon_l\Bar{\upsilon}\pa_l(-2\qk\dmu) \V_l\bigr)
=\Real\bigl( \overline{2\qk\upsilon\pa_l\Bar{\upsilon}_l} \V_l\bigr)\dmu.\qedhere
\end{align*}
\end{proof}

\begin{proof}[Proof of Theorem~\ref{thm:kod submanifold}.]
We first prove the claim about $\leftsp(\X,E)$.
We will apply Lemma~\ref{lem:implicit-variant} to the map
\begin{align}\label{immersion kodaira simple}
\Theta:\pot{E}\times H^0(\X,\sob{1,p}{E})&\to H^0(\X,\forms{0,1}{}\ban{p}{E}), &
(\V,\upsilon)&\mapsto((\delbar{E}-\V)\upsilon)
\end{align}
whose zero set is $\leftsp(\X,E)$.
 We first verify for any given pair $(\V,\upsilon) \in \leftsp(\X,E)$, that $\Spa{K} = \ker \Theta'(\V,\upsilon)$ is complemented.
 We do so via Lemma~\ref{lem:complemented2} and show that the restriction $p_{\pot{E}}|_{\Spa{K}}$ of the projection $p_{\pot{E}}: (\var \V, \var \upsilon) \mapsto \var \V$ is a Fredholm operator.
 The partial derivatives of $\Theta$ are $\tfrac{\partial \Theta(\V,\upsilon)}{\partial \V}(\var \V) = -\var \V \,\upsilon$ and $\tfrac{\partial \Theta(\V,\upsilon)}{\partial \upsilon}(\var \upsilon) = (\delbar{E}-\V)\var\upsilon$.
 So $\Spa{K}$ is the subspace
\begin{equation}\label{eq:Theta'-kernel simple}
\bigl\{(\var \V,\var\upsilon)\;\big|\; \var \V\upsilon = (\delbar{E}-\V)\var \upsilon \bigr\}
\end{equation}
of $\pot{E}\times H^0(\X,\sob{1,p}{E})$.
The operator $\delbar{E}-\V$ is a Fredholm operator from $H^0(\X,\sob{1,p}{E})$ to $H^0(\X,\forms{0,1}{}\ban{p}{E})$, see Chapter~\ref{chapter:riemann roch}.
 Therefore $\ker p_{\pot{E}}|_{\Spa{K}}$ is the kernel of the Fredholm operator
\begin{align}\label{eq:W2-L2-map simple}
H^0(\X,\sob{1,p}{E})&\to H^0(\X,\forms{0,1}{}\ban{p}{E}),&\var\upsilon&\mapsto(\delbar{E}-\V)\var\upsilon
\end{align}
and therefore finite-dimensional.
We now show that the codimension of the image of $p_{\pot{E}}|_{\Spa{K}}$ is less or equal to the finite dimension $k$ of the cokernel of~\eqref{eq:W2-L2-map simple}.
The image of $p_{\pot{E}}|_{\Spa{K}}$ contains all $\var \V \in \pot{E}$ which are mapped by  $\tfrac{\partial \Theta(\V,\upsilon, \phi)}{\partial \V}$ into the image of the Fredholm operator~\eqref{eq:W2-L2-map simple}.
The latter image is the intersection of the kernels of $k$ linear forms $\alpha_1,\dotsc,\alpha_d$, and the image of $p_{\pot{E}}|_{\Spa{K}}$ the intersection of the kernels of the linear forms $\alpha_1 \circ  \tfrac{\partial \Theta(\V,\upsilon, \phi)}{\partial \V},\dotsc,\alpha_d \circ  \tfrac{\partial \Theta(\V,\upsilon, \phi)}{\partial \V}$.
Hence the codimension of the image of $p_{\pot{E}}|_{\Spa{K}}$ is at most $k$.
This completes the proof that $p_{\pot{E}}|_{\Spa{K}}$ is Fredholm.
Therefore Lemma~\ref{lem:complemented2} applies and $\Spa{K}$ is complemented.

In a second step we verify for any pair $(\V,\upsilon)\in\leftsp(\X,E)$ the other assumption in Lemma~\ref{lem:implicit-variant}, namely that $\Theta'(\V,\upsilon)$ is surjective.
The image of $\Theta'(\V,\upsilon)$ is the sum of the images of $\tfrac{\partial \Theta(\V,\upsilon,\phi)}{\partial \V}$ and~\eqref{eq:W2-L2-map simple}.
Hence it suffices to show that the linear map 
\begin{align}\label{eq:deltaW-map simple}
\tfrac{\partial \Theta(\V,\upsilon)}{\partial \V}: \pot{E}&\to H^0(\X,\forms{0,1}{}\ban{p}{E}),&\var \V&\mapsto -\var \V \, \upsilon
\end{align}
is surjective onto the finite-dimensional cokernel of~\eqref{eq:W2-L2-map simple}.
In finite dimensions the surjectivity follows from the injectivity of the dual map.
In Chapter~\ref{chapter:riemann roch} we defined Hermitian pairings as the integrals in Lemma~\ref{lem:11 pairing properties}(iii) of the $(1,1)$\=/forms $\lh\cdot,\cdot\rh$, such that $-\delbar{KE^{-1}} +\V\sd$ acting on $H^0(\X,\ban{p/(p-1)}{KE^{-1}})$ is the adjoint operator of $\delbar{E}-\V$.
The dual of the cokernel of~\eqref{eq:W2-L2-map simple} is the kernel of the adjoint operator.
By Theorem~\ref{cauchy formula}(i)$\Leftrightarrow$(iv), the kernel of the adjoint operator is equal to $H^0(\X,\Q{KE^{-1},\V\sd})$.
The Hermitian pairing of any $\upsilon\pa\in H^0(\X,\Q{KE^{-1},\V\sd})$ in this kernel, with the image of $\var\V\in\pot{E}$ under the map~\eqref{eq:deltaW-map simple} is equal to $-\int_{\X}\lh\upsilon\pa,\var \V\upsilon\rh$.

Since $\upsilon$ is non-zero and invertible, the vanishing of this pairing with the elements of the image of~\eqref{eq:deltaW-map simple} implies $\upsilon\pa=0$.
Hence $\Theta'(\V,\upsilon)$ is surjective and Lemma~\ref{lem:implicit-variant} applies and shows that $\leftsp(\X,E)$ is a real Banach submanifold of $\pot{E}\times H^0(\X,\sob{1,p}{E})$.

For the case of $\kodaira(\X,E)$ the proof is similar.
 This time we apply Lemma~\ref{lem:implicit-variant} to the map
\begin{equation}\label{immersion kodaira}
\begin{aligned}
\Theta:\pot{E}\times H^0(\X,\sob{1,p}{E})^{\times 2}&\to H^0(\X,\forms{0,1}{}\ban{p}{E})^{\times2}\\(\V,\upsilon,\phi)&\mapsto((\delbar{E}-\V)\upsilon,(\delbar{E}-\V)\phi)
\end{aligned}\end{equation}
whose zero set is $\kodaira(\X,E)$.
We need to verify for any given Kodaira triple $(\V,\upsilon,\phi) \in \kodaira(\X,E)$, that $\Spa{K} = \ker \Theta'(\V,\upsilon, \phi)$ is complemented and $\Theta'(\V,\upsilon,\phi)$ is surjective.
$\Spa{K}$ is the subspace
\begin{equation}\label{eq:Theta'-kernel}
\bigl\{(\var \V,\var\upsilon,\var\phi)\;\big|\;(\var \V\upsilon,\var \V\phi)=((\delbar{E}-\V)\var \upsilon,(\delbar{E}-\V)\var\phi)\bigr\}
\end{equation}
and $\ker p_{\pot{E}}|_{\Spa{K}}$ is the kernel of the Fredholm operator
\begin{align}\label{eq:W2-L2-map}
H^0(\X,\sob{1,p}{E})^{\times 2}&\to H^0(\X,\forms{0,1}{}\ban{p}{E})^{\times 2},&(\var\upsilon, \var\phi)&\mapsto\bigl((\delbar{E}-\V)\var\upsilon, (\delbar{E}-\V)\var\phi\bigr)
\end{align}
and therefore finite-dimensional.
Again the codimension of the image of $p_{\pot{E}}|_{\Spa{K}}$ is less or equal to the finite dimension of the cokernel of~\eqref{eq:W2-L2-map}.
Therefore Lemma~\ref{lem:complemented2} applies and $\Spa{K}$ is complemented.
The image of $\Theta'(\V,\upsilon,\phi)$ is the sum of the images of $\tfrac{\partial \Theta(\V,\upsilon,\phi)}{\partial \V}$ and~\eqref{eq:W2-L2-map}, it suffices to show that the linear map 
\begin{align}\label{eq:deltaW-map}
-\frac{\partial \Theta(\V,\upsilon,\phi)}{\partial \V}: \pot{E}&\to H^0(\X,\forms{0,1}{}\ban{p}{E})^{\times2},&\var \V&\mapsto (\var \V \, \upsilon,\var \V \,\phi)
\end{align}
is surjective onto the finite-dimensional cokernel of~\eqref{eq:W2-L2-map}.
Again we conclude the surjectivity from the injectivity of the dual map.
The dual of the cokernel of~\eqref{eq:W2-L2-map} is the kernel of the adjoint operator, whose kernel is equal to $(\upsilon\pa,\phi\pa) \in H^0(\X,\Q{KE^{-1},\V\sd})^{\times 2}$ by Theorem~\ref{cauchy formula}(i)$\Leftrightarrow$(iv).
The real part of the Hermitian pairing of any $(\upsilon\pa,\phi\pa)$ in this kernel with $(\var \V \upsilon, \var \V \phi)$ is 
\[
\int_\X\Real\left(\lh\upsilon\pa,\var \V\upsilon\rh+\lh\phi\pa,\var \V \phi\rh\right)=\langle\ls\upsilon\pa,\upsilon\rs+\ls\phi\pa,\phi\rs,\var \V\rangle.
\]
Hence the linear map
\begin{align}\label{eq:eta-xi-pairing}
H^0(\X,\ban{p/(p-1)}{KE^{-1}})^{\times 2}&\to\pot{E},&
(\upsilon\pa,\phi\pa) &\mapsto \ls\upsilon\pa,\upsilon\rs + \ls\phi\pa,\phi\rs
\end{align}
is the dual of~\eqref{eq:deltaW-map} with respect to this real pairing and the inner product of $\pot{E}$.
Therefore $\Theta'(\V,\upsilon,\phi)$ is surjective if and only if the map 
\begin{align}\label{eq:singularity1}
H^0(\X,\Q{KE^{-1},\V\sd})^{\times 2}&\to\pot{E},&(\upsilon\pa,\phi\pa)&\mapsto \ls\upsilon\pa,\upsilon\rs + \ls\phi\pa,\phi\rs
\end{align}
is injective.
But writing this out in components shows
\[
\ls\upsilon\pa,\upsilon\rs + \ls\phi\pa,\phi\rs = 0
\quad \Leftrightarrow \quad
\upsilon\pa_l \Bar{\upsilon}_l + \phi\pa_l \Bar{\phi}_l = 0,
\]
which is exactly~\eqref{eq:dual map condition}.
Lemma~\ref{lem:strongly isothermic correspondence} applied to this situation now tells us at $(\V,\upsilon,\phi)\in\kodaira(\X,E)$ that $F=\upsilon^{-1}\phi$ strongly isothermic if and only if there exist a non-trivial element of the kernel of~\eqref{eq:singularity1}.
As we have assumed that $F$ is not strongly isothermic, then~\eqref{eq:singularity1} is injective and therefore the derivative of $\Theta$ is surjective.
Hence Lemma~\ref{lem:implicit-variant} applies and shows that $\kodaira(\X,E)$ is there locally a real Banach submanifold of $\pot{E}\times H^0(\X,\sob{1,p}{E})^{\times2}$.
\end{proof}

It remains to consider the strongly isothermic triples $(\V,\chi,\psi)\in\kodaira(\X,E)$.
If the kernel of~\eqref{eq:singularity1} is higher dimensional (the totally umbilic case), then by Lemma~\ref{thm:isothermic local characterization} we have map $F$ to a round sphere or plane.
We shall see in the next chapter that the derivative of the Willmore functional vanishes on the tangent space of $\kodaira(\X,E)$ at $F$ (as defined below).
It is not clear for a minimizer of the Willmore functional with a one-dimensional kernel of~\eqref{eq:singularity1} whether the derivative of the Willmore functional vanishes on the tangent space.

However, it will turn out that for these surfaces the tangent cone still spans the tangent space.
There are several definitions of the tangent cone in the literature.
The tangent cone defined here is contained in the tangent cone by most common definitions (for example, see~\cite[Ch.~II, \S 1.5]{Shafarevich} or~\cite[Ch.~7, \S 1, Definition 1G]{Wh}).  
In the current situation, the following definition applies with $\kodaira(\X,E)$ as the zero level set of $\Theta$ from~\eqref{immersion kodaira}.

\begin{definition}\label{D:kodaira-tangentspace-tangentcone}
Let $\Spa{A}$ and $\Spa{B}$ be Banach spaces, and $f: U \to \Spa{B}$ a $C^1$\=/map on an open subset $U \subset \Spa{A}$.
For any $a_0\in U$ the tangent space of the level set $f^{-1}[\{f(a_0)\}]$ at $a_0$ is defined as the kernel of $f'(a_0)$.
The tangent cone at $a_0$ is the cone of tangent vectors
\[\left\{\dot{\gamma}(0)\mid\gamma\in C^1((-\epsilon,\epsilon),\Spa{B})\cap C((-\epsilon,\epsilon),f^{-1}[\{f(a_0)\}])\text{ with }\gamma(0)=a_0\right\}.\]
\end{definition}
We characterize triples in $\kodaira(\X,E)$ with a kernel of~\eqref{eq:singularity1} of dimension $\leq1$ in the following theorem, which is proved in the remainder of this chapter.
\begin{theorem}
\label{thm:kodaira isothermic characterization}
\index{Isothermicity!Singularity of $\kodaira(\X,E)$}
On a compact Riemann surface let $(\V,\upsilon,\phi)\in\kodaira(\X,E)$ have a kernel of~\eqref{eq:singularity1} of real dimension $\le1$.
Then at $(\V,\upsilon,\phi)$ the tangent cone spans the tangent space and there is a smooth path in $\kodaira(\X,E)$ through $(\V,\upsilon,\phi)$ that can meet the set of strongly isothermic surfaces only at $(\V,\upsilon,\phi)$.
Moreover, the following are equivalent:
\begin{enumeratethm}
\item\label{case a} The admissible map $F = \upsilon^{-1}\,\phi$ is strongly isothermic.
\item\label{case b} The kernel of the map~\eqref{eq:singularity1} has real dimension one.
\item\label{case c} The tangent cone of $\kodaira(\X,E)$ at $(\V,\upsilon,\phi)$ is a proper subset of the tangent space.
\item\label{case d} $\kodaira(\X,E)$ is not a submanifold of $\pot{E}\times H^0(\X,\sob{1,p}{E})^{\times2}$ at $(\V,\upsilon,\phi)$.
\end{enumeratethm}
\end{theorem}
The idea of the proof is to describe $\kodaira(\X,E)$ in the neighborhood of such strongly isothermic surfaces as a level set of a single real function $\lambda$ on a Banach manifold, because
Lemma~\ref{lem:tangent cone spans} shows that at critical points of $\lambda$ with indefinite Hessian the tangent cone is a proper subset of the tangent space but still spans it.
In Lemma~\ref{lem:isothermic1} we construct the function $\lambda$ whose level set is $\kodaira(\X,E)$ and compute its derivative and Hessian~\eqref{eqn:kod second derivative}.
In order to show that the Hessian is indefinite, we blow-up the surface at a generic point to a plane and use Lemma~\ref{lem:sobolev local behavior} to control the convergence of the blow-up limit.
In Lemma~\ref{second variation plane} we define an indefinite bilinear form for the plane.
Finally, Lemma~\ref{lem:indefinite} makes a perturbation argument that the Hessian of $\lambda$ is also indefinite.


\begin{lemma}
\label{lem:tangent cone spans}
Let $\Spa{A}$ be a Banach space, $\lambda : U \to \mathbb{R}$ a smooth map on an open subset $U \subset \Spa{A}$ containing $0$, and $\gamma : (-\varepsilon,\varepsilon) \to U$ a $C^2$\=/path with $\gamma(0)=0$.
 Suppose that  $\lambda(0)=0=\lambda'(0)$.
 Then 
\begin{align*}
\left. \frac{d^2}{dt^2}\right|_{t=0} \lambda(\gamma(t)) = \lambda''(0)(\gamma'(0),\gamma'(0)).
\end{align*}
Suppose further that $\lambda''(0)$ is indefinite.
 Then $\Spa{A}$ is spanned by the tangent cone of $\lambda^{-1}[\{0\}]$ at $0$.

Furthermore, there exists $\gamma \in C^\infty((-\epsilon, \epsilon),\Spa{U})\cap C((-\epsilon,\epsilon),\lambda^{-1}[\{0\}])$ with $\gamma(0)=0$ and $\lambda'(\gamma(t)) \neq 0$ for all $t \in (-\epsilon, \epsilon) \setminus \{0\}$.  
\end{lemma}
\begin{proof}
For the first statement we compute directly
\begin{align*}
\frac{d}{dt} \lambda(\gamma(t)) &= \lambda'(\gamma(t))(\gamma'(t)), \\
\frac{d^2}{dt^2} \lambda(\gamma(t)) &= \lambda''(\gamma(t))(\gamma'(t),\gamma'(t)) + \lambda'(\gamma(t))(\gamma''(t)).
\end{align*}
The last term vanishes at $t=0$ by assumption, proving the formula for the second derivative.

For the remaining statements, we note that the tangent cone is contained in the set of isotropic vectors of $\lambda''(0)$ and therefore a strict subset of $\Spa{A}$.
Indeed, the second derivative of the equation $\lambda(\gamma(t))=0$ shows that $\lambda''(0)(\dot{\gamma}(0),\dot{\gamma}(0))=0$ for curves $\gamma\in C^1((-\epsilon,\epsilon),\Spa{U})$ with $\gamma(0)=0$.
Conversely we will show that an isotropic vector $v\in\Spa{A}$ of $\lambda''(0)$ belongs to the tangent cone, if there exits $w\in\Spa{A}$ with $\lambda''(0)(v,w)\not=0$.
We call such $v$ \emph{non-degenerate isotropic}.
Given any such $\upsilon\in\Spa{A}$ the function
\begin{align*}
\Lambda:(-\var,\var)\times(-\var,\var)&\to\mathbb{R}&(t,s)&\mapsto\begin{cases}2t^{-2}\lambda(tv+tsw)&\text{for }t\not=0\\\lambda''(0)(v+sw,v+sw)&\text{for }t=0\end{cases}
\end{align*}
is for small $\var>0$ smooth with $\Lambda(0,0)=0$ and $\frac{\partial\Lambda(0,0)}{\partial s}=2\lambda''(0)(v,w)\not=0$.
Due to the implicit function theorem there exists a smooth function $t\mapsto s(t)$ on a small interval $t\in(-\epsilon,\epsilon)$ with $s(0)=0$ and $\Lambda(t,s(t))=0$ and therefore also $\lambda(\gamma(t))=0$ for $\gamma(t)=tv+ts(t)w$.
The derivative $\dot{\gamma}(0)=v$ belongs to the tangent cone.
Moreover, note that since  $\tfrac{\partial \Lambda}{\partial s} \neq 0$ at $(0,0)$, this also holds at $\gamma(t)$ for small $t$.
By definition of $\Lambda$ this implies that $\lambda'(\gamma(t))\neq 0$ for small, non-zero $t$.
So the existence of a non-degenerate isotropic vector $v\in\Spa{A}$ would imply the last statement of the lemma.

The only property of $\lambda''$ thus far unused is its indefiniteness.
We now use this property to construct a set of non-degenerate isotropic vectors, which we will prove span $\Spa{A}$.
By indefiniteness there exist $x,y\in\Spa{A}$ with
\begin{align}\label{light cone coordinates}
\lambda''(0)(x,x)& =1 \;,& \lambda''(0)(y,y) & = -1 \;,& \lambda''(0)(x,y) = \lambda''(0)(y,x)&\in\mathbb{R}.
\end{align}
The replacement of $y$ by $(y-\lambda''(0)(x,y)x)/\sqrt{1+(\lambda''(0)(x,y))^2}$ preserves the second equation in~\eqref{light cone coordinates} and annihilates $\lambda''(0)(x,y)$.
The new vectors $x\pm y$ are non-degenerate isotropic.
By the previous step we know that they belong to the tangent cone, which therefore span $\mathbb{R}x+\mathbb{R}y$.
For any $u \in \Spa{A}$, we calculate
\begin{gather*}
\lambda''(0)(u+tx+sy,u+t'x+s'y)=tt'+a(t+t')+b(s+s')-ss'+c
\quad\text{with}\\
a=\lambda''(0)(u,x),\quad
b=\lambda''(0)(u,y),\quad
\text{and}\quad 
c= \lambda''(0)(u,u).
\end{gather*}
So long as $(t,s) \neq (-a,b)$ we can find $(t',s')$ such that this Hessian is non-zero.
Further we see that $u + tx + sy$ is isotropic if and only if $(t+a)^2-(s-b)^2=a^2-b^2-c$.
Since $\mathbb{R}^2\setminus\{(0,0)\}\to\mathbb{R},(\Tilde{t},\Tilde{s})\mapsto\Tilde{t}^2-\Tilde{s}^2$ is surjective, there exist $(t,s)\in\mathbb{R}^2$ such that $u + tx + sy$ is non-degenerate isotropic.
Due to the second step, $u$ itself belongs to the span of the tangent cone in addition to $tx+sy$.
\end{proof}

\noindent We now define locally a function $\lambda$ with three properties.
By the first property the level set is $\kodaira(\X,E)$.
Properties (i) and (ii) together imply that the tangent spaces of $\Spa{U}$ and $\kodaira(\X,E)$ coincide at $(\V,\upsilon,\phi)$.
Due to the foregoing lemma, property (iii) further ensures that the tangent cone of $\kodaira(\X,E)$ spans its tangent space.

\begin{lemma}\label{lem:isothermic1}
Let for a given $(\V,\upsilon,\phi)\in\kodaira(\X,E)$ on a compact Riemann surface $X$ the map~\eqref{eq:singularity1} have a real one-dimensional kernel.
Then there exist an open subset $O\subset\pot{E}\times H^0(\X,\sob{1,p}{E})^{\times2}$, a submanifold $\Spa{U} \ni (\V,\upsilon,\phi)$ of $O$ and a smooth function $\lambda: \Spa{U} \to \mathbb{R}$ with the following properties:
\begin{enumeratethm}
\item $\kodaira(\X,E) \cap O = \lambda^{-1}[\{0\}]$,
\item $\lambda'(\V,\upsilon,\phi)=0$,
\item $\lambda''(\V,\upsilon,\phi)$ is indefinite.
\end{enumeratethm}
\end{lemma}
\begin{proof}
Let $(\upsilon\pa,\phi\pa)$ span the real one-dimensional kernel of~\eqref{eq:singularity1}.
We choose
\begin{align}\label{dual vector}
(\alpha,\beta)\in H^0(\X,\forms{0,1}{}\ban{p}{E})^{\times2}
\quad\text{with}\quad
\int_\X\Real\left(\lh\upsilon\pa,\alpha\rh+\lh\phi\pa,\beta\rh\right)\ne0,
\end{align}
in the dual space of the kernel of~\eqref{eq:singularity1} with respect to the Hermitian pairing defined in Lemma~\ref{lem:11 pairing properties}(iii).
The quotient $H^0(\X,\forms{0,1}{}\ban{p}{E})^{\times2}/\mathbb{R}(\alpha,\beta)$ becomes a Banach space with the norm given by the infimum of the norms of all representatives.
We now show that in a neighborhood of the strongly isothermic triple $(\V,\upsilon,\phi)$ the zero set $\Spa{U}$ of the composition $\Tilde{\Theta}$ of $\Theta$~\eqref{immersion kodaira} with the natural projection onto the quotient $H^0(\X,\forms{0,1}{}\ban{p}{E})^{\times2}/\mathbb{R}(\alpha,\beta)$ is a submanifold of $\pot{E}\times H^0(\X,\sob{1,p}{E})^{\times2}$.
The proof is a modification of the proof of Theorem~\ref{thm:kod submanifold}.
We first modify the arguments which show that $\Theta'$ is surjective at triple that are not strongly isothermic.
The modification is thus: It suffices to show that the image of the linear map~\eqref{eq:deltaW-map} together with $(\alpha,\beta)$ spans the finite-dimensional cokernel of~\eqref{eq:W2-L2-map}.
Since the dual of a quotient space is the orthogonal complement of the space in the denominator, it is enough to show that the restriction of~\eqref{eq:eta-xi-pairing} to the orthogonal complement of $(\alpha,\beta)$ is injective.
This follows from~\eqref{dual vector}, since the intersection of the orthogonal complement of $(\alpha,\beta)$ in $\bigl(H^0(\X,\ban{p/(p-1)}{KE^{-1}})\bigr)^{\times2}$ with the kernel of~\eqref{eq:singularity1} is trivial.
Recall from the proof of Theorem~\ref{thm:kod submanifold} that $\Theta'(\V,\upsilon,\phi)$ is surjective if and only if the map~\eqref{eq:singularity1} is injective.
So in the present situation, $\Theta'(\V,\upsilon,\phi)$ is not surjective and the cokernel of $\Theta'(\V,\upsilon,\phi)$ is spanned by $(\alpha,\beta)$.
Therefore we have 
\begin{equation}\label{eq:ker-theta-tildetheta}
\ker(\Theta'(\V,\upsilon,\phi)) = (\Theta'(\V,\upsilon,\phi))^{-1}[\mathbb{R}(\alpha,\beta)] = \ker(\Tilde{\Theta}'(\V,\upsilon,\phi)) .
\end{equation}
Moreover $\ker(\Tilde{\Theta}'(\V,\upsilon,\phi))$ is together with $\ker {\Theta}'(\V,\upsilon,\phi)$ complemented as was shown in the proof of Theorem~\ref{thm:kod submanifold}.
Again Lemma~\ref{lem:implicit-variant} shows that $\Spa{U}$ is a submanifold of $\pot{E}\times H^0(\X,\sob{1,p}{E})^{\times2}$ locally at $(\V,\upsilon,\phi)$.

By definition of $\Spa{U}$ and due to~\eqref{dual vector} $\kodaira(\X,E)\cap\Spa{U}$ is the zero set of $\lambda:
 \Spa{U} \to \mathbb{R}$
\begin{align}
\label{eqn:kod lambda}
&\lambda:&
(\V',\upsilon',\phi') &\mapsto \int_\X \Real\left(\lh\upsilon\pa,(\delbar{E} - \V')\upsilon'\rh+ \lh\phi\pa,(\delbar{E} - \V')\phi'\rh \right).
\end{align}
Due to Equation~\eqref{eq:ker-theta-tildetheta}, the tangent space of the Banach manifold $\Spa{U}$ at $(\V,\upsilon,\phi)$ is equal to the tangent space of $\kodaira(\X,E)$ as defined in Definition~\ref{D:kodaira-tangentspace-tangentcone} at this point.
The derivative of $\lambda$ along a tangent vector $(\var \V, \var \upsilon, \var \phi)$ of $\Spa{U}$ at $(\V,\upsilon,\phi)$ is equal to the integral of the real part of
\[
- \lh \upsilon\pa,\var \V \upsilon \rh - \lh \phi\pa, \var \V \phi \rh + \lh \upsilon\pa, (\delbar{E}-\V)\var \upsilon \rh + \lh \phi\pa, (\delbar{E}-\V)\var \phi \rh  . 
\]
The last two terms vanish because $(\upsilon\pa,\phi\pa)$ belongs to the kernel of the adjoint operator of $\delbar{E}-\V$, and the sum of the first two terms vanishes because $(\upsilon\pa,\phi\pa)$ belongs to the kernel of~\eqref{eq:singularity1} and therefore also to the kernel of~\eqref{eq:eta-xi-pairing}.
Hence the first derivative of $\lambda$ at $(\V,\upsilon,\phi)$ vanishes.

At $(\V,\upsilon,\phi)$ the Hessian $\lambda''$ in the directions $(\var \V,\var\upsilon,\var\phi)$ and $(\var \V',\var\upsilon',\var\phi')$ in $T_{(\V,\upsilon,\phi)}\Spa{U}$ is equal to
\begin{gather*}
\labelthis{eqn:kod second derivative}
\lambda''(\V,\upsilon,\phi)\big((\var \V,\var\upsilon,\var\phi),(\var \V',\var\upsilon',\var\phi')\big) \\
\begin{aligned}
&=-\int_\X\hspace{-1mm}\Real\left(\lh\upsilon\pa,\var \V\var\upsilon'\rh + \lh\phi\pa,\var \V\var\phi'\rh + \lh\upsilon\pa,\var \V'\var\upsilon\rh + \lh\phi\pa,\var \V'\var\phi\rh\right) \\
&=-\langle\ls\upsilon\pa,\var\upsilon'\rs+\ls\phi\pa,\var\phi'\rs,\var \V\rangle-\langle\ls\upsilon\pa,\var\upsilon\rs+\ls\phi\pa,\var\phi\rs,\var \V'\rangle.
\end{aligned}
\end{gather*}
In Lemma~\ref{lem:indefinite} we show that this bilinear form is indefinite.
\end{proof}

Let us now explain how the blow-up will be performed.
Choose a local coordinate $z$ taking values in the open subset $\Omega\subset\mathbb{C}$.
By blow-up we mean the limit as $r\to 0^+$ of functions $\xi_r(z_\infty) := \xi(z_0+r z_\infty)$, where $\xi(z)$ is a function on $\Omega$.
We evaluate the limit of the family $\xi_r$ on the unit disc $z_\infty \in \mathbb{D}$, which corresponds for sufficiently small $r$ to $\Set{B}_r:=B(z_0,r)\subset\Omega$.
We restrict to a bounded subset of the exceptional fiber, because otherwise the constant functions would not be integrable.
Indeed, the blow-up of a continuous function is easily seen to be a constant function, so it is necessary to include them.

In our work generalizing quaternionic holomorphicity to Sobolev functions, we have mostly taken the point of view that local rather than pointwise properties are central, for example a root is defined by having $(z-z_0)^{-1}\xi$ belong to $\banach{q}$ rather than $\xi(z_0) = 0$.
But in our blow-up method we could not completely avoid the particular behavior of functions at the point $z_0$.
Control of the function at a point is the essence of the Lebesgue differentiation theorem~\cite[Chapter~1, \S5.7]{St}: for $\xi \in \banach{p}$ at almost every point $z_0$ the oscillation is small:
\[\lim_{r\to 0} \frac{1}{\mu(\Set{B}_r)}\int_{\Set{B}_r} |\xi(z) - \xi(z_0)|^p \;\dmu(z) = 0.\]
Such points are called Lebesgue or $\banach{p}$\=/points of $\xi$.
The `mass' that a function can have in a neighborhood of a Lebesgue point is restricted, since
\begin{align*}
r^{-2/p} \|\xi\|_{\banach{p}(\Set{B}_r)}
&\leq r^{-2/p} \|\xi - \xi(z_0)\|_{\banach{p}(\Set{B}_r)} + r^{-2/p} \|\xi(z_0)\|_{\banach{p}(\Set{B}_r)} \\
&\to 0 + \pi^{1/p}|\xi(z_0)|.
\end{align*}
In this paper we interpret the limit in the Lebesgue differentiation theorem as a type of blow-up and for each $r>0$ make the coordinate change $\xi_r(z_\infty) = \xi(z_0 + r z_\infty)$, so that $z \in \Set{B}_r$ is rescaled to $z_\infty \in \mathbb{D}$.
Under this rescaling $r^{-2/p}\|\xi\|_{\banach{p}(\Set{B}_r)}$ becomes $\|\xi_r\|_{\banach{p}(\mathbb{D})}$.
Thus $z_0$ is an $\banach{p}$\=/point of $\xi$ if and only if
\[\lim_{r\to 0} \|\xi_r - \xi(z_0)\|_{\banach{p}(\mathbb{D})} = 0.\]
We also use the results of~\cite{Zi}, which extend these ideas in the case of Sobolev functions.
The following lemma summarizes the tools that we need.
\begin{lemma}
\label{lem:sobolev local behavior}
Fix $1<p<2$ and $p^\ast=\frac{2-p}{2p}$.
Let $\xi\in\sobolev{1,p}(\Omega,\qat)$ and for any $z_0\in \Omega$ define $\xi_r(z_\infty) = \xi(z_0 + r z_\infty)$ and $P_r(z_\infty) = \xi(z_0) + rz_\infty \partial\xi(z_0) + r\overline{z_\infty} \barpartial\xi(z_0)$.
Almost all $z_0\in\Omega$ are simultaneously $\banach{p^\ast}$\=/points of $\xi$ and $\banach{p}$\=/points of $\partial\xi$ and $\barpartial \xi$ such that
\[\lim_{r\to 0} r^{-1}\| \xi_r - P_r \|_{\banach{p^\ast}(\mathbb{D})} = 0\]
holds.
If two $\xi$'s have these properties at a point $z_0$, so do their right-$\qat$\=/linear combinations.
Suppose further that $\xi$ is $\V$\=/holomorphic for $\V \in \banach{2}(\Omega,\qat^-)$.
Then such a point $z_0$ is a root of $\xi$ if and only if $\xi(z_0) = 0$.
If it has a root then
\[\lim_{r\to 0} \|r^{-1}\xi_r - z_\infty\partial\xi(z_0)\|_{\banach{p^\ast}(\mathbb{D})} = 0.\]
\end{lemma}

\begin{proof} 
In~\cite[Theorem~3.4.2]{Zi} a generalization of Lebesgue points for Sobolev functions is given.
The statement given is that $\xi \in \sobolev{1,p}$ is well approximated in almost everywhere by its first order Taylor polynomial:
\[\lim_{r \to 0} r^{-1-2/p}\| \xi - P \|_{\banach{p}(\Set{B}_r)} = 0\]
as $r \to 0$.
In our rescaled point of view this is $r^{-1}\| \xi_r - P_r \|_{\banach{p}(\mathbb{D})} \to 0$, which is slightly weaker than the first statement of the lemma.
Since finitely many conditions that hold almost everywhere hold simultaneously almost everywhere, we can additionally have that $z_0$ is a Lebesgue point for $\xi$ considered as an element of $\banach{p^\ast}(\Omega)$ and $\partial\xi, \barpartial \xi \in \banach{p}(\Omega)$.
This allows for an improved estimate.
We have that
\[r^{-1}\| \partial(\xi_r - P_r) \|_{\banach{p}(\mathbb{D})}=r^{-1}\| r(\partial\xi)_r - r\partial\xi(z_0) \|_{\banach{p}(\mathbb{D})}\to 0,\]
and likewise for $\barpartial(\xi_r - P_r)$.
Together we can say that
\[r^{-1}\| \xi_r - P_r \|_{\sobolev{1,p}(\mathbb{D})} \to 0.\]
Using the Sobolev inequality, we get
\[r^{-1}\| \xi_r - P_r \|_{\banach{p^\ast}(\mathbb{D})}\leq C\, r^{-1}\| \xi_r - P_r \|_{\sobolev{1,p}(\mathbb{D})}\to 0.\]
Because the rescaling operation $\xi_r$ and the Taylor approximation $P_r$ are both linear in the function $\xi$, an easy triangle inequality gives the required limits for any function in the right-$\qat$ span of two well-approximated functions.


Now we assume that $\xi$ is $\V$\=/holomorphic.
We recall Wente's inequality~\ref{wente}, which states that $|\xi|$ is continuous and that it is zero if and only if $\xi$ has a root in the sense of Definition~\ref{order of roots}.
Suppose that $\xi$ has a root at $z_0$.
The $\V$\=/holomorphicity implies that $\barpartial \xi(z_0) =\V(z_0)\xi(z_0) = 0$ and so the Taylor polynomial is just $P_r = rz_\infty \partial\xi(z_0)$.
Applying the first part of the lemma gives
\[
\lim_{r \to 0} \|r^{-1}\xi_r - z_\infty\partial\xi(z_0)\|_{\banach{p^\ast}(\mathbb{D})} = 0.
\qedhere
\]
\end{proof}
In the following lemma we construct an indefinite bilinear form $\lambda_\infty''$ for the plane, which will be shown to be a perturbation of $\lambda''$.
Despite the double prime notation, $\lambda_\infty''$ is defined as a bilinear form  and we do not care whether it is the Hessian of some function.
\begin{lemma}\label{second variation plane}
Let $\upsilon_\infty,\phi_\infty,\upsilon\pa_\infty,\phi\pa_\infty\in\qat\setminus\{0\}$ with $\upsilon_\infty\Bar{\upsilon}\pa_\infty=-\phi_\infty\Bar{\phi}\pa_\infty\in\qj\mathbb{C}$.
For $f\in C^\infty_0(\mathbb{D},\mathbb{C})$ the variation $d\Bar{z}_\infty\, \var \V_\infty=d\Bar{z}_\infty\,\barpartial^2 \qj f$ can be supplemented with variations with compact support $\var \upsilon_\infty$ and $\var \phi_\infty$ to a tangent vector $(dz_\infty\var \V_\infty,\var\upsilon_\infty,\var\phi_\infty)$ of $\kodaira(\mathbb{D},\unity)$ at the triple $(0,\upsilon_\infty,z_\infty\phi_\infty)$ corresponding to the flat plane $z_\infty\mapsto\upsilon_\infty^{-1}z_\infty\phi_\infty$ on $z_\infty\in\mathbb{D}$.
The corresponding~\eqref{eqn:kod second derivative} with $\upsilon=\upsilon_\infty$, $\phi=z_\infty\phi_\infty$, $\upsilon\pa=z_\infty\upsilon\pa_\infty$, $\phi\pa=\phi\pa_\infty$, which we denote $\lambda_\infty''$, 
is indefinite.
\end{lemma}
\begin{proof}
Recall $\kodaira(\mathbb{D},\unity)$ describes admissible maps $\mathbb{D}\to \qat$ which are quotients of two $\V$\=/holomorphic sections $\upsilon,\phi$ of the trivial $\qat$\=/line bundle $\unity_\qat$ on the unit disc $\mathbb{D}$.
This space contains the triple $(\V,\upsilon,\phi) = (0,\upsilon_\infty, z_\infty \phi_\infty)$, which corresponds to the admissible map $z_\infty\mapsto\upsilon_\infty^{-1}z_\infty\phi_\infty$ and describes a flat plane, as in Example~\ref{eg:plane kodaira}.
This plane has a dual plane in the sense of Definition~\ref{def:dual map} described in terms of Lemma~\ref{lem:strongly isothermic correspondence} by $(z_\infty\upsilon\pa_\infty,\phi\pa_\infty)$: Since $\upsilon_\infty\Bar{\upsilon}\pa_\infty=-\phi_\infty\Bar{\phi}\pa_\infty\in\qj\mathbb{C}$ they satisfy~\eqref{eq:dual map condition}
\[
-\phi\pa_\infty(\overline{z_\infty\phi_{\infty}}) 
= -\phi\pa_\infty\,\Bar{\phi}_{\infty} \Bar{z}_\infty 
= \upsilon\pa_\infty \, \Bar{\upsilon}_\infty\, \Bar{z}_\infty
= (z_\infty \,\upsilon\pa_\infty) \Bar{\upsilon}_\infty.
\]
We may change the coordinate $z_\infty$ such that these planes are described by constants with $\upsilon_\infty\Bar{\upsilon}\pa_\infty=-\phi_\infty\Bar{\phi}\pa_\infty\in\qk\mathbb{R}$.
For $(d\Bar{z}\,\var \V_\infty, \var \upsilon_\infty,\var\phi_\infty)$ to be tangential to $\kodaira(\mathbb{D},\unity)$ where $\var \V_\infty=\barpartial^2\qj f$ with $f\in C^\infty_0(\mathbb{D},\mathbb{C})$, the variations $\var\upsilon_\infty$ and $\var\phi_\infty$ have to obey (see~\eqref{eq:Theta'-kernel}):
\begin{align*}
\barpartial^2\qj f\upsilon_\infty&=\barpartial \var\upsilon_\infty&
\barpartial^2\qj fz_\infty\phi_\infty&=\barpartial \var\phi_\infty.
\end{align*}
The solutions are $\var\upsilon_\infty=\barpartial \qj f \upsilon_\infty$ and $\var\phi_\infty=(\barpartial \qj fz_\infty-\qj f)\phi_\infty$.
Note that $\var \V_\infty$, $\var \upsilon_\infty$ and $\var \phi_\infty$ are functions with compact support, in contrast to the constants $\upsilon_\infty$, $\phi_\infty$, $\upsilon\pa_\infty$ and $\phi\pa_\infty$.
From Equation~\eqref{eqn:kod second derivative} we thus obtain
\begin{align*}
&-\tfrac{1}{2}\lambda_\infty''((d\Bar{z}_\infty\,\var \V_\infty,\var\upsilon_\infty,\var\phi_\infty),(d\Bar{z}_\infty\,\var \V_\infty,\var\upsilon_\infty,\var\phi_\infty)) \\
&= \!\int_{\mathbb{D}}\!\Real\left(\lh z_\infty\upsilon\pa_\infty, d\Bar{z}_\infty\barpartial^2 (\qj f)\barpartial(\qj f)\upsilon_\infty\rh \!+\! \lh \phi\pa_\infty, d\Bar{z}_\infty\barpartial^2 (\qj f)(\barpartial (\qj f) z_\infty \!-\! \qj f)\phi_\infty\rh\right)  \\
&= 2 \int_{\mathbb{D}} \Real \bigl( \Bar{\upsilon}\pa_\infty \Bar{z}_\infty \qk  \barpartial^2 \Bar{f} {\partial} f \upsilon_\infty + \Bar{\phi}\pa_\infty \qk \barpartial^2  \Bar{f}(\partial f z_\infty  - f)\phi_\infty \bigr) \;\dmu  \\
&= 2 \int_{\mathbb{D}} \Real \bigl( ( \upsilon_\infty\Bar{\upsilon}\pa_\infty \Bar{z}_\infty + z_\infty \phi_\infty\Bar{\phi}\pa_\infty) \qk \barpartial^2 \Bar{f}\partial f - \phi_\infty\Bar{\phi}\pa_\infty \qk \barpartial^2 \Bar{f} f \bigr) \;\dmu  \\
&= -2\phi_\infty\Bar{\phi}\pa_\infty \qk \int_{\mathbb{D}} \Real \bigl( \barpartial^2 \Bar{f} f \bigr) \;\dmu.
\end{align*}
Finally, we should integrate this by parts and write real and imaginary parts $f = u + \qi w$
\begin{align*}
&= 4\phi_\infty\Bar{\phi}\pa_\infty \qk \int_{\mathbb{D}} \Real \bigl( \barpartial(\barpartial \Bar{f}\,f) - \barpartial \Bar{f}\, \barpartial f \bigr) \;\dmu
= - 4\phi_\infty\Bar{\phi}\pa_\infty \qk \int_{\mathbb{D}} \Real \bigl( \barpartial \Bar{f}\, \barpartial f \bigr) \;\dmu \\
&= -\phi_\infty\Bar{\phi}\pa_\infty \qk \int_{\mathbb{D}}(\tfrac{\partial u}{\partial x})^2 - (\tfrac{\partial u}{\partial y})^2 + (\tfrac{\partial w}{\partial x})^2 - (\tfrac{\partial w}{\partial y})^2 \dmu.
\end{align*}
If this is positive for $f(x,y)$ then it is negative for $f(y,x)$.
Thus $\lambda_\infty''$ is indefinite.
\end{proof}
\noindent Now we can combine Lemma~\ref{lem:sobolev local behavior} and Lemma~\ref{second variation plane}
with a perturbation argument to show that $\lambda''$ is indefinite, as promised above.
\begin{lemma}\label{lem:indefinite}
Let $(\V,\upsilon,\phi)\in\kodaira(\X,E)$ be strongly isothermic.
Then the second derivative $\lambda''$ of~\eqref{eqn:kod lambda} is indefinite at $(\V,\upsilon,\phi)$.
\end{lemma}
\begin{proof}  
As previously discussed, our strategy is to blow up $X$ at a point to a plane and apply a perturbation argument.
Let us explain why this is advantageous.
The subset of $\banach{2}$ functions that are Kodaira potentials has finite-codimension, so there are many possible variations.
But the Kodaira potential of the plane is $0$, and thus its variations are far easier to understand.
The steps of the proof are as follows.
The first phase of the proof constructs a suitable blow-up of $X$ to a plane.
This allows us to introduce the variations $(d\Bar{z}_\infty\,\var \V_\infty,\var\upsilon_\infty,\var\phi_\infty)$ from Lemma~\ref{second variation plane} and pull them back to variations $(\var \Tilde{\V}_r,\var\Tilde{\upsilon}_r,\var\Tilde{\phi}_r)$ on $X$.
These pulled-back variations are however not tangent to $\Spa{U}$ and so the second phase of the proof concerns perturbing these variations to variations $(\var {\V}_r,\var{\upsilon}_r,\var{\phi}_r)$ on $X$ that are tangent to $\Spa{U}$.
The final phase of this proof is to show that $(\var {\V}_r,\var{\upsilon}_r,\var{\phi}_r)$ converge suitably in the blow-up limit $r\downarrow 0$ to the $(\var \Tilde{\V}_r,\var\Tilde{\upsilon}_r,\var\Tilde{\phi}_r)$ and therefore $\lambda''(\var {\V}_r,\var{\upsilon}_r,\var{\phi}_r)$ converges to $\lambda''_\infty(d\Bar{z}_\infty\,\var \V_\infty,\var\upsilon_\infty,\var\phi_\infty) = \pm 1$.
This shows that $\lambda''$ is indefinite.

In our discussion of the space of Kodaira representations we have fixed a line bundle $E$, necessitating the use of potentials $\V$ with possibly non-zero $\V^+$.
However in this lemma we are dealing with a particular surface so we can choose $E$ such that $\V^+=0$.
This simplifying observation means that we can invoke Wente's inequality~\ref{wente} and the improved Darboux transformation from Chapter~\ref{chapter:kodaira}.
In particular, given a $\V$\=/holomorphic section $\upsilon$ and a point $z_0$ that is not a root of $\upsilon$ there exists a neighborhood of $z_0$ where $(\partial\upsilon)\upsilon^{-1}$ is square-integrable.

Fix $1<p<2$.
The Sobolev conjugate of $p$ is $p^\ast = \frac{2-p}{2p}$ and the choice $p' = \frac{p}{p-1}$ is suitable for Hölder's inequality since $1 = (p')^{-1} + 2^{-1} + (p^\ast)^{-1}$.
This is frequently also used in the form $2/p' + 2/p^\ast = 1$.
Let the kernel of~\eqref{eq:singularity1} be spanned by $(\upsilon\pa,\phi\pa)$.

Choose a chart $z$ and identify an open subset of $\X$ with an open set $\Omega \subset \mathbb{C}$.
 Let's find a suitable point $z_0$ of $\Omega$ to blowup.
 We choose $z_0$ such that certain properties hold.
 These properties each hold almost everywhere on $\Omega$, so they hold simultaneously almost everywhere.
Firstly, we require that $dF$, $\upsilon$, $\phi\pa$ have no roots at $z_0$ (in the sense of Definition~\ref{order of roots}).
Let $z_0$ be a point where Lemma~\ref{lem:sobolev local behavior} holds for the sections $\upsilon,\phi,\upsilon\pa$, and $\phi\pa$.
Additionally, we choose $z_0$ such that it is Lebesgue point of the potentials $\V$ and $\B+\U^-:= -(\partial\upsilon)\upsilon^{-1}$ (compare Theorem~\ref{thm:darboux}(1)(b)).
Moreover, due to Lemma~\ref{lem:change of order} we can control the derivative of $\phi$ using the $\U^-$\=/holomorphic function $\psi = (\partial +\B+\U)\phi$ from Theorem~\ref{thm:darboux}(1)(e).
Choose $z_0$ so that Lemma~\ref{lem:sobolev local behavior} also applies to $\psi$.
By our hypothesis that $dF = \lp(\overline{\qj \upsilon})^{-1},\psi\rp$ has no root at $z_0$, neither does $\psi$.
Later in the proof we will carry out a computation with respect to a basis of the image of $H^0(\X,\Q{K^{-1}E^{-1},\V\sd})^{\times 2}$ in $\pot{E}$ under~\eqref{eq:singularity1}, so let us choose a (finite) set of  elements that map to an orthonormal basis of the image.
Similar to $\psi$, we further require $z_0$ be Lebesgue points but not roots of the chosen basis elements of $H^0(\X,\Q{K^{-1}E^{-1},\V\sd})$.

If we consider $F$ as a map into $\mathbb{HP}^1$ by $[\upsilon:\phi]$, the quaternionic Möbius transformation~\eqref{eq:quaternionic moebius 1} of the target space
\begin{align*}
\upsilon&\mapsto\upsilon,&\phi&\mapsto \phi-\upsilon \, F(z_0),&\phi\pa&\mapsto\phi\pa,&\upsilon\pa&\mapsto\upsilon\pa+\phi\pa\, \overline{F(z_0)},&F&\mapsto F-F(z_0)
\end{align*}
corresponds to a translation of the surfaces $F$ and $-\Bar{F}$ (the dual surface is unchanged, but it is only determined up to a translation anyway).
These translations preserve the conditions that $F$ is strongly isothermic, the sections are holomorphic and that Lemma~\ref{lem:sobolev local behavior} applies.
This particular translation ensures that $\phi$, $\upsilon\pa$ and $F$ have a root at $z_0$, since $F(z_0) = \upsilon(z_0)^{-1}\phi(z_0) = - \overline{\phi\pa(z_0)^{-1}\upsilon\pa(z_0)}$ makes them zero at $z_0$ and we can use Lemma~\ref{lem:sobolev local behavior} (compare to the proof of Lemma~\ref{step 5}).
Notice that the function $\psi$ is unchanged by this translation because $(\partial+\B+\U^-)\upsilon = 0$.

Now we give the definitions of blow-ups of the Kodaira data of $F$ and of $(\upsilon\pa,\phi\pa)$, which will allow us to also control $\lambda''$ in the limit.
\index{Blow-up}
The idea, as explained in the paragraph before Lemma~\ref{lem:sobolev local behavior}, is that for each $r > 0$ we restrict functions to the disc $\Set{B}_r=B(z_0,r)\subset\Omega$ and rescale this disc to the unit disc $\mathbb{D}$ through the coordinate changes $z_\infty = \frac{z-z_0}{r}$.
By our choice of $z_0$, the functions $\upsilon_r(z_\infty) = \upsilon(z_0+rz_\infty)$ and $\phi\pa_r(z_\infty)=\phi\pa(z_0+rz_\infty)$ converge  in $\banach{p^\ast}(\mathbb{D},\qat)$ and $\banach{p'}(\mathbb{D},\qat)$ to non-zero constants $\upsilon_\infty = \upsilon(z_0) \in \qat$ and $\phi\pa_\infty = \phi\pa(z_0) \in \qat$, respectively.

Next note that $\phi$ is $\V$\=/holomorphic and has a root.
Due to Lemma~\ref{lem:sobolev local behavior} we know that $\phi_r(z_\infty) = r^{-1}\phi(z_0 + rz_\infty)$ converges in $\banach{p^\ast}(\mathbb{D},\qat)$ to $z_\infty\phi_\infty = z_\infty \partial\psi(z_0)$.
Be advised that this subscript notation does not match the notation in that lemma, but including the power $r^{-1}$ is the more convenient convention here.
The function $\psi_r(z_\infty) = \psi(z_0 + rz_\infty)$ also converges to $\phi_\infty$ in $\banach{p^\ast}$ since $\psi(z_0) = \partial\phi(z_0) + (B(z_0) + \U^-(z_0))\phi(z_0) = \phi_\infty + 0$.
As $\psi$ does not have a root, $\phi_\infty$ is non-zero.
Likewise $\upsilon\pa_r(z_\infty) = r^{-1}\upsilon\pa(z_0 + rz_\infty)$ converges in $\banach{p'}(\mathbb{D},\qat)$ to $z_\infty\upsilon\pa_\infty$.

That $(\upsilon\pa,\phi\pa)$ lies in the kernel of~\eqref{eq:singularity1} implies locally $\upsilon\Bar{\upsilon}\pa+\phi\Bar{\phi}\pa=0$.
Therefore we have 
\begin{equation}\label{sequence1}
\upsilon_\infty\Bar{\upsilon}\pa_\infty\Bar{z}_\infty
=\lim_{r\downarrow0}\upsilon_r\Bar{\upsilon}\pa_r
=-\lim_{r\downarrow0}\phi_r\Bar{\phi}\pa_r=-z_\infty\phi_\infty\Bar{\phi}\pa_\infty.
\end{equation}
We see that this blow-up of the Kodaira data not only defines an admissible map to the plane $z_\infty\mapsto\upsilon_\infty^{-1}z_\infty\phi_\infty$, but also induces a blow up of the dual surface and the data of Lemma~\ref{lem:strongly isothermic correspondence}.
Taking $z_\infty\in\{1,\qi\}$ shows that $\upsilon_\infty\Bar{\upsilon}\pa_\infty= - \phi_\infty\Bar{\phi}\pa_\infty \in \qj \mathbb{C}$.
This completes the blow-up of $X$ at $z_0$ to a plane.

Lemma~\ref{second variation plane} applies to planes and dual planes of this form.
In that lemma we constructed variations with compact support $(d\Bar{z}_\infty\,\var \V_\infty,\var\upsilon_\infty,\var\phi_\infty)$
in the tangent space of $\kodaira(\mathbb{D},\unity)$ at $(0,\upsilon_\infty,z_\infty\phi_\infty)$, parameterized by a function $f \in C^\infty_0(\mathbb{D},\mathbb{C})$.
Among these variations were two with
\[
\lambda_\infty''((d\Bar{z}_\infty\,\var \V_\infty,\var\upsilon_\infty,\var\phi_\infty),(d\Bar{z}_\infty\,\var \V_\infty,\var\upsilon_\infty,\var\phi_\infty))=\pm 1
\]
respectively.
From these variations we want to construct variations on $X$.
For small $r$ we define for $z \in \Set{B}_r$ the functions
\begin{align*}
\var \Tilde{\V}_r(z) 
&= d\Bar{z}\, r^{-2} (\barpartial^2\qj f)(\tfrac{z-z_0}{r}) \\\var \Tilde{\upsilon}_r(z) 
&= r^{-1} (\barpartial \qj f)(\tfrac{z-z_0}{r})\upsilon(z)\\\var \Tilde{\phi}_r(z) 
&= r^{-1}(\barpartial \qj f)(\tfrac{z-z_0}{r}) \phi(z) - \qj f(\tfrac{z-z_0}{r}) \psi(z)
\end{align*}
and extend them by zero on $\X\setminus \Set{B}_r$ to elements $(\var \Tilde{\V}_r,\var\Tilde{\upsilon}_r,\var\Tilde{\phi}_r)$ of the Banach space $\pot{E}\times H^0(\X,\sob{1,p}{E})^{\times2}$.

In general, these elements do not belong to the tangent space of the submanifold $\Spa{U}$ from Proposition~\ref{lem:isothermic1} at $(\V,\upsilon,\phi)$, but are rather useful intermediaries.
We will perturb them to obtain a family of variations $(\var {\V}_r,\var{\upsilon}_r,\var{\phi}_r)_{r\in(0,\varepsilon)}$ that are tangential to $\Spa{U}$ at $(\V,\upsilon,\phi)$.
Due to Lemma~\ref{lem:isothermic1}, this means that they are tangent to $\kodaira(\X,E)$.
The idea behind the definition of these intermediaries is to take the variations from Lemma~\ref{second variation plane} and replace $\upsilon_\infty, z_\infty\phi_\infty,\phi_\infty$ with the sections $\upsilon, \phi, \psi$; the later converge to the former at $z_0$ in a suitable sense.
The choice of leading factors of $r$ was not clear to us at the outset, but was motivated by working backwards through the argument with the following considerations.
First, the scaling of $\var \Tilde{\V}_r$ is determined by $\var {\V}_r$, since one is the perturbation of the other.
Second, the scaling of $\var {\V}_r$ determines the scaling of the variations $\var{\upsilon}_r,\var{\phi}_r$ by the fact that they must be tangent to $\Spa{U}$, c.f.\ Equation~\eqref{eq:Theta'-kernel}.
This fixes the differences between the powers of $r$ (one and two respectively).
Finally, although a `natural' scaling of $r^{-1}$ for $\var\Tilde{\V}_r$ would ensure that its $\pot{E}$\=/norm is the same as the $\banach{2}(\mathbb{D},\qat)$\=/norm of $\var \V_\infty$, its norm must in fact diverge in the limit $r\downarrow0$ like $r^{-1}$, otherwise the second derivative $\lambda''$ in this direction tends to zero.
Indeed, consider $\lambda''(\V,\upsilon,\phi)\big((\var \Tilde{\V}_r,\var\Tilde{\upsilon}_r,\var\Tilde{\phi}_r),(\var \Tilde{\V}_r,\var\Tilde{\upsilon}_r,\var\Tilde{\phi}_r)\big)$, namely
\begin{align*}
-2\int_{\X}\Real\left(\lh\upsilon\pa,\var \Tilde{\V}_r \var\Tilde{\upsilon}_r\rh+\lh\phi\pa,\var \Tilde{\V}_r \var\Tilde{\phi}_r\rh\right).
\end{align*}
The definition of the $(1,1)$\=/bracket is given in Definition~\ref{def:11 pairing} as a $(1,1)$\=/form and under the change of variables $z = z_0 + rz_\infty$ we have that 
$dz \wedge d\Bar{z} = r^2 dz_\infty \wedge d\Bar{z}_\infty$
Taking each term in turn,
\begin{align*}
\int_{\X} \lh\upsilon\pa,\var \Tilde{\V}_r \var\Tilde{\upsilon}_r\rh
&=\int_{\Set{B}_r}\Bar{\upsilon}\pa (-\qj) dz \wedge d\Bar{z} r^{-2} \var \V_\infty(\tfrac{z-z_0}{r}) r^{-1} \var\upsilon_\infty(\tfrac{z-z_0}{r})\upsilon_\infty^{-1}\upsilon \;\dmu \\
&= \int_{\mathbb{D}} \overline{\upsilon\pa_r} (-\qj) dz_\infty \wedge d\Bar{z}_\infty \var \V_\infty \var\upsilon_\infty \upsilon_\infty^{-1}\upsilon_r \; \dmu_\infty.
\end{align*}
The middle factors are bounded smooth functions and so in particular they are $\banach{2}(\mathbb{D},\qat)$.
Overall therefore the integral converges to $\int_{\mathbb{D}} \lh z_\infty \upsilon\pa_\infty,d\Bar{z}_\infty\var \V_\infty \var\upsilon_\infty\rh$.
Likewise for the second term
\begin{align*}
\int_{\X} \lh\phi\pa,\var \Tilde{\V}_r \var\Tilde{\phi}_r\rh
&= \int_{\mathbb{D}} \overline{\phi\pa_r} (-\qj) dz_\infty \wedge d\Bar{z}_\infty \var \V_\infty [ (\barpartial \qj f) \phi_r - \qj f \psi_r ]
\end{align*}
we have convergence to $\int_{\mathbb{D}} \lh\phi\pa_\infty,d\Bar{z}_\infty\var \V_\infty \var\phi_\infty\rh$.
Together these prove the convergence $\lambda''$ in direction of the intermediaries $(\var \Tilde{\V}_r,\var\Tilde{\upsilon}_r,\var\Tilde{\phi}_r)$ to $\lambda''_\infty$ in direction of $(d\Bar{z}_\infty\var \V_\infty,\var\upsilon_\infty,\var\phi_\infty)$ in the blow-up limit $r\downarrow0$.

While we are still thinking about the blow-up, let us gather some facts on which we will later rely.
 First observe that the intermediaries typically do not converge in the limit $r\downarrow0$, but
\begin{align*}
r^{2/p'}\|\var \Tilde{\upsilon}_r\|_{\banach{p^*}(X)}
&= r^{2/p' - 1 + 2/p^*} \| \var \upsilon_\infty \upsilon_\infty^{-1}\upsilon_r \|_{\banach{p^*}(\mathbb{D})} 
\leq \| \var \upsilon_\infty \|_\infty |\upsilon_\infty|^{-1} \|\upsilon_r\|_{\banach{p^*}(\mathbb{D})}
\end{align*}
is bounded, as is
\begin{align*}
r\|\var \Tilde{\V}_r\|_{2}
&= r^{-1} \|\var \V_\infty(\tfrac{z-z_0}{r})\|_{\banach2(\Set{B}_r)}
\leq r^{-1} \|\var \V_\infty\|_\infty (\pi r^2)^{1/2}
\end{align*}
and
\[r^{-2/p^\ast}\|\var \Tilde{\phi}_r\|_{\banach{p^\ast}(X)}\leq \|\barpartial \qj f\|_\infty \|\phi_r\|_{\banach{p^\ast}(\mathbb{D})} + \|\qj f\|_\infty \|\psi_r\|_{\banach{p^\ast}(\mathbb{D})}.\]
It follows from this and the fact that $\upsilon\pa$ has a root at $z_0$ that the products 
\begin{align*}
&\|\upsilon\pa\|_{\banach{p'}(\Set{B}_r)} \|\var\Tilde{\V}_r\|_{2}
= r^{2/p'} \times \|\upsilon\pa_r\|_{\banach{p'}(\mathbb{D})} \times r\|\var\Tilde{\V}_r\|_{2}, \\
&\|\upsilon\pa\|_{\banach{p'}(\Set{B}_r)} \|\var \Tilde{\upsilon}_r\|_{\banach{p^\ast}(X)}
= r \times \|\upsilon\pa_r\|_{\banach{p'}(\mathbb{D})} \times r^{2/p'} \|\var \Tilde{\upsilon}_r\|_{\banach{p^\ast}(X)}
\end{align*}
both tend to zero with $r\downarrow0$.
Another common situation is an $\banach{2}$\=/potential multiplied with $\var\Tilde{\upsilon}_r$.
For example,
\begin{align}
\label{eq:W varupsilon Lp limit}
\|\V \var\Tilde{\upsilon}_r\|_{\banach{p}(X,\qat)}
&\leq r^{2/p^\ast} \times r^{-1}\|\V\|_{\banach{2}(B(0,r),\qat)} \times r^{2/p'}\|\var\Tilde{\upsilon}_r\|_{\banach{p^\ast}(B(0,r),\qat)}
\end{align}
converges to zero because $z_0$ is a Lebesgue point of $\V$, which as discussed before Lemma~\ref{lem:sobolev local behavior} means that $r^{-1}\|\V\|_{\banach{2}(B(0,r),\qat)}$ is bounded.

We now begin the perturbation phase of our argument.
We will first use a suitable projection to ensure that $\var \V_r$ is mapped by $-\tfrac{\partial \Theta}{\partial \V}$ into the image of the Fredholm operator $(\delbar{E}-\V) \times (\delbar{E}-\V)$.
This is a necessary condition on $\var \V_r$ in order for $(\var \upsilon_r, \var \phi_r)$ to exist.
Second, we will then choose $(\var \upsilon_r, \var \phi_r) \in \bigl( (\delbar{E}-\V)^{-1}  \times (\delbar{E}-\V)^{-1} \bigr) (-\tfrac{\partial \Theta}{\partial \V}(\var \V_r))$.
Equation~\eqref{eq:Theta'-kernel} then implies that $(\var \V_r, \var \upsilon_r, \var \phi_r)$ is tangential to $\Spa{U}$.
Any pair $(\Tilde{\upsilon}\pa,\Tilde{\phi}\pa)\in H^0(\X,\Q{K^{-1}E^{-1},\V\sd})^{\times 2}$ gives an element of $\pot{E}$, which can be paired with $\var\Tilde{\V}$:
\begin{align*}
\langle\ls\Tilde{\upsilon}\pa,\upsilon\rs+\ls\Tilde{\phi}\pa,\phi\rs,\var\Tilde{\V}\rangle
=\int_\X \Real \left(\lh\Tilde{\upsilon}\pa,\var{\Tilde{\V}} \upsilon\rh+\lh\Tilde{\phi}\pa,\var{\Tilde{\V}} \phi\rh\right).
\end{align*}
Let $\Spa{C}$ be the set of these elements in $\pot{E}$.
Choose $\var \V_r$ to be the image of $\var \Tilde{\V}_r$ under the orthogonal projection onto the orthogonal complement of $\Spa{C}$.
Recall that $(\Tilde{\upsilon}\pa,\Tilde{\phi}\pa)\in H^0(\X,\Q{K^{-1}E^{-1},\V\sd})^{\times 2}$ is the dual of the cokernel of the operator~\eqref{eq:W2-L2-map}, see the argument below Equation~\eqref{eq:deltaW-map}.
Hence $(\var \V_r\upsilon, \var \V_r\phi)$ belongs to the image of the operator~\eqref{eq:W2-L2-map}.
Next we show that this perturbation $\var \V_r - \var\Tilde{\V}_r \in C$ is in some sense small; namely that $r^{-1}\|\var \V_r - \var\Tilde{\V}_r\|_{2}$ is bounded for small $r$.
For given $(\Tilde{\upsilon}\pa,\Tilde{\phi}\pa)\in H^0(\X,\Q{K^{-1}E^{-1},\V\sd})^{\times 2}$ let $\Tilde{\upsilon}\pa_r(z_\infty) = \Tilde{\upsilon}\pa(z_0+rz_\infty)$ and $\Tilde{\phi}\pa_r(z_\infty) = \Tilde{\phi}\pa(z_0+rz_\infty)$.
In particular, for the basis chosen at the outset of this proof we know that that $\Tilde{\upsilon}\pa_r$ and  $\Tilde{\phi}\pa_r$ converge to nonzero constants $\Tilde{\upsilon}\pa_\infty$ and $\Tilde{\phi}\pa_\infty$ in $\banach{p'}(\mathbb{D})$.
From the product rule for paired holomorphic sections Lemma~\ref{lem:11 pairing properties}(iv)
\begin{align*}
d\lp\Tilde{\upsilon}\pa,\var\Tilde{\upsilon}_r\rp 
&=-0+\lh \Tilde{\upsilon}\pa, (\delbar{E} - \V)\var\Tilde{\upsilon}_r \rh \\
&=\lh\Tilde{\upsilon}\pa, \var\Tilde{\V}_r \upsilon - r^{-1}(\barpartial \qj f)(\tfrac{z-z_0}{r})(\B+\U^-)\upsilon - \V r^{-1}(\barpartial \qj f)(\tfrac{z-z_0}{r})\upsilon \rh
\end{align*}
Using this formula to integrate by parts gives
\begin{align*}
r^{-1}&\left|\int_{\X} \lh\Tilde{\upsilon}\pa,\var{\Tilde{\V}_r} \upsilon\rh \right|
\leq 2r^{-2} \|\barpartial \qj f\|_\infty \int_{\Set{B}_r} |\Tilde{\upsilon}\pa| (|\B + \U^-|+|\V|) |\upsilon| \;\dmu\\
&\leq 2\|\barpartial \qj f\|_\infty  r^{-1}(\|\B + \U^-\|_{\banach{2}}+\|\V\|_{\banach{2}}) r^{-2/p'}\|\Tilde{\upsilon}\pa\|_{\banach{p'}} r^{-2/p^\ast}\|\upsilon\|_{\banach{p^\ast}}
\end{align*}
(the omitted domains of the Lebesgue spaces in this equation are all $\Set{B}_r$).
The right hand side is bounded because $z_0$ is a Lebesgue point of $\B$, $\U$, $\V$, $\Tilde{\upsilon}\pa$, and $\upsilon$.
For the other term we do not even have to integrate by parts due to the root of $\phi$:
\begin{align*}
r^{-1}\int_{\X}\lh\Tilde{\phi}\pa,\var{\Tilde{\V}_r}\phi\rh
&= r^{-1}\int_{\mathbb{D}}\lh\Tilde{\phi}\pa_r,d\Bar{z}_\infty r^{-2}\var \V_\infty r\phi_r\rh \,r^2 \\
&\to \int_{\mathbb{D}}\lh\Tilde{\phi}\pa_\infty, d\Bar{z}_\infty\var \V_\infty z_\infty\phi_\infty\rh.
\end{align*}
This implies that for any $\varepsilon > 0$
\begin{align}
\label{eq:tildeWr convergence}
\lim_{r\downarrow0} r^{-1+\varepsilon} \|\var \V_r - \var \Tilde{\V}_r\|_{2}
&=0.
\end{align}
It remains to construct the variations $\var\upsilon_r$ and $\var\phi_r$.
First however we must fix an inverse of $\delbar{E} - \V$.
Choose a complementary space $\Spa{A}\subset H^0(\X,\sob{1,p}{E})$ of the finite-dimensional subspace $H^0(\X,\Q{E,\V})$, and a complementary subspace $\Spa{B}\subset H^0(\X,\ban{p}{KE})$ of the image of the Fredholm operator $\delbar{E}-\V: H^0(\X,\sob{1,p}{E})\to H^0(\X,\ban{p}{KE})$.
Then the canonical map $\Spa{A}\to H^0(\X,\sob{1,p}{E})/H^0(\X,\Q{E,\V})$ is an isomorphism and the restriction of $H^0(\X,\ban{p}{KE})\to H^0(\X,\ban{p}{KE})/\Spa{B}$ to the image of $\delbar{E}-\V$ is also an isomorphism.
The quotient $H^0(\X,\ban{p}{KE})/\Spa{B}$ becomes a Banach space with the usual norm, namely the infimum of the norm on an equivalence class.
With these choices, $\delbar{E}-\V$ induces an isomorphism of Banach spaces $D: \Spa{A}\to H^0(\X,\ban{p}{KE})/\Spa{B}$.

For $(\var \V_r, \var \upsilon_r, \var \phi_r)$ to belong to the tangent space of $\Spa{U}$ means that $(\delbar{E}-\V)\var \upsilon_r = \var \V_r \upsilon$ and $(\delbar{E}-\V)\var \phi_r = \var \V_r \phi$ by Equation~\eqref{eq:Theta'-kernel}.
By our construction of $\var \V_r$, we know that $(\var \V_r \upsilon, \var \V_r \phi)$ does belong to the image of $(\delbar{E} - \V)\times(\delbar{E} - \V)$.
Hence we define
\begin{align*}
\var\upsilon_r := D^{-1}(\var \V_r \upsilon), \qquad
\var\phi_r := D^{-1}(\var \V_r \phi).
\end{align*}
However we want not only that $\var \upsilon_r$ and $\var \phi_r$ exist, but that they are perturbations of $\var \Tilde{\upsilon}_r$ and $\var \Tilde{\phi}_r$.
Observe 
\begin{align*}
(\delbar{E} - \V)\left( \var\upsilon_r - \var\Tilde{\upsilon}_r \right)
&= \var \V_r \upsilon - d\Bar{z}\,r^{-2} (\barpartial^2\qj f)(\tfrac{z-z_0}{r})\upsilon - r^{-1}(\barpartial \qj f)(\tfrac{z-z_0}{r})dz\partial\upsilon + \V \var\Tilde{\upsilon}_r \\
&= (\var \V_r - \var \Tilde{\V}_r) \upsilon + \var\Tilde{\upsilon}_r \upsilon^{-1}(\B+\U^-)\upsilon + \V \var\Tilde{\upsilon}_r.
\end{align*}
The first term converges to zero in $\banach{p}(\X,\qat)$ since $\var \V_r - \var \Tilde{\V}_r$ converges to zero in $\pot{E}$ from~\eqref{eq:tildeWr convergence} and $\upsilon$ is $\banach{p^\ast}(\X,\qat)$.
The following terms are all similar to one another, since in magnitude they are essentially $\var\Tilde{\upsilon}_r$ multiplied with an $\banach{2}$\=/potential.
But this situation was examined in~\eqref{eq:W varupsilon Lp limit} and found to converge to zero.
Together this shows that $D( \var\upsilon_r - \var\Tilde{\upsilon}_r )$ tends to zero in $\banach{p}(\X,\qat)$ and hence that $\var\upsilon_r - \var\Tilde{\upsilon}_r$ converges to $0$ in $\sobolev{1,p}(X)$.
By Sobolev embedding we also have convergence to zero in $\banach{p^\ast}(\X,\qat)$.

For the other tangent vector we prove that the perturbation is even smaller, which will be used in the final step.
Consider
\begin{align*}
r^{-2/p^\ast}&(\delbar{E}-\V)(\var\phi_r - \var\Tilde{\phi}_r)
= r^{-2/p^\ast}(\var \V_r - \var \Tilde{\V}_r) \phi \\
&- r^{-2/p^\ast}(\barpartial \qj f)(\tfrac{z-z_0}{r})r^{-1}dz(\partial \phi - \psi) 
- r^{-2/p^\ast}\qj f(\tfrac{z-z_0}{r})dz\partial \psi 
+ r^{-2/p^\ast}\V \var\Tilde{\phi}_r
\end{align*}
in the $\banach{p}(\X)$\=/norm.
The first term is bounded by $r^{-2/p^\ast}\|\var \V_r - \var \Tilde{\V}_r\|_{2} \|\phi\|_{\banach{p^\ast}(\X)}$.
Due to~\eqref{eq:tildeWr convergence}, the difference of variations of the potential can absorb a factor of $r^{-2/p^\ast}$ while still tending to zero.
For the second term, we recall that by definition $\psi = (\partial + \B + \U^-)\phi$ and bound it by
\[
\|\barpartial \qj f\|_{\banach{\infty}(\mathbb{D})} \times r^{-1}\|\B+\U^-\|_{\banach{2}(\Set{B}_r)} \times r^{-2/p^\ast}\|\phi\|_{\banach{p^\ast}(\Set{B}_r)}.\]
Because $r^{-1}\|\B+\U^-\|_{\banach{2}(\Set{B}_r)}$ is bounded, the convergence of $r^{-2/p^\ast}\|\phi\|_{\banach{p^\ast}(\Set{B}_r)}$ to zero is sufficient to see that the second term also vanishes as $r\downarrow0$.
For the third term we have the estimate $r^{-2/p^\ast} \|\qj f\|_\infty \|\partial \psi\|_{\banach{p}(\Set{B}_r)}$.
If we write the scaling factor as $r^{-2/p^\ast} = r^{1-2/p}$ and use the fact that $\partial\psi$ has an $\banach{p}$\=/point at $z_0$, we see that this term is vanishing.
Finally, the fourth term is handled by the observation that $r^{-2/p^\ast}\|\var \Tilde{\phi}_r\|_{\banach{p^\ast}(X)}$ is bounded and $\|\V\|_{\banach{2}(\Set{B}_r)}$ is vanishing.
Altogether this shows that
\begin{align}\label{eq:varphi difference strong}
\lim_{r\to 0} r^{-2/p^\ast}\|\var\phi_r - \var\Tilde{\phi}_r\|_{\banach{p^\ast}(\X)} = 0.
\end{align}
To conclude, we have constructed $(\var \V_r, \var \upsilon_r, \var \phi_r)$ in the tangent space of $\Spa{U}$ such that they converge to $(\var \Tilde{\V}_r, \var \Tilde{\upsilon}_r, \var \Tilde{\phi}_r)$ in $\pot{E}\times H^0(\X,\sob{1,p}{E})^{\times2}$.

Finally we must show that the second derivative $\lambda''$ at $(\V,\upsilon,\phi)$ in the direction $(\var \V_r, \var \upsilon_r, \var \phi_r)$ converges to $\lambda''$ in the direction $(\var \Tilde{\V}_r, \var \Tilde{\upsilon}_r, \var \Tilde{\phi}_r)$.
We bound the difference by twice
\begin{align}
\label{eq:lambda'' difference}
\int_{\X} \left| \lh\upsilon\pa,\var \V_r \var\upsilon_r - \var \Tilde{\V}_r \var\Tilde{\upsilon}_r\rh \right| + \int_{\X} \left| \lh\phi\pa,\var \V_r \var\phi_r - \var \Tilde{\V}_r \var\Tilde{\phi}_r\rh \right|
\end{align}
The first integral is bounded by
\begin{align*}
2&\|\upsilon\pa\|_{\banach{p'}(\X)} \|\var \V_r - \var\Tilde{\V}_r\|_{\banach{2}(\X)} \|\var\upsilon_r - \var\Tilde{\upsilon}_r \|_{\banach{p^\ast}(\X)} \\
&+ 2\|\upsilon\pa\|_{\banach{p'}(\Set{B}_r)} \|\var \V_r - \var\Tilde{\V}_r\|_{\banach{2}(\X)} \|\var\Tilde{\upsilon}_r \|_{\banach{p^\ast}(\Set{B}_r)} \\
&+ 2\|\upsilon\pa\|_{\banach{p'}(\Set{B}_r)} \|\var\Tilde{\V}_r\|_{\banach{2}(\Set{B}_r)} \|\var\upsilon_r - \var\Tilde{\upsilon}_r \|_{\banach{p^\ast}(\X)}.
\end{align*}
In the first term, $\|\upsilon\pa\|_{\banach{p'}(\X)}$ is finite and the other factors tend to zero as $r\downarrow0$.
In the second term, we already know that both $\|\upsilon\pa\|_{\banach{p'}(\Set{B}_r)} \|\var \Tilde{\upsilon}_r\|_{\banach{p^\ast}(X)}$ and $\|\var \V_r - \var\Tilde{\V}_r\|_{\banach{2}(\X)}$ tend to zero.
For the third term, likewise $\|\upsilon\pa\|_{\banach{p'}(\Set{B}_r)} \|\var\Tilde{\V}_r\|_{\banach{2}(\Set{B}_r)}$ and $\|\var\upsilon_r - \var\Tilde{\upsilon}_r \|_{\banach{p^\ast}(\X)}$ tend to zero.

We split the second integral of~\eqref{eq:lambda'' difference} into an estimate of the same shape, but with $\phi\pa$ in the place of $\upsilon\pa$ and $\phi$ in the place of $\upsilon$.
The first term of the estimate tends to zero for the same reasons as before.
For the other two terms, the important difference is that $\phi\pa$ does not have a root at $z_0$ as $\upsilon\pa$ does, so we need different arguments.
For the second term $2\|\phi\pa\|_{\banach{p'}(\Set{B}_r)} \|\var \V_r - \var\Tilde{\V}_r\|_{\banach{2}(\X)} \|\var\Tilde{\phi}_r \|_{\banach{p^\ast}(\Set{B}_r)}$ we note that both the difference of the variations of the potentials and $\|\var\Tilde{\phi}_r \|_{\banach{p^\ast}(\Set{B}_r)}$ are tending to zero.
So the only difficulty is in the third term.
Observe
\begin{align*}
2\|\phi\pa\|_{\banach{p'}(\Set{B}_r)} & \|\var\Tilde{\V}_r\|_{\banach{2}(\X)} \|\var\phi_r - \var\Tilde{\phi}_r \|_{\banach{p^\ast}(\X)} \\
&= 2r^{-2/p'}\|\phi\pa\|_{\banach{p'}(\Set{B}_r)} \times r\|\var\Tilde{\V}_r\|_{\banach{2}(\X)} \times r^{-2/p^\ast}\|\var\phi_r - \var\Tilde{\phi}_r \|_{\banach{p^\ast}(\X)}.
\end{align*}
With the given scaling factors, $r^{-2/p'}\|\phi\pa\|_{\banach{p'}(\Set{B}_r)}$ and $r\|\var\Tilde{\V}_r\|_{\banach{2}(\X)}$ are bounded in $r$, while we have seen in~\eqref{eq:varphi difference strong} that $r^{-2/p^\ast}\|\var\phi_r - \var\Tilde{\phi}_r \|_{\banach{p^\ast}(\X)}$ tends to zero.
Therefore the value of second derivative in the direction of these tangent vectors approaches the value of the second derivative in the direction of the intermediaries, which itself approaches the value of $\lambda''_\infty$, namely $1$ or $-1$ depending on the choice of $f$.
This proves that $\lambda''(\V,\upsilon,\phi)$ is indefinite.
\end{proof}


\begin{proof}[Proof of Theorem~\ref{thm:kodaira isothermic characterization}]
We first prove the equivalence of~(i)-(iv).
The equivalence of (i) and (ii) is simply Lemma~\ref{lem:strongly isothermic correspondence}.
Lemmas~\ref{lem:tangent cone spans} and~\ref{lem:isothermic1} together show that (ii) implies (iii).
Statement~(iv) is a direct consequence of~(iii).
Finally, the contrapositive of Theorem~\ref{thm:kod submanifold} is that (iv) implies (i).

It remains to show that the tangent cone spans the tangent space and there is a smooth path in $\kodaira(\X,E)$ that meets the set of strongly isothermic surfaces in at most one point.
At not strongly isothermic $(\V,\upsilon,\phi)$ the equivalence of~(i) and~(iii) implies the equality of the tangent cone and the tangent space of $\kodaira(\X,E)$.
Theorem~\ref{thm:kod submanifold} implies that the not strongly isothermic $(\V,\upsilon,\phi)$ constitute an open subset of $\kodaira(\X,E)$, and so there certainly exists a path that does not intersect the set of strongly isothermic surfaces.
At strongly isothermic $(\V,\upsilon,\phi)$ the combination of Lemma~\ref{lem:tangent cone spans} and Lemma~\ref{lem:isothermic1} shows that the tangent cone spans the tangent space.
For the path $\gamma$, take the path given by Lemma~\ref{lem:tangent cone spans}.
Because $\lambda'(\gamma(t))$ is surjective for $t \neq 0$ it follows that $\kodaira(\X,E)$ is a submanifold at $\gamma(t)$ and therefore these are not strongly isothermic surfaces.
\end{proof}

\chapter{Constrained Willmore Kodaira Triples}
\label{chapter:constrained 1}
In this chapter we describe critical points of the Willmore functional on the space of all admissible maps from a compact Riemann surface $\X$ to $\mathbb{R}^4$.
Such surfaces are called constrained Willmore surfaces~\cite{Pinkall1987}.
The term `constrained' here refers to the fact that we have fixed $\X$ as a Riemann surface and not just a smooth surface and that $F : \X \to \qat$ is a conformal map; we are constrained to a particular conformal class.
The main result is Theorem~\ref{constrained Willmore 1}, which establishes a quaternionic analysis criterion to be constrained Willmore and proves these surfaces are analytic and branched conformal.
This theorem is followed by a natural generalization, which invites the definition of $n$\=/Willmore surfaces.

The Willmore functional of an admissible map was introduced in Chapter~\ref{chapter:isothermic}.
In particular, recall the Equation~\eqref{eq:willmore energy kodaira} $\willmore(F) = 4\pi\deg(E') + 4\|\V'\|_2^2$ for the Willmore energy in terms of the Kodaira data $(E',\V',\upsilon',\phi')$ of $F$ with $\V'\in \pot{E'}^-$ (see Theorem~\ref{thm:kodaira}).
But in the present context, as discussed following Definition~\ref{def:kodaira}, the potential $\V$ of $(\V,\upsilon,\phi) \in \kodaira(\X,E)$ is not necessarily anti-commuting and so this formula may not apply.
However we can first use Corollary~\ref{gauge holomorphic structure} to transform $(\V,\upsilon,\phi) \in \kodaira(\X,E)$ into $(\V',\upsilon',\phi')\in\kodaira(\X,E')$ with $|\V'| = |\V^-|$.
Now it follows from Equation~\eqref{eq:willmore energy kodaira} that 
\begin{equation}\label{eq:willmore-F}
\index{Willmore energy}
\willmore(F) = 4\pi\deg(E) + 4\|\V^-\|_2^2 
\end{equation}
for the admissible map $F=\upsilon^{-1}\phi$ for any $(\V,\upsilon,\phi) \in \kodaira(\X,E)$.
Now that we have computed $\willmore(F)$ in terms of a triple of $\kodaira(\X,E)$, we investigate the critical points of the Willmore functional as a map on $\kodaira(\X,E)$.
The following lemma was explained to the third author by Ernst Kuwert~\cite{KS}.
\begin{lemma}\label{lem:dWillmore-space-cone}
The derivative of $(\V,\upsilon,\phi)\mapsto\willmore(\upsilon^{-1}\phi)$ vanishes on the tangent space of $\kodaira(\X,E)$ at some $(\V,\upsilon,\phi)$ if and only if it vanishes on the corresponding tangent cone.
\end{lemma}
\begin{proof}                   %
Let $(\V,\upsilon,\phi) \in \kodaira(\X,E)$ be a representative of $F=\upsilon^{-1}\phi$.
By Theorem~\ref{thm:isothermic local characterization} there are two cases: Either (i) $\V=0$ holds, or (ii) the kernel of~\eqref{eq:singularity1} has real dimension $\le1$.
These cases are not disjoint, for example the round sphere has $\V = 0$ but is not strongly isothermic (the kernel of~\eqref{eq:singularity1} is trivial).

In case (i) Equation~\eqref{eq:willmore-F} shows that the derivative of the Willmore functional vanishes on the tangent space of  $\pot{E}\times H^0(\X,\sob{1,p}{E})^{\times 2}$ at $(\V,\upsilon,\phi)$, and therefore in particular both on the tangent space and the tangent cone of $\kodaira(\X,E)$ at $(\V,\upsilon,\phi)$.

In case (ii) the tangent cone spans the tangent space by Theorem~\ref{thm:kodaira isothermic characterization}.
\end{proof}
\begin{definition}
\index{Constrained Willmore}
A Kodaira triple $(\V,\upsilon,\phi)\in\kodaira(\X,E)$ on a compact Riemann surface $\X$ is called \emph{constrained Willmore} if and only if either of the two equivalent conditions of Lemma~\ref{lem:dWillmore-space-cone} holds.
\end{definition}

There exists an Euler-Lagrange formulation for Willmore surfaces and a modification with Lagrange multipliers for constrained Willmore surfaces (compare \cite[Theorem~14 and 15]{BPP} and references therein).
There is a difficulty with using this formulation for weakly conformal maps that is exemplified by the following example: The catenoid is clearly a Willmore surface, and because Willmore surfaces are invariant with respect to inversions~\cite[Theorem~13.14]{Pinkall2024}, the inverted catenoid is also Willmore.
In Example~\ref{eg:inverted catenoid} we compactified the inverted catenoid and saw that the resulting surface was admissible.
However it is not Willmore, as we will show in Example~\ref{eg:inverted catenoid willmore} below.
This means that there exist compact weakly conformal surfaces which are Willmore only after restriction to the complement of a finite subset.
For the Euler-Lagrange formulation, this means that the Euler-Lagrange equation is satisfied only on the complement of a set of measure zero.
But in Sobolev spaces, solutions of differential equations always coincide with their restrictions to complements of sets of measure zero.
In Theorem~\ref{constrained Willmore 1} we treat constrained Willmore surfaces from the point of view of quaternionic analysis.
A straightforward application of this theorem solves the described problem.
Indeed, as we will see in Example~\ref{eg:inverted catenoid willmore}, the difference between the inverted catenoid and its compactification is reflected by sections $\upsilon\pa$, $\phi\pa$ which are holomorphic on the complement of the two points added during the compactification, but which have poles at those points.
\begin{theorem}
\label{constrained Willmore 1}
On a compact Riemann surface a triple $(\V,\upsilon,\phi)\in\kodaira(\X,E)$ is constrained Willmore if and only there exists $(\upsilon\pa,\phi\pa)\in(H^0(\X,\Q{KE^{-1},\V\sd}))^{\times2}$ which are mapped by~\eqref{eq:singularity1} onto the part $\V^-$ of $\V$ in $\pot{E}^-$.
The corresponding admissible maps are analytic and branched conformal.
\end{theorem}
\begin{proof}
Any triple $(\V,\upsilon,\phi)\in\kodaira(\X,E)$ with $\V^-=0$ is constrained Willmore and the condition in the first statement holds with $\upsilon\pa=\phi\pa=0$.
Otherwise, a triple $(\var \V,\var\upsilon,\var\phi)$ belongs to the tangent space of $\kodaira(\X,E)$ at $(\V,\upsilon,\phi)$, if and only if $-\tfrac{\partial \Theta(\V,\upsilon,\phi)}{\partial \V}$ (see~\eqref{eq:deltaW-map}) maps $\var \V$ into the image of $(\var\upsilon,\var\phi)$ with respect to the Fredholm operator~\eqref{eq:W2-L2-map}.
Hence for a variation $\var \V\in\pot{E}$ there exist $(\var\upsilon,\var\phi)$ such that the triple $(\var \V, \var \upsilon, \var \phi)$ belongs to the tangent space of $\kodaira(\X,E)$ at $(\V,\upsilon,\phi)$ if and only if $-\tfrac{\partial \Theta(\V,\upsilon,\phi)}{\partial \V}$ maps $\var \V$ onto the image of the Fredholm operator~\eqref{eq:W2-L2-map}, or equivalently onto the orthogonal complement of the kernel of the dual operator.
By definition of the bilinear map $\ls\cdot,\cdot\rs$, $-\tfrac{\partial \Theta(\V,\upsilon,\phi)}{\partial \V}$ maps $\var \V$ into the orthogonal complement of this kernel if and only if $\var \V$ is orthogonal to the image of the map~\eqref{eq:singularity1}, whose domain is exactly this kernel.
Due to Equation~\eqref{eq:willmore-F} the variation of the Willmore functional is equal to
\begin{align*}
\var\willmore&=8\langle\var \V^- , \V^-\rangle =8\langle\var\V,\V^-\rangle 
\quad\text{with}\quad
\var \V = \var \V^+ + \var \V^-. 
\end{align*}
Hence $F$ is constrained Willmore if and only if $\V^-$ belongs to the image of~\eqref{eq:singularity1}, which shows the first statement of the theorem.

It remains to show that $F$ is analytic, since it would then follow that $F$ is branched conformal. 
By Corollary~\ref{gauge holomorphic structure} we can again transform the triple $(\V,\upsilon,\phi) \in \kodaira{(\X,E)}$ into another triple $(\V'^-,\upsilon',\phi') \in \kodaira{(\X,E')}$ with $F=\upsilon^{-1}\phi=\upsilon'^{-1}\phi'$. 
Therefore we may assume without loss of generality that the commuting part of the potential vanishes, and $\V=\V^-=\ls\upsilon\pa,\upsilon\rs + \ls\phi\pa,\phi\rs$ holds. 
We separate out the proof of analyticity to Lemma~\ref{lem:solutions are analytic}, so that we can prove it in a slightly more general situation.
\end{proof}

From our point of view the image of the map~\eqref{eq:singularity1} describes the constraints on $\var\V$ to be part of a tangent vector along $\kodaira(\X,E)$ at $(\V,\upsilon,\phi)$. 
The existence of the elements $\upsilon\pa,\phi\pa\in H^0(\X,\Q{KE^{-1},\V\sd})$ with $\V^-=\ls\upsilon\pa,\upsilon\rs+\ls\phi\pa,\phi\rs\in\pot{E}^-$ is equivalent to the condition that the derivative of the Willmore functional vanishes on the subset of $\pot{E}\times H^0(\X,\sob{1,p}{E})^{\times2}$ whose elements represent admissible maps. 
In the study of variational analysis with constraints such elements are called Lagrange multipliers. 
For this reason we call $(\upsilon\pa,\phi\pa)$ the Lagrange multipliers of the constrained Willmore Kodaira triple $(\V,\upsilon,\phi)$.

The Möbius transformation of a constrained Willmore surface is again constrained Willmore.
This follows from the observation in Chapter~\ref{chapter:isothermic} that the Willmore energy is invariant under Möbius transformation.
However, we can also be explicit in terms of the Lagrange multipliers.
Let tilde denote the data corresponding to the transformed surface.
For the transformed surface to be constrained Willmore we must have
\[
\V^-
= \Tilde{\V}^-
= \ls \Tilde{\upsilon}\pa, \Tilde{\upsilon} \rs + \ls \Tilde{\phi}\pa, \Tilde{\phi} \rs,
\]
but using Lemma~\ref{lem:01 pairing properties}(ii) 
\begin{align*}
\ls \Tilde{\upsilon}\pa, \Tilde{\upsilon} \rs + \ls \Tilde{\phi}\pa, \Tilde{\phi} \rs
&= \ls \Tilde{\upsilon}\pa, \phi\gamma + \upsilon\delta \rs + \ls \Tilde{\phi}\pa, \phi\alpha + \upsilon\beta \rs \\
&= \ls \Tilde{\phi}\pa\Bar{\beta} + \Tilde{\upsilon}\pa\Bar{\delta}, \upsilon \rs + \ls \Tilde{\phi}\pa\Bar{\alpha} + \Tilde{\upsilon}\pa\Bar{\gamma}, \phi \rs.
\end{align*}
Compare to the transformation of the data of dual surfaces following Lemma~\ref{lem:isothermic local duals}.
Thus the new Lagrange multipliers are the solution to the linear system
\begin{equation}
\label{eq:Lagrange mutlipliers mobius}
\Tilde{\phi}\pa\Bar{\alpha} + \Tilde{\upsilon}\pa\Bar{\gamma} = \phi\pa,
\qquad
\Tilde{\phi}\pa\Bar{\beta} + \Tilde{\upsilon}\pa\Bar{\delta} = \upsilon\pa.
\end{equation}

\begin{example}[Catenoid]
\label{eg:inverted catenoid willmore}
\index{Catenoid!Inverted catenoid}
\index{Constrained Willmore!Inverted catenoid}
We now give the promised example of the catenoid and its compactified inversion.
We have the Kodaira data of the catenoid from Example~\ref{eg:catenoid kodaira}.
The Lagrange multipliers are $\upsilon\pa = -\frac{1}{2}z^{-1}\qj(z + \qk)(1 + |z|^2)^{-1}$ and $\phi\pa = 0$.
Both are $\V\sd$\=/holomorphic, since 
\begin{gather*}
\barpartial \upsilon\pa
= -z^{-1}\qj \partial \left[ (1 + |z|^2)^{-1}(z + \qk) \right]
= -z^{-1}\qj (1 -\Bar{z} \qk ) (1 + |z|^2)^{-2} \\
(\barpartial \upsilon\pa)(\upsilon\pa)^{-1}
= z^{-1}\qj (1 -\Bar{z} \qk ) (z + \qk)^{-1} (-\qj) z (1 + |z|^2)^{-1}
= \qk z \Bar{z}^{-1} (1 + |z|^2)^{-1} 
= \V\sd.
\end{gather*}
We also check that $\V$ really is in the image of~\eqref{eq:singularity1}:
\begin{align*}
\ls \upsilon\pa,\upsilon \rs + \ls \phi\pa,\phi \rs
&= -2\qk \tfrac{1}{2} z^{-1}\qj(z + \qk)(1 + |z|^2)^{-1} (\overline{1 - z^{-1}\qk})^{-1} \qj - 0 \\
&= 2\qk z^{-1} \Bar{z} (1 + |z|^2)^{-1} 
= \V.
\end{align*}
For the inverted catenoid, we solve Equations~\eqref{eq:Lagrange mutlipliers mobius} to get
\[
\Tilde{\upsilon}\pa = \tfrac{1}{2}\upsilon\pa \qi
, \qquad
\Tilde{\phi}\pa = \tfrac{1}{2}\upsilon\pa.
\]
This shows that the (non-compactified) inverted catenoid is also constrained Willmore.

To show the compactified inverted catenoid is not constrained Willmore, it is not enough to show that the above Lagrange multipliers do not work; we must consider all possible multipliers.
Any other solution $(\upsilon\pa,\phi\pa)$ mapping under~\eqref{eq:singularity1} to $\V^-$ must differ from the above solution by an element of the kernel of~\eqref{eq:singularity1}.
For the catenoid, this was computed in Example~\ref{eg:catenoid isothermic kernel}.
This also transforms to the inverted catenoid using Equations~\eqref{eq:Lagrange mutlipliers mobius}, giving the general solution
\begin{align*}
\Tilde{\upsilon}\pa 
&= \tfrac{1}{2}\upsilon\pa \qi + \tfrac{1}{2}z^{-2}(\upsilon + \phi\qi)\alpha 
&
\Tilde{\phi}\pa 
&= \tfrac{1}{2}\upsilon\pa + \tfrac{1}{2}z^{-2}(\upsilon\qi + \phi)\alpha \\
&= \tfrac{1}{2}\upsilon\pa \qi + \tfrac{1}{2}z^{-2}\Tilde{\phi}\alpha 
&
&= \tfrac{1}{2}\upsilon\pa + \tfrac{1}{2}z^{-2}\Tilde{\upsilon}\alpha
\end{align*}
for any $\alpha \in \qat$.

Now we come to the point: there is no choice of $\alpha$ such that these sections extend to holomorphic sections on the compactified inverted catenoid.
Hence the compactified inverted catenoid is not constrained Willmore.
In Example~\ref{eg:inverted catenoid} we showed that $\Tilde{\upsilon}$ and $\Tilde{\phi}$ extend continuously and are non-zero at $z=0$.
On the other hand,
\[
\lim_{z \to 0} z \upsilon\pa
= -\tfrac{1}{2}\qj(0 + \qk)(1 + 0)^{-1}
= -\tfrac{1}{2}\qi.
\]
Therefore if $\alpha \neq 0$ then $\Tilde{\upsilon}\pa$ has a double pole at $z= 0$, but if $\alpha = 0$ then it still has a simple pole.
Any extension of the Lagrange multipliers to are meromorphic sections.
\end{example}

We close this chapter with the completion of Lemma~\ref{constrained Willmore 1}, giving analyticity if the potential has a special form.
Just as we generalized the Darboux transformation to handle more than two sections, so too can we generalize to potentials that are a sum of more than two $\ls\cdot,\cdot\rs$\=/pairings.
The $d=1$ case will be considered in Example~\ref{example:d=1}.

\begin{lemma}
\label{lem:solutions are analytic}
On a (not necessarily compact) Riemann $\X$ surface let the following data be given: $E$ a holomorphic $\mathbb{C}$\=/line bundle, $\V\in\pot{E}^-$, $\phi_1,\ldots,\phi_d\in H^0(\X,\Q{E,\V})$ and $\phi_1\pa,\ldots,\phi_d\pa\in H^0(\X,\Q{KE^{-1},\V\sd})$ such that $\V$ is equal to $\V=\sum_{m=1}^d\ls\phi_m\pa,\phi_m\rs$. Then $\phi_1,\ldots,\phi_d$ and $\phi_1\pa,\ldots,\phi_d\pa$ are analytic.
\end{lemma}
\begin{proof}
We first prove that $\phi_1,\ldots,\phi_d$ and $\phi_1\pa,\ldots,\phi_d\pa$ are smooth. 
Since the claim is local we may restrict to open subsets $\Omega$ of $\X$ (and omit subscripts) on which $E$ and $KE^{-1}$ are trivial and the operator $\Op{I}_{\Omega,\V}$ of~\eqref{eq:resolvent} exists. 
For all $m\in\mathbb{N}_+$ and $1<p<2$ the right-inverse  $\Op{I}_{\Omega} = \Op{I}_{\Omega,0}$ of $\barpartial$ is a bounded operator from $\sobolev{m-1,p}(\Omega,\qat)$ onto $\sobolev{m,p}(\Omega,\qat)$ (compare~\cite[Chapter~V]{St}).
Therefore the local representatives of elements in the kernel of $\barpartial-\V$ and $\barpartial-\V\sd$ are contained in $\bigcap_{p<\infty}\sobolev{m,p}\loc(\Omega,\qat)$, if $\V$ belongs to $\bigcap_{p<\infty}\sobolev{m-1,p}\loc(\Omega,\qat)$. By the definition~\eqref{def:01 pairing} of the pairing $\ls\,\cdot\,,\,\cdot\,\rs$, the equation $\V=\sum_{m=1}^d\ls\phi_m\pa,\phi_m\rs$ represents the local representatives of $\V$ as a sum of products of $\qk$ with the local representatives of the $\phi\pa$'s and the local representatives of the $\Bar{\phi}$'s. It now follows inductively from the Sobolev product formula that the local functions of the $\phi$'s and the $\phi\pa$'s belong to $\bigcap_{p<\infty}\sobolev{m,p}\loc(\Omega,\qat)$ for every $m \geq 1$ and are therefore smooth.

We now apply the theorem~\cite[Theorem~C]{Mo} to show analyticity of the sections $\phi_1,\ldots,\phi_d$ and $\phi_1\pa,\ldots,\phi_d\pa$.
They obey the equations
\begin{align*}
\left(\barpartial-\sum\nolimits_{l=1}^d\ls\phi_l\pa,\phi_l\rs\right)\phi_1&=0,&\ldots&&\left(\barpartial-\sum\nolimits_{l=1}^d\ls\phi_l\pa,\phi_l\rs\right)\phi_d&=0.\\
\left(\barpartial-\sum\nolimits_{l=1}^d\ls\phi_l,\phi_l\pa\rs\right)\phi_1\pa&=0,&\ldots&&\left(\barpartial-\sum\nolimits_{l=1}^d\ls\phi_l,\phi_l\pa\rs\right)\phi_d\pa&=0.
\end{align*}
This is a non-linear analytic elliptic system (compare~\cite[p.\ 532]{DoNi}).
Indeed, the leading order terms of these equations are of the form $\barpartial(\phi_1,\dots,\phi_d,\phi_1\pa,\dots,\phi_d\pa)=0$.
Because $\partial$ and $\barpartial$ are elliptic differential operators, the ellipticity follows.
So by~\cite[Theorem~C]{Mo} any solution $(\phi_1,\dots,\phi_d,\phi_1\pa,\dots,\phi_d\pa)$ is analytic.
\end{proof}

\chapter{Strongly Isothermic Weierstraß Triples}
\label{chapter:isothermic weierstrass triples}
In this chapter we describe the set of admissible maps in terms of the Weierstraß representation.
More precisely we introduce a space of Weierstraß triples and show that this space is locally a Banach manifold if and only if the corresponding admissible map is not strongly isothermic.
This result is analogous to the Theorem~\ref{thm:kod submanifold}, which treats the corresponding Kodaira triples.
One may ask why it is important to investigate the strongly isothermic Weierstraß triples, when we already have an understanding of strongly isothermic Kodaira triples.
The profit of the investigation of isothermic Weierstraß data will be an understanding of the singularities of the space of admissible maps into the imaginary quaternions.
As a consequence almost all results of this and the following chapter have a Part~(A) that deals with the general case and a Part~(B) that deals with admissible maps into the imaginary quaternions.

We did investigate whether the results from strongly isothermic Kodaira triples could be carried over by way of a diffeomorphism.
This was possible for potentials in $\banach{q}$ for $q > 2$ but for general admissible maps the difficulties we encountered make us open to the possibility that $\kodaira(\X,E)$ and $\wei(\X,E)$ may endow the space of admissible maps with different differentiable structures.

As in Chapter~\ref{chapter:isothermic kodaira triples}, we consider non-necessarily anti-commuting potentials $\U$ in the Weierstraß representation~\ref{thm:weierstrass}, in order to describe the admissible maps with fixed degree $\deg(E')$ as a subset of a single Banach space.
The important difference to the Kodaira representation is the need to consider the periods of the $1$\=/form $\lp \chi,\psi \rp$.
This leads us to consider in Definition~\ref{def:wei} a space $\prewei(\X,E)$ of holomorphic sections and then the subset $\wei(\X,E)$ on which the $1$\=/form is exact (and thus the Weierstraß representation can be integrated to an admissible map on $\X$).
Definitions and Lemmas~\ref{def:period extension}--\ref{lem:wei-necessary condition} provide the tools to handle periods.
With these tools we are able to follow the path set in Chapter~\ref{chapter:isothermic kodaira triples}.
Where possible, we have tried to formulate the results of this chapter with similar wording to the corresponding statements about $\kodaira(\X,E)$.

\begin{definition}\label{def:wei}
\index{Weierstraß representation}
Let $E$ be a holomorphic line bundle on a Riemann surface $\X$.
The space of triples $(\U,\chi,\psi)\in\pot{KE}\times H^0(\X,\sob{1,p}{E^{-1}})\times H^0(\X,\sob{1,p}{KE})$ that obey 
\begin{gather}\label{eq:con 1}
\chi\in H^0(\X,\Q{E^{-1},\U\sd})
\quad\text{without roots and}\quad
\psi\in H^0(\X,\Q{KE,\U})
\end{gather}
is denoted by $\prewei(\X,E)$.
Recall~\eqref{eq:def period map}.
We define $\wei(\X,E)$ as the kernel of the period map
\begin{align}\label{eq:periods}
\index{Period map}
\prewei(\X,E)&\to\Hom(H_1(\X,\mathbb{Z}),\qat),&(\U, \chi,\psi)&\mapsto \per\lp\chi,\psi\rp.
\end{align}
\end{definition}
So let us now introduce the analog of $\wei(\X,E)$ whose elements correspond to admissible maps into $\Imag\qat$.
We characterized in Chapter~\ref{chapter:3-space} when an element of $\wei(\X,E)$ describes an admissible map into $\Imag\qat$ in terms of a holomorphic section $g$ of $KE^2$.
The existence of such a section is therefore a necessary condition on $E$.
By fixing $g$ we consider only triples with fixed root divisor of $dF$.
\begin{definition}\label{def:wei 2}
For any holomorphic $\mathbb{C}$\=/line bundle $E$ on the Riemann surface $\X$ and a holomorphic section $g\ne0$ of $KE^2$ let $\prewei(\X,E,g)$ be the subspace of $(\U,\chi)\in\pot{KE}\times H^0(\X,\sob{1,p}{E^{-1}})$ such that $(\U,\chi,g\chi)\in\prewei(\X,E)$.
Moreover, let $\wei(\X,E,g)$ be the kernel of the following map, which is induced by the period map~\eqref{eq:periods}
\begin{align}
\label{eq:periods imagsubspace}
\prewei(\X,E,g)&\to\Hom(H_1(\X,\mathbb{Z}),\Imag\qat),&(\U, \chi)\mapsto \per\lp\chi,g\chi\rp.
\end{align}
\end{definition}
If we assume in the sequel the existence of an element of $\prewei(\X,E,g)$ or of $\wei(\X,E,g)$ we also implicitly assume the existence of a non-trivial holomorphic section $g$ of $KE^2$, which is a condition on $E$.

Note that since $\U\sd$ is a potential for $\chi$ the equation $(\delbar{KE} -\U)(g\chi)=0$ is equivalent to $\U=g\U\sd g^{-1}$ (compare with the discussion preceding Lemma~\ref{lem:admissible imaginary}).
Moreover, if we decompose $\U=\U^++\U^-$, then the involution $\U\mapsto g\U\sd g^{-1}$ acts as $g(\U^++\U^-)\sd g^{-1}=-\U^++\frac{g}{\Bar{g}}(\U^-)\sd$.
We denote the invariant potentials by
\begin{gather}
\pot{KE,g} := \{\U\in\pot{KE}^-\mid\U=\tfrac{g}{\Bar{g}}\U\sd\}
\label{def:real potential}
\end{gather}
(compare Equation~\eqref{eq:imaginary 1}).
In particular, for a holomorphic section $g\ne 0$ of $KE^2$ we have $\prewei(\X,E,g)=\{(\U,\chi)\in\pot{KE,g}\times H^0(\X,\sob{1,p}{E^{-1}})\mid(\delbar{E^{-1}} -\U\sd)\chi=0\}$.
Finally, due to Lemma~\ref{lem:10 pairing properties}(i) and (iii) $\lp\chi,g\chi\rp=\lp g\chi,\chi\rp=-\overline{\lp\chi,g\chi\rp}$, the periods of this form indeed belong to $\Imag\qat$.

\begin{example}[Catenoid]
\label{eg:catenoid belongs to wei}
\index{Catenoid!Weierstraß representation}
As we have seen in Example~\ref{eg:catenoid weierstrass}, the Weierstraß data for the catenoid $F: \X = \mathbb{C}^\ast \to \qat$ is $E = \X \times \mathbb{C}$, $U \equiv 0$, and $\chi = \psi = 1+z^{-1}\qk$.
We see that $g\equiv 1$.
By definition of the Weierstraß representation, $dF = \lp \chi, \psi \rp$, the $1$\=/form is exact.
Hence $(0, \chi, \psi)$ belongs to $\wei(\mathbb{C}^\ast, E)$ and $(0,\chi)$ belongs to $\wei(\mathbb{C}^\ast, E, 1)$.
\end{example}

In our convention the holomorphic section $\chi$ of a Weierstraß triple is assumed to have no roots.
This could be always achieved by choosing the holomorphic $\mathbb{C}$\=/line bundle $E$ properly.
For this reason the roots of the derivative of an admissible map into $\Imag\qat$ are determined by $g$.
For given $g$ each of the spinors $\chi$ and $\psi$ determines the other.
Therefore $\wei(\X,E,g)$ is the fixed point set of an involution acting on a subspace of $\wei(\X,E)$ whose spinors have the appropriate roots:
\begin{lemma}
\label{lem:real wei fixed point}
Let $E$ be a holomorphic $\mathbb{C}$\=/line bundle on the Riemann surface $\X$ and $g\ne0$ a holomorphic section of $KE^2$ with root divisor $D$.
Then $\{(\U,\chi,g\chi)\mid(\U,\chi)\in\prewei(\X,E,g)\}$ and $\{(\U,\chi,g\chi)\mid(\U,\chi)\in\wei(\X,E,g)\}$ are respectively the fixed point sets of the involution 
\begin{align}\label{involution for imaginary maps}
(\U,\chi,\psi)\mapsto(g\U\sd g^{-1},g^{-1}\psi,g\chi).
\end{align}
on the subspaces of $\prewei(\X,E)$ and $\wei(\X,E)$ consisting of the elements $(\U,\chi,\psi)$ such that the root divisor of $\psi$ is $D$.
\end{lemma}
\begin{proof}
By assumption on the root divisor of $\psi$, the transformed section $g^{-1}\psi$ is indeed a global section of $\Q{(KE^2)^{-1}KE,g^{-1}\U g}=\Q{E^{-1},(g\U\sd g^{-1})\sd}$ without roots.
And $g\chi$ is indeed a global section of $\Q{KE^2E^{-1},g\U\sd g^{-1}}=\Q{KE,g\U\sd g^{-1}}$ with root divisor $D$.
It remains to show the invariance of the kernel of~\eqref{eq:periods imagsubspace}.
This follows by Lemma~\ref{lem:10 pairing properties}(i) and (iii) since $\lp g^{-1}\psi,g\chi\rp=\lp\psi,\chi\rp=-\overline{\lp\chi,\psi\rp}$.
\end{proof}
We shall now transfer Theorem~\ref{thm:kod submanifold} and Theorem~\ref{thm:kodaira isothermic characterization} to the description of admissible maps by the Weierstraß representation.
In the course of our research for this book we attempted to show that the Darboux transformation which transforms the Kodaira triples into Weierstraß triples induces a $C^1$\=/diffeomorphism between $\kodaira(\X,E)$ and $\wei(\X,E)$.
However, for general $\banach{2}\loc$\=/potentials we did not succeed.
We present here a separate investigation of strongly isothermic Weierstraß triples.
Since a Weierstraß triple has to obey two different conditions~\eqref{eq:con 1} and~\eqref{eq:periods} we have to prove that neither condition gives rise to singularities away from the strongly isothermic maps.
The first condition~\eqref{eq:con 1} can be treated by the same methods as in the proof of Theorem~\ref{thm:kod submanifold}.
For the second condition~\eqref{eq:periods} we need new ideas.
By Lemma~\ref{lem:fredholm}, for paired holomorphic sections~\eqref{eq:con 1} the quaternionic-valued $1$\=/form $\lp \chi, \psi \rp$ is closed.
This form is the derivative $dF$ of an admissible map $F:\X\to\qat$ if and only if $(\chi,\psi)$ belongs to the kernel of~\eqref{eq:periods}.
On a compact Riemann surface of genus $\genus$ this gives $2\genus$ quaternionic conditions.
We shall investigate the derivatives of these conditions.
For this purpose we need to introduce the tangent space of $\wei(\X,E)$, because Definition~\eqref{D:kodaira-tangentspace-tangentcone} does not yet apply to $\wei(\X,E)$.
But it does apply to $\prewei(\X,E)$ and $\prewei(\X,E,g)$ with the following maps
\begin{gather}
\begin{aligned}\pot{KE}\!\times\!H^0(\X,\sob{1,p}{E^{-1}})\!\times\!H^0(\X,\sob{1,p}{KE})&\to H^0(\X,\forms{0,1}{}\!\ban{p}{E^{-1}})\!\times\!H^0(\X,\forms{0,1}{}\!\ban{p}{KE}),\hspace{-2mm}\\
(\U,\chi,\psi)&\mapsto\big((\delbar{E^{-1}} -\U\sd)\chi,(\delbar{KE}-\U)\psi\big)\;.
\end{aligned}\label{eq:prewei map}\\
\begin{aligned}\pot{KE,g}\times H^0(\X,\sob{1,p}{E^{-1}})&\to H^0(\X,\forms{0,1}{}\ban{p}{E^{-1}}),&(\U,\chi)&\mapsto(\delbar{E^{-1}} -\U\sd)\chi.
\end{aligned}\label{eq:prewei map 2}
\end{gather}
The period map~\eqref{eq:periods} on the other hand is only defined on $\prewei(\X,E)$ in Definition~\ref{def:wei}, because the form $\lp \chi, \psi \rp$ needs to be closed in order to define a map on homology.
Analogously for $\wei(\X,E,g)$ the map~\eqref{eq:periods imagsubspace} is only defined on $\prewei(\X,E,g)$.

\begin{definition}
\label{def:period extension}
Fix a canonical basis $A_1,\ldots,A_g,B_1,\ldots,B_g$ of $H_1(\X,\mathbb{Z})$ and further for cycle $\gamma$ of this basis fix a smooth $1$\=/form $\eta_{\gamma}$ as in~\eqref{eq:def line integral}.
The extension of~\eqref{eq:periods} is defined to be 
\begin{gather}\begin{aligned}
\pot{KE}\times H^0(\X,\sob{1,p}{E^{-1}})\times H^0(\X,\sob{1,p}{KE})
&\to\Hom(H_1(\X,\mathbb{Z}),\qat)\\
(\U,\chi,\psi)&\mapsto\left(\gamma\mapsto\int_\X\eta_\gamma\wedge \lp\chi,\psi\rp \right)
\label{eq:periods extension map}
\end{aligned}\end{gather}
We define the derivative of~\eqref{eq:periods} as the restriction of the derivative of~\eqref{eq:periods extension map} to the tangent space of $\prewei(\X,E)$.
The cartesian product of~\eqref{eq:prewei map} with~\eqref{eq:periods extension map} gives a smooth map from $\pot{KE}\times H^0(\X,\sob{1,p}{E^{-1}})\times H^0(\X,\sob{1,p}{KE})$ to $H^0(\X,\forms{0,1}{}\ban{p}{E^{-1}})\times H^0(\X,\forms{0,1}{}\ban{p}{KE})\times\Hom(H_1(\X,\mathbb{Z}),\qat)$ whose kernel is $\wei(\X,E)$.
Analogously, the extension of~\eqref{eq:periods imagsubspace} is~\eqref{def:period extension} composed with $(\U,\chi) \mapsto (\U,\chi,g\chi)$ and $\wei(\X,E,g)$ is the kernel of the cartesian product of~\eqref{eq:prewei map 2} with this extension.
However, especially in situations where $\prewei(\X,E)$ and $\prewei(\X,E,g)$ are not manifolds, it will be useful to extend the period map to the full space $\pot{KE}\times H^0(\X,\sob{1,p}{E^{-1}})\times H^0(\X,\sob{1,p}{KE})$ and $\pot{KE}\times H^0(\X,\sob{1,p}{E^{-1}})$, respectively.
\end{definition}

Note that this extension depends on the particular choices of cycles and representing forms $\eta_\gamma$.
In the case that this is evaluated at a point of $\prewei(\X,E)$ then the form $\lp\chi,\psi\rp$ is closed and property~\eqref{eq:def line integral} makes the period map independent of the choices.
Hence~\eqref{eq:periods extension map} is indeed an extension of~\eqref{eq:periods}.
Now Definition~\ref{D:kodaira-tangentspace-tangentcone} applies to $\wei(\X,E)$ and $\wei(\X,E,g)$ and the tangent space is the kernel of the derivative of these map.
Because of the following lemma, we refer to both~\eqref{eq:periods} and its extension~\eqref{eq:periods extension map}, as well as~\eqref{eq:periods imagsubspace} and its extension, as the period map.
\begin{lemma}
\label{lem:period welldefined}
\phantom{boo}
\begin{enumeratethm}[label={\upshape(\Alph*)}]
\item For any $(\var \U,\var\chi,\var\psi)$ in the tangent space of $\prewei(\X,E)$ at $(\U,\chi,\psi)$ the form $\lp\var\chi,\psi\rp + \lp\chi,\var\psi\rp$ is closed.
As a consequence the derivative of~\eqref{eq:periods} at $(\U,\chi,\psi)$ and its kernel, which is the tangent space of $\wei(\X,E)$, are independent of the choice of the extension.
The derivative of~\eqref{eq:periods} is given by
\begin{equation}
\label{eq:periods derivative}
(\var \U,\var\chi,\var\psi) \mapsto \per \big(\lp\var\chi,\psi\rp + \lp\chi,\var\psi\rp\big) .
\end{equation}
\item For any $(\var \U,\var\chi)$ in the tangent space of $\prewei(\X,E,g)$ at $(\U,\chi)$ the form $\lp\var\chi,g\chi\rp+\lp\chi,g\var\chi\rp$ is closed.
As a consequence the derivative of~\eqref{eq:periods imagsubspace} at $(\U,\chi)$ and its kernel, which is the tangent space of $\wei(\X,E,g)$, are independent of the choice of the extension.
The derivative of~\eqref{eq:periods imagsubspace} is given by
\begin{equation}
\label{eq:periods imagsubspace derivative}
(\var \U,\var\chi) \mapsto \per \big(\lp\var\chi,g\chi\rp+\lp\chi,g\var\chi\rp\big) .
\end{equation}
\end{enumeratethm}
\end{lemma}
\begin{proof}
We only prove~(A).
The proof of (B) is similar.
The given element $(\var \U,\var\chi,\var\psi)$ of the tangent space of $\prewei(\X,E)$ obeys
\begin{align}
(\delbar{E^{-1}}-\U\sd)\var\chi&=\var\U\sd\chi,&(\delbar{KE}-\U)\var\psi&=\var\U\psi.
\end{align}
Thus with the aid of Lemma~\ref{lem:11 pairing properties}(iv) we calculate\hspace{5mm}
\begin{align*}d\left(\lp\var\chi,\psi\rp+\lp\chi,\var\psi\rp \right)
&=-\overline{\lh\psi,(\delbar{E^{-1}}-\U\sd)\var\chi\rh}+\lh\var\chi,(\delbar{KE}-\U)\psi\rh\\&\hspace{4mm}-\overline{\lh\var\psi,(\delbar{E^{-1}}-\U\sd)\chi\rh}+\lh\chi,(\delbar{KE}-\U)\var\psi\rh\\
&=-\overline{\lh\psi,\var\U\sd\chi\rh}+\lh\chi,\var \U \psi\rh= 0,
\end{align*}
with the last equality following from the equation in Lemma~\ref{lem:11 pairing properties}(iv) on the right-hand side.
The derivative~\eqref{eq:periods extension map} in the direction of $(\var \U, \var\chi, \var\psi)$ is
\[\left(\gamma\mapsto\int_\X \eta_{\gamma}\wedge \left( \lp\var\chi,\psi\rp + \lp\chi,\var\psi\rp \right) \right)\in\Hom(H_1(\X,\mathbb{Z}),\qat).\]
Since the form is closed, this is independent of the choice of the $\eta_{\gamma}$.
\end{proof}

Now that we have set up the basic concept of the space of Weierstraß triples and its tangent space, let us introduce a tool that will help us to control variations of the potential that preserve the periods of the admissible map.
This control is the subject of Lemma~\ref{lem:wei-necessary condition}.
It will allow us to prove our first result about the space of Weierstraß triples in Lemma~\ref{lem:wei-isothermic} that the space of not strongly isothermic Weierstraß triples is a Banach submanifold.

But before we launch into this examination, we observe that we can apply the same trick as for $\kodaira(\X,E)$ to simplify the potential of the triple.
We apply Corollary~\ref{gauge holomorphic structure} simultaneously to both sections of a given triple $(\U,\chi,\psi)\in\prewei(\X,E)$ as described in Chapter~\ref{chapter:riemann roch}.
This regauging transforms the quadruple $(E, \U,\chi,\psi)$ into another quadruple $(E',\U',\chi',\psi')$ such that $\U'\in\pot{KE'}^-$, $\chi'\in H^0(\X,\Q{(E')^{-1},{\U'}\sd})$, $\psi'\in H^0(\X,\Q{KE',\U'})$ and $\lp\chi',\psi'\rp =\lp\chi,\psi\rp$.
For the sake of simplicity we omit the primes from the notations of these objects in the sequel.

Let $\Sh{U}\sd=\Q{E^{-1},\U\sd}$ and $\Sh{V}\sd=\Q{KE^{-1},\V\sd}$ denote the Serre dual sheaves of $\Sh{U}=\Q{KE,\U}$ and $\Sh{V}=\Q{E,\V}$, respectively.
Due to the local Darboux transformation~\ref{thm:darboux}(1)(b) and Corollary~\ref{cor:global darboux} any section $\chi \in \Sh{U}\sd$ has a $(1,0)$\=/derivative $\partial_l\chi=(\B_l+ \Bar{\V}_l)\chi_l$, for some $\mathbb{C}$\=/valued $\B_l$ and anti-commuting potential $\V \in \pot{E}^-$.
Hence $\chi$ determines $\V$ and also a section $\upsilon=-(\overline{\qj \chi})^{-1}\in H^0(\X,\Sh{V})$ without roots of the dual bundle $E_\qat$, see Theorem~\ref{thm:darboux}(2)(c).
We obtain the following two period maps:
\begin{align}\label{eq:period 1}
H^0(\X,\Sh{V}\sd)&\to\Hom(H_1(\X,\mathbb{Z}),\qat),&\upsilon\pa&\mapsto \per\lp \upsilon,\upsilon\pa\rp, \\
\label{eq:period 2}
H^0(\X,\Sh{U})&\to\Hom(H_1(\X,\mathbb{Z}),\qat),&\chi\pa&\mapsto \per\lp\chi,\chi\pa\rp.
\end{align}
We now wish to prove a reciprocity law between the two integrands.
The essential calculation is contained in the following lemma.
\begin{lemma}
\label{lem:period pairing}
\index{Pairing!Periods}
The exists a symplectic pairing (non-degenerate, anti-symmetric, real) on $\Hom(H_1(\X,\mathbb{Z}),\qat)$, which for any choice of canonical basis $A_1,\ldots,A_g$, $B_1,\ldots,B_g$ of $H_1(\X,\mathbb{Z})$ may be computed with the formula:
\begin{equation}
\label{eq:pairing}
\langle\beta,\alpha\rangle := \sum_{i=1}^g\Real\big(\Bar{\beta}(A_i)\alpha(B_i)\!-\!\Bar{\beta}(B_i)\alpha(A_i)\big),
\end{equation}
Let $\omega$ and $\eta$ be two closed $\qat$\=/valued $1$\=/forms on $\X$.
Then 
\begin{align}
\label{eq:wedge gives period pairing}
\int_{\X} \Real \Bar{\omega} \wedge \eta
= \langle\per\omega, \per\eta\rangle.
\end{align}
\end{lemma}
\begin{proof}
Consider the intersection pairing on $H_1(\X,\mathbb{Z})$.
Its definition is standard, and can be found in \cite[Section~III.1]{FK} among other sources.
It induces a pairing on $H_1(\X,\mathbb{Z})^*$, which in turn can be extended linearly to a $\qat\otimes_{\mathbb{Z}}\qat$\=/valued pairing on $\Hom(H_1(\X,\mathbb{Z}),\qat)$.
Now composing with the inner product on $\qat$ gives desired pairing.
In particular, it is anti-symmetric and non-degenerate.
As a formula, if $\gamma_i$ is a basis of $H_1(\X,\mathbb{Z})$ then
\[
\langle \beta, \alpha \rangle
= \sum_{i,j} \Real \overline{\beta(\gamma_i)} \alpha(\gamma_j) \iota(\gamma_i,\gamma_j).
\]
If we have a canonical basis, then the intersection form takes a special form (indeed, this is the definition of a canonical basis), and the above formula becomes Equation~\eqref{eq:pairing}:
\[
\langle \beta, \alpha \rangle
= \Real \overline{(\beta(A_1), \ldots, \beta(B_g))}
\begin{pmatrix}
0 & I_g \\ - I_g & 0
\end{pmatrix}
\begin{pmatrix}
\alpha(A_1) \\ \vdots \\ \alpha(B_g)
\end{pmatrix}.
\]

For the second formula, we use the standard argument for the Riemann bilinear relations found in the remark following~\cite[III.2.3~Proposition]{FK}.
Represent $X$ as a $4\genus$\=/sided polygon $\Hat{X}$ inside the universal cover $\Tilde{X}$, whose boundary represents the canonical basis $A_1,\ldots,A_g,B_1,\ldots,B_g$ of $H_1(\X,\mathbb{Z})$, oriented correctly with respect to the orientation of $\X$.
Let $G$ be a primitive of $\omega$ on $\Tilde{\X}$, that is $dG = \omega$.
Then as a first step we transform the integral over $\X$ to an integral around $\boundary \hat{\X}$ in $\Tilde{\X}$ using Stokes' theorem:
\begin{align*}
\int_{\X} \Real \Bar{\omega} \wedge \eta
= \Real \int_{\Hat{\X}} d( \Bar{G} \eta)
= \Real \int_{\boundary \Hat{\X}} \Bar{G}\eta.
\end{align*}
The difference of the values of $G$ at two boundary points of $\Hat{\X}$ over the same point of $\X$ differ by the corresponding period of $\omega$.
By breaking the integral over the boundary into a sum of integrals over each side of the polygon and pairing them together, we get a sum where each term is either of the form
\[
\Real \left( \int_{A_i} \Bar{G}\eta + \int_{-A_i} \Bar{G}\eta \right)
= \Real \int_{A_i} \left(-\per \Bar{\omega}(B_i)\right) \eta
= - \Real \left(\overline{\per \omega(B_i)} \per \eta(A_i) \right)
\]
or of the form
\[
\Real \left( \int_{B_i} \Bar{G}\eta + \int_{-B_i} \Bar{G}\eta \right)
= \Real \int_{B_i} \per \Bar{\omega}(A_i) \eta
= \Real \left(\overline{\per \omega(A_i)} \per \eta(B_i) \right).
\]
Taking the sum gives the right hand side of Equation~\eqref{eq:wedge gives period pairing}.
\end{proof}

As applied to~\eqref{eq:period 1} and~\eqref{eq:period 2}
\begin{align*}
\langle \per\lp \upsilon,\upsilon\pa\rp, \per\lp \chi,\chi\pa\rp \rangle
&= \Real \int_\X \overline{\lp \upsilon,\upsilon\pa\rp} \wedge \lp \chi,\chi\pa\rp \\
&= \Real \int_\X \overline{\upsilon\pa} (-\qj dz) \upsilon \wedge \upsilon^{-1}\qj^{-1}\qj dz \chi\pa 
= 0.
\end{align*}
This shows that the images of~\eqref{eq:period 1} and~\eqref{eq:period 2} are orthogonal with respect to the period pairing.
Let us now show that these images are the orthogonal complements of each other.
\begin{lemma}\label{pairing}
Let $\U\in\pot{KE}$ on the compact Riemann surface $\X$, let $\Sh{U}=\Q{KE,\U}$ and let $\Sh{U}\sd = \Q{E^{-1},\U\sd}$ be the Serre dual sheaf, see Chapter~\ref{chapter:riemann roch}.
Let $\chi\in H^0(\X, \Sh{U}\sd)$ have no roots, and let $\Sh{V} = \Q{E,\V}$ be the sheaf of holomorphic sections of the $\qat$\=/line bundle with the section $\upsilon=-(\overline{\qj \chi})^{-1} \in H^0(\X, \Sh{V})$ and $\Sh{V}\sd=\Q{KE^{-1},\V\sd}$ the corresponding Serre dual sheaf.
Then the images of~\eqref{eq:period 1} and~\eqref{eq:period 2} are the orthogonal complements of each other with respect to the period pairing~\eqref{eq:pairing}.
\end{lemma}
\begin{proof}
We have proved above the orthogonality for the case $\V^+=0$.
For non-zero $\V^+$ we now show that we can apply Corollary~\ref{gauge holomorphic structure} simultaneously to the sheaves $\Sh{U}$, $\Sh{U}\sd$, $\Sh{V}$ and $\Sh{V}\sd$.
Let $\exp(h_l)$ denote on the open cover $\{\SO_l\}$ the local complex functions which in the proof of this Corollary transform the local functions $\upsilon_l$ of the sections $\upsilon$ of $\Sh{V}$ into the local functions $\upsilon'_l=\exp(h_l)\upsilon_l$ of the sections $\upsilon'$ of $\Sh{V'}$ with $\V'\in\pot{E'}^-$.
Consequently the local functions of $\chi=\big(\overline{\qj\upsilon})^{-1}$ and $\chi'=\big(\overline{\qj\upsilon'})^{-1}$ of the sections $\chi$ and $\chi'$ of $\Sh{U}\sd$ and $\Sh{U'}\sd$, respectively, are transformed by $\exp(-h_l)$:
\[
\chi'_l=\big(\overline{\qj\upsilon'_l}\big)^{-1}=\big(\overline{\qj \exp(h_l)\upsilon_l}\big)^{-1}=\big(\overline{\exp(\Bar{h}_l)\qj\upsilon_l}\big)^{-1}=\big(\overline{\qj\upsilon_l}\exp(h_l)\big)^{-1}=\exp(-h_l)\chi_l.
\]
By Theorem~\ref{thm:darboux}(1)(b) the commuting part $({\U'}\sd)^+$ of the potential ${\U'}\sd$ of ${\Sh{U}'}\sd$ vanishes and the commuting part of the potential $\U\sd$ of $\Sh{U}\sd$ is equal to $(\U\sd)^+=-\V^+=0$.
In particular, this application changes the corresponding holomorphic line bundles $E$ and $E^{-1}$ by inverse local functions into line bundles $E'$ and $E'^{-1}$ which are the inverse of each other.
Because the transformed ${\Sh{V'}}\sd$ and $\Sh{U}'$ are supposed to be the Serre dual of $\Sh{V}'$ and ${\Sh{U}'}\sd$ the local functions $\upsilon\pa_l$, ${\upsilon'}\pa_l$, $\psi_l$ and $\psi'_l$ of the corresponding Serre dual sheaves $\Sh{V}\sd$, ${\Sh{V}'}\sd$, $\Sh{U}$ and $\Sh{U}'$ transform as
\begin{align*}
{\upsilon'}\pa_l&=\exp(-h_l)\upsilon\pa_l,&\psi'_l&=\exp(h_l)\psi_l,
\end{align*}
respectively.
Due to Lemma~\ref{lem:10 pairing properties}(i) this transformation preserves the $1$\=/forms
\begin{align*}\lp\chi,\chi\pa\rp&=\lp\chi',{\chi'}\pa\rp,&\lp\upsilon,\upsilon\pa\rp&=\lp\upsilon',{\upsilon'}\pa\rp.
\end{align*}
Therefore the orthogonality follows in the general case from the special case $\V^+=0$.

It remains to prove in the special case $\V^+=0=\U^+$ that the quaternionic dimensions of the images of the maps~\eqref{eq:period 1} and~\eqref{eq:period 2} add up to $2\genus$.
The non-vanishing sections $\chi$ and $\upsilon=-(\overline{\qj\chi})^{-1}$ define connections on $E^{-1}_\qat$ and $E_\qat$, whose $(1,0)$\=/parts are respectively
\begin{align*}
\partial_l-\B_l+\V_l\sd&:H^0(\X,\Sh{U}\sd)\to H^0(\X,\Sh{V}\sd),&
\partial_l+\B_l+\U_l&:H^0(\X,\Sh{V})\to H^0(\X,\Sh{U}).
\end{align*}
Due to Corollary~\ref{cor:global darboux} (specifically the global version of Theorem~\ref{thm:darboux}(1)(e)), the images of these maps are the kernels of~\eqref{eq:period 1} and~\eqref{eq:period 2}.
We recall the zero curvature equations (see Equation~\eqref{eq:zero curvature 2}) of connections defined by $\chi$ and $\upsilon$:
\begin{align*}
(\barpartial_l + \V\sd_l)(\partial_l-\B_l+\V_l\sd)&=(\partial_l - \B_l + \U_l\sd)(\barpartial_l-\U_l\sd),\\
(\barpartial_l-\U_l)(\partial_l+\B_l+\U_l)&=(\partial_l+\B_l+\V_l)(\barpartial_l-\V_l).
\end{align*}
The map $\partial_l+\B_l-\V_l\sd$ has the one-dimensional kernel $\chi_l\qat$ and the map $\partial_l-\B_l-\U_l$ has the one-dimensional kernel $\upsilon_l\qat$.
The dimensions of the corresponding images are then respectively $\dim H^0(\X,\Sh{U}\sd) - 1$ and $\dim H^0(\X,\Sh{V}) -1$.
Hence the sum of the dimensions of the images of~\eqref{eq:period 1} and~\eqref{eq:period 2} is equal to
\[\left[\dim H^0(\X,\Sh{V}\sd) - \dim H^0(\X,\Sh{U}\sd) + 1 \right]
+ \left[ \dim H^0(\X,\Sh{U}) - \dim H^0(\X,\Sh{V}) + 1 \right]\]
Due to Serre duality~\ref{Serre duality} and the Riemann-Roch Theorem~\ref{riemann roch} this number is equal to $2-\deg(E)+\deg(KE)=2\genus$.
\end{proof}
We shall need an analogous statement for $\psi$, but in a slightly more general situation since $\psi$ might have roots.
For given $(\U,\chi,\psi)\in\prewei(\X,E)$ we use the notation of the foregoing lemma, so $\psi\in H^0(\X,\Sh{U})$.
With the divisor $D=(\psi)$ it is a non-vanishing section of the sheaf $\Sh{U}(-D)=\Q{KE(-D),\U(-D)}$ (see Lemma~\ref{quotient dimension}) with the Serre dual sheaf $\Sh{U}\sd(D)=\Q{(E(-D))^{-1},\U(-D)\sd}$.
By Theorem~\ref{thm:darboux} and Corollary~\ref{gauge holomorphic structure} $\zeta=-(\overline{\qj\psi})^{-1}$ is a section of the sheaf $\Sh{V}'=\Q{(KE)^{-1}(D),\V'}$ with potential $\V'\in\pot{(KE)^{-1}(D)}$. 
Let ${\Sh{V}'}\sd=\Q{K^2E(-D),{\V'}\sd}$ denote the Serre dual sheaf of $\Sh{V}'$. Note that the commuting part of $\V'$ is equal to ${\V'}^+=(\U\sd)^+=-\U^+=-\V^+$, since by~\eqref{eq:dchi} the commuting part of $\U$ is equal to $\U^+=\V^+$.
%
%
In this case the forgoing Lemma applies to the sheaves $\Sh{U}\sd(D)$, $\Sh{U}(-D)$, $\Sh{V}'$ and ${\Sh{V}'}\sd$ and shows that images of the following maps are the orthogonal complements of each other with respect to the period pairing~\eqref{eq:pairing}:
\begin{gather}\begin{aligned}\label{eq:period 3a}
H^0(\X,{\Sh{V}'}\sd)&\to\Hom(H_1(\X,\mathbb{Z}),\qat),&
\zeta\pa &\mapsto \per \lp\zeta,\zeta\pa\rp,\\
H^0(\X,\Sh{U}\sd(D))&\to\Hom(H_1(\X,\mathbb{Z}),\qat),&
\psi\pa &\mapsto \per \lp\psi,\psi\pa\rp.
\end{aligned}\end{gather}

We want to show that an analogous statement holds with $\Sh{U}\sd$ instead of $\Sh{U}\sd(D)$ and a singular holomorphic sheaf $\Sh{S}$ in the place of ${\Sh{V}'}\sd$.
For this purpose we apply Lemma~\ref{quotient dimension} and conceive the section $\zeta=-(\overline{\qj\psi})^{-1}$ of $\Sh{V}'$ without roots as a meromorphic section of $\Sh{V}'(-D)=\Q{(KE)^{-1},\V'(-D)}$.
For any germ $\zeta\pa$ of a meromorphic section of the Serre dual sheaf $\Sh{V}'(-D)\sd=\Q{K^2E,\V'(-D)\sd}$ at some point $x\in X$ there exists an open neighborhood $\Set{O}$ of $x$, such that $\zeta\pa\in H^0(\Set{O}\setminus\{x\},\Sh{V}'(-D)\sd)$.
In particular, the $1$\=/form $\lp\zeta,\zeta\pa\rp$ is closed on $\Set{O}\setminus\{x\}$.
We define $\Sh{S}$ as the unique subsheaf of the meromorphic sections of $\Sh{V'(-D)\sd}$, whose stalks $\Sh{S}_x$ at any $x\in X$ contain all germs $\zeta\pa$ of meromorphic sections of $\Sh{V}'(-D)\sd$ with poles of at most first order, such that the integral of $\lp\zeta,\zeta\pa\rp$ along a small circle around $x$ vanishes. 

For such a section $\zeta\pa\in\Sh{S}_x$ at a point $x \in \X\setminus D$ there exists locally a function $H$ defined through the Weierstraß representation $dH = \lp\zeta,\zeta\pa\rp$.
From its Kodaira representation we conclude $\psi H$ is $\U$\=/holomorphic.
At a point $x \in D = (\psi)$ by our definition of $\Sh{S}$ we can integrate $dH = \lp\zeta,\zeta\pa\rp$ to a function $H$ on $\Set{O}\setminus\{x\}$.
Furthermore, since $\zeta\pa$ has a pole of at most first order at $x$, the derivative of the function $H$ can have at $x$ a pole of at most one order higher than the pole of the meromorphic section $\zeta=-(\overline{\qj\psi})^{-1}$ of $\Sh{U}\sd$, which is equal to order of the root of $\psi$ at $x$. 
Since by Lemma~\ref{lem:change of order} the order of any pole of $H$ is one less than the corresponding pole order of the derivative $dH$, it follows that $\psi H$ is a holomorphic section of $\Sh{U}_x$.
In summary, $\psi H$ belongs to $\Sh{U}$ at all points of $\X$.

We can extend this argument to show that $\Sh{S}$ is a singular holomorphic sheaf in the sense of Definition~\ref{def:singular sheaf}.
Applying Lemma~\ref{lem:change of order} again shows that for any $x\in X\setminus D$ and any stalk $\chi\pa$ of a meromorphic section of $\Sh{V}'(-D)\sd$ at $x$ the derivative of $H=\psi^{-1}\chi\pa$ can have poles of all orders with exception of first order poles. 
Then Lemma~\ref{quotient dimension} implies the stalks $\Sh{S}_x$ of $\Sh{S}$ coincide with $\Sh{V}'(-D)\sd$ at all points of $\X\setminus D$. 
Furthermore, at the roots of $\psi$ they are the subsheaf of the meromorphic sections of $\Sh{S}_x$ with poles of at most first order at $x$ such that both stalks $\Sh{S}_x$ and $\Sh{V}'(-D)\sd_x$ have codimension one in the subspace $\Sh{S}_x+\Sh{V}'(-D)\sd_x$ of the stalk of meromorphic sections of $\Sh{V}'(-D)\sd_x$ with at most first order poles at $x$.
So $\Sh{S}$ is a singular holomorphic sheaf with underlying sheaf $\Sh{V}'(-D)\sd$ and $\deg(\Sh{S})=\deg(K^2E)$.

We now define the analogues of~\eqref{eq:period 3a}. 
Pay heed to the change of ordering of $\psi$ and $\psi\pa$.
This is done to avoid an excessive number of conjugations later on.
\begin{align}\label{eq:period 3}
H^0(\X,\Sh{S})&\to\Hom(H_1(\X,\mathbb{Z}),\qat),&\zeta\pa&\mapsto \per \lp\zeta\pa,\zeta\rp,\\
H^0(\X,\Sh{U}\sd)&\to\Hom(H_1(\X,\mathbb{Z}),\qat),&\psi\pa&\mapsto \per\lp\psi\pa,\psi\rp.\label{eq:period 4}
\end{align}
\begin{lemma}\label{pairing 2}
On a compact Riemann surface for all $(\U,\chi,\psi)\in\prewei(\X,E)$ the images of the two period maps~\eqref{eq:period 3} and~\eqref{eq:period 4} are the orthogonal complements of each other with respect to the period pairing~\eqref{eq:pairing}.
\end{lemma}

\begin{proof}
First note that the proof of Lemma~\ref{lem:period pairing} only requires that the differential locally has a primitive.
Therefore this lemma also applies to $dH=\lp\zeta,\zeta\pa\rp$ for $\zeta\pa\in H^0(\X,\Sh{S})$, even though it may have poles.
Then
\begin{align*}
\langle \per \lp\zeta\pa,\zeta\rp, \per \lp\psi\pa,\psi \rp \rangle
&= - \Real \int_X \lp\psi\pa,\psi \rp \wedge \overline{\lp\zeta\pa,\zeta\rp} \\ 
&= - \Real \int_X \overline{\psi\pa} \qj dz \psi \wedge (-\psi^{-1}\qj^{-1})(-\qj dz)\overline{\zeta\pa}
= 0.
\end{align*}
The other arguments of the proof of Lemma~\ref{pairing} also extend to the present situation; the Riemann-Roch theorem~\ref{thm:riemann roch singular} and Serre duality~\ref{thm:serre duality singular} imply that the dimensions of the images of the maps~\eqref{eq:period 3} and~\eqref{eq:period 4} add up to $2\genus$.
\end{proof}

As we did in Chapter~\ref{chapter:isothermic kodaira triples} we will control variations of the potential with a map analogous to~\ref{eq:singularity1}, namely
\[
(\chi G, \psi H) 
\mapsto \ls\chi G,\psi\rs+\ls\chi,\psi H\rs 
= \ls\chi (G + \Bar{H}),\psi\rs,
\]
using Lemma~\ref{lem:01 pairing properties}(ii) for the equality.
An important difference is that $G$ and $H$ are in general defined on $\Tilde{X}$, not $\X$ itself.
The potential is single-valued on $\X$ (and therefore an element of $\pot{KE}$) if and only if $G + \Bar{H}$ is, which in turn is equivalent to the period condition in~\ref{single valued} below.
Therefore we need to introduce a space of admissible maps on the universal covering $\Tilde{X}$ associated to the period maps~\eqref{eq:period 1} and~\eqref{eq:period 3} but for which the above potential is single valued on $\X$.
\begin{definition}
\label{def:extended pairs}
\phantom{boo}
\begin{enumeratethm}[label={\upshape(\Alph*)}]
\item
For given $(\U,\chi,\psi)\in\prewei(\X,E)$ on a Riemann surface $\X$ let $\Sh{U}$, $\Sh{U}\sd$, $\Sh{V}\sd$ and $\Sh{S}$ denote the four sheaves introduced in the Lemmas~\ref{pairing} and~\ref{pairing 2}.
Any $(\upsilon\pa,\zeta\pa)\in H^0(\X,\Sh{V}\sd)\times H^0(\X,\Sh{S})$ defines an admissible map $G$ on the universal covering $\Tilde{X}$ of $\X$ with $dG=\lp\upsilon,\upsilon\pa\rp$ and an admissible map $H$ on the complement of the roots of the pullback of $\psi$ to $\Tilde{X}$ with $dH=\lp\zeta,\zeta\pa\rp$.
By the definition of $\Sh{S}$, we know that $(\chi G, \psi H)$ are sections of the pullbacks of the sheaves $\Sh{U}\sd$ and $\Sh{U}$ to $\Tilde{X}$ respectively, in particular also at the roots of $\psi$.
We define 
\begin{align}
\label{single valued}
\con(\X,\chi,\psi)
:= \left\{ (\chi G,\psi H) \mid \per dG = - \per d\Bar{H} \right\} / \{ (\chi q,-\psi\Bar{q}) \mid q \in \mathbb{H}\}
\end{align}

\item
For given $(\U,\chi)\in\prewei(\X,E,g)$ on a Riemann surface $\X$, let $\Sh{U}$ denote the sheaf introduced in the Lemmas~\ref{pairing}-\ref{pairing 2}.
If $\V\in\pot{E}^-$ denotes as in Lemma~\ref{pairing} the potential of $\upsilon=-(\overline{\qj\chi})^{-1} \in H^0(\X,\Sh{V})$, then $\zeta=-(\overline{\qj\psi})^{-1}=-g^{-1}(\overline{\qj\chi})^{-1}=g^{-1}\upsilon$ is a meromorphic section of $\Q{(KE)^{-1},\qj\frac{g}{\Bar{g}}\V}$ with the potential $g^{-1}\V g=\frac{\Bar{g}}{g}\V\in\pot{(KE)^{-1}}^-$ (compare Equation~\eqref{eq:potential transformation}).
For the corresponding sheaf $\Sh{S}$ any $\zeta\pa\in H^0(\X,\Sh{S})$ defines on the complement of the preimage of the roots of $g$ in $\Tilde{\X}$ an admissible map $H$ with $dH=\lp\zeta,\zeta\pa\rp$ and $g\chi H\in H^0(\Tilde{\X},\Sh{U})$.
Due to Lemma~\ref{lem:01 pairing properties}(i) and (iii) $\ls\chi H,g\chi\rs=g|g|^{-2}\ls g\chi H,\chi\rs|g|^2g^{-1}=g{\ls\chi,g\chi H\rs}\sd g^{-1}$ and defines locally a $\banach{2}$\=/potential in $\pot{KE}$.
Moreover, the potential $\ls\chi H,g\chi\rs+\ls\chi,g\chi H\rs$ is single valued on $\X$ (and therefore an element of $\pot{KE,g}$) if and only if
\begin{gather}
\label{single valued 2}
\per dH (\gamma) = \int_\gamma dH
\in\Imag\qat, 
\quad\text{for all }\gamma\in H_1(\X,\mathbb{Z}).
\end{gather}
\normalsize
We define 
\begin{align}
\con(\X,\chi,g)
:= \left\{ g\chi H\in H^0(\Tilde{\X},\Sh{U}) \mid dH \text{ obeys~\eqref{single valued 2}} \right\}/
\{g\chi q\mid q \in \Imag\mathbb{H}\}.
\end{align}
\end{enumeratethm}
\end{definition}
The key feature of the spaces $\con(\X,\chi,\psi)$ and $\con(\X,\chi,g)$ is that they extend the spaces $(\psi\pa,\chi\pa)\in H^0(\X,\Sh{U}\sd)\times H^0(\X,\Sh{U})$ and $\chi\pa\in H^0(\X,\Sh{U})$ to elements of $H^0(\Tilde{\X},\Sh{U}\sd)\times H^0(\Tilde{\X},\Sh{U})$ and $H^0(\Tilde{\X},\Sh{U})$, such that the pairings $\ls\psi\pa,\psi\rs+\ls\chi,\chi\pa\rs$ and $\ls g^{-1}\chi\pa,g\chi\rs+\ls\chi,\chi\pa\rs$ remain single valued on $\X$, respectively.
Indeed, due to the global quaternionic Darboux transformation~\ref{cor:global darboux} such elements are of the form $(\chi G,\psi H)$ and $g\chi H$ as defined above, respectively.
They are determined by $(\upsilon\pa,\zeta\pa)$ and $\zeta\pa$ only up to $\{(\chi q,-\psi\Bar{q})\mid q\in \mathbb{H}\}$ and $\{g\chi q\mid q\in\Imag\qat\}$, respectively.

The spaces $\con(\X,\chi,\psi)$ and $\con(\X,\chi,g)$ are related to constraints on variations of the potential $\var \U$.
The constraint that would be represented by $(\chi,\psi)$ and $g\chi$ is trivial, which is why we quotient it out from the definitions of $\con(\X,\chi,\psi)$ and $\con(\X\chi,g)$, respectively.
However, for simplicity we will denote the elements of $\con(\X,\chi,\psi)$ and $\con(\X,\chi,g)$ by $(\chi G,\psi H)$ and $g\chi H$, and not by equivalence classes, with the understanding that our arguments do not depend on the choice of representative.
We reserve the symbols $(\psi\pa,\chi\pa)$ and $\chi\pa$ for when $(\chi G,\psi H)$ and $g\chi H$ are single valued on $X$, respectively.

We now derive for given $(\U,\chi,\psi)\in\wei(\X,E)$ a condition on $\var \U\in\pot{KE}$ which is equivalent to the existence of $(\var\chi,\var\psi)$ such that $(\var \U,\var\chi,\var\psi)$ belongs to the tangent space of $\wei(\X,E)$ at $(\U,\chi,\psi)$; and for given $(\U,\chi)\in\wei(\X,E,g)$ a condition on $\var \U\in\pot{KE,g}$ which is equivalent to the existence of $\var\chi$ such that $(\var \U,\var\chi)$ belongs to the tangent space of $\wei(\X,E,g)$ at $(\U,\chi)$.
Hence such $(\var \U,\var\chi,\var\psi)$ belong to the kernel of~\eqref{eq:periods derivative}, which is derivative of~\eqref{eq:periods}; and such $(\var \U,\var\chi)$ to kernel of~\eqref{eq:periods imagsubspace derivative}, which is the derivative of~\eqref{eq:periods imagsubspace}.
\begin{lemma}
\label{lem:wei-necessary condition}
\phantom{boo}
\begin{enumeratethm}[label={\upshape(\Alph*)}]
\item
Given $(\U,\chi,\psi)\in\wei(\X,E)$ any $\var \U\in\pot{KE}$ is the first component of a triple $(\var \U,\var\chi,\var\psi)$ which is tangent to $\wei(\X,E)$ at $(\U,\chi,\psi)$ if and only if $\var \U$ is orthogonal to the image of
\begin{align}\label{eq:singularity2b}
\con(\X,\chi,\psi)&\to\pot{KE},&(\chi G,\psi H)&\mapsto\ls\chi G,\psi\rs+\ls\chi,\psi H\rs.
\end{align}
The following map maps the kernel of~\eqref{eq:singularity2b} onto the orthogonal complement of the image of~\eqref{eq:periods derivative}:
\begin{align}\label{eq:translation period}
\con(\X,\chi,\psi)&\to\Hom(H_1(\X,\mathbb{Z}),\qat),
&(\chi G,\psi H) &\mapsto \per dG=- \per d\Bar{H}.
\end{align}
The codimension of the image of~\eqref{eq:periods derivative} is equal to the codimension of the kernel of the following map in the kernel of~\eqref{eq:singularity2b}:
\begin{gather}
\begin{aligned}
\label{eq:singularity2a}
H^0(\X,\Sh{U}\sd)\times H^0(\X,\Sh{U})/\{(\chi q,-\psi\Bar{q})\mid q\in\mathbb{H}\}&\to\pot{KE},\\
(\psi\pa,\chi\pa)&\mapsto\ls\psi\pa,\psi\rs+\ls\chi,\chi\pa\rs.
\end{aligned}
\end{gather}
\item
Given $(\U,\chi)\in\wei(\X,E,g)$ any $\var \U\in\pot{KE,g}$~\eqref{def:real potential} is the first component of a pair $(\var \U,\var\chi)$ which is tangent to $\wei(\X,E,g)$ at $(\U,\chi)$ if and only if $\var \U$ is orthogonal to the image of
\begin{align}\label{eq:singularity2c}
\con(\X,\chi,g)&\to\pot{KE,g},&g \chi H&\mapsto\ls\chi H,g\chi\rs+\ls\chi,g\chi H\rs.
\end{align}
The following map maps the kernel of~\eqref{eq:singularity2c} onto the orthogonal complement of the image of~\eqref{eq:periods imagsubspace derivative}:
\begin{align}\label{eq:translation period 2}
\con(\X,\chi,g)&\to\Hom(H_1(\X,\mathbb{Z}),\Imag\qat),&g \chi H&\mapsto \per dH.
\end{align}
The codimension of the image of~\eqref{eq:periods imagsubspace derivative} is equal to the codimension of the kernel of the following map in the kernel of~\eqref{eq:singularity2c}:
\begin{gather}
\label{eq:singularity2d}
\begin{aligned}
H^0(\X,\Sh{U})/\{g\chi q\mid q \in \Imag\mathbb{H}\}&\to\pot{KE,g},\\
\chi\pa&\mapsto\ls g^{-1}\chi\pa,g\chi\rs+\ls\chi,\chi\pa\rs.
\end{aligned}
\end{gather}
\end{enumeratethm}
\end{lemma}
\begin{proof}
\proofstep{(A):} Observe that the three maps~\eqref{eq:singularity2b}-\eqref{eq:singularity2a} are well-defined, since the addition of $ (\chi q, -\psi \Bar{q})$ to $(\psi\pa,\chi\pa)$ for constant $q\in \mathbb{H}$ does not change their values.

First we deal with the condition to be tangent to $\prewei(\X,E)$.
For any given $(\U,\chi,\psi)$ and $\var \U\in\pot{KE}$ there exist $(\var\chi,\var\psi)$ such that $(\var \U,\var\chi,\var\psi)$ belongs to the kernel of the derivative of~\eqref{eq:prewei map}, if and only if $(\var \U\sd\chi,\var \U\psi)$ belongs to the image of the Fredholm operator
\begin{gather}\label{fredholm weierstrass}\begin{aligned}
H^0(\X,\sob{1,p}{E^{-1}})\times H^0(\X,\sob{1,p}{KE})&\to H^0(\X,\forms{0,1}{}\ban{p}{E^{-1}})\times H^0(\X,\forms{0,1}{}\ban{p}{KE}),\\(\var\chi',\var\psi')&\mapsto\left((\delbar{E^{-1}} -\U\sd)\var\chi',(\delbar{KE}-\U)\var\psi'\right).
\end{aligned}\end{gather}
This image is the orthogonal complement of the kernel of the adjoint operator.
Hence, by definition of $\ls\cdot,\cdot\rs$, this is equivalent to the condition that $\var \U\in\pot{KE}$ is orthogonal to the image of~\eqref{eq:singularity2a}.

Such a $\var \U\in\pot{KE}$ determines a tangent vector $(\var \U, \var\chi, \var\psi)$ to $\prewei(\X,E)$ up to the addition of $(0, \var\chi'', \var\psi'')$ where $(\var\chi'',\var \psi'')$ belongs to the kernel of~\eqref{fredholm weierstrass}, namely $(\var\chi',\var\psi') \in H^0(\X,\Sh{U}\sd)\times H^0(\X,\Sh{U})$.
Thus at $(\U,\chi,\psi)\in\prewei(\X,E)$ the image of~\eqref{eq:periods derivative}, the derivative of the period map,  contains the image of
\begin{gather}\label{period dimension}\begin{aligned}
H^0(\X,\Sh{U}\sd)\times H^0(\X,\Sh{U})&\to\Hom(H_1(\X,\mathbb{Z}),\qat),\\
(\var\chi'',\var\psi'')&\mapsto \per \big(\lp\var\chi'',\psi\rp+\lp\chi,\var\psi''\rp\big).
\end{aligned}\end{gather}
This means that for given $\var \U \in \pot{KE}$ there exist a tangent vector to $\wei(\X,E)$ whose first component is $\var\U$ if and only if there exists a triple of the form $(\var \U, \var\chi,\var\psi)$ in the tangent space to $\prewei(\X,E)$ that is mapped by~\eqref{eq:periods derivative} to the image of~\eqref{period dimension}, since then $(\var \U, \var\chi-\var\chi'',\var\psi-\var\psi'')$ is tangent to $\wei(\X,E)$.
The existence of one such triple is equivalent to the condition that all such triples $(\var \U, \var\chi',\var\psi')$ in the tangent space to $\prewei(\X,E)$ are mapped to the image of~\eqref{period dimension} under~\eqref{eq:periods derivative}, because the difference of two such triples is of the form $(0,\var\chi'',\var \psi'')$ for some $(\var\chi'',\var\psi'')$ in the domain of~\eqref{period dimension}.

The image of~\eqref{period dimension} is the sum of the images of~\eqref{eq:period 2} and~\eqref{eq:period 4}.
Due to Lemmas~\ref{pairing} and~\ref{pairing 2} this is the sum of the orthogonal complements of the images of~\eqref{eq:period 1} and~\eqref{eq:period 3}.
By a simple vector space argument, this is then the orthogonal complement of the intersection of the images of~\eqref{eq:period 1} and~\eqref{eq:period 3}, which is by Definition~\ref{def:extended pairs} of~$\con(\X,\chi,\psi)$ the image of~\eqref{eq:translation period}.

To summarize the argument up to this point, given $\var \U \in \pot{KE}$ there exist $(\var\chi,\var\psi)$ such that $(\var \U,\var\chi,\var\psi)$ is tangent to $\wei(\X,E)$ at $(\U,\chi,\psi)$ if and only if
\begin{enumerate}
\item $\var \U$ is orthogonal to the image of~\eqref{eq:singularity2a}, and
\item the image of~\eqref{eq:translation period} is orthogonal to~\eqref{eq:periods derivative} for any tangent vector to $\prewei(\X,E)$ whose first first component is $\var \U$.
\end{enumerate}
So the first statement follows, if the conditions~(i) and (ii) on $\var \U\in\pot{KE}$ are equivalent to $\var \U$ being orthogonal to the image of~\eqref{eq:singularity2b}.
Since such $\var \U$ obey condition~(i), it suffices to show for any $\var \U$ obeying condition~(i) the equivalence of condition~(ii) and $\var \U$ being orthogonal to the image of~\eqref{eq:singularity2b}.
So we may assume that $(\var \U, \var \chi, \var\psi)$ is tangent to $\prewei(\X,E)$.
In particular, $(\delbar{E^{-1}} -\U\sd)\var\chi=\var \U\sd\chi$ and $(\delbar{KE}-\U)\var\psi=\var \U\psi$.
Lemma~\ref{lem:11 pairing properties}(iv) yields
\begin{multline*}
d\big(\lp\chi G,\var\psi\rp\! - \!\lp\var\chi,\psi H\rp\big)\\
\begin{aligned}
&=-\overline{\lh\var\psi,(\delbar{E^{-1}}-\U\sd)(\chi G)\rh}+\lh\chi G,(\delbar{KE}-\U)\var\psi\rh\\
&\hspace{4mm}+\overline{\lh\psi H,(\delbar{E^{-1}}-\U\sd)\var\chi\rh}-\lh\var\chi,(\delbar{KE}-\U)(\psi H)\rh\\
&= 0 + \lh\chi G,\var \U\psi\rh + \lh\chi,\var \U\psi H\rh + 0.
\end{aligned}
\end{multline*}
On the other hand
\begin{align*}
&\Real d\big[\lp\chi G,\var\psi\rp - \lp\var\chi,\psi H\rp\big]
= \Real d\big[\Bar{G}\lp\chi,\var\psi\rp - H\lp\var\chi,\psi \rp\big] \\
&= \Real d\Big[\Bar{G}\big( \lp\chi,\var\psi\rp+\lp\var\chi,\psi \rp \big)\Big] - \Real d\Big[ (\Bar{G} + H) \lp\var\chi,\psi \rp \Big] \\
&= \Real d\Bar{G} \wedge \big( \lp\chi,\var\psi\rp+\lp\var\chi,\psi \rp \big) - \Real d\Big[ (\Bar{G} + H) \lp\var\chi,\psi \rp \Big].
\end{align*}
Together we have
\begin{align*}
\labelthis{eq:period restriction obstruction}
&\left\langle\ls\chi G,\psi\rs+\ls\chi,\psi H\rs, \var \U\right\rangle 
=\int_{\X}\Real\big(\lh\chi G,\var \U\psi\rh+\lh\chi,\var \U\psi H\rh\big)\\
&=\int_{\X}\Real d\Bar{G} \wedge \big( \lp\chi,\var\psi\rp+\lp\var\chi,\psi \rp \big) - \int_{\X} \Real d\Big[ (\Bar{G} + H) \lp\var\chi,\psi \rp \Big]\\
&= \left\langle \per dG, \per\big(\lp\chi,\var\psi\rp+\lp\var\chi,\psi \rp \big) \right\rangle - 0.
\end{align*}
To explain the final line of the calculation, to the first term we have applied Lemma~\ref{lem:period pairing} and to the second term we have applied Stokes' theorem, noting that $(\Bar{G} + H) \lp\var\chi,\psi \rp$ is a $1$\=/form on $\X$, since $\Bar{G} + H$ is single-valued on $\X$ even though $G$ and $H$ separately may not be.
This shows that condition (ii) is equivalent to $\var \U$ being orthogonal to the image of~\eqref{eq:singularity2b}, proving the first statement.

To prove the second statement, 
we continue from the fact that the orthogonal complement of the image of~\eqref{period dimension} is the image of~\eqref{eq:translation period}.
Because~\eqref{period dimension} is the restriction of~\eqref{eq:periods derivative} to vectors of the form $(0,\var\chi,\var\psi)$, the orthogonal complement of~\eqref{eq:periods derivative} is a subset of the image of~\eqref{eq:translation period}.
We denote this subset by $\Spa{I}^\perp$.
Equation~\eqref{eq:period restriction obstruction} was derived under the assumption that $(\var\U,\var\chi,\var\psi)$ was tangent to $\prewei(\X,E)$.
So $\Spa{I}^\perp$ is characterized by the property that the preimage with respect to~\eqref{eq:translation period} is mapped by~\eqref{eq:singularity2b} into the orthogonal complement of all $\var \U\in\pot{KE}$ which can be supplemented to an element $(\var \U,\var\chi,\var\psi)$ of the tangent space of $\prewei(X,E)$ at $(\U,\chi,\psi)$.
As stated in the condition~(i) in the proof of the first statement this orthogonal complement is the image of~\eqref{eq:singularity2a}.
The preimage of this image with respect to~\eqref{eq:singularity2b} is the sum of the domain of~\eqref{eq:singularity2a} and the kernel of~\eqref{eq:singularity2b}.
Since~\eqref{eq:translation period} vanishes on the domain of~\eqref{eq:singularity2a}, the image $\Spa{I}^\perp$ of this sum with respect to~\eqref{eq:translation period} is equal to the image of the kernel of~\eqref{eq:singularity2b} with respect to~\eqref{eq:translation period}.
This proves the second statement.

From the second statement we get that the codimension of~\eqref{eq:periods derivative} is equal to the dimension of the image of restriction of~\eqref{eq:translation period} to the kernel of~\eqref{eq:singularity2b}.
We remark that the kernel of~\eqref{eq:singularity2a} is the intersection of the kernel of~\eqref{eq:singularity2b} and the kernel of~\eqref{eq:translation period}, so the image of the restriction of~\eqref{eq:translation period} to the kernel of~\eqref{eq:singularity2b} is isomorphic to the quotient of the kernel of~\eqref{eq:singularity2b} by the kernel of~\eqref{eq:singularity2a}.
This shows the third statement.

\proofstep{(B):} Observe that the three maps~\eqref{eq:singularity2c}-\eqref{eq:singularity2d} are well-defined, since the addition of $ g\chi q$ to $g\chi H$ for constant $q\in \mathbb{H}$ does not change their values.

By the following calculation using Lemma~\ref{lem:01 pairing properties}, the maps~\eqref{eq:singularity2c} and~\eqref{eq:singularity2d} (with $\chi\pa=g\chi H$) take indeed values in $\pot{KE,g}$~\eqref{def:real potential}:
\begin{align*}
g{\ls\chi H,g\chi\rs}\sd g^{-1}
= |g|^2 \ls \chi, \chi H\rs g^{-1}
= \ls \chi, \chi H\rs \Bar{g}
= \ls \chi, g \chi H\rs.
\end{align*}
Observe the following equivalence: for a given $\U\in\pot{KE}$ and $\chi\in H^0(\X,\Sh{U}\sd)$ the potential $\U$ belongs to $\pot{KE,g}$ if and only $g\chi\in H^0(\X,\Sh{U})$.
Hence for $\var{\U}\in\pot{KE,g}$ we only need to require that $-\var\U\sd \chi$ belongs to the image of $\delbar{E^{-1}} -\U\sd$.
By the argument in Part~(A), this is equivalent to $\var \U$ being orthogonal to $\{\ls\chi,\chi\pa\rs\mid\chi\pa\in H^0(\X,\Sh{U})\}$.
The involution of $\pot{KE}$ is unitary in the inner product~\eqref{eq:inner product} since
\[
\Real \overline{g\U\sd g^{-1}}\, g\W\sd g^{-1}
= \Real \Bar{g}^{-1} \qk\U\qk |g|^2 \qk\Bar{\W}\qk g^{-1}
= \Real \U\Bar{\W}
= \Real \Bar{\U}\W.
\]
Thus $\var \U\in\pot{KE,g}$ is orthogonal to the second term on the right hand side of
\[
2\ls\chi,\chi\pa\rs
= \big(\ls\chi,\chi\pa\rs + \ls g^{-1}\chi\pa,g\chi\rs\big) + \big(\ls\chi,\chi\pa\rs - \ls g^{-1}\chi\pa,g\chi\rs\big),
\]
which is anti-symmetric with respect to this involution.
Hence $\var \U\in\pot{KE,g}$ is the first component of a pair which is tangent to $\prewei(\X,E,g)$ at $(\U,\chi)$, if and only if it is orthogonal to the image of~\eqref{eq:singularity2d}.

The analogue of~\eqref{period dimension} is
\begin{gather}\label{period dimension imaginary}\begin{aligned}
H^0(\X,\Sh{U}\sd)&\to\Hom(H_1(\X,\mathbb{Z}),\Imag\qat),\\
\var\chi&\mapsto\left(\gamma\mapsto\int_\gamma \lp\var\chi,g\chi\rp+\lp\chi,g\var\chi\rp \right).
\end{aligned}\end{gather}
This means that for given $\var \U \in \pot{KE,g}$ there a tangent vector to $\wei(\X,E,g)$ whose first component is $\var U$ if and only if there exist an element $(\var \U, \var\chi)$ of the tangent space to $\prewei(\X,E,g)$ that is mapped by~\eqref{eq:periods imagsubspace derivative} to the image of~\eqref{period dimension imaginary}, since the sum of such a pair with an appropriate element in the domain of~\eqref{period dimension imaginary} is tangent to $\wei(\X,E,g)$.
Again the existence of one such pair is equivalent to the condition that~\eqref{eq:periods imagsubspace derivative} maps all such $(\var \U, \var\chi)$ in the tangent space to $\prewei(\X,E,g)$ to the image of~\eqref{period dimension imaginary}, because the difference of two such pairs belongs to $\{0\}\times H^0(\X,\Sh{U}\sd)$.

By $KE^2E^{-1}=KE$ and $\U\in\pot{KE,g}$, there is a map $H^0(\X,\Sh{U}\sd)\to H^0(\X,\Sh{U})$, $\var\chi\mapsto g\var\chi$.
Due to $\lp\var\chi,g\chi\rp=\lp g\var\chi,\chi\rp$ and since $\lp\cdot,\cdot\rp$ is anti-Hermitian, the map~\eqref{period dimension imaginary} indeed takes values in $\Hom(H_1(\X,\mathbb{Z}),\Imag\qat)$.
Moreover, its image consists of the imaginary parts of the image of~\eqref{eq:period 4} with $\psi=g\chi$.
The period pairing~\eqref{eq:pairing} is invariant under $\qat$\=/conjugation of both entries so restricts to a non-degenerate pairing on $\Hom(H_1(\X,\mathbb{Z}),\Imag\qat)$.

Now we claim that the images of~\eqref{eq:translation period 2} and~\eqref{period dimension imaginary} are the orthogonal complements of each other.
By definition of $\con(\X,\chi,g)$ the image of~\eqref{eq:translation period 2} is the intersection of the image of~\eqref{eq:period 3} with $\Hom(H_1(\X,\mathbb{Z}),\Imag\qat)$.
Due to Lemma~\ref{pairing 2} the image of~\eqref{eq:period 3} is the orthogonal complement of the image of~\eqref{eq:period 4}.
By writing the periods in~\eqref{period dimension imaginary} as
\[
\int_\gamma \lp\var\chi,g\chi\rp+\lp\chi,g\var\chi\rp 
= 
2 \int_\gamma \lp\chi,g\var\chi\rp + \int_\gamma \left(\lp\var\chi,g\chi\rp-\lp\chi,g\var\chi\rp \right)
\]
the claim follows since the second integral on the right is real and therefore orthogonal to $\Hom(H_1(\X,\mathbb{Z}),\Imag\qat)$.
Hence in this case $\var \U\in\pot{KE,g}$ is the first component of a pair $(\var \U,\var\chi)$ in the tangent space of $\wei(\X,E,g)$ at $(\U,\chi)$ if and only if
\begin{enumerate}
\item $\var \U$ is orthogonal to the image of~\eqref{eq:singularity2d}, and
\item the intersection image of~\eqref{eq:period 3} and $\Hom(H_1(\X,\mathbb{Z}),\Imag\qat)$ is orthogonal to~\eqref{eq:periods imagsubspace derivative} for any tangent vector to $\prewei(\X,E,g)$ whose first first component is $\var \U$.
\end{enumerate}
The map~\eqref{eq:translation period 2} is the restriction of~\eqref{eq:translation period} to the subspace $\{(\U,\chi,g\chi)\mid(\U,\chi)\in\prewei(\X,E,g)\}\subset\prewei(\X,E)$.
Similarly, the map~\eqref{eq:singularity2c} is the restriction of~\eqref{eq:singularity2b} to the subspace $\{(\chi H,g\chi H)\mid g\chi H\in\con(\X,\chi,g)\}\subset\con(\X,\chi,g\chi)$.
We replace in the remaining part of the proof in Part~(A) the image of~\eqref{eq:translation period} by the image of~\eqref{eq:translation period 2}, the map~\eqref{eq:periods} by~\eqref{eq:periods imagsubspace}, and the map~\eqref{eq:singularity2a} by the map~\eqref{eq:singularity2d}.
Then the modification of the argument for the first statement in Part~(A) proves the first statement in Part~(B).

For the proof of the second statement we recall that the restriction of~\eqref{eq:periods imagsubspace derivative} to vectors of the form $(0,\var\chi)$ has the same image as the analogue~\eqref{period dimension imaginary} of~\eqref{period dimension}.
Moreover, the orthogonal complement of this image is the image of~\eqref{eq:translation period 2}.
The modifications of the arguments for the second and the third statements in Part~(A) prove the analogous statements in Part~(B).
\end{proof}

\begin{lemma}
\label{lem:wei-isothermic}
Let $X$ be a (not necessarily compact) Riemann surface.
\begin{enumeratethm}[label={\upshape(\Alph*)}]
\item
Let $(\U,\chi,\psi)\in \wei(\X,E)$ correspond to an admissible map $F$.
\begin{enumerate}[label={\upshape(\roman*)}]
\item $(\chi G, \psi H)$ belongs to the kernel of~\eqref{eq:singularity2b} if and only if $G = -\Bar{H}$.
\item $F$ is strongly isothermic if and only if~\eqref{eq:singularity2b} is not injective.
\item If the kernel of~\eqref{eq:singularity2b} has dimension $>1$, then $(\delbar{E}-\U^+)(\overline{\qj\chi})^{-1}=0$.
\end{enumerate}

\item
Let $(\U,\chi)\in\wei(\X,E,g)$ correspond to an admissible map $F$ into $\Imag\qat$.
\begin{enumerate}[label={\upshape(\roman*)}]
\item $(g\chi H)$ belongs to the kernel of~\eqref{eq:singularity2c} if and only if $H=-\Bar{H}$.
\item $F$ is strongly isothermic if and only if~\eqref{eq:singularity2c} is not injective.
\item If the kernel of~\eqref{eq:singularity2c} has dimension $>1$, then $\delbar{E}(\overline{\qj\chi})^{-1}=0$.
\end{enumerate}
\end{enumeratethm}
\end{lemma}
\begin{proof}
As in Theorem~\ref{thm:kodaira isothermic characterization} and by applying the definition of $\ls\cdot,\cdot\rs$, any $(\chi G,\psi H)\in\con(\X,\chi,\psi)$ obeys 
\[
0 = \ls\chi G,\psi\rs + \ls\chi,\psi H\rs= \ls\chi (G + \Bar{H}),\psi\rs
\]
if and only if $G = -\Bar{H}$ on the complement of the pullback of the roots of $\psi$.
Analogously any $g\chi H\in\con(\X,\chi,g)$ obeys 
\[
0 = \ls\chi H,g\chi\rs + \ls\chi,g\chi H\rs = \ls\chi(H + \Bar{H}),g\chi \rs
\]
if and only if $H=-\Bar{H}$.
This proves~(i) in Parts~(A) and~(B).

Now we show~(ii) in Part~(A).
Unlike the proof of Theorem~\ref{thm:kod submanifold}, which used Lemma~\ref{lem:strongly isothermic correspondence} to characterize strongly isothermic, in this case the admissible map $G$ itself turns out to be the dual map to $F$.
Take a non-trivial pair $(\psi\pa, \chi\pa) = (\chi G, \psi H)$ in the kernel of~\eqref{eq:singularity2b}.
Then we have Kodaira representations $G=\chi^{-1}\psi\pa$ and $H=\psi^{-1}\chi\pa$ on the complement of the isolated roots of $\psi$ on $\Tilde{\X}$.
Non-triviality of the pair means that $(\psi\pa,\chi\pa)\not\in\{(\chi q,-\psi\Bar{q})\mid q\in\mathbb{H}\}$ respectively, so neither $G$ nor $H$ are not constant, which is part of the definition of admissible.
The left normals of $G$ and $H$ are $\chi^{-1}\qi\chi$ and $\psi^{-1}\qi\psi$ respectively, compare~\eqref{eq:weierstrass-NR}.
But because the pair belongs to the kernel of~\eqref{eq:singularity2b}, by (i) we know $G = -\Bar{H}$, and Example~\ref{eg: - bar F} tells us that the right normal of $G$ is $-\psi^{-1}\qi\psi$.
So the normals of $G$ are the negatives of the normals of $F$.
Although $G$ is defined on $\tilde{\X}$, its derivative $dG$ is single-valued on $\X$.
Hence $F$ is strongly isothermic (Definition~\ref{def:strongly isothermic}).

Vice versa, if $F$ is strongly isothermic, then there exists an admissible map $G$ on $\Tilde{\X}$ with left normal $\chi^{-1}\qi\chi$ and right normal $-\psi^{-1}\qi\psi$.
This implies that $G$ is the quotient $\chi^{-1}\psi\pa$ for another holomorphic section $\psi\pa\in H^0(\Tilde{\X},\tQ{E^{-1},\U\sd})$, where $\tQ{E^{-1},\U\sd}$ is the pullback to $\Tilde{\X}$ of $\Q{E^{-1},\U\sd}$.
Likewise $-\Bar{G}$ is the quotient $\psi^{-1}\chi\pa$ for another holomorphic section $\chi\pa\in H^0(\Tilde{X}\setminus\{\psi=0\},\tQ{KE,\U})$ on the complement of the roots of $\psi$ on $\Tilde{\X}$.
Defining $H = -\Bar{G}$, this pair $(\psi\pa,\chi\pa) = (\chi G, \psi H)$ belongs to $H^0(\Tilde{X}\setminus\{\psi=0\},\tQ{KE,\U}) \times H^0(\Tilde{X},\tQ{E^{-1},\U\sd})$ and the kernel of~\eqref{eq:singularity2b}.
This completes the proof of~(ii) in Part~(A).

For Part~(B)(ii) the argument is analogous to Part~(A).
The left and right normals of $F$ are $-\chi^{-1}\qi\chi$ and $(g\chi)^{-1}\qi(g\chi) = \chi^{-1}\qi\chi$. 
For a non-trivial $\chi\pa$ we have $H = (g\chi)^{-1} \chi\pa$ on the complement of the roots of $g$.
Its left normal is $\chi^{-1}\qi\chi$.
If $\chi\pa$ belongs to the kernel of~\eqref{eq:singularity2c}, then $H = -\Bar{H}$ and so the right normal of $H$ is $-\chi^{-1}\qi\chi$ (again by Example~\ref{eg: - bar F}).
Hence $F$ is strongly isothermic.
Conversely if $F$ is strongly isothermic, then there exists an admissible map $H$ on $\Tilde{\X}$ with left normal $\chi^{-1}\qi\chi$ and right normal $-\chi^{-1}\qi\chi$.
This implies that $H$ is the quotient $(g\chi)^{-1}\chi\pa$ for another holomorphic section $\chi\pa\in H^0(\Tilde{X}\setminus\{\psi=0\},\tQ{KE,\U})$ on the complement of the roots of $g$ on $\Tilde{\X}$.
Further, since $F$ is valued in $\Imag\qat$, we can choose $H$ to be valued in $\Imag\qat$ too, so that $H = -\Bar{H}$.
Hence $\chi\pa = g\chi H$ is in the kernel of~\eqref{eq:singularity2c}, completing the proof of~(ii) in Part~(B).

For the final statement~(iii) in Part~(A), suppose that there are two $(\upsilon\pa_1,\zeta\pa_1)$, $(\upsilon\pa_2,\zeta\pa_2)$ in $H^0(\X,\Sh{V}\sd)\times H^0(\X,\Sh{S})$ with admissible maps $dG_l=\lp(\overline{\qj\chi})^{-1},\upsilon\pa_l\rp$.
If the right normals $(\upsilon\pa_l)^{-1}\qi \upsilon\pa_l$~\eqref{eq:weierstrass-NR} of both $G_l$ are equal to $-\psi^{-1}\qi\psi$, then the arguments in the proof of Theorem~\ref{thm:isothermic local characterization} shows the vanishing of the part of the potential of $(\overline{\qj\chi})^{-1}$ that anti-commutes with the complex structure.
Suppose in Part~(B) that there are two $\zeta\pa_1,\zeta\pa_2\in H^0(\X,\Sh{S})$ with admissible maps $dH_l=\ls(\overline{\qj\chi})^{-1},g^{-1}\zeta\pa_l\rs$ having the same right normals, then this implies $\delbar{E}(\overline{\qj\chi})^{-1}=0$, since additionally by~\eqref{def:real potential} the commuting part of the potential $\delbar{E}(\overline{\qj\chi})^{-1}\overline{\qj\chi}$ of $(\overline{\qj\chi})^{-1}$ is zero.
\end{proof}

If the Weierstraß triple $(\U,\chi,\psi)\in\wei(\X,E)$ corresponds to a strongly isothermic map $F$ then we know that~\eqref{eq:singularity2b} has a non-trivial kernel.
Moreover if $(\V,\upsilon,\phi)\in\kodaira(\X,E)$ is the Kodaira triple corresponding to this strongly isothermic map $F$, then the real dimensions of the maps~\eqref{eq:singularity1} and~\eqref{eq:singularity2b} are equal.
This is because, by Lemma~\ref{lem:strongly isothermic correspondence} and the proof of~(ii) in the foregoing lemma, the elements of the kernels of both maps~\eqref{eq:singularity1} and~\eqref{eq:singularity2b} are in one-to-one correspondence to dual maps to $F$ on $\Tilde{\X}$.

\begin{example}[Catenoid]
\label{eg:catenoid isothermic weierstrass}
\index{Catenoid!Isothermicity}
For the catenoid we are in $\mathbb{R}^3$ so we can apply Part~(B).
In fact, because $\U \equiv 0$ there is a simple trick to produce an element of the kernel of~\eqref{eq:singularity2c}, namely $\chi\pa = \qi g \chi = \qi (1 + z^{-1}\qk)$.
The corresponding dual surface is
\begin{align*}
H
= (g\chi)^{-1} \chi^\ast
= \chi^{-1}\qi\chi
= (1 + z^{-1}\qk)^{-1} \qi (1 + z^{-1}\qk)
= (z + \qk)^{-1} \qi (z + \qk),
\end{align*}
which is the round sphere.
This is perhaps more enlightening than Example~\ref{eg:catenoid isothermic kernel}, which used the Kodaira representation.
\end{example}

In the following lemma, we give a sufficient criterion for $\prewei(\X,E)$ and $\prewei(\X,E,g)$ to be a Banach submanifold, and we give an additional sufficient criterion for $\wei(\X,E)$ and $\wei(\X,E,g)$ to be a Banach submanifold of $\prewei(\X,E)$ and $\prewei(\X,E,g)$, respectively.
Both criteria together show that $\wei(\X,E)$ and $\wei(\X,E,g)$ are Banach manifolds at elements corresponding to not strongly isothermic maps, in analogy to $\kodaira(\X,E)$.

\begin{lemma}\label{lem:wei-non-isothermic}
Let $\X$ be a compact Riemann surface.
\begin{enumeratethm}[label={\upshape(\Alph*)}]
\item 
At $(\U,\chi,\psi)$ with injective~\eqref{eq:singularity2a} $\prewei(\X,E)$ is a Banach submanifold of $\pot{KE}\times\sob{1,p}{E^{-1}}(\X)\times\sob{1,p}{KE}(\X)$.
At $(\U,\chi,\psi)$ that are not strongly isothermic the map~\eqref{eq:periods extension map} is submersive and $\wei(\X,E)$ a Banach manifold.

\item
At $(\U,\chi)\in\prewei(\X,E,g)$ with injective~\eqref{eq:singularity2d} $\prewei(\X,E,g)$ is a Banach submanifold of $\pot{KE,g}\times\sob{1,p}{KE}(\X)$.
At $(\U,\chi)$ that are not strongly isothermic the extension of~\eqref{eq:periods imagsubspace} is submersive and $\wei(\X,E,g)$ a Banach manifold.
\end{enumeratethm}
\end{lemma}
\begin{proof}
First we show that the subset given in the first statements of the parts~(A) and~(B) respectively is a Banach submanifold.
For this purpose we apply the arguments in the proof of Theorem~\ref{thm:kod submanifold}, where the map $\Theta$~\eqref{immersion kodaira} is replaced by the map $\Psi$~\eqref{eq:prewei map} in Part~(A) and by the map $\Phi$~\eqref{eq:prewei map 2} in Part~(B).
Then in Part~(A) $\tfrac{\partial \Theta}{\partial \upsilon}$ is replaced by $\tfrac{\partial \Psi}{\partial \chi}:\var\chi\mapsto(\delbar{E^{-1}} -\U\sd)\var\chi$ and $\tfrac{\partial \Theta}{\partial\phi}$ is replaced by $\tfrac{\partial \Psi}{\partial\psi}:\var\psi\mapsto(\delbar{KE}-\U)\var\psi$.
In Part~(B) the combined derivatives $(\frac{\partial\Theta}{\partial\upsilon},\frac{\partial\Theta}{\partial\phi})$ are replaced by $\frac{\partial\Phi}{\partial\chi}$.
In Part~(A) the kernel $\Spa{K}$ is replaced by the subspace
\begin{align*}
\bigl\{(\var \U,\var\chi,\var\psi)\;\big|\;(\var\U\sd\chi,\var \U\psi)&=((\delbar{E^{-1}} -\U\sd)\var \chi,(\delbar{KE}-\U)\var\psi)\bigr\}
\end{align*}
of $\pot{KE}\times\sob{1,p}{E^{-1}}(\X)\times\sob{1,p}{KE}(\X)$ and in Part~(B) by the subspace
\begin{align*}
\{(\var \U,\var\chi)\mid\var\U\sd \chi=(\delbar{E^{-1}} -\U\sd)\var\chi\Longleftrightarrow\var \U g\chi=(\delbar{KE}-\U)g\var\chi\}
\end{align*}
of $\pot{KE,g}\times\sob{1,p}{E^{-1}}(\X)$.
The Fredholm operator~\eqref{eq:W2-L2-map} is replaced in Part~(A) by
\begin{align*}
\sob{1,p}{E^{-1}}(\X)\times\sob{1,p}{KE}(\X)&\to\forms{0,1}{}\ban{p}{E^{-1}}(\X)\times \forms{0,1}{}\ban{p}{KE}(\X)\\(\var\chi,\var\psi)&\mapsto\bigl((\delbar{E^{-1}} -\U\sd)\var\chi\,,\,(\delbar{KE}-\U)\var\psi\bigr)
\end{align*}
and in Part~(B) by 
\[
\sob{1,p}{E^{-1}}(\X)\to \forms{0,1}{}\ban{p}{E^{-1}}(\X),\;\var\chi\mapsto(\delbar{E^{-1}} -\U\sd)\var\chi.
\]

Finally, $-\tfrac{\partial\Theta}{\partial V}$ (compare~\eqref{eq:deltaW-map}) is in Part~(A) replaced by the map $-\tfrac{\partial \Psi}{\partial \U}:\var \U\mapsto(\var\U\sd\,\chi,\var \U\,\psi)$, such that the pairing of the image with $(\psi\pa,\chi\pa)\in H^0(\X,\Sh{U}\sd)\times H^0(\X,\Sh{U})$ is $\langle\ls\psi\pa,\psi\rs + \ls\chi,\chi\pa\rs, \var \U\rangle$ and vanishes exactly on the subspace $\{(\chi q, -\psi \Bar{q})\mid q\in\qat\}$.
In Part~(B) $-\tfrac{\partial\Theta}{\partial V}$ is replaced by $-\tfrac{\partial \Phi}{\partial \U}:\var \U\mapsto \var\U\sd$, such that the pairing of the image with $\chi\pa\in H^0(\X,\Sh{U})$ is $\langle\ls g^{-1}\chi\pa,g\chi\rs + \ls\chi,\chi\pa\rs, \var \U\rangle$ and vanishes exactly on the subspace $\{g\chi q \mid q\in\qat\}$.
The arguments from the proof of Theorem~\ref{thm:kod submanifold} show that the derivatives of $\Psi$ and $\Phi$ are surjective if and only if~\eqref{eq:singularity2a} and~\eqref{eq:singularity2d} (the replacements of~\eqref{eq:singularity1}) are injective.
Again Lemma~\ref{lem:implicit-variant} applies and shows that the subsets of $\prewei(\X,E)$ and $\prewei(\X,E,g)$ given in the first statements are Banach submanifolds.

It remains to prove that at triples $(\U,\chi,\psi)\in\wei(\X,E)$ and pairs $(\U,\chi)\in\wei(\X,E,g)$ that are not strongly isothermic, the extended period maps~\eqref{eq:periods extension map} and the extension of~\eqref{eq:periods imagsubspace} are submersive.
This follows in both parts~(A) and~(B) from the second statements of Lemma~\ref{lem:wei-necessary condition}: For injective~\eqref{eq:singularity2b} and~\eqref{eq:singularity2c} the orthogonal complement of the image of~\eqref{eq:periods derivative} and~\eqref{eq:periods imagsubspace derivative} is trivial, respectively.
But these are the restrictions of the derivative of the extension~\eqref{eq:periods extension map} and of the extension of~\eqref{eq:periods imagsubspace}, and so the latter are also surjective.
\end{proof}
Since we have just shown that $\wei(\X,E)$ and $\wei(\X,E,g)$ are Banach submanifolds away from the strongly isothermic maps, the natural question is how these spaces look at the strongly isothermic maps.
The following theorem gives the answer.

\begin{theorem}\label{thm:weierstrass isothermic characterization}
\index{Isothermicity!Singularity of $\wei(\X,E)$}
Let $\X$ be a compact Riemann surface.
\begin{enumeratethm}[label={\upshape(\Alph*)}]
\item
Let $(\U,\chi,\psi)\in\wei(\X,E)$ have a kernel of~\eqref{eq:singularity2b} of real dimension $\le1$.
Then at $(\U,\chi,\psi)$ the tangent cone spans the tangent space and there is a smooth path in $\wei(\X,E)$ through $(\U,\chi,\psi)$ that meets the set of strongly isothermic surfaces only at $(\U,\chi,\psi)$.
Moreover, the following statements are equivalent: 
\begin{enumerate}[label={\upshape(\roman*)}]
\item The admissible map with $dF = \lp \chi,\psi \rp$ is strongly isothermic.
\item The kernel of the map~\eqref{eq:singularity2b} has real dimension one.
\item The tangent cone of $\wei(\X,E)$ at $(\U,\chi,\psi)$ is a proper subset of the tangent space.
\item The subset $\wei(\X,E)\subset\pot{KE}\times H^0(\X,\sob{1,p}{E^{-1}})\times H^0(\X,\sob{1,p}{KE})$ is not a submanifold at $(\U,\chi,\psi)$.
\end{enumerate}
\item
Let $(\U,\chi)\in\wei(\X,E,g)$ have a kernel of~\eqref{eq:singularity2c} of real dimension $\le1$.
Then at $(\U,\chi)$ the tangent cone spans the tangent space and there is a smooth path in $\wei(\X,E,g)$ through $(\U,\chi)$ that meets the set of strongly isothermic surfaces only at $(\U,\chi)$.
Moreover, the following statements are equivalent:
\begin{enumerate}[label={\upshape(\roman*)}]
\item The admissible map with $dF = \lp \chi,\psi \rp$ is strongly isothermic.
\item The kernel of the map~\eqref{eq:singularity2c} has real dimension one.
\item The tangent cone of $\wei(\X,E,g)$ at $(\U,\chi)$ is a proper subset of the tangent space.
\item $\wei(\X,E,g)$ is not a submanifold of $\pot{KE,g}\times H^0(\X,\sob{1,p}{E^{-1}})$ at $(\U,\chi)$.
\end{enumerate}
\end{enumeratethm}
\end{theorem}

This is completely analogous to Theorem~\ref{thm:kodaira isothermic characterization} and we approach the proof through the same method.
We show in Lemma~\ref{lem:weierstrass tangent cone} that $\wei(\X,E)$ is locally the zero set of a smooth real function $\lambda$ on a submanifold, whose derivative vanishes at the isothermic elements but whose second derivative is indefinite.
The indefiniteness is shown for a plane in Lemma~\ref{second variation plane weierstrass}, and for the general case through a blow-up and perturbation argument in Lemma~\ref{lem:indefinite weierstrass}.
The principal difference between the proof in the Kodaira case and the Weierstraß case is that $\wei(\X,E)$ is defined through two conditions; the holomorphicity of the sections and the period conditions.
Thus, roughly speaking, there are two chances for it to fail to be a manifold.
But because the dimension of the kernel of the map~\eqref{eq:singularity2b} has real dimension one, which represents the obstacles to being a manifold, only one condition or the other can fail, not both simultaneously.
This leads us to actually introduce two functions $\lambda_1$ and $\lambda_2$ depending on the obstacle.

\begin{lemma}
\label{lem:weierstrass tangent cone}
Let $\X$ be a compact Riemann surface.
\begin{enumeratethm}[label={\upshape(\Alph*)}]
\item
Suppose the map~\eqref{eq:singularity2b} of any $(\U,\chi,\psi)\in\wei(\X,E)$ has a real one-dimen\-sional kernel.
Then there exist an open subset $O$ of the Banach space $\pot{KE}\times\sob{1,p}{E^{-1}}(\X)\times\sob{1,p}{KE}(\X)$, a submanifold $\Spa{U}\ni(\U,\chi,\psi)$ of $O$ and a smooth function $\lambda:\Spa{U} \to \mathbb{R}$ with the following properties:
\begin{enumerate}[label={\upshape(\roman*)}]
\item $\wei(\X,E) \cap O = \lambda^{-1}[\{0\}]$,
\item $\lambda'(\U,\chi,\psi)=0$,
\item $\lambda''(\U,\chi,\psi)$ is indefinite.
\end{enumerate}

\item 
Suppose that at $(\U,\chi)\in\wei(\X,E,g)$ the kernel of the map~\eqref{eq:singularity2c} is one-dimensional.
Then there exists an open subset $O\subset\pot{KE}\times\sob{1,p}{E^{-1}}(\X)$, a submanifold $\Spa{U}\ni(\U,\chi)$ of $O$ and a smooth function $\lambda:\Spa{U} \to \mathbb{R}$ with $\wei(\X,E,g) \cap O = \lambda^{-1}[\{0\}]$ and properties (ii) and (iii) above.
\end{enumeratethm}
\end{lemma}

The Properties~(i) and (ii) together imply that the tangent spaces of $\Spa{U}$ and $\wei(\X,E)$ (respectively $\wei(\X,E,r)$) coincide at $(\U,\chi, \psi)$ (respectively $(\U,\chi)$).
By Lemma~\ref{lem:tangent cone spans} the Property~(iii) further ensures that the tangent cone spans the tangent space.

\begin{proof}
We prove (A) and comment throughout any additional steps to prove (B).
Let $(\chi G,\psi H)\in\con(\X,\chi,\psi)$ span the kernel of~\eqref{eq:singularity2b}.
There are two cases to consider: whether $(\psi\pa,\chi\pa)\in H^0(\X,\Sh{U}\sd)\times H^0(\X,\Sh{U})\subset\con(\X,\chi,\psi)$ or not.
If it does then the kernels of both maps~\eqref{eq:singularity2b} and~\eqref{eq:singularity2a} are one-dimensional.
If it doesn't, then the kernel of~\eqref{eq:singularity2a} is trivial.

For the first case, we apply similar arguments as in the proof of Lemma~\ref{lem:isothermic1} and mainly explain the differences.
As $(\psi\pa,\chi\pa)\in H^0(\X,\Sh{U}\sd)\times H^0(\X,\Sh{U})$, which is contained in $\con(\X,\chi,\psi)$, we choose $(\alpha, \beta) \in H^0(\X,\forms{0,1}{}\ban{p}{E^{-1}}) \times H^0(\X,\forms{0,1}{}\ban{p}{KE})$ with
\begin{equation}\label{eq:dual vector 2}
\int_\X \Real \left(\lh\psi\pa,\alpha\rh + \lh\chi\pa,\beta\rh \right) \neq 0.
\end{equation}
A modification of the first part of the proof of Lemma~\ref{lem:wei-non-isothermic} shows that the zero set of the composition of the map~\eqref{eq:prewei map} with the natural projection onto the quotient space $H^0(\X,\forms{0,1}{}\ban{p}{E^{-1}})\times H^0(\X,\forms{0,1}{}\ban{p}{KE})/\mathbb{R}(\alpha,\beta)$ is a submanifold $\Spa{U}_1'$ of an open neighborhood of $(\U,\chi,\psi)$.

We shall prove now that the tangent space of $\Spa{U}_1'$ at $(\U,\chi,\psi)$ coincides with the tangent space of $\prewei(\X,E)$.
By the chain rule, a vector lies in $T_{(\U,\chi,\psi)}\Spa{U}_1'$ if and only if it is mapped by the derivative of~\eqref{eq:prewei map} to the kernel of the projection, which is $\mathbb{R}(\alpha,\beta)$.
The pairing (from Lemma~\ref{lem:11 pairing properties}(iii)) of $(\psi\pa,\chi\pa)$ with the image of the derivative of~\eqref{eq:prewei map} along $(\var \U, \var \chi,\var \psi)$ is equal to 
\begin{align*}
&\Real \int_X \lh\psi\pa,-\var\U\psi + (\delbar{KE} - \U)\var\psi\rh
+ \lh\chi\pa, -\var\U\sd \chi + (\delbar{E^{-1}}-\U\sd)\var\chi\rh \\
&= - \langle \ls\psi\pa.\psi\rs + \langle\ls\chi, \chi\pa\rs, \var\U\rangle \\
&\qquad + \Real \int_X \overline{\lh\var\psi, (\delbar{E^{-1}}-\U\sd)\psi\pa\rh} + \overline{\lh\var\chi,(\delbar{KE} - \U)\chi\pa\rh}
= 0.
\end{align*}
In this calculation we have used the Lemma~\ref{lem:01 pairing properties}(iii) and~(iv), as well as integration by parts~\eqref{eq:integration by parts}.
The expression vanishes because $(\psi\pa,\chi\pa)$ belongs to the kernel of~\eqref{eq:singularity2a} and $(\psi\pa,\chi\pa) \in H^0(\X,\Sh{U}\sd)\times H^0(\X,\Sh{U})$ are holomorphic.
This shows that the image of the derivative of~\eqref{eq:prewei map} is contained in the orthogonal complement of $(\psi\pa,\chi\pa)$.
Our choice of $(\alpha,\beta)$ with~\eqref{eq:dual vector 2} means that it lies outside the image of the derivative of~\eqref{eq:prewei map}.
In other words, a vector is mapped to $\mathbb{R}(\alpha,\beta)$ by the derivative of~\eqref{eq:prewei map} if and only if it belongs to the kernel of the derivative of~\eqref{eq:prewei map}, which is tangent space of $\prewei(\X,E)$.

Since the kernels of both maps~\eqref{eq:singularity2a} and~\eqref{eq:singularity2b} are one-dimensional in this case, the third statement of Lemma~\ref{lem:wei-necessary condition} implies that the period map~\eqref{eq:periods extension map} is submersive on an open neighborhood of $(\U,\chi,\psi)$.
Therefore the zero set of the restriction of~\eqref{eq:periods extension map} to $\Spa{U}_1'$ defines a Banach submanifold $\Spa{U}_1$ of $\pot{KE}\times H^0(\X,\sob{1,p}{E^{-1}})\times H^0(\X,\sob{1,p}{KE})$.
Here we see the idea discussed before the lemma that only one obstacle can come into play at a time; $\prewei(\X,E)$ might not be a submanifold at this point since~\eqref{eq:prewei map} is not a submersion, but given a suitable submanifold $\Spa{U}_1'$ containing $\prewei(\X,E)$ there is no difficulty finding a further submanifold $\Spa{U}_1$ whose Weierstraß data has vanishing periods.

Now we define the function $\lambda_1:\Spa{U}_1\to\mathbb{R}$ in analogy to~\eqref{eqn:kod lambda}.
Recall $(\psi\pa,\chi\pa) = (\chi G, \psi H) \in \con(\X,\chi,\psi)$.
Let $\lambda_1$ at $(\U', \chi',\psi') \in \Spa{U}_1$ be the pairing of $(\chi' G,\psi' H)$ with the image of~\eqref{eq:prewei map}.
In symbols
\begin{gather}
\label{eq:def lambda weierstrass 1a}
\lambda_1:(\U'\hspace{-1.5mm},\chi'\hspace{-1.5mm},\psi')\mapsto
\Real\int_\X \!\big\lh\chi'G,\hspace{-.5mm}(\delbar{KE}\!-\!\U')\psi'\big\rh\! + \!\big\lh\psi'H,\hspace{-.5mm}(\delbar{E^{-1}} \!-\!{\U'}\sd)\chi'\big\rh
\end{gather}
We argue now that $\lambda_1$ has Property (i): the zero level set of this function is $\wei(\X,E)$.
The construction of $\U_1$ guarantees that we are in the kernel of~\eqref{eq:periods extension map}.
On $\U_1$, the image of~\eqref{eq:prewei map} is a multiple of $(\alpha,\beta)$.
By the choice~\eqref{eq:dual vector 2}, the pairing of $(\chi' G, \psi' H)$ with $(\alpha,\beta)$ is nonzero at $(\U,\chi,\psi)$ and so non-zero for sufficiently small $\U_1'$ and $\U_1$.
Therefore $\lambda_1$ vanishes if and only if~\eqref{eq:prewei map} does.
Hence $\lambda_1^{-1}[\{0\}]$ is a neighborhood of $\wei(\X,E)$.

Before we compute the derivatives of $\lambda_1$ so that we can prove Properties~(ii) and~(iii), let us bring $\lambda_1$ into a simpler form.
This is not strictly necessary, but shortens the arguments, especially in the second case below.
The reason such a simplification is possible is that we are in the strongly isothermic case, due to Lemma~\ref{lem:wei-isothermic}, so we know that $G = - \Bar{H}$.
Therefore, using Lemma~\ref{lem:11 pairing properties}(iv), we have
\begin{align*}
\lambda_1
&= \Real\int_\X \Bar{G}\big\lh\chi',\hspace{-.5mm}(\delbar{KE}\!-\!\U')\psi'\big\rh\! + \overline{\Bar{H}\big\lh\psi',\hspace{-.5mm}(\delbar{E^{-1}} \!-\!{\U'}\sd)\chi'\big\rh} \\
&= \Real\int_\X \Bar{G} \,d\lp\chi',\psi'\rp.
\end{align*}
Thus the first derivative is
\begin{align*}
\lambda_1'(\U,\chi,\psi)(\var \U,\var\chi,\var\psi)
&= \Real\int_\X \Bar{G} \,d\Big[\lp\var\chi,\psi\rp + \lp\chi,\var\psi\rp\Big].
\end{align*}
This is zero, Property (ii), because above we have argued that $(\var \U,\var\chi,\var\psi)$ is tangent to $\prewei(\X,E)$ and thus Lemma~\ref{lem:period welldefined} tells us that the form in the square bracket is closed.

We apply Lemma~\ref{lem:tangent cone spans} to obtain the second derivative,
\begin{align*}
\lambda_1''(\U,\chi,\psi)\big((\var \U,\var\chi,\var\psi),(\var \U,\var\chi,\var\psi)\big)
&= \Real\int_\X \Bar{G} \,2d\lp\var\chi,\var\psi\rp.
\end{align*}
It will be useful to have a form of this equation that doesn't involve $d$.
Applying Lemma~\ref{lem:11 pairing properties}(iv) again gives
\begin{align*}
\lambda_1''
&= 2 \Real\int_\X \Bar{G}\big[ \lh\var\chi,\hspace{-.5mm}(\delbar{KE}\!-\!\U)\var\psi\rh\! - \overline{\lh\var\psi,\hspace{-.5mm}(\delbar{E^{-1}} \!-\!\U\sd)\var\chi\rh} \big] \\
&= - 2 \Real\int_\X \lh\var\chi,\var\U\psi H\rh + \lh\chi G,\var\U\var\psi\rh.
\labelthis{eqn:lambda1''}
\end{align*}
We will prove that this is indefinite, Property (iii), in Lemma~\ref{lem:indefinite weierstrass}.





To prove Part~(B) for the analogous case, namely that~\eqref{eq:singularity2c} has a one dimensional kernel spanned by $\chi\pa = g\chi H$, we only need to choose $\beta \in H^0(\X,\forms{0,1}{}\ban{p}{KE})$ such that the pairing with $\chi\pa$ is non-zero.
The argument defining the submanifolds $\Spa{U}'_1$ and $\Spa{U}_1$ shows that on $\Spa{U}_1$ the image of~\eqref{eq:prewei map 2} lies in the span of $\beta$ and~\eqref{eq:periods imagsubspace} vanishes.
Because
\[
\lh\chi' H,(\delbar{KE} \!-\!\U')(g\chi')\rh 
= \lh \chi' H, g(\delbar{E^{-1}}-{\U'}\sd)\chi' \rh
= \lh g\chi' H, (\delbar{E^{-1}}-{\U'}\sd)\chi' \rh
\]
the definition of $\lambda_1(\U',\chi',g\chi')$ reduces to twice the pairing of the image of~\eqref{eq:prewei map 2} with $\chi\pa$.
Thus $\lambda_1$ is zero if and only if $(\U',\chi') \in \wei(\X,E,g)$.
Similar formulas hold for $\lambda_1'$ and $\lambda_1''$ with the standard substitutions.
In particular
\begin{align*}
\lambda_1''(\U,\chi)\big((\var \U, \var\chi),(\var \U, \var\chi)\big)
= - 4\Real \int_\X \lh \var\chi, \var\U g\chi H\rh.
\end{align*}

In the second case $(\chi G, \psi H)\in\con(\X,\chi,\psi)\setminus\big(H^0(\X,\Sh{U}\sd)\times H^0(\X,\Sh{U})\big)$ so~\eqref{eq:singularity2a} is injective whereas~\eqref{eq:singularity2b} has a one-dimensional kernel.
Hence by Lemma~\ref{lem:wei-non-isothermic} $\prewei(\X,E)$ is locally at $(\U,\chi,\psi)$ a Banach submanifold of $\pot{KE}\times\sob{1,p}{E^{-1}}(\X)\times\sob{1,p}{KE}(\X)$.
By the second statement of Lemma~\ref{lem:wei-necessary condition}, $\per dG$ spans the orthogonal complement of~\eqref{eq:periods derivative} at $(\U,\chi,\psi)$.
Let $P$ be any projection from $\Hom(H_1(\X,\mathbb{Z}),\qat)$ onto the orthogonal complement of $\per dG$.
Concretely, we may choose any $\beta \in \Hom(H_1(\X,\mathbb{Z}),\qat)$ that gives $1$ when paired with $\per dG$. 
Then $P_\beta(\alpha) = \alpha - \langle\per dG,\alpha\rangle\beta$ is such a projection.
Then there exists an open neighborhood $\Spa{U}_2'$ of $(\U,\chi,\psi)$ in $\prewei(\X,E)$ such that the composition of~\eqref{eq:periods} with $P_\beta$ is a submersion from this neighborhood onto the orthogonal complement of $\per dG$.
By the implicit function theorem the preimage $\Spa{U}_2$ of $\{0\}$ is a submanifold of $\Spa{U}_2'$.
Applying Lemma~\ref{lem:period pairing}, $\lambda_2 : \Spa{U}_2 \to \mathbb{R}$ is given by the following formula
\begin{align}\label{def lambda 2}
(\U',\chi',\psi')&\!\mapsto
\langle \per dG, \per \lp\chi',\psi'\rp \rangle
= \Real \int_\X d\Bar{G}\wedge \lp\chi',\psi'\rp.
\end{align}
The first derivative at $(\U,\chi,\psi)$ is 
\begin{align*}
\lambda_2'(\U,\chi,\psi)(\var \U,\var\chi,\var\psi)
&= \langle \per dG, \per \Big[\lp\var\chi,\psi\rp + \lp\chi,\var\psi\rp\Big],
\end{align*}
which vanishes since $\per dG$ spans the orthogonal complement of~\eqref{eq:periods derivative} at $(\U,\chi,\psi)$.
Therefore we have established Property~(ii) in this case.

For the second derivative at $(\U,\chi,\psi)$ we have
\begin{align*}
\lambda_2''(\U,\chi,\psi)\big((\var \U,\var\chi,\var\psi),(\var \U,\var\chi,\var\psi)\big)
&= \Real\int_\X d\Bar{G} \wedge 2\lp\var\chi,\var\psi\rp.
\end{align*}
As in the first case, we would like to calculate the $d$ away.
This requires an integration by parts, but $G$ is not single-valued on $\X$.
Let $\hat{X}$ be a fundamental domain of $\X$ inside its universal cover whose boundary is composed of curves realizing a canonical basis of $H_1(\X,\mathbb{Z})$, compare Lemma~\ref{lem:period pairing}.
Then
\begin{align*}
\lambda_2''
&= 2\Real\int_{\hat{\X}} d\Big[\Bar{G} \lp\var\chi,\var\psi\rp \Big] - \Bar{G} d \lp\var\chi,\var\psi\rp \\
&= 2\Real\int_{\hat{\X}} d\Big[\Bar{G} \lp\var\chi,\var\psi\rp \Big] + \lh\var\chi,\var\U\psi H\rh + \lh\chi G,\var\U\var\psi\rh.
\labelthis{eqn:lambda2''}
\end{align*}
Notice that this is essentially Equation~\eqref{eqn:lambda1''} with a sign change and an additional `boundary term'.
We will prove Property (iii) in Lemma~\ref{lem:indefinite weierstrass} as well.

Finally, to prove Part~(B) for this second case we observe that $\prewei(\X,E,g)$ is a submanifold and we can define $\beta$ and $\Spa{U}_2$ with~\eqref{eq:periods imagsubspace} and~\eqref{eq:translation period 2} in the place of~\eqref{eq:periods} and~\eqref{eq:translation period} because the relevant facts were proved in Lemma~\ref{lem:wei-necessary condition}(B) and are identical to Part~(A).
The rest follows likewise.
In particular, the second derivative $\lambda_2''(\U,\chi)\big((\var\U,\var\chi),(\var\U,\var\chi)\big)$ has the formula
\begin{equation}
\label{eqn:lambda2'' real}
2 \Real\int_{\Hat{\X}}d \left[ -H\lp\var \chi , g\var\chi \rp\right] + 4\Real \int_\X \lh \var\chi, \var\U g\chi H\rh.
\qedhere
\end{equation}
\end{proof}
In the following lemma we define an indefinite bilinear form $\lambda_\infty''$ for the plane, which will be shown to be a limiting case of $\lambda''$.
We only consider $\lambda_\infty''$ as a bilinear form and do not care whether it is the Hessian of some function.
\begin{lemma}
\label{second variation plane weierstrass}
\phantom{boo}
\begin{enumeratethm}[label={\upshape(\Alph*)}]
\item 
Let $\chi_\infty, \psi_\infty, \psi\pa_\infty,\chi\pa_\infty\in\qat\setminus\{0\}$ with $\psi\pa_\infty\Bar{\psi}_\infty=-\chi_\infty\Bar{\phi}\pa_\infty\in\qj\mathbb{C}\subset\qat$.
Consider the triple $(0,\chi_\infty,\psi_\infty) \in \wei(\mathbb{D},\unity)$ corresponding to the flat plane $z_\infty\mapsto \Bar{\chi}_\infty\qj z_\infty\psi_\infty$ on $z_\infty\in\mathbb{D}$.
For $f\in C^\infty_0(\mathbb{D},\mathbb{C})$ the variation $d\Bar{z}_\infty\var \U_\infty = d\Bar{z}_\infty\qj\barpartial \partial f$ can be supplemented with variations with compact support $\var \chi_\infty$ and $\var \psi_\infty$ to a tangent vector $(d\Bar{z}_\infty\var \U_\infty,\var\chi_\infty,\var\psi_\infty)$ of $\wei(\mathbb{D},\unity)$ at the triple $(0,\chi_\infty,\psi_\infty)$.
The following quadratic form is indefinite:
\begin{align*}
\lambda_\infty'' 
&:= \Real \int_{\mathbb{D}} \lh z_\infty \psi\pa_\infty,d\Bar{z}_\infty\var \U_\infty\var\psi_\infty\rh
+ \lh z_\infty\chi\pa_\infty,(d\Bar{z}_\infty\var \U_\infty)\sd \qk \var\chi_\infty\rh
\end{align*}

\item
For $g \in \mathbb{C}\setminus \{0\}$ consider the special case with $\psi_\infty = g \chi_\infty$ and $\psi\pa_\infty = g^{-1}\chi\pa_\infty$.
Choose $f$ to be valued in $g^{-1}\qi\mathbb{R}$.
Then the variations constructed in Part~(A) additionally obey $d\Bar{z}_\infty\var \U_\infty \in \pot{\unity,g}$, $\var\psi_\infty = g\var\chi_\infty$, and $(d\Bar{z}_\infty\var \U_\infty, \var \chi_\infty)$ is tangent to $\wei(\mathbb{D},\unity,g)$ at $(0,\chi_\infty)$.
The quadratic form $\lambda''_\infty$ is also indefinite on this restricted set of $f$.
\end{enumeratethm}
\end{lemma}
\begin{proof}
Because the disc $\mathbb{D}$ is simply connected, there is no difference between $\wei(\mathbb{D},\unity)$ and $\prewei(\mathbb{D},\unity)$.
Recall $\prewei(\mathbb{D},\unity)$ describes admissible maps $\mathbb{D}\to \qat$ whose derivative $dF = \lp\chi,\psi\rp$ is a pairing of holomorphic sections of the trivial $\qat$\=/line bundle $\unity_\qat$ on the unit disc $\mathbb{D}$.
This space contains the triple $(\U,\chi,\psi) = (0,\chi_\infty, \psi_\infty)$ of constant sections, which corresponds to the admissible map $z_\infty\mapsto \Bar{\chi}_\infty\qj z_\infty\psi_\infty$ and describes a flat plane.
By Lemma~\ref{lem:wei-isothermic} we can describe the dual map by a non-trivial element of the kernel of~\eqref{eq:singularity2b}.
This is equivalent to $\qk (\psi\pa_\infty\Bar{\psi}_\infty + \chi_\infty \Bar{\chi}\pa_\infty) = 0$, which holds for the constants $\psi\pa_\infty,\chi\pa_\infty$ in the statement of this lemma.
We may change the coordinate $z_\infty$ such that $\psi\pa_\infty\Bar{\psi}_\infty=-\chi_\infty\Bar{\chi}\pa_\infty\in\qk\mathbb{R}$.
For $(d\Bar{z}_\infty\var \U_\infty, \var \chi_\infty,\var\psi_\infty)$ to be tangential to $\prewei(\mathbb{D},\unity)$ the variations $\var\chi_\infty$ and $\var\psi_\infty$ must lie in the kernel of the derivative of Equation~\eqref{eq:prewei map}.
Letting $d\Bar{z}_\infty\var \U_\infty=d\Bar{z}_\infty\qj \barpartial \partial f$ with $f\in C^\infty_0(\mathbb{D},\mathbb{C})$ this yields
\begin{align*}
\qj\barpartial \partial f\psi_\infty&=\barpartial \var\psi_\infty&
\qk\overline{\qj\barpartial \partial f}\qk \chi_\infty&=\barpartial \var\chi_\infty\\
\var\psi_\infty&=\qj\barpartial f \psi_\infty &
\var\chi_\infty&=-\qj \barpartial \Bar{f}\chi_\infty,
\end{align*}
since the variations $\var \U_\infty, \var\chi_\infty, \var\psi_\infty$ are functions with compact support (not constants) and so these $\barpartial$\=/equations have unique solutions.

Substituting these into $\lambda_\infty''$, we obtain
\begin{align*}
\lambda''_\infty
&= \Real \int_{\mathbb{D}} \lh z_\infty \psi\pa_\infty,d\Bar{z}_\infty\var \U_\infty\var\psi_\infty\rh + \lh z_\infty\chi\pa_\infty,(d\Bar{z}_\infty\var\U_\infty)\sd \var\chi_\infty\rh \\
&= \Real \int_{\mathbb{D}} \Bar{\psi}\pa_\infty\Bar{z}_\infty (-2\qk) \qj \barpartial \partial f \, \qj\barpartial f \psi_\infty 
+ \Bar{\chi}\pa_\infty \Bar{z}_\infty (-2\qk) \,\qk\overline{\qj \barpartial \partial f} \qk (-\qj \barpartial \Bar{f}\chi_\infty) \;\dmu \\
&= \Real \int_{\mathbb{D}} \psi_\infty\Bar{\psi}\pa_\infty 2\qk z_\infty \partial \barpartial \Bar{f} \, \barpartial f 
+ \chi_\infty\Bar{\chi}\pa_\infty (-2\qk) z_\infty \qk \barpartial \partial \Bar{f} \qk \barpartial \Bar{f} \;\dmu \\
&= \Real \int_{\mathbb{D}} (-\psi\pa_\infty\Bar{\psi}_\infty) 2\qk z_\infty \partial \barpartial \Bar{f} \, \barpartial f 
+ \chi_\infty\Bar{\chi}\pa_\infty 2\qk z_\infty \partial \barpartial f \barpartial \Bar{f} \;\dmu \\
&= 2\chi_\infty\Bar{\chi}\pa_\infty \qk \Real \int_{\mathbb{D}} z_\infty \partial \barpartial \Bar{f} \, \barpartial f 
+ z_\infty \partial \barpartial f \barpartial \Bar{f} \;\dmu \\
&= 2\chi_\infty\Bar{\chi}\pa_\infty \qk \Real \int_{\mathbb{D}} z_\infty \partial \barpartial \Bar{f} \, \barpartial f 
+ \partial( z_\infty \barpartial f \barpartial \Bar{f}) - \barpartial f \barpartial \Bar{f} - z_\infty \barpartial f \partial\barpartial \Bar{f} \;\dmu \\
&= -2\chi_\infty\Bar{\chi}\pa_\infty \qk \int_{\mathbb{D}} \Real \left( \barpartial \Bar{f} \, \barpartial f \right)\;\dmu.
\end{align*}
As in Lemma~\ref{second variation plane}, writing $f$ as real and imaginary parts shows that $\lambda_\infty''$ is indefinite.

For Part~(B), our choice of $f$ is such that $gf$ is imaginary valued.
Checking~\eqref{def:real potential}
\[
\frac{g}{\Bar{g}} \var\U_\infty\sd
= \frac{g}{\Bar{g}} \qk \overline{\qj\barpartial \partial f} \qk
= -\frac{g}{\Bar{g}} \qj\barpartial \partial \Bar{f}
= - \qj\frac{1}{g} \barpartial \partial \overline{gf}
= \qj\barpartial \partial f
= \var\U_\infty,
\]
(we have omitted $d\Bar{z}_\infty$)
shows that $\var \U_\infty \in \pot{\unity,g}$.
If we write $f = g^{-1} \qi u$ for a real valued function $u$ then the quadratic form further simplifies to 
\[
\lambda''_\infty
=-2\chi_\infty\Bar{\chi}\pa_\infty \qk \int_{\mathbb{D}} \Real \left( \barpartial u \right)^2\;\dmu
=-\frac{1}{2}\chi_\infty\Bar{\chi}\pa_\infty \qk \int_{\mathbb{D}} \left( \partial_x u \right)^2 - \left( \partial_y u \right)^2\;\dmu,
\]
which is indefinite, even in this special case.
\end{proof}

\begin{lemma}\label{lem:indefinite weierstrass}
Let $(\U,\chi,\psi)\in\wei(\X,E)$ or $(\U,\chi)\in\wei(\X,E,g)$ be strongly isothermic and the kernel of~\eqref{eq:singularity2b} respectively~\eqref{eq:singularity2c} be one dimensional.
Then the second derivatives $\lambda_m''$ from the proof of Lemma~\ref{lem:weierstrass tangent cone} are indefinite.
\end{lemma}
\begin{proof}
The proof of this lemma closely follows the proof of Lemma~\ref{lem:indefinite}.
We give the complete argument for $\wei(\X,E)$ and then afterwards explain the modifications necessary for $\wei(\X,E,g)$.
We can regauge so that $\U = \U^-$.
We also fix $1<p<2$ and set $p^* = \frac{2-p}{2p}$ and $p' = \frac{p}{p-1}$.
Let $(\chi G, \psi H)$ span the kernel of~\eqref{eq:singularity2b} or~\eqref{eq:singularity2c}.
In particular $G = -\Bar{H}$ by Lemma~\ref{lem:wei-isothermic}.
At the beginning of this chapter we introduced smooth $1$\=/forms $\psi\pa_\gamma$ with Property~\eqref{eq:def line integral}.
They were chosen so that the union of their supports do not cover $\X$.
Choose a chart $z$ in the complement of the supports and use it to identify an open subset of $\X$ with an open set $\Omega \subset \mathbb{C}$.
We choose $z_0 \in \Omega$ such that certain properties hold.
These properties each hold almost everywhere on $\Omega$, so they hold simultaneously almost everywhere.
Firstly, we require that $\chi,\psi$ and $d\Bar{G} = -dH$ have no roots at $z_0$.
By avoiding the roots of the sections, we avoid the singular points of the sheaf $\Sh{S}$ and $(\chi G, \psi H)$ are locally pullbacks of holomorphic sections to $\Tilde{X}$.
We write $(\psi\pa,\chi\pa) = (\chi G, \psi H)$ in a neighborhood of $z_0$.
Let $z_0$ be a point where Lemma~\ref{lem:sobolev local behavior} holds for the sections $\chi,\psi,\psi\pa,\chi\pa$.
Additionally, we choose $z_0$ such that it is Lebesgue point of the potentials $\U^-$, $(\partial\chi)\chi^{-1}$ and $(\partial\psi)\psi^{-1}$.
Later in the proof we will carry out a computation with respect to a basis of the image of $\con(\X,\chi,\psi)$ in $\pot{KE}$ under~\eqref{eq:singularity2b}, so let us choose a (finite) set of elements that map to an orthonormal basis of the image.
We further require $z_0$ be Lebesgue points but not roots of these chosen basis elements.

With the point $z_0$ chosen, we are ready to define the blow-ups.
\index{Blow-up}
First however, we note that because $G$ and $H$ are determined only up to a constant, we can assume without loss of generality that $G$ and $H$ have a root at $z_0$.
This is achieved concretely by the transformation $\psi\pa \mapsto \psi\pa - \chi G(z_0)$ and $\chi\pa \mapsto \chi\pa - \psi H(z_0)$, which does not change that $z_0$ is a Lebesgue point of $\psi\pa$ or $\chi\pa$.
By blow-up, recall that we mean that for each $r > 0$ we restrict functions to the disc $\Set{B}_r=B(z_0,r)\subset\Omega$ and rescale this disc to the unit disc $\mathbb{D}$ through the coordinate changes $z_\infty = \frac{z-z_0}{r}$.
In the limit as $r\downarrow0$ the functions $\chi_r(z_\infty) = \chi(z_0+rz_\infty)$ and $\psi_r(z_\infty)=\psi(z_0+rz_\infty)$ converge  in $\banach{p^*}(\mathbb{D},\qat)$ to non-zero constants $\chi_\infty,\psi_\infty \in \qat$ respectively.
The sections $\psi\pa$ and $\chi\pa$ have a root at $z_0$ and so we apply the last statement of Lemma~\ref{lem:sobolev local behavior} to them.
In the limit as $r\downarrow0$ the functions $\psi\pa_r(z_\infty) = r^{-1}\psi\pa(z_0 + rz_\infty)$ and $\chi\pa_r(z_\infty) = r^{-1}\chi\pa(z_0 + rz_\infty)$ converge in $\banach{p'}(\mathbb{D},\qat)$ to $z_\infty\psi\pa_\infty$ and $z_\infty\chi\pa_\infty$, respectively.
Furthermore, we know that $dG$ and $dH$ have no root at $z_0$, from which we may conclude that $\psi\pa_\infty$ and $\chi\pa_\infty$ are non-zero.
Taking the limit of the local relation $\psi\Bar{\psi}\pa + \chi\pa\Bar{\chi}=0$ gives $\psi_\infty \Bar{\psi}\pa_\infty \Bar{z}_\infty = - z_\infty \chi\pa_\infty \Bar{\chi}_\infty$.
Considering $z_\infty = 1,\qi$ in turn shows that $\psi_\infty \Bar{\psi}\pa_\infty = - \chi\pa_\infty \Bar{\chi}_\infty \in \qj \mathbb{C}$.
Hence the blow-up defines an admissible map to a plane $z_\infty\mapsto \Bar{\chi}_\infty \qj z_\infty \psi_\infty$ and of the dual surface to the dual plane.

Therefore we can deploy Lemma~\ref{second variation plane weierstrass} to give us two variations with compact support $(d\Bar{z}_\infty\var \U_\infty,\var\chi_\infty,\var\psi_\infty)$ in the tangent space of $\wei(\mathbb{D},\unity)$ at $(0,\chi_\infty,\psi_\infty)$ with $\lambda_\infty'' = \pm 1$ respectively.
For small $r>0$, $\Set{B}_r$ is contained in $\Omega$.
For such $r$ we pull-back the variations on the blow-up to obtain
\begin{align*}
\var \Tilde{\U}_r(z) &:= d\Bar{z}\, r^{-2}\,\var \U_\infty(\tfrac{z-z_0}{r}) \;, \\
\var \Tilde{\chi}_r(z) := -r^{-1} \qj (\barpartial \Bar{f})(\tfrac{z-z_0}{r})\chi(z) \;, &\qquad
\var \Tilde{\psi}_r(z) := r^{-1} \qj(\barpartial f)(\tfrac{z-z_0}{r}) \psi(z)
\end{align*}
for $z \in \Set{B}_r$ and zero on the complement of $\Set{B}_r$.
These elements $(\var \Tilde{\U}_r,\var\Tilde{\chi}_r,\var\Tilde{\psi}_r)$ belong to the Banach space $\in \pot{KE}\times H^0(\X,\sob{1,p}{E^{-1}})\times H^0(\X,\sob{1,p}{KE})$ but do not belong to the tangent space of the submanifold $\Spa{U}_1$ or $\Spa{U}_2$ from Lemma~\ref{lem:weierstrass tangent cone} at $(\U,\chi,\psi)$.

The next part of the argument is completely analogous to Lemma~\ref{lem:indefinite}: we must project $\Tilde{\var \U}_r$ onto the orthogonal complement in $\pot{KE}$ of the image of the map~\eqref{eq:singularity2b}.
The result of this projection is $\var \U_r$.
Due to Lemma~\ref{lem:wei-necessary condition} this ensures that there exists a tangent vector $(\var{\U}_r,\var{\chi}_r,\var{\psi}_r)_r$.
Specifically, $(\var\U\sd_r\chi, \var \U_r\psi)$ belongs to the image of the operator~\eqref{fredholm weierstrass}.
We fix an inverse $D^{-1}$ to~\eqref{fredholm weierstrass} in the same way as we fix an inverse of $\delbar{E} - \V$ in Lemma~\ref{lem:indefinite}.
To give a few more details,~\eqref{fredholm weierstrass} is a pair consisting of $\delbar{E^{-1}}-\U\sd$ and $\delbar{KE} -\U$, and for each of these Fredholm operators we can find complementary subspaces to kernels and images.
The complementary subspaces of the kernel and image of the pair operator are the just the cartesian product of these.
Restricting to the complement of the kernel and composing with the projection onto the image makes $D$ an isomorphism.
Thus we obtain $(\var{\chi}_r,\var{\psi}_r) = D^{-1}(\var\U\sd_r\chi, \var \U_r\psi)$ in $H^0(\X,\sob{1,p}{E^{-1}})\times H^0(\X,\sob{1,p}{KE})$.

We need to estimate the size of the difference $(\var \Tilde{\U}_r,\var\Tilde{\chi}_r,\var\Tilde{\psi}_r)-(\var{\U}_r,\var{\chi}_r,\var{\psi}_r)$ in the limit $r\downarrow 0$.
The difference $\var{\U}_r - \Tilde{\var \U}_r$ can be controlled by the inner product
\begin{align*}
\left\langle\ls\Tilde{\upsilon}\pa,\psi\rs + \ls\chi,\Tilde{\phi}\pa\rs, \var \Tilde{\U}_r\right\rangle
=\Real\int_{\Set{B}_r}\lh\Tilde{\upsilon}\pa,\var\Tilde{\U}_r\psi\rh+\lh\Tilde{\phi}\pa,\var\Tilde{\U}_r\sd\chi\rh,
\end{align*}
where $(\Tilde{\upsilon}\pa,\Tilde{\phi}\pa)$ is the restriction to $\Set{B}_r$ of an element of $\con(\X,\chi,\psi)$, specifically an element chosen at the outset of the proof as part of a basis of the image of~\eqref{eq:singularity2b}.
Consider the blowups $\Tilde{\upsilon}\pa_r(z_\infty) = \Tilde{\upsilon}\pa(z_0+rz_\infty)$ and $\Tilde{\phi}\pa_r(z_\infty) = \Tilde{\phi}\pa(z_0+rz_\infty)$.
We now compute
\begin{align*}
\int_{\Set{B}_r} \lh\Tilde{\upsilon}\pa,\var\Tilde{\U}_r \psi\rh 
&= \int_{\mathbb{D}} \lh\Tilde{\upsilon}\pa_r,r^{-2}d\Bar{z}_\infty\var \U_\infty \psi_r\rh \;r^2
\to \int_{\mathbb{D}} \lh\Tilde{\upsilon}\pa_\infty, d\Bar{z}_\infty\var \U_\infty \psi_\infty \rh,
\end{align*}
as $r \downarrow 0$.
Since $\var \U_\infty$ is the derivative of a function of compact support, the integral  on the right is zero.
The same argument applies to $\lh\Tilde{\phi}\pa,\var\Tilde{\U}_r\sd\chi\rh$.
Together we conclude
\begin{align*}
\lim_{r\downarrow0} \|\var \Tilde{\U}_r - \var \U_r\|_{\pot{KE}}& =0 .
\end{align*}
Next we bound $\var\chi_r - \var\Tilde{\chi}_r$.
Consider
\begin{align*}
(\barpartial &- \U\sd)\left( \var\chi_r - \var\Tilde{\chi}_r \right) \\
&=\var\U\sd_r\chi + d\Bar{z}\,r^{-2} \qj (\partial\barpartial \Bar{f})(\tfrac{z-z_0}{r})\chi + r^{-1} \qj (\barpartial \Bar{f})(\tfrac{z-z_0}{r}) dz\partial\chi + \U\sd \var\Tilde{\chi}_r \\
&= \qk (\var \U_r - \var \Tilde{\U}_r) \qk \chi - \var\Tilde{\chi}_r \chi^{-1}dz \partial\chi + \U\sd \var\Tilde{\chi}_r.
\end{align*}
The first term converges to zero in $\banach{p}(\X,\qat)$ since we have just seen that $\var \U_r-\var\Tilde{\U}_r$ converges to zero in $\pot{KE}$ and $\chi$ is $\banach{p^*}(\X,\qat)$.
The remaining terms also vanish in the limit, since in magnitude they are essentially $\var\Tilde{\chi}_r$ multiplied with an $\banach{2}$\=/potential and this was shown to converge to zero in~\eqref{eq:W varupsilon Lp limit} under the assumption that $z_0$ is a Lebesgue point of the potential.
The same argument applies to $\var\psi_r$; unlike in the Kodaira case, the convergence properties of both variations are the same.
Together this shows that $D( \var\chi_r - \var\Tilde{\chi}_r, \var\psi_r - \var\Tilde{\psi}_r)$ tends to zero in $\banach{p}(\X,\qat)^{\times 2}$ and hence that $( \var\chi_r - \var\Tilde{\chi}_r, \var\psi_r - \var\Tilde{\psi}_r)$ converges to $0$ in $\sobolev{1,p}(X)^{\times 2}$.
By Sobolev embedding we also have convergence to zero in $\banach{p^*}(\X,\qat)^{\times 2}$.

We are now in a position to show that $\lambda_m''$ is indefinite by computing the limit of its value in the direction of $(\var{\U}_r,\var{\chi}_r,\var{\psi}_r)$.
From the formulas~\eqref{eqn:lambda1''} and\eqref{eqn:lambda2''} in Lemma~\ref{lem:weierstrass tangent cone} we have 
\begin{align*}
-\tfrac{1}{2}\lambda_1''
&= \Real\int_\X \lh\var\chi,\var\U\psi H\rh + \lh\chi G,\var\U\var\psi\rh, \\
\tfrac{1}{2}\lambda_2''
&= \Real\int_{\hat{\X}} d\Big[\Bar{G} \lp\var\chi,\var\psi\rp \Big] + \lh\var\chi,\var\U\psi H\rh + \lh\chi G,\var\U\var\psi\rh.
\end{align*}
We handle the first term of $\lambda_2''$, the `boundary' term, separately, because it is unlike the others and in fact tends to zero.
As in Lemma~\ref{lem:period pairing}, it is a finite sum of terms of the form
\[
\Real \left[ -\per d\Bar{G}(B_i) \int_{A_i} \lp\var\chi,\var\psi\rp \right]
\quad\text{or}\quad
\Real \left[ \per d\Bar{G}(A_i) \int_{B_i} \lp\var\chi,\var\psi\rp \right].
\]
So we should consider the value of 
\[\int_\gamma \lp \var \chi_r, \var \psi_r \rp
= \int_\X \eta_\gamma \wedge \lp \var \chi_r, \var \psi_r \rp
= \int_\X \eta_\gamma \wedge \lp \var \chi_r - \var\Tilde{\chi}_r, \var \psi_r - \var \Tilde{\psi}_r \rp\]
for $\gamma$ an element of the chosen basis of $H_1(\X,\mathbb{Z})$, compare to Equation~\eqref{eq:def line integral}.
The second equality holds due to the choice of $z_0$: the support of $(\var \Tilde{\U}_r,\var\Tilde{\chi}_r,\var\Tilde{\psi}_r)$ has empty intersection with the support of $\eta_\gamma$.
But we see that this integral is bounded by $\|\eta_\gamma\|_\infty \|\var\chi_r - \var\Tilde{\chi}_r \|_{\banach{p^*}(\X)} \|\var\psi_r - \var\Tilde{\psi}_r \|_{\banach{p^*}(\X)}$, an expression that tends to zero.

Thus to understand the limiting values of $\lambda_m''$ we can put aside the boundary term.
Consider the difference of the other terms and the corresponding expressions with $(\var \Tilde{\U}_r,\var\Tilde{\chi}_r,\var\Tilde{\psi}_r)$ as the variations.
An upper bound for $\lh\psi\pa,\var\U_r\var\psi_r\rh-\lh\psi\pa,\var\Tilde{\U}_r \var\Tilde{\psi}_r\rh$ in the $\banach{1}(\X,\qat)$\=/norm is
\begin{align*}
&\|\psi\pa\|_{\banach{p'}(\X)} \|\var \U_r - \var\Tilde{\U}_r\|_{\banach{2}(\X)} \|\var\psi_r - \var\Tilde{\psi}_r \|_{\banach{p^*}(\X)} \\
&+ \|\psi\pa\|_{\banach{p'}(\Set{B}_r)} \|\var \U_r - \var\Tilde{\U}_r\|_{\banach{2}(\X)} \|\var\Tilde{\psi}_r \|_{\banach{p^*}(\Set{B}_r)} \\
&+ \|\psi\pa\|_{\banach{p'}(\Set{B}_r)} \|\var\Tilde{\U}_r\|_{\banach{2}(\Set{B}_r)} \|\var\psi_r - \var\Tilde{\psi}_r \|_{\banach{p^*}(\X)}.
\end{align*}
The first term vanishes in the limit by the smallness of the perturbation.
In the remaining terms, we can use the fact that $\psi\pa$ has a root to show convergence to zero.
To take the last term as an example
\[
\|\psi\pa\|_{\banach{p'}(\Set{B}_r)} \|\var\Tilde{\U}_r\|_{\banach{2}(\Set{B}_r)}
= r^{1+ 2/p'} \|\psi\pa_r\|_{\banach{p'}(\mathbb{D})} r^{-2+1}\|\var \U_\infty\|_{\banach{2}(\mathbb{D})}
\to 0.
\]
The same calculation holds for $\lh\psi\pa,\var \U_r \var\psi_r \rh - \lh \psi\pa,\var \Tilde{\U}_r \var\Tilde{\psi}_r \rh$, and also using the $\banach{1}(\hat{\X},\qat)$\=/norm.
This shows that $\lambda_m''$ in the direction $(\var{\U}_r,\var{\chi}_r,\var{\psi}_r)$ have the same value in the limit as in the direction $(\var \Tilde{\U}_r,\var\Tilde{\chi}_r,\var\Tilde{\psi}_r)$.
Thus is remains to compute
\begin{align*}
&\lim_{r\downarrow0} \Real \int_{\Set{B}_r} \lh\psi\pa,\var \Tilde{\U}_r \var\Tilde{\psi}_r \rh + \lh\chi\pa,\var\Tilde{\U}\sd_r\var\Tilde{\chi}_r\rh \\
& =\lim_{r\downarrow0} \Real \int_{\mathbb{D}} \lh\psi\pa_r,\var \U_\infty \var\psi_\infty\psi_\infty^{-1}\psi_r\rh
+ \lh\chi\pa_r,\var\U\sd_\infty\var\chi_\infty\chi_\infty^{-1}\chi_r\rh \\
&= \Real \int_{\mathbb{D}} \lh z_\infty\psi\pa_\infty,\var \U_\infty \var\psi_\infty \rh
+ \lh z_\infty\chi\pa_\infty,\var\U\sd_\infty\var\chi_\infty\rh
= \lambda_\infty''.
\end{align*}
Since in the limit $\lambda_m''$ must approach $\pm 1$ there exist tangent vectors on which it is positive and negative: it is indefinite.

We now explain the changes to the argument for the $\wei(\X,E,g)$ situation.
The condition $\U \in \pot{KE,g}$ already forces $\U = \U^-$, and as usual we set $\psi = g\chi$ and $G = H$.
This gives the relation $\psi\pa = g^{-1} \chi\pa$.
The choice of the point $z_0$ proceeds largely the same.
Since $z_0$ is not a root of $\psi$ it is not a root of $g$.
The only change is that instead of elements in $\con(\X,\chi,\psi)$, we choose elements of $\con(\X,\chi,g)$ that are mapped to an orthogonal basis of the image of~\eqref{eq:singularity2c} in $\pot{KE,g}$.
Note that the transformation of $\psi\pa$ and $\chi\pa$, which gives them a root, preserves the relation $\psi\pa = g^{-1} \chi\pa$.
For the blow-up, since $g$ is holomorphic we have $\psi_\infty = g(z_0)\chi_\infty$ and $\psi\pa_\infty = g^{-1}(z_0) \chi\pa_\infty$.

From Lemma~\ref{second variation plane weierstrass}(B) we obtain variations $(d\Bar{z}_\infty\var \U_\infty,\var \chi_\infty, \var \psi_\infty)$.
The same formulas as above define $(\var \Tilde{\U}_r,\var\Tilde{\chi}_r,\var\Tilde{\psi}_r)$, which implies $\var \Tilde{\U}_r \in \pot{KE,g}$ and $\var\Tilde{\psi}_r = g\var\Tilde{\chi}_r$.
However the definition of $(\var{\U}_r,\var{\chi}_r)$ is necessarily different.
We project $\var \Tilde{\U}_r$ onto the orthogonal complement in $\pot{KE,g}$ of the image of~\eqref{eq:singularity2c}.
The first claim of Lemma~\ref{lem:wei-necessary condition}(B) tells us that there exists a tangent vector $(\var \U_r, \var \chi_r)$ to $\wei(\X,E,g)$ and we can obtain $\var \chi_r$ through an inverse of $\delbar{E^{-1}}-\U\sd$ as above.
The proof of the convergence of $\var \U_r - \var \Tilde{\U}_r$ and $\var \chi_r - \var \Tilde{\chi}_r$ to zero is unchanged.

Finally, to prove the convergence of $\lambda_m''$ to $\lambda_\infty''$ we observe the only difference is the presence of $g \var\chi_r$ in the `boundary' term $d\left[-H\lp\var \chi_r, g\var\chi_r\rp\right]$ of~\eqref{eqn:lambda2'' real}.
The easiest way to handle this is to define $\var \psi_r = g \var \chi_r$ and calculate
\begin{align*}
\|\var\psi_r - \var\Tilde{\psi}_r \|_{\banach{p^*}(\X)}
= \|g (\var\chi_r - \var\Tilde{\chi}_r) \|_{\banach{p^*}(\X)}
\leq \|g\|_\infty \|\var\chi_r - \var\Tilde{\chi}_r \|_{\banach{p^*}(\X)},
\end{align*}
where $\|g\|_\infty$ means the maxima of suprema of $g_l$ as functions on a cover $\{\SO_l\}$ that trivializes $KE^2$.
A reader may be concerned that although $g$ is holomorphic on each chart of the cover, it may be unbounded.
This possibility is anticipated by Remark~\ref{rem:special cover}: we know that $g_l$ extends to the boundary of $\SO_l$.
Then the right hand side of the above tends to zero with $r$.
Now we can apply the argument for Part~(A) to this term.
\end{proof}

\begin{proof}[Proof of Theorem~\ref{thm:weierstrass isothermic characterization}]
This proof is essentially the same as the proof of Theorem~6.8.
To prove the (A) part of this theorem, we use the (A) parts of the lemmas above, and similarly for Part~(B), and so we omit the (A) and (B) labels in this argument.
The equivalence of (i) and (ii) is just Lemma~\ref{lem:wei-isothermic}(ii).
This time it's Lemma~\ref{lem:weierstrass tangent cone} that provides the function $\lambda$ and allows us to apply Lemma~\ref{lem:tangent cone spans} to show that (ii) implies (iii).
Lemma~\ref{lem:tangent cone spans} also proves that the tangent cone spans the tangent space and the existence of a smooth path.
Finally~(iv) is a direct consequence of~(iii) and the contrapositive of Lemma~\ref{lem:wei-non-isothermic} is that (iv) implies (i).
\end{proof}

\chapter{Constrained Willmore Weierstraß Triples}
\label{chapter:constrained weierstrass}
The main goal of this chapter is the description of critical points of the Willmore functional on the space of all admissible maps from a compact Riemann surface $\X$ to $\Imag\qat$.
In Chapter~\ref{chapter:weierstrass}, we reformulated the condition that an admissible map takes values in $\Imag\qat$ by using the quaternionic Weierstraß representation~\ref{thm:weierstrass}.
This is the only efficient way we know of to understand such admissible maps.
Hence we transfer the previous characterization in Theorem~\ref{constrained Willmore 1} of constrained Willmore surfaces in terms of the quaternionic Kodaira representation to an analogous characterization in Theorem~\ref{constrained Willmore 2} in terms of the quaternionic Weierstraß representation.

As in the beginning of Chapter~\ref{chapter:constrained 1} we can apply Corollary~\ref{gauge holomorphic structure} to modify Equation~\eqref{eq:willmore energy weierstrass} to $\U \in \pot{KE}$: An admissible map that corresponds to $(\U,\chi,\psi)\in\wei(\X,E)$ or $(\U,\chi)\in\wei(\X,E,g)$ has Willmore energy $\willmore=4\|\U^-\|_2^2$. 
Hence at points where $\wei(\X,E)$ or $\wei(\X,E,g)$ is a manifold, the characterization of constrained Willmore Weierstra{\ss} triples depends on the obstructions on the variations $\var\U$ of the Weierstra{\ss} potentials that ensure the existence of variations $(\var\chi,\var\psi)$ or variations $\var\chi$ such that $(\var\U,\var\chi,\var\psi)$ is tangent to $\wei(\X,E)$ or $(\var\U,\var\chi)$ is tangent to $\wei(\X,E,g)$, respectively. Thus the smooth case follows from the characterization of these obstructions in Lemma~\ref{lem:wei-necessary condition}.

For non-smooth points of $\wei(\X,E)$ and $\wei(\X,E,g)$, Lemma~\ref{lem:wei-isothermic} describes two different cases. 
Our first task, Lemma~\ref{lem:dWillmore-space-cone 2}, is the treatment of the case with higher-dimensional kernels of the maps~\eqref{eq:singularity2b} and~\eqref{eq:singularity2c}. 
In this case, as in the analogous case in Chapter~\ref{chapter:constrained 1}, the derivative of the Willmore functional vanishes on the whole tangent space. 
The second case, that of one-dimensional kernels, may be handled by an extension of the smooth case, since the tangent cone spans the tangent space. We close this chapter with a proof that an admissible map has constrained Willmore Kodaira data if and only if it has constrained Willmore Weierstra{\ss} data.

\begin{definition}
\index{Constrained Willmore}
The Weierstraß data $(\U,\chi,\psi)\in\wei(\X,E)$ or $(\U,\chi)\in\wei(\X,E,g)$ are called \emph{constrained Willmore} if and only if the derivative of $\willmore$ vanishes on the corresponding tangent cone or tangent space.
\end{definition}

The equivalence of vanishing on the tangent space or on the tangent cone is not yet completely shown.
If the kernel of~\eqref{eq:singularity2b} is at most one than Theorem~\ref{thm:weierstrass isothermic characterization} the tangent cone spans the corresponding tangent space.
Lemma~\ref{lem:wei-isothermic} (iii) characterizes the case that the kernel has dimension greater than one.
In this case, the derivative of the Willmore energy always vanishes on the tangent space.

\begin{lemma}\label{lem:dWillmore-space-cone 2}
Let $\X$ be a (not necessarily compact) Riemann surface.
\begin{enumeratethm}[label={\upshape(\Alph*)}]
\item 
If $(\U,\chi,\psi)\in\wei(\X,E)$ obeys $(\delbar{E}-\U^+)(\overline{\qj\chi})^{-1}=0$, then $H^0(\X,\Q{KE,\U})$ contains an element $\chi\pa$ with $\U^-=\ls\chi,\chi\pa\rs$.
In particular, the derivative of $\willmore$ vanishes on the tangent space of $\wei(\X,E)$ at $(\U,\chi,\psi)$.
\item
If $(\U,\chi)\in\wei(\X,E,g)$ obeys $\delbar{E}(\overline{\qj\chi})^{-1}=0$, then $H^0(\X,\Q{E^{-1},\U\sd})$ contains an element $\psi\pa$ with $\U=\ls\chi,g\psi\pa\rs+\ls\psi\pa,g\chi\rs$.
In particular, the derivative of $\willmore$ vanishes on the tangent space of $\wei(\X,E,g)$ at $(\U,\chi)$.
\end{enumeratethm}
\end{lemma}
\begin{proof} (A): An application of Corollary~\ref{gauge holomorphic structure} transforms the triple $(\U,\chi,\psi)$ into another triple of the form $(\U',\chi',\psi')\in\wei(\X,E')$ with $\U'\in\pot{KE'}^-$.
To simplify notation we omit the primes and replace $(\U',\chi',\psi')$ by $(\U,\chi,\psi)\in\kodaira(\X,E)$.
Next we shall define $\chi\pa$ such that $\ls\chi,\chi\pa\rs=\U$.
As in the proof of Corollary~\ref{cor:global darboux} we cover $\X$ by domains of charts $z_l:\Set{O}_l\to\Omega_l$ such that $E$ is represented by the holomorphic cocycle $f_{ml}$ with respect to the open cover $(\Set{O}_l)_l$.
Let $(\V,\upsilon,\phi)\in\kodaira(\X,E)$ be the image of the inverse Darboux transform~\ref{cor:global darboux} of $(\U,\chi,\psi)$.
This means $\upsilon=-(\overline{\qj\chi})^{-1}$ and $\phi=\upsilon F$.
Theorem~\ref{thm:darboux} implies $\V=(\delbar{E} \upsilon) \upsilon^{-1}$.
Our assumption on $\chi$ means that this vanishes.
Therefore $\upsilon$, $\phi$ and $\U^-$ are represented locally by holomorphic functions and therefore analytic.
Consequently, the triple $(\U,\chi,\psi)$ is also analytic.

Locally on $\Set{O}_l$ the sections $\upsilon$ and $\chi$ are represented by the holomorphic functions $\upsilon_l$ and $\chi_l$.
The potentials $\B$ and $\U$ are represented by
\begin{align*}
\B_l&=\tfrac{1}{2}\big(\qi\partial_l\upsilon_l\upsilon_l^{-1}\qi-\partial_l\upsilon_l\upsilon_l^{-1}\big),&\U^-_l&=-\tfrac{1}{2}\big(\qi\partial_l\upsilon_l\upsilon_l^{-1}\qi+\partial_l\upsilon_l\upsilon_l^{-1}\big).
\end{align*}
From Definition~\ref{def:01 pairing} the $1$\=/form $\ls\chi,\chi\pa\rs$ is on $\Set{O}_l$ represented by $2\qk\chi_l\Bar{\chi}\pa_l$ where $\chi\pa_l$ represents locally $\chi\pa$.
In order to have $\U^-=\ls\chi,\chi\pa\rs$ we define
\[
2\chi\pa_l=\overline{\chi_l^{-1}\qk^{-1}\U^-_l}=\Bar{\U}^-_l\qk\Bar{\chi}_l^{-1}=-\U^-_l\qi\qj\Bar{\chi}^{-1}=-\U^-_l\qi(\overline{\qj\chi})^{-1}=-\U^-_l\qi\upsilon_l.
\]
Due to~\eqref{eq:potential transformation}, these local functions transforms as
\[
2\chi\pa_m=-\U^-_m\qi\upsilon_m=-f_{ml}\tfrac{dz_l}{dz_m}\U_l f_{ml}^{-1}\qi f_{ml}\upsilon_l=2f_{ml}\tfrac{dz_l}{dz_m}\chi\pa_l.
\]
Hence these local functions represent a global section of $KE_\qat$.
The flatness of the connection induced by the non-vanishing section $\upsilon$ of $E_\qat$ evidenced in Equation~\eqref{eq:zero curvature 2} implies $\barpartial_l \U^-_l=\U^-_l B_l$.
Together with Theorem~\ref{thm:darboux}(b) we obtain
\[
2\barpartial_l \chi\pa_l=-\barpartial_l \U_l\qi\upsilon_l-\U^-_l\qi\partial_l\upsilon_l=-\U^-_l\B_l\qi\upsilon_l+\U^-_l\qi(\B_l+\U^-_l)\upsilon_l=-\U^-_l\U^-_l\qi\upsilon_l=\U^-_l2\chi\pa_l.
\]
Therefore the local functions $\chi\pa_l:\Set{O}_l\to\qat$ indeed define a global section $\chi\pa\in H^0(\X,\Q{KE,\U^-})$ such that $\U^-=\ls\chi,\chi\pa\rs$.
Now Lemma~\ref{lem:wei-necessary condition} implies that the component $\var \U$ of any element $(\var \U,\var\chi,\var\psi)$ of the tangent space of $\wei(\X,E)$ at $(\U^-,\chi,\psi)$ is orthogonal to $\U^-$.
In particular $\var \willmore = 8 \langle \U^-,\var \U\rangle$, the derivative of $\willmore$, vanishes on this tangent space.
This proves Part~(A).

For Part~(B) we define $\psi\pa_l=-\frac{1}{4}g_l^{-1}\U^-_l\qi\upsilon_l$.
This has the property that $4\chi_l\Bar{\psi}\pa_l=-\qj\qi\Bar{\U}^-_l\Bar{g}_l^{-1}=-\qk\U^-_l\Bar{g}_l^{-1}$, and because of the reality condition~\eqref{involution for imaginary maps} of the potential $\overline{\qk\U^-_l\Bar{g}_l^{-1}}=g_l^{-1}g_l\qk\U^-_l\qk g_l^{-1}\Bar{\qk}=\qk\U^-_l\Bar{g}_l^{-1}$, we can see that $\chi_l\Bar{\psi}\pa_l$ is real.
This implies
\[
2\qk \chi_l \overline{g_l\psi\pa_l} +2\qk \psi\pa_l \overline{g \chi_l}
=2\qk ( \chi_l \Bar{\psi}\pa_l + \psi\pa_l \Bar{\chi}_l ) \Bar{g}_l
=4\qk\chi_l\Bar{\psi}\pa_l\Bar{g}_l
= \U^-_l .
\]
Observe that $2 g_l \psi\pa_l$ is the formula for $\chi\pa_l$ from Part~(A).
Therefore $\psi\pa$ is a global section of $E^{-1}_\qat$ and it is $\U^-$\=/holomorphic.
We again apply Lemma~\ref{lem:wei-necessary condition}, except we use Part~(B).
This says that $\var \U$ is the first component of a tangent vector to $\wei(\X,W,g)$ if and only if it is orthogonal to the image of 8.23, which we have shown to include $\U^-$.
Therefore the derivative of $\willmore$ vanishes at this point.
\end{proof}

In the following theorem we give a characterization in terms of additional data $(\chi G, \psi H)$, which could be viewed as analogous to Lagrange multipliers, as discussed following Theorem~\ref{constrained Willmore 1}.

\begin{theorem}
\label{constrained Willmore 2}
Let $\X$ be a compact Riemann surface.
\begin{enumeratethm}[label={\upshape(\Alph*)}]
\item
The Weierstraß data $(\U,\chi,\psi)\in\wei(\X,E)$ is constrained Willmore, if and only if $\U^-$ is the image of some $(\chi G,\psi H)\in\con(\X,\chi,\psi)$ with respect to the map~\eqref{eq:singularity2b}.
In this case the admissible map is analytic.

\item
The Weierstraß data $(\U,\chi)\in\wei(\X,E,g)$ is constrained Willmore, if and only if $\U^-$ is the image of some $g\chi H\in\con(\X,\chi,g)$ with respect to the map~\eqref{eq:singularity2c}.
In this case the admissible map is analytic.
\end{enumeratethm}
\end{theorem}
\begin{proof}
The proof is similar to the proof of Theorem~\ref{constrained Willmore 1}.
In the case described in Lemma~\ref{lem:dWillmore-space-cone 2}, we constructed the $\U^-$ explicitly with $(\chi G, \psi H) = (\psi\pa,0)$ and $g\chi H = \chi\pa$.
If on the other hand the kernel of~\eqref{eq:singularity2b} or~\eqref{eq:singularity2c} is at most one dimensional, then by Theorem~\ref{thm:weierstrass isothermic characterization} the tangent cones spans the tangent space.
Then the Weierstraß data is constrained Willmore if and only the derivative $\willmore$ vanishes on the tangent space: $\U^-$ is orthogonal to all first components $\var \U$ of tangent vectors of $\wei(\X,E)$ or $\wei(\X,E,g)$.
By Lemma~\ref{lem:wei-necessary condition} this holds if and only if $\U^-$ is in the image of~\eqref{eq:singularity2b} or~\eqref{eq:singularity2c}.
This shows that $\U^-$ has the required form.
The smoothness and analyticity now follows from Lemma~\ref{lem:solutions are analytic} by a local argument with an elliptic system in the same manner as in Theorem~\ref{constrained Willmore 1}.
\end{proof}
From our point of view the image of the maps~\eqref{eq:singularity2b} and~\eqref{eq:singularity2c} describes the constraints on $\var\U$ to be part of a tangent vector along $\wei(\X,E)$ at $(\U,\chi,\psi)$ or along $\wei(\X,E,g)$ at $(\U,\chi)$, respectively. The existence of the elements $(\chi G,\psi H)\in \con(\X,\chi,,\psi)$ with $\U^-=\ls\chi G,\psi\rs+\ls\chi,\phi H\rs$ and of $g\chi H\in\con(\X,\chi,g)$ with $\U=\ls\chi G,g\chi\rs+\ls\chi,g\chi H\rs$, is equivalent to the condition that the derivative of the Willmore functions vanishes on the subset of $\pot{KE}\times H^0(\X,\sob{1,p}{E^{-1}})\times H^0(\X,\sob{1,p}{KE})$ and of $\pot{KE,g}\times H^0(\X,\sob{1,p}{E^{-1}})$, whose elements represent admissible maps, respectively. In the variational analysis with constraints such elements are called Lagrange multipliers. For this reason, just as in the Kodaira case, we call $(\chi G,\psi H)\in\con(\X,\chi,\psi)$ and $g\chi H\in\con(\X,\chi,g)$ the Lagrange multipliers of the constrained Willmore Weierstrass data $(\U,\chi,\psi)$ and $(\U,\chi)$, respectively.

The two representations of admissible maps in terms of Kodaira data and in terms of Weierstraß data has lead to two different definitions of constrained Willmore admissible maps. 
Theorem~\ref{constrained kodaira and weierstrass} shows that these two definitions are equivalent, which means that Kodaira data are constrained Willmore if and only if the corresponding Weierstraß data are constrained Willmore. 
As a preparation we shall first prove an extended version of this statement. 
At the end of Chapter~\ref{chapter:weierstrass} we generalized the global Darboux transformation to a transformation between more general data $H_l\times H\pa_l$ of finite dimensional spaces of holomorphic sections of two paired holomorphic $\qat$\=/line bundles. 
Furthermore, Lemma~\ref{lem:solutions are analytic} contains a natural definition when such data are constrained Willmore. 
The following proposition shows that on a simply connected Riemann surfaces $\X$ the generalized global Darboux transformation and its inverse interchanges constrained Willmore data $H_l\times H\pa_l$ with constrained Willmore data $H_{l-1}\times H\pa_{l-1}$. 
We will not enter into discussion of constrained Willmore data of this form on non-simply connected Riemann surfaces because the period problem seems to be even more involved than the period problem of the Weierstraß data. 
However, the use the following proposition in the proof of Theorem~\ref{constrained kodaira and weierstrass} gives the first hints how one might define constrained Willmore data of the form $H_l\times H\pa_l$.

\begin{proposition}
\label{prop:constraint transformation}
\index{Darboux transformation}
Let $\X$ be a (not necessarily compact) Riemann surface. 
Suppose that $(E,\V,\upsilon)$ is the data describing a left normal on $\X$ as in Theorem~\ref{thm:kodaira normal}. Furthermore, let $\chi=(\overline{\qj\upsilon})^{-1}\in H^0(\X,\Q{E^{-1},\U\sd})$ and the corresponding $\U\sd\in\pot{E^{-1}}^-$ be as in Corollary~\ref{cor:global darboux}.
Finally, take sections $\phi_1,\ldots,\phi_{d-1}\in H^0(\X,\Q{E,\V})$ and $\psi\pa_1,\ldots,\psi\pa_{d-1}\in H^0(\X,\Q{E^{-1},\U\sd})$.
Then the following two conditions are equivalent:
\begin{enumeratethm}
\item[(i)] There exists $\upsilon\pa\in H^0(\X,\Q{KE^{-1},\V\sd})$ such that the potential $\V$ is equal to
\[\V=\ls\upsilon\pa,\upsilon\rs-\sum\limits_{m=1}^{d-1}\big\ls(\partial_l-(\partial_l\chi)\chi^{-1})\psi_m\pa,\phi_m\big\rs.\]
\item[(ii)] There exists $\chi\pa\in H^0(\X,\Q{KE,\U})$ such that the potential $\U\sd$ is equal to
\[\U\sd=\ls\chi\pa,\chi\rs+\sum\limits_{m=1}^{d-1}\big\ls(\partial_l-(\partial_l \upsilon)\upsilon^{-1})\phi_m,\psi\pa_m\big\rs.\]
\end{enumeratethm}
\end{proposition}
\begin{proof}
The analyticity of $\upsilon$ is equivalent to the analyticity of $\chi$. Furthermore, the global Darboux transformation Corollary~\ref{cor:global darboux} implies $d(\chi^{-1}\psi_m\pa)=-\big\lp\upsilon,(\partial_l-(\partial_l\chi)\chi^{-1})\psi_m\pa\big\rp$ and $d(\upsilon^{-1}\phi_m)=\big\lp\chi,(\partial_l-(\partial_l\upsilon)\upsilon^{-1})\phi_m\big\rp$. Therefore, due to Lemma~\ref{lem:solutions are analytic}, each of the conditions~(i) and (ii) implies that all the functions $\upsilon,\chi,\V,\U\sd,\phi_1,\ldots,\phi_{d-1}$ and $\psi_1\pa,\ldots,\psi_{d-1}\pa$ are analytic. Thus we may assume this.

The equations in (i) and (ii) determines $\upsilon\pa$ and $\chi\pa$ in terms of $\upsilon^{-1}, \chi^{-1},\V,\U\sd$ and $\phi_1,\ldots,\phi_{d-1},\psi\pa_1,\ldots,\psi\pa_{d-1}$ respectively. 
We shall derive two expressions in terms of these data whose vanishing is equivalent to $\upsilon\pa$ being $\V\pa$\=/holomorphic and to $\chi\pa$ being $\U$\=/holomorphic, respectively.
Afterwards we show that these two expressions are equal.
The whole argument is purely local.
So the proof only works with the local representatives on $\SO_l$.

We first calculate the local representative of $d\big(\ls\upsilon\pa,\upsilon\rs\big)$:
\begin{multline*}
2d\big(d\Bar{z}_l\qk\upsilon\pa_l\Bar{\upsilon}_l\big)=-2d\Bar{z}_l\qk\wedge\big(d\Bar{z}_l\barpartial_l\upsilon\pa_l\Bar{\upsilon}_l+\upsilon\pa_l\overline{d\upsilon_l}\big)\\
=2d\Bar{z}_l\left(\qi\wedge dz_l\qj\barpartial_l\upsilon\pa_l(\upsilon\pa_l)^{-1}\qk^2\upsilon\pa_l\Bar{\upsilon}_l+2\qk\upsilon\pa_l\wedge\overline{\big(dz_l(\B_l+\U_l)-d\Bar{z}_l\V_l\big)\upsilon_l}\right)\\
=\qj\barpartial_l\upsilon\pa_l(\upsilon\pa_l)^{-1}d\Bar{z}_l\qi\qk\wedge\big(2d\Bar{z}_l\qk\upsilon\pa_l\Bar{\upsilon}_l\big)+\big(2d\Bar{z}_l\qk\upsilon\pa_l\Bar{\upsilon}_l\big)\wedge\big(d\Bar{z}_l(\Bar{\B}_l+\V_l)-dz_l\U_l\big).
\end{multline*}
Here we used that $d\Bar{z}_l\qi\wedge dz_l=-2\dmu_l$ is real and commutes with all quaternions.
Hence $\upsilon\pa$ is $\V\sd$\=/holomorphic if and only if the local representative of $\ls\upsilon\pa,\upsilon\rs$ belongs to the kernel of the following operator
\begin{gather*}
\omega\mapsto d\omega+d\Bar{z}_l\V_l\wedge\omega-\omega\wedge\big(d\Bar{z}_l(\Bar{\B}_l+\V_l)-dz_l\U_l\big).
\end{gather*}
We first apply this operator to the left hand side of the equation in~(i):
\begin{align*}
d\big(d\Bar{z}_l\V_l\big)+d\Bar{z}_l\V_l\wedge d\Bar{z}_l\V_l-d\Bar{z}_l\V_l\wedge\big(d\Bar{z}_l(\Bar{\B}_l+\V_l)-dz_l\U_l\big)&=dz_l\wedge d\Bar{z}_l\big(\partial_l\V_l+\V_l\Bar{\B}_l\big).
\end{align*}
For the other terms we omit the factor $2$ and the subscript $m$ (so that $\phi_l$ means $(\phi_m)_l$, etc) and use~\eqref{eq:dchi}:
\begin{multline*}
-d\big(d\Bar{z}_l\qk(\partial_l-\B_l+\qj\V_l\qj)\psi_l\pa\Bar{\phi}_l\big)+d\Bar{z}_l\V_l\wedge d\Bar{z}_l\qk(\partial_l-\B_l+\qj\V_l\qj)\psi_l\pa\Bar{\phi}_l\\d\Bar{z}_l\qk(\partial_l-\B_l+\qj\V_l\qj)\psi_l\pa\Bar{\phi}_l\wedge\big(d\Bar{z}_l(\Bar{\B}_l+\V_l)-dz_l\U_l\big)\\
=d\Bar{z}_l\big(\qk\wedge d\Bar{z}_l\qj\V_l\qj-\V_l\wedge d\Bar{z}_l\qk\big)(\partial_l-\B_l+\qj\V_l\qj)\psi_l\pa\Bar{\phi}_l\hspace{40mm}\\d\Bar{z}_l\qk(\partial_l-\B_l+\qj\V_l\qj)\psi_l\pa\wedge\Big(\overline{(dz_l\partial_l+d\Bar{z}_l\V_l)\phi_l}+\Bar{\phi}_l\big(\overline{dz_l(\B_l+\U_l)-d\Bar{z}_l\V_l}\big)\Big)\\
=d\Bar{z}_l\qk(\partial_l-\B_l+\qj\V_l\qj)\psi_l\pa\wedge\overline{dz_l(\partial_l+\B_l+\U_l)\phi_l}.\hspace{45mm}
\end{multline*}
Therefore $\upsilon\pa$ is $\V\sd$\=/holomorphic if and only if the following expression vanishes
\begin{gather}\label{eq:upsilon pa holomorphic}
dz_l\wedge d\Bar{z}_l\big(\partial_l\V_l\!+\!\V_l\Bar{\B}_l\big)-
2\sum_{m=1}^{d-1}d\Bar{z}_l\qk(\partial_l\!-\!\B_l\!+\!\qj\V_l\qj)\psi_{m,l}\pa\wedge\overline{dz_l(\partial_l\!+\!\B_l\!+\!\U_l)\phi_{m,l}}.
\end{gather}
For the derivation of the expression, whose vanishing is equivalent to $\chi\pa$ being $\U$\=/holomorphic, we transform the $1$\=/forms on both sides of the equation in (ii) by $\omega\mapsto\qk\Bar{\omega}\qk$.
The local representative of $\qk\overline{\ls\chi\pa,\chi\rs}\qk$ has the exterior derivative
\begin{multline*}
d\big(2\qk\chi_l\Bar{\chi}_l\pa d\Bar{z}_l\big)=2\qk\big((d\chi_l)\Bar{\chi}_l\pa+\chi_l\overline{d\Bar{z}_l\barpartial_l\chi_l\pa}\big)\wedge d\Bar{z}_l\\
=2\qk\big(dz_l(\B_l-\qj\V_l\qj)+d\Bar{z}_l\qk\U_l\qk\big)\wedge\chi_l\Bar{\chi}_l\pa d\Bar{z}_l+2\qk\chi_l\chi_l\pa\overline{d\Bar{z}_l\barpartial_l\chi_l\pa(\chi_l\pa)^{-1}}\wedge d\Bar{z}_l\\
=\big(d\Bar{z}_l(\Bar{\B}_l-\V_l)-dz_l\U_l\big)\wedge\big(2\qk\chi_l\Bar{\chi}_l\pa d\Bar{z}_l\big)+\big(2\qk\chi_l\Bar{\chi}_l\pa d\Bar{z}_l\big)\wedge dz_l\qi\overline{\barpartial_l\chi_l\pa(\chi_l\pa)^{-1}}\qi.\hspace{2mm}
\end{multline*}
Hence $\chi\pa$ is $\U$\=/holomorphic if and only if the local representative of $\qk\overline{\ls\chi\pa,\chi\rs}\qk$ is annihilated by
\begin{gather*}
\omega\mapsto d\omega-\big(d\Bar{z}_l(\Bar{\B}_l-\V_l)-dz_l\U_l\big)\wedge\omega-\omega\wedge dz_l\U_l.
\end{gather*}
Again we start with the term on the left hand side of the equation in (ii):
\begin{align*}
d\big(\U_l d\Bar{z}_l\big)-\big(d\Bar{z}_l(\Bar{\B}_l-\V_l)-dz_l\U_l\big)\wedge\U_l d\Bar{z}_l-\U_l d\Bar{z}_l\wedge dz_l\U_l&=d\Bar{z}_l\wedge dz_l(\barpartial_l-\Bar{B}_l)\U_l.
\end{align*}
For the other terms we again omit the factor $2$ and the subscript $m$:
\begin{multline*}
d\big(\qk\psi\pa_l\overline{(\partial_l+\B_l+\U_l)\phi_l}d\Bar{z}_l\big)-\big(d\Bar{z}_l(\Bar{\B}_l-\V_l)-dz_l\U_l\big)\wedge\qk\psi\pa_l\overline{(\partial_l +\B_l+\U_l)\phi_l}d\Bar{z}_l\\-\qk\psi\pa_l\overline{(\partial_l +\B_l+\U_l)\phi_l}d\Bar{z}_l\wedge dz_l\U_l\\
=\qk\Big((dz_l\partial_l -d\Bar{z}_l\qk\U_l\qk)-dz_l(\B_l+\qk\V_l\qk)+d\Bar{z}_l\qk\U_l\qk\Big)\psi_l\pa\wedge\overline{(\partial_l +\B_l+\U_l)\phi_l}d\Bar{z}_l\\
+\qk\psi_l\pa\Big(\overline{d\Bar{z}_l\U_l(\partial_l +\B_l+\U_l)\phi_l}\wedge d\Bar{z}_l-\overline{(\partial_l +\B_l+\U_l)\phi_l}d\Bar{z}_l\wedge dz_l\U_l\Big)\\
=d\Bar{z}_l\qk\big(\partial_l -\B_l-\qk\V_l\qk)\psi_l\pa\wedge\overline{dz_l(\partial_l +\B_l+\U_l)\phi_l}.\hspace{43mm}
\end{multline*}
Since $\qj\V_l\qj=\qk\qi\V_l\qj=-\qk\V_l\qi\qj=-\qk\V_l\qk$, $\chi\pa$ is $\U$\=/holomorphic if and only if the following expression vanishes:
\begin{gather}\label{eq:chi pa holomorphic}
d\Bar{z}_l\wedge dz_l\big(\barpartial_l\U_l\!-\!\Bar{B}_l\U_l\big)-2\sum_{m=1}^{d-1}d\Bar{z}_l\qk\big(\partial_l \!-\!\B_l\!+\!\qj\V_l\qj)\psi_{m,l}\pa\wedge\overline{dz_l(\partial_l \!+\!\B_l\!+\!\U_l)\phi_{m,l}}.
\end{gather}
With $\V^+=0=\U^+$ the $\qat^-$\=/part of the zero curvature equation~\eqref{eq:zero curvature 2} for the connection $(\partial_l +\B_l+\U_l)+(\barpartial_l-\V_l)$ induced by $\upsilon$ yields the equation
\[
\barpartial_l\U_l-\U_l\B_l=-\partial_l \V_l-\B_l\V_l.\hspace{-2mm}
\]
Therefore the expressions~\eqref{eq:upsilon pa holomorphic} and~\eqref{eq:chi pa holomorphic} are equal.
\end{proof}

\begin{example}[Catenoid]
\index{Catenoid!Inverted catenoid}
\index{Constrained Willmore!Inverted catenoid}
Let us explain how we used this proposition to construct the Lagrange multipliers in Example~\ref{eg:inverted catenoid willmore}.
First, the catenoid is constrained Willmore in the Weierstraß sense because $\U = 0$ and we can simply take $(\chi\pa,\psi\pa) = (0,0)$.
The relationship between the Lagrange multipliers $(\upsilon\pa,\phi\pa)$ and $(\chi\pa,\psi\pa)$ is $\phi\pa = -(\partial - (\partial\chi)\chi^{-1})\psi\pa$.
In particular, this tells us that $\phi\pa = 0$.
Now that we know $\V$, $\upsilon$, $\phi$, and $\phi\pa$ we can solve for $\upsilon\pa$.
The proposition guarantees that the result, $\upsilon\pa = - \tfrac{1}{2}z^{-1}\qj(z + \qk)(1+|z|^2)^{-1}$, will be $\V\sd$\=/holomorphic.

Because the Weierstraß potential $\U$ is not preserved under general Möbius transformation, Equation~\eqref{eq:mobius potentials}, there is not an obvious formula for how $(\chi\pa,\psi\pa)$ transform.
One method of producing such a formula would be to first find the equivalent $(\upsilon\pa,\phi\pa)$, transform these according to~\eqref{eq:Lagrange mutlipliers mobius}, and then invert the Darboux transformation to find $(\Tilde{\chi}\pa,\Tilde{\psi}\pa)$.
\end{example}

Now that we can transform the Lagrange multipliers between the Kodaira and Weierstraß cases, we are ready to prove that constrained Willmore in the two cases are equivalent.

\begin{theorem}
\label{constrained kodaira and weierstrass}
\index{Constrained Willmore}
Let $\X$ be a compact Riemann surface.
\begin{enumeratethm}[label={\upshape(\Alph*)}]
\item
An admissible map has constrained Willmore Kodaira data if and only if it has constrained Willmore Weierstraß data.
\item
An admissible map $F: \X \to \Imag\qat$ without roots of $dF$ has constrained Willmore Kodaira data if and only if it has constrained Willmore Weierstraß data $(\U,\chi)\in\wei(\X,E,g)$.
\end{enumeratethm}
\end{theorem}
Note that the statement in Part~(B) is in general not true for admissible maps $F:\X\to\Imag\qat$, if $dF$ has roots: The corresponding data $(\U,\chi,\psi)\in\wei(\X,E)$ have deformations which do not preserve the roots of $dF$, whereas all deformations of the data $(\U,\chi)\in\wei(\X,E,g)$ preserve the roots of $dF$ which are the roots of $g$. For this reason not all constraints in $\con(\X,\chi,g)$ on $\var\U$ in $\wei(\X,E,g)$ induce a corresponding constraint in $\con(\X,\chi,\psi)$ of $\var\U$ in $\wei(\X,E)$. In fact, the derivatives of the functions $G$ of the elements $(\chi G,\psi H)\in\con(\X,\chi,\psi)$ are defined in Definition~\ref{def:extended pairs}(A) as $dG=\lp\upsilon,\upsilon\pa\rp$ with $\upsilon\pa\in H^0(\X,\Sh{V}\sd)$. However, the derivatives of the functions $H$ of the elements $g\chi H\in\con(\X,E,g)$ are in Definition~\ref{def:extended pairs}(B) defined as $dH=\lp\zeta,\zeta\pa\rp=\lp g^{-1}\upsilon,\zeta\pa\rp$ with $\zeta\pa\in H^0(\X,\Sh{S})$ and correspond to the elements $\upsilon\pa=g^{-1}\zeta\pa$ in the bigger space $H^0(\X,g^{-1}\Sh{S})$. In particular, not all Lagrange multipliers of $(\U,\chi)\in\wei(\X,E,g)$ induce a Lagrange multiplier of $(U,\chi,g\chi)\in\wei(\X,E)$.
\begin{proof}
We prove Part~(A) first. 
Theorems~\ref{constrained Willmore 1} and~\ref{constrained Willmore 2} show that the constrained Willmore Kodaira triples and the constrained Willmore Weierstraß triples are analytic. 
In a first step we apply Corollary~\ref{gauge holomorphic structure} simultaneously to the Kodaira triple $(\V,\upsilon,\phi)$ and the Weierstraß triple $(\U,\chi,\psi)$ such that $\V^+=\U^+$ becomes zero. 
Furthermore, the global and local sections of $\Q{KE^{-1},\V\sd},\Q{E^{-1}\U\sd}$ are transformed like $\chi$ and the global and local sections of $\Q{KE,\U}$ are transformed like $\upsilon,\phi$ and $\psi$. 
In this way the transformed sections are mapped by~\eqref{eq:singularity1},~\eqref{eq:singularity2b} and~\eqref{eq:singularity2c} to the transformed potentials.

Now we show that for a given constrained Willmore Kodaira triple $(\V,\upsilon,\phi)$ with $\V^+=0$ the corresponding Weierstraß triple is also constrained Willmore. 
On any simply connected subset $\SO_l$ of $\X$ there exists $\psi\pa_l$ such that $-(\partial_l-(\partial_l\chi_l)\chi^{-1}_l)\psi_l\pa=\phi_l\pa$ by the inverse Darboux transformation.
Due to the foregoing proposition these $\psi_l\pa$ can be supplemented by $\chi_l\pa$ to Lagrange multipliers of the Weierstraß triple on $\SO_l$, which therefore are constrained Willmore. 
It remains to show that these Lagrange multipliers $(\chi_l\pa,\psi_l\pa)$ define global elements of $\con(\X,\chi,\psi)$. 
By unique continuation we may assume that $\psi\pa$ and $\upsilon\pa$ extend to global $\V\sd$\=/holomorphic sections on the universal covering $\Tilde{\X}$. 
Furthermore $d(\chi^{-1}\psi\pa)=-\ls\upsilon,(\partial_l-(\partial_l\chi)\chi^{-1})\psi\pa\rs=\ls\upsilon,\phi\pa\rs$ is well defined on $\X$. 
This means that the difference $\psi_l\pa-\psi_m\pa$ of the local representatives $\psi_l\pa$ on the intersection of two such subsets $\SO_l\cap\SO_m$ is locally constant. 
Since $\chi\pa$ is uniquely determined by $\U\sd,\chi,\psi$ and $\psi\pa$, the differences $\chi_l\pa-\chi_m\pa$ of the local representatives of $\chi\pa$ on $\SO_l\cap\SO_m$ are also constant. 
Since $\U\sd$ is single-valued on $\X$ the local representatives $\chi_l\pa$ and $\psi_l\pa$ indeed meet the condition described in Definition~\ref{def:extended pairs} and define an element $(\chi G,\psi H)\in \con(\X,\chi,\psi)$.

For the reverse we also apply locally on $\SO_l$ the forgoing proposition. 
For given constrained Willmore Weierstraß triple $(\U,\chi,\psi)$ with Lagrange multipliers $(\chi G,\psi H)\in\con(\X,\chi,\psi)$ the differences $\chi_l G_l-\chi_m G_m$ of the local representatives of $\chi G$ are constant. 
Consequently the local representatives $\phi_l\pa=-(\partial_l-(\partial_l\chi_l)\chi_l^{-1})(\chi_l G_l)$ defines a global $\phi\pa\in H^0(\X,\Q{KE^{-1},\V\sd})$. 
Since $\upsilon\pa$ is uniquely determined by $\V,\upsilon,\phi$ and $\phi\pa$, it is also single valued on $\X$. 
Hence the corresponding Kodaira triple $(\V,\upsilon,\phi)$ triple is also constrained Willmore with Lagrange multipliers $(\upsilon\pa,\phi\pa)$.

For the proof of Part~(B) let $(\U,\chi)\in\wei(\X,E,g)$ be the Weierstraß data of an admissible map $F:\X\to\Imag\qat$ without roots of $dF$. In this case $g$ and $\psi=g\chi$ have no roots and the sheaf $\Sh{S}$ is the Serre dual sheaf ${\Sh{V}'}\sd=\Q{K^2E,g\V\sd g^{-1}}=\Q{gKE^{-1},g\V\sd g^{-1}}=g\Sh{V}\sd$ of the sheaf $\Sh{V}'=\Q{(KE)^{-1},\V'}=\Q{g^{-1}E,g^{-1}\V g}=g^{-1}\Sh{V}$ with the unique potential $g^{-1}\V g\in\pot{(KE)^{-1}}^-=\pot{g^{-1}E}^-$ such that $\zeta=-(\overline{\qj\psi})^{-1}=-g^{-1}(\overline{\qj\chi})^{-1}=g^{-1}\upsilon$ is $g^{-1}\V g$\=/holomorphic without roots. If $g\chi H\in\con(\X,\chi,g)$ are the corresponding Lagrange multiplier in Theorem~\ref{constrained Willmore 2}(B), then $(\chi G,\psi H)=(\chi H,g\chi H)$ are Lagrange multipliers in $\con(\X,\chi,\psi)$. This shows that the Weierstraß data $(\U,\chi,\psi)\in\wei(\X,E)$ are constrained Willmore. Due to Part~(A) the corresponding Kodaira data are also constrained Willmore.
\end{proof}

\chapter[Families of flat connections for constrained Willmore maps]{Families of flat connections for constrained Willmore immersions}
\label{ch:Flat families}
The Berlin group in~\cite[Lemma~6.2]{FLPP} characterized Willmore surfaces in terms of a family of flat connections.
Bohle~\cite{Boh} generalized this characterization to constrained Willmore surfaces in $\mathbb{R}^4$.
In this chapter we explain the relation of this characterization to our characterizations in Theorems~\ref{constrained Willmore 1} and~\ref{constrained Willmore 2}.
As preparation we shall now explain in reasonable generality that for any section $\upsilon$ of the frame bundle of a complex $\qat$\=/vector bundle the corresponding flat connection extends to a family of flat connections if and only if there exists a holomorphic section $\upsilon\pa$ of the paired holomorphic $\qat$\=/vector bundle, such that the potential of the unique holomorphic structure that contains $\upsilon$ in the kernel is equal to the part of $\ls\upsilon\pa,\upsilon\rs$ that anti-commutes with the complex structure.
We shall see in Theorem~\ref{constrained Willmore 3} that any admissible constrained Willmore map $F$ is of this form, if $dF$ has no roots.

In order to construct the flat family we need to employ complex $\qat$\=/vector bundles of rank $r > 1$.
Although until now we have worked with line bundles, the definitions of Chapter~\ref{chapter:prelim} apply to any rank, with one exception.
The natural definition of the real positive definite inner product~\eqref{inner product} of the potential of a holomorphic structure has no obvious extension to the corresponding potentials acting on complex $\qat$\=/vector bundles.
We shall use a different inner product, which is used by the Berlin group~\cite[Section~6.1]{BFLPP}:
\index{Pairing!Potentials of higher rank}
\begin{gather}\label{inner product 2}
\langle \V,\V'\rangle_l=\int_{O_l}\Real\big(\tr(\ast d\Bar{z}_l\V_l\wedge d\Bar{z}_l\V'_l)\big)
\end{gather}
Unfortunately this pairing is only non-degenerate on the subspace $\pot{E}^-$. 
However, with respect to this pairing, the space $\pot{E}^+$ is dual to the space of potentials that are locally on $\SO_l$ of the form $\V'_l dz_l$ for $\V'_l \in \banach{2}(\SO_l, \qat^-)$.
Since $d\Bar{z}\V_l=\V_l dz_l$ for $d\Bar{z}_l\V_l\in\pot{E}^-$ the dual space of $\pot{E}$ is the space of all potentials of the form $\V'_l dz$ with $\V'_l\in\banach{2}(\SO_l,\qat)$.

Because of this change, the notion of paired bundles must also be modified.
If all complex conjugations of sections and potentials are replaced by the composition of complex conjugation with transposition, then the Sections~\ref{chapter:weierstrass} and~\ref{chapter:riemann roch} carry over to $\mathbb{C}$\=/vector bundles $E$.
In particular, the pairings $\lp\cdot,\cdot\rp$ and $\lh\cdot,\cdot\rh$ carry over and define analogous pairings for underlying vector bundles $E$.

The main statements in this chapter are purely local, hence we shall mostly restrict to trivial vector $\qat$\=/bundles on an open subset $\Omega\subset\mathbb{C}$.
We will omit subscripts when we take this point of view.

\begin{lemma}\label{le:flat family}
Let $\X$ be a not-necessarily compact Riemann surface and $E$ a holomorphic $\mathbb{C}$\=/vector bundle of rank $r$. 
Moreover let $\upsilon$ be a smooth global section of the frame bundle corresponding to $E_\qat$ and
\index{Potential}
\[
\big(\partial + B + \U^-\big)\;+\;\big(\delbar{E}-\V\big)
\]
the unique flat connection such that $\upsilon$ is parallel. 
Then the family of connections
\begin{gather}\label{flat family 3}
\big(\partial + \B +\B'(1-\lambda)+\U^-\big)\;+\;\big(\delbar{E}-\V^+-\V^-\lambda\big).
\end{gather}
is flat for all unimodular complex $\lambda$, if and only if there exists $r$ holomorphic sections $\upsilon\pa=(\upsilon\pa_1,\ldots,\upsilon\pa_r)$ of the paired $\qat$\=/vector bundle which obey $\B'_l dz_l+d\Bar{z}_l\V^-_l=(\B'_l+\V^-_l)dz_l=-2\upsilon_l(\Bar{\upsilon}\pa)^T_l\qk dz_l$ on the sets $\SO_l$.
\end{lemma}
\begin{remark}
The smoothness of $\upsilon$ is not a further assumption but a consequence of the existence of the section $\upsilon\pa$: First an application of a more general version of Corollary~\ref{gauge holomorphic structure} which applies to connections on $\qat$\=/vector bundles gives $\V^+=0$. A simultaneous application will preserve the relation between the potentials $\B'$ and $\V$ and the sections $\upsilon$ and $\upsilon\pa$. Afterwards an elliptic bootstrapping argument analogous to the argument in the proof of Lemma~\ref{lem:solutions are analytic} yields that the transformed sections $\upsilon$ and $\upsilon\pa$ are smooth.
\end{remark}
\begin{proof}
The statement and the proof is purely local. Hence we may assume that $\X$ is an open subset $\Omega\subset\mathbb{C}$ and that all sections are functions on $\Omega$. In particular, we shall omit the index $l$. Let us recall that the section $\upsilon$ has the following derivative on $\Omega$:
\begin{align*}
d\upsilon&=\big(d\Bar{z}(\V)-dz(\B+\U^-)\big)\upsilon.
\end{align*}
The flatness of the connection induced by $\upsilon$ is equivalent to (compare~\eqref{eq:zero curvature 2}):
\index{Zero curvature equation}
\[
(\barpartial -V^+-\U^-)(\partial +\B+\U^-)=(\partial +\B+\V^-)(\barpartial -\V).
\]
The decomposition into the parts which commute and anti-commute with the complex structure gives
\begin{gather}\label{flat family 2}\begin{aligned}
(\barpartial-\V^+)\B-(\U^-)^2+(\partial +\B)\V^++(\V^-)^2&=0,\\
(\barpartial-\V^+)\U^--\U^-\B+(\partial +\B)\V^-+\V^-\V^+&=0.
\end{aligned}\end{gather}
We first assume that there exists $\upsilon\pa$ and $\B'$ which obey
\begin{align*}
(\B'+\V^-)dz &=-2\upsilon(\Bar{\upsilon}\pa)^T\qk dz, &
d\upsilon\pa &=dz\partial \upsilon\pa+d\Bar{z}\qk\Bar{\V}^T\qk\upsilon\pa
\end{align*}
Since $\upsilon$ is a section of the frame bundle, the first equation uniquely determines the solution $\upsilon\pa$.
The second equation means that $\upsilon\pa$ is a holomorphic section of the paired bundle.
Now we calculate the exterior derivative of $\B'dz+d\Bar{z}\V^-$:
\begin{align*}
&\big(\barpartial\B'+\partial \V^-)d\Bar{z}\wedge dz
=d\big(\B'dz+d\Bar{z}\V^-\big)
=-2d\big(\upsilon(\Bar{\upsilon}\pa)^T\qk dz\big)\\
&=-2d\upsilon\wedge(\Bar{\upsilon}\pa)^T\qk dz-2\upsilon(\Bar{\upsilon}\pa)^T\qk\V\qk dz\wedge\qk dz\\
&=\big(d\Bar{z}\V-dz(B+\U^-)\big)\wedge(\B'+\V^-)dz+2\upsilon(\Bar{\upsilon}\pa)^T\qk\V d\Bar{z}\wedge dz\\
&=\Big(\big(\V^+-\U^-\big)\B'+\big(\V^--\B\big)\V^--\big(\B'+\V^-\big)\V\Big)d\Bar{z}\wedge dz\\
&=\Big( \V^+\B'-\B'\V^+-\big(\U^-\B'+\B\V^-+\V^-\V^++\B'\V^-\big)\Big)d\Bar{z}\wedge dz.
\end{align*}
The decomposition into the parts which commute and anti-commute with the complex structure gives
\begin{align}\label{flat family 1}
\barpartial \B'&=\V^+\B'-\B'\V^+,&\partial \V^-&=-\big(\U^-\B'+(\B+\B')\V^-+\V^-\V^+\big).
\end{align}
Conversely, if there exists a $\B'$ such that Equation\eqref{flat family 1} is satisfied and we define $\upsilon\pa$ through $(\B'+\V^-)dz=-2\upsilon(\Bar{\upsilon}\pa)^T\qk dz$.
Then reusing the above calculation we see that the exterior derivative of $-2\upsilon(\Bar{\upsilon}\pa)^T\qk dz$ involves besides $d\upsilon$ only the $\barpartial$ derivatives of $\upsilon\pa$.
Thus Equation\eqref{flat family 1} is satisfied only if $\upsilon\pa$ is a holomorphic section of the paired bundle, that is $\barpartial\upsilon\pa =\qk\Bar{\V}^T\qk\upsilon\pa$.

It remains to prove that Equation\eqref{flat family 1} holds if and only if the family~\eqref{flat family 3} is flat for all unimodular $\lambda$.
For fixed $\lambda\in\mathbb{C}$ the flatness of~\eqref{flat family 3} is equivalent to
\begin{align*}
(\barpartial-\V^+)\B-(\U^-)^2+(\partial +\B)\V^++(\V^-)^2\Bar{\lambda}\lambda&=\big((\barpartial-\V^+)\B'+\B'\V^+\big)(\lambda-1),\\
(\barpartial-\Bar{\V}^+)\U^--\U^-\B+(\partial +\Bar{\B})\V^-+\V^-\V^+&=\\
=\big((\partial +\B)\V^-+\V^-\V^++&\U^-\B'\big)(1-\lambda)+\B'\V^-(\Bar{\lambda}-1)\lambda.
\end{align*}
If $\lambda$ is unimodular, then $\Bar{\lambda}\lambda = 1$ and $(\Bar{\lambda}-1)\lambda = 1-\lambda$.
In that case the equations~\eqref{flat family 2} and~\eqref{flat family 1} imply that the above pair of equations holds, since both sides vanish in both equations.
Conversely, if the family of connections~\eqref{flat family 3} is for all unimodular complex $\lambda$ flat, then together with~\eqref{flat family 2} this implies that Equations~\eqref{flat family 1} hold.
\end{proof}
\begin{example}\label{example:d=1}
\index{Sinh-Gordon equation}
In this example we transform the rank-one case of the flat family of connections~\eqref{flat family 3} into a flat family of connections on a $\mathbb{C}$\=/vector bundle of rank two, which takes values in $\mathrm{su}(2,\mathbb{C})$ for unimodular parameters $\lambda$. As in the proof of the preceding Lemma we restrict to the case of a trivial bundle on an open subset $\Omega\subset\mathbb{C}$. A simultaneous application of Corollary~\ref{gauge holomorphic structure} to the sections $\upsilon$ and $\upsilon\pa$ yields $\V^+=0$ and $\U=\U^-$ in~\eqref{flat family 3}.
The $\qat^+$ functions $\B$ and $\B'$ obey $\B\qj=\qj\Bar{\B}$ and $\B'\qj=\qj\Bar{\B}'$, and the products $\qj\V$ and $\qj\U$ of $\qj$ with the $\qat^-$ functions $\U$ and $\V$ are $\mathbb{C}$\=/valued functions obeying $\overline{\qj\V}=-\Bar{\V}\qj=\V\qj$ and $\overline{\qj\U}=\U\qj$. By conceiving $\qat$ as the two-dimensional $\mathbb{C}$\=/vector space of row vectors $\qat=\big(\begin{smallmatrix}\mathbb{C}\\\qj\mathbb{C}\end{smallmatrix}\big)$, the left multiplication with the potentials $\B$, $\B'$, $\V$ and $\U$ corresponds to the left multiplication with the following complex $2\times 2$ matrices:
\begin{align*}
\begin{pmatrix}\B&0\\0&\Bar{\B}\end{pmatrix}&,&\begin{pmatrix}\B'&0\\0&\Bar{\B}'\end{pmatrix}&,&\begin{pmatrix}0&\overline{\qj\V}\\-\qj\V&0\end{pmatrix}&,&\begin{pmatrix}0&\overline{\qj\U}\\-\qj\U&0\end{pmatrix}.
\end{align*}
The real part of the potentials in~\eqref{flat family 3} viewed as $\qat$\=/valued $1$\=/forms is equal to
\begin{equation}
\label{eq:gauge real part}
\tfrac{1}{2}\Big(dz\big(B+\B'(1-\lambda)\big)+\overline{dz(B+\B'\big(1-\lambda)\big)}\Big).
\end{equation}
A gauge transformation using a real frame that is parallel with respect to the connection whose potentials are equal to~\eqref{eq:gauge real part}, preserves the flatness of the connection. Furthermore, this gauge transformation annihilates the real part and preserves the imaginary part of the potentials in~\eqref{flat family 3}, since they are invariant with respect to the conjugation by a real frame. If we replace in the gauge transformed connection $\Bar{\lambda}$ by $\lambda^{-1}$, then the corresponding potentials become equal to
\begin{align*}
&dz\begin{pmatrix}\frac{1}{2}(\B+\B')&\overline{\qj\U}\\0&-\frac{1}{2}(\B+\B')\end{pmatrix}-d\Bar{z}\begin{pmatrix}\frac{1}{2}(\Bar{\B}+\Bar{\B}')&0\\\qj\U&-\frac{1}{2}(\Bar{\B}+\Bar{\B}')\end{pmatrix}\\
+&dz\begin{pmatrix}-\frac{1}{2}\B'&0\\dz\qj\V&\frac{1}{2}\B'\end{pmatrix}\lambda-d\Bar{z}\begin{pmatrix}-\frac{1}{2}\Bar{\B}'&\overline{\qj\V}\\0&\frac{1}{2}\Bar{\B}'\end{pmatrix}\lambda^{-1}.
\end{align*}
The coefficient of $\lambda$ is called Higgs field and is the anti-Hermitian conjugate of the coefficient of $\lambda^{-1}$. 
For unimodular $\lambda$ it is a $\mathrm{su}(2,\mathbb{C})$\=/connection which coincides for such $\lambda$ with the family with $\Imag\qat$\=/valued potentials. Hence it is flat for such $\lambda$, and by $\mathbb{C}$\=/linearity for all $\lambda\in\mathbb{C}$. 
It is exactly of the form considered in~\cite[Section~1]{Hi}. Hitchin shows that for non-nilpotent Higgs fields it describes non conformal harmonic maps into the three sphere. 
The formula shows that the Higgs field is nilpotent if and only if $\B'=0$. 
This case describes conformal harmonic maps to the three sphere and corresponds to the integrable system of the $\sinh$\=/Gordon equation.

Let us now establish for $\V^+=0=\U^+=\B'$ the relation to the $\sinh$\=/Gordon equation. 
First the Equation~\eqref{flat family 1} gives $\partial\V=-\B\V$, then Equation~\eqref{flat family 1} yields $\barpartial\U=\U\B$, and these equations together imply $\barpartial(\U\V)=(\barpartial\U)\V+\U\partial\V=\U\B\V-\U\B\V=0$.
Therefore $(dz\U)(d\Bar{z}\V) = dz^2 \U \V$ is a holomorphic quadratic differential.
Away from its roots we can choose the coordinate $z$ such that it is simply $dz^2$. 
On this domain the functions $\U$ and $\V$ have no roots. 
Now we make the ansatz $\V=\qj\exp(\frac{1}{2}u)=\exp(\frac{1}{2}\Bar{u})\qj$ and $\U=-\exp(-\frac{1}{2}u)\qj$ with a complex function $u$. 
Equations~\eqref{flat family 2}-\eqref{flat family 1} imply
\begin{align*}
\barpartial\partial\Bar{u}
&=2\barpartial\big(\partial\V\V^{-1}\big)
=2\barpartial\B
=2\big(\U^2-\V^2\big)\\
&=\qj\exp(-\tfrac{1}{2}u)\qj\exp(-\tfrac{1}{2}u)-\exp(\tfrac{1}{2}u)\qj\exp(\tfrac{1}{2}u)\qj
=\sinh(\Real(u)).
\end{align*}
In particular, $\Real(u)$ is a solution of $\triangle\Real(u)=4\sinh(\Real(u))$, which is the $\sinh$\=/Gordon equation. Furthermore, $\Imag(u)$ is a harmonic function, which on simply connected domains is the imaginary part of a holomorphic function $f$. If we change the trivialization of the holomorphic line bundle $E$ by the holomorphic invertible function $\exp(-\frac{1}{2}f)$, then, due to~\eqref{eq:potential transformation}, the new potentials $\V'$ and $\U'$ are equal to
\begin{align*}
\exp(-\tfrac{1}{2}f)\V\exp(\tfrac{1}{2}f)
&=\qj\exp(\tfrac{1}{2}(f-\Bar{f})+u)=\qj\exp(\Real(u))\quad\text{and}\\
\exp(-\tfrac{1}{2}f)\U\exp(\tfrac{1}{2}f)&=\exp(\tfrac{1}{2}(\Bar{f}-f)-u)\qj=\exp(-\Real(u))\qj,
\end{align*}
respectively. 
Hence with respect to a suitable coordinate $z$ and after a suitable change of the line bundle $E$ the function $-\qj\V$ is positive and the logarithm $u=\ln((-\qj\V)^2)$ solves the $\sinh$\=/Gordon equation $\triangle u = 4\sinh(\Real(u))$.
\end{example}
In the rest of this chapter we describe the holomorphic rank two $\qat$\=/bundle that was introduced by the Berlin group.
The holomorphic structure has the important property that it is the $\barpartial$\=/part of a flat connection which can be extended to a family of flat connections in case of Willmore and constrained Willmore surfaces~\cite{Boh}.
We examine whether this point of view can be incorporated into our framework and extended to non-immersed admissible maps. Again it suffices to do the calculation locally on an open subset $\Omega\subset\mathbb{C}$.

Since the Berlin group only considers smooth conformal immersions there is no need to bring the holomorphic $\qat$\=/line bundle of the left normal into our preferred form with trivial complex structure and non-trivial flat connection. 
In order to conceive the conformal maps as maps into the projective quaternionic space, they are identified with $\qat$\=/subbundles of the trivial rank two $\qat$\=/vector bundle, and the complex structure defined by the left normal is extended to a complex structure $S$ of this rank two vector bundle (\cite[Section~4.1]{BFLPP}). 
In~\cite[Chapter~5]{BFLPP} natural conditions are derived on this extension which makes it unique. Our considerations are based on the explicit calculation of $S$ in~\cite[Section~7.2]{BFLPP}. 
In our notation with reversed signs of $\ast$ and $R$ the formula~\cite[(7.3)]{BFLPP} takes the form~\eqref{eq:mean curvature 2}. 
The comparison of this formula with the formulas~\cite[(7.11)-(7.12)]{BFLPP} with reversed signs of $\ast$ and $R$ leads to the following formula for $S$ in terms of $F$, $N$, $R$ and the mean curvature $H$:
\[
S=\begin{pmatrix}1&F\\0&1\end{pmatrix}\begin{pmatrix}N&0\\\Bar{H}N&R\end{pmatrix}\begin{pmatrix}1&-F\\0&1\end{pmatrix}.
\]
Let us now insert $\upsilon=-(\overline{\qj\chi})^{-1}=-\qj\Bar{\chi}^{-1}$, $F=\upsilon^{-1}\phi$ and $H=2\chi^{-1}\qj\U\Bar{\psi}^{-1}$ (from Equation~\eqref{eq:mean curvature}) in order to express $S$ in terms of the local functions of $\upsilon$, $\phi$, $\psi$ and $\U$:
\begin{align}
S&=\begin{pmatrix}1&\upsilon^{-1}\phi\\0&1\end{pmatrix}\begin{pmatrix}\upsilon^{-1}\qi\upsilon&0\\-2\psi^{-1}\U\upsilon\upsilon^{-1}\qi\upsilon&\psi^{-1}\qi\psi\end{pmatrix}\begin{pmatrix}1&-\upsilon^{-1}\phi\\0&1\end{pmatrix}\nonumber\\
&=\begin{pmatrix}1&\upsilon^{-1}\phi\\0&1\end{pmatrix}\hspace{-1.5mm}\begin{pmatrix}\upsilon^{-1}&0\\0&\psi^{-1}\end{pmatrix}\hspace{-1.5mm}\begin{pmatrix}1&0\\-\U&1\end{pmatrix}\hspace{-1.5mm}\begin{pmatrix}\qi&0\\0&\qi\end{pmatrix}\hspace{-1.5mm}\begin{pmatrix}1&0\\\U&1\end{pmatrix}\hspace{-1.5mm}\begin{pmatrix}\upsilon&0\\0&\psi\end{pmatrix}\hspace{-1.5mm}\begin{pmatrix}1&-\upsilon^{-1}\phi\\0&1\end{pmatrix}\hspace{-2mm}\nonumber\\
&=\Upsilon^{-1}\begin{pmatrix}\qi&0\\0&\qi\end{pmatrix}\Upsilon\quad\text{with}\quad\Upsilon=\begin{pmatrix}\upsilon&-\phi\\\U\upsilon&\psi-\U\phi\end{pmatrix}.\label{eq:intertwiner}
\end{align}
Therefore $\Upsilon$ transforms the complex structure $S$ into the trivial complex structure on the rank-two vector bundle, which is given by left multiplication with $\qi$.
We consider this matrix as an intertwiner from the rank-two $\qat$\=/vector bundle with complex structure $S$ and trivial holomorphic structure (i.e.\ the $\barpartial$ part of the trivial flat connection induced by the derivative $d$) onto a rank-two $\qat$\=/vector bundle with trivial complex structure given by left multiplication with $\qi$ and non-trivial holomorphic structure.
We see that the assumption that $dF$ has no roots is necessary due to the presence of $\psi^{-1}$ in the formula for $S$.
If $dF$ has roots then $S$ is not a smooth complex structure.
Therefore the characterization of Willmore and constrained Willmore maps in terms of a family of flat connections on this rank-two vector bundle can not be applied to admissible maps that are not immersions.

To complete the picture (returning to the assumption that $dF$ has no roots) we calculate the potentials of the non-trivial holomorphic structure, as determined by the property that the intertwiner $\Upsilon$ belongs to its kernel (compare to Lemma~\ref{lem:local holomorphic structure}):
\begin{multline*}
\barpartial \left[\begin{pmatrix}1&0\\\U&1\end{pmatrix}\begin{pmatrix}\upsilon&-\phi\\0&\psi\end{pmatrix}\right]\begin{pmatrix}\upsilon&-\phi\\0&\psi\end{pmatrix}^{-1}\begin{pmatrix}1&0\\\U&1\end{pmatrix}^{-1}\\
\begin{aligned}
&=\left[\begin{pmatrix}\barpartial&0\\\barpartial\U+\U\partial&\barpartial \end{pmatrix}\begin{pmatrix}\upsilon&-\phi\\0&\psi\end{pmatrix}\right]\begin{pmatrix}\upsilon^{-1}&\upsilon^{-1}\phi\psi^{-1}\\0&\psi^{-1}\end{pmatrix}\begin{pmatrix}1&0\\-\U&1\end{pmatrix}\\
&=\begin{pmatrix}\V\upsilon&-\V\phi\\
\barpartial\U\upsilon+\U\partial\upsilon&-\barpartial\U\phi-\U\partial\phi+\U\end{pmatrix}\begin{pmatrix}\upsilon^{-1}&\upsilon^{-1}\phi\psi^{-1}\\0&\psi^{-1}\end{pmatrix}\begin{pmatrix}1&0\\-\U&1\end{pmatrix}\\
&=\begin{pmatrix}\V&0\\
\barpartial\U+\U\partial\upsilon\upsilon^{-1}&\U(1-\upsilon\partial F\psi^{-1})\end{pmatrix}
\begin{pmatrix}1&0\\-\U&1\end{pmatrix}.\end{aligned}
\end{multline*}
We insert the formulas in Theorem~\ref{thm:darboux} for $\partial\upsilon$ and $\partial\phi$ and obtain
\[
\upsilon\partial F\psi^{-1}=-\partial\upsilon\upsilon^{-1}\phi\psi^{-1}+\upsilon\partial\phi\psi^{-1}=(\B+\U)\phi\psi^{-1}+(\psi-(\B+\U)\phi)\psi^{-1}=1.
\]
So the $(0,1)$\=/potential of $\Upsilon$ simplifies to $\hspace{3mm}\displaystyle{
\barpartial\Upsilon\Upsilon^{-1}=\begin{pmatrix}\V&0\\-\U^2-\U\B+\barpartial\U&0\end{pmatrix}}$.

In order to calculate the $(1,0)$\=/potential we use besides the potentials $\V$, $\B$ and $\U$ the $(1,0)$\=/potentials $\B'$ and $\V'$ such that $d\psi\psi^{-1}=-dz(\B'+\V')+d\Bar{z}\U$ holds (compare to~\eqref{eq:dchi}). This implies
\begin{multline*}
\partial\left[\begin{pmatrix}1&0\\\U&1\end{pmatrix}\begin{pmatrix}\upsilon&-\phi\\0&\psi\end{pmatrix}\right]\begin{pmatrix}\upsilon&-\phi\\0&\psi\end{pmatrix}^{-1}\begin{pmatrix}1&0\\\U&1\end{pmatrix}^{-1}\\
\begin{aligned}&=\left[\begin{pmatrix}\partial&0\\\partial \U+\U\barpartial&\partial\end{pmatrix}\begin{pmatrix}\upsilon&-\phi\\0&\psi\end{pmatrix}\right]\begin{pmatrix}\upsilon^{-1}&\upsilon^{-1}\phi\psi^{-1}\\0&\psi^{-1}\end{pmatrix}\begin{pmatrix}1&0\\-\U&1\end{pmatrix}\\
&=\begin{pmatrix}-(\B+\U)\upsilon&-\psi+(\B+\U)\phi\\(\partial\U+\U\V)\upsilon&-(\partial\U+\U\V)\phi-(\B'+\V')\psi\end{pmatrix}\begin{pmatrix}\upsilon^{-1}&\upsilon^{-1}\phi\psi^{-1}\\0&\psi^{-1}\end{pmatrix}\begin{pmatrix}1&0\\-\U&1\end{pmatrix}\hspace{-1mm}\\
&=\begin{pmatrix}-\B-\U&-1\\\partial \U\!+\!\U\V&-\B'\!-\!\V'\end{pmatrix}\begin{pmatrix}1&0\\-\U&1\end{pmatrix}=\begin{pmatrix}-\B&-1\\\U\V+\V'\U+\partial\U+\B'\U&-\B'-\V'\end{pmatrix}.
\end{aligned}
\end{multline*}
Observe that the entries of these potentials are $\banach{2}$, except for the bottom left entry.
For constrained Willmore immersions, we know that the potentials are smooth.

If we replace in Lemma~\ref{le:flat family} the $(0,1)$\=/potential $\V^+$ and the $(1,0)$\=/potential $\B$, which both commute with the complex structure, and the $(0,1)$\=/potential $\V^-$ and $(1,0)$\=/potential $\U^-$, which both anti-commute with the complex structure, by the corresponding potentials of $\Upsilon$ as listed below, then the Equations~\eqref{flat family 2} describe the flatness of the connections induced by this frame $\Upsilon$:
\begin{gather}\label{2x2potentials}
\begin{pmatrix}0&0\\-\U^2&0\end{pmatrix},\begin{pmatrix}\B&1\\-\U\V\!-\!\V'\U&\B'\end{pmatrix},\begin{pmatrix}\V&0\\ \barpartial\U\!-\!\U\B&0\end{pmatrix},\begin{pmatrix}0&0\\-\partial\U\!-\!\B'\U&\V'\end{pmatrix}.
\end{gather}
Hence the flatness of the connection induced by the frame $\Upsilon$ is equivalent to the Equations~\eqref{flat family 2} with the former replacements. They take the following form:
\begin{gather*}
\left(\barpartial+\begin{pmatrix}0&0\\\U^2&0\end{pmatrix}\right)\begin{pmatrix}\B&1\\-\U\V-\V'\U&\B'\end{pmatrix}-\begin{pmatrix}0&0\\-\partial\U-\B'\U&\V'\end{pmatrix}^2\\
=\left(\partial+\begin{pmatrix}\B&1\\-\U\V-\V'\U&\B'\end{pmatrix}\right)\begin{pmatrix}0&0\\\U^2&0\end{pmatrix}-\begin{pmatrix}\V&0\\ \barpartial \U-\U\B&0\end{pmatrix}^2,\\
\left(\barpartial+\begin{pmatrix}0&0\\\U^2&0\end{pmatrix}\right)\begin{pmatrix}0&0\\-\partial\U-\B'\U&\V'\end{pmatrix}-\begin{pmatrix}0&0\\
-\partial\U-\B'\U&\V'\end{pmatrix}\begin{pmatrix}\B&1\\-\U\V-\V'\U&\B'\end{pmatrix}\\
=-\left(\partial+\begin{pmatrix}\B&1\\-\U\V-\V'\U&\B'\end{pmatrix}\right)\begin{pmatrix}\V&0\\ \barpartial \U-\U\B&0\end{pmatrix}-\begin{pmatrix}\V&0\\ \barpartial \U-\U\B&0\end{pmatrix}\begin{pmatrix}0&0\\\U^2&0\end{pmatrix}.
\end{gather*}
The corresponding equations on the three entries with the exception of the bottom left entries are quite simple and take the following form
\begin{gather}\label{eq:extended flatness}\begin{aligned}
\barpartial \B-\U^2+\V^2&=0,&0&=0,&\barpartial\B'+\U^2-{\V'}^2&=0,\\
\partial \V+\B\V+\barpartial \U-\U\B&=0,&0&=0,&\barpartial\V'+\partial\U+\B'\U-\V'\B'&=0.
\end{aligned}\end{gather}
Both sides of the equations on the bottom left entries take the following form:
\begin{align*}
-\barpartial(\U\V+\V'\U)+\U^2\B+\V'(\partial \U+\B'\U)&=
\partial(\U^2)+\B'\U^2-(\barpartial\U-\U\B)\V,\\
(\partial\U\!+\!\B'\U)\B\!-\!\barpartial(\partial\U\!+\!\B'\U)\!+\!\V'(\U\V\!+\!\V'\U)&=(\partial\!+\!\B')(\U\B\!-\!\barpartial\U)\!+\!(\U\V\!+\!\V'\U)\V.
\end{align*}
These equations are equivalent to
\begin{align*}
-\U(\partial \V-\U\B+\barpartial \U+\B\V)&=(\barpartial\V'-\V'\B'+\partial\U+\B'\U)\U,\\
-(\barpartial\B'+\U^2-\V'^2)\U&=\U(\barpartial \B-\U^2+\V^2).
\end{align*}
This shows that the flatness follows from the four non-trivial equations~\eqref{eq:extended flatness}.
They are equivalent to the flatness of the connections which are induced by the non-vanishing section $\upsilon$ of $E_\qat$ and $\psi$ of $KE_\qat$, respectively.
Furthermore, these four equations also imply that with the exception of the bottom left entries both sides of the former flatness equations vanish separately. This is also true for both sides of the bottom left entry of the second flatness equation, if and only if in addition to the four non-trivial equations in~\eqref{eq:extended flatness} the following equation holds:
\[
-\barpartial\partial \U-\B'\barpartial \U+\partial\U\B+\B'\U\B+\U^3+\V'\U\V=0.
\]
In this case it is possible to multiply either of the anti-commuting parts of the matrix potentials in~\eqref{2x2potentials} by a unimodular complex scalar to produce a family of flat connections.
This is carried out in~\cite[Lemma~6.2]{FLPP} and corresponds to Willmore immersions.

A generalization of this condition to constrained Willmore immersions was derived by Bohle in~\cite[Section~2.4]{Boh}.
Let us now formulate Bohle's approach in the present setting: 
The constrained Weierstraß triples $(\U,\chi,\psi)$ without roots of $\chi$ and $\psi$ are characterized by a solution $\B'':\Omega\to\qat^+$ of
\begin{align}\label{eq:flatness of family}
\bar\partial\B''&=0,&\partial\barpartial\U\!+\!\B'\barpartial\U\!-\!\B\partial\U\!-\!\B'\U\B\!-\!\U^3\!-\!\V'\U\V&=-\B''\V\!-\!\V'\B''.
\end{align}
We shall prove that the existence of such a $\B''$ characterizes constrained Willmore maps. First we claim that for given solutions $\U$, $\V$, $\B$, $\V'$ and $\B'$ of~\eqref{eq:extended flatness}, the Equations~\eqref{eq:flatness of family} on $\B''$ are exactly the Equations~\eqref{flat family 1} with the potentials~\eqref{2x2potentials} acting on the rank two complex $\qat$\=/vector bundle, whose complex structure is induced by $S$. In fact, the intertwiner $\Upsilon$~\eqref{eq:intertwiner} transforms this vector bundle with trivial flat connection into the rank two $\qat$\=/vector bundle with complex structure given by left multiplication with $\qi$ and flat connection whose the rank two potentials $\V^+$, $\B$, $\V^-$ and $\U^-$ in~\eqref{flat family 3} are equal to the potentials listed in~\eqref{2x2potentials}. If the additional rank two potential $\B'$ in~\eqref{flat family 3} has the form $\big(\begin{smallmatrix}0&0\\\B''&0\end{smallmatrix}\big)$ with $\B'':\Omega\to\qat^-$, then the Equations~\eqref{flat family 1} take the form
\begin{gather*}
\bar\partial\begin{pmatrix}0&0\\\B''&0\end{pmatrix}=\begin{pmatrix}0&0\\B''&0\end{pmatrix}\begin{pmatrix}0&0\\\U^2&0\end{pmatrix}-\begin{pmatrix}0&0\\\U^2&0\end{pmatrix}\begin{pmatrix}0&0\\B''&0\end{pmatrix},\\
\left(\partial+\begin{pmatrix}\B&1\\\B''\!-\!\U\V\!-\!\V'\U&\B'\end{pmatrix}\!\right)\begin{pmatrix}\V&0\\ \barpartial\U\!-\!\U\B&0\end{pmatrix}=\begin{pmatrix}0&0\\\partial\U\!+\!\B'\U&-\V'\end{pmatrix}\begin{pmatrix}0&0\\\B''&0\end{pmatrix}.
\end{gather*}
The first equation is equivalent to $\barpartial\B''=0$ and the second to
\begin{align*}
\partial\V+\B\V+\barpartial\U-\U\B&=0,&0&=0,\\\partial\barpartial\U-\partial(\U\B)+(\B''-\U\V-\V'\U)\V+\B'(\barpartial\U-\U\B)&=-\V'\B'',&0&=0. 
\end{align*}
All together these equations are equivalent to~\eqref{eq:flatness of family}, if we assume that the Equations~\eqref{eq:extended flatness} hold.
In particular, the family~\eqref{flat family 3} takes the form
\begin{gather}\label{flat family 4}\begin{aligned}
&\partial+\begin{pmatrix}\B&1\\-\U\V-\V'\U&\B'\end{pmatrix}+\begin{pmatrix}0&0\\\B''&0\end{pmatrix}(1-\lambda)-\begin{pmatrix}0&0\\\partial\U+\B'\U&-\V'\end{pmatrix}\\+\;&\barpartial+\begin{pmatrix}0&0\\\U^2&0\end{pmatrix}-\begin{pmatrix}\V&0\\\bar\partial\U-\U\B&0\end{pmatrix}\lambda.
\end{aligned}\end{gather}
Let us state the characterization of constrained Willmore Weierstraß triples $(\U,\chi,\psi)$ without roots of $\chi$ and $\psi$ in a manner which does not refer to the additional constraints $\chi\pa$ and $\psi\pa$:
\begin{theorem}\label{constrained Willmore 3}
Let $(\V,\upsilon,\phi)$ be on a compact Riemann surface $\X$ the Kodaira data of an admissible map $F=\upsilon^{-1}\phi$ with $\V\in\pot{E}^-$ and without roots of $dF$. Then $F$ is constrained Willmore if and only if there exists in addition to the potentials $\B_l+\U_l=-(\partial_l\upsilon_l)\upsilon_l^{-1}$ and $\B'_l+\V'_l=-(\partial_l\psi_l)\psi_l^{-1}$ of the Weierstra{\ss} section $\psi_l=(\partial_l+\B_l+\U_l)\phi_l$ a solution $\B''$ of~\eqref{eq:flatness of family}. This condition is locally equivalent to the flatness of the connection~\eqref{flat family 4} for all unimodular complex $\lambda$.
\end{theorem}
\begin{proof}
Because the statement and the proof is local, we may represent all involved sections by local functions on an open subset $\Omega\subset\mathbb{C}$. Let us first show that the constrained Willmore Kodaira triples $(\V,\upsilon,\phi$) with $\V\in\pot{E}^-$, which are characterized in Theorem~\ref{constrained Willmore 1} obey the condition in the theorem. Due to Theorem~\ref{constrained Willmore 1} there exists solutions $\upsilon\pa$ and $\phi\pa$ of
\begin{align*}
(\barpartial-\qk\Bar{\V}\qk)\upsilon\pa&=0,&(\barpartial-\qk\Bar{\V}\qk)\phi\pa&=0,&d\Bar{z}\V&=-2d\Bar{z}\qk(\upsilon\pa\Bar{\upsilon}+\phi\pa\Bar{\phi}).
\end{align*}
The exterior derivative of the $1$\=/form $-d\Bar{z}\V$ is
\begin{align*}
d\Bar{z}\wedge dz\partial\V& =-d\big(d\Bar{z}\V\big)
=2d\Bar{z}\qk\wedge d\big(-\upsilon\pa\Bar{\upsilon}-\phi\pa\Bar{\phi}\big) \\
&=-2d\Bar{z}\qk\wedge d\Bar{z}\qk\Bar{\V}\qk\big(\upsilon\pa\Bar{\upsilon}+\phi\pa\Bar{\phi}\big)-2d\Bar{z}\qk\wedge\big(\upsilon\pa d\Bar{\upsilon}+\phi\pa d\Bar{\phi}\big)\\
&=-2d\Bar{z}\wedge dz\V\qk\big(\upsilon\pa\Bar{\upsilon}+\phi\pa\Bar{\phi}\big)-2d\Bar{z}\qk\left(\upsilon\pa\Bar{\upsilon}\wedge\overline{d\upsilon\upsilon^{-1}}+\phi\pa\Bar{\phi}\wedge\overline{d\phi\phi^{-1}}\right)\\
&=d\Bar{z}\wedge dz\V^2-2d\Bar{z}\qk\left(\big(\upsilon\pa\Bar{\upsilon}+\phi\pa\Bar{\phi}\big)\wedge\big(\overline{d\Bar{z}\V-dz(\B+\U)}\big)+\phi\pa\Bar{\psi}\wedge d\Bar{z}\right)\\
&=d\Bar{z}\wedge dz\V^2+d\Bar{z}\V\wedge\big(-\V dz-\Bar{\B}d\Bar{z}+\U d\Bar{z}\big)-2d\Bar{z}\qk\phi\pa\Bar{\psi}\wedge d\Bar{z}\\
&=d\Bar{z}\wedge dz\left(\V^2-\V^2-\B\V-2\qk\big(\phi\pa\Bar{\psi}\big)^+\right)
\end{align*}
This implies $2\qk\big(\phi\pa\Bar{\psi}\big)^+=-\partial\V-\B\V$.
By inserting the left bottom equation of~\eqref{eq:extended flatness} we obtain
\[
2\qk\big(\phi\pa\Bar{\psi}\big)^+=\barpartial \U-\U\B.
\]
The first of the following pairs of equations is equivalent to the second:
\begin{align*}
\V&=-2\qk\big(\upsilon\pa\Bar{\upsilon}+\phi\pa\Bar{\phi}\big),&\barpartial\U-\U\B&=2\qk\big(\phi\pa\Bar{\psi}\big)^+,\\
\V&=-2\big(\upsilon\Bar{\upsilon}\pa+\phi\Bar{\phi}\pa\big)\qk,&\barpartial\U-\U\B&=2\big(\psi\Bar{\phi}\pa\big)^+\qk.
\end{align*}
Therefore the former intertwiner $\Upsilon$~\eqref{eq:intertwiner} obeys
\begin{gather*}
-2\begin{pmatrix}\upsilon&-\phi\\\U\upsilon&\psi-\U\phi\end{pmatrix}\overline{\begin{pmatrix}\upsilon\pa&-\phi\pa\\0&0\end{pmatrix}}^T\hspace{-1mm}\qk dz=-2\begin{pmatrix}1&0\\\U&1\end{pmatrix}\begin{pmatrix}\upsilon&-\phi\\0&\psi\end{pmatrix}\begin{pmatrix}\Bar{\upsilon}\pa&0\\-\Bar{\phi}\pa&0\end{pmatrix}\qk dz\\
=-2\begin{pmatrix}1&0\\\U&1\end{pmatrix}\begin{pmatrix}\upsilon\Bar{\upsilon}\pa+\phi\Bar{\phi}\pa&0\\-\psi\Bar{\phi}\pa&0\end{pmatrix}\qk dz=\begin{pmatrix}\V&0\\\U\V+2\big(\psi\Bar{\phi}\pa\big)^-\qk+\barpartial \U-\U\B&0\end{pmatrix}dz.
\end{gather*}
In particular, this section of the corresponding frame bundle together with the two sections $\big(\begin{smallmatrix}\upsilon\pa\\0\end{smallmatrix}\big)$ and  $\big(\begin{smallmatrix}-\phi\pa\\0\end{smallmatrix}\big)$ of the paired rank-two $\qat$\=/vector bundle obeys one of the two equivalent conditions in Lemma~\ref{le:flat family} with the rank two potential $\B'$ defined as $(\begin{smallmatrix}0&0\\\B''&0\end{smallmatrix})$ with $\B''=\U\V+2(\psi\Bar{\phi})^-\qk$.
Now the condition in the theorem follows from that Lemma, since for the flat connection with parallel $\Upsilon$ the Equations~\eqref{flat family 1} take the form~\eqref{eq:flatness of family}.
Furthermore, due to Lemma~\ref{le:flat family} this equation is equivalent to the flatness of the family~\eqref{eq:flatness of family} for all complex unimodular $\lambda$.

To prove the converse we assume that for a given Kodaira triple $(\V,\upsilon,\phi)$ of an admissible immersion with $\V\in\pot{E}^-$ and with potentials $\B+\U=-\partial\upsilon\upsilon^{-1}$ and $\B'+\V'=-\partial\psi\psi^{-1}$ there exists a solution $\B''$ of~\eqref{eq:flatness of family}.
We have seen above that in this case the family~\eqref{flat family 4} is flat for all complex unimodular $\lambda$.
Therefore Lemma~\ref{le:flat family} implies that the sum of potentials $\big(\begin{smallmatrix}0&0\\\B''&0\end{smallmatrix}\big)+\big(\begin{smallmatrix}\V&0\\\partial \U-\B\U&0\end{smallmatrix}\big)$ can be expressed as a pairing of the intertwiner $\Upsilon$~\eqref{eq:intertwiner} with two holomorphic sections of the corresponding paired holomorphic vector bundle.
Since the second column of this sum vanishes, these sections are of the form $\big(\begin{smallmatrix}\upsilon\pa\\0\end{smallmatrix}\big)$ and  $\big(\begin{smallmatrix}-\phi\pa\\0\end{smallmatrix}\big)$, respectively, with unique holomorphic sections $\upsilon\pa$ and $\phi\pa$ paired with $\upsilon$ and $\phi$.
Since the top left entry of this sum is $\V$, these holomorphic sections $(\upsilon\pa,\phi\pa)$ obey the condition in Theorem~\ref{constrained Willmore 1} and $(\V,\upsilon,\phi)$ is constrained Willmore.
\end{proof}

We close this chapter and the book with a prospective view. 
At the end of Chapter~\ref{chapter:weierstrass} we described how the Darboux transformation should be generalized to quaternionic $d$\=/dimensional spaces of holomorphic sections of paired holomorphic $\qat$\=/line bundles. 
Furthermore, in Lemma~\ref{lem:solutions are analytic} we considered an analogue to constrained Willmore Kodaira triples or Weierstra{\ss} triples for tuples with more than two holomorphic sections of paired holomorphic $\qat$\=/line bundles. 
This book mainly investigated the case $d=2$, which corresponds to admissible maps. In the Example~\ref{example:d=1} we have seen that for $d=1$ the analogue to the constrained Willmore case is exactly the description of harmonic maps to the $3$\=/sphere investigated by Hitchin in~\cite{Hi}. 
Therefore we expect that for $d>2$ there also exist interesting geometric interpretations of these nonlinear elliptic equations.

\backmatter
\let\addcontentslineOriginal\addcontentsline
\renewcommand{\addcontentsline}[3]{}


\bibliographystyle{alpha}
\bibliography{embdededbib.bib}%

\newcommand{\etalchar}[1]{$^{#1}$}
\begin{thebibliography}{BHPP97}

\bibitem[ABR01]{ABR}
Sheldon Axler, Paul Bourdon, and Wade Ramey.
\newblock {\em Harmonic {{Function Theory}}}.
\newblock Number 137 in Graduate Texts in Mathematics. Springer, New York,
  2001.

\bibitem[AF08]{Ad}
Robert~A. Adams and John J.~F. Fournier.
\newblock {\em Sobolev Spaces}.
\newblock Number 140 in Pure and Applied Mathematics. Academic Press,
  Amsterdam, 2008.

\bibitem[Ahl79]{Ah}
Lars~V. Ahlfors.
\newblock {\em Complex Analysis: An Introduction to the Theory of Analytic
  Functions of One Complex Variable}.
\newblock International Series in Pure and Applied Mathematics. McGraw-Hill,
  New York, 1979.

\bibitem[Aub98]{Au}
Thierry Aubin.
\newblock {\em Some {{Nonlinear Problems}} in {{Riemannian Geometry}}}.
\newblock Springer {{Monographs}} in {{Mathematics}}. Springer, Berlin, 1998.

\bibitem[Aud94]{Aud}
Mich{\`e}le Audin.
\newblock Symplectic and almost complex manifolds.
\newblock In Mich{\`e}le Audin and Jacques Lafontaine, editors, {\em
  Holomorphic Curves in Symplectic Geometry}, pages 41--74. Birkh{\"a}user
  Basel, Basel, 1994.

\bibitem[BB13]{BB}
David~D. Bleecker and Bernhelm Boo{\ss}.
\newblock {\em Index Theory with Applications to Mathematics and Physics}.
\newblock International press, Somerville, 2013.

\bibitem[BE53]{MOT}
Harry Bateman and Arthur Erd{\'e}lyi.
\newblock {\em Higher Transcendental Functions. {{Vol}} 2}.
\newblock McGraw-Hill, New York, NY, 1953.

\bibitem[BFL{\etalchar{+}}02]{BFLPP}
Francis~E. Burstall, Dirk Ferus, Katrin Leschke, Franz Pedit, and Ulrich
  Pinkall.
\newblock {\em Conformal {{Geometry}} of {{Surfaces}} in {{S4}} and
  {{Quaternions}}}.
\newblock Number 1772 in Lecture {{Notes}} in {{Mathematics}}. Springer,
  Berlin, 2002.

\bibitem[BHPP97]{Burstall1997a}
F.~Burstall, U.~{Hertrich-Jeromin}, F.~Pedit, and U.~Pinkall.
\newblock Curved flats and isothermic surfaces.
\newblock {\em Mathematische Zeitschrift}, 225(2):199--209, June 1997.

\bibitem[BK03]{BK}
Matthias Bauer and Ernst Kuwert.
\newblock Existence of minimizing {{Willmore}} surfaces of prescribed genus.
\newblock {\em International Mathematics Research Notices}, 2003(10):553, 2003.

\bibitem[BLPP12]{BLPP}
Christoph Bohle, Katrin Leschke, Franz Pedit, and Ulrich Pinkall.
\newblock Conformal maps from a 2-torus to the 4-sphere.
\newblock {\em Journal f{\"u}r die reine und angewandte Mathematik (Crelles
  Journal)}, 2012(671):1--30, 2012.

\bibitem[Boh10]{Boh}
Christoph Bohle.
\newblock Constrained {{Willmore}} tori in the 4-sphere.
\newblock {\em Journal of Differential Geometry}, 86(1):71--132, 2010.

\bibitem[Boh12]{Bohle2012a}
Christoph Bohle.
\newblock Constant mean curvature tori as stationary solutions to the
  {{Davey}}--{{Stewartson}} equation.
\newblock {\em Mathematische Zeitschrift}, 271(1-2):489--498, June 2012.

\bibitem[Bou09]{Bo}
Nicolas Bourbaki.
\newblock {\em Elements of Mathematics. {{Algebra I}}}.
\newblock Springer, Berlin, 2009.

\bibitem[BPP02]{Burstall2002b}
Francis Burstall, Franz Pedit, and Ulrich Pinkall.
\newblock Schwarzian derivatives and flows of surfaces.
\newblock In {\em Contemporary {{Mathematics}}}, volume 308, pages 39--61.
  American Mathematical Society, Providence, Rhode Island, 2002.

\bibitem[BPP08]{BPP}
Christoph Bohle, G.~Paul Peters, and Ulrich Pinkall.
\newblock Constrained {{Willmore}} surfaces.
\newblock {\em Calculus of Variations and Partial Differential Equations},
  32(2):263--277, 2008.

\bibitem[Bry84]{Bryant1984}
Robert~L. Bryant.
\newblock A duality theorem for {{Willmore}} surfaces.
\newblock {\em Journal of Differential Geometry}, 20(1):23--53, 1984.

\bibitem[BS88]{BS}
Colin Bennett and Robert~C. Sharpley.
\newblock {\em Interpolation of Operators}.
\newblock Number 129 in Pure and Applied Mathematics. Academic Press, Boston,
  1988.

\bibitem[BS18]{BSa}
Theo B{\"u}hler and D.~Salamon.
\newblock {\em Functional Analysis}.
\newblock Number volume 191 in Graduate Studies in Mathematics. American
  Mathematical Society, Providence, Rhode Island, 2018.

\bibitem[Car39]{Ca}
Torsten Carleman.
\newblock Sur un probleme d'unicite pour les sytemes d'equations aux derivees
  partielles a doux variables independantes.
\newblock {\em Arkiv f{\"o}r Matematik, Astronomi och Fysik}, 26(17):1--9,
  1939.

\bibitem[CGS95]{Cieslinski1995}
Jan Cie{\'s}li{\'n}ski, Piotr Goldstein, and Antoni Sym.
\newblock Isothermic surfaces in {{E3}} as soliton surfaces.
\newblock {\em Physics Letters A}, 205(1):37--43, 1995.

\bibitem[Con95]{Co2}
John~B. Conway.
\newblock {\em Functions of {{One Complex Variable II}}}.
\newblock Number 159 in Graduate Texts in Mathematics. Springer, New York,
  1995.

\bibitem[CP98]{CP}
Michael Cwikel and Evgeniy Pustylnik.
\newblock Sobolev type embeddings in the limiting case.
\newblock {\em The Journal of Fourier Analysis and Applications},
  4(4-5):433--446, 1998.

\bibitem[DN55]{DoNi}
Avron Douglis and Louis Nirenberg.
\newblock Interior estimates for elliptic systems of partial differential
  equations.
\newblock {\em Communications on Pure and Applied Mathematics}, 8(4):503--538,
  1955.

\bibitem[dR84]{dR}
Georges de~Rham.
\newblock {\em Differentiable Manifolds: Forms, Currents, Harmonic Forms}.
\newblock Number 266 in Die {{Grundlehren}} Der Mathematischen
  {{Wissenschaften}} in {{Einzeldarstellungen}}. Springer, Berlin, 1984.

\bibitem[Eis09]{Ei}
Luther~Pfahler Eisenhart.
\newblock {\em A Treatise on the Differential Geometry of Curves and Surfaces}.
\newblock {Ginn and Company}, Boston, MA, 1909.

\bibitem[EK82]{EK}
Fritz Ehlers and Horst Kn{\"o}rrer.
\newblock An algebro-geometric interpretation of the
  {{B{\"a}cklund-transformation}} for the {{Korteweg-de Vries}} equation.
\newblock {\em Commentarii Mathematici Helvetici}, 57(1):1--10, 1982.

\bibitem[Fed96]{Federer1996}
Herbert Federer.
\newblock {\em Geometric {{Measure Theory}}}.
\newblock Classics in {{Mathematics}}. Springer Berlin Heidelberg, Berlin,
  Heidelberg, 1996.

\bibitem[FK92]{FK}
Hershel~M. Farkas and Irwin Kra.
\newblock {\em Riemann Surfaces}.
\newblock Number~71 in Graduate Texts in Mathematics. Springer, New York, 1992.

\bibitem[FLPP01]{FLPP}
D.~Ferus, K.~Leschke, F.~Pedit, and U.~Pinkall.
\newblock Quaternionic holomorphic geometry: {{Pl{\"u}cker}} formula, {{Dirac}}
  eigenvalue estimates and energy estimates of harmonic 2-tori.
\newblock {\em Inventiones Mathematicae}, 146(3):507--593, 2001.

\bibitem[For81]{Fo}
Otto Forster.
\newblock {\em {Lectures on Riemann surfaces}}.
\newblock Number~81 in Graduate Texts in Mathematics. Springer, New York, 1981.

\bibitem[Fri98]{Fr2}
Thomas Friedrich.
\newblock On the {{Spinor Representation}} of {{Surfaces}} in {{Euclidean}}
  3-{{Space}}.
\newblock {\em Journal of Geometry and Physics}, 28(1-2):143--157, 1998.

\bibitem[Fri00]{Fr1}
Thomas Friedrich.
\newblock {\em Dirac {{Operators}} in {{Riemannian Geometry}}}.
\newblock Number~25 in Graduate Studies in Mathematics. American Mathematical
  Society, Providence, Rhode Island, 2000.

\bibitem[Gas06]{Gasper2006}
George Gasper.
\newblock Formulas of the dirichlet-mehler type.
\newblock In {\em Fractional Calculus and Its Applications: {{Proceedings}} of
  the International Conference Held at the University of New Haven, June 1974},
  pages 207--215. Springer, 2006.

\bibitem[GH94]{GrHa}
Phillip Griffiths and Joe Harris.
\newblock {\em Principles of Algebraic Geometry}.
\newblock J. Wiley \& Sons, New York, NY, 1994.

\bibitem[GHZ05]{GHZ}
Chaohao Gu, Hesheng Hu, and Zixiang Zhou.
\newblock {\em Darboux {{Transformations}} in {{Integrable Systems}}}.
\newblock Springer, Dordrecht, 2005.

\bibitem[GR09]{GuRo}
Robert~C. Gunning and Hugo Rossi.
\newblock {\em Analytic Functions of Several Complex Variables}.
\newblock Prentice-Hall, Englewood Cliffs, NJ, 2009.

\bibitem[Heb96]{Heb}
Emmanuel Hebey.
\newblock {\em Sobolev Spaces on {{Riemannian}} Manifolds}.
\newblock Number 1635 in Lecture Notes in Mathematics. Springer, Berlin, 1996.

\bibitem[Hel01]{Heg}
Sigurdur Helgason.
\newblock {\em Differential {{Geometry}}, {{Lie Groups}}, and {{Symmetric
  Spaces}}}.
\newblock Number~34 in Graduate Studies in Mathematics. American Mathematical
  Society, Providence, Rhode Island, 2001.

\bibitem[H{\'e}l02]{Hel}
Fr{\'e}d{\'e}ric H{\'e}lein.
\newblock {\em Harmonic Maps, Conservation Laws, and Moving Frames}.
\newblock Number 150 in Cambridge Tracts in Mathematics. Cambridge University
  Press, Cambridge, 2002.

\bibitem[Hit90]{Hi}
Nigel~J. Hitchin.
\newblock Harmonic maps from a 2-torus to the 3-sphere.
\newblock {\em Journal of Differential Geometry}, 31:627--710, 1990.

\bibitem[HO83]{Hoffman1983}
David~A. Hoffman and Robert Osserman.
\newblock The {{Gauss}} map of surfaces in {{Rn}}.
\newblock {\em J. Differential Geom}, 18(4):733--754, 1983.

\bibitem[H{\"o}r64]{Ho}
Lars H{\"o}rmander.
\newblock {\em {Linear Partial Differential Operators}}.
\newblock Number 116 in {Grundlehren der mathematischen Wissenschaften}.
  Springer, Berlin, 1964.

\bibitem[HR79]{Hewitt1979}
Edwin Hewitt and Kenneth~A. Ross.
\newblock {\em Abstract {{Harmonic Analysis}}: {{Volume I}}, {{Structure}} of
  {{Topological Groups Integration Theory Group Representations}}}.
\newblock Number 115 in Die {{Grundlehren}} Der {{Mathematischen
  Wissenschaften}}. Springer, New York, 2 edition, 1979.

\bibitem[Huy05]{Hu}
Daniel Huybrechts.
\newblock {\em Complex {{Geometry}}}.
\newblock Universitext. Springer, Berlin, 2005.

\bibitem[Jer86]{Je}
David Jerison.
\newblock Carleman inequalities for the {{Dirac}} and {{Laplace}} operators and
  unique continuation.
\newblock {\em Advances in Mathematics}, 62(2):118--134, 1986.

\bibitem[Jos11]{Jo}
J{\"u}rgen Jost.
\newblock {\em Riemannian {{Geometry}} and {{Geometric Analysis}}}.
\newblock Universitext. Springer, Berlin, 2011.

\bibitem[Kau99]{Ka}
Robert Kaufman.
\newblock Sobolev spaces, dimension, and random series.
\newblock {\em Proceedings of the American Mathematical Society},
  128(2):427--431, 1999.

\bibitem[Kim95]{Ki2}
Yonne~Mi Kim.
\newblock Carleman inequalities for the {{Dirac}} operator and strong unique
  continuation.
\newblock {\em Proceedings of the American Mathematical Society},
  123(7):2103--2112, 1995.

\bibitem[KL12]{KL}
Ernst Kuwert and Yuxiang Li.
\newblock {$W^{2,2}$-conformal immersions of a closed Riemann surface into
  $\mathbb{R}^n$}.
\newblock {\em Communications in Analysis and Geometry}, 20(2):313–340, 2012.

\bibitem[Kna86]{Kn}
Antony~W. Knapp.
\newblock {\em Representation Theory of Semisimple Groups: An Overview Based on
  Examples (PMS-36)}.
\newblock Princeton University Press, revised edition, 1986.

\bibitem[Kon92]{Kon1}
B.~G. Konopelchenko.
\newblock {\em Introduction to {{Multidimensional Integrable Equations}}: {{The
  Inverse Spectral Transform}} in 2+1 {{Dimensions}}}.
\newblock Springer, Boston, 1992.

\bibitem[Kon96]{Kon2}
B.~G. Konopelchenko.
\newblock Induced {{Surfaces}} and {{Their Integrable Dynamics}}.
\newblock {\em Studies in Applied Mathematics}, 96(1):9--51, 1996.

\bibitem[KS13]{KS}
Ernst Kuwert and Reiner Sch{\"a}tzle.
\newblock Minimizers of the {{Willmore}} functional under fixed conformal
  class.
\newblock {\em Journal of Differential Geometry}, 93(3), 2013.

\bibitem[KZPS76]{KZPS}
M.~A. Krasnoselskii, P.~P. Zabreyko, E.~I. Pustylnik, and P.~E. Sobolevski.
\newblock {\em Integral {{Operators}} in {{Spaces}} of {{Summable Functions}}}.
\newblock Noordhoff International Publishing, Leyden, 1976.

\bibitem[LT73]{LT}
Joram Lindenstrauss and Lior Tzafriri.
\newblock {\em Classical {{Banach}} Spaces}.
\newblock Number 338 in Lecture Notes in Mathematics. Springer, Berlin, 1973.

\bibitem[LY82]{LY}
Peter Li and Shing-Tung Yau.
\newblock A new conformal invariant and its applications to the {{Willmore}}
  conjecture and the first eigenvalue of compact surfaces.
\newblock {\em Inventiones Mathematicae}, 69(2):269--291, 1982.

\bibitem[Man94]{Ma}
Niculae Mandache.
\newblock Some remarks concerning unique continuation for the {{Dirac}}
  operator.
\newblock {\em Letters in Mathematical Physics}, 31(2):85--92, 1994.

\bibitem[Mor58]{Mo}
Charles~B. Morrey.
\newblock On the {{Analyticity}} of the {{Solutions}} of {{Analytic Non-Linear
  Elliptic Systems}} of {{Partial Differential Equations}}: {{Part~I}}.
  {{Analyticity}} in the {{Interior}}.
\newblock {\em American Journal of Mathematics}, 80(1):198--218, 1958.

\bibitem[Mun75]{Mun}
James~Raymond Munkres.
\newblock {\em Topology: A First Course}.
\newblock Prentice-Hall, Englewood Cliffs, NJ, 1975.

\bibitem[Nar92]{Na}
Raghavan Narasimhan.
\newblock {\em Compact {{Riemann}} Surfaces}.
\newblock Lectures in Mathematics {{ETH Z{\"u}rich}}. Birkh{\"a}user Verlag,
  Basel ; Boston, 1992.

\bibitem[Nic21]{Ni}
Liviu~I. Nicolaescu.
\newblock {\em Lectures on the Geometry of Manifolds}.
\newblock World Scientific, New Jersey, 2021.

\bibitem[NN57]{NN}
A.~Newlander and L.~Nirenberg.
\newblock Complex {{Analytic Coordinates}} in {{Almost Complex Manifolds}}.
\newblock {\em The Annals of Mathematics}, 65(3):391, 1957.

\bibitem[O'N63]{O}
Richard O'Neil.
\newblock Convolution operators and {{L}}(p,q) spaces.
\newblock {\em Duke Mathematical Journal}, 30(1):129--142, 1963.

\bibitem[Pal65]{Pa}
Richard~S. Palais.
\newblock {\em Seminar on the {{Atiyah-Singer}} Index Theorem}.
\newblock Number~57 in Annals of {{Mathematics Studies}}. Princeton University
  Press, Princeton, NJ, 1965.

\bibitem[Pal88]{Pal}
Bennett Palmer.
\newblock Isothermic surfaces and the {{Gauss}} map.
\newblock {\em Proceedings of the American Mathematical Society},
  104(3):876--884, 1988.

\bibitem[PG24]{Pinkall2024}
Ulrich Pinkall and Oliver Gross.
\newblock {\em Differential {{Geometry}}: {{From Elastic Curves}} to {{Willmore
  Surfaces}}}.
\newblock Compact {{Textbooks}} in {{Mathematics}}. Springer International
  Publishing, Cham, 2024.

\bibitem[PP98]{PP}
Franz Pedit and Ulrich Pinkall.
\newblock Quaternionic analysis on {{Riemann}} surfaces and differential
  geometry.
\newblock In {\em Proceedings of the {{International Congress}} of
  {{Mathematicians}} ({{Berlin}}, 1998), Vol. II}, Documenta Mathematica, pages
  389--400, Berlin, 1998.

\bibitem[PS87]{Pinkall1987}
U.~Pinkall and I.~Sterling.
\newblock Willmore {{Surfaces}}.
\newblock {\em The Mathematical Intelligencer}, 9(2):38--43, January 1987.

\bibitem[Ric97]{Richter1997thesis}
J{\"o}rg Richter.
\newblock {\em Conformal Maps of a {{Riemannian}} Surface onto the Space of
  Quaternions}.
\newblock Doktors der {{Naturwissenschaften}} thesis, Technischen
  Universit{\"a}t Berlin, Berlin, 1997.

\bibitem[Riv08]{Ri}
Tristan Rivi{\`e}re.
\newblock Analysis aspects of {{Willmore}} surfaces.
\newblock {\em Inventiones mathematicae}, 174(1):1--45, 2008.

\bibitem[Riv14]{Ri2}
Tristan Rivi{\`e}re.
\newblock Variational principles for immersed surfaces with {{L2-bounded}}
  second fundamental form.
\newblock {\em Journal f{\"u}r die reine und angewandte Mathematik (Crelles
  Journal)}, 2014(695):41--98, 2014.

\bibitem[Rou14]{Rouviere2014}
Fran{\c c}ois Rouvi{\`e}re.
\newblock {\em Symmetric {{Spaces}} and the {{Kashiwara-Vergne Method}}},
  volume 2115 of {\em Lecture {{Notes}} in {{Mathematics}}}.
\newblock Springer International Publishing, Cham, 2014.

\bibitem[Roy88]{Ro2}
H.~L. Royden.
\newblock {\em Real Analysis}.
\newblock Macmillan, New York, 1988.

\bibitem[RS80]{RS1}
Michael Reed and Barry Simon.
\newblock {\em Methods of Modern Mathematical Physics. 1: {{Functional}}
  Analysis}.
\newblock Academic Press, New York, 1980.

\bibitem[RS00]{Reiter2000}
Hans Reiter and Jan~D Stegeman.
\newblock {\em Classical {{Harmonic Analysis}} and {{Locally Compact Groups}}}.
\newblock Oxford University Press, Oxford, 2000.

\bibitem[RS07]{RS2}
Michael Reed and Barry Simon.
\newblock {\em Methods of Modern Mathematical Physics. 2: {{Fourier}} Analysis,
  Self-Adjointness}.
\newblock Academic Press, San Diego, 2007.

\bibitem[Sch04]{Sch2}
Martin~Ulrich Schmidt.
\newblock Existence of minimizing willmore surfaces of prescribed conformal
  class, 2004.

\bibitem[Sha13]{Shafarevich}
Igor~R. Shafarevich.
\newblock {\em Basic {{Algebraic Geometry}} 1: {{Varieties}} in {{Projective
  Space}}}.
\newblock Springer, Berlin, 2013.

\bibitem[Sim86]{Si1}
Leon Simon.
\newblock Existence of {{Willmore}} surfaces.
\newblock In {\em Proceedings of the {{Centre}} for {{Mathematics}} and Its
  {{Applications}}}, volume~10, pages 187--216, Canberra, 1986.

\bibitem[Sim93]{Si2}
Leon Simon.
\newblock Existence of surfaces minimizing the {{Willmore}} functional.
\newblock {\em Communications in Analysis and Geometry}, 1(2):281--326, 1993.

\bibitem[Sog93]{So}
Christopher~D. Sogge.
\newblock {\em Fourier {{Integrals}} in {{Classical Analysis}}}.
\newblock Cambridge University Press, 1993.

\bibitem[Spi99]{Sp}
Michael Spivak.
\newblock {\em A Comprehensive Introduction to Differential Geometry},
  volume~IV.
\newblock Publish or Perish, Inc, Houston, TX, 1999.

\bibitem[Ste70]{St}
Elias~M. Stein.
\newblock {\em Singular Integrals and Differentiability Properties of
  Functions}.
\newblock Number~30 in Princeton Mathematical Series. Princeton University
  Press, Princeton, N.J, 1970.

\bibitem[SW72]{SW}
Elias~M. Stein and Guido Weiss.
\newblock {\em Introduction to {{Fourier Analysis}} on {{Euclidean Spaces}}}.
\newblock Princeton University Press, 1972.

\bibitem[Sze39]{Sz}
G{\'a}bor Szeg{\H o}.
\newblock {\em Orthogonal Polynomials}.
\newblock Number~23 in Colloquium Publications - {{American Mathematical
  Society}}. American Mathematical Society, Providence, 1939.

\bibitem[Tai97]{Ta1}
Iskander~A. Taimanov.
\newblock Modified {{Novikov--Veselov}} equation and differential geometry of
  surfaces.
\newblock In Victor Buchstaber and Serge Novikov, editors, {\em Solitons,
  {{Geometry}}, and {{Topology}}: {{On}} the {{Crossroad}}}, number 179 in
  American {{Mathematical Society Translations}}: {{Series}} 2, pages 133--151.
  American Mathematical Society, Providence, RI, 1997.

\bibitem[Tai98]{Ta2}
Iskander~A. Taimanov.
\newblock The {{Weierstrass}} representation of closed surfaces in
  {{$\mathbb{R}$3}}.
\newblock {\em Functional Analysis and Its Applications}, 32(4):258--267, 1998.

\bibitem[Tai06]{Ta3}
Iskander~A. Taimanov.
\newblock Two-dimensional {{Dirac}} operator and the theory of surfaces.
\newblock {\em Russian Mathematical Surveys}, 61(1):79--159, 2006.

\bibitem[Tho24]{Thomsen1924}
G.~Thomsen.
\newblock {Grundlagen der konformen fl{\"a}chentheorie}.
\newblock {\em Abhandlungen aus dem Mathematischen Seminar der Universit{\"a}t
  Hamburg}, 3(1):31--56, December 1924.

\bibitem[Wen69]{We}
Henry~C Wente.
\newblock An existence theorem for surfaces of constant mean curvature.
\newblock {\em Journal of Mathematical Analysis and Applications},
  26(2):318--344, 1969.

\bibitem[Whi72]{Wh}
Hassler Whitney.
\newblock {\em Complex Analytic Varieties}.
\newblock Addison-Wesley Publishing Co., Reading, Massachusetts, 1972.

\bibitem[Wil65]{Willmore1965}
Thomas~J Willmore.
\newblock Note on embedded surfaces.
\newblock {\em An. Sti. Univ.“Al. I. Cuza” Iasi Sect. I a Mat.(NS) B},
  11:493--496, 1965.

\bibitem[Wol93]{Wo}
Thomas~H. Wolff.
\newblock Recent work on sharp estimates in second-order elliptic unique
  continuation problems.
\newblock {\em Journal of Geometric Analysis}, 3(6):621--650, 1993.

\bibitem[Zie89]{Zi}
William~P. Ziemer.
\newblock {\em Weakly Differentiable Functions: {{Sobolev}} Spaces and
  Functions of Bounded Variation}.
\newblock Number 120 in Graduate Texts in Mathematics. Springer, New York,
  1989.

\end{thebibliography}

\newpage
\chapter*{List of Notation}


\noindent
\begin{tabularx}{\textwidth}{@{}l@{\quad}X@{}}
$\qi, \qj, \qk, \qat$ & Quaternions \\
$\Real, \Imag$ & Real and imaginary parts of a quaternion \\ 
$\X$ & Riemann surface \\
$z_l[\SO_l] = \Omega_l$ & Coordinates of $\X$ \\
$\mu(z)$ & Lebesgue measure in $z$-coordinate \\ 
$\SU, \SV$ & Open sets of $\X$ \\ 
$F : \X \to \qat$ & Admissible map \\
$N, R$ & Left and right normals of $F$ \\
$E$ & Holomorphic $\mathbb{C}$\=/line bundle \\
$E_\qat = E \otimes_\mathbb{C} \qat$ & Complex $\qat$\=/bundle \\
$K$ & Canonical bundle of $\X$ \\
$\pot{E}$ & $(0,1)$-Potentials of $E_\qat$ \\
$\pot{KE,g}$ & $(0,1)$-Potentials of $KE_\qat$ with  $\U=g\Bar{g}^{-1}\U\sd$ \\
$\willmore$ & Willmore energy \\

& \\
\multicolumn{2}{@{}l}{\textbf{Sheaves of sections}} \\
$\forms{p}{E}, \forms{p,q}{E}$ & Differential forms of degree $p$, bidegree $(p,q)$ valued in $E_\qat$ \\ 
$\ban{p}{E} = \sob{0,p}{E}$ & Lebesgue integrable sections  of $E_\qat$ \\ 
$\Sh{M}_{E,\V}$ & $\V$\=/meromorphic sections of $E_\qat$ \\
$\Q{E,\V}$ & $\V$\=/holomorphic sections  of $E_\qat$ \\ 
$\tQ{E,\V}$ & Pullback of $\Q{E,\V}$ to universal cover \\ 
$\sob{k,p}{E}$ & Sobolev sections  of $E_\qat$ \\ 
\multicolumn{2}{@{}l}{ 
$\Sh{V}=\Q{E,\V},\;
\Sh{U}=\Q{KE,\U},\;
\Sh{V}\sd=\Q{KE^{-1},\V\sd},\;
\Sh{U}\sd=\Q{E^{-1},\U\sd}$
} \\

& \\
\multicolumn{2}{@{}l}{\textbf{Operators and Pairings}} \\
$\barpartial_l, \partial_l$ & Wirtinger derivatives in $z_l$-coordinate \\ 
$\delbar{E}$ & Holomorphic structure of $E$ \\ 
$\dirac$ & Dirac operator \\ 
$\Op{I}_{\Omega}, \Op{I}_{\Omega,\V}$ & Right inverse of $\barpartial, \barpartial - \V$ on $\Omega$ \\
$\diracJ$ & Conjugation of $E_\qat$ \\ 
$\per \omega$ & Period map of $\omega$ \\ 
$\V = \V^+ + \V^-$ & Commuting and anti-commuting part \\ 
$\V\sd$ & Serre dual potential to $\V$ \\
$\lz \chi, \xi \rz$ & Pairing of $E_\qat$ and $E_\qat^{-1}$ \\ 
$\lp \chi, \psi \rp$ & Pairing of $E_\qat^{-1}$ and $KE_\qat$ valued in $\forms{1}{}$ \\ 
$\lh \xi, \omega \rh$ & Pairing of $KE^{-1}_\qat$ and $\forms{0,1}{E}$ valued in $\forms{2}{}$ \\ 
$\ls \upsilon\pa, \upsilon \rs$ & Pairing of $KE^{-1}_\qat$ and $E_\qat$ valued in $\pot{E}$
\end{tabularx}

\noindent
\begin{tabularx}{\textwidth}{@{}l@{\quad}X@{}}
\multicolumn{2}{@{}l}{\textbf{Sections used with a particular meaning}} \\
$\upsilon$ & Kodaira representation; Section of $\Q{E,\V}$ without roots \\
$\phi$ & Kodaira representation; Section of $\Q{E,\V}$ \\
\multicolumn{2}{@{}l}{ 
$d\Bar{z}\V - dz(\B + \U^-) = (d\upsilon)\upsilon^{-1},\; \U^+ = \V^+$
} \\
$\chi = \Bar{\upsilon}^{-1}$ & Weierstraß representation; Section of $\Q{E^{-1},\U\sd}$ \\
$\psi = (\partial + \B + \U^- )\phi$ & Weierstraß representation; Section of $\Q{KE,\U}$ \\
$\upsilon\pa$ & Section that can be $\ls\cdot,\cdot\rs$\=/paired with $\upsilon$ \\ 
$\var\V$ & Tangent vector to $\pot{E}$ \\ 
$\var\upsilon$ & Tangent vector to $H^0(\X,\Q{E,\V})$ \\

& \\
\multicolumn{2}{@{}l}{\textbf{Spaces}} \\
$\kodaira(\X,E)$ & Space of Kodaira representations \\ 
$\leftsp(\X,E)$ & Space of parameterizations of left normals \\ 
$\prewei(\X,E)$ & Space of Weierstraß representations without period conditions \\ 
$\wei(\X,E)$ & Space of Weierstraß representations \\ 
$\con(\X,\chi,\psi)$ & A space of admissible maps on the universal cover of $\X$ (Def~\ref{def:extended pairs}) \\ 
\end{tabularx}

\newcommand{\mypageout}[1]{\thepage}
\index{Index!Important terms|mypageout}
\begingroup
\renewcommand{\baselinestretch}{0.9}
\printindex
\endgroup

\end{document}